\pdfoutput=1

\newif\ifOPTcover
\OPTcovertrue

\newif\ifOPTbastard
\OPTbastardfalse

\def\OPTcovercolor{0.375,0,0,0.625} 
\def\OPTcovertextcolor{0,0,0,0.125}

\def\OPTlinkcolor{0,0.45,0}   

\def\OPTblue{blue}
\def\OPTred{red}
\def\OPTpurple{purple}

\def\OPTcolormodel{rgb}
\def\OPTcolxA{1., 0.85, 0.85}    
\def\OPTcolxB{1., 0.988462, 0.85}
\def\OPTcolxC{0.873077, 1., 0.85}
\def\OPTcolxD{0.85, 1., 0.965385}
\def\OPTcolxE{0.85, 0.896154, 1.}
\def\OPTcolxF{0.942308, 0.85, 1.}
\def\OPTcolxG{1., 0.85, 0.919231}
\def\OPTcolxH{1., 0.919231, 0.85}
\def\OPTcolxI{0.942308, 1., 0.85}
\def\OPTcolxJ{0.85, 1., 0.896154}
\def\OPTcolxK{0.85, 0.965385, 1.}
\def\OPTcolxL{0.873077, 0.85, 1.}
\def\OPTcolxM{1., 0.85, 0.988462}

\newcommand{\narrowbreak}{}
\newcommand{\narrowequation}[1]{$#1$}
\newenvironment{narrowmultline}{\csname equation\endcsname}{\csname endequation\endcsname}
\newenvironment{narrowmultline*}{\csname equation*\endcsname}{\csname endequation*\endcsname}

\def\OPTfontsize{11pt}        

\def\OPTpagesize{8.5in,11in}  
\def\OPTtopmargin{1in}        
\def\OPTbottommargin{1in}     
\def\OPTinnermargin{0.75in}   
\def\OPTbindingoffset{0.35in} 
\def\OPToutermargin{1.0in}    

\def\OPTcoverwidth{8.45in}    
\def\OPTcoverheight{10.95in}  
\def\OPTtopskip{43pt}         
\def\OPTbotskip{43pt}         
\def\OPTcovertitlefont{74.5}  
\def\OPTcovertitleskip{24pt}  
\def\OPTcoversubtitlefont{27} 
\def\OPTcoverauthorfont{22}   
\def\OPTcoverauthorskip{8pt}  

\def\OPTbacktitlefont{\LARGE} 
\def\OPTbackfont{\large}      

\def\OPTbastardtitlefont{22}
\def\OPTbastardsubtitlefont{15}
\def\OPTbastardwidth{0.40\textwidth} 
\def\OPTbastardtitleskip{8pt}

\def\OPTtitletitlefont{37}    
\def\OPTtitletitleskip{10pt}  
\def\OPTtitlesubtitlefont{25} 
\def\OPTtitlewidth{0.65\textwidth} 
\def\OPTtitleskip{30pt}       
\def\OPTtitleauthorfont{18}   
\def\OPTtitleauthorskip{6pt}  

\def\OPTpartfont{40}          
\def\OPTpartskip{20pt}        
\def\OPTparttitlefont{60}     

\def\OPTchapterfont{23}       
\def\OPTchapterskip{20pt}     
\def\OPTchaptertitlefont{35}  

\def\OPTprefacecols{3}        

\def\OPTwidow{\enlargethispage{\baselineskip}}

\def\OPTsmalltable{}

\newcommand{\OPTspherescolwidth}{30pt}

\newcommand{\OPTbibliographyfont}{\small}

\newcommand{\OPTindexfont}{\small}      
\newcommand{\OPTindexcolumnsep}{20.0pt} 

\def\OPThalftorus{torus-lores-bw}

\let\OPTfrontimage=\OPTlofrontimage
\let\OPTbackimage=\OPTlobackimage

\documentclass[\OPTfontsize]{book}

\PassOptionsToPackage{table}{xcolor}

\usepackage{etex} 

\usepackage[papersize={\OPTpagesize},
            twoside,
            includehead,
            top=\OPTtopmargin,
            bottom=\OPTbottommargin,
            inner=\OPTinnermargin,
            outer=\OPToutermargin,
            bindingoffset=\OPTbindingoffset]{geometry}

\usepackage[backref=page,
            colorlinks,
            citecolor=linkcolor,
            linkcolor=linkcolor,
            urlcolor=linkcolor,
            unicode,
            pdfauthor={Univalent Foundations Program},
            pdftitle={Homotopy Type Theory: Univalent Foundations of Mathematics},
            pdfsubject={Mathematics},
            pdfkeywords={type theory, homotopy theory, univalence axiom}]{hyperref}
\renewcommand{\backref}[1]{}
\renewcommand{\backrefalt}[4]{%
   \ifcase #1 %
   (No citations.)
   \or
   (Cited on page\ #2.)
   \else
   (Cited on pages\ #2.)
   \fi}

\usepackage{layout}

\usepackage[utf8]{inputenc}

\usepackage{pifont}

\usepackage{multicol}

\usepackage{mathpazo}
\usepackage[scaled=0.95]{helvet}
\usepackage{courier}
\usepackage{graphicx}
\usepackage{comment}

\usepackage{fancyhdr} 

\usepackage{nextpage} 

\usepackage{amssymb,amsmath,amsthm,stmaryrd,mathrsfs,wasysym}
\usepackage{enumitem,mathtools,xspace}
\usepackage{xstring} 

\usepackage{xcolor} 
\usepackage{wallpaper} 

\usepackage{booktabs} 
\usepackage{array} 

\definecolor{linkcolor}{rgb}{\OPTlinkcolor}
\usepackage{aliascnt}
\usepackage[capitalize]{cleveref}
\usepackage[all,2cell,cmtip]{xy}
\UseAllTwocells
\usepackage{braket} 

\usepackage{tikz}
\usetikzlibrary{decorations.pathmorphing,arrows}

\usepackage{etoolbox}           

\usepackage{mathpartir}         

\usepackage[numbered]{bookmark} 

\makeatletter
\newcommand{\sectionNotes}{\phantomsection\section*{Notes}\addcontentsline{toc}{section}{Notes}\markright{\textsc{\@chapapp{} \thechapter{} Notes}}}
\newcommand{\sectionExercises}[1]{\phantomsection\section*{Exercises}\addcontentsline{toc}{section}{Exercises}\markright{\textsc{\@chapapp{} \thechapter{} Exercises}}}
\makeatother

\newcommand{\jdeq}{\equiv}      
\let\judgeq\jdeq
\newcommand{\defeq}{\vcentcolon\equiv}  

\newcommand{\define}[1]{\textbf{#1}}

\newcommand{\Fin}{\ensuremath{\mathsf{Fin}}}
\newcommand{\fmax}{\ensuremath{\mathsf{fmax}}}

\def\prdsym{\textstyle\prod}
\makeatletter
\def\prd#1{\@ifnextchar\bgroup{\prd@parens{#1}}{\@ifnextchar\sm{\prd@parens{#1}\@eatsm}{\prd@noparens{#1}}}}
\def\prd@parens#1{\@ifnextchar\bgroup%
  {\mathchoice{\@dprd{#1}}{\@tprd{#1}}{\@tprd{#1}}{\@tprd{#1}}\prd@parens}%
  {\@ifnextchar\sm%
    {\mathchoice{\@dprd{#1}}{\@tprd{#1}}{\@tprd{#1}}{\@tprd{#1}}\@eatsm}%
    {\mathchoice{\@dprd{#1}}{\@tprd{#1}}{\@tprd{#1}}{\@tprd{#1}}}}}
\def\@eatsm\sm{\sm@parens}
\def\prd@noparens#1{\mathchoice{\@dprd@noparens{#1}}{\@tprd{#1}}{\@tprd{#1}}{\@tprd{#1}}}
\def\lprd#1{\@ifnextchar\bgroup{\@lprd{#1}\lprd}{\@@lprd{#1}}}
\def\@lprd#1{\mathchoice{{\textstyle\prod}}{\prod}{\prod}{\prod}({\textstyle #1})\;}
\def\@@lprd#1{\mathchoice{{\textstyle\prod}}{\prod}{\prod}{\prod}({\textstyle #1}),\ }
\def\tprd#1{\@tprd{#1}\@ifnextchar\bgroup{\tprd}{}}
\def\@tprd#1{\mathchoice{{\textstyle\prod_{(#1)}}}{\prod_{(#1)}}{\prod_{(#1)}}{\prod_{(#1)}}}
\def\dprd#1{\@dprd{#1}\@ifnextchar\bgroup{\dprd}{}}
\def\@dprd#1{\prod_{(#1)}\,}
\def\@dprd@noparens#1{\prod_{#1}\,}

\def\lam#1{{\lambda}\@lamarg#1:\@endlamarg\@ifnextchar\bgroup{.\,\lam}{.\,}}
\def\@lamarg#1:#2\@endlamarg{\if\relax\detokenize{#2}\relax #1\else\@lamvar{\@lameatcolon#2},#1\@endlamvar\fi}
\def\@lamvar#1,#2\@endlamvar{(#2\,{:}\,#1)}
\def\@lameatcolon#1:{#1}
\let\lamt\lam
\def\lamu#1{{\lambda}\@lamuarg#1:\@endlamuarg\@ifnextchar\bgroup{.\,\lamu}{.\,}}
\def\@lamuarg#1:#2\@endlamuarg{#1}

\def\fall#1{\forall (#1)\@ifnextchar\bgroup{.\,\fall}{.\,}}

\def\exis#1{\exists (#1)\@ifnextchar\bgroup{.\,\exis}{.\,}}

\def\smsym{\textstyle\sum}
\def\sm#1{\@ifnextchar\bgroup{\sm@parens{#1}}{\@ifnextchar\prd{\sm@parens{#1}\@eatprd}{\sm@noparens{#1}}}}
\def\sm@parens#1{\@ifnextchar\bgroup%
  {\mathchoice{\@dsm{#1}}{\@tsm{#1}}{\@tsm{#1}}{\@tsm{#1}}\sm@parens}%
  {\@ifnextchar\prd%
    {\mathchoice{\@dsm{#1}}{\@tsm{#1}}{\@tsm{#1}}{\@tsm{#1}}\@eatprd}%
    {\mathchoice{\@dsm{#1}}{\@tsm{#1}}{\@tsm{#1}}{\@tsm{#1}}}}}
\def\@eatprd\prd{\prd@parens}
\def\sm@noparens#1{\mathchoice{\@dsm@noparens{#1}}{\@tsm{#1}}{\@tsm{#1}}{\@tsm{#1}}}
\def\lsm#1{\@ifnextchar\bgroup{\@lsm{#1}\lsm}{\@@lsm{#1}}}
\def\@lsm#1{\mathchoice{{\textstyle\sum}}{\sum}{\sum}{\sum}({\textstyle #1})\;}
\def\@@lsm#1{\mathchoice{{\textstyle\sum}}{\sum}{\sum}{\sum}({\textstyle #1}),\ }
\def\tsm#1{\@tsm{#1}\@ifnextchar\bgroup{\tsm}{}}
\def\@tsm#1{\mathchoice{{\textstyle\sum_{(#1)}}}{\sum_{(#1)}}{\sum_{(#1)}}{\sum_{(#1)}}}
\def\dsm#1{\@dsm{#1}\@ifnextchar\bgroup{\dsm}{}}
\def\@dsm#1{\sum_{(#1)}\,}
\def\@dsm@noparens#1{\sum_{#1}\,}

\def\wtypesym{{\mathsf{W}}}
\def\wtype#1{\@ifnextchar\bgroup%
  {\mathchoice{\@twtype{#1}}{\@twtype{#1}}{\@twtype{#1}}{\@twtype{#1}}\wtype}%
  {\mathchoice{\@twtype{#1}}{\@twtype{#1}}{\@twtype{#1}}{\@twtype{#1}}}}
\def\lwtype#1{\@ifnextchar\bgroup{\@lwtype{#1}\lwtype}{\@@lwtype{#1}}}
\def\@lwtype#1{\mathchoice{{\textstyle\mathsf{W}}}{\mathsf{W}}{\mathsf{W}}{\mathsf{W}}({\textstyle #1})\;}
\def\@@lwtype#1{\mathchoice{{\textstyle\mathsf{W}}}{\mathsf{W}}{\mathsf{W}}{\mathsf{W}}({\textstyle #1}),\ }
\def\twtype#1{\@twtype{#1}\@ifnextchar\bgroup{\twtype}{}}
\def\@twtype#1{\mathchoice{{\textstyle\mathsf{W}_{(#1)}}}{\mathsf{W}_{(#1)}}{\mathsf{W}_{(#1)}}{\mathsf{W}_{(#1)}}}
\def\dwtype#1{\@dwtype{#1}\@ifnextchar\bgroup{\dwtype}{}}
\def\@dwtype#1{\mathsf{W}_{(#1)}\,}

\newcommand{\suppsym}{{\mathsf{sup}}}
\newcommand{\supp}{\ensuremath\suppsym\xspace}

\def\wtypeh#1{\@ifnextchar\bgroup%
  {\mathchoice{\@lwtypeh{#1}}{\@twtypeh{#1}}{\@twtypeh{#1}}{\@twtypeh{#1}}\wtypeh}%
  {\mathchoice{\@@lwtypeh{#1}}{\@twtypeh{#1}}{\@twtypeh{#1}}{\@twtypeh{#1}}}}
\def\lwtypeh#1{\@ifnextchar\bgroup{\@lwtypeh{#1}\lwtypeh}{\@@lwtypeh{#1}}}
\def\@lwtypeh#1{\mathchoice{{\textstyle\mathsf{W}^h}}{\mathsf{W}^h}{\mathsf{W}^h}{\mathsf{W}^h}({\textstyle #1})\;}
\def\@@lwtypeh#1{\mathchoice{{\textstyle\mathsf{W}^h}}{\mathsf{W}^h}{\mathsf{W}^h}{\mathsf{W}^h}({\textstyle #1}),\ }
\def\twtypeh#1{\@twtypeh{#1}\@ifnextchar\bgroup{\twtypeh}{}}
\def\@twtypeh#1{\mathchoice{{\textstyle\mathsf{W}^h_{(#1)}}}{\mathsf{W}^h_{(#1)}}{\mathsf{W}^h_{(#1)}}{\mathsf{W}^h_{(#1)}}}
\def\dwtypeh#1{\@dwtypeh{#1}\@ifnextchar\bgroup{\dwtypeh}{}}
\def\@dwtypeh#1{\mathsf{W}^h_{(#1)}\,}

\makeatother

\let\setof\Set    

\newcommand{\tup}[2]{(#1,#2)}
\newcommand{\proj}[1]{\ensuremath{\mathsf{pr}_{#1}}\xspace}
\newcommand{\fst}{\ensuremath{\proj1}\xspace}
\newcommand{\snd}{\ensuremath{\proj2}\xspace}
\newcommand{\ac}{\ensuremath{\mathsf{ac}}\xspace} 
\newcommand{\un}{\ensuremath{\mathsf{upun}}\xspace} 

\newcommand{\rec}[1]{\mathsf{rec}_{#1}}
\newcommand{\ind}[1]{\mathsf{ind}_{#1}}
\newcommand{\indid}[1]{\ind{=_{#1}}} 
\newcommand{\indidb}[1]{\ind{=_{#1}}'} 

\newcommand{\uppt}{\ensuremath{\mathsf{uppt}}\xspace}

\newcommand{\pairpath}{\ensuremath{\mathsf{pair}^{\mathord{=}}}\xspace}
\newcommand{\projpath}[1]{\ensuremath{\apfunc{\proj{#1}}}\xspace}

\newcommand{\pairr}[1]{{\mathopen{}(#1)\mathclose{}}}
\newcommand{\Pairr}[1]{{\mathopen{}\left(#1\right)\mathclose{}}}

\newcommand{\im}{\ensuremath{\mathsf{im}}} 

\newcommand{\leftwhisker}{\mathbin{{\ct}_{\ell}}}
\newcommand{\rightwhisker}{\mathbin{{\ct}_{r}}}
\newcommand{\hct}{\star}

\newcommand{\modal}{\ensuremath{\ocircle}}
\let\reflect\modal
\newcommand{\modaltype}{\ensuremath{\type_\modal}}
\newcommand{\mreturn}{\ensuremath{\eta}}
\let\project\mreturn
\newcommand{\ext}{\mathsf{ext}}
\renewcommand{\P}{\ensuremath{\type_{P}}\xspace}

\newcommand{\idsym}{{=}}
\newcommand{\id}[3][]{\ensuremath{#2 =_{#1} #3}\xspace}
\newcommand{\idtype}[3][]{\ensuremath{\mathsf{Id}_{#1}(#2,#3)}\xspace}
\newcommand{\idtypevar}[1]{\ensuremath{\mathsf{Id}_{#1}}\xspace}
\newcommand{\defid}{\coloneqq}

\newcommand{\dpath}[4]{#3 =^{#1}_{#2} #4}

\newcommand{\refl}[1]{\ensuremath{\mathsf{refl}_{#1}}\xspace}

\newcommand{\ct}{%
  \mathchoice{\mathbin{\raisebox{0.5ex}{$\displaystyle\centerdot$}}}%
             {\mathbin{\raisebox{0.5ex}{$\centerdot$}}}%
             {\mathbin{\raisebox{0.25ex}{$\scriptstyle\,\centerdot\,$}}}%
             {\mathbin{\raisebox{0.1ex}{$\scriptscriptstyle\,\centerdot\,$}}}
}

\newcommand{\opp}[1]{\mathord{{#1}^{-1}}}
\let\rev\opp

\newcommand{\trans}[2]{\ensuremath{{#1}_{*}\mathopen{}\left({#2}\right)\mathclose{}}\xspace}

\newcommand{\transf}[1]{\ensuremath{{#1}_{*}}\xspace} 
\newcommand{\transfib}[3]{\ensuremath{\mathsf{transport}^{#1}(#2,#3)\xspace}}
\newcommand{\Transfib}[3]{\ensuremath{\mathsf{transport}^{#1}\Big(#2,\, #3\Big)\xspace}}
\newcommand{\transfibf}[1]{\ensuremath{\mathsf{transport}^{#1}\xspace}}

\newcommand{\transtwo}[2]{\ensuremath{\mathsf{transport}^2\mathopen{}\left({#1},{#2}\right)\mathclose{}}\xspace}

\newcommand{\transconst}[3]{\ensuremath{\mathsf{transportconst}}^{#1}_{#2}(#3)\xspace}
\newcommand{\transconstf}{\ensuremath{\mathsf{transportconst}}\xspace}

\newcommand{\mapfunc}[1]{\ensuremath{\mathsf{ap}_{#1}}\xspace} 
\newcommand{\map}[2]{\ensuremath{{#1}\mathopen{}\left({#2}\right)\mathclose{}}\xspace}

\newcommand{\mapdepfunc}[1]{\ensuremath{\mathsf{apd}_{#1}}\xspace} 

\let\apfunc\mapfunc
\let\ap\map
\let\apdfunc\mapdepfunc
\let\apd\mapdep

\newcommand{\aptwofunc}[1]{\ensuremath{\mathsf{ap}^2_{#1}}\xspace}
\newcommand{\aptwo}[2]{\ensuremath{\aptwofunc{#1}\mathopen{}\left({#2}\right)\mathclose{}}\xspace}
\newcommand{\apdtwofunc}[1]{\ensuremath{\mathsf{apd}^2_{#1}}\xspace}
\newcommand{\apdtwo}[2]{\ensuremath{\apdtwofunc{#1}\mathopen{}\left(#2\right)\mathclose{}}\xspace}

\newcommand{\idfunc}[1][]{\ensuremath{\mathsf{id}_{#1}}\xspace}

\newcommand{\htpy}{\sim}

\newcommand{\bisim}{\sim}       
\newcommand{\eqr}{\sim}         

\newcommand{\eqv}[2]{\ensuremath{#1 \simeq #2}\xspace}
\newcommand{\eqvspaced}[2]{\ensuremath{#1 \;\simeq\; #2}\xspace}
\newcommand{\eqvsym}{\simeq}    
\newcommand{\texteqv}[2]{\ensuremath{\mathsf{Equiv}(#1,#2)}\xspace}
\newcommand{\isequiv}{\ensuremath{\mathsf{isequiv}}}
\newcommand{\qinv}{\ensuremath{\mathsf{qinv}}}
\newcommand{\ishae}{\ensuremath{\mathsf{ishae}}}
\newcommand{\linv}{\ensuremath{\mathsf{linv}}}
\newcommand{\rinv}{\ensuremath{\mathsf{rinv}}}
\newcommand{\biinv}{\ensuremath{\mathsf{biinv}}}
\newcommand{\lcoh}[3]{\mathsf{lcoh}_{#1}(#2,#3)}
\newcommand{\rcoh}[3]{\mathsf{rcoh}_{#1}(#2,#3)}
\newcommand{\hfib}[2]{{\mathsf{fib}}_{#1}(#2)}

\newcommand{\total}[1]{\ensuremath{\mathsf{total}(#1)}}

\newcommand{\UU}{\ensuremath{\mathcal{U}}\xspace}
\let\bbU\UU
\let\type\UU
\newcommand{\typele}[1]{\ensuremath{{#1}\text-\mathsf{Type}}\xspace}
\newcommand{\typeleU}[1]{\ensuremath{{#1}\text-\mathsf{Type}_\UU}\xspace}
\newcommand{\typelep}[1]{\ensuremath{{(#1)}\text-\mathsf{Type}}\xspace}

\let\ntype\typele
\let\ntypeU\typeleU

\renewcommand{\set}{\ensuremath{\mathsf{Set}}\xspace}
\newcommand{\setU}{\ensuremath{\mathsf{Set}_\UU}\xspace}
\newcommand{\prop}{\ensuremath{\mathsf{Prop}}\xspace}
\newcommand{\propU}{\ensuremath{\mathsf{Prop}_\UU}\xspace}
\newcommand{\pointed}[1]{\ensuremath{#1_\bullet}}

\newcommand{\card}{\ensuremath{\mathsf{Card}}\xspace}
\newcommand{\ord}{\ensuremath{\mathsf{Ord}}\xspace}
\newcommand{\ordsl}[2]{{#1}_{/#2}}

\newcommand{\ua}{\ensuremath{\mathsf{ua}}\xspace} 
\newcommand{\idtoeqv}{\ensuremath{\mathsf{idtoeqv}}\xspace}
\newcommand{\univalence}{\ensuremath{\mathsf{univalence}}\xspace} 

\newcommand{\iscontr}{\ensuremath{\mathsf{isContr}}}
\newcommand{\contr}{\ensuremath{\mathsf{contr}}} 
\newcommand{\isset}{\ensuremath{\mathsf{isSet}}}
\newcommand{\isprop}{\ensuremath{\mathsf{isProp}}}

\let\hfiber\hfib

\newcommand{\trunc}[2]{\mathopen{}\left\Vert #2\right\Vert_{#1}\mathclose{}}
\newcommand{\ttrunc}[2]{\bigl\Vert #2\bigr\Vert_{#1}}
\newcommand{\Trunc}[2]{\Bigl\Vert #2\Bigr\Vert_{#1}}
\newcommand{\truncf}[1]{\Vert \blank \Vert_{#1}}
\newcommand{\tproj}[3][]{\mathopen{}\left|#3\right|_{#2}^{#1}\mathclose{}}
\newcommand{\tprojf}[2][]{|\blank|_{#2}^{#1}}
\def\pizero{\trunc0}

\newcommand{\brck}[1]{\trunc{}{#1}}
\newcommand{\bbrck}[1]{\ttrunc{}{#1}}
\newcommand{\Brck}[1]{\Trunc{}{#1}}
\newcommand{\bproj}[1]{\tproj{}{#1}}
\newcommand{\bprojf}{\tprojf{}}

\newcommand{\Parens}[1]{\Bigl(#1\Bigr)}

\let\extendsmb\ext
\newcommand{\extend}[1]{\extendsmb(#1)}

\newcommand{\emptyt}{\ensuremath{\mathbf{0}}\xspace}

\newcommand{\unit}{\ensuremath{\mathbf{1}}\xspace}
\newcommand{\ttt}{\ensuremath{\star}\xspace}

\newcommand{\bool}{\ensuremath{\mathbf{2}}\xspace}
\newcommand{\btrue}{{1_{\bool}}}
\newcommand{\bfalse}{{0_{\bool}}}

\newcommand{\inlsym}{{\mathsf{inl}}}
\newcommand{\inrsym}{{\mathsf{inr}}}
\newcommand{\inl}{\ensuremath\inlsym\xspace}
\newcommand{\inr}{\ensuremath\inrsym\xspace}

\newcommand{\seg}{\ensuremath{\mathsf{seg}}\xspace}

\newcommand{\freegroup}[1]{F(#1)}
\newcommand{\freegroupx}[1]{F'(#1)} 

\newcommand{\glue}{\mathsf{glue}}

\newcommand{\Sn}{\mathbb{S}}
\newcommand{\base}{\ensuremath{\mathsf{base}}\xspace}
\newcommand{\lloop}{\ensuremath{\mathsf{loop}}\xspace}
\newcommand{\surf}{\ensuremath{\mathsf{surf}}\xspace}

\newcommand{\susp}{\Sigma}
\newcommand{\north}{\mathsf{N}}
\newcommand{\south}{\mathsf{S}}
\newcommand{\merid}{\mathsf{merid}}

\newcommand{\blank}{\mathord{\hspace{1pt}\text{--}\hspace{1pt}}}

\newcommand{\nameless}{\mathord{\hspace{1pt}\underline{\hspace{1ex}}\hspace{1pt}}}

\newcommand{\uset}{\ensuremath{\mathcal{S}et}\xspace}
\newcommand{\ucat}{\ensuremath{{\mathcal{C}at}}\xspace}
\newcommand{\urel}{\ensuremath{\mathcal{R}el}\xspace}
\newcommand{\uhilb}{\ensuremath{\mathcal{H}ilb}\xspace}
\newcommand{\utype}{\ensuremath{\mathcal{T}\!ype}\xspace}

\newbox\pbbox
\setbox\pbbox=\hbox{\xy \POS(65,0)\ar@{-} (0,0) \ar@{-} (65,65)\endxy}

\newcommand{\inv}[1]{{#1}^{-1}}
\newcommand{\idtoiso}{\ensuremath{\mathsf{idtoiso}}\xspace}
\newcommand{\isotoid}{\ensuremath{\mathsf{isotoid}}\xspace}
\newcommand{\op}{^{\mathrm{op}}}
\newcommand{\y}{\ensuremath{\mathbf{y}}\xspace}
\newcommand{\dgr}[1]{{#1}^{\dagger}}
\newcommand{\unitaryiso}{\mathrel{\cong^\dagger}}
\newcommand{\cteqv}[2]{\ensuremath{#1 \simeq #2}\xspace}

\newcommand{\N}{\ensuremath{\mathbb{N}}\xspace}
\let\nat\N
\newcommand{\natp}{\ensuremath{\nat'}\xspace} 

\newcommand{\zerop}{\ensuremath{0'}\xspace}   
\newcommand{\suc}{\mathsf{succ}}
\newcommand{\sucp}{\ensuremath{\suc'}\xspace} 
\newcommand{\add}{\mathsf{add}}
\newcommand{\ack}{\mathsf{ack}}
\newcommand{\ite}{\mathsf{iter}}
\newcommand{\assoc}{\mathsf{assoc}}
\newcommand{\dbl}{\ensuremath{\mathsf{double}}}
\newcommand{\dblp}{\ensuremath{\dbl'}\xspace} 

\newcommand{\lst}[1]{\mathsf{List}(#1)}
\newcommand{\nil}{\mathsf{nil}}
\newcommand{\cons}{\mathsf{cons}}

\newcommand{\vect}[2]{\ensuremath{\mathsf{Vec}_{#1}(#2)}\xspace}

\newcommand{\Z}{\ensuremath{\mathbb{Z}}\xspace}
\newcommand{\Zsuc}{\mathsf{succ}}
\newcommand{\Zpred}{\mathsf{pred}}

\newcommand{\Q}{\ensuremath{\mathbb{Q}}\xspace}

\newcommand{\funext}{\mathsf{funext}}
\newcommand{\happly}{\mathsf{happly}}

\newcommand{\code}{\ensuremath{\mathsf{code}}\xspace}
\newcommand{\encode}{\ensuremath{\mathsf{encode}}\xspace}
\newcommand{\decode}{\ensuremath{\mathsf{decode}}\xspace}

\newcommand{\function}[4]{\left\{\begin{array}{rcl}#1 &
      \longrightarrow & #2 \\ #3 & \longmapsto & #4 \end{array}\right.}

\newcommand{\cocone}[2]{\mathsf{cocone}_{#1}(#2)}
\newcommand{\composecocone}[2]{#1\circ#2}

\newcommand{\Ddiag}{\mathscr{D}}

\newcommand{\Map}{\mathsf{Map}}

\newcommand{\interval}{\ensuremath{I}\xspace}
\newcommand{\izero}{\ensuremath{0_{\interval}}\xspace}
\newcommand{\ione}{\ensuremath{1_{\interval}}\xspace}

\newcommand{\epi}{\ensuremath{\twoheadrightarrow}}
\newcommand{\mono}{\ensuremath{\rightarrowtail}}

\newcommand{\bin}{\ensuremath{\mathrel{\widetilde{\in}}}}

\newcommand{\semigroupstrsym}{\ensuremath{\mathsf{SemigroupStr}}}
\newcommand{\semigroupstr}[1]{\ensuremath{\mathsf{SemigroupStr}}(#1)}
\newcommand{\semigroup}[0]{\ensuremath{\mathsf{Semigroup}}}

\newcommand{\emptyctx}{\ensuremath{\cdot}}
\newcommand{\production}{\vcentcolon\vcentcolon=}
\newcommand{\conv}{\downarrow}
\newcommand{\ctx}{\ensuremath{\mathsf{ctx}}}
\newcommand{\wfctx}[1]{#1\ \ctx}
\newcommand{\oftp}[3]{#1 \vdash #2 : #3}
\newcommand{\jdeqtp}[4]{#1 \vdash #2 \jdeq #3 : #4}

\newcommand{\tmtp}[2]{#1 \mathord{:} #2}

\newcommand{\form}{\textsc{form}}
\newcommand{\intro}{\textsc{intro}}
\newcommand{\elim}{\textsc{elim}}
\newcommand{\comp}{\textsc{comp}}
\newcommand{\uniq}{\textsc{uniq}}
\newcommand{\Weak}{\mathsf{Wkg}}
\newcommand{\Vble}{\mathsf{Vble}}

\newcommand{\Subst}{\mathsf{Subst}}

\newcommand{\cc}{\mathsf{c}}
\newcommand{\pp}{\mathsf{p}}
\newcommand{\cct}{\widetilde{\mathsf{c}}}
\newcommand{\ppt}{\widetilde{\mathsf{p}}}
\newcommand{\Wtil}{\ensuremath{\widetilde{W}}\xspace}

\newcommand{\istype}[1]{\mathsf{is}\mbox{-}{#1}\mbox{-}\mathsf{type}}
\newcommand{\nplusone}{\ensuremath{(n+1)}}
\newcommand{\nminusone}{\ensuremath{(n-1)}}
\newcommand{\fact}{\mathsf{fact}}

\newcommand{\kbar}{\overline{k}} 

\newcommand{\natw}{\ensuremath{\mathbf{N^w}}\xspace}
\newcommand{\zerow}{\ensuremath{0^\mathbf{w}}\xspace}
\newcommand{\sucw}{\ensuremath{\mathbf{s^w}}\xspace}
\newcommand{\nalg}{\nat\mathsf{Alg}}
\newcommand{\nhom}{\nat\mathsf{Hom}}
\newcommand{\ishinitw}{\mathsf{isHinit}_{\mathsf{W}}}
\newcommand{\ishinitn}{\mathsf{isHinit}_\nat}
\newcommand{\w}{\mathsf{W}}
\newcommand{\walg}{\w\mathsf{Alg}}
\newcommand{\whom}{\w\mathsf{Hom}}

\newcommand{\RC}{\ensuremath{\mathbb{R}_\mathsf{c}}\xspace} 
\newcommand{\RD}{\ensuremath{\mathbb{R}_\mathsf{d}}\xspace} 
\newcommand{\R}{\ensuremath{\mathbb{R}}\xspace}           
\newcommand{\barRD}{\ensuremath{\bar{\mathbb{R}}_\mathsf{d}}\xspace} 

\newcommand{\close}[1]{\sim_{#1}} 
\newcommand{\closesym}{\mathord\sim}
\newcommand{\rclim}{\mathsf{lim}} 
\newcommand{\rcrat}{\mathsf{rat}} 
\newcommand{\rceq}{\mathsf{eq}_{\RC}} 
\newcommand{\CAP}{\mathcal{C}}    
\newcommand{\Qp}{\Q_{+}}
\newcommand{\apart}{\mathrel{\#}}  
\newcommand{\dcut}{\mathsf{isCut}}  
\newcommand{\cover}{\triangleleft} 
\newcommand{\intfam}[3]{(#2, \lam{#1} #3)} 

\newcommand{\bsim}{\frown}
\newcommand{\bbsim}{\smile}

\newcommand{\hapx}{\diamondsuit\approx}
\newcommand{\hapname}{\diamondsuit}
\newcommand{\hapxb}{\heartsuit\approx}
\newcommand{\hapbname}{\heartsuit}
\newcommand{\tap}[1]{\bullet\approx_{#1}\triangle}
\newcommand{\tapname}{\triangle}
\newcommand{\tapb}[1]{\bullet\approx_{#1}\square}
\newcommand{\tapbname}{\square}

\newcommand{\NO}{\ensuremath{\mathsf{No}}\xspace}
\newcommand{\surr}[2]{\{\,#1\,\big|\,#2\,\}}
\newcommand{\LL}{\mathcal{L}}
\newcommand{\RR}{\mathcal{R}}
\newcommand{\noeq}{\mathsf{eq}_{\NO}} 

\newcommand{\ble}{\trianglelefteqslant}
\newcommand{\blt}{\vartriangleleft}
\newcommand{\bble}{\sqsubseteq}
\newcommand{\bblt}{\sqsubset}

\newcommand{\hle}{\diamondsuit\preceq}
\newcommand{\hlt}{\diamondsuit\prec}
\newcommand{\hlname}{\diamondsuit}
\newcommand{\hleb}{\heartsuit\preceq}
\newcommand{\hltb}{\heartsuit\prec}
\newcommand{\hlbname}{\heartsuit}
\newcommand{\tle}{\triangle\preceq}
\newcommand{\tlt}{\triangle\prec}
\newcommand{\tlname}{\triangle}
\newcommand{\tleb}{\square\preceq}
\newcommand{\tltb}{\square\prec}
\newcommand{\tlbname}{\square}

\newcommand{\vset}{\mathsf{set}}  
\def\cd{\tproj0}
\newcommand{\inj}{\ensuremath{\mathsf{inj}}} 
\newcommand{\acc}{\ensuremath{\mathsf{acc}}} 

\newcommand{\power}[1]{\mathcal{P}(#1)} 
\newcommand{\powerp}[1]{\mathcal{P}_+(#1)} 

\makeatletter
\def\defthm#1#2#3{%
  \newaliascnt{#1}{thm}
  \newtheorem{#1}[#1]{#2}
  \aliascntresetthe{#1}
  \crefname{#1}{#2}{#3}}

\newtheorem{thm}{Theorem}[section]
\crefname{thm}{Theorem}{Theorems}
\defthm{cor}{Corollary}{Corollaries}
\defthm{lem}{Lemma}{Lemmas}
\defthm{axiom}{Axiom}{Axioms}
\theoremstyle{definition}
\defthm{defn}{Definition}{Definitions}
\theoremstyle{remark}
\defthm{rmk}{Remark}{Remarks}
\defthm{eg}{Example}{Examples}
\defthm{egs}{Examples}{Examples}
\defthm{notes}{Notes}{Notes}
\newtheorem{ex}{Exercise}[chapter]
\crefname{ex}{Exercise}{Exercises}

\crefformat{section}{\S#2#1#3}
\Crefformat{section}{Section~#2#1#3}
\crefrangeformat{section}{\S\S#3#1#4--#5#2#6}
\Crefrangeformat{section}{Sections~#3#1#4--#5#2#6}
\crefmultiformat{section}{\S\S#2#1#3}{ and~#2#1#3}{, #2#1#3}{ and~#2#1#3}
\Crefmultiformat{section}{Sections~#2#1#3}{ and~#2#1#3}{, #2#1#3}{ and~#2#1#3}
\crefrangemultiformat{section}{\S\S#3#1#4--#5#2#6}{ and~#3#1#4--#5#2#6}{, #3#1#4--#5#2#6}{ and~#3#1#4--#5#2#6}
\Crefrangemultiformat{section}{Sections~#3#1#4--#5#2#6}{ and~#3#1#4--#5#2#6}{, #3#1#4--#5#2#6}{ and~#3#1#4--#5#2#6}

\crefformat{appendix}{Appendix~#2#1#3}
\Crefformat{appendix}{Appendix~#2#1#3}
\crefrangeformat{appendix}{Appendices~#3#1#4--#5#2#6}
\Crefrangeformat{appendix}{Appendices~#3#1#4--#5#2#6}
\crefmultiformat{appendix}{Appendices~#2#1#3}{ and~#2#1#3}{, #2#1#3}{ and~#2#1#3}
\Crefmultiformat{appendix}{Appendices~#2#1#3}{ and~#2#1#3}{, #2#1#3}{ and~#2#1#3}
\crefrangemultiformat{appendix}{Appendices~#3#1#4--#5#2#6}{ and~#3#1#4--#5#2#6}{, #3#1#4--#5#2#6}{ and~#3#1#4--#5#2#6}
\Crefrangemultiformat{appendix}{Appendices~#3#1#4--#5#2#6}{ and~#3#1#4--#5#2#6}{, #3#1#4--#5#2#6}{ and~#3#1#4--#5#2#6}

\crefname{part}{Part}{Parts}

\crefname{figure}{Figure}{Figures}

\let\autoref\cref

\let\c@equation\c@thm
\numberwithin{equation}{section}

\setitemize[1]{leftmargin=2em}
\setenumerate[1]{leftmargin=*}

\setitemize[1]{itemsep=-0.2em}
\setenumerate[1]{itemsep=-0.2em}

\def\noteson{%
\gdef\note##1{\mbox{}\marginpar{\color{blue}\textasteriskcentered\ ##1}}}

\noteson

\newcommand{\Coq}{\textsc{Coq}\xspace}
\newcommand{\Agda}{\textsc{Agda}\xspace}
\newcommand{\NuPRL}{\textsc{NuPRL}\xspace}

\newcommand{\indexdef}[1]{\index{#1|defstyle}}   
\newcommand{\indexfoot}[1]{\index{#1|footstyle}} 
\newcommand{\indexsee}[2]{\index{#1|see{#2}}}    

\newcommand{\ZF}{Zermelo--Fraenkel}
\newcommand{\CZF}{Constructive \ZF{} Set Theory}

\newcommand{\LEM}[1]{\ensuremath{\mathsf{LEM}_{#1}}\xspace}
\newcommand{\choice}[1]{\ensuremath{\mathsf{AC}_{#1}}\xspace}

\newcommand{\mentalpause}{\medskip} 

\newcounter{symindex}
\newcommand{\symlabel}[1]{\refstepcounter{symindex}\label{#1}}

\usepackage{makeidx}
\makeindex

\renewcommand{\chaptermark}[1]{\markboth{\textsc{Chapter \thechapter. #1}}{}}
\renewcommand{\sectionmark}[1]{\markright{\textsc{\thesection\ #1}}}

\usepackage[raggedright]{titlesec}
\titleformat{\part}[display]{\fontsize{\OPTpartfont}{\OPTpartfont}\fontseries{m}\fontshape{sc}\selectfont}{\hfil\partname\ \Roman{part}}{\OPTpartskip}{\fontsize{\OPTparttitlefont}{\OPTparttitlefont}\fontseries{b}\fontshape{sc}\selectfont\hfil}
\titleformat{\chapter}[display]{\fontsize{\OPTchapterfont}{\OPTchapterfont}\fontseries{m}\fontshape{it}\selectfont}{\chaptertitlename\ \thechapter}{\OPTchapterskip}{\fontsize{\OPTchaptertitlefont}{\OPTchaptertitlefont}\fontseries{b}\fontshape{n}\selectfont}

\definecolor{covercolor}{cmyk}{\OPTcovercolor}
\definecolor{covertext}{cmyk}{\OPTcovertextcolor}
\pagecolor{white}
\nopagecolor

\begin{document}

\title{Homotopy Type Theory: Univalent Foundations of Mathematics}
\author{The Univalent Foundations Program}

\frontmatter 

\pagestyle{empty}


\newgeometry{noheadfoot,bindingoffset=-5pt,top=0pt,bottom=0pt,inner=0pt,outer=0pt}%
\ifOPTcover
\setcounter{page}{-1} 
\newlength{\coverheight}
\setlength{\coverheight}{\OPTcoverheight}
\newlength{\coverwidth}
\setlength{\coverwidth}{\OPTcoverwidth}
\newcommand{\frontpage}{
\begin{minipage}[b][\coverheight][c]{\coverwidth}
\hbox{}\hfill
\begin{minipage}[b][\coverheight][t]{0.83\coverwidth}
\color{covertext}
\vspace{\OPTtopskip}
{\fontsize{\OPTcovertitlefont}{\OPTcovertitlefont}\fontseries{b}\selectfont%
  \hfill Homotopy\par \hfill Type Theory}\par
\vspace*{\OPTcovertitleskip}
{\fontsize{\OPTcoversubtitlefont}{\OPTcoversubtitlefont}\fontshape{it}\selectfont
\hfill Univalent Foundations of Mathematics}

\vfill

{\fontsize{\OPTcoverauthorfont}{\OPTcoverauthorfont}\fontseries{b}\fontshape{sc}\selectfont

\hfill The Univalent Foundations Program
\par
\vspace*{\OPTcoverauthorskip}
\hfill Institute for Advanced Study\par

\vspace*{\OPTbotskip}
}

\end{minipage}\hfill\hbox{}
\end{minipage}
}

\ThisLLCornerWallPaper{1.1}{\OPTfrontimage}
\pagecolor{covercolor}
\frontpage
\newpage
\nopagecolor
\cleartooddpage
\else
\fi

\ifOPTbastard
\cleartooddpage
\hbox{}
\vspace{0.2\textwidth}
{\centering
\makebox[\OPTbastardwidth][s]{
\fontsize{\OPTbastardtitlefont}{\OPTbastardtitlefont}\fontshape{n}\selectfont%
\textbf{Homotopy Type Theory}}\par
\vspace*{\OPTbastardtitleskip}
\makebox[\OPTbastardwidth][s]{
\fontsize{\OPTbastardsubtitlefont}{\OPTbastardsubtitlefont}\fontshape{n}\selectfont%
\textit{Univalent Foundations of Mathematics}}\par
}
\else
\fi

\cleartooddpage
\hbox{}\vfill
{\centering
\makebox[\OPTtitlewidth][s]{\fontsize{\OPTtitletitlefont}{\OPTtitletitlefont}\fontseries{b}\selectfont%
Homotopy Type Theory}\par
\vspace*{\OPTtitletitleskip}
\makebox[\OPTtitlewidth][s]{\fontsize{\OPTtitlesubtitlefont}{\OPTtitlesubtitlefont}\fontshape{it}\selectfont%
Univalent Foundations of Mathematics}\par
\vspace*{\OPTtitleskip}
{\fontsize{\OPTtitleauthorfont}{\OPTtitleauthorfont}\fontshape{n}\selectfont%
The Univalent Foundations Program\par
\vspace*{\OPTtitleauthorskip}
Institute for Advanced Study\par
}
\vspace*{\OPTtitleskip}
\vspace*{\OPTtitleskip}
\includegraphics[width=\OPTtitlewidth]{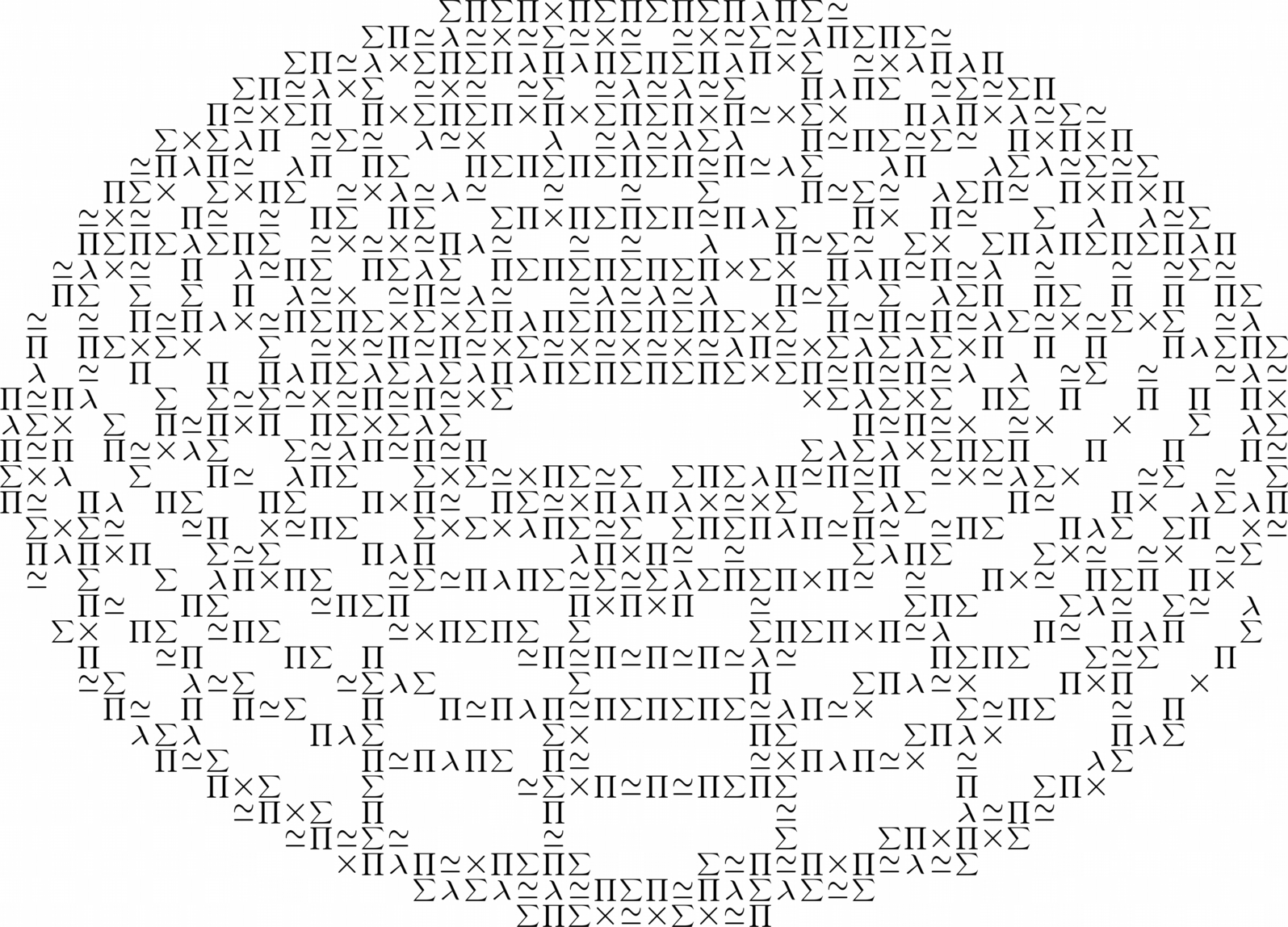}\par
}

\vfill
\hbox{}

\clearpage
\restoregeometry

\hbox{}
\vfill
\newcommand{\OPTversion}{first-edition-257-g5561b73
}
{\small
\noindent
\emph{``Homotopy Type Theory: Univalent Foundations of Mathematics''}\\
\copyright\ 2013 The Univalent Foundations Program

\medskip
\noindent
Book version: \texttt{\OPTversion}

\medskip
\noindent
MSC 2010 classification:
\texttt{03-02},
\texttt{55-02},
\texttt{03B15}

\bigskip
\footnotesize

\noindent
This work is licensed under the
\textbf{\emph{Creative Commons Attribution-ShareAlike 3.0 Unported License.}}
To view a copy of this license, visit
\href{http://creativecommons.org/licenses/by-sa/3.0/}{http://creativecommons.org/licenses/by-sa/3.0/}.

\bigskip

\noindent
This book is freely available at \href{http://homotopytypetheory.org/book/}{http://homotopytypetheory.org/book/}.

\bigskip

\noindent
\emph{\textbf{\small Acknowledgment}}

\medskip

\noindent
Apart from the generous support from the Institute for Advanced Study, some contributors
to the book were partially or fully supported by the following agencies and grants:
\begin{itemize}
\item Association of Members of the Institute for Advanced Study: a grant to the Institute for Advanced Study 
\item Agencija za raziskovalno dejavnost Republike Slovenije:  
\href{http://www.sicris.si/search/prg.aspx?id=6120}{P1--0294},
\href{http://www.sicris.si/search/prj.aspx?id=7109}{N1--0011}.

\item Air Force Office of Scientific Research:
  FA9550-11-1-0143, and 
  FA9550-12-1-0370.  
  {
    \setlength{\parskip}{0pt}
    \begin{quote}
      \noindent\scriptsize
      This material is based in part upon work supported by the AFOSR under the above awards.
      Any opinions, findings, and conclusions or recommendations expressed in this publication are those of the author(s) and do not necessarily reflect the views of the AFOSR.
    \end{quote}
  }

\item Engineering and Physical Sciences Research Council: 
   \href{http://gow.epsrc.ac.uk/NGBOViewGrant.aspx?GrantRef=EP/G034109/1}{EP/G034109/1}, 
   \href{http://gow.epsrc.ac.uk/NGBOViewGrant.aspx?GrantRef=EP/G03298X/1}{EP/G03298X/1}. 

\item European Union's 7th Framework Programme under grant agreement nr.\ 243847 (%
\href{http://wiki.portal.chalmers.se/cse/pmwiki.php/ForMath/ForMath/}{ForMath}). 

\item National Science Foundation:
  \href{http://www.nsf.gov/awardsearch/showAward.do?AwardNumber=1001191}{DMS-1001191}, 
  \href{http://www.nsf.gov/awardsearch/showAward.do?AwardNumber=1100938}{DMS-1100938}, 
  \href{http://www.nsf.gov/awardsearch/showAward.do?AwardNumber=1116703}{CCF-1116703}, 
  and 
  \href{http://www.nsf.gov/awardsearch/showAward.do?AwardNumber=1128155}{DMS-1128155}. 
  {
    \setlength{\itemsep}{0pt}
    \begin{quote}
      \noindent\scriptsize
      This material is based in part upon work supported by the
      National Science Foundation under the above awards.  Any opinions,
      findings, and conclusions or recommendations expressed in this
      material are those of the author(s) and do not necessarily reflect the
      views of the National Science Foundation.
    \end{quote}
  }
\item The Simonyi Fund: a grant to the Institute for Advanced Study          
  \end{itemize}

}
\cleartooddpage


\cleartooddpage

\hypertarget{preface}{}
\bookmark[dest=preface]{Preface}

\pagestyle{fancyplain}

\chapter*{Preface}
\label{cha:preface}

\subsection*{IAS Special Year on Univalent Foundations}

A Special Year on Univalent Foundations of Mathematics was held in 2012-13 at the Institute for Advanced Study, School of Mathematics, organized by Steve Awodey, Thierry Coquand, and Vladimir Voevodsky.  The following people were the official participants.

\begin{multicols}{\OPTprefacecols}{
\begin{itemize}
\item[] Peter Aczel
\item[] Benedikt Ahrens
\item[] Thorsten Altenkirch
\item[] Steve Awodey
\item[] Bruno Barras
\item[] Andrej Bauer
\item[] Yves Bertot
\item[] Marc Bezem
\item[] Thierry Coquand
\item[] Eric Finster
\item[] Daniel Grayson
\item[] Hugo Herbelin
\item[] Andr\'e Joyal
\item[] Dan Licata
\item[] Peter Lumsdaine
\item[] Assia Mahboubi
\item[] Per Martin-L\"of
\item[] Sergey Melikhov
\item[] Alvaro Pelayo
\item[] Andrew Polonsky
\item[] Michael Shulman
\item[] Matthieu Sozeau
\item[] Bas Spitters
\item[] Benno van den Berg
\item[] Vladimir Voevodsky
\item[] Michael Warren
\item[] Noam Zeilberger
\end{itemize}
}
\end{multicols}

\noindent There were also the following students, whose participation was no less valuable.

\begin{multicols}{\OPTprefacecols}{
\begin{itemize}
\item[] Carlo Angiuli
\item[] Anthony Bordg
\item[] Guillaume Brunerie
\item[] Chris Kapulkin
\item[] Egbert Rijke
\item[] Kristina Sojakova
\end{itemize}
}
\end{multicols}

\noindent In addition, there were the following short- and long-term visitors, including student visitors, whose contributions to the Special Year were also essential.

\begin{multicols}{\OPTprefacecols}{
\begin{itemize}
\item[] Jeremy Avigad
\item[] Cyril Cohen
\item[] Robert Constable
\item[] Pierre-Louis Curien
\item[] Peter Dybjer
\item[] Mart{\'\i}n Escard{\'o}
\item[] Kuen-Bang Hou
\item[] Nicola Gambino
\item[] Richard Garner
\item[] Georges Gonthier
\item[] Thomas Hales
\item[] Robert Harper
\item[] Martin Hofmann
\item[] Pieter Hofstra
\item[] Joachim Kock
\item[] Nicolai Kraus
\item[] Nuo Li
\item[] Zhaohui Luo
\item[] Michael Nahas
\item[] Erik Palmgren
\item[] Emily Riehl
\item[] Dana Scott
\item[] Philip Scott
\item[] Sergei Soloviev
\end{itemize}
}
\end{multicols}

\subsection*{About this book}

We did not set out to write a book. The present work has its origins in our collective attempts to develop a new style of ``informal type theory'' that can be read and understood by a human being, as a complement to a formal proof that can be checked by a machine.  
Univalent foundations is closely tied to the idea of a foundation of mathematics that can be implemented in a computer proof assistant.  Although such a formalization is not part of this book, much of the material presented here was actually done first in the fully formalized setting inside a proof assistant, and only later ``unformalized'' to arrive at the presentation you find before you --- a remarkable inversion of the usual state of affairs in formalized mathematics.  

Each of the above-named individuals contributed something to the Special Year --- and so to this book --- in the form of ideas, words, or deeds.  The spirit of collaboration that prevailed throughout the year was truly extraordinary. 

\mentalpause

Special thanks are due to the Institute for Advanced Study, without which this book would obviously never have come to be.  It proved  to be an ideal setting for the creation of this new branch of mathematics: stimulating, congenial, and supportive.  May some trace of this unique atmosphere linger in the pages of this book, and in the future development of this new field of study.

\bigskip

\begin{flushright}
The Univalent Foundations Program\\
Institute for Advanced Study\\
Princeton, April 2013
\end{flushright}


\cleartooddpage[\thispagestyle{empty}]

\hypertarget{toc}{}
\bookmark[dest=toc]{Table of Contents}

\setcounter{tocdepth}{1}        
\tableofcontents
\setcounter{tocdepth}{2}        
\cleartooddpage[\thispagestyle{empty}]

\mainmatter 

\pagestyle{fancyplain}

\chapter*{Introduction}
\markboth{\textsc{Introduction}}{}
\addcontentsline{toc}{chapter}{Introduction}
\setcounter{page}{1}
\pagenumbering{arabic}

\emph{Homotopy type theory} is a new branch of mathematics that combines aspects of several different fields in a surprising way. It is based on a recently discovered connection between \emph{homotopy theory} and \emph{type theory}.
Homotopy theory is an outgrowth of algebraic topology and homological algebra, with relationships to higher category theory; while type theory is a branch of mathematical logic and theoretical computer science.
Although the connections between the two are currently the focus of intense investigation, it is increasingly clear that they are just the beginning of a subject that will take more time and more hard work to fully understand.
It touches on topics as seemingly distant as the homotopy groups of spheres, the algorithms for type checking, and the definition of weak $\infty$-groupoids.

Homotopy type theory also brings new ideas into the very foundation of mathematics.
\index{foundations, univalent}%
On the one hand, there is Voevodsky's subtle and beautiful \emph{univalence axiom}. 
\index{univalence axiom}%
The univalence axiom implies, in particular, that isomorphic structures can be identified, a principle that mathematicians have been happily using on workdays, despite its incompatibility with the ``official'' doctrines of conventional foundations.
On the other hand, we have \emph{higher inductive types}, which provide direct, logical descriptions of some of the basic spaces and constructions of homotopy theory: spheres, cylinders, truncations, localizations, etc.
Both ideas are impossible to capture directly in classical set-theoretic foundations, but when combined in homotopy type theory, they permit an entirely new kind of ``logic of homotopy types''.
\index{foundations}%

This suggests a new conception of foundations of mathematics, with intrinsic homotopical content, an ``invariant'' conception of the objects of mathematics --- and convenient machine implementations, which can serve as a practical aid to the working mathematician.
This is the \emph{Univalent Foundations} program.
The present book is intended as a first systematic exposition of the basics of univalent foundations, and a collection of examples of this new style of reasoning --- but without requiring the reader to know or learn any formal logic, or to use any computer proof assistant.

\OPTwidow

We emphasize that homotopy type theory is a young field, and univalent foundations is very much a work in progress.
This book should be regarded as a ``snapshot'' of the state of the field at the time it was written, rather than a polished exposition of an established edifice.
As we will discuss briefly later, there are many aspects of homotopy type theory that are not yet fully understood --- but as of this writing, its broad outlines seem clear enough.
The ultimate theory will probably not look exactly like the one described in this book, but it will surely be \emph{at least} as capable and powerful; we therefore believe that univalent foundations will eventually become a viable alternative to set theory as the ``implicit foundation'' for the unformalized mathematics done by most mathematicians.

\subsection*{Type theory}

Type theory was originally invented by Bertrand Russell \cite{Russell:1908},\index{Russell, Bertrand} as a device for blocking the paradoxes in the logical foundations of mathematics  that were under investigation at the time. It was later developed as a rigorous formal system in its own right (under the name ``$\lambda$-calculus'') by Alonzo Church \cite{Church:1933cl,Church:1940tu,Church:1941tc}.  Although it is not generally regarded as the foundation for classical mathematics, set theory being more customary, type theory still has numerous applications, especially in computer science and the theory of programming languages~\cite{Pierce-TAPL}.
\index{programming}%
\index{type theory}%
\index{lambda-calculus@$\lambda$-calculus}%
Per Martin-L\"{o}f \cite{Martin-Lof-1972,Martin-Lof-1973,Martin-Lof-1979,martin-lof:bibliopolis}, among others,
developed a ``predicative'' modification of Church's type system, which is now usually called dependent, constructive, intuitionistic, or simply \emph{Martin\--L\"of type theory}. This is the basis of the system that we consider here; it was originally intended as a rigorous framework for the formalization of constructive mathematics.  In what follows, we will often use ``type theory'' to refer specifically to this system and similar ones, although type theory as a subject is much broader (see \cite{somma,kamar} for the history of type theory).

In type theory, unlike set theory, objects are classified using a primitive notion of \emph{type}, similar to the data-types used in programming languages.  These elaborately structured types can be used to express detailed specifications of the objects classified, giving rise to principles of reasoning about these objects.  To take a very simple example, the objects of a product type $A\times B$ are known to be of the form $\pairr{a,b}$, and so one automatically knows how to construct them and how to decompose them. Similarly, an object of function type $A\to B$ can be acquired from an object of type $B$ parametrized by objects of type $A$, and can be evaluated at an argument of type $A$.  This rigidly predictable behavior of all objects (as opposed to set theory's more liberal formation principles, allowing inhomogeneous sets) is one aspect of type theory that has led to its extensive use in verifying the correctness of computer programs.  The clear reasoning principles associated with the construction of types also form the basis of modern \emph{computer proof assistants},%
\index{proof!assistant}%
\indexsee{computer proof assistant}{proof assistant}
\index{mathematics!formalized}%
which are used for formalizing mathematics and verifying the correctness of formalized proofs.  We return to this aspect of type theory below.  

One problem in understanding type theory from a mathematical point of view, however, has always been that the basic concept of \emph{type} is unlike that of \emph{set} in ways that have been hard to make precise.  We believe that the new idea of regarding types, not as strange sets (perhaps constructed without using classical logic), but as spaces, viewed from the perspective of homotopy theory, is a significant step forward.  In particular, it solves the problem of understanding how the notion of equality of elements of a type differs from that of elements of a set.

In homotopy theory one is concerned with spaces
\index{topological!space}%
and continuous mappings between them, 
\index{function!continuous!in classical homotopy theory}%
up to homotopy.  A \emph{homotopy}
\index{homotopy!topological}%
between a pair of continuous maps $f : X \to Y$
and  $g : X\to Y$ is 
a continuous map $H : X \times [0, 1] \to Y$ satisfying
$H(x, 0) = f (x)$  and $H(x, 1) = g(x)$. The homotopy $H$ may be thought of as a ``continuous deformation'' of $f$ into $g$. The spaces $X$ and $Y$ are said to be \emph{homotopy equivalent},
\index{homotopy!equivalence!topological}%
$\eqv X Y$, if there are continuous maps going back and forth, the composites of which are homotopical to the respective identity mappings, i.e., if they are isomorphic ``up to homotopy''.  Homotopy equivalent spaces have the same algebraic invariants (e.g., homology, or the fundamental group), and are said to have the same \emph{homotopy type}.

\subsection*{Homotopy type theory}

Homotopy type theory (HoTT) interprets type theory from a homotopical perspective.
In homotopy type theory, we regard the types as ``spaces'' (as studied in homotopy theory) or higher groupoids, and the logical constructions (such as the product $A\times B$) as homotopy-invariant constructions on these spaces.
In this way, we are able to manipulate spaces directly without first having to develop point-set topology (or any combinatorial replacement for it, such as the theory of simplicial sets).
To briefly explain this perspective, consider first the basic concept of type theory, namely that
the \emph{term} $a$ is of \emph{type} $A$, which is written:
\[ a:A. \]
This expression is traditionally thought of as akin to:
\begin{center}
``$a$ is an element of the set $A$.''
\end{center}
However, in homotopy type theory we think of it instead as:
\begin{center}
``$a$ is a point of the space $A$.''
\end{center}
\index{continuity of functions in type theory@``continuity'' of functions in type theory}%
Similarly, every function $f : A\to B$ in type theory is regarded as a continuous map from the space $A$ to the space $B$.

We should stress that these ``spaces'' are treated purely homotopically, not topologically.
For instance, there is no notion of ``open subset'' of a type or of ``convergence'' of a sequence of elements of a type.
We only have ``homotopical'' notions, such as paths between points and homotopies between paths, which also make sense in other models of homotopy theory (such as simplicial sets).
Thus, it would be more accurate to say that we treat types as \emph{$\infty$-groupoids}\index{.infinity-groupoid@$\infty$-groupoid}; this is a name for the ``invariant objects'' of homotopy theory which can be presented by topological spaces,
\index{topological!space}%
simplicial sets, or any other model for homotopy theory.
However, it is convenient to sometimes use topological words such as ``space'' and ``path'', as long as we remember that other topological concepts are not applicable.

(It is tempting to also use the phrase \emph{homotopy type}
\index{homotopy!type}%
for these objects, suggesting the dual interpretation of ``a type (as in type theory) viewed homotopically'' and ``a space considered from the point of view of homotopy theory.''
The latter is a bit different from the classical meaning of ``homotopy type'' as an \emph{equivalence class} of spaces modulo homotopy equivalence, although it does preserve the meaning of phrases such as ``these two spaces have the same homotopy type''.)

The idea of interpreting types as structured objects, rather than sets, has a long pedigree, and is known to clarify various mysterious aspects of type theory.
For instance, interpreting types as sheaves helps explain the intuitionistic nature of type-theoretic logic, while interpreting them as partial equivalence relations or ``domains'' helps explain its computational aspects.
It also implies that we can use type-theoretic reasoning to study the structured objects, leading to the rich field of categorical logic.
The homotopical interpretation fits this same pattern: it clarifies the nature of \emph{identity} (or equality) in type theory, and allows us to use type-theoretic reasoning in the study of homotopy theory.

The key new idea of the homotopy interpretation is that the logical notion of identity $a = b$ of two objects $a, b: A$ of the same type $A$ can be understood as the existence of a path $p : a \leadsto b$ from point $a$ to point $b$ in the space $A$.
This also means that two functions $f, g: A\to B$ can be identified if they are homotopic, since a homotopy is just a (continuous) family of paths $p_x: f(x) \leadsto g(x)$ in $B$, one for each $x:A$.
In type theory, for every type $A$ there is a (formerly somewhat mysterious) type $\idtypevar{A}$ of identifications of two objects of $A$; in homotopy type theory, this is just the \emph{path space} $A^I$ of all continuous maps $I\to A$ from the unit interval.
\index{unit!interval}%
\index{interval!topological unit}%
\index{path!topological}%
\index{topological!path}%
In this way, a term $p : \idtype[A]{a}{b}$ represents a path $p : a \leadsto b$ in $A$. 

The idea of homotopy type theory arose around 2006 in independent work by Awodey and Warren~\cite{AW} and Voevodsky~\cite{VV}, but it was inspired by 
Hofmann and Streicher's earlier groupoid interpretation~\cite{hs:gpd-typethy}.
Indeed, higher-dimensional category theory (particularly the theory of weak $\infty$-groupoids) is now known to be intimately connected to homotopy theory, as proposed by Grothendieck and now being studied intensely by mathematicians of both sorts.
The original semantic models of Awodey--Warren and Voevodsky use well-known notions and techniques from homotopy theory which are now also in use in higher category theory, such as  Quillen model categories and Kan\index{Kan complex} simplicial sets\index{simplicial!sets}.
\index{Quillen model category}%
\index{model category}%

Voevodsky recognized that the simplicial interpretation of type theory satisfies a further crucial property, dubbed \emph{univalence}, which had not previously been considered in type theory (although Church's principle of extensionality for propositions turns out to be a very special case of it).
Adding univalence to type theory in the form of a new axiom has far-reaching consequences, many of which are natural, simplifying and compelling.
The univalence axiom also further strengthens the homotopical view of type theory, since it holds in the simplicial model and other related models, while failing under the view of types as sets.

\subsection*{Univalent foundations}

Very briefly, the basic idea of the univalence axiom can be explained as follows.
In type theory, one can have a universe $\UU$, the terms of which are themselves types, $A : \UU$, etc.
Those types that are terms of $\UU$ are commonly called \emph{small} types.
\index{type!small}%
\index{small!type}%
Like any type, $\UU$ has an identity type $\idtypevar{\UU}$, which expresses the identity relation $A = B$ between small types.
Thinking of types as spaces, $\UU$ is a space, the points of which are spaces; to understand its identity type, we must ask, what is a path $p : A \leadsto B$ between spaces in $\UU$?
The univalence axiom says that such paths correspond to homotopy equivalences $\eqv A B$, (roughly) as explained above.
A bit more precisely, given any (small) types $A$ and $B$, in addition to the primitive type $\idtype[\UU]AB$ of identifications of $A$ with $B$, there is the defined type $\texteqv AB$ of equivalences from $A$ to $B$.
Since the identity map on any object is an equivalence, there is a canonical map,
\[\idtype[\UU]AB\to\texteqv AB.\]
The univalence axiom states that this map is itself an equivalence.
At the risk of oversimplifying, we can state this succinctly as follows:

\begin{description}\index{univalence axiom}%
\item[Univalence Axiom:]  $\eqvspaced{(A = B)}{(\eqv A B)}$.
\end{description}
In other words, identity is equivalent to equivalence. \index{identity}%
In particular, one may say that ``equivalent types are identical''.
However, this phrase is somewhat misleading, since it may sound like a sort of ``skeletality'' condition which \emph{collapses} the notion of equivalence to coincide with identity, whereas in fact univalence is about \emph{expanding} the notion of identity so as to coincide with the (unchanged) notion of equivalence.

From the homotopical point of view, univalence implies that spaces of the same homotopy type are connected by a path in the universe $\UU$, in accord with the intuition of a classifying space for (small) spaces.
From the logical point of view, however, it is a radically new idea: it says that isomorphic things can be identified!  Mathematicians are of course used to identifying isomorphic structures in practice, but they generally do so by ``abuse of notation''\index{abuse!of notation}, or some other informal device, knowing that the objects involved are not ``really'' identical.  But in this new foundational scheme, such structures can be formally identified, in the logical sense that every property or construction involving one also applies to the other. Indeed, the identification is now made explicit, and properties and constructions can be systematically transported along it.  Moreover, the different ways in which such identifications may be made themselves form a structure that one can (and should!)\ take into account.

Thus in sum, for points $A$ and $B$ of the universe $\UU$ (i.e., small types), the univalence axiom identifies the following three notions:
\begin{itemize}
\item (logical) an identification $p:A=B$ of $A$ and $B$
\item (topological) a path $p:A \leadsto B$ from $A$ to $B$ in $\UU$
\item (homotopical) an equivalence $p:\eqv A B$ between $A$ and $B$.
\end{itemize}

\subsection*{Higher inductive types}\index{type!higher inductive}%

One of the classical advantages of type theory is its simple and effective techniques for working with inductively defined structures.
The simplest nontrivial inductively defined structure is the natural numbers, which is inductively generated by zero and the successor function.
From this statement one can algorithmically\index{algorithm} extract the principle of mathematical induction, which characterizes the natural numbers.
More general inductive definitions encompass lists and well-founded trees of all sorts, each of which is characterized by a corresponding ``induction principle''.
This includes most data structures used in certain programming languages; hence the usefulness of type theory in formal reasoning about the latter.
If conceived in a very general sense, inductive definitions also include examples such as a disjoint union $A+B$, which may be regarded as ``inductively'' generated by the two injections $A\to A+B$ and $B\to A+B$.
The ``induction principle'' in this case is ``proof by case analysis'', which characterizes the disjoint union.

In homotopy theory, it is natural to consider also ``inductively defined spaces'' which are generated not merely by a collection of \emph{points}, but also by collections of \emph{paths} and higher paths.
Classically, such spaces are called \emph{CW complexes}.
\index{CW complex}%
For instance, the circle $S^1$ is generated by a single point and a single path from that point to itself.
Similarly, the 2-sphere $S^2$ is generated by a single point $b$ and a single two-dimensional path from the constant path at $b$ to itself, while the torus $T^2$ is generated by a single point, two paths $p$ and $q$ from that point to itself, and a two-dimensional path from $p\ct q$ to $q\ct p$.

By using the identification of paths with identities in homotopy type theory, these sort of ``inductively defined spaces'' can be characterized in type theory by ``induction principles'', entirely analogously to classical examples such as the natural numbers and the disjoint union.
The resulting \emph{higher inductive types}
\index{type!higher inductive}%
give a direct ``logical'' way to reason about familiar spaces such as spheres, which (in combination with univalence) can be used to perform familiar arguments from homotopy theory, such as calculating homotopy groups of spheres, in a purely formal way.
The resulting proofs are a marriage of classical homotopy-theoretic ideas with classical type-theoretic ones, yielding new insight into both disciplines.

Moreover, this is only the tip of the iceberg: many abstract constructions from homotopy theory, such as homotopy colimits, suspensions, Postnikov towers, localization, completion, and spectrification, can also be expressed as higher inductive types.
Many of these are classically constructed using Quillen's ``small object argument'', which can be regarded as a finite way of algorithmically describing an infinite CW complex presentation\index{presentation!of a space as a CW complex} of a space, just as ``zero and successor'' is a finite algorithmic\index{algorithm} description of the infinite set of natural numbers.
Spaces produced by the small object argument are infamously complicated and difficult to understand; the type-theoretic approach is potentially much simpler, bypassing the need for any explicit construction by giving direct access to the appropriate ``induction principle''.
Thus, the combination of univalence and higher inductive types suggests the possibility of a revolution, of sorts, in the practice of homotopy theory.

\subsection*{Sets in univalent foundations}

\index{set|(}%

We have claimed that univalent foundations can eventually serve as a foundation for ``all'' of mathematics, but so far we have discussed 
only homotopy theory.  Of course, there are many specific examples of the use of type theory without the new homotopy type theory features to formalize mathematics,
\index{mathematics!formalized}%
\index{theorem!Feit--Thompson}%
\index{theorem!odd-order}%
\index{Feit--Thompson theorem}%
\index{odd-order theorem}%
such as the recent formalization of the Feit--Thompson odd-order theorem in \Coq~\cite{gonthier}.

But the traditional view is that mathematics is founded on set theory, in the sense that all mathematical objects and constructions can be coded into a theory such as Zermelo--Fraenkel set theory (ZF).
\index{set theory!Zermelo--Fraenkel}%
\indexsee{Zermelo-Fraenkel set theory}{set theory}%
\indexsee{ZF}{set theory}%
\indexsee{ZFC}{set theory}%
However, it is well-established by now that for most mathematics outside of set theory proper, the intricate hierarchical membership structure of sets in ZF is really unnecessary: a more ``structural'' theory, such as Lawvere's\index{Lawvere} Elementary Theory of the Category of Sets~\cite{lawvere:etcs-long}, suffices.
\index{Elementary Theory of the Category of Sets}%

In univalent foundations, the basic objects are ``homotopy types'' rather than sets, but we can \emph{define} a class of types which behave like sets.
Homotopically, these can be thought of as spaces in which every connected component is contractible, i.e.\ those which are homotopy equivalent to a discrete space.
\index{discrete!space}%
It is a theorem  that the category of such ``sets'' satisfies Lawvere's\index{Lawvere} axioms (or related ones, depending on the details of the theory).
Thus, any sort of mathematics that can be represented in an ETCS-like theory (which, experience suggests, is essentially all of mathematics) can equally well be represented in univalent foundations.  

This supports the claim that univalent foundations is at least as good as existing foundations of mathematics.
A mathematician working in univalent foundations can build structures out of sets in a familiar way, with more general homotopy types waiting in the foundational background until there is need of them.
For this reason, most of the applications in this book have been chosen to be areas where univalent foundations has something \emph{new} to contribute that distinguishes it from existing foundational systems.

Unsurprisingly, homotopy theory and category theory are two of these, but perhaps less obvious is that univalent foundations has something new and interesting to offer even in subjects such as set theory and real analysis.
For instance, the univalence axiom allows us to identify isomorphic structures, while higher inductive types allow direct descriptions of objects by their universal properties.
Thus we can generally avoid resorting to arbitrarily chosen representatives or transfinite iterative constructions.
In fact, even the objects of study in ZF set theory can be characterized, inside the sets of univalent foundations, by such an inductive universal property.

\index{set|)}%

\subsection*{Informal type theory}

\index{mathematics!formalized|(defstyle}%
\index{informal type theory|(defstyle}%
\index{type theory!informal|(defstyle}%
\index{type theory!formal|(}%
One difficulty often encountered by the classical mathematician when faced with learning about type theory is that it is usually presented as a fully or partially formalized deductive system.
This style, which is very useful for proof-theoretic investigations, is not particularly convenient for use in applied, informal reasoning.
Nor is it even familiar to most working mathematicians, even those who might be interested in foundations of mathematics.
One objective of the present work is to develop an informal style of doing mathematics in univalent foundations that is at once rigorous and precise, but is also closer to the language and style of presentation of everyday mathematics.

In present-day mathematics, one usually constructs and reasons about mathematical objects in a way that could in principle, one presumes, be formalized in a system of elementary set theory, such as ZFC --- at least given enough ingenuity and patience.
For the most part, one does not even need to be aware of this possibility, since it largely coincides with the condition that a proof be ``fully rigorous'' (in the sense that all mathematicians have come to understand intuitively through education and experience).
But one does need to learn to be careful about a few aspects of ``informal set theory'': the use of collections too large or inchoate to be sets; the axiom of choice and its equivalents; even (for undergraduates) the method of proof by contradiction; and so on.
Adopting a new foundational system such as homotopy type theory as the \emph{implicit formal basis} of informal reasoning will require adjusting some of one's instincts and practices.
The present text is intended to serve as an example of this ``new kind of mathematics'', which is still informal, but could now in principle be formalized in homotopy type theory, rather than ZFC, again given enough ingenuity and patience.

It is worth emphasizing that, in this new system, such formalization can have real practical benefits.
The formal system of type theory is suited to computer systems and has been implemented in existing proof assistants.
\index{proof!assistant}%
A proof assistant is a computer program which guides the user in construction of a fully formal proof, only allowing valid steps of reasoning.
It also provides some degree of automation, can search libraries for existing theorems, and can even extract numerical algorithms\index{algorithm} \index{extraction of algorithms} from the resulting (constructive) proofs.

We believe that this aspect of the univalent foundations program distinguishes it from other approaches to foundations, potentially providing a new practical utility for the working mathematician.
Indeed, proof assistants based on older type theories have already been used to formalize substantial mathematical proofs, such as the four-color theorem\index{theorem!four-color} \index{four-color theorem} and the Feit--Thompson theorem.
Computer implementations of univalent foundations are presently works in progress (like the theory itself).
\index{proof!assistant}%
However, even its currently available implementations (which are mostly small modifications to existing proof assistants such as \Coq and 
\Agda) have already demonstrated their worth, not only in the formalization of known proofs, but in the discovery of new ones.
Indeed, many of the proofs described in this book were actually \emph{first} done in a fully formalized form in a proof assistant, and are only now being ``unformalized'' for the first time --- a reversal of the usual relation between formal and informal mathematics.

One can imagine a not-too-distant future when it will be possible for mathematicians to verify the correctness of their own papers by working within the system of univalent foundations, formalized in a proof assistant, and that doing so will become as natural as typesetting their own papers in \TeX.
In principle, this could be equally true for any other foundational system, but we believe it to be more practically attainable using univalent foundations, as witnessed by the present work and its formal counterpart.

\index{type theory!formal|)}%
\index{informal type theory|)}%
\index{type theory!informal|)}%
\index{mathematics!formalized|)}%

\subsection*{Constructivity} 

\index{mathematics!constructive|(}%

One of the most striking differences between classical\index{mathematics!classical} foundations and type theory is the idea of \emph{proof relevance}, according to which mathematical statements, and even their proofs, become first-class mathematical objects.
In type theory, we represent mathematical statements by types, which can be regarded simultaneously as both mathematical constructions and mathematical assertions, a conception also known as \emph{propositions as types}.
\index{proposition!as types}%
Accordingly, we can regard a term $a : A$ as both an element of the type $A$ (or in homotopy type theory, a point of the space $A$), and at the same time, a proof of the proposition $A$.
To take an example, suppose we have sets $A$ and $B$ (discrete spaces),
\index{discrete!space}%
and consider the statement ``$A$ is isomorphic to $B$.''
In type theory, this can be rendered as:
\begin{narrowmultline*}
  \mathsf{Iso}(A,B) \defeq \narrowbreak
  \sm{f : A\to B}{g : B\to A}\Big(\big(\tprd{x:A} g(f(x)) = x\big) \times \big(\tprd{y:B}\, f(g(y)) = y\big)\Big).
\end{narrowmultline*}
Reading the type constructors $\Sigma, \Pi, \times$  here  as ``there exists'', ``for all'', and ``and'' respectively yields the usual formulation of ``$A$ and $B$ are isomorphic''; on the other hand, reading them as sums and products yields the \emph{type of all isomorphisms} between $A$ and $B$!  To prove that $A$ and $B$ are isomorphic, one  constructs a proof $p : \mathsf{Iso}(A,B)$, which is therefore the same  as constructing an isomorphism between $A$ and $B$, i.e., exhibiting a pair of functions $f, g$ together with \emph{proofs} that their composites are the respective identity maps.  The latter proofs, in turn, are nothing but homotopies of the appropriate sorts.  In this way, \emph{proving a proposition is the same as constructing an element of some particular type.}

In particular, to prove a statement of the form ``$A$ and $B$'' is just to prove $A$ and to prove $B$, i.e., to give an element of the type $A\times B$.
And to prove that $A$ implies $B$ is just to find an element of $A\to B$, i.e.\ a function from $A$ to $B$ (determining a mapping of proofs of $A$ to proofs of $B$).
This ``constructive'' conception (for more on which, see~\cite{kolmogorov,TroelstraI,TroelstraII}) is what gives type theory its good 
computational character.
For instance, every proof that something exists carries with it enough information to actually find such an object; and from a proof that  ``$A$ or $B$'' holds, one can extract either a proof that $A$ holds or one that $B$ holds.
Thus, from every proof we can automatically extract an algorithm;\index{algorithm} \index{extraction of algorithms} this can be very useful in applications to computer programming.

However, this conception of logic does behave in ways that are unfamiliar to most mathematicians.
\index{axiom!of choice}%
On one hand, a naive translation of the \emph{axiom of choice} yields a statement that we can simply prove.
Essentially, this notion of ``there exists'' is strong enough to ensure that, by showing that for every $x: A$ there exists a $y:B$ such that $R(x,y)$, we automatically construct a function $f : A\to B$ such that, for all $x:A$, we have $R(x, f(x))$.

\index{excluded middle}%
\index{univalence axiom}%
On the other hand, this notion of ``or'' is so strong that a naive translation of the \emph{law of excluded middle} is inconsistent with the univalence axiom.
For if we assume ``for all $A$, either $A$ or not $A$'', then since proving ``$A$'' means exhibiting an element of it, we would have a uniform way of selecting an element from every nonempty type --- a sort of Hilbertian choice operator.
However, univalence implies that the element of $A$ selected by such a choice operator must be invariant under all self-equivalences of $A$, since these are identified with self-identities and every operation must respect identity.
But clearly some types have automorphisms with no fixed points, e.g.\ we can swap the elements of a two-element type.
\index{automorphism!fixed-point-free}%

Thus, the logic of ``proposition as types'' suggested by traditional type theory is not the ``classical'' logic familiar to most mathematicians.
But it is also different from the logic sometimes called ``intuitionistic'', which may lack \emph{both} the law of excluded middle and the axiom of choice. For present purposes, it may be called \emph{constructive logic}
\index{logic!constructive vs classical}%
(but one should be aware that the terms  ``intuitionistic'' and ``constructive'' are often used differently).

The computational advantages of constructive logic imply that we should not  discard it lightly; but for some purposes in classical mathematics, its non-classical character can be problematic.
Many mathematicians are, of course, accustomed to rely on the law of excluded middle; while the ``axiom of choice'' that is available in constructive logic looks superficially similar to its classical namesake, but does not have all of its strong consequences.
Fortunately, homotopy type theory gives a finer analysis of this situation,  allowing various different kinds of logic to coexist and intermix.

The new insight that makes this possible is that the system of all types, just like spaces in classical homotopy theory, is ``stratified'' according to the dimensions in which their higher homotopy structure exists or collapses.
In particular, Voevodsky has found a purely type-theoretic definition of \emph{homotopy $n$-types}, corresponding to spaces with no nontrivial homotopy information above dimension $n$.
(The $0$-types are the ``sets'' mentioned previously as satisfying Lawvere's axioms\index{Lawvere}.)
Moreover, with higher inductive types, we can universally ``truncate'' a type into an $n$-type; in classical homotopy theory this would be its $n^{\mathrm{th}}$ Postnikov\index{Postnikov tower} section.\index{n-type@$n$-type}

With these notions in hand, the homotopy $(-1)$-types, which we call \emph{(mere) propositions}, support a logic that is much more like traditional ``intuitionistic'' logic.
(Classically, every $(-1)$-type is empty or contractible; we interpret these possibilities as the truth values ``false'' and ``true'' respectively.)
The ``$(-1)$-truncated axiom of choice'' is not automatically true, but is a strong assumption with the same sorts of consequences as its counterpart in classical\index{mathematics!classical} set theory.
Similarly, the ``$(-1)$-truncated law of excluded middle'' may be assumed, with many of the same consequences as in classical mathematics.
Thus, the homotopical perspective reveals that classical and constructive logic can coexist, as endpoints of a spectrum of different systems, with an infinite number of possibilities in between (the homotopy $n$-types for $-1 < n < \infty$).
We may speak of ``\LEM{n}'' and ``\choice{n}'', with $\choice{\infty}$ being provable and \LEM{\infty} inconsistent with univalence, while $\choice{-1}$ and $\LEM{-1}$ are the versions familiar to classical mathematicians (hence in most cases it is appropriate to assume the subscript $(-1)$ when none is given).  Indeed, one can even have useful systems in which only \emph{certain} types satisfy such further ``classical'' principles, while types in general remain ``constructive.''\index{excluded middle}\index{axiom!of choice}

It is worth emphasizing that univalent foundations does not \emph{require} the use of constructive or intuitionistic logic.\index{logic!intuitionistic}\index{logic!constructive} %
Most of classical mathematics which depends on the law of excluded middle and the axiom of choice can be performed in univalent foundations, simply by assuming that these two principles hold (in their proper, $(-1)$-truncated, form).
However, type theory does encourage avoiding these principles when they are unnecessary, for several reasons.

First of all, every mathematician knows that a theorem is more powerful when proven using fewer assumptions, since it applies to more examples.
The situation with \choice{} and \LEM{} is no different:
type theory admits many interesting ``nonstandard'' models, such as in sheaf toposes,\index{topos} where classicality principles such as \choice{} and \LEM{} tend to fail.
Homotopy type theory admits similar models in higher toposes, such as are studied in~\cite{ToenVezzosi02,Rezk05,lurie:higher-topoi}.
Thus, if we avoid using these principles, the theorems we prove will be valid internally to all such models.

Secondly, one of the additional virtues of type theory is its computable character.
In addition to being a foundation for mathematics, type theory is a formal theory of computation, and can be treated as a powerful programming language.
\index{programming}%
From this perspective, the rules of the system cannot be chosen arbitrarily the way set-theoretic axioms can: there must be a harmony between them which allows all proofs to be ``executed'' as programs.
We do not yet fully understand the new principles introduced by homotopy type theory, such as univalence and higher inductive types, from
this point of view, but the basic outlines are emerging; see, for example,~\cite{lh:canonicity}.
It has been known for a long time, however, that principles such as \choice{} and \LEM{} are fundamentally antithetical to computability, since they assert baldly that certain things exist without giving any way to compute them.
Thus, avoiding them is necessary to maintain the character of type theory as a theory of computation.

Fortunately, constructive reasoning is not as hard as it may seem.
In some cases, simply by rephrasing some definitions, a theorem can be made constructive and its proof more elegant.
Moreover, in univalent foundations this seems to happen more often.
For instance:
\begin{enumerate}
\item In set-theoretic foundations, at various points in homotopy theory and category theory one needs the axiom of choice to perform transfinite constructions.
  But with higher inductive types, we can encode these constructions directly and constructively.
  In particular, none of the ``synthetic'' homotopy theory in \autoref{cha:homotopy} requires \LEM{} or \choice{}.
\item In set-theoretic foundations, the statement ``every fully faithful and essentially surjective functor is an equivalence of categories'' is equiv\-a\-lent to the axiom of choice.
  But with the univalence axiom, it is just \emph{true}; see \autoref{cha:category-theory}.
\item In set theory, various circumlocutions are required to obtain notions of ``cardinal number'' and ``ordinal number'' which canonically represent isomorphism classes of sets and well-ordered sets, respectively --- possibly involving the axiom of choice or the axiom of foundation.
  But with univalence and higher inductive types, we can obtain such representatives directly by truncating the universe; see \autoref{cha:set-math}.
\item In set-theoretic foundations, the definition of the real numbers as equivalence classes of Cauchy sequences requires either the law of excluded middle or the axiom of (countable) choice to be well-behaved.
  But with higher inductive types, we can give a version of this definition which is well-behaved and avoids any choice principles; see \autoref{cha:real-numbers}.
\end{enumerate}
Of course, these simplifications could as well be taken as evidence that the new methods will not, ultimately, prove to be really constructive.  However, we emphasize again that the reader does not have to care, or worry, about constructivity in order to read this book.  The point is that in all of the above examples, the version of the theory we give has independent advantages, whether or not \LEM{} and \choice{} are assumed to be available.  Constructivity, if attained, will be an added bonus.\index{constructivity}%

Given this discussion of adding new principles such as univalence, higher inductive types, \choice{}, and \LEM{}, one may wonder whether the resulting system remains consistent.
(One of the original virtues of type theory, relative to set theory, was that it can be seen to be consistent by proof-theoretic means).
As with any foundational system, consistency\index{consistency} is a relative question: ``consistent with respect to what?''
The short answer is that all of the constructions and axioms considered in this book have a model in the category of Kan\index{Kan complex} complexes, due to Voevodsky~\cite{klv:ssetmodel} (see~\cite{ls:hits} for higher inductive types).
Thus, they are known to be consistent relative to ZFC (with as many inaccessible cardinals
\index{inaccessible cardinal}\index{consistency}%
as we need nested univalent universes).
Giving a more traditionally type-theoretic account of this consistency is work in progress (see,
e.g.,~\cite{lh:canonicity,coquand2012constructive}).

We summarize the different points of view of the type-theoretic operations in \autoref{tab:pov}.

\begin{table}[htb]
  \centering
  \OPTsmalltable
 \begin{tabular}{lllll}
    \toprule
       Types && Logic & Sets & Homotopy\\ \addlinespace[2pt]
    \midrule
       $A$ && proposition & set & space\\ \addlinespace[2pt]
       $a:A$ && proof & element & point \\ \addlinespace[2pt]
       $B(x)$ && predicate & family of sets & fibration \\ \addlinespace[2pt]
       $b(x) : B(x)$ && conditional proof & family of elements & section\\ \addlinespace[2pt]
       $\emptyt, \unit$ && $\bot, \top$ & $\emptyset, \{ \emptyset \}$ & $\emptyset, *$\\ \addlinespace[2pt]
       $A + B$ && $A\vee B$ & disjoint union & coproduct\\ \addlinespace[2pt]
       $A\times B$ && $A\wedge B$ & set of pairs & product space\\ \addlinespace[2pt]
       $A\to B$ && $A\Rightarrow B$ & set of functions & function space\\ \addlinespace[2pt]
       $\sm{x:A}B(x)$ &&  $\exists_{x:A}B(x)$ & disjoint sum & total space\\ \addlinespace[2pt]
       $\prd{x:A}B(x)$ &&  $\forall_{x:A}B(x)$ & product & space of sections\\ \addlinespace[2pt]
       $\mathsf{Id}_{A}$ && equality $=$ & $\setof{\pairr{x,x} | x\in A}$ & path space $A^I$ \\ \addlinespace[2pt]
    \bottomrule
  \end{tabular}
  \caption{Comparing points of view on type-theoretic operations}\label{tab:pov}
\end{table}

\index{mathematics!constructive|)}%

\subsection*{Open problems} 

\index{open!problem|(}%

For those interested in contributing to this new branch of mathematics, it may be encouraging to know that there are many interesting open questions.

\index{univalence axiom!constructivity of}%
Perhaps the most pressing of them is the ``constructivity'' of the Univalence Axiom, posed by Voevodsky in \cite{Universe-poly}.
The basic system of type theory follows the structure of Gentzen's natural deduction. Logical connectives are defined by their introduction rules, and have elimination rules justified by computation rules. Following this pattern, and using Tait's computability method, originally designed to analyse G\"odel's Dialectica interpretation, one can show the property of \emph{normalization} for type theory. This in turn implies important properties such as decidability of type-checking (a crucial property since type-checking corresponds to proof-checking, and one can argue that we should be able to ``recognize a proof when we see one''), and the so-called ``canonicity\index{canonicity} property'' that any closed term of the type of natural numbers reduces to a numeral. This last property, and the uniform structure of introduction/elimination rules, are lost when one extends type theory with an axiom, such as the axiom of function extensionality, or the univalence axiom. Voevodsky has formulated a precise mathematical conjecture connected to this question of canonicity for type theory extended with the axiom of Univalence: given a closed term of the type of natural numbers, is it always possible to find a numeral and a proof that this term is equal to this numeral, where this proof of equality may itself use the univalence axiom? More generally, an important issue is whether it is possible to provide a constructive justification of the univalence axiom.
What about if one adds other homotopically motivated constructions, like higher inductive types?
These questions remain open at the present time, although methods are currently being developed to try to find answers.

Another basic issue is the difficulty of working with types, such as the natural numbers, that are essentially sets (i.e., discrete spaces),
\index{discrete!space}%
containing only trivial paths.
At present, homotopy type theory can really only characterize spaces up to homotopy equivalence, which means that these ``discrete spaces'' may only be \emph{homotopy equivalent} to discrete spaces.
Type-theoretically, this means there are many paths that are equal to reflexivity, but not \emph{judgmentally} equal to it (see \cref{sec:types-vs-sets} for the meaning of ``judgmentally'').
While this homotopy-invariance has advantages, these ``meaningless'' identity terms do introduce needless complications into arguments and constructions, so it would be convenient to have a systematic way of eliminating or collapsing them.

A more specialized, but no less important, problem is the relation between homotopy type theory and the research on \emph{higher toposes}%
\index{.infinity1-topos@$(\infty,1)$-topos}
currently happening at the intersection of higher category theory and homotopy theory.
There is a growing conviction among those familiar with both subjects that they are intimately connected.
For instance, the notion of a univalent universe should coincide with that of an object classifier, while higher inductive types should be an ``elementary'' reflection of local presentability.
More generally, homotopy type theory should be the ``internal language'' of $(\infty,1)$-toposes, just as intuitionistic higher-order logic is the internal language of ordinary 1-toposes.
Despite this general consensus, however, details remain to be worked out --- in particular, questions of coherence and strictness remain to be addressed  --- and doing so will undoubtedly lead to further insights into both concepts.

\index{mathematics!formalized}%
But by far the largest field of work to be done is in the ongoing formalization of everyday mathematics in this new system.
Recent successes in formalizing some facts from basic homotopy theory and category theory have been encouraging; some of these are described in \cref{cha:homotopy,cha:category-theory}.
Obviously, however, much work remains to be done.

\index{open!problem|)}%

The homotopy type theory community maintains a web site and group blog at \url{http://homotopytypetheory.org}, as well as a discussion email list.
Newcomers are always welcome!

\subsection*{How to read this book}

This book is divided into two parts.
\autoref{part:foundations}, ``Foundations'', develops the fundamental concepts of homotopy type theory.
This is the mathematical foundation on which the development of specific subjects is built, and which is required for the understanding of the univalent foundations approach. To a programmer, this is ``library code''.
Since univalent foundations is a new and different kind of mathematics, its basic notions take some getting used to; thus \autoref{part:foundations} is fairly extensive.

\autoref{part:mathematics}, ``Mathematics'', consists of four chapters that build on the basic notions of \autoref{part:foundations} to exhibit some of the new things we can do with univalent foundations in four different areas of mathematics: homotopy theory (\autoref{cha:homotopy}), category theory (\autoref{cha:category-theory}), set theory (\autoref{cha:set-math}), and real analysis (\autoref{cha:real-numbers}).
The chapters in \autoref{part:mathematics} are more or less independent of each other, although occasionally one will use a lemma proven in another.

A reader who wants to seriously understand univalent foundations, and be able to work in it, will eventually have to read and understand most of \autoref{part:foundations}.
However, a reader who just wants to get a taste of univalent foundations and what it can do may understandably balk at having to work through over 200 pages before getting to the ``meat'' in \autoref{part:mathematics}.
Fortunately, not all of \autoref{part:foundations} is necessary in order to read the chapters in \autoref{part:mathematics}.
Each chapter in \autoref{part:mathematics} begins with a brief overview of its subject, what univalent foundations has to contribute to it, and the necessary background from \autoref{part:foundations}, so the courageous reader can turn immediately to the appropriate chapter for their favorite subject.
For those who want to understand one or more chapters in \autoref{part:mathematics} more deeply than this, but are not ready to read all of \autoref{part:foundations}, we provide here a brief summary of \autoref{part:foundations}, with remarks about which parts are necessary for which chapters in \autoref{part:mathematics}.

\autoref{cha:typetheory} is about the basic notions of type theory, prior to any homotopical interpretation.
A reader who is familiar with Martin-L\"of type theory can quickly skim it to pick up the particulars of the theory we are using.
However, readers without experience in type theory will need to read \autoref{cha:typetheory}, as there are many subtle differences between type theory and other foundations such as set theory.

\autoref{cha:basics} introduces the homotopical viewpoint on type theory, along with the basic notions supporting this view, and describes the homotopical behavior of each component of the type theory from \autoref{cha:typetheory}.
It also introduces the \emph{univalence axiom} (\autoref{sec:compute-universe}) --- the first of the two basic innovations of homotopy type theory.
Thus, it is quite basic and we encourage everyone to read it, especially \crefrange{sec:equality}{sec:basics-equivalences}.

\autoref{cha:logic} describes how we represent logic in homotopy type theory, and its connection to classical logic as well as to constructive and intuitionistic logic.
Here we define the law of excluded middle, the axiom of choice, and the axiom of propositional resizing (although, for the most part, we do not need to assume any of these in the rest of the book), as well as the \emph{propositional truncation} which is essential for representing traditional logic.
This chapter is essential background for \autoref{cha:set-math,cha:real-numbers}, less important for \autoref{cha:category-theory}, and not so necessary for \autoref{cha:homotopy}.

\autoref{cha:equivalences,cha:induction} study two special topics in detail: equivalences (and related notions) and generalized inductive definitions.
While these are important subjects in their own rights and provide a deeper understanding of homotopy type theory, for the most part they are not necessary for \autoref{part:mathematics}.
Only a few lemmas from \autoref{cha:equivalences} are used here and there, while the general discussions in \autoref{sec:bool-nat,sec:strictly-positive,sec:generalizations} are helpful for providing the intuition required for \autoref{cha:hits}.
The generalized sorts of inductive definition discussed in \autoref{sec:generalizations} are also used in a few places in \autoref{cha:set-math,cha:real-numbers}.

\autoref{cha:hits} introduces the second basic innovation of homotopy type theory --- \emph{higher inductive types} --- with many examples.
Higher inductive types are the primary object of study in \autoref{cha:homotopy}, and some particular ones play important roles in \autoref{cha:set-math,cha:real-numbers}.
They are not so necessary for \autoref{cha:category-theory}, although one example is used in \autoref{sec:rezk}.

Finally, \autoref{cha:hlevels} discusses homotopy $n$-types and related notions such as $n$-connected types.
These notions are important for \autoref{cha:homotopy}, but not so important in the rest of \autoref{part:mathematics}, although the case $n=-1$ of some of the lemmas are used in \autoref{sec:piw-pretopos}.

This completes \autoref{part:foundations}.
As mentioned above, \autoref{part:mathematics} consists of four largely unrelated chapters, each describing what univalent foundations has to offer to a particular subject.

Of the chapters in \autoref{part:mathematics}, \autoref{cha:homotopy} (Homotopy theory) is perhaps the most radical.
Univalent foundations has a very different ``synthetic'' approach to homotopy theory in which homotopy types are the basic objects (namely, the types) rather than being constructed using topological spaces or some other set-theoretic model.
This enables new styles of proof for classical theorems in algebraic topology, of which we present a sampling, from $\pi_1(\Sn^1)=\Z$ to the Freudenthal suspension theorem.

In \autoref{cha:category-theory} (Category theory), we develop some basic (1-)category theory, adhering to the principle of the univalence axiom that \emph{equality is isomorphism}.
This has the pleasant effect of ensuring that all definitions and constructions are automatically invariant under equivalence of categories: indeed, equivalent categories are equal just as equivalent types are equal.
(It also has connections to higher category theory and higher topos theory.)

\autoref{cha:set-math} (Set theory) studies sets in univalent foundations.
The category of sets has its usual properties, hence provides a foundation for any mathematics that doesn't need homotopical or higher-categorical structures.
We also observe that univalence makes cardinal and ordinal numbers a bit more pleasant, and that higher inductive types yield a cumulative hierarchy satisfying the usual axioms of Zermelo--Fraenkel set theory.

In \autoref{cha:real-numbers} (Real numbers), we summarize the construction of Dedekind real numbers, and then observe that higher inductive types allow a definition of Cauchy real numbers that avoids some associated problems in constructive mathematics.
Then we sketch a similar approach to Conway's surreal numbers.

Each chapter in this book ends with a Notes section, which collects historical comments, references to the literature, and attributions of results, to the extent possible.
We have also included Exercises at the end of each chapter, to assist the reader in gaining familiarity with doing mathematics in univalent foundations.

Finally, recall that this book was written as a massively collaborative effort by a large number of people.
We have done our best to achieve consistency in terminology and notation, and to put the mathematics in a linear sequence that flows logically, but it is very likely that some imperfections remain.
We ask the reader's forgiveness for any such infelicities, and welcome suggestions for improvement of the next edition.


\part{Foundations}
\label{part:foundations}

\chapter{Type theory}
\label{cha:typetheory}

\section{Type theory versus set theory}
\label{sec:types-vs-sets}
\label{sec:axioms}

\index{type theory}
Homotopy type theory is (among other things) a foundational language for mathematics, i.e., an alternative to Zermelo--Fraenkel\index{set theory!Zermelo--Fraenkel} set theory.
However, it behaves differently from set theory in several important ways, and that can take some getting used to.
Explaining these differences carefully requires us to be more formal here than we will be in the rest of the book.
As stated in the introduction, our goal is to write type theory \emph{informally}; but for a mathematician accustomed to set theory, more precision at the beginning can help avoid some common misconceptions and mistakes.

We note that a set-theoretic foundation has two ``layers'': the deductive system of first-order logic,\index{first-order!logic} and, formulated inside this system, the axioms of a particular theory, such as ZFC.
Thus, set theory is not only about sets, but rather about the interplay between sets (the objects of the second layer) and propositions (the objects of the first layer).

By contrast, type theory is its own deductive system: it need not be formulated inside any superstructure, such as first-order logic.
Instead of the two basic notions of set theory, sets and propositions, type theory has one basic notion: \emph{types}.
Propositions (statements which we can prove, disprove, assume, negate, and so on\footnote{Confusingly, it is also a common practice (dating 
back to Euclid) to use the word ``proposition'' synonymously with ``theorem''.
  We will confine ourselves to the logician's usage, according to which a \emph{proposition} is a statement \emph{susceptible to} proof, whereas a \emph{theorem}\indexfoot{theorem} (or ``lemma''\indexfoot{lemma} or ``corollary''\indexfoot{corollary}) is such a statement that \emph{has been} proven.
Thus ``$0=1$'' and its negation ``$\neg(0=1)$'' are both propositions, but only the latter is a theorem.}) are identified with particular types, via the correspondence shown in \autoref{tab:pov} on page~\pageref{tab:pov}.
Thus, the mathematical activity of \emph{proving a theorem} is identified with a special case of the mathematical activity of \emph{constructing an object}---in this case, an inhabitant of a type that represents a proposition.

\index{deductive system}%
This leads us to another difference between type theory and set theory, but to explain it we must say a little about deductive systems in general.
Informally, a deductive system is a collection of \define{rules}
\indexdef{rule}%
for deriving things called \define{judgments}.
\indexdef{judgment}%
If we think of a deductive system as a formal game,
\index{game!deductive system as}%
then the judgments are the ``positions'' in the game which we reach by following the game rules.
We can also think of a deductive system as a sort of algebraic theory, in which case the judgments are the elements (like the elements of a group) and the deductive rules are the operations (like the group multiplication).
From a logical point of view, the judgments can be considered to be the ``external'' statements, living in the metatheory, as opposed to the ``internal'' statements of the theory itself.

In the deductive system of first-order logic (on which set theory is based), there is only one kind of judgment: that a given proposition has a proof.
That is, each proposition $A$ gives rise to a judgment ``$A$ has a proof'', and all judgments are of this form.
A rule of first-order logic such as ``from $A$ and $B$ infer $A\wedge B$'' is actually a rule of ``proof construction'' which says that given the judgments ``$A$ has a proof'' and ``$B$ has a proof'', we may deduce that ``$A\wedge B$ has a proof''.
Note that the judgment ``$A$ has a proof'' exists at a different level from the \emph{proposition} $A$ itself, which is an internal statement of the theory.

The basic judgment of type theory, analogous to ``$A$ has a proof'', is written ``$a:A$'' and pronounced as ``the term $a$ has type $A$'', or more loosely ``$a$ is an element of $A$'' (or, in homotopy type theory, ``$a$ is a point of $A$'').
\indexdef{term}%
\indexdef{element}%
\indexdef{point!of a type}%
When $A$ is a type representing a proposition, then $a$ may be called a \emph{witness}\index{witness!to the truth of a proposition} to the provability of $A$, or \emph{evidence}\index{evidence, of the truth of a proposition} of the truth of $A$ (or even a \emph{proof}\index{proof} of $A$, but we will try to avoid this confusing terminology).
In this case, the judgment $a:A$ is derivable in type theory (for some $a$) precisely when the analogous judgment ``$A$ has a proof'' is derivable in first-order logic (modulo differences in the axioms assumed and in the encoding of mathematics, as we will discuss throughout the book).
 
On the other hand, if the type $A$ is being treated more like a set than like a proposition (although as we will see, the distinction can become blurry), then ``$a:A$'' may be regarded as analogous to the set-theoretic statement ``$a\in A$''.
However, there is an essential difference in that ``$a:A$'' is a \emph{judgment} whereas ``$a\in A$'' is a \emph{proposition}.
In particular, when working internally in type theory, we cannot make statements such as ``if $a:A$ then it is not the case that $b:B$'', nor can we ``disprove'' the judgment ``$a:A$''.

A good way to think about this is that in set theory, ``membership'' is a relation which may or may not hold between two pre-existing objects ``$a$'' and ``$A$'', while in type theory we cannot talk about an element ``$a$'' in isolation: every element \emph{by its very nature} is an element of some type, and that type is (generally speaking) uniquely determined.
Thus, when we say informally ``let $x$ be a natural number'', in set theory this is shorthand for ``let $x$ be a thing and assume that $x\in\nat$'', whereas in type theory ``let $x:\nat$'' is an atomic statement: we cannot introduce a variable without specifying its type.\index{membership}

At first glance, this may seem an uncomfortable restriction, but it is arguably closer to the intuitive mathematical meaning of ``let $x$ be a natural number''.
In practice, it seems that whenever we actually \emph{need} ``$a\in A$'' to be a proposition rather than a judgment, there is always an ambient set $B$ of which $a$ is known to be an element and $A$ is known to be a subset.
This situation is also easy to represent in type theory, by taking $a$ to be an element of the type $B$, and $A$ to be a predicate on $B$; see \autoref{subsec:prop-subsets}.

A last difference between type theory and set theory is the treatment of equality.
The familiar notion of equality in mathematics is a proposition: e.g.\ we can disprove an equality or assume an equality as a hypothesis.
Since in type theory, propositions are types, this means that equality is a type: for elements $a,b:A$ (that is, both $a:A$ and $b:A$) we have a type ``$\id[A]ab$''.
(In \emph{homotopy} type theory, of course, this equality proposition can behave in unfamiliar ways: see \autoref{sec:identity-types,cha:basics}, and the rest of the book).
When $\id[A]ab$ is inhabited, we say that $a$ and $b$ are \define{(propositionally) equal}.
\index{propositional!equality}%
\index{equality!propositional}%

However, in type theory there is also a need for an equality \emph{judgment}, existing at the same level as the judgment ``$x:A$''.\index{judgment}
\symlabel{defn:judgmental-equality}%
This is called \define{judgmental equality}
\indexdef{equality!judgmental}%
\indexdef{judgmental equality}%
or \define{definitional equality},
\indexdef{equality!definitional}%
\indexsee{definitional equality}{equality, definitional}%
and we write it as $a\jdeq b : A$ or simply $a \jdeq b$.
It is helpful to think of this as meaning ``equal by definition''.
For instance, if we define a function $f:\nat\to\nat$ by the equation $f(x)=x^2$, then the expression $f(3)$ is equal to $3^2$ \emph{by definition}.
Inside the theory, it does not make sense to negate or assume an equality-by-definition; we cannot say ``if $x$ is equal to $y$ by definition, then $z$ is not equal to $w$ by definition''.
Whether or not two expressions are equal by definition is just a matter of expanding out the definitions; in particular, it is algorithmically\index{algorithm} decidable (though the algorithm is necessarily meta-theoretic, not internal to the theory).\index{decidable!definitional equality}

As type theory becomes more complicated, judgmental equality can get more subtle than this, but it is a good intuition to start from.
Alternatively, if we regard a deductive system as an algebraic theory, then judgmental equality is simply the equality in that theory, analogous to the equality between elements of a group---the only potential for confusion is that there is \emph{also} an object \emph{inside} the deductive system of type theory (namely the type ``$a=b$'') which behaves internally as a notion of ``equality''.

The reason we \emph{want} a judgmental notion of equality is so that it can control the other form of judgment, ``$a:A$''.
For instance, suppose we have given a proof that $3^2=9$, i.e.\ we have derived the judgment $p:(3^2=9)$ for some $p$.
Then the same witness $p$ ought to count as a proof that $f(3)=9$, since $f(3)$ is $3^2$ \emph{by definition}.
The best way to represent this is with a rule saying that given the judgments $a:A$ and $A\jdeq B$, we may derive the judgment $a:B$.

Thus, for us, type theory will be a deductive system based on two forms of judgment:
\begin{center}
\medskip
\begin{tabular}{cl}
  \toprule
  Judgment & Meaning\\
  \midrule
  $a : A$       & ``$a$ is an object of type $A$''\\
  $a \jdeq b : A$ & ``$a$ and $b$ are definitionally equal objects of type $A$''\\
  \bottomrule
\end{tabular}
\medskip
\end{center}
\symlabel{defn:defeq}%
When introducing a definitional equality, i.e., defining one thing to be equal to another, we will use the symbol ``$\defeq$''.
Thus, the above definition of the function $f$ would be written as $f(x)\defeq x^2$.

Because judgments cannot be put together into more complicated statements, the symbols ``$:$'' and ``$\jdeq$'' bind more loosely than anything else.%
\footnote{In formalized\indexfoot{mathematics!formalized} type theory, commas and turnstiles can bind even more loosely.
  For instance, $x:A,y:B\vdash c:C$ is parsed as $((x:A),(y:B))\vdash (c:C)$.
  However, in this book we refrain from such notation until \autoref{cha:rules}.}
Thus, for instance, ``$p:\id{x}{y}$'' should be parsed as ``$p:(\id{x}{y})$'', which makes sense since ``$\id{x}{y}$'' is a type, and not as ``$\id{(p:x)}{y}$'', which is senseless since ``$p:x$'' is a judgment and cannot be equal to anything.
Similarly, ``$A\jdeq \id{x}{y}$'' can only be parsed as ``$A\jdeq(\id{x}{y})$'', although in extreme cases such as this, one ought to add parentheses anyway to aid reading comprehension.
Moreover, later on we will fall into the common notation of chaining together equalities --- e.g.\ writing $a=b=c=d$ to mean ``$a=b$ and $b=c$ and $c=d$, hence $a=d$'' --- and we will also include judgmental equalities in such chains.
Context usually suffices to make the intent clear.

This is perhaps also an appropriate place to mention that the common mathematical notation ``$f:A\to B$'', expressing the fact that $f$ is a function from $A$ to $B$, can be regarded as a typing judgment, since we use ``$A\to B$'' as notation for the type of functions from $A$ to $B$ (as is standard practice in type theory; see \autoref{sec:pi-types}).

\index{assumption|(defstyle}%
Judgments may depend on \emph{assumptions} of the form $x:A$, where $x$ is a variable
\indexdef{variable}%
and $A$ is a type.
For example, we may construct an object $m + n : \nat$ under the assumptions that $m,n : \nat$.
Another example is that assuming $A$ is a type, $x,y : A$, and $p : \id[A]{x}{y}$, we may construct an element $p^{-1} : \id[A]{y}{x}$.
The collection of all such assumptions is called the \define{context};%
\index{context}
from a topological point of view it may be thought of as a ``parameter\index{parameter!space} space''.
In fact, technically the context must be an ordered list of assumptions, since later assumptions may depend on previous ones: the assumption $x:A$ can only be made \emph{after} the assumptions of any variables appearing in the type $A$.

If the type $A$ in an assumption $x:A$ represents a proposition, then the assumption is a type-theoretic version of a \emph{hypothesis}:
\indexdef{hypothesis}%
we assume that the proposition $A$ holds.
When types are regarded as propositions, we may omit the names of their proofs.
Thus, in the second example above we may instead say that assuming $\id[A]{x}{y}$, we can prove $\id[A]{y}{x}$.
However, since we are doing ``proof-relevant'' mathematics,
\index{mathematics!proof-relevant}%
we will frequently refer back to proofs as objects.
In the example above, for instance, we may want to establish that $p^{-1}$ together with the proofs of transitivity and reflexivity behave like a groupoid; see \autoref{cha:basics}.

Note that under this meaning of the word \emph{assumption}, we can assume a propositional equality (by assuming a variable $p:x=y$), but we cannot assume a judgmental equality $x\jdeq y$, since it is not a type that can have an element.
However, we can do something else which looks kind of like assuming a judgmental equality: if we have a type or an element which involves a variable $x:A$, then we can \emph{substitute} any particular element $a:A$ for $x$ to obtain a more specific type or element.
We will sometimes use language like ``now assume $x\jdeq a$'' to refer to this process of substitution, even though it is not an \emph{assumption} in the technical sense introduced above.
\index{assumption|)}%

By the same token, we cannot \emph{prove} a judgmental equality either, since it is not a type in which we can exhibit a witness.
Nevertheless, we will sometimes state judgmental equalities as part of a theorem, e.g.\ ``there exists $f:A\to B$ such that $f(x)\jdeq y$''.
This should be regarded as the making of two separate judgments: first we make the judgment $f:A\to B$ for some element $f$, then we make the additional judgment that $f(x)\jdeq y$.

In the rest of this chapter, we attempt to give an informal presentation of type theory, sufficient for the purposes of this book; we give a more formal account in \autoref{cha:rules}.
Aside from some fairly obvious rules (such as the fact that judgmentally equal things can always be substituted\index{substitution} for each other), the rules of type theory can be grouped into \emph{type formers}.
Each type former consists of a way to construct types (possibly making use of previously constructed types), together with rules for the construction and behavior of elements of that type.
In most cases, these rules follow a fairly predictable pattern, but we will not attempt to make this precise here; see however the beginning of \autoref{sec:finite-product-types} and also \autoref{cha:induction}.\index{type theory!informal}

\index{axiom!versus rules}%
\index{rule!versus axioms}%
An important aspect of the type theory presented in this chapter is that it consists entirely of \emph{rules}, without any \emph{axioms}.
In the description of deductive systems in terms of judgments, the \emph{rules} are what allow us to conclude one judgment from a collection of others, while the \emph{axioms} are the judgments we are given at the outset.
If we think of a deductive system as a formal game, then the rules are the rules of the game, while the axioms are the starting position.
And if we think of a deductive system as an algebraic theory, then the rules are the operations of the theory, while the axioms are the \emph{generators} for some particular free model of that theory.

In set theory, the only rules are the rules of first-order logic (such as the rule allowing us to deduce ``$A\wedge B$ has a proof'' from ``$A$ has a proof'' and ``$B$ has a proof''): all the information about the behavior of sets is contained in the axioms.
By contrast, in type theory, it is usually the \emph{rules} which contain all the information, with no axioms being necessary.
For instance, in \autoref{sec:finite-product-types} we will see that there is a rule allowing us to deduce the judgment ``$(a,b):A\times B$'' from ``$a:A$'' and ``$b:B$'', whereas in set theory the analogous statement would be (a consequence of) the pairing axiom.

The advantage of formulating type theory using only rules is that rules are ``procedural''.
In particular, this property is what makes possible (though it does not automatically ensure) the good computational properties of type theory, such as ``canonicity''.\index{canonicity}
However, while this style works for traditional type theories, we do not yet understand how to formulate everything we need for \emph{homotopy} type theory in this way.
In particular, in \autoref{sec:compute-pi,sec:compute-universe,cha:hits} we will have to augment the rules of type theory presented in this chapter by introducing additional axioms, notably the \emph{univalence axiom}.
In this chapter, however, we confine ourselves to a traditional rule-based type theory.

\section{Function types}
\label{sec:function-types}

\index{type!function|(defstyle}%
\indexsee{function type}{type, function}%
Given types $A$ and $B$, we can construct the type $A \to B$ of \define{functions}
\index{function|(defstyle}%
\indexsee{map}{function}%
\indexsee{mapping}{function}%
with domain $A$ and codomain $B$.
We also sometimes refer to functions as \define{maps}.
\index{domain!of a function}%
\index{codomain, of a function}%
\index{function!domain of}%
\index{function!codomain of}%
\index{functional relation}%
Unlike in set theory, functions are not defined as
functional relations; rather they are a primitive concept in type theory.
We explain the function type by prescribing what we can do with functions, 
how to construct them and what equalities they induce.

Given a function $f : A \to B$ and an element of the domain $a : A$, we
can \define{apply}
\indexdef{application!of function}%
\indexdef{function!application}%
\indexsee{evaluation}{application, of a function}
the function to obtain an element of the codomain $B$,
denoted $f(a)$ and called the \define{value} of $f$ at $a$.
\indexdef{value!of a function}%
It is common in type theory to omit the parentheses\index{parentheses} and denote $f(a)$ simply by $f\,a$, and we will sometimes do this as well.

But how can we construct elements of $A \to B$? There are two equivalent ways:
either by direct definition or by using
$\lambda$-abstraction. Introducing a function by definition
\indexdef{definition!of function, direct}%
means that
we introduce a function by giving it a name --- let's say, $f$ --- and saying
we define $f : A \to B$ by giving an equation
\begin{equation}
  \label{eq:expldef}
  f(x) \defeq \Phi
\end{equation}
where $x$ is a variable
\index{variable}%
and $\Phi$ is an expression which may use $x$.
In order for this to be valid, we have to check that $\Phi : B$ assuming $x:A$.

Now we can compute $f(a)$ by replacing the variable $x$ in $\Phi$ with
$a$. As an example, consider the function $f : \nat \to \nat$ which is
defined by $f(x) \defeq x+x$.  (We will define $\nat$ and $+$ in \autoref{sec:inductive-types}.)
Then $f(2)$ is judgmentally equal to $2+2$.

If we don't want to introduce a name for the function, we can use
\define{$\lambda$-abstraction}.
\index{lambda abstraction@$\lambda$-abstraction|defstyle}%
\indexsee{function!lambda abstraction@$\lambda$-abstraction}{$\lambda$-abstraction}%
\indexsee{abstraction!lambda-@$\lambda$-}{$\lambda$-abstraction}%
Given an expression $\Phi$ of type $B$ which may use $x:A$, as above, we write $\lam{x:A} \Phi$ to indicate the same function defined by~\eqref{eq:expldef}.
Thus, we have
\[ (\lamt{x:A}\Phi) : A \to B. \]
For the example in the previous paragraph, we have the typing judgment
\[ (\lam{x:\nat}x+x) : \nat \to \nat. \]
As another example, for any types $A$ and $B$ and any element $y:B$, we have a \define{constant function}
\indexdef{constant!function}%
\indexdef{function!constant}%
$(\lam{x:A} y): A\to B$.

We generally omit the type of the variable $x$ in a $\lambda$-abstraction and write $\lam{x}\Phi$, since the typing $x:A$ is inferable from the judgment that the function $\lam x \Phi$ has type $A\to B$.
By convention, the ``scope''
\indexdef{variable!scope of}%
\indexdef{scope}%
of the variable binding ``$\lam{x}$'' is the entire rest of the expression, unless delimited with parentheses\index{parentheses}.
Thus, for instance, $\lam{x} x+x$ should be parsed as $\lam{x} (x+x)$, not as $(\lam{x}x)+x$ (which would, in this case, be ill-typed anyway).

Another equivalent notation is
\symlabel{mapsto}%
\[ (x \mapsto \Phi) : A \to B. \]
\symlabel{blank}%
We may also sometimes use a blank ``$\blank$'' in the expression $\Phi$ in place of a variable, to denote an implicit $\lambda$-abstraction.
For instance, $g(x,\blank)$ is another way to write $\lam{y} g(x,y)$.

Now a $\lambda$-abstraction is a function, so we can apply it to an argument $a:A$.
We then have the following \define{computation rule}\indexdef{computation rule!for function types}\footnote{Use of this equality is often referred to as \define{$\beta$-conversion}
\indexsee{beta-conversion@$\beta $-conversion}{$\beta$-reduction}%
\indexsee{conversion!beta@$\beta $-}{$\beta$-reduction}%
or \define{$\beta$-reduction}.%
\index{beta-reduction@$\beta $-reduction|footstyle}%
\indexsee{reduction!beta@$\beta $-}{$\beta$-reduction}%
}, which is a definitional equality:
\[(\lamu{x:A}\Phi)(a) \jdeq \Phi'\]
where $\Phi'$ is the
expression $\Phi$ in which all occurrences of $x$ have been replaced by $a$.
Continuing the above example, we have
\[ (\lamu{x:\nat}x+x)(2) \jdeq 2+2. \]
Note that from any function $f:A\to B$, we can construct a lambda abstraction function $\lam{x} f(x)$.
Since this is by definition ``the function that applies $f$ to its argument'' we consider it to be definitionally equal to $f$:\footnote{Use of this equality is often referred to as \define{$\eta$-conversion}
\indexsee{eta-conversion@$\eta $-conversion}{$\eta$-expansion}%
\indexsee{conversion!eta@$\eta $-}{$\eta$-expansion}%
or \define{$\eta$-expansion.
\index{eta-expansion@$\eta $-expansion|footstyle}%
\indexsee{expansion, eta-@expansion, $\eta $-}{$\eta$-expansion}%
}}
\[ f \jdeq (\lam{x} f(x)). \]
This equality is the \define{uniqueness principle for function types}\indexdef{uniqueness!principle!for function types}, because it shows that $f$ is uniquely determined by its values.

The introduction of functions by definitions with explicit parameters can be reduced
to simple definitions by using $\lambda$-abstraction: i.e., we can read 
a definition of $f: A\to B$ by
\[ f(x) \defeq \Phi \]
as 
\[ f \defeq \lamu{x:A}\Phi.\]

When doing calculations involving variables, we have to be 
careful when replacing a variable with an expression that also involves
variables, because we want to preserve the binding structure of
expressions. By the \emph{binding structure}\indexdef{binding structure} we mean the
invisible link generated by binders such as $\lambda$, $\Pi$ and
$\Sigma$ (the latter we are going to meet soon) between the place where the variable is introduced and where it is used. As an example, consider $f : \nat \to (\nat \to \nat)$
defined as 
\[ f(x) \defeq \lamu{y:\nat} x + y. \] 
Now if we have assumed somewhere that $y : \nat$, then what is $f(y)$? It would be wrong to just naively replace $x$ by $y$ everywhere in the expression ``$\lam{y}x+y$'' defining $f(x)$, obtaining $\lamu{y:\nat} y + y$, because this means that $y$ gets \define{captured}.
\indexdef{capture, of a variable}%
\indexdef{variable!captured}%
Previously, the substituted\index{substitution} $y$ was referring to our assumption, but now it is referring to the argument of the $\lambda$-abstraction. Hence, this naive substitution would destroy the binding structure, allowing us to perform calculations which are semantically unsound.

But what \emph{is} $f(y)$ in this example? Note that bound (or ``dummy'')
variables
\indexdef{variable!bound}%
\indexdef{variable!dummy}%
\indexsee{bound variable}{variable, bound}%
\indexsee{dummy variable}{variable, bound}%
such as $y$ in the expression $\lamu{y:\nat} x + y$
have only a local meaning, and can be consistently replaced by any
other variable, preserving the binding structure. Indeed, $\lamu{y:\nat} x + y$ is declared to be judgmentally equal\footnote{Use of this equality is often referred to as \define{$\alpha$-conversion.
\indexfoot{alpha-conversion@$\alpha $-conversion}
\indexsee{conversion!alpha@$\alpha$-}{$\alpha$-conversion}
}} to
$\lamu{z:\nat} x + z$.  It follows that 
$f(y)$ is judgmentally equal to  $\lamu{z:\nat} y + z$, and that answers our question.  (Instead of $z$,
any variable distinct from $y$ could have been used, yielding an equal result.)

Of course, this should all be familiar to any mathematician: it is the same phenomenon as the fact that if $f(x) \defeq \int_1^2 \frac{dt}{x-t}$, then $f(t)$ is not $\int_1^2 \frac{dt}{t-t}$ but rather $\int_1^2 \frac{ds}{t-s}$.
A $\lambda$-abstraction binds a dummy variable in exactly the same way that an integral does.

We have seen how to define functions in one variable. One
way to define functions in several variables would be to use the
cartesian product, which will be introduced later; a function with
parameters $A$ and $B$ and results in $C$ would be given the type 
$f : A \times B \to C$. However, there is another choice that avoids
using product types, which is called \define{currying}
\indexdef{currying}%
\indexdef{function!currying of}%
(after the mathematician Haskell Curry).
\index{programming}%

The idea of currying is to represent a function of two inputs $a:A$ and $b:B$ as a function which takes \emph{one} input $a:A$ and returns \emph{another function}, which then takes a second input $b:B$ and returns the result.
That is, we consider two-variable functions to belong to an iterated function type, $f : A \to (B \to C)$.
We may also write this without the parentheses\index{parentheses}, as $f : A \to B \to C$, with
associativity\index{associativity!of function types} to the right as the default convention.  Then given $a : A$ and $b : B$,
we can apply $f$ to $a$ and then apply the result to $b$, obtaining
$f(a)(b) : C$. To avoid the proliferation of parentheses, we allow ourselves to
write $f(a)(b)$ as $f(a,b)$ even though there are no products
involved.
When omitting parentheses around function arguments entirely, we write $f\,a\,b$ for $(f\,a)\,b$, with the default associativity now being to the left so that $f$ is applied to its arguments in the correct order.

Our notation for definitions with explicit parameters extends to
this situation: we can define a named function $f : A \to B \to C$ by
giving an equation
\[ f(x,y) \defeq \Phi\]
where $\Phi:C$ assuming $x:A$ and $y:B$. Using $\lambda$-abstraction\index{lambda abstraction@$\lambda$-abstraction} this
corresponds to
\[ f \defeq \lamu{x:A}{y:B} \Phi, \]
which may also be written as 
\[ f \defeq x \mapsto y \mapsto \Phi. \] 
We can also implicitly abstract over multiple variables by writing multiple blanks, e.g.\ $g(\blank,\blank)$ means $\lam{x}{y} g(x,y)$.
Currying a function of three or more arguments is a straightforward extension of what we have just described.
 
\index{type!function|)}%
\index{function|)}%

\section{Universes and families}
\label{sec:universes}

So far, we have been using the expression ``$A$ is a type'' informally. We
are going to make this more precise by introducing \define{universes}.
\index{type!universe|(defstyle}%
\indexsee{universe}{type, universe}%
A universe is a type whose elements are types. As in naive set theory,
we might wish for a universe of all types $\UU_\infty$ including itself
(that is, with $\UU_\infty : \UU_\infty$).
However, as in set
theory, this is unsound, i.e.\ we can deduce from it that every type,
including the empty type representing the proposition False (see \autoref{sec:coproduct-types}), is inhabited.
For instance, using a
representation of sets as trees, we can directly encode Russell's
paradox\index{paradox} \cite{coquand:paradox}.

To avoid the paradox we introduce a hierarchy of universes
\indexsee{hierarchy!of universes}{type, universe}%
\[ \UU_0 : \UU_1 : \UU_2 : \cdots \]
where every universe $\UU_i$ is an element of the next universe
$\UU_{i+1}$. Moreover, we assume that our universes are
\define{cumulative},
\indexdef{type!universe!cumulative}%
\indexdef{cumulative!universes}%
that is that all the elements of the $i^{\mathrm{th}}$
universe are also elements of the $(i+1)^{\mathrm{st}}$ universe, i.e.\ if
$A:\UU_i$ then also $A:\UU_{i+1}$.
This is convenient, but has the slightly unpleasant consequence that elements no longer have unique types, and is a bit tricky in other ways that need not concern us here; see the Notes.

When we say that $A$ is a type, we mean that it inhabits some universe
$\UU_i$. We usually want to avoid mentioning the level
\indexdef{universe level}%
\indexsee{level}{universe level or $n$-type}%
\indexsee{type!universe!level}{universe level}%
$i$ explicitly,
and just assume that levels can be assigned in a consistent way; thus we
may write $A:\UU$ omitting the level. This way we can even write
$\UU:\UU$, which can be read as $\UU_i:\UU_{i+1}$, having left the
indices implicit.  Writing universes in this style is referred to as
\define{typical ambiguity}.
\indexdef{typical ambiguity}%
It is convenient but a bit dangerous, since it allows us to write valid-looking proofs that reproduce the paradoxes of self-reference.
If there is any doubt about whether an argument is correct, the way to check it is to try to assign levels consistently to all universes appearing in it.
When some universe \UU is assumed, we may refer to types belonging to \UU as \define{small types}.
\indexdef{small!type}%
\indexdef{type!small}%

To model a collection of types varying over a given type $A$, we use functions $B : A \to \UU$  whose
codomain is a universe. These functions are called
\define{families of types} (or sometimes \emph{dependent types});
\indexsee{family!of types}{type, family of}%
\indexdef{type!family of}%
\indexsee{type!dependent}{type, family of}%
\indexsee{dependent!type}{type, family of}%
they correspond to families of sets as used in
set theory.

\symlabel{fin}%
An example of a type family is the family of finite sets $\Fin
: \nat \to \UU$, where $\Fin(n)$ is a type with exactly $n$ elements.
(We cannot \emph{define} the family $\Fin$ yet --- indeed, we have not even introduced its domain $\nat$ yet --- but we will be able to soon; see \autoref{ex:fin}.)
We may denote the elements of $\Fin(n)$ by $0_n,1_n,\dots,(n-1)_n$, with subscripts to emphasize that the elements of $\Fin(n)$ are different from those of $\Fin(m)$ if $n$ is different from $m$, and all are different from the ordinary natural numbers (which we will introduce in \autoref{sec:inductive-types}).
\index{finite!sets, family of}%

A more trivial (but very important) example of a type family is the \define{constant} type family
\indexdef{constant!type family}%
\indexdef{type!family of!constant}%
at a type $B:\UU$, which is of course the constant function $(\lam{x:A} B):A\to\UU$.

As a \emph{non}-example, in our version of type theory there is no type family ``$\lam{i:\nat} \UU_i$''.
Indeed, there is no universe large enough to be its codomain.
Moreover, we do not even identify the indices $i$ of the universes $\UU_i$ with the natural numbers \nat of type theory (the latter to be introduced in \autoref{sec:inductive-types}).

\index{type!universe|)}%

\section{Dependent function types (\texorpdfstring{$\Pi$}{Π}-types)}
\label{sec:pi-types}

\index{type!dependent function|(defstyle}%
\index{function!dependent|(defstyle}%
\indexsee{dependent!function}{function, dependent}%
\indexsee{type!Pi-@$\Pi$-}{type, dependent function}%
\indexsee{Pi-type@$\Pi$-type}{type, dependent function}%
In type theory we often use a more general version of function
types, called a \define{$\Pi$-type} or \define{dependent function type}. The elements of
a $\Pi$-type are functions
whose codomain type can vary depending on the
element of the domain to which the function is applied, called \define{dependent functions}. The name ``$\Pi$-type''
is used because this type can also be regarded as the  cartesian
product over a given type.

Given a type $A:\UU$ and a family $B:A \to \UU$, we may construct
the type of dependent functions $\prd{x:A}B(x) : \UU$.
There are many alternative notations for this type, such as
\[ \tprd{x:A} B(x) \qquad \dprd{x:A}B(x) \qquad \lprd{x:A} B(x). \]
If $B$ is a constant family, then the dependent product type is the ordinary function type:
\[\tprd{x:A} B \jdeq (A \to B).\]
Indeed, all the constructions of $\Pi$-types are generalizations of the corresponding constructions on ordinary function types.

\indexdef{definition!of function, direct}%
We can introduce dependent functions by explicit definitions: to
define $f : \prd{x:A}B(x)$, where $f$ is the name of a dependent function to be
defined, we need an expression $\Phi : B(x)$ possibly involving the variable $x:A$,
\index{variable}%
and we write
\[ f(x) \defeq \Phi \qquad \mbox{for $x:A$}.\]
Alternatively, we can use \define{$\lambda$-abstraction}%
\index{lambda abstraction@$\lambda$-abstraction|defstyle}%
\begin{equation}
  \label{eq:lambda-abstraction}
  \lamu{x:A} \Phi \ :\ \prd{x:A} B(x).
\end{equation}
\indexdef{application!of dependent function}%
\indexdef{function!dependent!application}%
As with non-dependent functions, we can \define{apply} a dependent function $f : \prd{x:A}B(x)$ to an argument $a:A$ to obtain an element $f(a):B(a)$.
The equalities are the same as for the ordinary function type, i.e.\ we have the computation rule
\index{computation rule!for dependent function types}%
given $a:A$ we have $f(a) \jdeq \Phi'$ and  
$(\lamu{x:A} \Phi)(a) \jdeq \Phi'$, where $\Phi' $ is obtained by replacing all
occurrences of $x$ in $\Phi$ by $a$ (avoiding variable capture, as always).
Similarly, we have the uniqueness principle $f\jdeq (\lam{x} f(x))$ for any $f:\prd{x:A} B(x)$.
\index{uniqueness!principle!for dependent function types}%

As an example, recall from \autoref{sec:universes} that there is a type family $\Fin:\nat\to\UU$ whose values are the standard finite sets, with elements $0_n,1_n,\dots,(n-1)_n : \Fin(n)$.
There is then a dependent function $\fmax : \prd{n:\nat} \Fin(n+1)$
which returns the ``largest'' element of each nonempty finite type, $\fmax(n) \defeq n_{n+1}$.
\index{finite!sets, family of}%
As was the case for $\Fin$ itself, we cannot define $\fmax$ yet, but we will be able to soon; see \autoref{ex:fin}.

Another important class of dependent function types, which we can define now, are functions which are \define{polymorphic}
\indexdef{function!polymorphic}%
\indexdef{polymorphic function}%
over a given universe.
A polymorphic function is one which takes a type as one of its arguments, and then acts on elements of that type (or other types constructed from it).
\symlabel{idfunc}%
\indexdef{function!identity}%
\indexdef{identity!function}%
An example is the polymorphic identity function $\idfunc : \prd{A:\UU} A \to A$, which we define by $\idfunc{} \defeq \lam{A:\type}{x:A} x$.

We sometimes write some arguments of a dependent function as subscripts.
For instance, we might equivalently define the polymorphic identity function by $\idfunc[A](x) \defeq x$.
Moreover, if an argument can be inferred from context, we may omit it altogether.
For instance, if $a:A$, then writing $\idfunc(a)$ is unambiguous, since $\idfunc$ must mean $\idfunc[A]$ in order for it to be applicable to $a$.

Another, less trivial, example of a polymorphic function is the ``swap'' operation that switches the order of the arguments of a (curried) two-argument function:
\[ \mathsf{swap} : \prd{A:\UU}{B:\UU}{C:\UU} (A\to B\to C) \to (B\to A \to C)
\]
We can define this by
\[ \mathsf{swap}(A,B,C,g) \defeq \lam{b}{a} g(a)(b). \]
We might also equivalently write the type arguments as subscripts:
\[ \mathsf{swap}_{A,B,C}(g)(b,a) \defeq g(a,b). \]

Note that as we did for ordinary functions, we use currying to define dependent functions with
several arguments (such as $\mathsf{swap}$). However, in the dependent case the second domain may
depend on the first one, and the codomain may depend on both. That is,
given $A:\UU$ and type families $B : A \to \UU$ and $C : \prd{x:A}B(x) \to \UU$, we may construct
the type $\prd{x:A}{y : B(x)} C(x,y)$ of functions with two
arguments.
(Like $\lambda$-abstractions, $\Pi$s automatically scope\index{scope} over the rest of the expression unless delimited; thus $C : \prd{x:A}B(x) \to \UU$ means $C : \prd{x:A}(B(x) \to \UU)$.)
In the case when $B$ is constant and equal to $A$, we may condense the notation and write $\prd{x,y:A}$; for instance, the type of $\mathsf{swap}$ could also be written as
\[ \mathsf{swap} : \prd{A,B,C:\UU} (A\to B\to C) \to (B\to A \to C).
\]
Finally, given $f:\prd{x:A}{y : B(x)} C(x,y)$ and arguments $a:A$ and $b:B(a)$, we have $f(a)(b) : C(a,b)$, which,
as before, we write as $f(a,b) : C(a,b)$.

\index{type!dependent function|)}%
\index{function!dependent|)}%

\section{Product types}
\label{sec:finite-product-types}

Given types $A,B:\UU$ we introduce the type $A\times B:\UU$, which we call their \define{cartesian product}.
\indexsee{cartesian product}{type, product}%
\indexsee{type!cartesian product}{type, product}%
\index{type!product|(defstyle}%
\indexsee{product!of types}{type, product}%
We also introduce a nullary product type, called the \define{unit type} $\unit : \UU$.
\indexsee{nullary!product}{type, unit}%
\indexsee{unit!type}{type, unit}%
\index{type!unit|(defstyle}%
We intend the elements of $A\times B$ to be pairs $\tup{a}{b} : A \times B$, where $a:A$ and $b:B$, and the only element of $\unit$ to be some particular object $\ttt : \unit$.
\indexdef{pair!ordered}%
However, unlike in set theory, where we define ordered pairs to be particular sets and then collect them all together into the cartesian product, in type theory, ordered pairs are a primitive concept, as are functions.   

\begin{rmk}\label{rmk:introducing-new-concepts}
  There is a general pattern for introduction of a new kind of type in type theory, and because products are our second example following this pattern,\footnote{The description of universes above is an exception.} it is worth emphasizing the general form:
  To specify a type, we specify:
  \begin{enumerate}
  \item how to form new types of this kind, via \define{formation rules}.
    \indexdef{formation rule}%
    \index{rule!formation}%
(For example, we can form the function type $A \to B$ when $A$ is a type and when $B$ is a type. We can form the dependent function type $\prd{x:A} B(x)$ when $A$ is a type and $B(x)$ is a type for $x:A$.)

  \item how to construct elements of that type.  
    These are called the type's \define{constructors} or \define{introduction rules}.
    \indexdef{constructor}%
    \indexdef{rule!introduction}%
    \indexdef{introduction rule}%
    (For example, a function type has one constructor, $\lambda$-abstraction.
    Recall that a direct definition like $f(x)\defeq 2x$ can equivalently be phrased
    as a $\lambda$-abstraction $f\defeq \lam{x} 2x$.)

  \item how to use elements of that type.  
    These are called the type's \define{eliminators} or \define{elimination rules}.
    \indexsee{rule!elimination}{eliminator}%
    \indexsee{elimination rule}{eliminator}%
    \indexdef{eliminator}%
    (For example, the function type has one eliminator, namely function application.)

  \item 
    a \define{computation rule}\indexdef{computation rule}\footnote{also referred to as \define{$\beta$-reduction}\index{beta-reduction@$\beta $-reduction|footstyle}}, which expresses how an eliminator acts on a constructor.
(For example, for functions, the computation rule states that $(\lamu{x:A}\Phi)(a)$ is judgmentally equal to the substitution of $a$ for $x$ in $\Phi$.)

  \item 
    an optional \define{uniqueness principle}\indexdef{uniqueness!principle}\footnote{also referred to as \define{$\eta$-expansion}\index{eta-expansion@$\eta $-expansion|footstyle}}, which expresses
uniqueness of maps into or out of that type.  
For some types, the uniqueness principle characterizes maps into the type, by stating that 
every element of the type is uniquely determined by the results of applying eliminators to it, and can be reconstructed from those results by applying a constructor---thus expressing how constructors act on eliminators, dually to the computation rule.  
(For example, for functions, the uniqueness principle says that any function $f$ is judgmentally equal to the ``expanded'' function $\lamu{x} f(x)$, and thus is uniquely determined by its values.)
For other types, the uniqueness principle says that every map (function) \emph{from} that type is uniquely determined by some data. (An example is the coproduct type introduced in \cref{sec:coproduct-types}, whose uniqueness principle is mentioned in \cref{sec:universal-properties}.)  
    
    When the uniqueness principle is not taken as a rule of judgmental equality, it is often nevertheless provable as a \emph{propositional} equality from the other rules for the type.
    In this case we call it a \define{propositional uniqueness principle}.
    \indexdef{uniqueness!principle, propositional}%
    \indexsee{propositional!uniqueness principle}{uniqueness principle, propositional}%
    (In later chapters we will also occasionally encounter \emph{propositional computation rules}.)
    \indexdef{computation rule!propositional}%
  \end{enumerate}
The inference rules in \autoref{sec:syntax-more-formally} are organized and named accordingly; see, for example, \autoref{sec:more-formal-pi}, where each possibility is realized.
\end{rmk}

The way to construct pairs is obvious: given $a:A$ and $b:B$, we may form $(a,b):A\times B$.
Similarly, there is a unique way to construct elements of $\unit$, namely we have $\ttt:\unit$.
We expect that ``every element of $A\times B$ is a pair'', which is the uniqueness principle for products; we do not assert this as a rule of type theory, but we will prove it later on as a propositional equality.

Now, how can we \emph{use} pairs, i.e.\ how can we define functions out of a product type?
Let us first consider the definition of a non-dependent function $f : A\times B \to C$.
Since we intend the only elements of $A\times B$ to be pairs, we expect to be able to define such a function by prescribing the result
when $f$ is applied to a pair $\tup{a}{b}$.
We can prescribe these results by providing a function $g : A \to B \to C$.
Thus, we introduce a new rule (the elimination rule for products), which says that for any such $g$, we can define a function $f : A\times B \to C$ by
\[ f(\tup{a}{b}) \defeq g(a)(b). \]
We avoid writing $g(a,b)$ here, in order to emphasize that $g$ is not a function on a product.
(However, later on in the book we will often write $g(a,b)$ both for functions on a product and for curried functions of two variables.)
This defining equation is the computation rule for product types\index{computation rule!for product types}.

Note that in set theory, we would justify the above definition of $f$ by the fact that every element of $A\times B$ is a pair, so that it suffices to define $f$ on pairs.
By contrast, type theory reverses the situation: we assume that a function on $A\times B$ is well-defined as soon as we specify its values on tuples, and from this (or more precisely, from its more general version for dependent functions, below) we will be able to \emph{prove} that every element of $A\times B$ is a pair.
From a category-theoretic perspective, we can say that we define the product $A\times B$ to be left adjoint to the ``exponential'' $B\to C$, which we have already introduced.

As an example, we can derive the \define{projection}
\indexsee{function!projection}{projection}%
\indexsee{component, of a pair}{projection}%
\indexdef{projection!from cartesian product type}%
functions
\symlabel{defn:proj}%
\begin{align*}
  \fst & :  A \times B \to A \\
  \snd & :  A \times B \to B
\end{align*}
with the defining equations 
\begin{align*}
  \fst(\tup{a}{b}) & \defeq  a \\
  \snd(\tup{a}{b}) & \defeq  b.
\end{align*}
\symlabel{defn:recursor-times}%
Rather than invoking this principle of function definition every time we want to define a function, an alternative approach is to invoke it once, in a universal case, and then simply apply the resulting function in all other cases.
That is, we may define a function of type
\begin{equation}
  \rec{A\times B} : \prd{C:\UU}(A \to B \to C) \to A \times B \to C
\end{equation}
with the defining equation
\[\rec{A\times B}(C,g,\tup{a}{b}) \defeq g(a)(b). \]
Then instead of defining functions such as $\fst$ and $\snd$ directly by a defining equation, we could  define
\begin{align*}
  \fst &\defeq \rec{A\times B}(A, \lam{a}{b} a)\\
  \snd &\defeq \rec{A\times B}(B, \lam{a}{b} b).
\end{align*}
We refer to the function $\rec{A\times B}$ as the \define{recursor}
\indexsee{recursor}{recursion principle}%
for product types.  The name ``recursor'' is a bit unfortunate here, since no recursion is taking place.  It comes from the fact that product types are a degenerate example of a general framework for inductive types, and for types such as the natural numbers, the recursor will actually be recursive.  We may also speak of the \define{recursion principle} for cartesian products, meaning the fact that we can define a function $f:A\times B\to C$ as above by giving its value on pairs.
\index{recursion principle!for cartesian product}%

We leave it as a simple exercise to show that the recursor can be
derived from the projections and vice versa.

\symlabel{defn:recursor-unit}%
We also have a recursor for the unit type:
\[\rec{\unit} : \prd{C:\UU}C \to \unit \to C\]
with the defining equation
\[ \rec{\unit}(C,c,\ttt) \defeq c. \]
Although we include it to maintain the pattern of type definitions, the recursor for $\unit$ is completely useless,
because we could have defined such a function directly
by simply ignoring the argument of type $\unit$.

To be able to define \emph{dependent} functions over the product type, we have
to generalize the recursor. Given $C: A \times B \to \UU$, we may
define a function $f : \prd{x : A \times B} C(x)$ by providing a
function
\narrowequation{
 g : \prd{x:A}\prd{y:B} C(\tup{x}{y})
}
with defining equation
\[ f(\tup x y) \defeq g(x)(y). \] 
For example, in this way we can prove the propositional uniqueness principle, which says that every element of $A\times B$ is equal to a pair.
\index{uniqueness!principle, propositional!for product types}%
Specifically, we can construct a function
\[ \uppt : \prd{x:A \times B} (\id[A\times B]{\tup{\fst {(x)}}{\snd {(x)}}}{x}). \]
Here we are using the identity type, which we are going to introduce below in \autoref{sec:identity-types}.
However, all we need to know now is that there is a reflexivity element $\refl{x} : \id[A]{x}{x}$ for any $x:A$.
Given this, we can define
\[ \uppt(\tup{a}{b}) \defeq \refl{\tup{a}{b}}. \]
This construction works, because in the case that $x \defeq \tup{a}{b}$ we can 
calculate 
\[ \tup{\fst(\tup{a}{b})}{\snd{(\tup{a}{b})}} \jdeq \tup{a}{b} \]
using the defining equations for the projections. Therefore,
\[ \refl{\tup{a}{b}} : \id{\tup{\fst(\tup{a}{b})}{\snd{(\tup{a}{b})}}}{\tup{a}{b}} \]
is well-typed, since both sides of the equality are judgmentally equal.

More generally, the ability to define dependent functions in this way means that to prove a property for all elements of a product, it is enough 
to prove it for its canonical elements, the tuples.
When we come to inductive types such as the natural numbers, the analogous property will be the ability to write proofs by induction.
Thus, if we do as we did above and apply this principle once in the universal case, we call the resulting function \define{induction} for product types: given $A,B : \UU$ we have
\symlabel{defn:induction-times}%
\[ \ind{A\times B} : \prd{C:A \times B \to \UU}
\Parens{\prd{x:A}{y:B} C(\tup{x}{y})} \to \prd{x:A \times B} C(x) \]
with the defining equation 
\[ \ind{A\times B}(C,g,\tup{a}{b}) \defeq g(a)(b). \]
Similarly, we may speak of a dependent function defined on pairs being obtained from the \define{induction principle}
\index{induction principle}%
\index{induction principle!for product}%
of the cartesian product.
It is easy to see that the recursor is just the special case of induction
in the case that the family $C$ is constant.
Because induction describes how to use an element of the product type, induction is also called the \define{(dependent) eliminator},
\indexsee{eliminator!of inductive type!dependent}{induction principle}%
and recursion the \define{non-dependent eliminator}.
\indexsee{eliminator!of inductive type!non-dependent}{recursion principle}%
\indexsee{non-dependent eliminator}{recursion principle}%
\indexsee{dependent eliminator}{induction principle}%


Induction for the unit type turns out to be more useful than the
recursor: 
\symlabel{defn:induction-unit}%
\[ \ind{\unit} : \prd{C:\unit \to \UU} C(\ttt) \to \prd{x:\unit}C(x)\]
with the defining equation
\[ \ind{\unit}(C,c,\ttt) \defeq c. \]
Induction enables us to prove the propositional uniqueness principle for $\unit$, which asserts that its only inhabitant is $\ttt$.
That is, we can construct
\[\un : \prd{x:\unit} \id{x}{\ttt} \]
by using the defining equations
\[\un(\ttt) \defeq \refl{\ttt} \]
or equivalently by using induction:
\[\un \defeq \ind{\unit}(\lamu{x:\unit} \id{x}{\ttt},\refl{\ttt}). \]

\index{type!product|)}%
\index{type!unit|)}%

\section{Dependent pair types (\texorpdfstring{$\Sigma$}{Σ}-types)}
\label{sec:sigma-types}

\index{type!dependent pair|(defstyle}%
\indexsee{type!dependent sum}{type, dependent pair}%
\indexsee{type!Sigma-@$\Sigma$-}{type, dependent pair}%
\indexsee{Sigma-type@$\Sigma$-type}{type, dependent pair}%
\indexsee{sum!dependent}{type, dependent pair}%

Just as we generalized function types (\autoref{sec:function-types}) to dependent function types (\autoref{sec:pi-types}), it is often useful to generalize the product types from \autoref{sec:finite-product-types} to allow the type of
the second component of a pair to vary depending on the choice of the first
component. This is called a \define{dependent pair type}, or \define{$\Sigma$-type}, because in set theory it
corresponds to an indexed sum (in the sense of coproduct or
disjoint union) over a given type.

Given a type $A:\UU$ and a family $B : A \to \UU$, the dependent
pair type is written as $\sm{x:A} B(x) : \UU$.
Alternative notations are 
\[ \tsm{x:A} B(x) \hspace{2cm} \dsm{x:A}B(x) \hspace{2cm} \lsm{x:A} B(x). \]
Like other binding constructs such as $\lambda$-abstractions and $\Pi$s, $\Sigma$s automatically scope\index{scope} over the rest of the expression unless delimited, so e.g.\ $\sm{x:A} B(x) \times C(x)$ means $\sm{x:A} (B(x) \times C(x))$.

\symlabel{defn:dependent-pair}%
\indexdef{pair!dependent}%
The way to construct elements of a dependent pair type is by pairing: we have
$\tup{a}{b} : \sm{x:A} B(x)$ given $a:A$ and $b:B(a)$.
If $B$ is constant, then the dependent pair type is the
ordinary cartesian product type:
\[ \Parens{\sm{x:A} B} \jdeq (A \times B).\]
All the constructions on $\Sigma$-types arise as straightforward generalizations of the ones for product types, with dependent functions often replacing non-dependent ones.

For instance, the recursion principle%
\index{recursion principle!for dependent pair type}
says that to define a non-dependent function out of a $\Sigma$-type
$f : (\sm{x:A} B(x)) \to C$, we provide a function 
$g : \prd{x:A} B(x) \to C$, and then we can define $f$ via the defining
equation
\[ f(\tup{a}{b}) \defeq g(a)(b). \]
\indexdef{projection!from dependent pair type}%
For instance, we can derive the first projection from a $\Sigma$-type:
\symlabel{defn:dependent-proj1}%
\begin{equation*}
  \fst : \Parens{\sm{x : A}B(x)} \to A.
\end{equation*}
by the defining equation
\begin{equation*}
  \fst(\tup{a}{b}) \defeq a.
\end{equation*}
However, since the type of the second component of a pair
\narrowequation{
  (a,b):\sm{x:A} B(x)
}
is $B(a)$, the second projection must be a \emph{dependent} function, whose type involves the first projection function:
\symlabel{defn:dependent-proj2}%
\[ \snd : \prd{p:\sm{x : A}B(x)}B(\fst(p)). \]
Thus we need the \emph{induction} principle%
\index{induction principle!for dependent pair type}
for $\Sigma$-types (the ``dependent eliminator'').
This says that to construct a dependent function out of a $\Sigma$-type into a family $C : (\sm{x:A} B(x)) \to \UU$, we need a function
\[ g : \prd{a:A}{b:B(a)} C(\tup{a}{b}). \]
We can then derive a function 
\[ f : \prd{p : \sm{x:A}B(x)} C(p) \]
with  defining equation\index{computation rule!for dependent pair type}
\[ f(\tup{a}{b}) \defeq g(a)(b).\]
Applying this with $C(p)\defeq B(\fst(p))$, we can define
\narrowequation{
\snd : \prd{p:\sm{x : A}B(x)}B(\fst(p))
}
with the obvious equation
\[ \snd(\tup{a}{b})  \defeq  b. \]
To convince ourselves that this is correct, we note that $B (\fst(\tup{a}{b})) \jdeq B(a)$, using the defining equation for $\fst$, and
indeed $b : B(a)$.

We can package the recursion and induction principles into the recursor for $\Sigma$:
\symlabel{defn:recursor-sm}%
\[ \rec{\sm{x:A}B(x)} : \dprd{C:\UU}\Parens{\tprd{x:A} B(x) \to C} \to
\Parens{\tsm{x:A}B(x)} \to C \]
with the defining equation
\[ \rec{\sm{x:A}B(x)}(C,g,\tup{a}{b}) \defeq g(a)(b) \]
and the corresponding induction operator:
\symlabel{defn:induction-sm}%
\begin{narrowmultline*}
  \ind{\sm{x:A}B(x)} : \narrowbreak
    \dprd{C:(\sm{x:A} B(x)) \to \UU}
    \Parens{\tprd{a:A}{b:B(a)} C(\tup{a}{b})}
    \to \dprd{p : \sm{x:A}B(x)} C(p)
\end{narrowmultline*}
with the defining equation 
\[ \ind{\sm{x:A}B(x)}(C,g,\tup{a}{b}) \defeq g(a)(b). \]
As before, the recursor is the special case of induction
when the family $C$ is constant.

As a further example, consider the following principle, where $A$ and $B$ are types and $R:A\to B\to \UU$.
\[ \ac : \Parens{\tprd{x:A} \tsm{y :B} R(x,y)} \to
\Parens{\tsm{f:A\to B} \tprd{x:A} R(x,f(x))}
\]
We may regard $R$ as a ``proof-relevant relation''
\index{mathematics!proof-relevant}%
between $A$ and $B$, with $R(a,b)$ the type of witnesses for relatedness of $a:A$ and $b:B$.
Then $\ac$ says intuitively that if we have a dependent function $g$ assigning to every $a:A$ a dependent pair $(b,r)$ where $b:B$ and $r:R(a,b)$, then we have a function $f:A\to B$ and a dependent function assigning to every $a:A$ a witness that $R(a,f(a))$.
Our intuition tells us that we can just split up the values of $g$ into their components.
Indeed, using the projections we have just defined, we can define:
\[ \ac(g) \defeq \Parens{\lamu{x:A} \fst(g(x)),\, \lamu{x:A} \snd(g(x))}. \]
To verify that this is well-typed, note that if $g:\prd{x:A} \sm{y :B} R(x,y)$, we have
\begin{align*}
\lamu{x:A} \fst(g(x)) &: A \to  B, \\
\lamu{x:A} \snd(g(x)) &: \tprd{x:A} R(x,\fst(g(x))).
\end{align*}
Moreover, the type $\prd{x:A} R(x,\fst(g(x)))$ is the result of substituting the function $\lamu{x:A} \fst(g(x))$ for $f$ in the family being summed over in the co\-do\-main of \ac:
\[ \tprd{x:A} R(x,\fst(g(x))) \jdeq
\Parens{\lamu{f:A\to B} \tprd{x:A} R(x,f(x))}\big(\lamu{x:A} \fst(g(x))\big). \]
Thus, we have
\[ \Parens{\lamu{x:A} \fst(g(x)),\, \lamu{x:A} \snd(g(x))} : \tsm{f:A\to B} \tprd{x:A} R(x,f(x))\]
as required.

If we read $\Pi$ as ``for all'' and $\Sigma$ as ``there exists'', then the type of the function $\ac$ expresses:
\emph{if for all $x:A$ there is a $y:B$ such that $R(x,y)$, then there is a function $f : A \to B$ such that for all $x:A$ we have $R(x,f(x))$}.
Since this sounds like a version of the axiom of choice, the function \ac has traditionally been called the \define{type-theoretic axiom of choice}, and as we have just shown, it can be proved directly from the rules of type theory, rather than having to be taken as an axiom.
\index{axiom!of choice!type-theoretic}%
However, note that no choice is actually involved, since the choices have already been given to us in the premise: all we have to do is take it apart into two functions: one representing the choice and the other its correctness.
In \autoref{sec:axiom-choice} we will give another formulation of an ``axiom of choice'' which is closer to the usual one.

Dependent pair types are often used to define types of mathematical structures, which commonly consist of several dependent pieces of data.
To take a simple example, suppose we want to define a \define{magma}\indexdef{magma} to be a type $A$ together with a binary operation $m:A\to A\to A$.
The precise meaning of the phrase ``together with''\index{together with} (and the synonymous ``equipped with''\index{equipped with}) is that ``a magma'' is a \emph{pair} $(A,m)$ consisting of a type $A:\UU$ and an operation $m:A\to A\to A$.
Since the type $A\to A\to A$ of the second component $m$ of this pair depends on its first component $A$, such pairs belong to a dependent pair type.
Thus, the definition ``a magma is a type $A$ together with a binary operation $m:A\to A\to A$'' should be read as defining \emph{the type of magmas} to be
\[ \mathsf{Magma} \defeq \sm{A:\UU} (A\to A\to A). \]
Given a magma, we extract its underlying type (its ``carrier''\index{carrier}) with the first projection $\proj1$, and its operation with the second projection $\proj2$.
Of course, structures built from more than two pieces of data require iterated pair types, which may be only partially dependent; for instance the type of pointed magmas (magmas $(A,m)$ equipped with a basepoint $e:A$) is
\[ \mathsf{PointedMagma} \defeq \sm{A:\UU} (A\to A\to A) \times A. \]
We generally also want to impose axioms on such a structure, e.g.\ to make a pointed magma into a monoid or a group.
This can also be done using $\Sigma$-types; see \autoref{sec:pat}.

In the rest of the book, we will sometimes make definitions of this sort explicit, but eventually we trust the reader to translate them from English into $\Sigma$-types.
We also generally follow the common mathematical practice of using the same letter for a structure of this sort and for its carrier (which amounts to leaving the appropriate projection function implicit in the notation): that is, we will speak of a magma $A$ with its operation $m:A\to A\to A$.

Note that the canonical elements of $\mathsf{PointedMagma}$ are of the form $(A,(m,e))$ where $A:\UU$, $m:A\to A\to A$, and $e:A$.
Because of the frequency with which iterated $\Sigma$-types of this sort arise, we use the usual notation of ordered triples, quadruples and so on to stand for nested pairs (possibly dependent) associating to the right.
That is, we have $(x,y,z) \defeq (x,(y,z))$ and $(x,y,z,w)\defeq (x,(y,(z,w)))$, etc.

\index{type!dependent pair|)}%

\section{Coproduct types}
\label{sec:coproduct-types}

Given $A,B:\UU$, we introduce their \define{coproduct} type $A+B:\UU$.
\indexsee{coproduct}{type, coproduct}%
\index{type!coproduct|(defstyle}%
\indexsee{disjoint!sum}{type, coproduct}%
\indexsee{disjoint!union}{type, coproduct}%
\indexsee{sum!disjoint}{type, coproduct}%
\indexsee{union!disjoint}{type, coproduct}%
This corresponds to the \emph{disjoint union} in set theory, and we may also use that name for it.
In type theory, as was the case with functions and products, the coproduct must be a fundamental construction, since there is no previously given notion of ``union of types''.
We also introduce a
nullary version: the \define{empty type $\emptyt:\UU$}.
\indexsee{nullary!coproduct}{type, empty}%
\indexsee{empty type}{type, empty}%
\index{type!empty|(defstyle}%

There are two ways to construct elements of $A+B$, either as $\inl(a) : A+B$ for $a:A$, or as
$\inr(b):A+B$ for $b:B$. There are no ways to construct elements of the empty type. 

\index{recursion principle!for coproduct}
To construct a non-dependent function $f : A+B \to C$, we need 
functions $g_0 : A \to C$ and $g_1 : B \to C$. Then $f$ is defined
via the defining equations
\begin{align*}
  f(\inl(a)) &\defeq g_0(a), \\
  f(\inr(b)) &\defeq g_1(b).
\end{align*}
That is, the function $f$ is defined by \define{case analysis}.
\indexdef{case analysis}%
As before, we can derive the recursor:
\symlabel{defn:recursor-plus}%
\[ \rec{A+B} : \dprd{C:\UU}(A \to C) \to (B\to C) \to A+B \to C\]
with the defining equations
\begin{align*}
\rec{A+B}(C,g_0,g_1,\inl(a)) &\defeq g_0(a), \\
\rec{A+B}(C,g_0,g_1,\inr(b)) &\defeq g_1(b).
\end{align*}

\index{recursion principle!for empty type}
We can always construct a function $f : \emptyt \to C$ without
having to give any defining equations, because there are no elements of \emptyt on which to define $f$.
Thus, the recursor for $\emptyt$ is
\symlabel{defn:recursor-emptyt}%
\[\rec{\emptyt} : \tprd{C:\UU} \emptyt \to C,\]
which constructs the canonical function from the empty type to any other type.
Logically, it corresponds to the principle \textit{ex falso quodlibet}.
\index{ex falso quodlibet@\textit{ex falso quodlibet}}

\index{induction principle!for coproduct}
To construct a dependent function $f:\prd{x:A+B}C(x)$ out of a coproduct, we assume as given the family 
$C: (A + B) \to \UU$, and 
require 
\begin{align*}
  g_0 &: \prd{a:A} C(\inl(a)), \\
  g_1 &: \prd{b:B} C(\inr(b)).
\end{align*}
This yields $f$ with the defining equations:\index{computation rule!for coproduct type}
\begin{align*}
  f(\inl(a)) &\defeq g_0(a), \\
  f(\inr(b)) &\defeq g_1(b).
\end{align*}
We package this scheme into the induction principle for coproducts:
\symlabel{defn:induction-plus}%
\begin{narrowmultline*}
  \ind{A+B} :
  \dprd{C: (A + B) \to \UU}
  \Parens{\tprd{a:A} C(\inl(a))} \to \narrowbreak
  \Parens{\tprd{b:B} C(\inr(b))} \to \tprd{x:A+B}C(x). 
\end{narrowmultline*}
As before, the recursor arises in the case that the family $C$ is
constant. 

\index{induction principle!for empty type}
The induction principle for the empty type
\symlabel{defn:induction-emptyt}%
\[ \ind{\emptyt} : \prd{C:\emptyt \to \UU}{z:\emptyt} C(z) \]
gives us a way to define a trivial dependent function out of the
empty type. 

\index{type!coproduct|)}%
\index{type!empty|)}%

\section{The type of booleans}
\label{sec:type-booleans}

\indexsee{boolean!type of}{type of booleans}%
\index{type!of booleans|(defstyle}%
The type of booleans $\bool:\UU$ is intended to have exactly two elements 
$\bfalse,\btrue : \bool$. It is clear that we could construct this
type out of coproduct
\index{type!coproduct}%
and unit\index{type!unit} types as $\unit + \unit$. However,
since it is used frequently, we give the explicit rules here.
Indeed, we are going to observe that we can also go the other way
and derive binary coproducts from $\Sigma$-types and $\bool$.

\index{recursion principle!for type of booleans}
To derive a function $f : \bool \to C$ we need $c_0,c_1 : C$ and
add the defining equations
\begin{align*}
  f(\bfalse) &\defeq c_0, \\
  f(\btrue)  &\defeq c_1.
\end{align*}
The recursor corresponds to the if-then-else construct in
functional programming:
\symlabel{defn:recursor-bool}%
\[ \rec{\bool} : \prd{C:\UU}  C \to C \to \bool \to C \]
with the defining equations
\begin{align*}
  \rec{\bool}(C,c_0,c_1,\bfalse) &\defeq c_0, \\
  \rec{\bool}(C,c_0,c_1,\btrue)  &\defeq c_1.
\end{align*}

\index{induction principle!for type of booleans}
Given $C : \bool \to \UU$, to derive a dependent function 
$f : \prd{x:\bool}C(x)$ we need $c_0:C(\bfalse)$ and $c_1 : C(\btrue)$, in which case we can give the defining equations
\begin{align*}
  f(\bfalse) &\defeq c_0, \\
  f(\btrue)  &\defeq c_1.
\end{align*}
We package this up into the induction principle
\symlabel{defn:induction-bool}%
\[ \ind{\bool} : \dprd{C:\bool \to \UU}  C(\bfalse) \to C(\btrue)
\to \tprd{x:\bool} C(x) \]
with the defining equations
\begin{align*}
  \ind{\bool}(C,c_0,c_1,\bfalse) &\defeq c_0, \\
  \ind{\bool}(C,c_0,c_1,\btrue)  &\defeq c_1.
\end{align*}

As an example, using the induction principle we can prove that, as we expect, every element of $\bool$ is either $\btrue$ or $\bfalse$.
As before, we use the equality types which we have not yet introduced, but we need only the fact that everything is equal to itself by $\refl{x}:x=x$.

\begin{thm}\label{thm:allbool-trueorfalse}
  We have
  \[ \prd{x:\bool}(x=\bfalse)+(x=\btrue). \]
\end{thm}
\begin{proof}
  We use the induction principle with $C(x) \defeq (x=\bfalse)+(x=\btrue)$.
  The two inputs are $\inl(\refl{\bfalse}) : C(\bfalse)$ and $\inr(\refl{\btrue}):C(\btrue)$.
\end{proof}

We have remarked that $\Sigma$-types can be regarded as analogous to indexed disjoint unions, while coproducts are binary disjoint unions.
It is natural to expect that a binary disjoint union $A+B$ could be constructed as an indexed one over the two-element type \bool.
For this we need a type family $P:\bool\to\type$ such that $P(\bfalse)\jdeq A$ and $P(\btrue)\jdeq B$.
Indeed, we can obtain such a family precisely by the recursion principle for $\bool$.
\index{type!family of}%
(The ability to define \emph{type families} by induction and recursion, using the fact that the universe $\UU$ is itself a type, is a subtle and important aspect of type theory.)
Thus, we could have defined
\index{type!coproduct}%
\[ A + B \defeq \sm{x:\bool} \rec{\bool}(\UU,A,B,x). \]
with
\begin{align*}
  \inl(a) &\defeq \tup{\bfalse}{a}, \\
  \inr(b) &\defeq \tup{\btrue}{b}.
\end{align*}
We leave it as an exercise to derive the induction principle of a coproduct type from this definition.
(See also \autoref{ex:sum-via-bool,sec:appetizer-univalence}.)

We can apply the same idea to products and $\Pi$-types: we could have defined
\[ A \times B \defeq \prd{x:\bool}\rec{\bool}(\UU,A,B,x) \]
Pairs could then be constructed using induction for \bool:
\[ \tup{a}{b} \defeq \ind{\bool}(\rec{\bool}(\UU,A,B),a,b) \]
while the projections are straightforward applications
\begin{align*}
  \fst(p) &\defeq p(\bfalse), \\
  \snd(p) &\defeq p(\btrue).
\end{align*}
The derivation of the induction principle for binary products defined in this way is a bit more involved, and requires function extensionality, which we will introduce in \autoref{sec:compute-pi}.
Moreover, we do not get the same judgmental equalities; see \autoref{ex:prod-via-bool}.
This is a recurrent issue when encoding one type as another; we will return to it in \autoref{sec:htpy-inductive}. 

We may occasionally refer to the elements $\bfalse$ and $\btrue$ of $\bool$ as ``false'' and ``true'' respectively.
However, note that unlike in classical\index{mathematics!classical} mathematics, we do not use elements of $\bool$ as truth values
\index{value!truth}%
or as propositions.
(Instead we identify propositions with types; see \autoref{sec:pat}.)
In particular, the type $A \to \bool$ is not generally the power set\index{power set} of $A$; it represents only the ``decidable'' subsets of $A$ (see \autoref{cha:logic}).
\index{decidable!subset}%

\index{type!of booleans|)}%

\section{The natural numbers}
\label{sec:inductive-types}

\indexsee{type!of natural numbers}{natural numbers}%
\index{natural numbers|(defstyle}%
\indexsee{number!natural}{natural numbers}%
The rules we have introduced so far do not allow us to construct any infinite types.
The simplest infinite type we can think of (and one which is of course also extremely useful) is the type $\nat : \UU$ of natural numbers.
The elements of $\nat$ are constructed using $0 : \nat$\indexdef{zero} and the successor\indexdef{successor} operation $\suc : \nat \to \nat$.
When denoting natural numbers, we adopt the usual decimal notation $1 \defeq \suc(0)$, $2 \defeq \suc(1)$, $3 \defeq \suc(2)$, \dots.

The essential property of the natural numbers is that we can define functions by recursion and perform proofs by induction --- where now the words ``recursion'' and ``induction'' have a more familiar meaning.
\index{recursion principle!for natural numbers}%
To construct a non-dependent function $f : \nat \to C$ out of the natural numbers by recursion, it is enough to provide a starting point $c_0 : C$ and a ``next step'' function $c_s : \nat \to C \to C$.
This gives rise to $f$ with the defining equations\index{computation rule!for natural numbers}
\begin{align*}
  f(0) &\defeq c_0, \\
  f(\suc(n)) &\defeq c_s(n,f(n)).
\end{align*}
We say that $f$ is defined by \define{primitive recursion}.
\indexdef{primitive!recursion}%
\indexdef{recursion!primitive}%

As an example, we look at how to define a function on natural numbers which doubles its argument.
In this case we have $C\defeq \nat$.
We first need to supply the value of $\dbl(0)$, which is easy: we put $c_0 \defeq 0$.
Next, to compute the value of $\dbl(\suc(n))$ for a natural number $n$, we first compute the value of $\dbl(n)$ and then perform the successor operation twice.
This is captured by the recurrence\index{recurrence} $c_s(n,y) \defeq \suc(\suc(y))$.
Note that the second argument $y$ of $c_s$ stands for the result of the \emph{recursive call}\index{recursive call} $\dbl(n)$.

Defining $\dbl:\nat\to\nat$ by primitive recursion in this way, therefore, we obtain the defining equations:
\begin{align*}
  \dbl(0) &\defeq 0\\
  \dbl(\suc(n)) &\defeq \suc(\suc(\dbl(n))).
\end{align*}
This indeed has the correct computational behavior: for example, we have 
\begin{align*}
  \dbl(2) &\jdeq \dbl(\suc(\suc(0)))\\
  & \jdeq c_s(\suc(0), \dbl(\suc(0))) \\
                 & \jdeq \suc(\suc(\dbl(\suc(0)))) \\
                 & \jdeq \suc(\suc(c_s(0,\dbl(0)))) \\
                 & \jdeq \suc(\suc(\suc(\suc(\dbl(0))))) \\
                 & \jdeq \suc(\suc(\suc(\suc(c_0)))) \\
                 & \jdeq \suc(\suc(\suc(\suc(0))))\\
                 &\jdeq 4.
\end{align*}
We can define multi-variable functions by primitive recursion as well, by currying and allowing $C$ to be a function type.
\indexdef{addition!of natural numbers}
For example, we define addition $\add : \nat \to \nat \to \nat$ with $C \defeq \nat \to \nat$ and the following ``starting point'' and ``next step'' data:
\begin{align*}
  c_0 & : \nat \to \nat \\
  c_0 (n) & \defeq n \\
  c_s & : \nat \to (\nat \to \nat) \to (\nat \to \nat) \\
  c_s(m,g)(n) & \defeq \suc(g(n)).
\end{align*}
We thus obtain $\add : \nat \to \nat \to \nat$ satisfying the definitional equalities
\begin{align*}
  \add(0,n) &\jdeq n \\
  \add(\suc(m),n) &\jdeq \suc(\add(m,n)). 
\end{align*}
As usual, we write $\add(m,n)$ as $m+n$.
The reader is invited to verify that $2+2\jdeq 4$.

As in previous cases, we can package the principle of primitive recursion into a recursor:
\[\rec{\nat}  : \dprd{C:\UU} C \to (\nat \to C \to C) \to \nat \to C \]
with the defining equations
\symlabel{defn:recursor-nat}%
\begin{align*}
\rec{\nat}(C,c_0,c_s,0)  &\defeq c_0, \\
\rec{\nat}(C,c_0,c_s,\suc(n)) &\defeq c_s(n,\rec{\nat}(C,c_0,c_s,n)).
\end{align*}
Using $\rec{\nat}$ we can present $\dbl$ and $\add$ as follows:
\begin{align}
\dbl &\defeq \rec\nat\big(\nat,\, 0,\, \lamu{n:\nat}{y:\nat} \suc(\suc(y))\big) \label{eq:dbl-as-rec}\\
\add &\defeq \rec{\nat}\big(\nat \to \nat,\, \lamu{n:\nat} n,\, \lamu{n:\nat}{g:\nat \to \nat}{m :\nat} \suc(g(m))\big).
\end{align}
Of course, all functions definable only using the primitive recursion principle will be \emph{computable}.
(The presence of higher function types --- that is, functions with other functions as arguments --- does, however, mean we can define more than the usual primitive recursive functions; see e.g.~\autoref{ex:ackermann}.)
This is appropriate in constructive mathematics;
\index{mathematics!constructive}%
in \autoref{sec:intuitionism,sec:axiom-choice} we will see how to augment type theory so that we can define more general mathematical functions.

\index{induction principle!for natural numbers}
We now follow the same approach as for other types, generalizing primitive recursion to dependent functions to obtain an \emph{induction principle}.
Thus, assume as given a family $C : \nat \to \UU$, an element $c_0 : C(0)$, and a function $c_s : \prd{n:\nat} C(n) \to C(\suc(n))$; then we can construct $f : \prd{n:\nat} C(n)$ with the defining equations:\index{computation rule!for natural numbers}
\begin{align*}
  f(0) &\defeq c_0, \\
  f(\suc(n)) &\defeq c_s(n,f(n)).
\end{align*}
We can also package this into a single function
\symlabel{defn:induction-nat}%
\[\ind{\nat}  : \dprd{C:\nat\to \UU} C(0) \to \Parens{\tprd{n : \nat} C(n) \to C(\suc(n))} \to \tprd{n : \nat} C(n) \]
with the defining equations
\begin{align*}
\ind{\nat}(C,c_0,c_s,0)  &\defeq c_0, \\
\ind{\nat}(C,c_0,c_s,\suc(n)) &\defeq c_s(n,\ind{\nat}(C,c_0,c_s,n)).
\end{align*}
Here we finally see the connection to the classical notion of proof by induction.
Recall that in type theory we represent propositions by types, and proving a proposition by inhabiting the corresponding type.
In particular, a \emph{property} of natural numbers is represented by a family of types $P:\nat\to\type$.
From this point of view, the above induction principle says that if we can prove $P(0)$, and if for any $n$ we can prove $P(\suc(n))$ assuming $P(n)$, then we have $P(n)$ for all $n$.
This is, of course, exactly the usual principle of proof by induction on natural numbers.

\index{associativity!of addition!of natural numbers}
As an example, consider how we might represent an explicit proof that $+$ is associative.
(We will not actually write out proofs in this style, but it serves as a useful example for understanding how induction is represented formally in type theory.)
To derive
\[\assoc : \prd{i,j,k:\nat} \id{i + (j + k)}{(i + j) + k}, \]
it is sufficient to supply
\[ \assoc_0 :  \prd{j,k:\nat} \id{0 + (j + k)}{(0+ j) + k} \]
and
\begin{narrowmultline*}
  \assoc_s  : \prd{i:\nat} \left(\prd{j,k:\nat} \id{i + (j + k)}{(i + j) + k}\right)
   \narrowbreak
   \to \prd{j,k:\nat} \id{\suc(i) + (j + k)}{(\suc(i) + j) + k}.
\end{narrowmultline*}
To derive $\assoc_0$, recall that $0+n \jdeq n$, and hence  $0 + (j + k) \jdeq j+k \jdeq (0+ j) + k$.
Hence we can just set
\[ \assoc_0(j,k) \defeq \refl{j+k}. \]
For $\assoc_s$, recall that the definition of $+$ gives $\suc(m)+n \jdeq \suc(m+n)$, and hence 
\begin{align*}
   \suc(i) + (j + k)  &\jdeq \suc(i+(j+k)) \qquad\text{and}\\
   (\suc(i)+j)+k &\jdeq \suc((i+j)+k).
\end{align*}
Thus, the output type of $\assoc_s$ is equivalently $\id{\suc(i+(j+k))}{\suc((i+j)+k)}$.
But its input (the ``inductive hypothesis'')
\index{hypothesis!inductive}%
\index{inductive!hypothesis}%
yields $\id{i+(j+k)}{(i+j)+k}$, so it suffices to invoke the fact that if two natural numbers are equal, then so are their successors.
(We will prove this obvious fact in \autoref{lem:map}, using the induction principle of identity types.)
We call this latter fact
$\apfunc{\suc} : 
(\id[\nat]{m}{n}) \to (\id[\nat]{\suc(m)}{\suc(n)})$, so we can define
\[\assoc_s(i,h,j,k) \defeq \apfunc{\suc}( 
h(j,k)). \]
Putting these together with $\ind{\nat}$, we obtain a proof of associativity.

\index{natural numbers|)}%

\section{Pattern matching and recursion}
\label{sec:pattern-matching}

\index{pattern matching|(defstyle}%
\indexsee{matching}{pattern matching}%
\index{definition!by pattern matching|(}%
The natural numbers introduce an additional subtlety over the types considered up until now.
In the case of coproducts, for instance, we could define a function $f:A+B\to C$ either with the recursor:
\[ f \defeq \rec{A+B}(C, g_0, g_1) \]
or by giving the defining equations:
\begin{align*}
  f(\inl(a)) &\defeq g_0(a)\\
  f(\inr(b)) &\defeq g_1(b).
\end{align*}
To go from the former expression of $f$ to the latter, we simply use the computation rules for the recursor.
Conversely, given any defining equations
\begin{align*}
  f(\inl(a)) &\defeq \Phi_0\\
  f(\inr(b)) &\defeq \Phi_1
\end{align*}
where $\Phi_0$ and $\Phi_1$ are expressions that may involve the variables
\index{variable}%
$a$ and $b$ respectively, we can express these equations equivalently in terms of the recursor by using $\lambda$-abstraction\index{lambda abstraction@$\lambda$-abstraction}:
\[ f\defeq \rec{A+B}(C, \lam{a} \Phi_0, \lam{b} \Phi_1).\]
In the case of the natural numbers, however, the ``defining equations'' of a function such as $\dbl$:
\begin{align}
  \dbl(0) &\defeq 0 \label{eq:dbl0}\\
  \dbl(\suc(n)) &\defeq \suc(\suc(\dbl(n)))\label{eq:dblsuc}
\end{align}
involve \emph{the function $\dbl$ itself} on the right-hand side.
However, we would still like to be able to give these equations, rather than~\eqref{eq:dbl-as-rec}, as the definition of \dbl, since they are much more convenient and readable.
The solution is to read the expression ``$\dbl(n)$'' on the right-hand side of~\eqref{eq:dblsuc} as standing in for the result of the recursive call, which in a definition of the form $\dbl\defeq \rec{\nat}(\nat,c_0,c_s)$ would be the second argument of $c_s$.

More generally, if we have a ``definition'' of a function $f:\nat\to C$ such as
\begin{align*}
  f(0) &\defeq \Phi_0\\
  f(\suc(n)) &\defeq \Phi_s
\end{align*}
where $\Phi_0$ is an expression of type $C$, and $\Phi_s$ is an expression of type $C$ which may involve the variable $n$ and also the symbol ``$f(n)$'', we may translate it to a definition
\[ f \defeq \rec{\nat}(C,\,\Phi_0,\,\lam{n}{r} \Phi_s') \]
where $\Phi_s'$ is obtained from $\Phi_s$ by replacing all occurrences of ``$f(n)$'' by the new variable $r$.

This style of defining functions by recursion (or, more generally, dependent functions by induction) is so convenient that we frequently adopt it.
It is called definition by \define{pattern matching}.
Of course, it is very similar to how a computer programmer may define a recursive function with a body that literally contains recursive calls to itself.
However, unlike the programmer, we are restricted in what sort of recursive calls we can make: in order for such a definition to be re-expressible using the recursion principle, the function $f$ being defined can only appear in the body of $f(\suc(n))$ as part of the composite symbol ``$f(n)$''.
Otherwise, we could write nonsense functions such as
\begin{align*}
  f(0)&\defeq 0\\
  f(\suc(n)) &\defeq f(\suc(\suc(n))).
\end{align*}
If a programmer wrote such a function, it would simply call itself forever on any positive input, going into an infinite loop and never returning a value.
In mathematics, however, to be worthy of the name, a \emph{function} must always associate a unique output value to every input value, so this would be unacceptable.

This point will be even more important when we introduce more complicated inductive types in \autoref{cha:induction,cha:hits,cha:real-numbers}.
Whenever we introduce a new kind of inductive definition, we always begin by deriving its induction principle.
Only then do we introduce an appropriate sort of ``pattern matching'' which can be justified as a shorthand for the induction principle.

\index{pattern matching|)}%
\index{definition!by pattern matching|)}%

\section{Propositions as types}
\label{sec:pat}

\index{proposition!as types|(defstyle}%
\index{logic!propositions as types|(}%
As mentioned in the introduction, to show that a proposition is true in type theory corresponds to exhibiting an element of the type corresponding to that proposition.
\index{evidence, of the truth of a proposition}%
\index{witness!to the truth of a proposition}%
\index{proof|(}
We regard the elements of this type as \emph{evidence} or \emph{witnesses} that the proposition is true. (They are sometimes even called \emph{proofs}, but this terminology can be misleading, so we generally avoid it.)
In general, however, we will not construct witnesses explicitly; instead we present the proofs in ordinary mathematical prose, in such a way that they could be translated into an element of a type.
This is no different from reasoning in classical set theory, where we don't expect to see an explicit derivation using the rules of predicate logic and the axioms of set theory.

However, the type-theoretic perspective on proofs is nevertheless different in important ways.
The basic principle of the logic of type theory is that a proposition is not merely true or false, but rather can be seen as the collection of all possible witnesses of its truth.
Under this conception, proofs are not just the means by which mathematics is communicated, but rather are mathematical objects in their own right, on a par with more familiar objects such as numbers, mappings, groups, and so on.
Thus, since types classify the available mathematical objects and govern how they interact, propositions are nothing but special  types --- namely, types whose elements are proofs.

\index{propositional!logic}%
\index{logic!propositional}%
The basic observation which makes this identification feasible is that we have the following natural correspondence between \emph{logical} operations on propositions, expressed in English, and \emph{type-theoretic} operations on their corresponding types of witnesses.
\index{false}%
\index{true}%
\index{conjunction}%
\index{disjunction}%
\index{implication}%
\begin{center}
\medskip
\begin{tabular}{ll}
  \toprule
  English & Type Theory\\
  \midrule
  True & $\unit$ \\
  False & $\emptyt$ \\
  $A$ and $B$ & $A \times B$ \\
  $A$ or $B$ & $A + B$ \\
  If $A$ then $B$ & $A \to B$ \\
  $A$ if and only if $B$ & $(A \to B) \times (B \to A)$ \\
  Not $A$ &  $A \to \emptyt$ \\
  \bottomrule
\end{tabular}
\medskip
\end{center}

The point of the correspondence is that in each case, the rules for constructing and using elements of the type on the right correspond to the rules for reasoning about the proposition on the left.
For instance, the basic way to prove a statement of the form ``$A$ and $B$'' is to prove $A$ and also prove $B$, while the basic way to construct an element of $A\times B$ is as a pair $(a,b)$, where $a$ is an element (or witness) of $A$ and $b$ is an element (or witness) of $B$.
And if we want to use ``$A$ and $B$'' to prove something else, we are free to use both $A$ and $B$ in doing so, analogously to how the induction principle for $A\times B$ allows us to construct a function out of it by using elements of $A$ and of $B$.

Similarly, the basic way to prove an implication\index{implication} ``if $A$ then $B$'' is to assume $A$ and prove $B$, while the basic way to construct an element of $A\to B$ is to give an expression which denotes an element (witness) of $B$ which may involve an unspecified variable element (witness) of type $A$.
And the basic way to use an implication ``if $A$ then $B$'' is deduce $B$ if we know $A$, analogously to how we can apply a function $f:A\to B$ to an element of $A$ to produce an element of $B$.
We strongly encourage the reader to do the exercise of verifying that the rules governing the other type constructors translate sensibly into logic.

Of special note is that the empty type $\emptyt$ corresponds to falsity.\index{false}
When speaking logically, we refer to an inhabitant of $\emptyt$ as a \define{contradiction}:
\indexdef{contradiction}%
thus there is no way to prove a contradiction,%
\footnote{More precisely, there is no \emph{basic} way to prove a contradiction, i.e.\ \emptyt has no constructors.
If our type theory were inconsistent, then there would be some more complicated way to construct an element of \emptyt.}
while from a contradiction anything can be derived.
We also define the \define{negation}
\indexdef{negation}%
of a type $A$ as
\begin{equation*}
  \neg A \ \defeq\ A \to \emptyt.
\end{equation*}
Thus, a witness of $\neg A$ is a function $A \to \emptyt$, which we may construct by assuming $x : A$ and deriving an element of~$\emptyt$.
\index{proof!by contradiction}%
\index{logic!constructive vs classical}
Note that although the logic we obtain is ``constructive'', as discussed in the introduction, this sort of ``proof by contradiction'' (assume $A$ and derive a contradiction, concluding $\neg A$) is perfectly valid constructively: it is simply invoking the \emph{meaning} of ``negation''.
The sort of ``proof by contradiction'' which is disallowed is to assume $\neg A$ and derive a contradiction as a way of proving $A$.
Constructively, such an argument would only allow us to conclude $\neg\neg A$, and the reader can verify that there is no obvious way to get from $\neg\neg A$ (that is, from $(A\to \emptyt)\to\emptyt$) to $A$.

\mentalpause

The above translation of logical connectives into type-forming operations is referred to as \define{propositions as types}: it gives us a way to translate propositions and their proofs, written in English, into types and their elements.
For example, suppose we want to prove the following tautology (one of ``de Morgan's laws''):
\index{law!de Morgan's|(}%
\index{de Morgan's laws|(}%
\begin{equation}\label{eq:tautology1}
  \text{\emph{``If not $A$ and not $B$, then not ($A$ or $B$)''}.}
\end{equation}
An ordinary English proof of this fact might go as follows.
\begin{quote}
  Suppose not $A$ and not $B$, and also suppose $A$ or $B$; we will derive a contradiction.
  There are two cases.
  If $A$ holds, then since not $A$, we have a contradiction.
  Similarly, if $B$ holds, then since not $B$, we also have a contradiction.
  Thus we have a contradiction in either case, so not ($A$ or $B$).
\end{quote}
Now, the type corresponding to our tautology~\eqref{eq:tautology1}, according to the rules given above, is
\begin{equation}\label{eq:tautology2}
  (A\to \emptyt) \times (B\to\emptyt) \to (A+B\to\emptyt)
\end{equation}
so we should be able to translate the above proof into an element of this type.

As an example of how such a translation works, let us describe how a mathematician reading the above English proof might simultaneously construct, in his or her head, an element of~\eqref{eq:tautology2}.
The introductory phrase ``Suppose not $A$ and not $B$'' translates into defining a function, with an implicit application of the recursion principle for the cartesian product in its domain $(A\to\emptyt)\times (B\to\emptyt)$.
This introduces unnamed variables
\index{variable}%
(hypotheses)
\index{hypothesis}%
of types $A\to\emptyt$ and $B\to\emptyt$.
When translating into type theory, we have to give these variables names; let us call them $x$ and $y$.
At this point our partial definition of an element of~\eqref{eq:tautology2} can be written as
\[ f((x,y)) \defeq\; \Box\;:A+B\to\emptyt \]
with a ``hole'' $\Box$ of type $A+B\to\emptyt$ indicating what remains to be done.
(We could equivalently write $f \defeq \rec{(A\to\emptyt)\times (B\to\emptyt)}(A+B\to\emptyt,\lam{x}{y} \Box)$, using the recursor instead of pattern matching.)
The next phrase ``also suppose $A$ or $B$; we will derive a contradiction'' indicates filling this hole by a function definition, introducing another unnamed hypothesis $z:A+B$, leading to the proof state:
\[ f((x,y))(z) \defeq \;\Box\; :\emptyt \]
Now saying ``there are two cases'' indicates a case split, i.e.\ an application of the recursion principle for the coproduct $A+B$.
If we write this using the recursor, it would be
\[ f((x,y))(z) \defeq \rec{A+B}(\emptyt,\lam{a} \Box,\lam{b}\Box,z) \]
while if we write it using pattern matching, it would be
\begin{align*}
  f((x,y))(\inl(a)) &\defeq \;\Box\;:\emptyt\\
  f((x,y))(\inr(b)) &\defeq \;\Box\;:\emptyt.
\end{align*}
Note that in both cases we now have two ``holes'' of type $\emptyt$ to fill in, corresponding to the two cases where we have to derive a contradiction.
Finally, the conclusion of a contradiction from $a:A$ and $x:A\to\emptyt$ is simply application of the function $x$ to $a$, and similarly in the other case.
\index{application!of hypothesis or theorem}%
(Note the convenient coincidence of the phrase ``applying a function'' with that of ``applying a hypothesis'' or theorem.)
Thus our eventual definition is
\begin{align*}
  f((x,y))(\inl(a)) &\defeq x(a)\\
  f((x,y))(\inr(b)) &\defeq y(b).
\end{align*}

As an exercise, you should verify 
the converse tautology \emph{``If not ($A$ or $B$), then  (not $A$) and (not $B$)}'' by exhibiting an element of 
\[ ((A + B) \to \emptyt) \to (A \to \emptyt) \times (B \to \emptyt), \]
for any types $A$ and $B$, using the rules we have just introduced.

\index{logic!classical vs constructive|(}
However, not all classical\index{mathematics!classical} tautologies hold under this interpretation.
For example, the rule 
\emph{``If not ($A$ and $B$), then (not $A$) or (not $B$)''} is not valid: we cannot, in general, construct an element of the corresponding type
\[ ((A \times B) \to \emptyt) \to (A \to \emptyt) + (B \to \emptyt).\]
This reflects the fact that the ``natural'' propositions-as-types logic of type theory is \emph{constructive}.
This means that it does not include certain classical principles, such as the law of excluded middle (\LEM{})\index{excluded middle}
or proof by contradiction,\index{proof!by contradiction}
and others which depend on them, such as this instance of de Morgan's law.
\index{law!de Morgan's|)}%
\index{de Morgan's laws|)}%

Philosophically, constructive logic is so-called because it confines itself to constructions that can be carried out \emph{effectively}, which is to say those with a computational meaning.
Without being too precise, this means there is some sort of algorithm\index{algorithm} specifying, step-by-step, how to build an object (and, as a special case, how to see that a theorem is true).
This requires omission of \LEM{}, since there is no \emph{effective}\index{effective!procedure} procedure for deciding whether a proposition is true or false.

The constructivity of type-theoretic logic means it has an intrinsic computational meaning, which is of interest to computer scientists.
It also means that type theory provides \emph{axiomatic freedom}.\index{axiomatic freedom}
For example, while by default there is no construction witnessing \LEM{}, the logic is still compatible with the existence of one (see \autoref{sec:intuitionism}).
Thus, because type theory does not \emph{deny} \LEM{}, we may consistently add it as an assumption, and work conventionally without restriction.
In this respect, type theory enriches, rather than constrains, conventional mathematical practice.

We encourage the reader who is unfamiliar with constructive logic to work through some more examples as a means of getting familiar with it.
See \autoref{ex:tautologies,ex:not-not-lem} for some suggestions.
\index{logic!classical vs constructive|)}

\mentalpause

So far we have discussed only propositional logic.
\index{quantifier}%
\index{quantifier!existential}%
\index{quantifier!universal}%
\index{predicate!logic}%
\index{logic!predicate}%
Now we consider \emph{predicate} logic, where in addition to logical connectives like ``and'' and ``or'' we have quantifiers ``there exists'' and ``for all''.
In this case, types play a dual role: they serve as propositions and also as types in the conventional sense, i.e., domains we quantify over.
A predicate over a type $A$ is represented as a family $P : A \to \UU$, assigning to every element $a : A$ a type $P(a)$ corresponding to the proposition that $P$ holds for $a$. We now extend the above translation with an explanation of the quantifiers:
\begin{center}
  \medskip
  \begin{tabular}{ll}
    \toprule
    English & Type Theory\\
    \midrule
    For all $x:A$, $P(x)$ holds & $\prd{x:A} P(x)$ \\
    There exists $x:A$ such that $P(x)$ & $\sm{x:A}$ $P(x)$ \\
    \bottomrule
  \end{tabular}
  \medskip
\end{center}
As before, we can show that tautologies of (constructive) predicate logic translate into inhabited types.
For example, \emph{If for all $x:A$, $P(x)$ and $Q(x)$ then (for all $x:A$, $P(x)$) and (for all $x:A$, $Q(x)$)} translates to
\[ (\tprd{x:A} P(x) \times Q(x)) \to (\tprd{x:A} P(x)) \times (\tprd{x:A} Q(x)). \]
An informal proof of this tautology might go as follows:
\begin{quote}
  Suppose for all $x$, $P(x)$ and $Q(x)$.
  First, we suppose given $x$ and prove $P(x)$.
  By assumption, we have $P(x)$ and $Q(x)$, and hence we have $P(x)$.
  Second, we suppose given $x$ and prove $Q(x)$.
  Again by assumption, we have $P(x)$ and $Q(x)$, and hence we have $Q(x)$.
\end{quote}
The first sentence begins defining an implication as a function, by introducing a witness for its hypothesis:\index{hypothesis}
\[ f(p) \defeq \;\Box\; : (\tprd{x:A} P(x)) \times (\tprd{x:A} Q(x)). \]
At this point there is an implicit use of the pairing constructor to produce an element of a product type, which is somewhat signposted in this example by the words ``first'' and ``second'':
\[ f(p) \defeq \Big( \;\Box\; : \tprd{x:A} P(x) \;,\; \Box\; : \tprd{x:A}Q(x) \;\Big). \]
The phrase ``we suppose given $x$ and prove $P(x)$'' now indicates defining a \emph{dependent} function in the usual way, introducing a variable
\index{variable}%
for its input.
Since this is inside a pairing constructor, it is natural to write it as a $\lambda$-abstraction\index{lambda abstraction@$\lambda$-abstraction}:
\[ f(p) \defeq \Big( \; \lam{x} \;\big(\Box\; : P(x)\big) \;,\; \Box\; : \tprd{x:A}Q(x) \;\Big). \]
Now ``we have $P(x)$ and $Q(x)$'' invokes the hypothesis, obtaining $p(x) : P(x)\times Q(x)$, and ``hence we have $P(x)$'' implicitly applies the appropriate projection:
\[ f(p) \defeq \Big( \; \lam{x} \proj1(p(x))  \;,\; \Box\; : \tprd{x:A}Q(x) \;\Big). \]
The next two sentences fill the other hole in the obvious way:
\[ f(p) \defeq \Big( \; \lam{x} \proj1(p(x))  \;,\; \lam{x} \proj2(p(x)) \; \Big). \]
Of course, the English proofs we have been using as examples are much more verbose than those that mathematicians usually use in practice; they are more like the sort of language one uses in an ``introduction to proofs'' class.
The practicing mathematician has learned to fill in the gaps, so in practice we can omit plenty of details, and we will generally do so.
The criterion of validity for proofs, however, is always that they can be translated back into the construction of an element of the corresponding type.

\symlabel{leq-nat}%
As a more concrete example, consider how to define inequalities of natural numbers.
One natural definition is that $n\le m$ if there exists a $k:\nat$ such that $n+k=m$.
(This uses again the identity types that we will introduce in the next section, but we will not need very much about them.)
Under the propositions-as-types translation, this would yield:
\[ (n\le m) \defeq \sm{k:\nat} (\id{n+k}{m}). \]
The reader is invited to prove the familiar properties of $\le$ from this definition.
For strict inequality, there are a couple of natural choices, such as
\[ (n<m) \defeq \sm{k:\nat} (\id{n+\suc(k)}{m}) \]
or
\[ (n<m) \defeq (n\le m) \times \neg(\id{n}{m}). \]
The former is more natural in constructive mathematics, but in this case it is actually equivalent to the latter, since $\nat$ has ``decidable equality'' (see \autoref{sec:intuitionism,prop:nat-is-set}).
\index{decidable!equality}%

The representation of propositions as types also allows us to incorporate axioms into the definition of types as mathematical structures using $\Sigma$-types, as discussed in \autoref{sec:sigma-types}.
For example, suppose we want to define a \define{semigroup}\index{semigroup} to be a type $A$ equipped with a binary operation $m:A\to A\to A$ (that is, a magma\index{magma}) and such that for all $x,y,z:A$ we have $m(x,m(y,z)) = m(m(x,y),z)$.
This latter proposition is represented by the type $\prd{x,y,z:A} m(x,m(y,z)) = m(m(x,y),z)$, so the type of semigroups is
\[ \semigroup \defeq \sm{A:\UU}{m:A\to A\to A} \prd{x,y,z:A} m(x,m(y,z)) = m(m(x,y),z). \]
From an inhabitant of this type we can extract the carrier $A$, the operation $m$, and a witness of the axiom, by applying appropriate projections.
We will return to this example in \autoref{sec:equality-of-structures}.

Note also that we can use the universes in type theory to represent ``higher order logic'' --- that is, we can quantify over all propositions or over all predicates.
For example, we can represent the proposition \emph{for all properties $P : A \to \UU$, if $P(a)$ then $P(b)$} as
\[ \prd{P : A \to \UU} P(a) \to P(b) \]
where $A : \UU$ and $a,b : A$.
However, \emph{a priori} this proposition lives in a different, higher, universe than the
propositions we are quantifying over; that is
\[ \Parens{\prd{P : A \to \UU_i} P(a) \to P(b)} : \UU_{i+1}. \]
We will return to this issue in \autoref{subsec:prop-subsets}.

\mentalpause

We have described here a ``proof-relevant''
\index{mathematics!proof-relevant}%
translation of propositions, where the proofs of disjunctions and existential statements carry some information.
For instance, if we have an inhabitant of $A+B$, regarded as a witness of ``$A$ or $B$'', then we know whether it came from $A$ or from $B$.
Similarly, if we have an inhabitant of $\sm{x:A} P(x)$, regarded as a witness of ``there exists $x:A$ such that $P(x)$'', then we know what the element $x$ is (it is the first projection of the given inhabitant).

As a consequence of the proof-relevant nature of this logic, we may have ``$A$ if and only if $B$'' (which, recall, means $(A\to B)\times (B\to A)$), and yet the types $A$ and $B$ exhibit different behavior.
For instance, it is easy to verify that ``$\mathbb{N}$ if and only if $\unit$'', and yet clearly $\mathbb{N}$ and $\unit$ differ in important ways.
The statement ``$\mathbb{N}$ if and only if $\unit$'' tells us only that \emph{when regarded as a mere proposition}, the type $\mathbb{N}$ represents the same proposition as $\unit$ (in this case, the true proposition).
We sometimes express ``$A$ if and only if $B$'' by saying that $A$ and $B$ are \define{logically equivalent}.
\indexdef{logical equivalence}%
\indexdef{equivalence!logical}%
This is to be distinguished from the stronger notion of \emph{equivalence of types} to be introduced in \autoref{sec:basics-equivalences,cha:equivalences}:
although $\mathbb{N}$ and $\unit$ are logically equivalent, they are not equivalent types.

In \autoref{cha:logic} we will introduce a class of types called ``mere propositions'' for which equivalence and logical equivalence coincide.
Using these types, we will introduce a modification to the above-described logic that is sometimes appropriate, in which the additional information contained in disjunctions and existentials is discarded.

Finally, we note that the propositions-as-types correspondence can be viewed in reverse, allowing us to regard any type $A$ as a proposition, which we prove by exhibiting an element of $A$.
Sometimes we will state this proposition as ``$A$ is \define{inhabited}.''
\indexdef{inhabited type}%
\indexsee{type!inhabited}{inhabited type}%
That is, when we say that $A$ is inhabited, we mean that we have given a (particular) element of $A$, but that we are choosing not to give a name to that element.
Similarly, to say that $A$ is \emph{not inhabited} is the same as to give an element of $\neg A$.
In particular, the empty type $\emptyt$ is obviously not inhabited, since $\neg \emptyt \jdeq (\emptyt \to \emptyt)$ is inhabited by $\idfunc[\emptyt]$.\footnote{This should not be confused with the statement that type theory is consistent, which is the \emph{meta-theoretic} claim that it is not possible to obtain an element of $\emptyt$ by following the rules of type theory.\indexfoot{consistency}}

\index{proof|)}%
\index{proposition!as types|)}%
\index{logic!propositions as types|)}%

\section{Identity types}
\label{sec:identity-types}

\index{type!identity|(defstyle}%
\indexsee{identity!type}{type, identity}%
\indexsee{type!equality}{type, identity}%
\indexsee{equality!type}{type, identity}%
While the previous constructions can be seen as generalizations of
standard set theoretic constructions, our way of handling identity  seems to be
specific to type theory.
According to the propositions-as-types conception, the \emph{proposition} that two elements of the same type $a,b:A$ are equal must correspond to some \emph{type}.
Since this proposition depends on what $a$ and $b$ are, these \define{equality types} or \define{identity types} must be type families dependent on two copies of $A$.

We may write the family as $\idtypevar{A}:A\to A\to\type$, so that $\idtype[A]ab$ is the type representing the proposition of equality between $a$ and $b$.
Once we are familiar with propositions-as-types, however, it is convenient to also use the standard equality symbol for this; thus ``$\id{a}{b}$'' will also be a notation for the \emph{type} $\idtype[A]ab$ corresponding to the proposition that $a$ equals $b$.
For clarity, we may also write ``$\id[A]{a}{b}$'' to specify the type $A$.
If we have an element of $\id[A]{a}{b}$, we may say that $a$ and $b$ are \define{equal}, or sometimes \define{propositionally equal} if we want to emphasize that this is different from the judgmental equality $a\jdeq b$ discussed in \autoref{sec:types-vs-sets}.
\indexdef{equality!propositional}%
\indexdef{propositional!equality}%

Just as we remarked in \autoref{sec:pat} that the propositions-as-types versions of ``or'' and ``there exists'' can include more information than just the fact that the proposition is true, nothing prevents  the type $\id{a}{b}$ from also including more information.
Indeed, this is the cornerstone of the homotopical interpretation, where we regard witnesses of $\id{a}{b}$ as \emph{paths}\indexdef{path} or \emph{equivalences} between $a$ and $b$ in the space $A$.  Just as there can be more than one path between two points of a space, there can be more than one witness that two objects are equal.  Put differently, we may regard $\id{a}{b}$ as the type of \emph{identifications}\indexdef{identification} of $a$ and $b$, and there may be many different ways in which $a$ and $b$ can be identified.
We will return to the interpretation in \autoref{cha:basics}; for now we focus on the basic rules for the identity type.
Just like all the other types considered in this chapter, it will have rules for formation, introduction, elimination, and computation, which behave formally in exactly the same way.

The formation rule says that given a type $A:\UU$ and two elements $a,b:A$, we can form the type $(\id[A]{a}{b}):\UU$ in the same universe.
The basic way to construct an element of $\id{a}{b}$ is to know that $a$ and $b$ are the same.
Thus, the introduction rule is a dependent function
\[\refl{} : \prd{a:A} (\id[A]{a}{a})\]
called \define{reflexivity},
\indexdef{reflexivity!of equality}%
which says that every element of $A$ is equal to itself (in a specified way).  We regard $\refl{a}$ as being the
constant path\indexdef{path!constant}\indexsee{loop!constant}{path, constant}
at the point $a$.

In particular, this means that if $a$ and $b$ are \emph{judgmentally} equal, $a\jdeq b$, then we also have an element $\refl{a} : \id[A]{a}{b}$.
This is well-typed because $a\jdeq b$ means that also the type $\id[A]{a}{b}$ is judgmentally equal to $\id[A]{a}{a}$, which is the type of $\refl{a}$.

The induction principle (i.e.\ the elimination rule) for the identity types is one of the most subtle parts of type theory, and crucial to the homotopy interpretation.
We begin by considering an important consequence of it, the principle that ``equals may be substituted for equals,'' as expressed by the following:
\index{indiscernability of identicals}%
\index{equals may be substituted for equals}%
\begin{description}
\item[Indiscernability of identicals:]
For every family 
\[
C : A \to \UU
\]
there is a function
\[
f : \prd{x,y:A}{p:\id[A] x y} C(x) \to C(y)
\]
such that
\[
f(x,x,\refl{x}) \defeq \idfunc[C(x)].
\]
\end{description}
This says that every family of types $C$ respects equality, in the sense that applying $C$ to \emph{equal} elements of $A$ also results in a function between the resulting types. The displayed equality states that the function associated to reflexivity is the identity function (and we shall see that, in general, the function $f(x,y,p): C(x) \to C(y)$ is always an equivalence of types).

Indiscernability of identicals can be regarded as a recursion principle for the identity type, analogous to those given for booleans and natural numbers above.
Just as $\rec{\nat}$ gives a specified map $\nat\to C$ for any other type $C$ of a certain sort, indiscernability of identicals gives a specified map from $\id[A] x y$ to certain other reflexive, binary relations on $A$, namely those of the form $C(x) \to C(y)$ for some unary predicate $C(x)$.
We could also formulate a more general recursion principle with respect to reflexive relations of the more general form $C(x,y)$.
However, in order to fully characterize the identity type, we must generalize this recursion principle to an induction principle, which not only considers maps out of $\id[A] x y$ but also families over it.
Put differently, we consider not only allowing equals to be substituted for equals, but also taking into account the evidence $p$ for the equality.
    
\subsection{Path induction}

\index{generation!of a type, inductive|(}
The induction principle for the identity type is called \define{path induction},
\index{path!induction|(}%
\index{induction principle!for identity type|(}%
in view of the homotopical interpretation to be explained in  the introduction to \autoref{cha:basics}.  It can be seen as stating that the family of identity types is freely generated by the elements of the form $\refl{x}: \id{x}{x}$.

\begin{description}
\item[Path induction:] 
  Given a family 
  \[ C : \prd{x,y:A} (\id[A]{x}{y}) \to \UU \]
  and a function
  \[ c :  \prd{x:A} C(x,x,\refl{x}),\]
  there is a function
  \[ f : \prd{x,y:A}{p:\id[A]{x}{y}} C(x,y,p) \]
  such that 
  \[ f(x,x,\refl{x}) \defeq c(x). \]
\end{description}

Note that just like the induction principles for products, coproducts, natural numbers, and so on, path induction allows us to define \emph{specified} functions which exhibit appropriate computational behavior.
Just as we have \emph{the} function $f:\nat\to C$ defined by recursion from $c_0:C$ and $c_s:\nat \to C \to C$, which moreover satisfies $f(0)\jdeq c_0$ and $f(\suc(n))\jdeq c_s(n,f(n))$, we have \emph{the} function $f : \dprd{x,y:A}{p:\id[A]{x}{y}} C(x,y,p)$ defined by path induction from $c :  \prd{x:A} C(x,x,\refl{x})$, which moreover satisfies $f(x,x,\refl{x}) \jdeq c(x)$.

To understand the meaning of this principle, consider first the simpler case when $C$
does not depend on $p$.  Then we have $C:A\to A\to \UU$, which we may
regard as a predicate depending on two elements of $A$.  We are
interested in knowing when the proposition $C(x,y)$ holds for some pair
of elements $x,y:A$.  In this case, the hypothesis of path induction
says that we know $C(x,x)$ holds for all $x:A$, i.e.\ that if we
evaluate $C$ at the pair $x, x$, we get a true proposition --- so $C$ is
a reflexive relation.  The conclusion then tells us that $C(x,y)$ holds
whenever $\id{x}{y}$.  This is exactly the more general recursion principle
for reflexive relations mentioned above.

The general, inductive form of the rule allows $C$ to also depend on the witness $p:\id{x}{y}$ to the identity between $x$ and $y$.  In the premise, we not only replace $x, y$ by $x,x$, but also simultaneously replace $p$ by reflexivity: to prove a property for all elements $x,y$ and paths $p : \id{x}{y}$ between them, it suffices to consider all the cases where the elements are $x,x$ and the path is $\refl{x}: \id{x}{x}$.  If we were viewing types just as sets, it would be unclear what this buys us, but since there may be many different identifications $p : \id{x}{y}$ between $x$ and $y$, it makes sense to keep track of them in considering families over the type $\id[A]{x}{y}$.
In \autoref{cha:basics} we will see that this is very important to the homotopy interpretation.

If we package up path induction into a single function, it takes the form:
\symlabel{defn:induction-ML-id}%
\begin{narrowmultline*}
  \indid{A} :  \dprd{C : \prd{x,y:A} (\id[A]{x}{y}) \to \UU}
  \Parens{\tprd{x:A} C(x,x,\refl{x})} \to 
  \narrowbreak
  \dprd{x,y:A}{p:\id[A]{x}{y}}   C(x,y,p)
\end{narrowmultline*}
with the equality\index{computation rule!for identity types}
\[ \indid{A}(C,c,x,x,\refl{x}) \defeq c(x). \]
The function $ \indid{A}$ is traditionally called $J$.
\indexsee{J@$J$}{induction principle for identity type}%
We leave it as an easy exercise to show that indiscernability of identicals follows from path induction.  

\mentalpause

Given a proof $p : \id{a}{b}$,
path induction requires us to replace \emph{both} $a$ and $b$ with the same unknown element $x$; thus in order to define an element of a family
$C$, for all pairs of elements of $A$, it suffices to define it on the diagonal.
In some proofs, however, it is simpler to make use of an equation $p : \id{a}{b}$ by replacing all occurrences of $b$ with $a$ (or vice versa), because it is sometimes easier to do the remainder of the proof for the specific element $a$ mentioned in the equality than for a general unknown $x$.  This motivates a second induction principle for identity types, which says that the family of types $\id[A]{a}{x}$ is generated by the element $\refl{a} : \id{a}{a}$.  As we show below, this second principle is equivalent to the first; it is just sometimes a more convenient formulation.

\index{path!induction based}%
\index{induction principle!for identity type!based}%
\begin{description}
\item[Based path induction:] 
  Fix an element $a:A$, and suppose given a family
  \[ C : \prd{x:A} (\id[A]{a}{x}) \to \UU \]
  and an element
  \[ c : C(a,\refl{a}). \]
  Then we obtain a function
  \[ f : \prd{x:A}{p:\id{a}{x}} C(x,p) \]
  such that
  \[ f(a,\refl{a}) \defeq c.\]
\end{description}

Here, $C(x,p)$ is a family of types, where $x$ is an element of $A$ and $p$ is an element of the identity type $\id[A]{a}{x}$, for fixed $a$ in $A$. The based path induction principle says that to define an element of this family for all $x$ and $p$, it suffices to consider
just the case where $x$ is $a$ and $p$ is $\refl{a} : \id{a}{a}$.

Packaged as a function, based path induction becomes:
\symlabel{defn:induction-PM-id}%
\begin{align*}
  \indidb{A} :  \dprd{a:A}{C : \prd{x:A} (\id[A]{a}{x}) \to \UU}
  C(a,\refl{a}) \to \dprd{x:A}{p : \id[A]{a}{x}} C(x,p) 
\end{align*}
with the equality
\[ \indidb{A}(a,C,c,a,\refl{a}) \defeq c. \]

Below, we show that path induction and based path induction are equivalent.  Because of this, we will sometimes be sloppy and also refer to based path induction simply as ``path induction,'' relying on the reader to infer which principle is meant from the form of the proof.  

\begin{rmk}
Intuitively, the induction principle for the natural numbers expresses the fact that the only natural numbers are $0$ and $\suc(n)$, so if we prove a property for these cases, then we have proved it for all natural numbers.  Applying this same reading to path induction, we might loosely say that path induction expresses the fact that the only path is \refl{}, so if we prove a property for reflexivity, then we have proved it for all paths.  However, this reading is quite confusing in the context of the homotopy interpretation of paths, where there may be many different ways in which two elements $a$ and $b$ can be identified, and therefore many different elements of the identity type!  How can there be many different paths, but at the same time we have an induction principle asserting that the only path is reflexivity?

The key observation is that it is not the identity \emph{type} that is inductively defined, but the identity \emph{family}.
In particular, path induction says that the \emph{family} of types $(\id[A]{x}{y})$, as $x,y$ vary over all elements of $A$, is inductively defined by the elements of the form $\refl{x}$.
This means that to give an element of any other family $C(x,y,p)$ dependent on a \emph{generic} element $(x,y,p)$ of the identity family, it suffices to consider the cases of the form $(x,x,\refl{x})$.
In the homotopy interpretation, this says that the type of triples $(x,y,p)$, where $x$ and $y$ are the endpoints of the path $p$ (in other words, the $\Sigma$-type $\sm{x,y:A}(\id{x}{y})$), is inductively generated by the constant loops\index{path!constant} at each point $x$.
In homotopy theory, the space corresponding to $\sm{x,y:A}(\id{x}{y})$ is the \emph{free path space} --- the space of paths in $A$ whose endpoints may vary --- and it is in fact the case that any point of this space is homotopic to the constant loop at some point, since we can simply retract one of its endpoints along the given path.

Similarly, based path induction says that for a fixed $a:A$, the \emph{family} of types $(\id[A]{a}{y})$, as $y$ varies over all elements of $A$, is inductively defined by the element $\refl{a}$.
Thus, to give an element of any other family $C(y,p)$ dependent on a generic element $(y,p)$ of this family, it suffices to consider the case $(a,\refl{a})$.
Homotopically, this expresses the fact that the space of paths starting at some chosen point (the \emph{based path space} at that point, which type-theoretically is $\sm{y:A} (\id{a}{y})$) is contractible to the constant loop on the chosen point.
Note that according to propositions-as-types, the type $\sm{y:A}(\id{a}{y})$ can be regarded as ``the type of all elements of $A$ which are equal to $a$'', a type-theoretic version of the ``singleton\index{type!singleton} subset'' $\{a\}$.

Neither of these two principles provides a way to give an element of a family $C(p)$ where $p$ has \emph{two fixed endpoints} $a$ and $b$.  In particular, for a family $C(p : \id[A]{a}{a})$ dependent on a loop, we \emph{cannot} apply path induction and consider only the case for $C(\refl{a})$, and consequently, we cannot prove
that all loops are reflexivity.  Thus, inductively defining the identity family does not prohibit non-reflexivity paths in specific instances of the identity type.
In other words, a path $p:\id{x}{x}$ may  be not equal to reflexivity as an element of $(\id{x}{x})$, but the pair $(x,p)$ will nevertheless be equal to the pair $(x,\refl{x})$ as elements of $\sm{y:A}(\id{x}{y})$.
\end{rmk}

\index{path!induction|)}%
\index{induction principle!for identity type|)}%
\index{generation!of a type, inductive|)}

\subsection{Equivalence of path induction and based path induction}

The two induction principles for the identity type introduced above are equivalent.
It is easy to see that path induction follows from based path induction principle.
Indeed, let us assume the premises of path induction:
\begin{align*}
C &: \prd{x,y:A}(\id[A]{x}{y}) \to \UU,\\
c &: \prd{x:A} C(x,x,\refl{x}).
\end{align*}
Now, given an element $x:A$, we can instantiate both of the above, obtaining
\begin{align*}
C' &: \prd{y:A} (\id[A]{x}{y}) \to \UU,  \\
C' &\defeq C(x), \\
c' &: C'(x,\refl{x}), \\
c' &\defeq c(x).
\end{align*}
Clearly, $C'$ and $c'$ match the premises of based path induction and hence we can construct 
\begin{equation*}
  g : \prd{y:A}{p : \id{x}{y}} C'(y,p)
\end{equation*}
with the defining equality
\[ g(x,\refl{x}) \defeq c'.\]
Now we observe that $g$'s codomain is equal to $C(x,y,p)$.
Thus, discharging our assumption $x:A$, we can derive a function 
\[ f : \prd{x,y:A}{p : \id[A]{x}{y}} C(x,y,p) \]
with the required judgmental equality $f(x,x,\refl{x}) \judgeq g(x,\refl{x}) \defeq c' \defeq c(x)$.

Another proof of this fact is to observe that any such $f$ can be obtained as an instance of $\indidb{A}$
so it suffices to define $\indid{A}$ in terms of $\indidb{A}$ as
\[ \indid{A}(C,c,x,y,p) \defeq \indidb{A}(x,C(x),c(x),y,p). \]

The other direction is a bit trickier; it is not clear how we can use a particular instance of path induction to derive a particular instance of
based path induction. What we can do instead is to construct one instance of path induction which shows 
all possible instantiations of based path induction at once.
Define
\begin{align*}
D &: \prd{x,y:A} (\id[A]{x}{y}) \to \UU, \\
D(x,y,p) &\defeq \prd{C : \prd{z:A} (\id[A]{x}{z}) \to \UU} C(x,\refl{x}) \to C(y,p).
\end{align*}
Then we can construct the function
\begin{align*}
d &: \prd{x : A} D(x,x,\refl{x}), \\
d &\defeq \lamu{x:A}\lamu{C:\prd{z:A}{p : \id[A]{x}{z}} \UU}\lam{c:C(x,\refl{x})} c
\end{align*}
and hence using path induction obtain
\[ f : \prd{x,y:A}{p:\id[A]{x}{y}} D(x,y,p) \]
with $f(x,x,\refl{x}) \defeq d(x)$. Unfolding the definition of $D$, we can expand the type of $f$:
\[ f : \prd{x,y:A}{p:\id[A]{x}{y}}{C : \prd{z:A} (\id[A]{x}{z}) \to \UU} C(x,\refl{x}) \to C(y,p). \]
Now given $x:A$ and $p:\id[A]{a}{x}$, we can derive the conclusion of based path induction:
\[ f(a,x,p,C,c) : C(x,p). \]
Notice that we also obtain the correct definitional equality.

Another proof is to observe that any use of based path induction is an instance of $\indidb{A}$  and to define
\begin{narrowmultline*}
\indidb{A}(a,C,c,x,p) \defeq \narrowbreak
\indid{A}
  \begin{aligned}[t]
    \big(
    &\big(\lamu{x,y:A}{p:\id[A]{x}{y}} \tprd{C : \prd{z:A} (\id[A]{x}{z}) \to \UU} C(x,\refl{x}) \to C(y,p) \big),\\
    &(\lamu{x:A}{C:\prd{z:A} (\id[A]{x}{z}) \to \UU}{d:C(x,\refl{x})} d),
     a, x, p, C, c \big) 
   \end{aligned}
\end{narrowmultline*}

Note that the construction given above uses universes. That is, if we want to
model $\indidb{A}$ with $C : \prd{x:A} (\id[A]{a}{x}) \to \UU_i$, we need
to use $\indid{A}$ with 
\[ D:\prd{x,y:A} (\id[A]{x}{y}) \to \UU_{i+1} \]
since $D$ quantifies over all $C$ of the given type. While this is
compatible with our definition of universes, it is also possible to
derive $\indidb{A}$ without using universes: we can show that $\indid{A}$ entails \autoref{lem:transport,thm:contr-paths}, and that these two principles imply $\indidb{A}$ directly.
We leave the details to the reader as \autoref{ex:pm-to-ml}.

We can use either of the foregoing formulations of identity types
to establish that equality is an equivalence relation, that every function preserves equality and that every family respects equality. We leave the details to the next chapter, where this will be derived  and explained in the context of homotopy type theory.

\subsection{Disequality}
\label{sec:disequality}

Finally, let us also say something about \define{disequality},
\indexdef{disequality}%
which is negation of equality:%
\footnote{We use ``inequality''
  to refer to $<$ and $\leq$. Also, note that this is negation of the \emph{propositional} identity type.
Of course, it makes no sense to negate judgmental equality $\jdeq$, because judgments are not subject to logical operations.}
\begin{equation*}
  (x \neq_A y) \ \defeq\ \lnot (\id[A]{x}{y}).
\end{equation*}
If $x\neq y$, we say that $x$ and $y$ are \define{unequal}
\indexdef{unequal}%
or \define{not equal}.
Just like negation, disequality plays a less important role here than it does in classical\index{mathematics!classical}
mathematics. For example, we cannot prove that two things are equal by proving that they
are not unequal: that would be an application of the classical law of double negation, see \autoref{sec:intuitionism}.

Sometimes it is useful to phrase disequality in a positive way. For example,
in~\autoref{RD-inverse-apart-0} we shall prove that a real number $x$ has an inverse if,
and only if, its distance from~$0$ is positive, which is a stronger requirement than $x
\neq 0$.

\index{type!identity|)}%

\sectionNotes

The type theory presented here is a version of Martin-L\"{o}f's intuitionistic type 
theory~\cite{Martin-Lof-1972,Martin-Lof-1973,Martin-Lof-1979,martin-lof:bibliopolis}, which itself is based on and influenced 
by the foundational work of Brouwer \cite{beeson}, Heyting~\cite{heyting1966intuitionism}, Scott~\cite{scott70}, de 
Bruijn~\cite{deBruijn-1973}, Howard~\cite{howard:pat}, Tait~\cite{Tait-1966,Tait-1968}, and Lawvere~\cite{lawvere:adjinfound}\index{Lawvere}.
\index{proof!assistant}%
Three principal variants of Martin-L\"{o}f's type theory underlie the \NuPRL \cite{constable+86nuprl-book}, \Coq~\cite{Coq}, and 
\Agda \cite{norell2007towards} computer implementations of type theory.  The theory given here differs from these formulations in a number 
of respects, some of which are critical to the homotopy interpretation, while others are technical conveniences or involve concepts that 
have not yet been studied in the homotopical setting.

\index{type theory!intensional}%
\index{type theory!extensional}%
\index{intensional type theory}%
\index{extensional!type theory}%
Most significantly, the type theory described here is derived from the \emph{intensional} version of Martin-L\"{o}f's type 
theory~\cite{Martin-Lof-1973}, rather than the \emph{extensional} version~\cite{Martin-Lof-1979}.  Whereas the extensional theory makes no 
distinction between judgmental and propositional equality, the intensional theory regards judgmental equality as purely definitional, and 
admits a much broader, proof-relevant interpretation of the identity type that is central to the homotopy interpretation.  From the 
homotopical perspective, extensional type theory confines itself to homotopically discrete sets (see \autoref{sec:basics-sets}), whereas the 
intensional theory admits types with higher-dimensional structure.  The \NuPRL system~\cite{constable+86nuprl-book} is extensional, whereas 
both \Coq~\cite{Coq} and \Agda~\cite{norell2007towards} are intensional.  Among intensional type theories, there are a number of variants 
that differ in the structure of identity proofs.  The most liberal interpretation, on which we rely here, admits a \emph{proof-relevant} 
interpretation of equality, whereas more restricted variants impose restrictions such as \emph{uniqueness of identity proofs 
  (UIP)}~\cite{Streicher93},
\indexsee{UIP}{uniqueness of identity proofs}%
\index{uniqueness!of identity proofs}%
stating that any two proofs of equality are judgmentally equal, and \emph{Axiom K}~\cite{Streicher93},
\index{axiom!Streicher's Axiom K}
stating that 
the only proof of equality is reflexivity (up to judgmental equality).  These additional requirements may be selectively imposed in the \Coq 
and \Agda\ systems.


Another point of variation among intensional theories is the strength of judgmental equality, particularly as regards objects of function type.  Here we include the uniqueness principle\index{uniqueness!principle} ($\eta$-conversion) $f \jdeq \lam{x} f(x)$, as a principle of judgmental equality.  This principle is used, for example, in \autoref{sec:univalence-implies-funext}, to show that univalence implies propositional function extensionality.  Uniqueness principles are sometimes considered for other types.
For instance, the uniqueness principle\index{uniqueness!principle!for product types} for cartesian products would be a judgmental version of the propositional equality $\uppt$ which we constructed in \autoref{sec:finite-product-types}, saying that $u \jdeq (\proj1(u),\proj2(u))$.
This and the corresponding version for dependent pairs would be reasonable choices (which we did not make), but we cannot include all such rules, because the corresponding uniqueness principle for identity types would trivialize all the higher homotopical structure.  So we are \emph{forced} to leave it out, and the question then becomes where to draw the line. With regards to inductive types, we discuss these points further in~\autoref{sec:htpy-inductive}.

It is important for our purposes that (propositional) equality of functions is taken to be \emph{extensional} (in a different sense than that used above!).
This is not a consequence of the rules in this chapter; it will be expressed by \autoref{axiom:funext}.
\index{function extensionality}%
This decision is significant for our purposes, because it specifies that equality of functions is as expected in mathematics.  Although we include \autoref{axiom:funext} as an axiom, it may be derived from the univalence axiom and the uniqueness principle for functions\index{uniqueness!principle!for function types} (see \autoref{sec:univalence-implies-funext}), as well as from the existence of an interval type (see \autoref{thm:interval-funext}).

Regarding inductive types such as products, $\Sigma$-types, coproducts, natural numbers, and so on (see \autoref{cha:induction}), there are additional choices regarding precisely how to  formulate induction and recursion.
\index{pattern matching}%
Formally, one may describe type theory by taking either \emph{pattern matching} or \emph{induction principles} as basic and deriving the other; see \autoref{cha:rules}.
However, pattern matching in general is not yet well understood from the homotopical perspective (in particular, ``nested'' or ``deep'' pattern matching is difficult to make general sense of for higher inductive types).
Moreover, it can be dangerous unless sufficient care is taken: for instance, the form of pattern matching implemented by default in \Agda
\index{proof!assistant!Agda@\textsc{Agda}}%
allows proving Axiom K.
\index{axiom!Streicher's Axiom K}%
For these reasons, we have chosen to regard the induction principle as the basic property of an inductive definition, with pattern matching justified in terms of induction.

\index{proof!assistant!Coq@\textsc{Coq}}%
Unlike the type theory of \Coq, we do not include a primitive type of propositions.  Instead, as discussed in \autoref{sec:pat}, we embrace 
the \emph{propositions-as-types (PAT)} principle, identifying propositions with types.
This was suggested originally by de Bruijn~\cite{deBruijn-1973}, Howard~\cite{howard:pat}, Tait~\cite{Tait-1968}, and Martin-L\"{o}f~\cite{Martin-Lof-1972}.
(Our decision is explained more fully in \autoref{subsec:pat?,subsec:hprops}.)

We do, however, include a full cumulative hierarchy of universes, so that the type formation and equality judgments become instances of the membership and equality judgments for a universe.
As a convenience, we regard objects of a universe as types, rather than as codes for types; in the terminology of \cite{martin-lof:bibliopolis}, this means we use ``Russell-style universes'' rather than ``Tarski-style universes''.
\index{type!universe!Tarski-style}%
\index{type!universe!Russell-style}%
An alternative would be to use Tarski-style universes, with an explicit coercion\index{coercion, universe-raising} function required to make an element $A:\UU$ of a universe into a type $\mathsf{El}(A)$, and just say that the coercion is omitted when working informally.

We also treat the universe hierarchy as cumulative, in that every type in $\UU_i$ is also in $\UU_j$ for each $j\geq i$.
There are different ways to implement cumulativity formally: the simplest is just to include a rule that if $A:\UU_i$ then $A:\UU_j$.
However, this has the annoying consequence that for a type family $B:A\to \UU_i$ we cannot conclude $B:A\to\UU_j$, although we can conclude $\lam{a} B(a) : A\to\UU_j$.
A more sophisticated approach that solves this problem is to introduce a judgmental subtyping relation $<:$ generated by $\UU_i<:\UU_j$, but this makes the type theory more complicated to study.
Another alternative would be to include an explicit coercion function $\uparrow : \UU_i \to \UU_j$, which could be omitted when working informally.

It is also not necessary that the universes be indexed by natural numbers and linearly ordered.
For some purposes, it is more appropriate to assume only that every universe is an element of some larger universe, together with a ``directedness'' property that any two universes are jointly contained in some larger one.
There are many other possible variations, such as including a universe ``$\UU_{\omega}$'' that contains all $\UU_i$ (or even higher ``large cardinal'' type universes), or by internalizing the hierarchy into a type family $\lam{i} \UU_i$.
The latter is in fact done in \Agda.

The path induction principle for identity types was formulated by Martin-L\"{o}f~\cite{Martin-Lof-1972}.
The based path induction rule in the setting of Martin-L\"of type theory is due to Paulin-Mohring \cite{Moh93}; it can be seen as an intensional generalization of the concept of ``pointwise functionality''\index{pointwise!functionality} for hypothetical judgments from \NuPRL~\cite[Section~8.1]{constable+86nuprl-book}.
The fact that Martin-L\"of's rule implies Paulin-Mohring's was proved by Streicher using Axiom K (see~\autoref{sec:hedberg}), by Altenkirch and Goguen as in \autoref{sec:identity-types}, and finally by Hofmann without universes (as in \autoref{ex:pm-to-ml}); see~\cite[\S1.3 and Addendum]{Streicher93}.

\sectionExercises

\begin{ex}\label{ex:composition}
  Given functions $f:A\to B$ and $g:B\to C$, define their \define{composite}
  \indexdef{composition!of functions}%
  \indexdef{function!composition}%
 $g\circ f:A\to C$.
  \index{associativity!of function composition}%
  Show that we have $h \circ (g\circ f) \jdeq (h\circ g)\circ f$.
\end{ex}

\begin{ex}
  Derive the recursion principle for products $\rec{A\times B} $ using only the projections, and verify that the definitional equalities are valid.
  Do the same for $\Sigma$-types.
\end{ex}

\begin{ex}
  Derive the induction principle for products $\ind{A\times B}$, using only the projections and the propositional uniqueness principle $\uppt$.
  Verify that the definitional equalities are valid.
  Generalize $\uppt$ to $\Sigma$-types, and do the same for $\Sigma$-types.
  \emph{(This requires concepts from \autoref{cha:basics}.)}
\end{ex}

\begin{ex}\label{ex:iterator}
\index{iterator!for natural numbers}
Assuming as given only the \emph{iterator} for natural numbers
\[\ite : \prd{C:\UU} C \to (C \to C) \to \nat \to C \]
with the defining equations
\begin{align*}
\ite(C,c_0,c_s,0)  &\defeq c_0, \\
\ite(C,c_0,c_s,\suc(n)) &\defeq c_s(\ite(C,c_0,c_s,n))
\end{align*}
derive the recursor $\rec{\nat}$.
\end{ex}

\begin{ex}\label{ex:sum-via-bool}
\index{type!coproduct}%
Show that if we define $A + B \defeq \sm{x:\bool} \rec{\bool}(\UU,A,B,x)$, then we can give a definition of $\ind{A+B}$ for which the definitional equalities stated in \autoref{sec:coproduct-types} hold.
\end{ex}

\begin{ex}\label{ex:prod-via-bool}
\index{type!product}%
Show that if we define $A \times B \defeq \prd{x:\bool}\rec{\bool}(\UU,A,B,x)$, then we can give a definition of  $\ind{A\times B}$ for which the definitional equalities stated in \autoref{sec:finite-product-types} hold propositionally (i.e.\ using equality types).
\emph{(This requires the function extensionality axiom, which is introduced in \autoref{sec:compute-pi}.)}
\end{ex}

\begin{ex}\label{ex:pm-to-ml}
Give an alternative derivation of $\indidb{A}$ from $\indid{A}$ which avoids the use of universes.
  \emph{(This is easiest using concepts from later chapters.)}
\end{ex}

\begin{ex}
  \index{multiplication!of natural numbers}%
  Define multiplication and exponentiation using $\rec{\nat}$.
  Verify that $(\nat,+,0,\times,1)$ is a semiring\index{semiring} using only $\ind{\nat}$.  
\end{ex}

\begin{ex}\label{ex:fin}
  \index{finite!sets, family of}%
  Define the type family $\Fin : \nat \to \UU$ mentioned at the end of \autoref{sec:universes}, and the dependent function $\fmax : \prd{n:\nat} \Fin(n+1)$ mentioned in \autoref{sec:pi-types}.
\end{ex}

\begin{ex}\label{ex:ackermann}
  \indexdef{function!Ackermann}%
  \indexdef{Ackermann function}%
  Show that the Ackermann function $\ack : \nat \to \nat \to \nat$ is definable using only $\rec{\nat}$ satisfying the following equations:
  \begin{align*}
    \ack(0,n) &\jdeq \suc(n), \\
    \ack(\suc(m),0) &\jdeq \ack(m,1), \\
    \ack(\suc(m),\suc(n)) &\jdeq \ack(m,\ack(\suc(m),n)).
  \end{align*}
\end{ex}

\begin{ex}\label{ex:neg-ldn}
  Show that for any type $A$, we have $\neg\neg\neg A \to \neg A$.
\end{ex}

\begin{ex}\label{ex:tautologies}
  Using the propositions as types interpretation, derive the following tautologies.
  \begin{enumerate}
  \item If $A$, then (if $B$ then $A$).
  \item If $A$, then not (not $A$).
  \item If (not $A$ or not $B$), then not ($A$ and $B$).
  \end{enumerate}
\end{ex}

\begin{ex}\label{ex:not-not-lem}
  Using propositions-as-types, derive the double negation of the principle of excluded middle, i.e.\ prove \emph{not (not ($P$ or not $P$))}.
\end{ex}

\begin{ex}\label{ex:without-K}
  Why do the induction principles for identity types not allow us to construct a function $f: \prd{x:A}{p:\id{x}{x}} (\id{p}{\refl{x}})$ with the defining equation
  \[ f(x,\refl{x}) \defeq \refl{\refl{x}} \quad ?\]
\end{ex}

\begin{ex}\label{ex:subtFromPathInd}
  Show that indiscernability of identicals follows from path induction.  
\end{ex}


\chapter{Homotopy type theory}
\label{cha:basics}

The central new idea in homotopy type theory is that types can be regarded as
spaces in homotopy theory, or higher-dimensional groupoids in category
theory.  

\index{classical!homotopy theory|(}
\index{higher category theory|(}
We begin with a brief summary of the connection between homotopy theory
and higher-dimensional category theory.  
In classical homotopy theory, a space $X$ is a set of points equipped
with a topology,
\indexsee{space!topological}{topological space}
\index{topological!space}
and a path between points $x$ and $y$ is represented by
a continuous map $p : [0,1] \to X$, where $p(0) = x$ and $p(1) = y$.
\index{path!topological}
\index{topological!path}
This function can be thought of as giving a point in $X$ at each
``moment in time''.  For many purposes, strict equality of paths
(meaning, pointwise equal functions) is too fine a notion.  For example,
one can define operations of path concatenation (if $p$ is a path from
$x$ to $y$ and $q$ is a path from $y$ to $z$, then the concatenation $p
\ct q$ is a path from $x$ to $z$) and inverses ($\opp p$ is a path
from $y$ to $x$).  However, there are natural equations between these
operations that do not hold for strict equality: for example, the path
$p \ct \opp p$ (which walks from $x$ to $y$, and then back along the
same route, as time goes from $0$ to $1$) is not strictly equal to the
identity path (which stays still at $x$ at all times).

The remedy is to consider a coarser notion of equality of paths called
\emph{homotopy}.
\index{homotopy!topological}
A homotopy between a pair of continuous maps $f :
X_1 \to X_2$ and $g : X_1\to X_2$ is a continuous map $H : X_1
\times [0, 1] \to X_2$ satisfying $H(x, 0) = f (x)$ and $H(x, 1) =
g(x)$.  In the specific case of paths $p$ and $q$ from $x$ to $y$, a homotopy is a
continuous map $H : [0,1] \times [0,1] \rightarrow X$
such that $H(s,0) = p(s)$ and $H(s,1) = q(s)$ for all $s\in [0,1]$.
In this case we require also that $H(0,t) = x$ and $H(1,t)=y$ for all $t\in [0,1]$,
so that for each $t$ the function $H(\blank,t)$ is again a path from $x$ to $y$;
a homotopy of this sort is said to be \emph{endpoint-preserving} or \emph{rel endpoints}.
Such a homotopy is the
image in $X$ of a square that fills in the space between $p$ and $q$,
which can be thought of as a ``continuous deformation'' between $p$ and
$q$, or a 2-dimensional \emph{path between paths}.\index{path!2-}

For example, because
$p \ct \opp p$ walks out and back along the same route, you know that
you can continuously shrink $p \ct \opp p$ down to the identity
path---it won't, for example, get snagged around a hole in the space.
Homotopy is an equivalence relation, and operations such as
concatenation, inverses, etc., respect it.  Moreover, the homotopy
equivalence classes of loops\index{loop} at some point $x_0$ (where two loops $p$
and $q$ are equated when there is a \emph{based} homotopy between them,
which is a homotopy $H$ as above that additionally satisfies $H(0,t) =
H(1,t) = x_0$ for all $t$) form a group called the \emph{fundamental
  group}.\index{fundamental!group}  This group is an \emph{algebraic invariant} of a space, which
can be used to investigate whether two spaces are \emph{homotopy
  equivalent} (there are continuous maps back and forth whose composites
are homotopic to the identity), because equivalent spaces have
isomorphic fundamental groups.

Because homotopies are themselves a kind of 2-dimensional path, there is
a natural notion of 3-dimensional \emph{homotopy between homotopies},\index{path!3-}
and then \emph{homotopy between homotopies between homotopies}, and so
on.  This infinite tower of points, path, homotopies, homotopies between
homotopies, \ldots, equipped with algebraic operations such as the
fundamental group, is an instance of an algebraic structure called a
(weak) \emph{$\infty$-groupoid}.  An $\infty$-groupoid\index{.infinity-groupoid@$\infty$-groupoid} consists of a
collection of objects, and then a collection of \emph{morphisms}\indexdef{morphism!in an .infinity-groupoid@in an $\infty$-groupoid} between
objects, and then \emph{morphisms between morphisms}, and so on,
equipped with some complex algebraic structure; a morphism at level $k$ is called a \define{$k$-morphism}\indexdef{k-morphism@$k$-morphism}.  Morphisms at each level
have identity, composition, and inverse operations, which are weak in
the sense that they satisfy the groupoid laws (associativity of
composition, identity is a unit for composition, inverses cancel) only
up to morphisms at the next level, and this weakness gives rise to
further structure. For example, because associativity of composition of
morphisms $p \ct (q \ct r) = (p \ct q) \ct r$ is itself a
higher-dimensional morphism, one needs an additional operation relating
various proofs of associativity: the various ways to reassociate $p \ct
(q \ct (r \ct s))$ into $((p \ct q) \ct r) \ct s$ give rise to Mac
Lane's pentagon\index{pentagon, Mac Lane}.  Weakness also creates non-trivial interactions between
levels.

Every topological space $X$ has a \emph{fundamental $\infty$-groupoid}
\index{.infinity-groupoid@$\infty$-groupoid!fundamental}
\index{fundamental!.infinity-groupoid@$\infty$-groupoid}
whose
$k$-mor\-ph\-isms are the $k$-dimen\-sional paths in $X$.  The weakness of the
$\infty$-group\-oid corresponds directly to the fact that paths form a
group only up to homotopy, with the $(k+1)$-paths serving as the
homotopies between the $k$-paths.  Moreover, the view of a space as an
$\infty$-groupoid preserves enough aspects of the space to do homotopy theory:
the fundamental $\infty$-groupoid construction is adjoint\index{adjoint!functor} to the
geometric\index{geometric realization} realization of an $\infty$-groupoid as a space, and this
adjunction preserves homotopy theory (this is called the \emph{homotopy
  hypothesis/theorem},
\index{hypothesis!homotopy}
\index{homotopy!hypothesis}
because whether it is a hypothesis or theorem
depends on how you define $\infty$-groupoid).  For example, you can
easily define the fundamental group of an $\infty$-groupoid, and if you
calculate the fundamental group of the fundamental $\infty$-groupoid of
a space, it will agree with the classical definition of fundamental
group of that space.  Because of this correspondence, homotopy theory
and higher-dimensional category theory are intimately related.

\index{classical!homotopy theory|)}%
\index{higher category theory|)}%

\mentalpause

Now, in homotopy type theory each type can be seen to have the structure
of an $\infty$-groupoid.  Recall that for any type $A$, and any $x,y:A$,
we have a identity type $\id[A]{x}{y}$, also written $\idtype[A]{x}{y}$
or just $x=y$.  Logically, we may think of elements of $x=y$ as evidence
that $x$ and $y$ are equal, or as identifications of $x$ with
$y$. Furthermore, type theory (unlike, say, first-order logic) allows us
to consider such elements of $\id[A]{x}{y}$ also as individuals which
may be the subjects of further propositions.  Therefore, we can
\emph{iterate} the identity type: we can form the type
$\id[{(\id[A]{x}{y})}]{p}{q}$ of identifications between
identifications $p,q$, and the type
$\id[{(\id[{(\id[A]{x}{y})}]{p}{q})}]{r}{s}$, and so on.  The structure
of this tower of identity types corresponds precisely to that of the
continuous paths and (higher) homotopies between them in a space, or an
$\infty$-groupoid.\index{.infinity-groupoid@$\infty$-groupoid}

Thus, we will frequently refer to an element $p : \id[A]{x}{y}$ as
a \define{path}
\index{path}
from $x$ to $y$; we call $x$ its \define{start point}
\indexdef{start point of a path}
\indexdef{path!start point of}
and $y$ its \define{end point}.
\indexdef{end point of a path}
\indexdef{path!end point of}
Two paths $p,q : \id[A]{x}{y}$ with the same start and end point are said to be \define{parallel},
\indexdef{parallel paths}
\indexdef{path!parallel}
in which case an element $r : \id[{(\id[A]{x}{y})}]{p}{q}$ can
be thought of as a homotopy, or a morphism between morphisms;
we will often refer to it as a \define{2-path}
\indexdef{path!2-}\indexsee{2-path}{path, 2-}%
or a \define{2-dimensional path}
\index{dimension!of paths}%
\indexsee{2-dimensional path}{path, 2-}\indexsee{path!2-dimensional}{path, 2-}%
Similarly, $\id[{(\id[{(\id[A]{x}{y})}]{p}{q})}]{r}{s}$ is the type of
\define{3-dimensional paths}
\indexdef{path!3-}\indexsee{3-path}{path, 3-}\indexsee{3-dimensional path}{path, 3-}\indexsee{path!3-dimensional}{path, 3-}%
between two parallel 2-dimensional paths, and so on.  If the
type $A$ is ``set-like'', such as \nat, these iterated identity types
will be uninteresting (see \autoref{sec:basics-sets}), but in the
general case they can model non-trivial homotopy types.


An important difference between homotopy type theory and classical homotopy theory is that homotopy type theory provides a \emph{synthetic}
\index{synthetic mathematics}%
\index{geometry, synthetic}%
\index{Euclid of Alexandria}%
description of spaces, in the following sense. Synthetic geometry is geometry in the style of Euclid~\cite{Euclid}: one starts from some basic notions (points and lines), constructions (a line connecting any two points), and axioms
(all right angles are equal), and deduces consequences logically.  This is in contrast with analytic
\index{analytic mathematics}%
geometry, where notions such as points and lines are represented concretely using cartesian coordinates in $\R^n$---lines are sets of points---and the basic constructions and axioms are derived from this representation.  While classical homotopy theory is analytic (spaces and paths are made of points), homotopy type theory is synthetic: points, paths, and paths between paths are basic, indivisible, primitive notions.

Moreover, one of the amazing things about homotopy type theory is that all of the basic constructions and axioms---all of the
higher groupoid structure----arises automatically from the induction
principle for identity types.
Recall from \autoref{sec:identity-types} that this says that if
\begin{itemize}
\item for every $x,y:A$ and every $p:\id[A]xy$ we have a type $D(x,y,p)$, and
\item for every $a:A$ we have an element $d(a):D(a,a,\refl a)$, 
\end{itemize}
then
\begin{itemize}
\item there exists an element $\indid{A}(D,d,x,y,p):D(x,y,p)$ for \emph{every} two elements $x,y:A$ and $p:\id[A]xy$, such that $\indid{A}(D,d,a,a,\refl a) \jdeq d(a)$.
\end{itemize}
In other words, given dependent functions
\begin{align*}
D & :\prd{x,y:A}{p:\id{x}{y}} \type\\
d & :\prd{a:A} D(a,a,\refl{a})
\end{align*}
there is a dependent function
\[\indid{A}(D,d):\prd{x,y:A}{p:\id{x}{y}} D(x,y,p)\]
such that 
\begin{equation}\label{eq:Jconv}
\indid{A}(D,d,a,a,\refl{a})\jdeq d(a)
\end{equation}
for every $a:A$.
Usually, every time we apply this induction rule we will either not care about the specific function being defined, or we will immediately give it a different name.

Informally, the induction principle for identity types says that if we want to construct an object (or prove a statement) which depends on an inhabitant $p:\id[A]xy$ of an identity type, then it suffices to perform the construction (or the proof) in the special case when $x$ and $y$ are the same (judgmentally) and $p$ is the reflexivity element $\refl{x}:x=x$ (judgmentally).
When writing informally, we may express this with a phrase such as ``by induction, it suffices to assume\dots''.
This reduction to the ``reflexivity case'' is analogous to the reduction to the ``base case'' and ``inductive step'' in an ordinary proof by induction on the natural numbers, and also to the ``left case'' and ``right case'' in a proof by case analysis on a disjoint union or disjunction.\index{induction principle!for identity type}%

The ``conversion rule''~\eqref{eq:Jconv} is less familiar in the context of proof by induction on natural numbers, but there is an analogous notion in the related concept of definition by recursion.
If a sequence\index{sequence} $(a_n)_{n\in \mathbb{N}}$ is defined by giving $a_0$ and specifying $a_{n+1}$ in terms of $a_n$, then in fact the $0^{\mathrm{th}}$ term of the resulting sequence \emph{is} the given one, and the given recurrence relation relating $a_{n+1}$ to $a_n$ holds for the resulting sequence.
(This may seem so obvious as to not be worth saying, but if we view a definition by recursion as an algorithm\index{algorithm} for calculating values of a sequence, then it is precisely the process of executing that algorithm.)
The rule~\eqref{eq:Jconv} is analogous: it says that if we define an object $f(p)$ for all $p:x=y$ by specifying what the value should be when $p$ is $\refl{x}:x=x$, then the value we specified is in fact the value of $f(\refl{x})$.

This induction principle endows each type with the structure of an $\infty$-groupoid\index{.infinity-groupoid@$\infty$-groupoid}, and each function between two types the structure of an $\infty$-functor\index{.infinity-functor@$\infty$-functor} between two such groupoids.  This is interesting from a mathematical point view, because it gives a new way to work with
$\infty$-groupoids.  It is interesting from a type-theoretic point view, because it reveals new operations that are associated with each type and function.  In the remainder of this chapter, we begin to explore this structure.

\section{Types are higher groupoids}
\label{sec:equality}

\index{type!identity|(}%
\index{path|(}%
\index{.infinity-groupoid@$\infty$-groupoid!structure of a type|(}%
We now derive from the induction principle the beginnings of the structure of a higher groupoid.
We begin with symmetry of equality, which, in topological language, means that ``paths can be reversed''.

\begin{lem}\label{lem:opp}
  For every type $A$ and every $x,y:A$ there is a function
  \begin{equation*}
    (x= y)\to(y= x)
  \end{equation*}
  denoted $p\mapsto \opp{p}$, such that $\opp{\refl{x}}\jdeq\refl{x}$ for each $x:A$.
  We call $\opp{p}$ the \define{inverse} of $p$.
  \indexdef{path!inverse}%
  \indexdef{inverse!of path}%
  \index{equality!symmetry of}
  \index{symmetry!of equality}%
\end{lem}

Since this is our first time stating something as a ``Lemma'' or ``Theorem'', let us pause to consider what that means.
Recall that propositions (statements susceptible to proof) are identified with types, whereas lemmas and theorems (statements that have been proven) are identified with \emph{inhabited} types.
Thus, the statement of a lemma or theorem should be translated into a type, as in \autoref{sec:pat}, and its proof translated into an inhabitant of that type.
According to the interpretation of the universal quantifier ``for every'', the type corresponding to \autoref{lem:opp} is
\[ \prd{A:\UU}{x,y:A} (x= y)\to(y= x). \]
The proof of \autoref{lem:opp} will consist of constructing an element of this type, i.e.\ deriving the judgment $f:\prd{A:\UU}{x,y:A} (x= y)\to(y= x)$ for some $f$.
We then introduce the notation $\opp{(\blank)}$ for this element $f$, in which the arguments $A$, $x$, and $y$ are omitted and inferred from context.
(As remarked in \autoref{sec:types-vs-sets}, the secondary statement ``$\opp{\refl{x}}\jdeq\refl{x}$ for each $x:A$'' should be regarded as a separate judgment.)

\begin{proof}[First proof]
  Assume given $A:\UU$, and
  let $D:\prd{x,y:A}{p:x= y} \type$ be the type family defined by $D(x,y,p)\defeq (y= x)$.
  In other words, $D$ is a function assigning to any $x,y:A$ and $p:x=y$ a type, namely the type $y=x$.
  Then we have an element
  \begin{equation*}
    d\defeq \lam{x} \refl{x}:\prd{x:A} D(x,x,\refl{x}).
  \end{equation*}
  Thus, the induction principle for identity types gives us an element
  \narrowequation{ \indid{A}(D,d,x,y,p): (y= x)}
  for each $p:(x= y)$.
  We can now define the desired function $\opp{(\blank)}$ to be $\lam{p} \indid{A}(D,d,x,y,p)$, i.e.\ we set $\opp{p} \defeq \indid{A}(D,d,x,y,p)$.
  The conversion rule~\eqref{eq:Jconv} gives $\opp{\refl{x}}\jdeq \refl{x}$, as required.
\end{proof}

We have written out this proof in a very formal style, which may be helpful while the induction rule on identity types is unfamiliar.
To be even more formal, we could say that \autoref{lem:opp} and its proof together consist of the judgment
\begin{narrowmultline*}
  \lam{A}{x}{y}{p} \indid{A}((\lam{x}{y}{p} (y=x)), (\lam{x} \refl{x}), x, y, p)
  \narrowbreak : \prd{A:\UU}{x,y:A} (x= y)\to(y= x)
\end{narrowmultline*}
(along with an additional equality judgment).
However, eventually we prefer to use more natural language, such as in the following equivalent proof.

\begin{proof}[Second proof]
  We want to construct, for each $x,y:A$ and $p:x=y$, an element $\opp{p}:y=x$.
  By induction, it suffices to do this in the case when $y$ is $x$ and $p$ is $\refl{x}$.
  But in this case, the type $x=y$ of $p$ and the type $y=x$ in which we are trying to construct $\opp{p}$ are both simply $x=x$.
  Thus, in the ``reflexivity case'', we can define $\opp{\refl{x}}$ to be simply $\refl{x}$.
  The general case then follows by the induction principle, and the conversion rule $\opp{\refl{x}}\jdeq\refl{x}$ is precisely the proof in the reflexivity case that we gave.
\end{proof}

We will write out the next few proofs in both styles, to help the reader become accustomed to the latter one.
Next we prove the transitivity of equality, or equivalently we ``concatenate paths''.

\begin{lem}\label{lem:concat}
  For every type $A$ and every $x,y,z:A$ there is a function
  \begin{equation*}
  (x= y) \to   (y= z)\to (x=  z)
  \end{equation*}
  written $p \mapsto q \mapsto p\ct q$, such that $\refl{x}\ct \refl{x}\jdeq \refl{x}$ for any $x:A$.
  We call $p\ct q$ the \define{concatenation} or \define{composite} of $p$ and $q$.
  \indexdef{path!concatenation}%
  \indexdef{path!composite}%
  \indexdef{concatenation of paths}%
  \indexdef{composition!of paths}%
  \index{equality!transitivity of}%
  \index{transitivity!of equality}%
\end{lem}

\begin{proof}[First proof]
  Let $D:\prd{x,y:A}{p:x=y} \type$ be the type family
  \begin{equation*}
    D(x,y,p)\defeq \prd{z:A}{q:y=z} (x=z).
  \end{equation*}
  Note that $D(x,x,\refl x) \jdeq \prd{z:A}{q:x=z} (x=z)$.
  Thus, in order to apply the induction principle for identity types to this $D$, we need a function of type
  \begin{equation}\label{eq:concatD}
    \prd{x:A} D(x,x,\refl{x})
  \end{equation}
  which is to say, of type
  \[ \prd{x,z:A}{q:x=z} (x=z). \]
  Now let $E:\prd{x,z:A}{q:x=z}\type$ be the type family $E(x,z,q)\defeq (x=z)$.
  Note that $E(x,x,\refl x) \jdeq (x=x)$.
  Thus, we have the function
  \begin{equation*}
    e(x) \defeq \refl{x} : E(x,x,\refl{x}).
  \end{equation*}
  By the induction principle for identity types applied to $E$, we obtain a function
  \begin{equation*}
    d(x,z,q) : \prd{x,z:A}{q:x=z} E(x,z,q).
  \end{equation*}
  But $E(x,z,q)\jdeq (x=z)$, so this is~\eqref{eq:concatD}.
  Thus, we can use this function $d$ and apply the induction principle for identity types to $D$, to obtain our desired function of type
  \begin{equation*}
    \prd{x,y,z:A}{q:y=z}{p:x=y} (x=z).
  \end{equation*}
  The conversion rules for the two induction principles give us $\refl{x}\ct \refl{x}\jdeq \refl{x}$ for any $x:A$.
\end{proof}

\begin{proof}[Second proof]
  We want to construct, for every $x,y,z:A$ and every $p:x=y$ and $q:y=z$, an element of $x=z$.
  By induction on $p$, it suffices to assume that $y$ is $x$ and $p$ is $\refl{x}$.
  In this case, the type $y=z$ of $q$ is $x=z$.
  Now by induction on $q$, it suffices to assume also that $z$ is $x$ and $q$ is $\refl{x}$.
  But in this case, $x=z$ is $x=x$, and we have $\refl{x}:(x=x)$.
\end{proof}

The reader may well feel that we have given an overly convoluted proof of this lemma.
In fact, we could stop after the induction on $p$, since at that point what we want to produce is an equality $x=z$, and we already have such an equality, namely $q$.
Why do we go on to do another induction on $q$?

The answer is that, as described in the introduction, we are doing \emph{proof-relevant} mathematics.
\index{mathematics!proof-relevant}%
When we prove a lemma, we are defining an inhabitant of some type, and it can matter what \emph{specific} element we defined in the course of the proof, not merely the type that that element inhabits (that is, the \emph{statement} of the lemma).
\autoref{lem:concat} has three obvious proofs: we could do induction over $p$, induction over $q$, or induction over both of them.
If we proved it three different ways, we would have three different elements of the same type.
It's not hard to show that these three elements are equal (see \autoref{ex:basics:concat}), but as they are not \emph{definitionally} equal, there can still be reasons to prefer one over another.

In the case of \autoref{lem:concat}, the difference hinges on the computation rule.
If we proved the lemma using a single induction over $p$, then we would end up with a computation rule of the form $\refl{y} \ct q \jdeq q$.
If we proved it with a single induction over $q$, we would have instead $p\ct\refl{x}\jdeq p$, while proving it with a double induction (as we did) gives only $\refl{x}\ct\refl{x} \jdeq \refl{x}$.

\index{mathematics!formalized}%
The asymmetrical computation rules can sometimes be convenient when doing formalized mathematics, as they allow the computer to simplify more things automatically.
However, in informal mathematics, and arguably even in the formalized case, it can be confusing to have a concatenation operation which behaves asymmetrically and to have to remember which side is the ``special'' one.
Treating both sides symmetrically makes for more robust proofs; this is why we have given the proof that we did.
(However, this is admittedly a stylistic choice.)

The table below summarizes the ``equality'', ``homotopical'', and ``higher-groupoid" points of view on what we have done so far.
\begin{center}
  \medskip
  \begin{tabular}{ccc}
    \toprule
    Equality & Homotopy & $\infty$-Groupoid\\
    \midrule
    reflexivity\index{equality!reflexivity of} & constant path & identity morphism\\
    symmetry\index{equality!symmetry of} & inversion of paths & inverse morphism\\
    transitivity\index{equality!transitivity of} & concatenation of paths & composition of morphisms\\
    \bottomrule
  \end{tabular}
  \medskip
\end{center}

In practice, transitivity is often applied to prove an equality by a chain of intermediate steps.
We will use the common notation for this such as $a=b=c=d$.
If the intermediate expressions are long, or we want to specify the witness of each equality, we may write
\begin{align*}
  a &= b & \text{(by $p$)}\\ &= c &\text{(by $q$)} \\ &= d &\text{(by $r$)}.
\end{align*}
In either case, the notation indicates construction of the element $(p\ct q)\ct r: (a=d)$.
(We choose left-associativity for concreteness, although in view of \autoref{thm:omg}\ref{item:omg4} below it makes litle difference.)
If it should happen that $b$ and $c$, say, are judgmentally equal, then we may write
\begin{align*}
  a &= b & \text{(by $p$)}\\ &\jdeq c \\ &= d &\text{(by $r$)}
\end{align*}
to indicate construction of $p\ct r : (a=d)$.

Now, because of proof-relevance, we can't stop after proving ``symmetry'' and ``transitivity'' of equality: we need to know that these \emph{operations} on equalities are well-behaved.
(This issue is invisible in set theory, where symmetry and transitivity are mere \emph{properties} of equality, rather than structure on
paths.)
From the homotopy-theoretic point of view, concatenation and inversion are just the ``first level'' of higher groupoid structure --- we also need coherence\index{coherence} laws on these operations, and analogous operations at higher dimensions.
For instance, we need to know that concatenation is \emph{associative}, and that inversion provides \emph{inverses} with respect to concatenation.

\begin{lem}\label{thm:omg}
  \index{associativity!of path concatenation}%
  \index{unit!law for path concatenation}%
  Suppose $A:\type$, that $x,y,z,w:A$ and that $p:x= y$ and $q:y = z$ and $r:z=w$.
  We have the following:
  \begin{enumerate}
  \item $p= p\ct \refl{y}$ and $p = \refl{x} \ct p$.\label{item:omg1}
  \item $\opp{p}\ct p=  \refl{y}$ and $p\ct \opp{p}= \refl{x}$.
  \item $\opp{(\opp{p})}= p$.
  \item $p\ct (q\ct r)=  (p\ct q)\ct r$.\label{item:omg4}
  \end{enumerate}
\end{lem}

Note, in particular, that \ref{item:omg1}--\ref{item:omg4} are themselves propositional equalities, living in the identity types \emph{of} identity types, such as $p=_{x=y}q$ for $p,q:x=y$.
Topologically, they are \emph{paths of paths}, i.e.\ homotopies.
It is a familiar fact in topology that when we concatenate a path $p$ with the reversed path $\opp p$, we don't literally obtain a constant path (which corresponds to the equality $\refl{}$ in type theory) --- instead we have a homotopy, or higher path, from $p\ct\opp p$ to the constant path.

\begin{proof}[Proof of~\autoref{thm:omg}]
  All the proofs use the induction principle for equalities.
  \begin{enumerate}
  \item \emph{First proof:} let $D:\prd{x,y:A}{p:x=y} \type$ be the type family given by 
    \begin{equation*}
      D(x,y,p)\defeq (p= p\ct \refl{y}).
    \end{equation*}
    Then $D(x,x,\refl{x})$ is $\refl{x}=\refl{x}\ct\refl{x}$.
    Since $\refl{x}\ct\refl{x}\jdeq\refl{x}$, it follows that $D(x,x,\refl{x})\jdeq (\refl{x}=\refl{x})$.
    Thus, there is a function
    \begin{equation*}
      d\defeq\lam{x} \refl{\refl{x}}:\prd{x:A} D(x,x,\refl{x}).
    \end{equation*}
    Now the induction principle for identity types gives an element $\indid{A}(D,d,p):(p= p\ct\refl{y})$ for each $p:x= y$.
    The other equality is proven similarly.

    \mentalpause

    \noindent
    \emph{Second proof:} by induction on $p$, it suffices to assume that $y$ is $x$ and that $p$ is $\refl x$.
    But in this case, we have $\refl{x}\ct\refl{x}\jdeq\refl{x}$.
  \item \emph{First proof:} let $D:\prd{x,y:A}{p:x=y} \type$ be the type family given by 
    \begin{equation*}
      D(x,y,p)\defeq (\opp{p}\ct p=  \refl{y}).
    \end{equation*}
    Then $D(x,x,\refl{x})$ is $\opp{\refl{x}}\ct\refl{x}=\refl{x}$.
    Since $\opp{\refl{x}}\jdeq\refl{x}$ and $\refl{x}\ct\refl{x}\jdeq\refl{x}$, we get that $D(x,x,\refl{x})\jdeq (\refl{x}=\refl{x})$.
    Hence we find the function
    \begin{equation*}
      d\defeq\lam{x} \refl{\refl{x}}:\prd{x:A} D(x,x,\refl{x}).
    \end{equation*}
    Now path induction gives an element $\indid{A}(D,d,p):\opp{p}\ct p=\refl{y}$ for each $p:x= y$ in $A$.
    The other equality is similar.

    \mentalpause

    \noindent \emph{Second proof} By induction, it suffices to assume $p$ is $\refl x$.
    But in this case, we have $\opp{p} \ct p \jdeq \opp{\refl x} \ct \refl x \jdeq \refl x$.

  \item \emph{First proof:} let $D:\prd{x,y:A}{p:x=y} \type$ be the type family given by
    \begin{equation*}
      D(x,y,p)\defeq (\opp{\opp{p}}= p).
    \end{equation*}
    Then $D(x,x,\refl{x})$ is the type $(\opp{\opp{\refl x}}=\refl{x})$.
    But since $\opp{\refl{x}}\jdeq \refl{x}$ for each $x:A$, we have $\opp{\opp{\refl{x}}}\jdeq \opp{\refl{x}} \jdeq\refl{x}$, and thus $D(x,x,\refl{x})\jdeq(\refl{x}=\refl{x})$.
    Hence we find the function
    \begin{equation*}
      d\defeq\lam{x} \refl{\refl{x}}:\prd{x:A} D(x,x,\refl{x}).
    \end{equation*}
    Now path induction gives an element $\indid{A}(D,d,p):\opp{\opp{p}}= p$ for each $p:x= y$.

    \mentalpause

    \noindent \emph{Second proof:} by induction, it suffices to assume $p$ is $\refl x$.
    But in this case, we have $\opp{\opp{p}}\jdeq \opp{\opp{\refl x}} \jdeq \refl x$.

  \item \emph{First proof:} let $D_1:\prd{x,y:A}{p:x=y} \type$ be the type family given by
    \begin{equation*}
      D_1(x,y,p)\defeq\prd{z,w:A}{q:y= z}{r:z= w} \big(p\ct (q\ct r)=  (p\ct q)\ct r\big).
    \end{equation*}
    Then $D_1(x,x,\refl{x})$ is
    \begin{equation*}
      \prd{z,w:A}{q:x= z}{r:z= w} \big(\refl{x}\ct(q\ct r)= (\refl{x}\ct q)\ct r\big).
    \end{equation*}
    To construct an element of this type, let $D_2:\prd{x,z:A}{q:x=z} \type$ be the type family
    \begin{equation*}
      D_2 (x,z,q) \defeq \prd{w:A}{r:z=w} \big(\refl{x}\ct(q\ct r)= (\refl{x}\ct q)\ct r\big).
    \end{equation*}
    Then $D_2(x,x,\refl{x})$ is
    \begin{equation*}
      \prd{w:A}{r:x=w} \big(\refl{x}\ct(\refl{x}\ct r)= (\refl{x}\ct \refl{x})\ct r\big).
    \end{equation*}
    To construct an element of \emph{this} type, let $D_3:\prd{x,w:A}{r:x=w} \type$ be the type family
    \begin{equation*}
      D_3(x,w,r) \defeq \big(\refl{x}\ct(\refl{x}\ct r)= (\refl{x}\ct \refl{x})\ct r\big).
    \end{equation*}
    Then $D_3(x,x,\refl{x})$ is
    \begin{equation*}
      \big(\refl{x}\ct(\refl{x}\ct \refl{x})= (\refl{x}\ct \refl{x})\ct \refl{x}\big)
    \end{equation*}
    which is definitionally equal to the type $(\refl{x} = \refl{x})$, and is therefore inhabited by $\refl{\refl{x}}$.
    Applying the path induction rule three times, therefore, we obtain an element of the overall desired type.

    \mentalpause

    \noindent \emph{Second proof:} by induction, it suffices to assume $p$, $q$, and $r$ are all $\refl x$.
    But in this case, we have
    \begin{align*}
      p\ct (q\ct r)
      &\jdeq \refl{x}\ct(\refl{x}\ct \refl{x})\\
      &\jdeq \refl{x}\\
      &\jdeq (\refl{x}\ct \refl x)\ct \refl x\\
      &\jdeq (p\ct q)\ct r.
    \end{align*}
    Thus, we have $\refl{\refl{x}}$ inhabiting this type. \qedhere
  \end{enumerate}
\end{proof}

\begin{rmk}
  There are other ways to define these higher paths.
  For instance, in \autoref{thm:omg}\ref{item:omg4} we might do induction only over one or two paths rather than all three.
  Each possibility will produce a \emph{definitionally} different proof, but they will all be equal to each other.
  Such an equality between any two particular proofs can, again, be proven by induction, reducing all the paths in question to reflexivities and then observing that both proofs reduce themselves to reflexivities.
\end{rmk}

In view of \autoref{thm:omg}\ref{item:omg4}, we will often write $p\ct q\ct r$ for $(p\ct q)\ct r$, and similarly $p\ct q\ct r \ct s$ for $((p\ct q)\ct r)\ct s$ and so on.
We choose left-associativity for definiteness, but it makes no real difference.
We generally trust the reader to insert instances of \autoref{thm:omg}\ref{item:omg4} to reassociate such expressions as necessary.

We are still not really done with the higher groupoid structure: the paths~\ref{item:omg1}--\ref{item:omg4} must also satisfy their own higher coherence\index{coherence} laws, which are themselves higher paths,
\index{associativity!of path concatenation!coherence of}%
\index{globular operad}%
\index{operad}%
\index{groupoid!higher}%
and so on ``all the way up to infinity'' (this can be made precise using e.g.\ the notion of a globular operad).
However, for most purposes it is unnecessary to make the whole infinite-dimensional structure explicit.
One of the nice things about homotopy type theory is that all of this structure can be \emph{proven} starting from only the inductive property of identity types, so we can make explicit as much or as little of it as we need.

In particular, in this book we will not need any of the complicated combinatorics involved in making precise notions such as ``coherent structure at all higher levels''.
In addition to ordinary paths, we will use paths of paths (i.e.\ elements of a type $p =_{x=_A y} q$), which as remarked previously we call \emph{2-paths}\index{path!2-} or \emph{2-dimensional paths}, and perhaps occasionally paths of paths of paths (i.e.\ elements of a type $r = _{p =_{x=_A y} q} s$), which we call \emph{3-paths}\index{path!3-} or \emph{3-dimensional paths}.
It is possible to define a general notion of \emph{$n$-dimensional path}
\indexdef{path!n-@$n$-}%
\indexsee{n-path@$n$-path}{path, $n$-}%
\indexsee{n-dimensional path@$n$-dimensional path}{path, $n$-}%
\indexsee{path!n-dimensional@$n$-dimensional}{path, $n$-}%
(see \autoref{ex:npaths}), but we will not need it.

We will, however, use one particularly important and simple case of higher paths, which is when the start and end points are the same.
In set theory, the proposition $a=a$ is entirely uninteresting, but in homotopy theory, paths from a point to itself are called \emph{loops}\index{loop} and carry lots of interesting higher structure.
Thus, given a type $A$ with a point $a:A$, we define its \define{loop space}
\index{loop space}%
$\Omega(A,a)$ to be the type $\id[A]{a}{a}$.
We may sometimes write simply $\Omega A$ if the point $a$ is understood from context.

Since any two elements of $\Omega A$ are paths with the same start and end points, they can be concatenated;
thus we have an operation $\Omega A\times \Omega A\to \Omega A$.
More generally, the higher groupoid structure of $A$ gives $\Omega A$ the analogous structure of a ``higher group''.

It can also be useful to consider the loop space\index{loop space!iterated}\index{iterated loop space} \emph{of} the loop space of $A$, which is the space of 2-dimensional loops on the identity loop at $a$.
This is written $\Omega^2(A,a)$ and represented in type theory by the type $\id[({\id[A]{a}{a}})]{\refl{a}}{\refl{a}}$.
While $\Omega^2(A,a)$, as a loop space, is again a ``higher group'', it now also has some additional structure resulting from the fact that its elements are 2-dimensional loops between 1-dimensional loops.  

\begin{thm}[Eckmann--Hilton]\label{thm:EckmannHilton}
  The composition operation on the second loop space
  \begin{equation*}
    \Omega^2(A)\times \Omega^2(A)\to \Omega^2(A)
  \end{equation*}
  is commutative: $\alpha\ct\beta = \beta\ct\alpha$, for any $\alpha, \beta:\Omega^2(A)$.
  \index{Eckmann--Hilton argument}%
\end{thm}

\begin{proof}
First, observe that the composition of $1$-loops $\Omega A\times \Omega A\to \Omega A$ induces an operation 
\[
\star : \Omega^2(A)\times \Omega^2(A)\to \Omega^2(A)
\]
as follows: consider elements $a, b, c : A$ and 1- and 2-paths,
\begin{align*}
  p &: a = b,       &       r &: b = c \\
  q &: a = b,       &       s &: b = c \\
  \alpha &: p = q,  &   \beta &: r = s
\end{align*}
as depicted in the following diagram (with paths drawn as arrows).
\[
 \xymatrix@+5em{
   {a} \rtwocell<10>^p_q{\alpha}
   &
   {b} \rtwocell<10>^r_s{\beta}
   &
   {c}
 }
\]
Composing the upper and lower 1-paths, respectively, we get two paths $p\ct r,\ q\ct s : a = c$, and there is then a ``horizontal composition''
\begin{equation*}
  \alpha\hct\beta : p\ct r = q\ct s
\end{equation*}
between them, defined as follows.
First, we define $\alpha \rightwhisker r : p\ct r = q\ct r$ by path induction on $r$, so that
\[ \alpha \rightwhisker \refl{b} \jdeq \opp{\mathsf{ru}_p} \ct \alpha \ct \mathsf{ru}_q \]
where $\mathsf{ru}_p : p = p \ct \refl{b}$ is the right unit law from \autoref{thm:omg}\ref{item:omg1}.
We could similarly define $\rightwhisker$ by induction on $\alpha$, or on all paths in sight, resulting in different judgmental equalities, but for present purposes the definition by induction on $r$ will make things simpler.
Similarly, we define $q\leftwhisker \beta : q\ct r = q\ct s$ by induction on $q$, so that
\[ \refl{b} \leftwhisker \beta \jdeq \opp{\mathsf{lu}_r} \ct \beta \ct \mathsf{lu}_s \]
where $\mathsf{lu}_r$ denotes the left unit law.
The operations $\leftwhisker$ and $\rightwhisker$ are called \define{whiskering}\indexdef{whiskering}.
Next, since $\alpha \rightwhisker r$ and $q\leftwhisker \beta$ are composable 2-paths, we can define the \define{horizontal composition}
\indexdef{horizontal composition!of paths}%
\indexdef{composition!of paths!horizontal}%
by:
\[
\alpha\hct\beta\ \defeq\ (\alpha\rightwhisker r) \ct (q\leftwhisker \beta).
\]
Now suppose that $a \jdeq  b \jdeq  c$, so that all the 1-paths $p$, $q$, $r$, and $s$ are elements of $\Omega(A,a)$, and assume moreover that $p\jdeq q \jdeq r \jdeq s\jdeq \refl{a}$, so that $\alpha:\refl{a} = \refl{a}$ and $\beta:\refl{a} = \refl{a}$ are composable in both orders.
In that case, we have
\begin{align*}
  \alpha\hct\beta
  &\jdeq (\alpha\rightwhisker\refl{a}) \ct (\refl{a}\leftwhisker \beta)\\
  &= \opp{\mathsf{ru}_{\refl{a}}} \ct \alpha \ct \mathsf{ru}_{\refl{a}} \ct \opp{\mathsf{lu}_{\refl a}} \ct \beta \ct \mathsf{lu}_{\refl{a}}\\
  &\jdeq \opp{\refl{\refl{a}}} \ct \alpha \ct \refl{\refl{a}} \ct \opp{\refl{\refl a}} \ct \beta \ct \refl{\refl{a}}\\
  &= \alpha \ct \beta.
\end{align*}
(Recall that $\mathsf{ru}_{\refl{a}} \jdeq \mathsf{lu}_{\refl{a}} \jdeq \refl{\refl{a}}$, by the computation rule for path induction.)
On the other hand, we can define another horizontal composition analogously by
\[
\alpha\hct'\beta\ \defeq\ (p\leftwhisker \beta)\ct (\alpha\rightwhisker s).
\]
and we similarly learn that
\[
\alpha\hct'\beta\ =\ (\refl{a}\leftwhisker \beta)\ct (\alpha\rightwhisker \refl{a}) = \beta\ct\alpha.
\]
\index{interchange law}%
But, in general, the two ways of defining horizontal composition agree, $\alpha\hct\beta = \alpha\hct'\beta$, as we can see by induction on $\alpha$ and $\beta$ and then on the two remaining 1-paths, to reduce everything to reflexivity.
Thus we have
\[\alpha \ct \beta = \alpha\hct\beta = \alpha\hct'\beta = \beta\ct\alpha.
\qedhere
\]
\end{proof}

The foregoing fact, which is known as the \emph{Eckmann--Hilton argument}, comes from classical homotopy theory, and indeed it is used in \autoref{cha:homotopy} below to show that the higher homotopy groups of a type are always abelian\index{group!abelian} groups.
The whiskering and horizontal composition operations defined in the proof are also a general part of the $\infty$-groupoid structure of types.
They satisfy their own laws (up to higher homotopy), such as
\[ \alpha \rightwhisker (p\ct q) = (\alpha \rightwhisker p) \rightwhisker q \]
and so on.
From now on, we trust the reader to apply path induction whenever needed to define further operations of this sort and verify their properties.

As this example suggests, the algebra of higher path types is much more intricate than just the groupoid-like structure at each level; the levels interact to give many further operations and laws, as in the study of iterated loop spaces in homotopy theory.
Indeed, as in classical homotopy theory, we can make the following general definitions:

\begin{defn} \label{def:pointedtype}
  A \define{pointed type}
  \indexsee{pointed!type}{type, pointed}%
  \indexdef{type!pointed}%
  $(A,a)$ is a type $A:\type$ together with a point $a:A$, called its \define{basepoint}.
  \indexdef{basepoint}%
  We write $\pointed{\type} \defeq \sm{A:\type} A$ for the type of pointed types in the universe $\type$.
\end{defn}

\begin{defn} \label{def:loopspace}
  Given a pointed type $(A,a)$, we define the \define{loop space}
  \indexdef{loop space}%
  of $(A,a)$ to be the following pointed type:
  \[\Omega(A,a)\defeq ((\id[A]aa),\refl a).\]
  An element of it will be called a \define{loop}\indexdef{loop} at $a$.
  For $n:\N$, the \define{$n$-fold iterated loop space} $\Omega^{n}(A,a)$
  \indexdef{loop space!iterated}%
  \indexsee{loop space!n-fold@$n$-fold}{loop space, iterated}%
  of a pointed type $(A,a)$ is defined recursively by:
  \begin{align*}
    \Omega^0(A,a)&\defeq(A,a)\\
    \Omega^{n+1}(A,a)&\defeq\Omega^n(\Omega(A,a)).
  \end{align*}
  An element of it will be called an \define{$n$-loop}
  \indexdef{loop!n-@$n$-}%
  \indexsee{n-loop@$n$-loop}{loop, $n$-}%
  or an \define{$n$-dimensional loop}
  \indexsee{loop!n-dimensional@$n$-dimensional}{loop, $n$-}%
  \indexsee{n-dimensional loop@$n$-dimensional loop}{loop, $n$-}%
  at $a$.
\end{defn}

We will return to iterated loop spaces in \autoref{cha:hlevels,cha:hits,cha:homotopy}.
\index{.infinity-groupoid@$\infty$-groupoid!structure of a type|)}%
\index{type!identity|)}
\index{path|)}%

\section{Functions are functors}
\label{sec:functors}

\index{function|(}%
\index{functoriality of functions in type theory@``functoriality'' of functions in type theory}%
Now we wish to establish that functions $f:A\to B$ behave functorially on paths.
In traditional type theory, this is equivalently the statement that functions respect equality.
\index{continuity of functions in type theory@``continuity'' of functions in type theory}%
Topologically, this corresponds to saying that every function is ``continuous'', i.e.\ preserves paths.

\begin{lem}\label{lem:map}
  Suppose that $f:A\to B$ is a function.
  Then for any $x,y:A$ there is an operation
  \begin{equation*}
    \apfunc f : (\id[A] x y) \to (\id[B] {f(x)} {f(y)}).
  \end{equation*}
  Moreover, for each $x:A$ we have $\apfunc{f}(\refl{x})\jdeq \refl{f(x)}$.
  \indexdef{application!of function to a path}%
  \indexdef{path!application of a function to}%
  \indexdef{function!application to a path of}%
  \indexdef{action!of a function on a path}%
\end{lem}

The notation $\apfunc f$ can be read either as the \underline{ap}plication of $f$ to a path, or as the \underline{a}ction on \underline{p}aths of $f$.

\begin{proof}[First proof]
  Let $D:\prd{x,y:A}{p:x=y}\type$ be the type family defined by
  \[D(x,y,p)\defeq (f(x)= f(y)).\]
  Then we have
  \begin{equation*}
    d\defeq\lam{x} \refl{f(x)}:\prd{x:A} D(x,x,\refl{x}).
  \end{equation*}
  By path induction, we obtain $\apfunc f : \prd{x,y:A}{p:x=y}(f(x)=g(x))$.
  The computation rule implies $\apfunc f({\refl{x}})\jdeq\refl{f(x)}$ for each $x:A$.
\end{proof}

\begin{proof}[Second proof]
  By induction, it suffices to assume $p$ is $\refl{x}$.
  In this case, we may define $\apfunc f(p) \defeq \refl{f(x)}:f(x)= f(x)$.
\end{proof}

We will often write $\apfunc f (p)$ as simply $\ap f p$.
This is strictly speaking ambiguous, but generally no confusion arises.
It matches the common convention in category theory of using the same symbol for the application of a functor to objects and to morphisms.

We note that $\apfunc{}$ behaves functorially, in all the ways that one might expect.

\begin{lem}\label{lem:ap-functor}
  For functions $f:A\to B$ and $g:B\to C$ and paths $p:\id[A]xy$ and $q:\id[A]yz$, we have:
  \begin{enumerate}
  \item $\apfunc f(p\ct q) = \apfunc f(p) \ct \apfunc f(q)$.\label{item:apfunctor-ct}
  \item $\apfunc f(\opp p) = \opp{\apfunc f (p)}$.\label{item:apfunctor-opp}
  \item $\apfunc g (\apfunc f(p)) = \apfunc{g\circ f} (p)$.\label{item:apfunctor-compose}
  \item $\apfunc {\idfunc[A]} (p) = p$.
  \end{enumerate}
\end{lem}
\begin{proof}
  Left to the reader.
\end{proof}
\index{function|)}%

As was the case for the equalities in \autoref{thm:omg}, those in \autoref{lem:ap-functor} are themselves paths, which satisfy their own coherence laws (which can be proved in the same way), and so on.

\section{Type families are fibrations}
\label{sec:fibrations}

\index{type!family of|(}%
\index{transport|(defstyle}%
Since \emph{dependently typed} functions are essential in type theory, we will also need a version of \autoref{lem:map} for these.
However, this is not quite so simple to state, because if $f:\prd{x:A} B(x)$ and $p:x=y$, then $f(x):B(x)$ and $f(y):B(y)$ are elements of distinct types, so that \emph{a priori} we cannot even ask whether they are equal.
The missing ingredient is that $p$ itself gives us a way to relate the types $B(x)$ and $B(y)$.

\begin{lem}[Transport]\label{lem:transport}
  Suppose that $P$ is a type family over $A$ and that $p:\id[A]xy$.
  Then there is a function $\transf{p}:P(x)\to P(y)$.
\end{lem}

\begin{proof}[First proof]
  Let $D:\prd{x,y:A}{p:\id{x}{y}} \type$ be the type family defined by
  \[D(x,y,p)\defeq P(x)\to P(y).\]
  Then we have the function
  \begin{equation*}
    d\defeq\lam{x} \idfunc[P(x)]:\prd{x:A} D(x,x,\refl{x}),
  \end{equation*}
  so that the induction principle gives us $\indid{A}(D,d,x,y,p):P(x)\to P(y)$ for $p:x= y$, which we define to be $\transf p$.
\end{proof}

\begin{proof}[Second proof]
  By induction, it suffices to assume $p$ is $\refl x$.
  But in this case, we can take $\transf{(\refl x)}:P(x)\to P(x)$ to be the identity function.
\end{proof}

Sometimes, it is necessary to notate the type family $P$ in which the transport operation happens.
In this case, we may write
\[\transfib P p \blank : P(x) \to P(y).\]

Recall that a type family $P$ over a type $A$ can be seen as a property of elements of $A$, which holds at $x$ in $A$ if $P(x)$ is inhabited.
Then the transportation lemma says that $P$ respects equality, in the sense that if $x$ is equal to $y$, then $P(x)$ holds if and only if $P(y)$ holds.
In fact, we will see later on that if $x=y$ then actually $P(x)$ and $P(y)$ are \emph{equivalent}.

Topologically, the transportation lemma can be viewed as a ``path lifting'' operation in a fibration.
\index{fibration}%
\indexdef{total!space}%
We think of a type family $P:A\to \type$ as a \emph{fibration} with base space $A$, with $P(x)$ being the fiber over $x$, and with $\sm{x:A}P(x)$ being the \define{total space} of the fibration, with first projection $\sm{x:A}P(x)\to A$.
The defining property of a fibration is that given a path $p:x=y$ in the base space $A$ and a point $u:P(x)$ in the fiber over $x$, we may lift the path $p$ to a path in the total space starting at $u$.
The point $\trans p u$ can be thought of as the other endpoint of this lifted path.
We can also define the path itself in type theory:

\begin{lem}[Path lifting property]\label{thm:path-lifting}
  \indexdef{path!lifting}%
  \indexdef{lifting!path}%
  Let $P:A\to\type$ be a type family over $A$ and assume we have $u:P(x)$ for some $x:A$.
  Then for any $p:x=y$, we have
  \begin{equation*}
    \mathsf{lift}(u,p):(x,u)=(y,\trans{p}{u})
  \end{equation*}
  in $\sm{x:A}P(x)$.
\end{lem}
\begin{proof}
  Left to the reader.
  We will prove a more general theorem in \autoref{sec:compute-sigma}.
\end{proof}

In classical homotopy theory, a fibration is defined as a map for which there \emph{exist} liftings of paths; while in contrast, we have just shown that in type theory, every type family comes with a \emph{specified} ``path-lifting function''.
This accords with the philosophy of constructive mathematics, according to which we cannot show that something exists except by exhibiting it.

\begin{rmk}
  Although we may think of a type family $P:A\to \type$ as like a fibration, it is generally not a good idea to say things like ``the fibration $P:A\to\type$'', since this sounds like we are talking about a fibration with base $\type$ and total space $A$.
  To repeat, when a type family $P:A\to \type$ is regarded as a fibration, the base is $A$ and the total space is $\sm{x:A} P(x)$.

  We may also occasionally use other topological terminology when speaking about type families.
  For instance, we may refer to a dependent function $f:\prd{x:A} P(x)$ as a \define{section}
  \indexdef{section!of a type family}%
  of the fibration $P$, and we may say that something happens \define{fiberwise}
  \indexdef{fiberwise}%
  if it happens for each $P(x)$.
  For instance, a section $f:\prd{x:A} P(x)$ shows that $P$ is ``fiberwise inhabited''.
\end{rmk}

\index{function!dependent|(}
Now we can prove the dependent version of \autoref{lem:map}.
The topological intuition is that given $f:\prd{x:A} P(x)$ and a path $p:\id[A]xy$, we ought to be able to apply $f$ to $p$ and obtain a path in the total space of $P$ which ``lies over'' $p$, as shown below.

\begin{center}
  \begin{tikzpicture}[yscale=.5,xscale=2]
    \draw (0,0) arc (-90:170:8ex) node[anchor=south east] {$A$} arc (170:270:8ex);
    \draw (0,6) arc (-90:170:8ex) node[anchor=south east] {$\sm{x:A} P(x)$} arc (170:270:8ex);
    \draw[->] (0,5.8) -- node[auto] {$\proj1$} (0,3.2);
    \node[circle,fill,inner sep=1pt,label=left:{$x$}] (b1) at (-.5,1.4) {};
    \node[circle,fill,inner sep=1pt,label=right:{$y$}] (b2) at (.5,1.4) {};
    \draw[decorate,decoration={snake,amplitude=1}] (b1) -- node[auto,swap] {$p$} (b2);
    \node[circle,fill,inner sep=1pt,label=left:{$f(x)$}] (b1) at (-.5,7.2) {};
    \node[circle,fill,inner sep=1pt,label=right:{$f(y)$}] (b2) at (.5,7.2) {};
    \draw[decorate,decoration={snake,amplitude=1}] (b1) -- node[auto] {$f(p)$} (b2);
  \end{tikzpicture}
\end{center}

We \emph{can} obtain such a thing from \autoref{lem:map}.
Given $f:\prd{x:A} P(x)$, we can define a non-dependent function $f':A\to \sm{x:A} P(x)$ by setting $f'(x)\defeq (x,f(x))$, and then consider $\ap{f'}{p} : f'(x) = f'(y)$.
However, it is not obvious from the type of such a path that it lies over a specific path in $A$ (in this case, $p$), which is sometimes important.

The solution is to use the transport lemma.
Since there is a canonical path from $u:P(x)$ to $\trans p u :P(y)$ which (at least intuitively) lies over $p$, any path from $u$ to $v:P(y)$ lying over $p$ should factor through this path, essentially uniquely, by a path from $\trans p u$ to $v$ lying entirely in the fiber $P(y)$.
Thus, up to equivalence, it makes sense to define ``a path from $u$ to $v$ lying over $p:x=y$'' to mean a path $\trans p u = v$ in $P(y)$.
And, indeed, we can show that dependent functions produce such paths.

\begin{lem}[Dependent map]\label{lem:mapdep}
  \indexdef{application!of dependent function to a path}%
  \indexdef{path!application of a dependent function to}%
  \indexdef{function!dependent!application to a path of}%
  \indexdef{action!of a dependent function on a path}%
  Suppose $f:\prd{x: A} P(x)$; then we have a map
  \[\apdfunc f : \prd{p:x=y}\big(\id[P(y)]{\trans p{f(x)}}{f(y)}\big).\]
\end{lem}

\begin{proof}[First proof]
  Let $D:\prd{x,y:A}{p:\id{x}{y}} \type$ be the type family defined by
  \begin{equation*}
    D(x,y,p)\defeq \trans p {f(x)}= f(y).
  \end{equation*}
  Then $D(x,x,\refl{x})$ is $\trans{(\refl{x})}{f(x)}= f(x)$.
  But since $\trans{(\refl{x})}{f(x)}\jdeq f(x)$, we get that $D(x,x,\refl{x})\jdeq (f(x)= f(x))$.
  Thus, we find the function
  \begin{equation*}
    d\defeq\lam{x} \refl{f(x)}:\prd{x:A} D(x,x,\refl{x})
  \end{equation*}
  and now path induction gives us $\apdfunc f(p):\trans p{f(x)}= f(y)$ for each $p:x= y$.
\end{proof}

\begin{proof}[Second proof]
  By induction, it suffices to assume $p$ is $\refl x$.
  But in this case, the desired equation is $\trans{(\refl{x})}{f(x)}\jdeq f(x)$, which holds judgmentally.
\end{proof}

We will refer generally to paths which ``lie over other paths'' in this sense as \emph{dependent paths}.
\indexsee{dependent!path}{path, dependent}%
\index{path!dependent}%
They will play an increasingly important role starting in \autoref{cha:hits}.
In \autoref{sec:computational} we will see that for a few particular kinds of type families, there are equivalent ways to represent the notion of dependent paths that are sometimes more convenient.

Now recall from \autoref{sec:pi-types} that a non-dependently typed function $f:A\to B$ is just the special case of a dependently typed function $f:\prd{x:A} P(x)$ when $P$ is a constant type family, $P(x) \defeq B$.
In this case, $\apdfunc{f}$ and $\apfunc{f}$ are closely related, because of the following lemma:

\begin{lem}\label{thm:trans-trivial}
  If $P:A\to\type$ is defined by $P(x) \defeq B$ for a fixed $B:\type$, then for any $x,y:A$ and $p:x=y$ and $b:B$ we have a path
  \[ \transconst Bpb : \transfib P p b = b. \]
\end{lem}
\begin{proof}[First proof]
  Fix a $b:B$, and let $D:\prd{x,y:A}{p:\id{x}{y}} \type$ be the type family defined by
  \[ D(x,y,p) \defeq (\transfib P p b = b). \]
  Then $D(x,x,\refl x)$ is $(\transfib P{\refl{x}}{b} = b)$, which is judgmentally equal to $(b=b)$ by the computation rule for transporting.
  Thus, we have the function
  \[ d \defeq \lam{x} \refl{b} : \prd{x:A} D(x,x,\refl x). \]
  Now path induction gives us an element of
  \narrowequation{
    \prd{x,y:A}{p:x=y}(\transfib P p b = b),}
  as desired.
\end{proof}
\begin{proof}[Second proof]
  By induction, it suffices to assume $y$ is $x$ and $p$ is $\refl x$.
  But $\transfib P {\refl x} b \jdeq b$, so in this case what we have to prove is $b=b$, and we have $\refl{b}$ for this.
\end{proof}

Thus, by concatenating with $\transconst B p b$, for any $x,y:A$ and $p:x=y$ and $f:A\to B$ we obtain functions
\begin{align}
  \big(f(x) = f(y)\big) &\to \big(\trans{p}{f(x)} = f(y)\big)\label{eq:ap-to-apd}
  \qquad\text{and} \\
  \big(\trans{p}{f(x)} = f(y)\big) &\to \big(f(x) = f(y)\big).\label{eq:apd-to-ap}
\end{align}
In fact, these functions are inverse equivalences (in the sense to be introduced in \autoref{sec:basics-equivalences}), and they relate $\apfunc f (p)$  to $\apdfunc f (p)$.

\begin{lem}\label{thm:apd-const}
  For $f:A\to B$ and $p:\id[A]xy$, we have
  \[ \apdfunc f(p) = \transconst B p{f(x)} \ct \apfunc f (p). \]
\end{lem}
\begin{proof}[First proof]
  Let $D:\prd{x,y:A}{p:\id xy} \type$ be the type family defined by
  \[ D(x,y,p) \defeq \big(\apdfunc f (p) = \transconst Bp{f(x)} \ct \apfunc f (p)\big). \]
  Thus, we have
  \[D(x,x,\refl x) \jdeq \big(\apdfunc f (\refl x) = \transconst B{\refl x}{f(x)} \ct \apfunc f ({\refl x})\big).\]
  But by definition, all three paths appearing in this type are $\refl{f(x)}$, so we have
  \[ \refl{\refl{f(x)}} : D(x,x,\refl x). \]
  Thus, path induction gives us an element of $\prd{x,y:A}{p:x=y} D(x,y,p)$, which is what we wanted.
\end{proof}
\begin{proof}[Second proof]
  By induction, it suffices to assume $y$ is $x$ and $p$ is $\refl x$.
  In this case, what we have to prove is $\refl{f(x)} = \refl{f(x)} \ct \refl{f(x)}$, which is true judgmentally.
\end{proof}

Because the types of $\apdfunc{f}$ and $\apfunc{f}$ are different, it is often clearer to use different notations for them.

\index{function!dependent|)}%

At this point, we hope the reader is starting to get a feel for proofs by induction on identity types.
From now on we stop giving both styles of proofs, allowing ourselves to use whatever is most clear and convenient (and often the second, more concise one).
Here are a few other useful lemmas about transport; we leave it to the reader to give the proofs (in either style).

\begin{lem}\label{thm:transport-concat}
  Given $P:A\to\type$ with $p:\id[A]xy$ and $q:\id[A]yz$ while $u:P(x)$, we have
  \[ \trans{q}{\trans{p}{u}} = \trans{(p\ct q)}{u}. \]
\end{lem}

\begin{lem}\label{thm:transport-compose}
  For a function $f:A\to B$ and a type family $P:B\to\type$, and any $p:\id[A]xy$ and $u:P(f(x))$, we have
  \[ \transfib{P\circ f}{p}{u} = \transfib{P}{\apfunc f(p)}{u}. \]
\end{lem}

\begin{lem}\label{thm:ap-transport}
  For $P,Q:A\to \type$ and a family of functions $f:\prd{x:A} P(x)\to Q(x)$, and any $p:\id[A]xy$ and $u:P(x)$, we have
  \[ \transfib{Q}{p}{f_x(u)} = f_y(\transfib{P}{p}{u}). \]
\end{lem}

\index{type!family of|)}%
\index{transport|)}

\section{Homotopies and equivalences}
\label{sec:basics-equivalences}

\index{homotopy|(defstyle}%

So far, we have seen how the identity type $\id[A]xy$ can be regarded as a type of \emph{identifications}, \emph{paths}, or \emph{equivalences} between two elements $x$ and $y$ of a type $A$.
Now we investigate the appropriate notions of ``identification'' or ``sameness'' between \emph{functions} and between \emph{types}.
In \autoref{sec:compute-pi,sec:compute-universe}, we will see that homotopy type theory allows us to identify these with instances of the identity type, but before we can do that we need to understand them in their own right.

Traditionally, we regard two functions as the same if they take equal values on all inputs.
Under the propositions-as-types interpretation, this suggests that two functions $f$ and $g$ (perhaps dependently typed) should be the same if the type $\prd{x:A} (f(x)=g(x))$ is inhabited.
Under the homotopical interpretation, this dependent function type consists of \emph{continuous} paths or \emph{functorial} equivalences, and thus may be regarded as the type of \emph{homotopies} or of \emph{natural isomorphisms}.\index{isomorphism!natural}%
We will adopt the topological terminology for this.

\begin{defn} \label{defn:homotopy}
  Let $f,g:\prd{x:A} P(x)$ be two sections of a type family $P:A\to\type$.
  A \define{homotopy}
  from $f$ to $g$ is a dependent function of type
  \begin{equation*}
    (f\htpy g) \defeq \prd{x:A} (f(x)=g(x)).
  \end{equation*}
\end{defn}

Note that a homotopy is not the same as an identification $(f=g)$.
However, in \autoref{sec:compute-pi} we will introduce an axiom making homotopies and identifications ``equivalent''.

The following proofs are left to the reader.

\begin{lem}\label{lem:homotopy-props}
  Homotopy is an equivalence relation on each function type $A\to B$.
  That is, we have elements of the types
  \begin{gather*}
    \prd{f:A\to B} (f\htpy f)\\
    \prd{f,g:A\to B} (f\htpy g) \to (g\htpy f)\\
    \prd{f,g,h:A\to B} (f\htpy g) \to (g\htpy h) \to (f\htpy h).
  \end{gather*}
\end{lem}


\index{functoriality of functions in type theory@``functoriality'' of functions in type theory}%
\index{continuity of functions in type theory@``continuity'' of functions in type theory}%
Just as functions in type theory are automatically ``functors'', homotopies are automatically
\index{naturality of homotopies@``naturality'' of homotopies}%
``natural transformations'', in the following sense.
Recall that for $f:A\to B$ and $p:\id[A]xy$, we may write $\ap f p$ to mean $\apfunc{f} (p)$.

\begin{lem}\label{lem:htpy-natural}
  Suppose $H:f\htpy g$ is a homotopy between functions $f,g:A\to B$ and let $p:\id[A]xy$.  Then we have
  \begin{equation*}
    H(x)\ct\ap{g}{p}=\ap{f}{p}\ct H(y).
  \end{equation*}
  We may also draw this as a commutative diagram:\index{diagram}
  \begin{align*}
    \xymatrix{
      f(x) \ar@{=}[r]^{\ap fp} \ar@{=}[d]_{H(x)} & f(y) \ar@{=}[d]^{H(y)} \\
      g(x) \ar@{=}[r]_{\ap gp} & g(y)
    }
  \end{align*}
\end{lem}
\begin{proof}
  By induction, we may assume $p$ is $\refl x$.
  Since $\apfunc{f}$ and $\apfunc g$ compute on reflexivity, in this case what we must show is
  \[ H(x) \ct \refl{g(x)} = \refl{f(x)} \ct H(x). \]
  But this follows since both sides are equal to $H(x)$.
\end{proof}

\begin{cor}\label{cor:hom-fg}
  Let $H : f \htpy \idfunc[A]$ be a homotopy, with $f : A \to A$. Then for any $x : A$ we have \[ H(f(x)) = \ap f{H(x)}. \]
\end{cor}
\noindent
Here $f(x)$ denotes the ordinary application of $f$ to $x$, while $\ap f{H(x)}$ denotes $\apfunc{f}(H(x))$.
\begin{proof}
By naturality of $H$, the following diagram of paths commutes:
\begin{align*}
\xymatrix@C=3pc{
ffx \ar@{=}[r]^-{\ap f{Hx}} \ar@{=}[d]_{H(fx)} & fx \ar@{=}[d]^{Hx} \\
fx \ar@{=}[r]_-{Hx} & x
}
\end{align*}
That is, $\ap f{H x} \ct H x = H(f x) \ct H x$.
We can now whisker by $\opp{(H x)}$ to cancel $H x$, obtaining
\[ \ap f{H x}
= \ap f{H x} \ct H x \ct \opp{(H x)}
= H(f x) \ct H x \ct \opp{(H x)}
= H(f x)
\]
as desired (with some associativity paths suppressed).
\end{proof}

Of course, like the functoriality of functions (\autoref{lem:ap-functor}), the equality in \autoref{lem:htpy-natural} is a path which satisfies its own coherence laws, and so on.

\index{homotopy|)}%

\index{equivalence|(}%
Moving on to types, from a traditional perspective one may say that a function $f:A\to B$ is an \emph{isomorphism} if there is a function $g:B\to A$ such that both composites $f\circ g$ and $g\circ f$ are pointwise equal to the identity, i.e.\ such that $f \circ g \htpy \idfunc[B]$ and $g\circ f \htpy \idfunc[A]$.
\indexsee{homotopy!equivalence}{equivalence}%
A homotopical perspective suggests that this should be called a \emph{homotopy equivalence}, and from a categorical one, it should be called an \emph{equivalence of (higher) groupoids}.
However, when doing proof-relevant mathematics,
\index{mathematics!proof-relevant}%
the corresponding type
\begin{equation}
  \sm{g:B\to A} \big((f \circ g \htpy \idfunc[B]) \times (g\circ f \htpy \idfunc[A])\big)\label{eq:qinvtype}
\end{equation}
is poorly behaved.
For instance, for a single function $f:A\to B$ there may be multiple unequal inhabitants of~\eqref{eq:qinvtype}.
(This is closely related to the observation in higher category theory that often one needs to consider \emph{adjoint} equivalences\index{adjoint!equivalence} rather than plain equivalences.)
For this reason, we give~\eqref{eq:qinvtype} the following historically accurate, but slightly de\-rog\-a\-to\-ry-sounding name instead.

\begin{defn}
  For a function $f:A\to B$, a \define{quasi-inverse}
  \indexdef{quasi-inverse}%
  \indexsee{function!quasi-inverse of}{quasi-inverse}%
  of $f$ is a triple $(g,\alpha,\beta)$ consisting of a function $g:B\to A$ and homotopies
$\alpha:f\circ g\htpy \idfunc[B]$ and $\beta:g\circ f\htpy \idfunc[A]$.
\end{defn}

\symlabel{qinv}
Thus,~\eqref{eq:qinvtype} is \emph{the type of quasi-inverses of $f$}; we may denote it by $\qinv(f)$.

\begin{eg}\label{eg:idequiv}
  \index{identity!function}%
  \index{function!identity}%
  The identity function $\idfunc[A]:A\to A$ has a quasi-inverse given by $\idfunc[A]$ itself, together with homotopies defined by $\alpha(y) \defeq \refl{y}$ and $\beta(x) \defeq \refl{x}$.
\end{eg}

\begin{eg}\label{eg:concatequiv}
  For any $p:\id[A]xy$ and $z:A$, the functions
  \begin{align*}
    (p\ct \blank)&:(\id[A]yz) \to (\id[A]xz) \qquad\text{and}\\
    (\blank \ct p)&:(\id[A]zx) \to (\id[A]zy)
  \end{align*}
  have quasi-inverses given by $(\opp p \ct \blank)$ and $(\blank \ct \opp p)$, respectively; see \autoref{ex:equiv-concat}.
\end{eg}

\begin{eg}\label{thm:transportequiv}
  For any $p:\id[A]xy$ and $P:A\to\type$, the function
  \[\transfib{P}{p}{\blank}:P(x) \to P(y)\]
  has a quasi-inverse given by $\transfib{P}{\opp p}{\blank}$; this follows from \autoref{thm:transport-concat}.
\end{eg}

\symlabel{basics-isequiv}
In general, we will only use the word \emph{isomorphism}
\index{isomorphism!of sets}
(and similar words such as \emph{bijection})
\index{bijection}
in the special case when the types $A$ and $B$ ``behave like sets'' (see \autoref{sec:basics-sets}).
In this case, the type~\eqref{eq:qinvtype} is unproblematic.
We will reserve the word \emph{equivalence} for an improved notion $\isequiv (f)$ with the following properties:%
\begin{enumerate}
\item For each $f:A\to B$ there is a function $\qinv(f) \to \isequiv (f)$.\label{item:be1}
\item Similarly, for each $f$ we have $\isequiv (f) \to \qinv(f)$; thus the two are logically equivalent (see \autoref{sec:pat}).\label{item:be2}
\item For any two inhabitants $e_1,e_2:\isequiv(f)$ we have $e_1=e_2$.\label{item:be3}
\end{enumerate}
In \autoref{cha:equivalences} we will see that there are many different definitions of $\isequiv(f)$ which satisfy these three properties, but that all of them are equivalent.
For now, to convince the reader that such things exist, we mention only the easiest such definition:
\begin{equation}\label{eq:isequiv-invertible}
  \isequiv(f) \;\defeq\;
  \Parens{\sm{g:B\to A} (f\circ g \htpy \idfunc[B])}
  \times
  \Parens{\sm{h:B\to A} (h\circ f \htpy \idfunc[A])}.
\end{equation}
We can show~\ref{item:be1} and~\ref{item:be2} for this definition now.
A function $\qinv(f) \to \isequiv (f)$ is easy to define by taking $(g,\alpha,\beta)$ to $(g,\alpha,g,\beta)$.
In the other direction, given $(g,\alpha,h,\beta)$, let $\gamma$ be the composite homotopy
\[ g \overset{\beta}{\htpy} h\circ f\circ g \overset{\alpha}{\htpy} h \]
and let $\beta':g\circ f\htpy \idfunc[A]$ be obtained from $\gamma$ and $\beta$.
Then $(g,\alpha,\beta'):\qinv(f)$.

Property~\ref{item:be3} for this definition is not too hard to prove either, but it requires identifying the identity types of cartesian products and dependent pair types, which we will discuss in \autoref{sec:compute-cartprod,sec:compute-sigma}.
Thus, we postpone it as well; see \autoref{sec:biinv}.
At this point, the main thing to take away is that there is a well-behaved type which we can pronounce as ``$f$ is an equivalence'', and that we can prove $f$ to be an equivalence by exhibiting a quasi-inverse to it.
In practice, this is the most common way to prove that a function is an equivalence.

In accord with the proof-relevant philosophy,
\index{mathematics!proof-relevant}%
\emph{an equivalence} from $A$ to $B$ is defined to be a function $f:A\to B$ together with an inhabitant of $\isequiv (f)$, i.e.\ a proof that it is an equivalence.
We write $(\eqv A B)$ for the type of equivalences from $A$ to $B$, i.e.\ the type
\begin{equation}\label{eq:eqv}
  (\eqv A B) \defeq \sm{f:A\to B} \isequiv(f).
\end{equation}
Property~\ref{item:be3} above will ensure that if two equivalences are equal as functions (that is, the underlying elements of $A\to B$ are equal), then they are also equal as equivalences (see \autoref{sec:compute-sigma}).
Thus, we often abuse notation by denoting an equivalence by the same letter as its underlying function.

We conclude by observing:

\begin{lem}\label{thm:equiv-eqrel}
  Type equivalence is an equivalence relation on \type.
  More specifically:
  \begin{enumerate}
  \item For any $A$, the identity function $\idfunc[A]$ is an equivalence; hence $\eqv A A$.
  \item For any $f:\eqv A B$, we have an equivalence $f^{-1} : \eqv B A$.
  \item For any $f:\eqv A B$ and $g:\eqv B C$, we have $g\circ f : \eqv A C$.
  \end{enumerate}
\end{lem}
\begin{proof}
  The identity function is clearly its own quasi-inverse; hence it is an equivalence.

  If $f:A\to B$ is an equivalence, then it has a quasi-inverse, say $f^{-1}:B\to A$.
  Then $f$ is also a quasi-inverse of $f^{-1}$, so $f^{-1}$ is an equivalence $B\to A$.

  Finally, given $f:\eqv A B$ and $g:\eqv B C$ with quasi-inverses $f^{-1}$ and $g^{-1}$, say, then for any $a:A$ we have $f^{-1} g^{-1} g f a = f^{-1} f a = a$, and for any $c:C$ we have $g f f^{-1} g^{-1} c = g g^{-1} c = c$.
  Thus $f^{-1} \circ g^{-1}$ is a quasi-inverse to $g\circ f$, hence the latter is an equivalence.
\end{proof}

\index{equivalence|)}%

\section{The higher groupoid structure of type formers}
\label{sec:computational}

In \autoref{cha:typetheory}, we introduced many ways to form new types: cartesian products, disjoint unions, dependent products, dependent sums, etc.
In \autoref{sec:equality,sec:functors,sec:fibrations}, we saw that \emph{all} types in homotopy type theory behave like spaces or higher groupoids.
Our goal in the rest of the chapter is to make explicit how this higher structure behaves in the case of the particular types defined in \autoref{cha:typetheory}.

It turns out that for many types $A$, the equality types $\id[A]xy$ can be characterized, up to equivalence, in terms of whatever data was used to construct $A$.
For example, if $A$ is a cartesian product $B\times C$, and $x\jdeq (b,c)$ and $y\jdeq(b',c')$, then we have an equivalence
\begin{equation}\label{eq:prodeqv}
  \eqv{\big((b,c)=(b',c')\big)}{\big((b=b')\times (c=c')\big)}.
\end{equation}
In more traditional language, two ordered pairs are equal just when their components are equal (but the equivalence~\eqref{eq:prodeqv} says rather more than this).
The higher structure of the identity types can also be expressed in terms of these equivalences; for instance, concatenating two equalities between pairs corresponds to pairwise concatenation.

Similarly, when a type family $P:A\to\type$ is built up fiberwise using the type forming rules from \autoref{cha:typetheory}, the operation $\transfib{P}{p}{\blank}$ can be characterized, up to homotopy, in terms of the corresponding operations on the data that went into $P$.
For instance, if $P(x) \jdeq B(x)\times C(x)$, then we have
\[\transfib{P}{p}{(b,c)} = \big(\transfib{B}{p}{b},\transfib{C}{p}{c}\big).\]

Finally, the type forming rules are also functorial, and if a function $f$ is built from this functoriality, then the operations $\apfunc f$ and $\apdfunc f$ can be computed based on the corresponding ones on the data going into $f$.
For instance, if $g:B\to B'$ and $h:C\to C'$ and we define $f:B\times C \to B'\times C'$ by $f(b,c)\defeq (g(b),h(c))$, then modulo the equivalence~\eqref{eq:prodeqv}, we can identify $\apfunc f$ with ``$(\apfunc g,\apfunc h)$''.

The next few sections (\crefrange{sec:compute-cartprod}{sec:compute-nat}) will be devoted to stating and proving theorems of this sort for all the basic type forming rules, with one section for each basic type former.
Here we encounter a certain apparent deficiency in currently available type theories;
as will become clear in later chapters, it would seem to be more convenient and intuitive if these characterizations of identity types, transport, and so on were \emph{judgmental}\index{judgmental equality} equalities.
However, in the theory presented in \autoref{cha:typetheory}, the identity types are defined uniformly for all types by their induction principle, so we cannot ``redefine'' them to be different things at different types.
Thus, the characterizations for particular types to be discussed in this chapter are, for the most part, \emph{theorems} which we have to discover and prove, if possible.

Actually, the type theory of \autoref{cha:typetheory} is insufficient to prove the desired theorems for two of the type formers: $\Pi$-types and universes.
For this reason, we are forced to introduce axioms into our type theory, in order to make those ``theorems'' true.
Type-theoretically, an \emph{axiom} (c.f.~\autoref{sec:axioms}) is an ``atomic'' element that is declared to inhabit some specified type, without there being any rules governing its behavior other than those pertaining to the type it inhabits.
\index{axiom!versus rules}%

\index{function extensionality}%
\indexsee{extensionality, of functions}{function extensionality}
\index{univalence axiom}%
The axiom for $\Pi$-types (\autoref{sec:compute-pi}) is familiar to type theorists: it is called \emph{function extensionality}, and states (roughly) that if two functions are homotopic in the sense of \autoref{sec:basics-equivalences}, then they are equal.
The axiom for universes (\autoref{sec:compute-universe}), however, is a new contribution of homotopy type theory due to Voevodsky: it is called the \emph{univalence axiom}, and states (roughly) that if two types are equivalent in the sense of \autoref{sec:basics-equivalences}, then they are equal.
We have already remarked on this axiom in the introduction; it will play a very important role in this book.%
\footnote{We have chosen to introduce these principles as axioms, but there are potentially other ways to formulate a type theory in which they hold.
  See the Notes to this chapter.}

It is important to note that not \emph{all} identity types can be ``determined'' by induction over the construction of types.
Counterexamples include most nontrivial higher inductive types (see \autoref{cha:hits,cha:homotopy}).
For instance, calculating the identity types of the types $\Sn^n$ (see \autoref{sec:circle}) is equivalent to calculating the higher homotopy groups of spheres, a deep and important field of research in algebraic topology.

\section{Cartesian product types}
\label{sec:compute-cartprod}

\index{type!product|(}%
Given types $A$ and $B$, consider the cartesian product type $A \times B$.  
For any elements $x,y:A\times B$ and a path $p:\id[A\times B]{x}{y}$, by functoriality we can extract paths $\ap{\proj1}p:\id[A]{\proj1(x)}{\proj1(y)}$ and $\ap{\proj2}p:\id[B]{\proj2(x)}{\proj2(y)}$.
Thus, we have a function
\begin{equation}\label{eq:path-prod}
  (\id[A\times B]{x}{y}) \to (\id[A]{\proj1(x)}{\proj1(y)}) \times (\id[B]{\proj2(x)}{\proj2(y)}).
\end{equation}

\begin{thm}\label{thm:path-prod}
  For any $x$ and $y$, the function~\eqref{eq:path-prod} is an equivalence.
\end{thm}

Read logically, this says that two pairs are equal if they are equal
componentwise.  Read category-theoretically, this says that the
morphisms in a product groupoid are pairs of morphisms.  Read
homotopy-theoretically, this says that the paths in a product
space are pairs of paths.

\begin{proof}
  We need a function in the other direction:
  \begin{equation}
    (\id[A]{\proj1(x)}{\proj1(y)}) \times (\id[B]{\proj2(x)}{\proj2(y)}) \to (\id[A\times B]{x}{y}). \label{eq:path-prod-inverse}
  \end{equation}
  By the induction rule for cartesian products, we may assume that $x$ and $y$ are both pairs, i.e.\ $x\jdeq (a,b)$ and $y\jdeq (a',b')$ for some $a,a':A$ and $b,b':B$.
  In this case, what we want is a function
  \begin{equation*}
    (\id[A]{a}{a'}) \times (\id[B]{b}{b'}) \to \big(\id[A\times B]{(a,b)}{(a',b')}\big).
  \end{equation*}
  Now by induction for the cartesian product in its domain, we may assume given $p:a=a'$ and $q:b=b'$.
  And by two path inductions, we may assume that $a\jdeq a'$ and $b\jdeq b'$ and both $p$ and $q$ are reflexivity.
  But in this case, we have $(a,b)\jdeq(a',b')$ and so we can take the output to also be reflexivity.

  It remains to prove that~\eqref{eq:path-prod-inverse} is quasi-inverse to~\eqref{eq:path-prod}.
  This is a simple sequence of inductions, but they have to be done in the right order.

  In one direction, let us start with $r:\id[A\times B]{x}{y}$.
  We first do a path induction on $r$ in order to assume that $x\jdeq y$ and $r$ is reflexivity.
  In this case, since $\apfunc{\proj1}$ and $\apfunc{\proj2}$ are defined by path induction,~\eqref{eq:path-prod} takes $r\jdeq \refl{x}$ to the pair $(\refl{\proj1x},\refl{\proj2x})$.
  Now by induction on $x$, we may assume $x\jdeq (a,b)$, so that this is $(\refl a, \refl b)$.
  Thus,~\eqref{eq:path-prod-inverse} takes it by definition to $\refl{(a,b)}$, which (under our current assumptions) is $r$.
  
  In the other direction, if we start with $s:(\id[A]{\proj1(x)}{\proj1(y)}) \times (\id[B]{\proj2(x)}{\proj2(y)})$, then we first do induction on $x$ and $y$ to assume that they are pairs $(a,b)$ and $(a',b')$, and then induction on $s:(\id[A]{a}{a'}) \times (\id[B]{b}{b'})$ to reduce it to a pair $(p,q)$ where $p:a=a'$ and $q:b=b'$.
  Now by induction on $p$ and $q$, we may assume they are reflexivities $\refl a$ and $\refl b$, in which case~\eqref{eq:path-prod-inverse} yields $\refl{(a,b)}$ and then~\eqref{eq:path-prod} returns us to $(\refl a,\refl b)\jdeq (p,q)\jdeq s$.
\end{proof}

In particular, we have shown that~\eqref{eq:path-prod} has an inverse~\eqref{eq:path-prod-inverse}, which we may denote by
\symlabel{defn:pairpath}
\[
\pairpath : (\id{\proj{1}(x)}{\proj{1}(y)}) \times (\id{\proj{2}(x)}{\proj{2}(y)}) \to (\id x y).
\]
Note that a special case of this yields the propositional uniqueness principle\index{uniqueness!principle, propositional!for product types} for products: $z = (\proj1(z),\proj2(z))$.

It can be helpful to view \pairpath as a \emph{constructor} or \emph{introduction rule} for $\id x y$, analogous to the ``pairing'' constructor of $A\times B$ itself, which introduces the pair $(a,b)$ given $a:A$ and $b:B$.
From this perspective, the two components of~\eqref{eq:path-prod}:
\begin{align*}
  \projpath{1} &: (\id{x}{y}) \to (\id{\proj{1}(x)}{\proj{1} (y)})\\
  \projpath{2} &: (\id{x}{y}) \to (\id{\proj{2}(x)}{\proj{2} (y)})
\end{align*}
are \emph{elimination} rules.
Similarly, the two homotopies which witness~\eqref{eq:path-prod-inverse} as quasi-inverse to~\eqref{eq:path-prod} consist, respectively, of \emph{propositional computation rules}:
\index{computation rule!propositional!for identities between pairs}%
\begin{align*}
  {\projpath{1}{(\pairpath(p, q)})}
  &= 
  {p} \qquad\text{for } p:\id{\proj{1} x}{\proj{1} y} \\
  {\projpath{2}{(\pairpath(p,q)})}
  &= 
  {q} \qquad\text{for } q:\id{\proj{2} x}{\proj{2} y}
\end{align*}
and a \emph{propositional uniqueness principle}:
\index{uniqueness!principle, propositional!for identities between pairs}%
\[
\id{r}{\pairpath(\projpath{1} (r), \projpath{2} (r)) }
\qquad\text{for } r : \id[A \times B] x y.
\]

We can also characterize the reflexivity, inverses, and composition of paths in $A\times B$ componentwise:
\begin{align*}
  {\refl{(z : A \times B)}}
  &= {\pairpath (\refl{\proj{1} z},\refl{\proj{2} z})} \\
  {\opp{p}}
  &= {\pairpath \big(\opp{\projpath{1} (p)},\, \opp{\projpath{2} (p)}\big)} \\
  {{p \ct q}}
  &= {\pairpath \big({\projpath{1} (p)} \ct {\projpath{1} (q)},\,{\projpath{2} (p)} \ct {\projpath{2} (q)}\big)}.
\end{align*}
The same is true for the rest of the higher groupoid structure considered in \autoref{sec:equality}.
All of these equations can be derived by using path induction on the given paths and then returning reflexivity.  

\index{transport!in product types}%
We now consider transport in a pointwise product of type families.
Given type families $ A, B : Z \to \type$, we abusively write $A\times B:Z\to \type$ for the type family defined by $(A\times B)(z) \defeq A(z) \times B(z)$.
Now given $p : \id[Z]{z}{w}$ and $x : A(z) \times B(z)$, we can transport $x$ along $p$ to obtain an element of $A(w)\times B(w)$.

\begin{thm}\label{thm:trans-prod}
  In the above situation, we have
  \[
  \id[A(y) \times B(y)]
  {\transfib{A\times B}px}
  {(\transfib{A}{p}{\proj{1}x}, \transfib{B}{p}{\proj{2}x})}.
  \]
\end{thm}
\begin{proof}
  By path induction, we may assume $p$ is reflexivity, in which case we have
  \begin{align*}
    \transfib{A\times B}px&\jdeq x\\
    \transfib{A}{p}{\proj{1}x}&\jdeq \proj1x\\
    \transfib{B}{p}{\proj{2}x}&\jdeq \proj2x.
  \end{align*}
  Thus, it remains to show $x = (\proj1 x, \proj2x)$.
  But this is the propositional uniqueness principle for product types, which, as we remarked above, follows from \autoref{thm:path-prod}.
\end{proof}

Finally, we consider the functoriality of $\apfunc{}$ under cartesian products.
Suppose given types $A,B,A',B'$ and functions $g:A\to A'$ and $h:B\to B'$; then we can define a function $f:A\times B\to A'\times B'$ by $f(x) \defeq (g(\proj1x),h(\proj2x))$.

\begin{thm}\label{thm:ap-prod}
  In the above situation, given $x,y:A\times B$ and $p:\proj1x=\proj1y$ and $q:\proj2x=\proj2y$, we have
  \[ \id[(f(x)=f(y))]{\ap{f}{\pairpath(p,q)}} {\pairpath(\ap{g}{p},\ap{h}{q})}. \]
\end{thm}
\begin{proof}
  Note first that the above equation is well-typed.
  On the one hand, since $\pairpath(p,q):x=y$ we have $\ap{f}{\pairpath(p,q)}:f(x)=f(y)$.
  On the other hand, since $\proj1(f(x))\jdeq g(\proj1x)$ and $\proj2(f(x))\jdeq h(\proj2x)$, we also have $\pairpath(\ap{g}{p},\ap{h}{q}):f(x)=f(y)$.

  Now, by induction, we may assume $x\jdeq(a,b)$ and $y\jdeq(a',b')$, in which case we have $p:a=a'$ and $q:b=b'$.
  Thus, by path induction, we may assume $p$ and $q$ are reflexivity, in which case the desired equation holds judgmentally.
\end{proof}

\index{type!product|)}%

\section{\texorpdfstring{$\Sigma$}{Σ}-types}
\label{sec:compute-sigma}

\index{type!dependent pair|(}%
Let $A$ be a type and $B:A\to\type$ a type family.
Recall that the $\Sigma$-type, or dependent pair type, $\sm{x:A} B(x)$ is a generalization of the cartesian product type.
Thus, we expect its higher groupoid structure to also be a generalization of the previous section.
In particular, its paths should be pairs of paths, but it takes a little thought to give the correct types of these paths.

Suppose that we have a path $p:w=w'$ in $\sm{x:A}P(x)$.
Then we get $\ap{\proj{1}}{p}:\proj{1}(w)=\proj{1}(w')$.
However, we cannot directly ask whether $\proj{2}(w)$ is identical to $\proj{2}(w')$ since they don't have to be in the same type.
But we can transport\index{transport} $\proj{2}(w)$ along the path $\ap{\proj{1}}{p}$, and this does give us an element of the same type as $\proj{2}(w')$.
By path induction, we do in fact obtain a path $\trans{\ap{\proj{1}}{p}}{\proj{2}(w)}=\proj{2}(w')$.

Recall from the discussion preceding \autoref{lem:mapdep} that
\narrowequation{
  \trans{\ap{\proj{1}}{p}}{\proj{2}(w)}=\proj{2}(w')
}
can be regarded as the type of paths from $\proj2(w)$ to $\proj2(w')$ which lie over the path $\ap{\proj1}{p}$ in $A$.
\index{fibration}%
\index{total!space}%
Thus, we are saying that a path $w=w'$ in the total space determines (and is determined by) a path $p:\proj1(w)=\proj1(w')$ in $A$ together with a path from $\proj2(w)$ to $\proj2(w')$ lying over $p$, which seems sensible.

\begin{rmk}
  Note that if we have $x:A$ and $u,v:P(x)$ such that $(x,u)=(x,v)$, it does not follow that $u=v$.
  All we can conclude is that there exists $p:x=x$ such that $\trans p u = v$.
  This is a well-known source of confusion for newcomers to type theory, but it makes sense from a topological viewpoint: the existence of a path $(x,u)=(x,v)$ in the total space of a fibration between two points that happen to lie in the same fiber does not imply the existence of a path $u=v$ lying entirely \emph{within} that fiber.
\end{rmk}

The next theorem states that we can also reverse this process.
Since it is a direct generalization of \autoref{thm:path-prod}, we will be more concise.

\begin{thm}\label{thm:path-sigma}
Suppose that $P:A\to\type$ is a type family over a type $A$ and let $w,w':\sm{x:A}P(x)$. Then there is an equivalence
\begin{equation*}
\eqvspaced{(w=w')}{\dsm{p:\proj{1}(w)=\proj{1}(w')} \trans{p}{\proj{2}(w)}=\proj{2}(w')}.
\end{equation*}
\end{thm}

\begin{proof}
We define for any $w,w':\sm{x:A}P(x)$, a function
\begin{equation*}
f: (w=w') \to \dsm{p:\proj{1}(w)=\proj{1}(w')} \trans{p}{\proj{2}(w)}=\proj{2}(w')
\end{equation*}
by path induction, with
\begin{equation*}
f(w,w,\refl{w})\defeq(\refl{\proj{1}(w)},\refl{\proj{2}(w)}).
\end{equation*}
We want to show that $f$ is an equivalence.

In the reverse direction, we define
\begin{narrowmultline*}
  g : \prd{w,w':\sm{x:A}P(x)} 
      \Parens{\sm{p:\proj{1}(w)=\proj{1}(w')}\trans{p}{\proj{2}(w)}=\proj{2}(w')}
      \to
      \narrowbreak
      (w=w')
\end{narrowmultline*}
by first inducting on $w$ and $w'$, which splits them into $(w_1,w_2)$ and
$(w_1',w_2')$ respectively, so it suffices to show 
\begin{equation*}
\Parens{\sm{p:w_1 = w_1'}\trans{p}{w_2}=w_2'} \to ((w_1,w_2)=(w_1',w_2')).
\end{equation*}
Next, given a pair $\sm{p:w_1 = w_1'}\trans{p}{w_2}=w_2'$, we can
use $\Sigma$-induction to get $p : w_1 = w_1'$ and $q :
\trans{p}{w_2}=w_2'$.  Inducting on $p$, we have $q :
\trans{\refl{}}{w_2}=w_2'$, and it suffices to show 
$(w_1,w_2)=(w_1,w_2')$.  But $\trans{\refl{}}{w_2} \jdeq w_2$, so
inducting on $q$ reduces to the goal to 
$(w_1,w_2)=(w_1,w_2)$, which we can prove with $\refl{(w_1,w_2)}$.  

Next we show that $f \circ g$ is the identity for all $w$, $w'$ and
$r$, where $r$ has type
\[\dsm{p:\proj{1}(w)=\proj{1}(w')} (\trans{p}{\proj{2}(w)}=\proj{2}(w')).\]
First, we break apart the pairs $w$, $w'$, and $r$ by pair induction, as in the
definition of $g$, and then use two path inductions to reduce both components
of $r$ to \refl{}.  Then it suffices to show that 
$f (g(\refl{},\refl{})) = \refl{}$, which is true by definition.

Similarly, to show that $g \circ f$ is the identity for all $w$, $w'$,
and $p : w = w'$, we can do path induction on $p$, and then induction to
split $w$, at which point it suffices to show that
$g(f (\refl{(w_1,w_2)})) = \refl{(w_1,w_2)}$, which is true by
definition.

Thus, $f$ has a quasi-inverse, and is therefore an equivalence.  
\end{proof}

As we did in the case of cartesian products, we can deduce a propositional uniqueness principle as a special case.

\begin{cor}\label{thm:eta-sigma}
  \index{uniqueness!principle, propositional!for dependent pair types}%
  For $z:\sm{x:A} P(x)$, we have $z = (\proj1(z),\proj2(z))$.
\end{cor}
\begin{proof}
  We have $\refl{\proj1(z)} : \proj1(z) = \proj1(\proj1(z),\proj2(z))$, so by \autoref{thm:path-sigma} it will suffice to exhibit a path $\trans{(\refl{\proj1(z)})}{\proj2(z)} = \proj2(\proj1(z),\proj2(z))$.
  But both sides are judgmentally equal to $\proj2(z)$.
\end{proof}

Like with binary cartesian products, we can think of 
the backward direction of \autoref{thm:path-sigma} as
an introduction form (\pairpath{}{}), the forward direction as
elimination forms (\projpath{1} and \projpath{2}), and the equivalence
as giving a propositional computation rule and uniqueness principle for these.

Note that the lifted path $\mathsf{lift}(u,p)$  of $p:x=y$ at $u:P(x)$ defined in \autoref{thm:path-lifting} may be identified with the special case of the introduction form
\[\pairpath(p,\refl{\trans p u}):(x,u) = (y,\trans p u).\]
\index{transport!in dependent pair types}%
This appears in the statement of action of transport on $\Sigma$-types, which is also a generalization of the action for binary cartesian products:

\begin{thm}\label{transport-Sigma}
  Suppose we have type families
  \begin{equation*}
    P:A\to\type
    \qquad\text{and}\qquad
    Q:\Parens{\sm{x:A} P(x)}\to\type.
  \end{equation*}
  Then we can construct the type family over $A$ defined by
  \begin{equation*}
    x \mapsto \sm{u:P(x)} Q(x,u).
  \end{equation*}
  For any path $p:x=y$ and any $(u,z):\sm{u:P(x)} Q(x,u)$ we have
  \begin{equation*}
    \trans{p}{u,z}=\big(\trans{p}{u},\,\trans{\pairpath(p,\refl{\trans pu})}{z}\big).
  \end{equation*}
\end{thm}

\begin{proof}
Immediate by path induction.
\end{proof}

We leave it to the reader to state and prove a generalization of
\autoref{thm:ap-prod} (see \autoref{ex:ap-sigma}), and to characterize
the reflexivity, inverses, and composition of $\Sigma$-types
componentwise.

\index{type!dependent pair|)}%

\section{The unit type}
\label{sec:compute-unit}

\index{type!unit|(}%
Trivial cases are sometimes important, so we mention briefly the case of the unit type~\unit.

\begin{thm}\label{thm:path-unit}
  For any $x,y:\unit$, we have $\eqv{(x=y)}{\unit}$.
\end{thm}
\begin{proof}
  A function $(x=y)\to\unit$ is easy to define by sending everything to \ttt.
  Conversely, for any $x,y:\unit$ we may assume by induction that $x\jdeq \ttt\jdeq y$.
  In this case we have $\refl{\ttt}:x=y$, yielding a constant function $\unit\to(x=y)$.

  To show that these are inverses, consider first an element $u:\unit$.
  We may assume that $u\jdeq\ttt$, but this is also the result of the composite $\unit \to (x=y)\to\unit$.

  On the other hand, suppose given $p:x=y$.
  By path induction, we may assume $x\jdeq y$ and $p$ is $\refl x$.
  We may then assume that $x$ is \ttt, in which case the composite $(x=y) \to \unit\to(x=y)$ takes $p$ to $\refl x$, i.e.\ to~$p$.
\end{proof}

In particular, any two elements of $\unit$ are equal.
We leave it to the reader to formulate this equivalence in terms of introduction, elimination, computation, and uniqueness rules.
\index{transport!in unit type}%
The transport lemma for \unit is simply the transport lemma for constant type families (\autoref{thm:trans-trivial}).

\index{type!unit|)}%

\section{\texorpdfstring{$\Pi$}{Π}-types and the function extensionality axiom}
\label{sec:compute-pi}

\index{type!dependent function|(}%
\index{type!function|(}%
\index{homotopy|(}%
Given a type $A$ and a type family $B : A \to \type$, consider the dependent function type $\prd{x:A}B(x)$.
We expect the type $f=g$ of paths from $f$ to $g$ in $\prd{x:A} B(x)$ to be equivalent to 
the type of pointwise paths:\index{pointwise!equality of functions}
\begin{equation}
  \eqvspaced{(\id{f}{g})}{\Parens{\prd{x:A} (\id[B(x)]{f(x)}{g(x)})}}.\label{eq:path-forall}
\end{equation}
From a traditional perspective, this would say that two functions which are equal at each point are equal as functions.
\index{continuity of functions in type theory@``continuity'' of functions in type theory}%
From a topological perspective, it would say that a path in a function space is the same as a continuous homotopy.
\index{functoriality of functions in type theory@``functoriality'' of functions in type theory}%
And from a categorical perspective, it would say that an isomorphism in a functor category is a natural family of isomorphisms.

Unlike the case in the previous sections, however, the basic type theory presented in \autoref{cha:typetheory} is insufficient to prove~\eqref{eq:path-forall}.
All we can say is that there is a certain function
\begin{equation}\label{eq:happly}
  \happly : (\id{f}{g}) \to \prd{x:A} (\id[B(x)]{f(x)}{g(x)})
\end{equation}
which is easily defined by path induction.
For the moment, therefore, we will assume:

\begin{axiom}[Function extensionality]\label{axiom:funext}
  \indexsee{axiom!function extensionality}{function extensionality}%
  \indexdef{function extensionality}%
  For any $A$, $B$, $f$, and $g$, the function~\eqref{eq:happly} is an equivalence.
\end{axiom}

We will see in later chapters that this axiom follows both from univalence (see \autoref{sec:compute-universe,sec:univalence-implies-funext}) and from an interval type (see \autoref{sec:interval}).

In particular, \autoref{axiom:funext} implies that~\eqref{eq:happly} has a quasi-inverse
\[
\funext : \Parens{\prd{x:A} (\id{f(x)}{g(x)})} \to {(\id{f}{g})}.
\]
This function is also referred to as ``function extensionality''.
As we did with $\pairpath$ in \autoref{sec:compute-cartprod}, we can regard $\funext$ as an \emph{introduction rule} for the type $\id f g$.
From this point of view, $\happly$ is the \emph{elimination rule}, while the homotopies witnessing $\funext$ as quasi-inverse to $\happly$ become a propositional computation rule\index{computation rule!propositional!for identities between functions}
\[
\id{\happly({\funext{(h)}},x)}{h(x)} \qquad\text{for }h:\prd{x:A} (\id{f(x)}{g(x)})
\]
and a propositional uniqueness principle\index{uniqueness!principle!for identities between functions}:
\[
\id{p}{\funext (x \mapsto \happly(p,{x}))} \qquad\text{for } p: (\id f g).
\]

We can also compute the identity, inverses, and composition in $\Pi$-types; they are simply given by pointwise operations:\index{pointwise!operations on functions}.
\begin{align*}
\refl{f} &= \funext(x \mapsto \refl{f(x)}) \\
\opp{\alpha} &= \funext (x \mapsto \opp{\happly (\alpha,x)})  \\
{\alpha} \ct \beta &= \funext (x \mapsto {\happly({\alpha},x) \ct \happly({\beta},x)}).
\end{align*}
The first of these equalities follows from the definition of $\happly$, while the second and third are easy path inductions.

Since the non-dependent function type $A\to B$ is a special case of the dependent function type $\prd{x:A} B(x)$ when $B$ is independent of $x$, everything we have said above applies in non-dependent cases as well.
\index{transport!in function types}%
The rules for transport, however, are somewhat simpler in the non-dependent case.
Given a type $X$, a path $p:\id[X]{x_1}{x_2}$, type families $A,B:X\to \type$, and a function $f : A(x_1) \to B(x_1)$,  we have
\begin{align}\label{eq:transport-arrow}
  \transfib{A\to B}{p}{f} &=
  \Big(x \mapsto \transfib{B}{p}{f(\transfib{A}{\opp p}{x})}\Big)
\end{align}
where $A\to B$ denotes abusively the type family $X\to \type$ defined by
\[(A\to B)(x) \defeq (A(x)\to B(x)).\]
In other words, when we transport a function $f:A(x_1)\to B(x_1)$ along a path $p:x_1=x_2$, we obtain the function $A(x_2)\to B(x_2)$ which transports its argument backwards along $p$ (in the type family $A$), applies $f$, and then transports the result forwards along $p$ (in the type family $B$).
This can be proven easily by path induction.

\index{transport!in dependent function types}%
Transporting dependent functions is similar, but more complicated.
Suppose given $X$ and $p$ as before, type families $A:X\to \type$ and $B:\prd{x:X} (A(x)\to\type)$, and also a dependent function $f : \prd{a:A(x_1)} B(x_1,a)$.
Then for $p:\id[A]{x_1}{x_2}$ and $a:A(x_2)$, we have
\begin{narrowmultline*}
  \transfib{\Pi_A(B)}{p}{f}(a) = \narrowbreak
  \Transfib{\widehat{B}}{\opp{(\pairpath(\opp{p},\refl{ \trans{\opp p}{a} }))}}{f(\transfib{A}{\opp p}{a})}
\end{narrowmultline*}
where $\Pi_A(B)$ and $\widehat{B}$ denote respectively the type families
\begin{equation}\label{eq:transport-arrow-families}
\begin{array}{rclcl}
\Pi_A(B) &\defeq& \big(x\mapsto \prd{a:A(x)} B(x,a) \big) &:& X\to \type\\
\widehat{B} &\defeq& \big(w \mapsto B(\proj1w,\proj2w) \big) &:& \big(\sm{x:X} A(x)\big) \to \type.
\end{array}
\end{equation}
If these formulas look a bit intimidating, don't worry about the details.
The basic idea is just the same as for the non-dependent function type: we transport the argument backwards, apply the function, and then transport the result forwards again.

Now recall that for a general type family $P:X\to\type$, in \autoref{sec:functors} we defined the type of \emph{dependent paths} over $p:\id[X]xy$ from $u:P(x)$ to $v:P(y)$ to be $\id[P(y)]{\trans{p}{u}}{v}$.
When $P$ is a family of function types, there is an equivalent way to represent this which is often more convenient.
\index{path!dependent!in function types}

\begin{lem}\label{thm:dpath-arrow}
  Given type families $A,B:X\to\type$ and $p:\id[X]xy$, and also $f:A(x)\to B(x)$ and $g:A(y)\to B(y)$, we have an equivalence
  \[ \eqvspaced{ \big(\trans{p}{f} = {g}\big) } { \prd{a:A(x)}  (\trans{p}{f(a)} = g(\trans{p}{a})) }. \]
  Moreover, if $q:\trans{p}{f} = {g}$ corresponds under this equivalence to $\widehat q$, then for $a:A(x)$, the path
  \[ \happly(q,\trans p a) : (\trans p f)(\trans p a) = g(\trans p a)\]
  is equal to the composite
  \begin{align*}
    (\trans p f)(\trans p a)
    &= \trans p {f (\trans {\opp p}{\trans p a})}
    \tag{by~\eqref{eq:transport-arrow}}\\
    &= \trans p {f(a)}\\
    &= g(\trans p a).
    \tag{by $\widehat{q}$}
  \end{align*}
\end{lem}
\begin{proof}
  By path induction, we may assume $p$ is reflexivity, in which case the desired equivalence reduces to function extensionality.
  The second statement then follows by the computation rule for function extensionality.
\end{proof}

As usual, the case of dependent functions is similar, but more complicated.
\index{path!dependent!in dependent function types}

\begin{lem}\label{thm:dpath-forall}
  Given type families $A:X\to\type$ and $B:\prd{x:X} A(x)\to\type$ and $p:\id[X]xy$, and also $f:\prd{a:A(x)} B(x,a)$ and $g:\prd{a:A(y)} B(y,a)$, we have an equivalence
  \[ \eqvspaced{ \big(\trans{p}{f} = {g}\big) } { \Parens{\prd{a:A(x)}  \transfib{\widehat{B}}{\pairpath(p,\refl{\trans pa})}{f(a)} = g(\trans{p}{a}) } } \]
  with $\widehat{B}$ as in~\eqref{eq:transport-arrow-families}.
\end{lem}

We leave it to the reader to prove this and to formulate a suitable computation rule.

\index{homotopy|)}%
\index{type!dependent function|)}%
\index{type!function|)}%

\section{Universes and the univalence axiom}
\label{sec:compute-universe}

\index{type!universe|(}%
\index{equivalence|(}%
Given two types $A$ and $B$, we may consider them as elements of some universe type \type, and thereby form the identity type $\id[\type]AB$.
As mentioned in the introduction, \emph{univalence} is the identification of $\id[\type]AB$ with the type $(\eqv AB)$ of equivalences from $A$ to $B$, which we described in \autoref{sec:basics-equivalences}.
We perform this identification by way of the following canonical function.

\begin{lem}
  For types $A,B:\type$, there is a certain function,
  \begin{equation}\label{eq:uidtoeqv}
    \idtoeqv : (\id[\type]AB) \to (\eqv A B),
  \end{equation}
  defined in the proof.
\end{lem}
\begin{proof}
  We could construct this directly by induction on equality, but the following description is more convenient.
  \index{identity!function}%
  \index{function!identity}%
  Note that the identity function $\idfunc[\type]:\type\to\type$ may be regarded as a type family indexed by the universe \type; it assigns to each type $X:\type$ the type $X$ itself.
  (When regarded as a fibration, its total space is the type $\sm{A:\type}A$ of ``pointed types''; see also \autoref{sec:object-classification}.)
  Thus, given a path $p:A =_\type B$, we have a transport\index{transport} function $\transf{p}:A \to B$.
  We claim that $\transf{p}$ is an equivalence.
  But by induction, it suffices to assume that $p$ is $\refl A$, in which case $\transf{p} \jdeq \idfunc[A]$, which is an equivalence by \autoref{eg:idequiv}.
  Thus, we can define $\idtoeqv(p)$ to be $\transf{p}$ (together with the above proof that it is an equivalence).
\end{proof}

We would like to say that \idtoeqv is an equivalence.
However, as with $\happly$ for function types, the type theory described in \autoref{cha:typetheory} is insufficient to guarantee this.
Thus, as we did for function extensionality, we formulate this property as an axiom: Voevodsky's \emph{univalence axiom}.

\begin{axiom}[Univalence]\label{axiom:univalence}
  \indexdef{univalence axiom}%
  \indexsee{axiom!univalence}{univalence axiom}%
  For any $A,B:\type$, the function~\eqref{eq:uidtoeqv} is an equivalence,
  \[
\eqv{(\id[\type]{A}{B})}{(\eqv A B)}.
\]
\end{axiom}

Technically, the univalence axiom is a statement about a particular universe type $\UU$.
If a universe $\UU$ satisfies this axiom, we say that it is \define{univalent}.
\indexdef{type!universe!univalent}%
\indexdef{univalent universe}%
Except when otherwise noted (e.g.\ in \autoref{sec:univalence-implies-funext}) we will assume that \emph{all} universes are univalent.

\begin{rmk}
  It is important for the univalence axiom that we defined $\eqv AB$ using a ``good'' version of $\isequiv$ as described in \autoref{sec:basics-equivalences}, rather than (say) as $\sm{f:A\to B} \qinv(f)$.
\end{rmk}

In particular, univalence means that \emph{equivalent types may be identified}.
As we did in previous sections, it is useful to break this equivalence into:
\symlabel{ua}
\begin{itemize}
\item An introduction rule for {(\id[\type]{A}{B})},
  \[
  \ua : ({\eqv A B}) \to (\id[\type]{A}{B}).
  \]
\item The elimination rule, which is $\idtoeqv$,
  \[
  \idtoeqv \jdeq \transfibf{X \mapsto X} : (\id[\type]{A}{B}) \to (\eqv A B).
  \]
\item The propositional computation rule\index{computation rule!propositional!for univalence},
  \[
  \transfib{X \mapsto X}{\ua(f)}{x} = f(x).
  \]
\item The propositional uniqueness principle: \index{uniqueness!principle, propositional!for univalence}
  for any $p : \id A B$,
  \[
  \id{p}{\ua(\transfibf{X \mapsto X}(p))}.
  \]
\end{itemize}
We can also identify the reflexivity, concatenation, and inverses of equalities in the universe with the corresponding operations on equivalences:
\begin{align*}
  \refl{A} &= \ua(\idfunc[A]) \\
  \ua(f) \ct \ua(g) &= \ua(g\circ f) \\
  \opp{\ua(f)} &= \ua(f^{-1}).
\end{align*}
The first of these follows because $\idfunc[A] = \idtoeqv(\refl{A})$ by definition of \idtoeqv, and \ua is the inverse of \idtoeqv.
For the second, if we define $p \defeq \ua(f)$ and $q\defeq \ua(g)$, then we have
\[ \ua(g\circ f) = \ua(\idtoeqv(q) \circ \idtoeqv(p)) = \ua(\idtoeqv(p\cdot q)) = p\cdot q\]
using \autoref{thm:transport-concat} and the definition of $\idtoeqv$.
The third is similar.

The following observation, which is a special case of \autoref{thm:transport-compose}, is often useful when applying the univalence axiom.

\begin{lem}\label{thm:transport-is-ap}
  For any type family $B:A\to\type$ and $x,y:A$ with a path $p:x=y$ and $u:B(x)$, we have
  \begin{align*}
    \transfib{B}{p}{u} &= \transfib{X\mapsto X}{\apfunc{B}(p)}{u}\\
    &= \idtoeqv(\apfunc{B}(p))(u).
  \end{align*}
\end{lem}

\index{equivalence|)}%
\index{type!universe|)}%

\section{Identity type}
\label{sec:compute-paths}

\index{type!identity|(}%
Just as the type \id[A]{a}{a'} is characterized up to isomorphism, with
a separate ``definition'' for each $A$, there is no simple
characterization of the type \id[{\id[A]{a}{a'}}]{p}{q} of paths between
paths $p,q : \id[A]{a}{a'}$.
However, our other general classes of theorems do extend to identity types, such as the fact that they respect equivalence.

\begin{thm}\label{thm:paths-respects-equiv}
  If $f : A \to B$ is an equivalence, then for all $a,a':A$, so is
  \[\apfunc{f} : (\id[A]{a}{a'}) \to (\id[B]{f(a)}{f(a')}).\]
\end{thm}
\begin{proof}
  Let $\opp f$ be a quasi-inverse of $f$, with homotopies
  \begin{equation*}
    \alpha:\prd{b:B} (f(\opp f(b))=b)
    \qquad\text{and}\qquad
    \beta:\prd{a:A} (\opp f(f(a)) = a).
  \end{equation*}
  The quasi-inverse of $\apfunc{f}$ is, essentially, \apfunc{\opp f}.
  However, the type of \apfunc{\opp f} is
  \[\apfunc{\opp f} : (\id{f(a)}{f(a')}) \to (\id{\opp f(f(a))}{\opp f(f(a'))}).\]
  Thus, in order to obtain an element of $\id[A]{a}{a'}$ we must concatenate with the paths $\opp{\beta(a)}$ and $\beta (a')$ on either side.
  To show that this gives a quasi-inverse of $\apfunc{f}$, on one hand we must show that for any $p:a=a'$ we have
  \[ \opp{\beta(a)} \ct \apfunc{\opp f}(\apfunc{f}(p)) \ct \beta(a') = p. \]
  This follows from the functoriality of $\apfunc{}$ on function composition and the naturality of homotopies, see \autoref{lem:ap-functor}\ref{item:apfunctor-compose} and \autoref{lem:htpy-natural}.
  On the other hand, we must show that for any $q:f(a)=f(a')$ we have
  \[ \apfunc{f}\big( \opp{\beta(a)} \ct \apfunc{\opp f}(q) \ct \beta(a') \big) = q. \]
  This follows in the same way, using also the functoriality of $\apfunc{}$ on path-concatenation and inverses, see \autoref{lem:ap-functor}\ref{item:apfunctor-ct} and~\ref{item:apfunctor-opp}.
\end{proof}

Thus, if for some type $A$ we have a full characterization of $\id[A]{a}{a'}$, the type $\id[{\id[A]{a}{a'}}]{p}{q}$ is determined as well.  
For example:
\begin{itemize}
\item Paths $p = q$, where $p,q : \id[A \times B]{w}{w'}$, are equivalent to pairs of paths
  \[\id[{\id[A]{\proj{1} w}{\proj{1} w'}}]{\projpath{1}{p}}{\projpath{1}{q}}
  \quad\text{and}\quad
  \id[{\id[B]{\proj{2} w}{\proj{2} w'}}]{\projpath{2}{p}}{\projpath{2}{q}}.
  \]
\item Paths $p = q$, where $p,q : \id[\prd{x:A} B(x)]{f}{g}$, are equivalent to homotopies
  \[\prd{x:A} (\id[f(x)=g(x)] {\happly(p)(x)}{\happly(q)(x)}).\]
\end{itemize}

\index{transport!in identity types}%
Next we consider transport in families of paths, i.e.\ transport in $C:A\to\type$ where each $C(x)$ is an identity type.
The simplest case is when $C(x)$ is a type of paths in $A$ itself, perhaps with one endpoint fixed.

\begin{lem}\label{cor:transport-path-prepost}
  For any $A$ and $a:A$, with $p:x_1=x_2$, we have
  \begin{align*}
    \transfib{x \mapsto (\id{a}{x})} {p} {q} &= q \ct p
    & &\text{for $q:a=x_1$,}\\
    \transfib{x \mapsto (\id{x}{a})} {p} {q} &= \opp {p} \ct q 
    & &\text{for $q:x_1=a$,}\\
    \transfib{x \mapsto (\id{x}{x})} {p} {q} &= \opp{p} \ct q \ct p
    & &\text{for $q:x_1=x_1$.}
  \end{align*}
\end{lem}
\begin{proof}
  Path induction on $p$, followed by the unit laws for composition.
\end{proof}

In other words, transporting with ${x \mapsto \id{c}{x}}$ is post-composition, and transporting with ${x \mapsto \id{x}{c}}$ is contravariant pre-composition.
These may be familiar as the functorial actions of the covariant and contravariant hom-functors $\hom(c, {\blank})$ and $\hom({\blank},c)$ in category theory.

Combining \autoref{cor:transport-path-prepost,thm:transport-compose}, we obtain a more general form:

\begin{thm}\label{thm:transport-path}
  For $f,g:A\to B$, with $p : \id[A]{a}{a'}$ and $q : \id[B]{f(a)}{g(a)}$, we have
  \begin{equation*}
    \id[f(a') = g(a')]{\transfib{x \mapsto \id[B]{f(x)}{g(x)}}{p}{q}}
    {\opp{(\apfunc{f}{p})} \ct q \ct \apfunc{g}{p}}.
  \end{equation*}
\end{thm}

Because $\apfunc{(x \mapsto x)}$ is the identity function and $\apfunc{(x \mapsto c)}$ (where $c$ is a constant) is \refl{c}, \autoref{cor:transport-path-prepost} is a special case.
A yet more general version is when $B$ can be a family of types indexed on $A$:

\begin{thm}\label{thm:transport-path2}
  Let $B : A \to \type$ and $f,g : \prd{x:A} B(x)$, with $p : \id[A]{a}{a'}$ and $q : \id[B(a)]{f(a)}{g(a)}$.
  Then we have
  \begin{equation*}
    \transfib{x \mapsto \id[B(x)]{f(x)}{g(x)}}{p}{q} = 
    \opp{(\apfunc{f}{p})} \ct \apdfunc{(\transfibf{A}{p})}(q) \ct \apfunc{g}{p}.
  \end{equation*}
\end{thm}

Finally, as in \autoref{sec:compute-pi}, for families of identity types there is another equivalent characterization of dependent paths.
\index{path!dependent!in identity types}

\begin{thm}\label{thm:dpath-path}
  For $p:\id[A]a{a'}$ with $q:a=a$ and $r:a'=a'$, we have
  \[ \eqvspaced{ \big(\transfib{x\mapsto (x=x)}{p}{q} = r \big) }{ \big( q \ct p = p \ct r \big). } \]
\end{thm}
\begin{proof}
  Path induction on $p$, followed by the fact that composing with the unit equalities $q\ct 1 = q$ and $r = 1\ct r$ is an equivalence.
\end{proof}

There are more general equivalences involving the application of functions, akin to \autoref{thm:transport-path,thm:transport-path2}.

\index{type!identity|)}%

\section{Coproducts}
\label{sec:compute-coprod}

\index{type!coproduct|(}%
\index{encode-decode method|(}%
So far, most of the type formers we have considered have been what are called \emph{negative}.
\index{type!negative}\index{negative!type}%
\index{polarity}%
Intuitively, this means that their elements are determined by their behavior under the elimination rules: a (dependent) pair is determined by its projections, and a (dependent) function is determined by its values.
The identity types of negative types can almost always be characterized straightforwardly, along with all of their higher structure, as we have done in \crefrange{sec:compute-cartprod}{sec:compute-pi}.
The universe is not exactly a negative type, but its identity types behave similarly: we have a straightforward characterization (univalence) and a description of the higher structure.
Identity types themselves, of course, are a special case.

We now consider our first example of a \emph{positive} type former.
\index{type!positive}\index{positive!type}%
Again informally, a positive type is one which is ``presented'' by certain constructors, with the universal property of a presentation\index{presentation!of a positive type by its constructors} being expressed by its elimination rule.
(Categorically speaking, a positive type has a ``mapping out'' universal property, while a negative type has a ``mapping in'' universal property.)
Because computing with presentations is, in general, an uncomputable problem, for positive types we cannot always expect a straightforward characterization of the identity type.
However, in many particular cases, a characterization or partial characterization does exist, and can be obtained by the general method that we introduce with this example.

(Technically, our chosen presentation of cartesian products and $\Sigma$-types is also positive.
However, because these types also admit a negative presentation which differs only slightly, their identity types have a direct characterization that does not require the method to be described here.)

Consider the coproduct type $A+B$, which is ``presented'' by the injections $\inl:A\to A+B$ and $\inr:B\to A+B$.
Intuitively, we expect that $A+B$ contains exact copies of $A$ and $B$ disjointly, so that we should have
\begin{align}
  {(\inl(a_1)=\inl(a_2))}&\eqvsym {(a_1=a_2)} \label{eq:inlinj}\\
  {(\inr(b_1)=\inr(b_2))}&\eqvsym {(b_1=b_2)}\\
  {(\inl(a)= \inr(b))} &\eqvsym {\emptyt}. \label{eq:inlrdj}
\end{align}
We prove this as follows.
Fix an element $a_0:A$; we will characterize the type family
\begin{equation}
  (x\mapsto (\inl(a_0)=x)) : A+B \to \type.\label{eq:sumcodefam}
\end{equation}
A similar argument would characterize the analogous family $x\mapsto (x = \inr(b_0))$ for any $b_0:B$.
Together, these characterizations imply~\eqref{eq:inlinj}--\eqref{eq:inlrdj}.

In order to characterize~\eqref{eq:sumcodefam}, we will define a type family $\code:A+B\to\type$ and show that $\prd{x:A+B} (\eqv{(\inl(a_0)=x)}{\code(x)})$.
Since we want to conclude~\eqref{eq:inlinj} from this, we should have $\code(\inl(a)) = (a_0=a)$, and since we also want to conclude~\eqref{eq:inlrdj}, we should have $\code (\inr(b)) = \emptyt$.
The essential insight is that we can use the recursion principle of $A+B$ to \emph{define} $\code:A+B\to\type$ by these two equations:
\begin{align*}
  \code(\inl(a)) &\defeq (a_0=a),\\
  \code(\inr(b)) &\defeq \emptyt.
\end{align*}
This is a very simple example of a proof technique that is used quite a
bit when doing homotopy theory in homotopy type theory; see
e.g.\ \autoref{sec:pi1-s1-intro,sec:general-encode-decode}.
We can now show:

\begin{thm}\label{thm:path-coprod}
  For all $x:A+B$ we have $\eqv{(\inl(a_0)=x)}{\code(x)}$.
\end{thm}
\begin{proof}
  The key to the following proof is that we do it for all points $x$ together, enabling us to use the elimination principle for the coproduct.
  We first define a function
  \[ \encode : \prd{x:A+B}{p:\inl(a_0)=x} \code(x) \]
  by transporting reflexivity along $p$:
  \[ \encode(x,p) \defeq \transfib{\code}{p}{\refl{a_0}}. \]
  Note that $\refl{a_0} : \code(\inl(a_0))$, since $\code(\inl(a_0))\jdeq (a_0=a_0)$ by definition of \code.
  Next, we define a function
  \[ \decode : \prd{x:A+B}{c:\code(x)} (\inl(a_0)=x). \]
  To define $\decode(x,c)$, we may first use the elimination principle of $A+B$ to divide into cases based on whether $x$ is of the form $\inl(a)$ or the form $\inr(b)$.

  In the first case, where $x\jdeq \inl(a)$, then $\code(x)\jdeq (a_0=a)$, so that $c$ is an identification between $a_0$ and $a$.
  Thus, $\apfunc{\inl}(c):(\inl(a_0)=\inl(a))$ so we can define this to be $\decode(\inl(a),c)$.

  In the second case, where $x\jdeq \inr(b)$, then $\code(x)\jdeq \emptyt$, so that $c$ inhabits the empty type.
  Thus, the elimination rule of $\emptyt$ yields a value for $\decode(\inr(b),c)$.

  This completes the definition of \decode; we now show that $\encode(x,{\blank})$ and $\decode(x,{\blank})$ are quasi-inverses for all $x$.
  On the one hand, suppose given $x:A+B$ and $p:\inl(a_0)=x$; we want to show
  \narrowequation{
    \decode(x,\encode(x,p)) = p.
  }
  But now by (based) path induction, it suffices to consider $x\jdeq\inl(a_0)$ and $p\jdeq \refl{\inl(a_0)}$:
  \begin{align*}
    \decode(x,\encode(x,p))
    &\jdeq \decode(\inl(a_0),\encode(\inl(a_0),\refl{\inl(a_0)}))\\
    &\jdeq \decode(\inl(a_0),\transfib{\code}{\refl{\inl(a_0)}}{\refl{a_0}})\\
    &\jdeq \decode(\inl(a_0),\refl{a_0})\\
    &\jdeq \ap{\inl}{\refl{a_0}}\\
    &\jdeq \refl{\inl(a_0)}\\
    &\jdeq p.
  \end{align*}
  On the other hand, let $x:A+B$ and $c:\code(x)$; we want to show $\encode(x,\decode(x,c))=c$.
  We may again divide into cases based on $x$.
  If $x\jdeq\inl(a)$, then $c:a_0=a$ and $\decode(x,c)\jdeq \apfunc{\inl}(c)$, so that
  \begin{align}
    \encode(x,\decode(x,c))
    &\jdeq \transfib{\code}{\apfunc{\inl}(c)}{\refl{a_0}}
    \notag\\
    &= \transfib{a\mapsto (a_0=a)}{c}{\refl{a_0}}
    \tag{by \autoref{thm:transport-compose}}\\
    &= \refl{a_0} \ct c
    \tag{by \autoref{cor:transport-path-prepost}}\\
    &= c. \notag
  \end{align}
  Finally, if $x\jdeq \inr(b)$, then $c:\emptyt$, so we may conclude anything we wish.
\end{proof}

\noindent
Of course, there is a corresponding theorem if we fix $b_0:B$ instead of $a_0:A$.

In particular, \autoref{thm:path-coprod} implies that for any $a : A$ and $b : B$ there are functions
\[ \encode(a, {\blank}) : (\inl(a_0)=\inl(a)) \to (a_0=a)\]
and
\[ \encode(b, {\blank}) : (\inl(a_0)=\inr(b)) \to \emptyt. \]
The second of these states
``$\inl(a_0)$ is not equal to $\inr(b)$'', i.e.\ the images of \inl and \inr are disjoint. The traditional reading of the first one, where identity types are viewed as propositions, is just injectivity of $\inl$.  The
full homotopical statement of \autoref{thm:path-coprod} gives more information: the types $\inl(a_0)=\inl(a)$ and
$a_0=a$ are actually equivalent, as are $\inr(b_0)=\inr(b)$ and $b_0=b$.

\begin{rmk}\label{rmk:true-neq-false}
In particular, since the two-element type $\bool$ is equivalent to $\unit+\unit$, we have $\bfalse\neq\btrue$.
\end{rmk}

This proof illustrates a general method for describing path spaces, which we will use often.  To characterize a path space, the first step is to define a comparison fibration ``$\code$'' that provides a more explicit description of the paths.  There are several different methods for proving that such a comparison fibration is equivalent to the paths (we show a few different proofs of the same result in \autoref{sec:pi1-s1-intro}).  The one we have used here is called the \define{encode-decode method}:
\indexdef{encode-decode method}
the key idea is to define $\decode$ generally for all instances of the fibration (i.e.\ as a function $\prd{x:A+B} \code(x) \to (\inl(a_0)=x)$), so that path induction can be used to analyze $\decode(x,\encode(x,p))$.  

\index{transport!in coproduct types}%
As usual, we can also characterize the action of transport in coproduct types.
Given a type~$X$, a path $p:\id[X]{x_1}{x_2}$, and type families $A,B:X\to\type$, we have
\begin{align*}
  \transfib{A+B}{p}{\inl(a)} &= \inl (\transfib{A}{p}{a}),\\
  \transfib{A+B}{p}{\inr(b)} &= \inr (\transfib{B}{p}{b}),
\end{align*}
where as usual, $A+B$ in the superscript denotes abusively the type family $x\mapsto A(x)+B(x)$.
The proof is an easy path induction.

\index{encode-decode method|)}%
\index{type!coproduct|)}%

\section{Natural numbers}
\label{sec:compute-nat}

\index{natural numbers|(}%
\index{encode-decode method|(}%
We use the encode-decode method to characterize the path space of the the natural numbers, which are also a positive type.
In this case, rather than fixing one endpoint, we characterize the two-sided path space all at once.
Thus, the codes for identities are a type family
\[\code:\N\to\N\to\type,\]
defined by double recursion over \N as follows:
\begin{align*}
  \code(0,0) &\defeq \unit\\
  \code(\suc(m),0) &\defeq \emptyt\\
  \code(0,\suc(n)) &\defeq \emptyt\\
  \code(\suc(m),\suc(n)) &\defeq \code(m,n).
\end{align*}
We also define by recursion a dependent function $r:\prd{n:\N} \code(n,n)$, with
\begin{align*}
  r(0) &\defeq \ttt\\
  r(\suc(n)) &\defeq r(n).
\end{align*}

\begin{thm}\label{thm:path-nat}
  For all $m,n:\N$ we have $\eqv{(m=n)}{\code(m,n)}$.
\end{thm}
\begin{proof}
  We define
  \[ \encode : \prd{m,n:\N} (m=n) \to \code(m,n) \]
  by transporting, $\encode(m,n,p) \defeq \transfib{\code(m,{\blank})}{p}{r(m)}$.
  And we define
  \[ \decode : \prd{m,n:\N} \code(m,n) \to (m=n) \]
  by double induction on $m,n$.
  When $m$ and $n$ are both $0$, we need a function $\unit \to (0=0)$, which we define to send everything to $\refl{0}$.
  When $m$ is a successor and $n$ is $0$ or vice versa, the domain $\code(m,n)$ is \emptyt, so the eliminator for \emptyt suffices.
  And when both are successors, we can define $\decode(\suc(m),\suc(n))$ to be the composite
  \begin{narrowmultline*}
    \code(\suc(m),\suc(n))\jdeq\code(m,n)
    \xrightarrow{\decode(m,n)} \narrowbreak
    (m=n)
    \xrightarrow{\apfunc{\suc}}
    (\suc(m)=\suc(n)).
  \end{narrowmultline*}
  Next we show that $\encode(m,n)$ and $\decode(m,n)$ are quasi-inverses for all $m,n$.

  On one hand, if we start with $p:m=n$, then by induction on $p$ it suffices to show
  \[\decode(n,n,\encode(n,n,\refl{n}))=\refl{n}.\]
  But $\encode(n,n,\refl{n}) \jdeq r(n)$, so it suffices to show that $\decode(n,n,r(n)) =\refl{n}$.
  We can prove this by induction on $n$.
  If $n\jdeq 0$, then $\decode(0,0,r(0)) =\refl{0}$ by definition of \decode.
  And in the case of a successor, by the inductive hypothesis we have $\decode(n,n,r(n)) = \refl{n}$, so it suffices to observe that $\apfunc{\suc}(\refl{n}) \jdeq \refl{\suc(n)}$.

  On the other hand, if we start with $c:\code(m,n)$, then we proceed by double induction on $m$ and $n$.
  If both are $0$, then $\decode(0,0,c) \jdeq \refl{0}$, while $\encode(0,0,\refl{0})\jdeq r(0) \jdeq \ttt$.
  Thus, it suffices to recall from \autoref{sec:compute-unit} that every inhabitant of $\unit$ is equal to \ttt.
  If $m$ is $0$ but $n$ is a successor, or vice versa, then $c:\emptyt$, so we are done.
  And in the case of two successors, we have
  \begin{multline*}
    \encode(\suc(m),\suc(n),\decode(\suc(m),\suc(n),c))\\
    \begin{aligned}
    &= \encode(\suc(m),\suc(n),\apfunc{\suc}(\decode(m,n,c)))\\
    &= \transfib{\code(\suc(m),{\blank})}{\apfunc{\suc}(\decode(m,n,c))}{r(\suc(m))}\\
    &= \transfib{\code(\suc(m),\suc({\blank}))}{\decode(m,n,c)}{r(\suc(m))}\\
    &= \transfib{\code(m,{\blank})}{\decode(m,n,c)}{r(m)}\\
    &= \encode(m,n,\decode(m,n,c))\\
    &= c
  \end{aligned}
  \end{multline*}
  using the inductive hypothesis.
\end{proof}

In particular, we have
\begin{equation}\label{eq:zero-not-succ}
  \encode(\suc(m),0) : (\suc(m)=0) \to \emptyt
\end{equation}
which shows that ``$0$ is not the successor of any natural number''.
We also have the composite
\begin{narrowmultline}\label{eq:suc-injective}
  (\suc(m)=\suc(n))
  \xrightarrow{\encode} \narrowbreak
  \code(\suc(m),\suc(n))
  \jdeq \code(m,n) \xrightarrow{\decode} (m=n)
\end{narrowmultline}
which shows that the function $\suc$ is injective.
\index{successor}%

We will study more general positive types in \autoref{cha:induction,cha:hits}.
In \autoref{cha:homotopy}, we will see that the same technique used here to characterize the identity types of coproducts and \nat can also be used to calculate homotopy groups of spheres.

\index{encode-decode method|)}%
\index{natural numbers|)}%

\section{Example: equality of structures}
\label{sec:equality-of-structures}

We now consider one example to illustrate the interaction between the groupoid structure on a type and the type
formers.  In the introduction we remarked that one of the
advantages of univalence is that two isomorphic things are interchangeable,
in the sense that every property or construction involving one also
applies to the other.  Common ``abuses of notation''\index{abuse!of notation} become formally
true.  Univalence itself says that equivalent types are equal, and
therefore interchangeable, which includes e.g.\  the common practice of identifying isomorphic sets.  Moreover, when we define other
mathematical objects as sets, or even general types, equipped with structure or properties, we
can derive the correct notion of equality for them from univalence.  We will illustrate this
point with a significant example in \cref{cha:category-theory}, where we
define the basic notions of category theory in such a way that equality
of categories is equivalence, equality of functors is natural
isomorphism, etc. See in particular \autoref{sec:sip}.
 In this section, we describe a very simple example, coming from algebra.

For simplicity, we use \emph{semigroups} as our example, where a
semigroup is a type equipped with an associative ``multiplication''
operation.  The same ideas apply to other algebraic structures, such as
monoids, groups, and rings.
Recall from \autoref{sec:sigma-types,sec:pat} that the definition of a kind of mathematical structure should be interpreted as defining the type of such structures as a certain iterated $\Sigma$-type.
In the case of semigroups this yields the following.

\begin{defn}
Given a type $A$, the type \semigroupstr{A} of \define{semigroup structures}
\indexdef{semigroup!structure}%
\index{structure!semigroup}%
\index{associativity!of semigroup operation}%
with carrier\index{carrier} $A$ is defined by
\[
\semigroupstr{A} \defeq \sm{m:A \to A \to A} \prd{x,y,z:A} m(x,m(y,z)) = m(m(x,y),z).
\]
A \define{semigroup}
\indexdef{semigroup}%
is a type together with such a structure:
\[
\semigroup \defeq \sm{A:\type} \semigroupstr A
\]
\end{defn}

\noindent 
In the next two sections, we describe two ways in which univalence makes
it easier to work with such semigroups.

\subsection{Lifting equivalences}

\index{lifting!equivalences}%
When working loosely, one might say that a bijection between sets $A$
and $B$ ``obviously'' induces an isomorphism between semigroup
structures on $A$ and semigroup structures on $B$.  With univalence,
this is indeed obvious, because given an equivalence between types $A$
and $B$, we can automatically derive a semigroup structure on $B$ from
one on $A$, and moreover show that this derivation is an equivalence of
semigroup structures.  The reason is that \semigroupstrsym\ is a family
of types, and therefore has an action on paths between types given by
$\mathsf{transport}$:
\[
\transfibf{\semigroupstrsym}{(\ua(e))} : \semigroupstr{A} \to \semigroupstr{B}.
\]
Moreover, this map is an equivalence, because 
$\transfibf{C}(\alpha)$ is always an equivalence with inverse 
$\transfibf{C}{(\opp \alpha)}$, see \cref{thm:transport-concat,thm:omg}.

While the univalence axiom\index{univalence axiom} ensures that this map exists, we need to use
facts about $\mathsf{transport}$ proven in the preceding sections to
calculate what it actually does. Let $(m,a)$ be a semigroup structure on
$A$, and we investigate the induced semigroup structure on $B$ given by
\[
\transfib{\semigroupstrsym}{\ua(e)}{(m,a)}.
\]
First, because
\semigroupstr{X} is defined to be a $\Sigma$-type, by
\cref{transport-Sigma},
\begin{narrowmultline}\label{eq:transport-semigroup-step1}
  \transfib{\semigroupstrsym}{\ua(e)}{(m,a)} = \narrowbreak
  \begin{aligned}[t]
    \big(&\transfib{X \mapsto (X \to X \to X)}{\ua(e)}{m}, \\
     &\transfib{(X,m) \mapsto \mathsf{Assoc}(X,m)}{(\pairpath(\ua(e),\refl{}))}{a}\big)
  \end{aligned}
\end{narrowmultline}
where $\mathsf{Assoc}(X,m)$ is the type $\prd{x,y,z:X} m(x,m(y,z)) = m(m(x,y),z)$.  
That is, the induced semigroup structure consists of an induced
multiplication operation on $B$
\begin{flalign*}
& m' : B \to B \to B \\
& m'(b_1,b_2) \defeq \transfib{X \mapsto (X \to X \to X)}{\ua(e)}{m}(b_1,b_2)
\end{flalign*}
together with an induced proof that $m'$ is associative.  By function
extensionality, it suffices to investigate the behavior of $m'$ when
applied to arguments $b_1,b_2 : B$. By applying
\eqref{eq:transport-arrow} twice, we have that $m'(b_1,b_2)$ is equal to
\begin{narrowmultline*}
  \transfibf{X \mapsto X}\big(
      \ua(e), \narrowbreak
      m(\transfib{X \mapsto X}{\opp{\ua(e)}}{b_1},
        \transfib{X \mapsto X}{\opp{\ua(e)}}{b_2}
       )
   \big).
\end{narrowmultline*}
Then, because $\ua$ is quasi-inverse to $\transfibf{X\mapsto X}$, this is equal to
\[
e(m(\opp{e}(b_1), \opp{e}(b_2))).
\]
Thus, given two elements of $B$, the induced multiplication $m'$ 
sends them to $A$ using the equivalence $e$, multiplies them in $A$, and
then brings the result back to $B$ by $e$, just as one would expect.

Moreover, though we do not show the proof, one can calculate that the
induced proof that $m'$ is associative (the second component of the pair
in \eqref{eq:transport-semigroup-step1}) is equal to a function sending
$b_1,b_2,b_3 : B$ to a path given by the following steps:
\begin{equation}
  \label{eq:transport-semigroup-assoc}
  \begin{aligned}
    m'(m'(b_1,b_2),b_3)
    &= e(m(\opp{e}(m'(b_1,b_2)),\opp{e}(b_3))) \\
    &= e(m(\opp{e}(e(m(\opp{e}(b_1),\opp{e}(b_2)))),\opp{e}(b_3))) \\
    &= e(m(m(\opp{e}(b_1),\opp{e}(b_2)),\opp{e}(b_3))) \\
    &= e(m(\opp{e}(b_1),m(\opp{e}(b_2),\opp{e}(b_3)))) \\
    &= e(m(\opp{e}(b_1),\opp{e}(e(m(\opp{e}(b_2),\opp{e}(b_3)))))) \\
    &= e(m(\opp{e}(b_1),\opp{e}(m'(b_2,b_3)))) \\
    &= m'(b_1,m'(b_2,b_3)).
\end{aligned}
\end{equation}
These steps use the proof $a$ that $m$ is associative and the inverse
laws for $e$.  From an algebra perspective, it may seem strange to
investigate the identity of a proof that an operation is associative,
but this makes sense if we think of $A$ and $B$ as general spaces, with
non-trivial homotopies between paths.  In \cref{cha:logic}, we will
introduce the notion of a \emph{set}, which is a type with only trivial
homotopies, and if we consider semigroup structures on sets, then any
two such associativity proofs are automatically equal.

\subsection{Equality of semigroups}

Using the equations for path spaces discussed in the previous sections,
we can investigate when two semigroups are equal. Given semigroups
$(A,m,a)$ and $(B,m',a')$, by \cref{thm:path-sigma}, the type of paths
\narrowequation{
  (A,m,a) =_\semigroup (B,m',a')
}
is equal to the type of pairs
\begin{align*}
p_1 &: A =_{\type} B \qquad\text{and}\\
p_2 &: \transfib{\semigroupstrsym}{p_1}{(m,a)}{(m',a')}.
\end{align*}
By univalence, $p_1$ is $\ua(e)$ for some equivalence $e$. By
\cref{thm:path-sigma}, function extensionality, and the above analysis of
transport in the type family $\semigroupstrsym$, $p_2$ is equivalent to a pair
of proofs, the first of which shows that
\begin{equation*} \label{eq:equality-semigroup-mult}
\prd{y_1,y_2:B} e(m(\opp{e}(y_1), \opp{e}(y_2))) = m'(y_1,y_2)
\end{equation*}
and the second of which shows that $a'$ is equal to the induced
associativity proof constructed from $a$ in
\eqref{eq:transport-semigroup-assoc}.  But by cancellation of inverses
\eqref{eq:equality-semigroup-mult} is equivalent to
\[
\prd{x_1,x_2:A} e(m(x_1, x_2)) = m'(e(x_2),e(x_2)).
\]
This says that $e$ commutes with the binary operation, in the sense
that it takes multiplication in $A$ (i.e.\ $m$) to multiplication in $B$
(i.e.\ $m'$).  A similar rearrangement is possible for the equation relating
$a$ and $a'$.  Thus, an equality of semigroups consists exactly of an
equivalence on the carrier types that commutes with the semigroup
structure.  

For general types, the proof of associativity is thought of as part of
the structure of a semigroup.  However, if we restrict to set-like types
(again, see \cref{cha:logic}), the
equation relating $a$ and $a'$ is trivially true.  Moreover, in this
case, an equivalence between sets is exactly a bijection.  Thus, we have
arrived at a standard definition of a \emph{semigroup isomorphism}:\index{isomorphism!semigroup} a
bijection on the carrier sets that preserves the multiplication
operation.  It is also possible to use the category-theoretic definition
of isomorphism, by defining a \emph{semigroup homomorphism}\index{homomorphism!semigroup} to be a map
that preserves the multiplication, and arrive at the conclusion that equality of
semigroups is the same as two mutually inverse homomorphisms; but we
will not show the details here; see \autoref{sec:sip}.

The conclusion is that, thanks to univalence, semigroups are equal
precisely when they are isomorphic as algebraic structures. As we will see in \autoref{sec:sip}, the
conclusion applies more generally: in homotopy type theory, all constructions of
mathematical structures automatically respect isomorphisms, without any
tedious proofs or abuse of notation.

\section{Universal properties}
\label{sec:universal-properties}

\index{universal!property|(}%
By combining the path computation rules described in the preceding sections, we can show that various type forming operations satisfy the expected universal properties, interpreted in a homotopical way as equivalences.
For instance, given types $X,A,B$, we have a function
\index{type!product}%
\begin{equation}\label{eq:prod-ump-map}
  (X\to A\times B) \to (X\to A)\times (X\to B)
\end{equation}
defined by $f \mapsto (\proj1 \circ f, \proj2\circ f)$.

\begin{thm}\label{thm:prod-ump}
  \index{universal!property!of cartesian product}%
  \eqref{eq:prod-ump-map} is an equivalence.
\end{thm}
\begin{proof}
  We define the quasi-inverse by sending $(g,h)$ to $\lam{x}(g(x),h(x))$.
  (Technically, we have used the induction principle for the cartesian product $(X\to A)\times (X\to B)$, to reduce to the case of a pair.
  From now on we will often apply this principle without explicit mention.)

  Now given $f:X\to A\times B$, the round-trip composite yields the function
  \begin{equation}
    \lam{x} (\proj1(f(x)),\proj2(f(x))).\label{eq:prod-ump-rt1}
  \end{equation}
  By \autoref{thm:path-prod}, for any $x:X$ we have $(\proj1(f(x)),\proj2(f(x))) = f(x)$.
  Thus, by function extensionality, the function~\eqref{eq:prod-ump-rt1} is equal to $f$.

  On the other hand, given $(g,h)$, the round-trip composite yields the pair $(\lam{x} g(x),\lam{x} h(x))$.
  By the uniqueness principle for functions, this is (judgmentally) equal to $(g,h)$.
\end{proof}

In fact, we also have a dependently typed version of this universal property.
Suppose given a type $X$ and type families $A,B:X\to \type$.
Then we have a function
\begin{equation}\label{eq:prod-umpd-map}
  \Parens{\prd{x:X} (A(x)\times B(x))} \to \Parens{\prd{x:X} A(x)} \times \Parens{\prd{x:X} B(x)}
\end{equation}
defined as before by $f \mapsto (\proj1 \circ f, \proj2\circ f)$.

\begin{thm}\label{thm:prod-umpd}
  \eqref{eq:prod-umpd-map} is an equivalence.
\end{thm}
\begin{proof}
  Left to the reader.
\end{proof}

Just as $\Sigma$-types are a generalization of cartesian products, they satisfy a generalized version of this universal property.
Jumping right to the dependently typed version, suppose we have a type $X$ and type families $A:X\to \type$ and $P:\prd{x:X} A(x)\to\type$.
Then we have a function
\index{type!dependent pair}%
\begin{equation}
  \label{eq:sigma-ump-map}
  \Parens{\prd{x:X}\dsm{a:A(x)} P(x,a)} \to
  \Parens{\sm{g:\prd{x:X} A(x)} \prd{x:X} P(x,g(x))}.
\end{equation}
Note that if we have $P(x,a) \defeq B(x)$ for some $B:X\to\type$, then~\eqref{eq:sigma-ump-map} reduces to~\eqref{eq:prod-umpd-map}.

\begin{thm}\label{thm:ttac}
  \index{universal!property!of dependent pair type}%
  \eqref{eq:sigma-ump-map} is an equivalence.
\end{thm}
\begin{proof}
  As before, we define a quasi-inverse to send $(g,h)$ to the function $\lam{x} (g(x),h(x))$.
  Now given $f:\prd{x:X} \sm{a:A(x)} P(x,a)$, the round-trip composite yields the function
  \begin{equation}
    \lam{x} (\proj1(f(x)),\proj2(f(x))).\label{eq:prod-ump-rt2}
  \end{equation}
  Now for any $x:X$, by \autoref{thm:eta-sigma} (the uniqueness principle for $\Sigma$-types) we have
  \begin{equation*}
    (\proj1(f(x)),\proj2(f(x))) = f(x).
  \end{equation*}
  Thus, by function extensionality,~\eqref{eq:prod-ump-rt2} is equal to $f$.
  On the other hand, given $(g,h)$, the round-trip composite yields $(\lam {x} g(x),\lam{x} h(x))$, which is judgmentally equal to $(g,h)$ as before.
\end{proof}

\index{axiom!of choice!type-theoretic}
This is noteworthy because the propositions-as-types interpretation of~\eqref{eq:sigma-ump-map} is ``the axiom of choice''.
If we read $\Sigma$ as ``there exists'' and $\Pi$ (sometimes) as ``for all'', we can pronounce:
\begin{itemize}
\item $\prd{x:X} \sm{a:A(x)} P(x,a)$ as ``for all $x:X$ there exists an $a:A(x)$ such that $P(x,a)$'', and
\item $\sm{g:\prd{x:X} A(x)} \prd{x:X} P(x,g(x))$ as ``there exists a choice function $g:\prd{x:X} A(x)$ such that for all $x:X$ we have $P(x,g(x))$''.
\end{itemize}
Thus, \autoref{thm:ttac} says that not only is the axiom of choice ``true'', its antecedent is actually equivalent to its conclusion.
(On the other hand, the classical\index{mathematics!classical} mathematician may find that~\eqref{eq:sigma-ump-map} does not carry the usual meaning of the axiom of choice, since we have already specified the values of $g$, and there are no choices left to be made.
We will return to this point in \autoref{sec:axiom-choice}.)

The above universal property for pair types is for ``mapping in'', which is familiar from the category-theoretic notion of products.
However, pair types also have a universal property for ``mapping out'', which may look less familiar.
In the case of cartesian products, the non-dependent version simply expresses
the cartesian closure adjunction\index{adjoint!functor}:
\[ \eqvspaced{\big((A\times B) \to C\big)}{\big(A\to (B\to C)\big)}.\]
The dependent version of this is formulated for a type family $C:A\times B\to \type$:
\[ \eqvspaced{\Parens{\prd{w:A\times B} C(w)}}{\Parens{\prd{x:A}{y:B} C(x,y)}}. \]
Here the left-to-right function is simply the induction principle for $A\times B$, while the right-to-left is evaluation at a pair.
We leave it to the reader to prove that these are quasi-inverses.
There is also a version for $\Sigma$-types:
\begin{equation}
  \eqvspaced{\Parens{\prd{w:\sm{x:A} B(x)} C(w)}}{\Parens{\prd{x:A}{y:B(x)} C(x,y)}}.\label{eq:sigma-lump}
\end{equation}
Again, the left-to-right function is the induction principle.

Some other induction principles are also part of universal properties of this sort.
For instance, path induction is the right-to-left direction of an equivalence as follows:
\index{type!identity}%
\index{universal!property!of identity type}%
\begin{equation}
  \label{eq:path-lump}
  \eqvspaced{\Parens{\prd{x:A}{p:a=x} B(x,p)}}{B(a,\refl a)}
\end{equation}
for any $a:A$ and type family $B:\prd{x:A} (a=x) \to\type$.
However, inductive types with recursion, such as the natural numbers, have more complicated universal properties; see \autoref{cha:induction}.

\index{type!limit}%
\index{type!colimit}%
\index{limit!of types}%
\index{colimit!of types}%
Since \autoref{thm:prod-ump} expresses the usual universal property of a cartesian product (in an appropriate homotopy-theoretic sense), the categorically inclined reader may well wonder about other limits and colimits of types.
In \autoref{ex:coprod-ump} we ask the reader to show that the coproduct type $A+B$ also has the expected universal property, and the nullary cases of $\unit$ (the terminal object) and $\emptyt$ (the initial object) are easy.
\index{type!empty}%
\index{type!unit}%
\indexsee{initial!type}{type, empty}%
\indexsee{terminal!type}{type, unit}%

\indexdef{pullback}%
For pullbacks, the expected explicit construction works: given $f:A\to C$ and $g:B\to C$, we define
\begin{equation}
  A\times_C B \defeq \sm{a:A}{b:B} (f(a)=g(b)).\label{eq:defn-pullback}
\end{equation}
In \autoref{ex:pullback} we ask the reader to verify this.
Some more general homotopy limits can be constructed in a similar way, but for colimits we will need a new ingredient; see \autoref{cha:hits}.

\index{universal!property|)}%

\sectionNotes

The definition of identity types, with their induction principle, is due to Martin-L\"of \cite{Martin-Lof-1972}.
\index{intensional type theory}%
\index{extensional!type theory}%
\index{type theory!intensional}%
\index{type theory!extensional}%
\index{reflection rule}%
As mentioned in the notes to \autoref{cha:typetheory}, our identity types are those that belong to \emph{intensional} type theory, by contrast with those of \emph{extensional} type theory which have an additional ``reflection rule'' saying that if $p:x=y$, then in fact $x\jdeq y$.
This reflection rule implies that all the higher groupoid structure collapses (see \autoref{ex:equality-reflection}), so for nontrivial homotopy we must use the intensional version. 
One may argue, however, that homotopy type theory is in another sense more ``extensional'' than traditional extensional type theory, because of the function extensionality and univalence axioms.  

The proofs of symmetry (inversion) and transitivity (concatenation) for equalities are well-known in type theory.
The fact that these make each type into a 1-groupoid (up to homotopy) was exploited in~\cite{hs:gpd-typethy} to give the first ``homotopy'' style semantics for type theory.  

The actual homotopical interpretation, with identity types as path spaces, and type families as fibrations, is due to \cite{AW}, who used the formalism of Quillen model categories.  An interpretation in (strict) $\infty$-groupoids\index{.infinity-groupoid@$\infty$-groupoid} was also given in the thesis \cite{mw:thesis}.
For a construction of \emph{all} the higher operations and coherences of an $\infty$-groupoid in type theory, see~\cite{pll:wkom-type} and~\cite{bg:type-wkom}.

\index{proof!assistant!Coq@\textsc{Coq}}%
Operations such as $\transfib{P}{p}{\blank}$ and $\apfunc{f}$, and one good notion of equivalence, were first studied extensively in type theory by Voevodsky, using the proof assistant \Coq.
Subsequently, many other equivalent definitions of equivalence have been found, which are compared in \autoref{cha:equivalences}.

The ``computational'' interpretation of identity types, transport, and so on described in \autoref{sec:computational} has been emphasized by~\cite{lh:canonicity}.
They also described a ``1-truncated'' type theory (see \autoref{cha:hlevels}) in which these rules are judgmental equalities.
The possibility of extending this to the full untruncated theory is a subject of current research.

\index{function extensionality}%
The naive form of function extensionality which says that ``if two functions are pointwise equal, then they are equal'' is a common axiom in type theory, going all the way back to \cite{PM2}.
Some stronger forms of function extensionality were considered in~\cite{garner:depprod}.
The version we have used, which identifies the identity types of function types up to equivalence, was first studied by Voevodsky, who also proved that it is implied by the naive version (and by univalence; see \autoref{sec:univalence-implies-funext}).

\index{univalence axiom}%
The univalence axiom is also due to Voevodsky.
It was originally motivated by semantic considerations; see~\cite{klv:ssetmodel}.

In the type theory we are using in this book, function extensionality and univalence have to be assumed as axioms, i.e.\ elements asserted to belong to some type but not constructed according to the rules for that type.
While serviceable, this has a few drawbacks.
For instance, type theory is formally better-behaved if we can base it entirely on rules rather than asserting axioms.
It is also sometimes inconvenient that the theorems of \crefrange{sec:compute-cartprod}{sec:compute-nat} are only propositional equalities (paths) or equivalences, since then we must explicitly mention whenever we pass back and forth across them.
One direction of current research in homotopy type theory is to describe a type system in which these rules are \emph{judgmental} equalities, solving both of these problems at once.
So far this has only been done in some simple cases, although preliminary results such as~\cite{lh:canonicity} are promising.
There are also other potential ways to introduce univalence and function extensionality into a type theory, such as having a sufficiently powerful notion of ``higher quotients'' or ``higher inductive-recursive types''.

The simple conclusions in \crefrange{sec:compute-coprod}{sec:compute-nat} such as ``$\inl$ and $\inr$ are injective and disjoint'' are well-known in type theory, and the construction of the function \encode is the usual way to prove them.
The more refined approach we have described, which characterizes the entire identity type of a positive type (up to equivalence), is a more recent development; see e.g.~\cite{ls:pi1s1}.

\index{axiom!of choice!type-theoretic}%
The type-theoretic axiom of choice~\eqref{eq:sigma-ump-map} was noticed in William Howard's original paper~\cite{howard:pat} on the propositions-as-types correspondence, and was studied further by Martin-L\"of with the introduction of his dependent type theory.  It is mentioned as a ``distributivity law'' in Bourbaki's set theory \cite{Bourbaki}.\index{Bourbaki}%

For a more comprehensive (and formalized) discussion of pullbacks and more general homotopy limits in homotopy type theory, see~\cite{AKL13}.
Limits of diagrams over directed graphs are the easiest general sort of limit to formalize; the problem with diagrams over categories (or more generally $(\infty,1)$-categories)
\index{.infinity1-category@$(\infty,1)$-category}%
\indexsee{category!.infinity1-@$(\infty,1)$-}{$(\infty,1)$-category}%
is that in general, infinitely many coherence conditions are involved in the notion of (homotopy coherent) diagram.\index{diagram}
Resolving this problem is an important open question\index{open!problem} in homotopy type theory.

\sectionExercises

\begin{ex}\label{ex:basics:concat}
  Show that the three obvious proofs of \autoref{lem:concat} are pairwise equal.
\end{ex}

\begin{ex}
  Show that the three equalities of proofs constructed in the previous exercise form a commutative triangle.
  In other words, if the three definitions of concatenation are denoted by $(p \mathbin{\ct_1} q)$, $(p\mathbin{\ct_2} q)$, and $(p\mathbin{\ct_3} q)$, then the concatenated equality
  \[(p\mathbin{\ct_1} q) = (p\mathbin{\ct_2} q) = (p\mathbin{\ct_3} q)\]
  is equal to the equality $(p\mathbin{\ct_1} q) = (p\mathbin{\ct_3} q)$.
\end{ex}

\begin{ex}
  Give a fourth, different, proof of \autoref{lem:concat}, and prove that it is equal to the others.
\end{ex}

\begin{ex}\label{ex:npaths}
  Define, by induction on $n$, a general notion of \define{$n$-dimensional path}\index{path!n-@$n$-} in a type $A$, simultaneously with the type of boundaries for such paths.
\end{ex}

\begin{ex}
  Prove that the functions~\eqref{eq:ap-to-apd} and~\eqref{eq:apd-to-ap} are inverse equivalences.
\end{ex}

\begin{ex}\label{ex:equiv-concat}
  Prove that if $p:x=y$, then the function $(p\ct \blank):(y=z) \to (x=z)$ is an equivalence.
\end{ex}

\begin{ex}\label{ex:ap-sigma}
  State and prove a generalization of \autoref{thm:ap-prod} from cartesian products to $\Sigma$-types.
\end{ex}

\begin{ex}
  State and prove an analogue of \autoref{thm:ap-prod} for coproducts.
\end{ex}

\begin{ex}\label{ex:coprod-ump}
  \index{universal!property!of coproduct}%
  Prove that coproducts have the expected universal property,
  \[ \eqv{(A+B \to X)}{(A\to X)\times (B\to X)}. \]
  Can you generalize this to an equivalence involving dependent functions?
\end{ex}

\begin{ex}\label{ex:sigma-assoc}
  Prove that $\Sigma$-types are ``associative'',
  \index{associativity!of Sigma-types@of $\Sigma$-types}%
  in that for any $A:\UU$ and families $B:A\to\UU$ and $C:(\sm{x:A} B(x))\to\UU$, we have
  \[\eqvspaced{\Parens{\sm{x:A}{y:B(x)} C(\pairr{x,y})}}{\Parens{\sm{p:\sm{x:A}B(x)} C(p)}}. \]
\end{ex}

\begin{ex}\label{ex:pullback}
  A (homotopy) \define{commutative square}
  \indexdef{commutative!square}%
  \begin{equation*}
  \vcenter{\xymatrix{
      P\ar[r]^h\ar[d]_k &
      A\ar[d]^f\\
      B\ar[r]_g &
      C
      }}
  \end{equation*}
  consists of functions $f$, $g$, $h$, and $k$ as shown, together with a path $f \circ h= g \circ k$.
  Note that this is exactly an element of the pullback $(P\to A) \times_{P\to C} (P\to B)$ as defined in~\eqref{eq:defn-pullback}.
  A commutative square is called a (homotopy) \define{pullback square}
  \indexdef{pullback}%
  if for any $X$, the induced map
  \[ (X\to P) \to (X\to A) \times_{(X\to C)} (X\to B) \]
  is an equivalence.
  Prove that the pullback $P \defeq A\times_C B$ defined in~\eqref{eq:defn-pullback} is the corner of a pullback square.
\end{ex}

\begin{ex}\label{ex:pullback-pasting}
  Suppose given two commutative squares
  \begin{equation*}
    \vcenter{\xymatrix{
        A\ar[r]\ar[d] &
        C\ar[r]\ar[d] &
        E\ar[d]\\
        B\ar[r] &
        D\ar[r] &
        F
      }}
  \end{equation*}
  and suppose that the right-hand square is a pullback square.
  Prove that the left-hand square is a pullback square if and only if the outer rectangle is a pullback square.
\end{ex}

\begin{ex}\label{ex:eqvboolbool}
  Show that $\eqv{(\eqv\bool\bool)}{\bool}$.
\end{ex}

\begin{ex}\label{ex:equality-reflection}
  Suppose we add to type theory the \emph{equality reflection rule} which says that if there is an element $p:x=y$, then in fact $x\jdeq y$.
  Prove that for any $p:x=x$ we have $p\jdeq \refl{x}$.
  (This implies that every type is a \emph{set} in the sense to be introduced in \autoref{sec:basics-sets}; see \autoref{sec:hedberg}.)
\end{ex}


\chapter{Sets and logic}
\label{cha:logic}

Type theory, formal or informal, is a collection of rules for manipulating types and their elements.
But when writing mathematics informally in natural language, we generally use familiar words, particularly logical connectives such as ``and'' and ``or'', and logical quantifiers such as ``for all'' and ``there exists''.
In contrast to set theory, type theory offers us more than one way to regard these English phrases as operations on types.
This potential ambiguity needs to be resolved, by setting out local or global conventions, by introducing new annotations to informal mathematics, or both.
This requires some getting used to, but is offset by the fact that because type theory permits this finer analysis of logic, we can represent mathematics more faithfully, with fewer ``abuses of language'' than in set-theoretic foundations.
In this chapter we will explain the issues involved, and justify the choices we have made.

\section{Sets and \texorpdfstring{$n$}{n}-types}
\label{sec:basics-sets}

\index{set|(defstyle}%

In order to explain the connection between the logic of type theory and the logic of set theory, it is helpful to have a notion of \emph{set} in type theory.
While types in general behave like spaces or higher groupoids, there is a subclass of them that behave more like the sets in a traditional set-theoretic system.
Categorically, we may consider \emph{discrete} groupoids, which are determined by a set of objects and only identity morphisms as higher morphisms; while topologically, we may consider spaces having the discrete topology.\index{discrete!space}
More generally, we may consider groupoids or spaces that are \emph{equivalent} to ones of this sort; since everything we do in type theory is up to homotopy, we can't expect to tell the difference.

Intuitively, we would expect a type to ``be a set'' in this sense if it has no higher homotopical information: any two parallel paths are equal (up to homotopy), and similarly for parallel higher paths at all dimensions.
Fortunately, because everything in homotopy type theory is automatically functorial/continuous,
\index{continuity of functions in type theory@``continuity'' of functions in type theory}%
\index{functoriality of functions in type theory@``functoriality'' of functions in type theory}%
it turns out to be sufficient to ask this at the bottom level.

\begin{defn}\label{defn:set}
  A type $A$ is a \define{set}
  if for all $x,y:A$ and all $p,q:x=y$, we have $p=q$.
\end{defn}

More precisely, the proposition $\isset(A)$ is defined to be the type
\[ \isset(A) \defeq \prd{x,y:A}{p,q:x=y} (p=q). \]

As mentioned in \autoref{sec:types-vs-sets},
the sets in homotopy type theory are not like the sets in ZF set theory, in that there is no global ``membership predicate'' $\in$.
They are more like the sets used in structural mathematics and in category theory, whose elements are ``abstract points'' to which we give structure with functions and relations.
This is all we need in order to use them as a foundational system for most set-based mathematics; we will see some examples in \autoref{cha:set-math}.

Which types are sets?
In \autoref{cha:hlevels} we will study a more general form of this question in depth, but for now we can observe some easy examples.

\begin{eg}
  The type \unit is a set.
  For by \autoref{thm:path-unit}, for any $x,y:\unit$ the type $(x=y)$ is equivalent to \unit.
  Since any two elements of \unit are equal, this implies that any two elements of $x=y$ are equal.
\end{eg}

\begin{eg}
  The type $\emptyt$ is a set, for given any $x,y:\emptyt$ we may deduce anything we like, by the induction principle of $\emptyt$.
\end{eg}

\begin{eg}\label{thm:nat-set}
  The type \nat of natural numbers is also a set.
  This follows from \autoref{thm:path-nat}, since all equality types $\id[\nat]xy$ are equivalent to either \unit or \emptyt, and any two inhabitants of \unit or \emptyt are equal.
  We will see another proof of this fact in \autoref{cha:hlevels}.
\end{eg}

Most of the type forming operations we have considered so far also preserve sets.

\begin{eg}\label{thm:isset-prod}
  If $A$ and $B$ are sets, then so is $A\times B$.
  For given $x,y:A\times B$ and $p,q:x=y$, by \autoref{thm:path-prod} we have $p= \pairpath(\projpath1(p),\projpath2(p))$ and $q= \pairpath(\projpath1(q),\projpath2(q))$.
  But $\projpath1(p)=\projpath1(q)$ since $A$ is a set, and $\projpath2(p)=\projpath2(q)$ since $B$ is a set; hence $p=q$.

  Similarly, if $A$ is a set and $B:A\to\type$ is such that each $B(x)$ is a set, then $\sm{x:A} B(x)$ is a set.
\end{eg}

\begin{eg}\label{thm:isset-forall}
  If $A$ is \emph{any} type and $B:A\to \type$ is such that each $B(x)$ is a set, then the type $\prd{x:A} B(x)$ is a set.
  For suppose $f,g:\prd{x:A} B(x)$ and $p,q:f=g$.
  By function extensionality, we have
  \begin{equation*}
    p = {\funext (x \mapsto \happly(p,x))}
    \quad\text{and}\quad
    q = {\funext (x \mapsto \happly(q,x))}.
  \end{equation*}
  But for any $x:A$, we have
  \begin{equation*}
   \happly(p,x):f(x)=g(x)
   \qquad\text{and}\qquad
   \happly(q,x):f(x)=g(x), 
  \end{equation*}
  so since $B(x)$ is a set we have $\happly(p,x) = \happly(q,x)$.
  Now using function extensionality again, the dependent functions $(x \mapsto \happly(p,x))$ and $(x \mapsto \happly(q,x))$ are equal, and hence (applying $\apfunc{\funext}$) so are~$p$ and~$q$.
\end{eg}

For more examples, see \autoref{ex:isset-coprod,ex:isset-sigma}.  For a more systematic investigation of the subsystem (category) of all sets in homotopy type theory, see~\autoref{cha:set-math}.

\index{n-type@$n$-type|(}%

Sets are just the first rung on a ladder of what are called \emph{homotopy $n$-types}.
The next rung consists of \emph{$1$-types}, which are analogous to $1$-groupoids in category theory.
The defining property of a set (which we may also call a \emph{$0$-type}) is that it has no non-trivial paths.
Similarly, the defining property of a $1$-type is that it has no non-trivial paths between paths:

\begin{defn}\label{defn:1type}
  A type $A$ is a \define{1-type}
  \indexdef{1-type}%
  if for all $x,y:A$ and $p,q:x=y$ and $r,s:p=q$, we have $r=s$.
\end{defn}

Similarly, we can define $2$-types, $3$-types, and so on.
We will define the general notion of $n$-type inductively in \autoref{cha:hlevels}, and study the relationships between them.

However, for now it is useful to have two facts in mind.
First, the levels are upward-closed: if $A$ is an $n$-type then $A$ is an $(n+1)$-type.
For example:

\begin{lem}\label{thm:isset-is1type}
  If $A$ is a set (that is, $\isset(A)$ is inhabited), then $A$ is a 1-type.
\end{lem}
\begin{proof}
  Suppose $f:\isset(A)$; then for any $x,y:A$ and $p,q:x=y$ we have $f(x,y,p,q):p=q$.
  Fix $x$, $y$, and $p$, and define $g: \prd{q:x=y} (p=q)$ by $g(q) \defeq f (x,y,p,q)$.
  Then for any $r:q=q'$, we have $\apdfunc{g}(r) : \trans{r}{g(q)} = g(q')$.
  By \autoref{cor:transport-path-prepost}, therefore, we have $g(q) \ct r = g(q')$.

  In particular, suppose given $x,y,p,q$ and $r,s:p=q$, as in \autoref{defn:1type}, and define $g$ as above.
  Then $g(p) \ct r = g(q)$ and also $g(p) \ct s = g(q)$, hence by cancellation $r=s$.
\end{proof}

Second, this stratification of types by level is not degenerate, in the
sense that not all types are sets:  

\begin{eg}\label{thm:type-is-not-a-set}
  \index{type!universe}%
  The universe \type is not a set.
  To prove this, it suffices to exhibit a type $A$ and a path $p:A=A$ which is not equal to $\refl A$.
  Take $A=\bool$, and let $f:A\to A$ be defined by $f(\bfalse)\defeq \btrue$ and $f(\btrue)\defeq \bfalse$.
  Then $f(f(x))=x$ for all $x$ (by an easy case analysis), so $f$ is an equivalence.
  Hence, by univalence, $f$ gives rise to a path $p:A=A$.

  If $p$ were equal to $\refl A$, then (again by univalence) $f$ would equal the identity function of $A$.
  But this would imply that $\bfalse=\btrue$, contradicting \autoref{rmk:true-neq-false}.
\end{eg}

In \autoref{cha:hits,cha:homotopy} we will show that for any $n$, there are types which are not $n$-types.

Note that $A$ is a 1-type exactly when for any $x,y:A$, the identity type $\id[A]xy$ is a set.
(Thus, \autoref{thm:isset-is1type} could equivalently be read as saying that the identity types of a set are also sets.)
This will be the basis of the recursive definition of $n$-types we will give in \autoref{cha:hlevels}.

We can also extend this characterization ``downwards'' from sets.
That is, a type $A$ is a set just when for any $x,y:A$, any two elements of $\id[A]xy$ are equal.
Since sets are equivalently 0-types, it is natural to call a type a \emph{$(-1)$-type} if it has this latter property (any two elements of it are equal).
Such types may be regarded as \emph{propositions in a narrow sense}, and their study is just what is usually called ``logic''; it will occupy us for the rest of this chapter.

\index{n-type@$n$-type|)}%
\index{set|)}%

\section{Propositions as types?}
\label{subsec:pat?}

\index{proposition!as types|(}%
\index{logic!constructive vs classical|(}%
\index{anger|(}%
Until now, we have been following the straightforward ``propositions as types'' philosophy described in \autoref{sec:pat}, according to which English phrases such as ``there exists an $x:A$ such that $P(x)$'' are interpreted by corresponding types such as $\sm{x:A} P(x)$, with the proof of a statement being regarded as judging some specific element to inhabit that type.
However, we have also seen some ways in which the ``logic'' resulting from this reading seems unfamiliar to a classical mathematician.
For instance, in \autoref{thm:ttac} we saw that the statement
\index{axiom!of choice!type-theoretic}%
\begin{equation}\label{eq:english-ac}
  \parbox{\textwidth-2cm}{``If for all $x:X$ there exists an $a:A(x)$ such that $P(x,a)$, then there exists a function $g:\prd{x:A} A(x)$ such that for all $x:X$ we have $P(x,g(x))$,''}
\end{equation}
which looks like the classical\index{mathematics!classical} \emph{axiom of choice}, is always true under this reading. This is a noteworthy, and often useful, feature of the propositions-as-types logic, but it also illustrates how significantly it differs from the classical interpretation of logic, under which the axiom of choice is not a logical truth, but an additional ``axiom''.

On the other hand, we can now also show that corresponding statements looking like the classical \emph{law of double negation} and \emph{law of excluded middle} are incompatible with the univalence axiom.
\index{univalence axiom}%

\begin{thm}\label{thm:not-dneg}
  \index{double negation, law of}%
  It is not the case that for all $A:\UU$ we have $\neg(\neg A) \to A$.
\end{thm}
\begin{proof}
  Recall that $\neg A \jdeq (A\to\emptyt)$.
  We also read ``it is not the case that \dots'' as the operator $\neg$.
  Thus, in order to prove this statement, it suffices to assume given some $f:\prd{A:\UU} (\neg\neg A \to A)$ and construct an element of \emptyt.

  The idea of the following proof is to observe that $f$, like any function in type theory, is ``continuous''.
  \index{continuity of functions in type theory@``continuity'' of functions in type theory}%
  \index{functoriality of functions in type theory@``functoriality'' of functions in type theory}%
  By univalence, this implies that $f$ is \emph{natural} with respect to equivalences of types.
  From this, and a fixed-point-free autoequivalence\index{automorphism!fixed-point-free}, we will be able to extract a contradiction.

  Let $e:\eqv\bool\bool$ be the equivalence defined by $e(\btrue)\defeq\bfalse$ and $e(\bfalse)\defeq\btrue$, as in \autoref{thm:type-is-not-a-set}.
  Let $p:\bool=\bool$ be the path corresponding to $e$ by univalence, i.e.\ $p\defeq \ua(e)$.
  Then we have $f(\bool) : \neg\neg\bool \to\bool$ and
  \[\apd f p : \transfib{A\mapsto (\neg\neg A \to A)}{p}{f(\bool)} = f(\bool).\]
  Hence, for any $u:\neg\neg\bool$, we have
  \[\happly(\apd f p,u) : \transfib{A\mapsto (\neg\neg A \to A)}{p}{f(\bool)}(u) = f(\bool)(u).\]

  Now by~\eqref{eq:transport-arrow}, transporting $f(\bool):\neg\neg\bool\to\bool$ along $p$ in the type family ${A\mapsto (\neg\neg A \to A)}$ is equal to the function which transports its argument along $\opp p$ in the type family $A\mapsto \neg\neg A$, applies $f(\bool)$, then transports the result along $p$ in the type family $A\mapsto A$:
  \begin{narrowmultline*}
    \transfib{A\mapsto (\neg\neg A \to A)}{p}{f(\bool)}(u) =
    \narrowbreak
    \transfib{A\mapsto A}{p}{f(\bool) (\transfib{A\mapsto \neg\neg
        A}{\opp{p}}{u})}
  \end{narrowmultline*}
  However, any two points $u,v:\neg\neg\bool$ are equal by function extensionality, since for any $x:\neg\bool$ we have $u(x):\emptyt$ and thus we can derive any conclusion, in particular $u(x)=v(x)$.
  Thus, we have $\transfib{A\mapsto \neg\neg A}{\opp{p}}{u} = u$, and so from $\happly(\apd f p,u)$ we obtain an equality
  \[ \transfib{A\mapsto A}{p}{f(\bool)(u)} = f(\bool)(u).\]
  Finally, as discussed in \autoref{sec:compute-universe}, transporting in the type family $A\mapsto A$ along the path $p\jdeq \ua(e)$ is equivalent to applying the equivalence $e$; thus we have
  \begin{equation}
    e(f(\bool)(u)) = f(\bool)(u).\label{eq:fpaut}
  \end{equation}
  However, we can also prove that
  \begin{equation}
    \prd{x:\bool} \neg(e(x)=x).\label{eq:fpfaut}
  \end{equation}
  This follows from a case analysis on $x$: both cases are immediate from the definition of $e$ and the fact that $\bfalse\neq\btrue$ (\autoref{rmk:true-neq-false}).
  Thus, applying~\eqref{eq:fpfaut} to $f(\bool)(u)$ and~\eqref{eq:fpaut}, we obtain an element of $\emptyt$.
\end{proof}

\begin{rmk}
  \index{choice operator}%
  \indexsee{operator!choice}{choice operator}%
  In particular, this implies that there can be no Hilbert-style ``choice operator'' which selects an element of every nonempty type.
  The point is that no such operator can be \emph{natural}, and under the univalence axiom, all functions acting on types must be natural with respect to equivalences.
\end{rmk}

\begin{rmk}
  It is, however, still the case that $\neg\neg\neg A \to \neg A$ for any $A$; see \autoref{ex:neg-ldn}.
\end{rmk}

\begin{cor}\label{thm:not-lem}
  \index{excluded middle}%
  It is not the case that for all $A:\UU$ we have $A+(\neg A)$.
\end{cor}
\begin{proof}
  Suppose we had $g:\prd{A:\UU} (A+(\neg A))$.
  We will show that then $\prd{A:\UU} (\neg\neg A \to A)$, so that we can apply \autoref{thm:not-dneg}.
  Thus, suppose $A:\UU$ and $u:\neg\neg A$; we want to construct an element of $A$.

  Now $g(A):A+(\neg A)$, so by case analysis, we may assume either $g(A)\jdeq \inl(a)$ for some $a:A$, or $g(A)\jdeq \inr(w)$ for some $w:\neg A$.
  In the first case, we have $a:A$, while in the second case we have $u(w):\emptyt$ and so we can obtain anything we wish (such as $A$).
  Thus, in both cases we have an element of $A$, as desired.
\end{proof}

Thus, if we want to assume the univalence axiom (which, of course, we do) and still leave ourselves the option of classical\index{mathematics!classical} reasoning (which is also desirable), we cannot use the unmodified propositions-as-types principle to interpret \emph{all} informal mathematical statements into type theory, since then the law of excluded middle would be false.
However, neither do we want to discard propositions-as-types entirely, because of its many good properties (such as simplicity, constructivity, and computability).
We now discuss a modification of propositions-as-types which resolves these problems; in \autoref{subsec:when-trunc} we will return to the question of which logic to use when.
\index{anger|)}%
\index{proposition!as types|)}%
\index{logic!constructive vs classical|)}%

\section{Mere propositions}
\label{subsec:hprops}

\index{logic!of mere propositions|(}%
\index{mere proposition|(defstyle}%
\indexsee{proposition!mere}{mere proposition}%
We have seen that the propositions-as-types logic has both good and bad properties.
Both have a common cause: when types are viewed as propositions, they can contain more information than mere truth or falsity, and all ``logical'' constructions on them must respect this additional information.
This suggests that we could obtain a more conventional logic by restricting attention to types that do \emph{not} contain any more information than a truth value, and only regarding these as logical propositions.

Such a type $A$ will be ``true'' if it is inhabited, and ``false'' if its inhabitation yields a contradiction (i.e.\ if $\neg A \jdeq (A\to\emptyt)$ is inhabited).
\index{inhabited type}%
What we want to avoid, in order to obtain a more traditional sort of logic, is treating as logical propositions those types for which giving an element of them gives more information than simply knowing that the type is inhabited.
For instance, if we are given an element of \bool, then we receive more information than the mere fact that \bool contains some element.
Indeed, we receive exactly \emph{one bit}
\index{bit}%
more information: we know \emph{which} element of \bool we were given.
By contrast, if we are given an element of \unit, then we receive no more information than the mere fact that \unit contains an element, since any two elements of \unit are equal to each other.
This suggests the following definition.

\begin{defn}\label{defn:isprop}
  A type $P$ is a \define{mere proposition}
  if for all $x,y:P$ we have $x=y$.
\end{defn}

Note that since we are still doing mathematics \emph{in} type theory, this is a definition \emph{in} type theory, which means it is a type --- or, rather, a type family.
Specifically, for any $P:\type$, the type $\isprop(P)$ is defined to be
\[ \isprop(P) \defeq \prd{x,y:P} (x=y). \]
Thus, to assert that ``$P$ is a mere proposition'' means to exhibit an inhabitant of $\isprop(P)$, which is a dependent function connecting any two elements of $P$ by a path.
The continuity/naturality of this function implies that not only are any two elements of $P$ equal, but $P$ contains no higher homotopy either.

\begin{lem}\label{thm:inhabprop-eqvunit}
  If $P$ is a mere proposition and $x_0:P$, then $\eqv P \unit$.
\end{lem}
\begin{proof}
  Define $f:P\to\unit$ by $f(x)\defeq \ttt$, and $g:\unit\to P$ by $g(u)\defeq x_0$.
  The claim follows from the next lemma, and the observation that \unit is a mere proposition by \autoref{thm:path-unit}.
\end{proof}

\begin{lem}\label{lem:equiv-iff-hprop}
  If $P$ and $Q$ are mere propositions such that $P\to Q$ and $Q\to P$, then $\eqv P Q$.
\end{lem}
\begin{proof}
  Suppose given $f:P\to Q$ and $g:Q\to P$.
  Then for any $x:P$, we have $g(f(x))=x$ since $P$ is a mere proposition.
  Similarly, for any $y:Q$ we have $f(g(y))=y$ since $Q$ is a mere proposition; thus $f$ and $g$ are quasi-inverses.
\end{proof}

That is, as promised in \autoref{sec:pat}, if two mere propositions are logically equivalent, then they are equivalent.

In homotopy theory, a space that is homotopy equivalent to \unit is said to be \emph{contractible}.
Thus, any mere proposition which is inhabited is contractible (see also \autoref{sec:contractibility}).
On the other hand, the uninhabited type \emptyt is also (vacuously) a mere proposition.
In classical\index{mathematics!classical} mathematics, at least, these are the only two possibilities.

Mere propositions are also called \emph{subterminal objects} (if thinking categorically), \emph{subsingletons} (if thinking set-theoretically), or \emph{h-propositions}.
\indexsee{object!subterminal}{mere proposition}%
\indexsee{subterminal object}{mere proposition}%
\indexsee{subsingleton}{mere proposition}%
\indexsee{h-proposition}{mere proposition}
The discussion in \autoref{sec:basics-sets} suggests we should also call them \emph{$(-1)$-types}; we will return to this in \autoref{cha:hlevels}.
The adjective ``mere'' emphasizes that although any type may be regarded as a proposition (which we prove by giving an inhabitant of it), a type that is a mere proposition cannot usefully be regarded as any \emph{more} than a proposition: there is no additional information contained in a witness of its truth.

Note that a type $A$ is a set if and only if for all $x,y:A$, the identity type $\id[A]xy$ is a mere proposition.
On the other hand, by copying and simplifying the proof of \autoref{thm:isset-is1type}, we have:

\begin{lem}\label{thm:prop-set}
  Every mere proposition is a set.
\end{lem}
\begin{proof}
  Suppose $f:\isprop(A)$; thus for all $x,y:A$ we have $f(x,y):x=y$.  Fix $x:A$
  and define $g(y)\defeq f(x,y)$.   Then for any $y,z:A$ and $p:y=z$ we have $\apd
  g p : \trans{p}{g(y)}={g(z)}$.  Hence by \autoref{cor:transport-path-prepost}, we have
  $g(y)\ct p = g(z)$, which is to say that $p=\opp{g(y)}\ct g(z)$.  Thus, for
  any $p,q:x=y$, we have $p = \opp{g(x)}\ct g(y) = q$.
\end{proof}

In particular, this implies:

\begin{lem}\label{thm:isprop-isprop}\label{thm:isprop-isset}
  For any type $A$, the types $\isprop(A)$ and $\isset(A)$ are mere propositions.
\end{lem}
\begin{proof}
  Suppose $f,g:\isprop(A)$.  By function extensionality, to show $f=g$ it
  suffices to show $f(x,y)=g(x,y)$ for any $x,y:A$.  But $f(x,y)$ and $g(x,y)$
  are both paths in $A$, and hence are equal because, by either $f$ or $g$, we
  have that $A$ is a mere proposition, and hence by \autoref{thm:prop-set} is a
  set.  Similarly, suppose $f,g:\isset(A)$, which is to say that for all
  $a,b:A$ and $p,q:a=b$, we have $f(a,b,p,q):p=q$ and $g(a,b,p,q):p=q$.  But by then since $A$ is a set (by
  either $f$ or $g$), and hence a 1-type, it follows that $f(a,b,p,q)=g(a,b,p,q)$; hence $f=g$ by
  function extensionality.
\end{proof}

We have seen one other example so far: condition~\ref{item:be3} in \autoref{sec:basics-equivalences} asserts that for any function $f$, the type $\isequiv (f)$ should be a mere proposition.

\index{logic!of mere propositions|)}%
\index{mere proposition|)}%

\section{Classical vs.\ intuitionistic logic}
\label{sec:intuitionism}

\index{logic!constructive vs classical|(}%
\index{denial|(}%
With the notion of mere proposition in hand, we can now give the proper formulation of the \define{law of excluded middle}
\indexdef{excluded middle}%
\indexsee{axiom!excluded middle}{excluded middle}%
\indexsee{law!of excluded middle}{excluded middle}%
in homotopy type theory:
\begin{equation}
  \label{eq:lem}
  \LEM{}\;\defeq\;
  \prd{A:\UU} \Big(\isprop(A) \to (A + \neg A)\Big).
\end{equation}
Similarly, the \define{law of double negation}
\indexdef{double negation, law of}%
\indexdef{axiom!double negation}%
\indexdef{law!of double negation}%
is
\begin{equation}
  \label{eq:ldn}
  \prd{A:\UU} \Big(\isprop(A) \to (\neg\neg A \to A)\Big).
\end{equation}
The two are also easily seen to be equivalent to each other---see \autoref{ex:lem-ldn}---so from now on we will generally speak only of \LEM{}.

This formulation of \LEM{} avoids the ``paradoxes'' of \autoref{thm:not-dneg,thm:not-lem}, since \bool is not a mere proposition.
In order to distinguish it from the more general propositions-as-types formulation, we rename the latter:
\symlabel{lem-infty}
\begin{equation*}
  \LEM\infty \defeq \prd{A:\UU} (A + \neg A).
\end{equation*}
For emphasis, the proper version~\eqref{eq:lem}
may be denoted $\LEM{-1}$;
see also \autoref{ex:lemnm}.
Although $\LEM{}$
is not a consequence of the basic type theory described in \autoref{cha:typetheory}, it may be consistently assumed as an axiom (unlike its $\infty$-counterpart).
For instance, we will assume it in \autoref{sec:wellorderings}.

However, it can be surprising how far we can get without using \LEM{}.
Quite often, a simple reformulation of a definition or theorem enables us to avoid invoking excluded middle.
While this takes a little getting used to sometimes, it is often worth the hassle, resulting in more elegant and more general proofs.
We discussed some of the benefits of this in the introduction.

For instance, in classical\index{mathematics!classical} mathematics, double negations are frequently used unnecessarily.
A very simple example is the common assumption that a set $A$ is ``nonempty'', which literally means it is \emph{not} the case that $A$ contains \emph{no} elements.
Almost always what is really meant is the positive assertion that $A$ \emph{does} contain at least one element, and by removing the double negation we make the statement less dependent on \LEM{}.
Recall that we say that a type $A$ is \emph{inhabited}
\index{inhabited type}%
when we assert $A$ itself as a proposition (i.e.\ we construct an element of $A$, usually unnamed).
Thus, often when translating a classical proof into constructive logic, we replace the word ``nonempty'' by ``inhabited'' (although sometimes we must replace it instead by ``merely inhabited''; see \autoref{subsec:prop-trunc}).

Similarly, it is not uncommon in classical mathematics to find unnecessary proofs by contradiction.
\index{proof!by contradiction}%
Of course, the classical form of proof by contradiction proceeds by way of the law of double negation: we assume $\neg A$ and derive a contradiction, thereby deducing $\neg \neg A$, and thus by double negation we obtain $A$.
However, often the derivation of a contradiction from $\neg A$ can be rephrased slightly so as to yield a direct proof of $A$, avoiding the need for \LEM{}.

It is also important to note that if the goal is to prove a \emph{negation}\index{negation}, then ``proof by contradiction'' does not involve \LEM{}.
In fact, since $\neg A$ is by definition the type $A\to\emptyt$, by definition to prove $\neg A$ is to prove a contradiction (\emptyt) under the assumption of $A$.
Similarly, the law of double negation does hold for negated propositions: $\neg\neg\neg A \to \neg A$.
With practice, one learns to distinguish more carefully between negated and non-negated propositions and to notice when \LEM{} is being used and when it is not.

Thus, contrary to how it may appear on the surface, doing mathematics ``constructively'' does not usually involve giving up important theorems, but rather finding the best way to state the definitions so as to make the important theorems constructively provable.
That is, we may freely use the \LEM{} when first investigating a subject, but once that subject is better understood, we can hope to refine its definitions and proofs so as to avoid that axiom.
This sort of observation is even more pronounced in \emph{homotopy} type theory, where the powerful tools of univalence and higher inductive types allow us to constructively attack many problems that traditionally would require classical\index{mathematics!classical} reasoning.
We will see several examples of this in \autoref{part:mathematics}.

It is also worth mentioning that even in constructive mathematics, the law of excluded middle can hold for \emph{some} propositions.
The name traditionally given to such propositions is \emph{decidable}.

\begin{defn}\label{defn:decidable-equality}
  \mbox{}
  \begin{enumerate}
  \item A type $A$ is called \define{decidable}
    \indexdef{decidable!type}%
    \indexdef{type!decidable}%
    if $A+\neg A$.
  \item Similarly, a type family $B:A\to \type$ is \define{decidable}
    \indexdef{decidable!type family}%
    \indexdef{type!family of!decidable}%
    if \narrowequation{\prd{a:A} (B(a)+\neg B(a)).}
  \item In particular, $A$ has \define{decidable equality}
    \indexdef{decidable!equality}%
    \indexsee{equality!decidable}{decidable equality}%
    if \narrowequation{\prd{a,b:A} ((a=b) + \neg(a=b)).}
  \end{enumerate}
\end{defn}

Thus, $\LEM{}$ is exactly the statement that all mere propositions are decidable, and hence so are all families of mere propositions.
In particular, $\LEM{}$ implies that all sets (in the sense of \autoref{sec:basics-sets}) have decidable equality.
Having decidable equality in this sense is very strong; see \autoref{thm:hedberg}.

\index{denial|)}%
\index{logic!constructive vs classical|)}%

\section{Subsets and propositional resizing}
\label{subsec:prop-subsets}

\index{mere proposition|(}%

As another example of the usefulness of mere propositions, we discuss subsets (and more generally subtypes).
Suppose $P:A\to\type$ is a type family, with each type $P(x)$ regarded as a proposition.
Then $P$ itself is a \emph{predicate} on $A$, or a \emph{property} of elements of $A$.

In set theory, whenever we have a predicate on $P$ on a set $A$, we may form the subset $\setof{x\in A | P(x)}$.
In type theory, the obvious analogue is the $\Sigma$-type $\sm{x:A} P(x)$.
An inhabitant of $\sm{x:A} P(x)$ is, of course, a pair $(x,p)$ where $x:A$ and $p$ is a proof of $P(x)$.
However, for general $P$, an element $a:A$ might give rise to more than one distinct element of $\sm{x:A} P(x)$, if the proposition $P(a)$ has more than one distinct proof.
This is counter to the usual intuition of a subset.
But if $P$ is a \emph{mere} proposition, then this cannot happen.

\begin{lem}\label{thm:path-subset}
  Suppose $P:A\to\type$ is a type family such that $P(x)$ is a mere proposition for all $x:A$.
  If $u,v:\sm{x:A} P(x)$ are such that $\proj1(u) = \proj1(v)$, then $u=v$.
\end{lem}
\begin{proof}
  Suppose $p:\proj1(u) = \proj1(v)$.
  By \autoref{thm:path-sigma}, to show $u=v$ it suffices to show $\trans{p}{\proj2(u)} = \proj2(v)$.
  But $\trans{p}{\proj2(u)}$ and $\proj2(v)$ are both elements of $P(\proj1(v))$, which is a mere proposition; hence they are equal.
\end{proof}

For instance, recall that in \autoref{sec:basics-equivalences} we defined
\[(\eqv A B) \;\defeq\; \sm{f:A\to B} \isequiv (f),\]
where each type $\isequiv (f)$ was supposed to be a mere proposition.
It follows that if two equivalences have equal underlying functions, then they are equal as equivalences.

\label{defn:setof}%
Henceforth, if $P:A\to \type$ is a family of mere propositions (i.e.\ each $P(x)$ is a mere proposition), we may write
\begin{equation}
  \label{eq:subset}
  \setof{x:A | P(x)}
\end{equation}
as an alternative notation for $\sm{x:A} P(x)$.
(There is no technical reason not to use this notation for arbitrary $P$ as well, but such usage could be confusing due to unintended connotations.)
If $A$ is a set, we call \eqref{eq:subset} a \define{subset}
\indexdef{subset}%
\indexdef{type!subset}%
of $A$; for general $A$ we might call it a \define{subtype}.
\indexdef{subtype}%
We may also refer to $P$ itself as a \emph{subset} or \emph{subtype} of $A$; this is actually more correct, since the type~\eqref{eq:subset} in isolation doesn't remember its relationship to $A$.

\symlabel{membership}
Given such a $P$ and $a:A$, we may write $a\in P$ or $a\in \setof{x:A | P(x)}$ to refer to the mere proposition $P(a)$.
If it holds, we may say that $a$ is a \define{member} of $P$.
\symlabel{subset}
Similarly, if $\setof{x:A | Q(x)}$ is another subset of $A$, then we say that $P$ is \define{contained}
\indexdef{containment!of subsets}%
\indexdef{inclusion!of subsets}%
in $Q$, and write $P\subseteq Q$, if we have $\prd{x:A}(P(x)\rightarrow Q(x))$.

As further examples of subtypes, we may define the ``subuniverses'' of sets and of mere propositions in a universe \UU:
\symlabel{setU}\symlabel{propU}
\begin{align*}
  \setU &\defeq \setof{A:\UU | \isset(A) },\\
  \propU &\defeq \setof{A:\UU | \isprop(A) }.
\end{align*}
An element of $\setU$ is a type $A:\UU$ together with evidence $s:\isset(A)$, and similarly for $\propU$.
\autoref{thm:path-subset} implies that $\id[\setU]{(A,s)}{(B,t)}$ is equivalent to $\id[\UU]AB$ (and hence to $\eqv AB$).
Thus, we will frequently abuse notation and write simply $A:\setU$ instead of $(A,s):\setU$.
We may also drop the subscript \UU if there is no need to specify the universe in question.

Recall that for any two universes $\UU_i$ and $\UU_{i+1}$, if $A:\UU_i$ then also $A:\UU_{i+1}$.
Thus, for any $(A,s):\set_{\UU_i}$ we also have $(A,s):\set_{\UU_{i+1}}$, and similarly for $\prop_{\UU_i}$, giving natural maps
\begin{align}
  \set_{\UU_i} &\to \set_{\UU_{i+1}}\label{eq:set-up},\\
  \prop_{\UU_i} &\to \prop_{\UU_{i+1}}.\label{eq:prop-up}
\end{align}
The map~\eqref{eq:set-up} cannot be an equivalence, since then we could reproduce the paradoxes\index{paradox} of self-reference that are familiar from Cantorian set theory.
However, although~\eqref{eq:prop-up} is not automatically an equivalence in the type theory we have presented so far, it is consistent to suppose that it is.
That is, we may consider adding to type theory the following axiom.

\begin{axiom}[Propositional resizing]
  \indexsee{axiom!propositional resizing}{propositional resizing}%
  \indexdef{propositional!resizing}%
  \indexsee{resizing!propositional}{propositional resizing}%
  The map $\prop_{\UU_i} \to \prop_{\UU_{i+1}}$ is an equivalence.
\end{axiom}

We refer to this axiom as \define{propositional resizing},
since it means that any mere proposition in the universe $\UU_{i+1}$ can be ``resized'' to an equivalent one in the smaller universe $\UU_i$.
It follows automatically if $\UU_{i+1}$ satisfies \LEM{} (see \autoref{ex:lem-impred}).
We will not assume this axiom in general, although in some places we will use it as an explicit hypothesis.
It is a form of \emph{impredicativity} for mere propositions, and by avoiding its use, the type theory is said to remain \emph{predicative}.
\indexsee{impredicativity!for mere propositions}{propositional resizing}%
\indexdef{mathematics!predicative}%

In practice, what we want most frequently is a slightly different statement: that a universe \UU under consideration contains a type which ``classifies all mere propositions''.
In other words, we want a type $\Omega:\UU$ together with an $\Omega$-indexed family of mere propositions, which contains every mere proposition up to equivalence.
This statement follows from propositional resizing as stated above if $\UU$ is not the smallest universe $\UU_0$, since then we can define $\Omega\defeq \prop_{\UU_0}$.

One use for impredicativity is to define power sets.
It is natural to define the \define{power set}\indexdef{power set} of a set $A$ to be $A\to\propU$; but in the absence of impredicativity, this definition depends
(even up to equivalence) on the choice of the universe \UU.
But with propositional resizing, we can define the power set to be
\symlabel{powerset}%
\[ \power A \defeq (A\to\Omega),\]
which is then independent of $\UU$.
See also \autoref{subsec:piw}.

\section{The logic of mere propositions}
\label{subsec:logic-hprop}

\index{logic!of mere propositions|(}%
We mentioned in \autoref{sec:types-vs-sets} that in contrast to type theory, which has only one basic notion (types), set-theoretic foundations have two basic notions: sets and propositions.
Thus, a classical\index{mathematics!classical} mathematician is accustomed to manipulating these two kinds of objects separately.

It is possible to recover a similar dichotomy in type theory, with the role of the set-theoretic propositions being played by the types (and type families) that are \emph{mere} propositions.
In many cases, the logical connectives and quantifiers can be represented in this logic by simply restricting the corresponding type-former to the mere propositions.
Of course, this requires knowing that the type-former in question preserves mere propositions.

\begin{eg}
  If $A$ and $B$ are mere propositions, so is $A\times B$.
  This is easy to show using the characterization of paths in products, just like \autoref{thm:isset-prod} but simpler.
  Thus, the connective ``and'' preserves mere propositions.
\end{eg}

\begin{eg}\label{thm:isprop-forall}
  If $A$ is any type and $B:A\to \type$ is such that for all $x:A$, the type $B(x)$ is a mere proposition, then $\prd{x:A} B(x)$ is a mere proposition.
  The proof is just like \autoref{thm:isset-forall} but simpler: given $f,g:\prd{x:A} B(x)$, for any $x:A$ we have $f(x)=g(x)$ since $B(x)$ is a mere proposition.
  But then by function extensionality, we have $f=g$.

  In particular, if $B$ is a mere proposition, then so is $A\to B$ regardless of what $A$ is.
  In even more particular, since \emptyt is a mere proposition, so is $\neg A \jdeq (A\to\emptyt)$.
  \index{quantifier!universal}%
  Thus, the connectives ``implies'' and ``not'' preserve mere propositions, as does the quantifier ``for all''.
\end{eg}

On the other hand, some type formers do not preserve mere propositions.
Even if $A$ and $B$ are mere propositions, $A+B$ will not in general be.
For instance, \unit is a mere proposition, but $\bool=\unit+\unit$ is not.
Logically speaking, $A+B$ is a ``purely constructive'' sort of ``or'': a witness of it contains the additional information of \emph{which} disjunct is true.
Sometimes this is very useful, but if we want a more classical sort of ``or'' that preserves mere propositions, we need a way to ``truncate'' this type into a mere proposition by forgetting this additional information.

\index{quantifier!existential}%
The same issue arises with the $\Sigma$-type $\sm{x:A} P(x)$.
This is a purely constructive interpretation of ``there exists an $x:A$ such that $P(x)$'' which remembers the witness $x$, and hence is not generally a mere proposition even if each type $P(x)$ is.
(Recall that we observed in \autoref{subsec:prop-subsets} that $\sm{x:A} P(x)$ can also be regarded as ``the subset of those $x:A$ such that $P(x)$''.)

\section{Propositional truncation}
\label{subsec:prop-trunc}

\index{truncation!propositional|(defstyle}%
\indexsee{type!squash}{truncation, propositional}%
\indexsee{squash type}{truncation, propositional}%
\indexsee{bracket type}{truncation, propositional}%
\indexsee{type!bracket}{truncation, propositional}%
The \emph{propositional truncation}, also called the \emph{$(-1)$-truncation}, \emph{bracket type}, or \emph{squash type}, is an additional type former which ``squashes'' or ``truncates'' a type down to a mere proposition, forgetting all information contained in inhabitants of that type other than their existence.

More precisely, for any type $A$, there is a type $\brck{A}$.
It has two constructors:
\begin{itemize}
\item For any $a:A$ we have $\bproj a : \brck A$.
\item For any $x,y:\brck A$, we have $x=y$.
\end{itemize}
The first constructor means that if $A$ is inhabited, so is $\brck A$.
The second ensures that $\brck A$ is a mere proposition; usually we leave the witness of this fact nameless.

\index{recursion principle!for truncation}%
The recursion principle of $\brck A$ says that:
\begin{itemize}
\item If $B$ is a mere proposition and we have $f:A\to B$, then there is an induced $g:\brck A \to B$ such that $g(\bproj a) \jdeq f(a)$ for all $a:A$.
\end{itemize}
In other words, any mere proposition which follows from (the inhabitedness of) $A$ already follows from $\brck A$.
Thus, $\brck A$, as a mere proposition, contains no more information than the inhabitedness of $A$.
(There is also an induction principle for $\brck A$, but it is not especially useful; see \autoref{ex:prop-trunc-ind}.)

In \autoref{ex:lem-brck,ex:impred-brck,sec:hittruncations} we will describe some ways to construct $\brck{A}$ in terms of more general things.
For now, we simply assume it as an additional rule alongside those of \autoref{cha:typetheory}.

With the propositional truncation, we can extend the ``logic of mere propositions'' to cover disjunction and the existential quantifier.
Specifically, $\brck{A+B}$ is a mere propositional version of ``$A$ or $B$'', which does not ``remember'' the information of which disjunct is true.

The recursion principle of truncation implies that we can still do a case analysis on $\brck{A+B}$ \emph{when attempting to prove a mere proposition}.
That is, suppose we have an assumption $u:\brck{A+B}$ and we are trying to prove a mere proposition $Q$.
In other words, we are trying to define an element of $\brck{A+B} \to Q$.
Since $Q$ is a mere proposition, by the recursion principle for propositional truncation, it suffices to construct a function $A+B\to Q$.
But now we can use case analysis on $A+B$.

Similarly, for a type family $P:A\to\type$, we can consider $\brck{\sm{x:A} P(x)}$, which is a mere propositional version of ``there exists an $x:A$ such that $P(x)$''.
As for disjunction, by combining the induction principles of truncation and $\Sigma$-types, if we have an assumption of type $\brck{\sm{x:A} P(x)}$, we may introduce new assumptions $x:A$ and $y:P(x)$ \emph{when attempting to prove a mere proposition}.
In other words, if we know that there exists some $x:A$ such that $P(x)$, but we don't have a particular such $x$ in hand, then we are free to make use of such an $x$ as long as we aren't trying to construct anything which might depend on the particular value of $x$.
Requiring the codomain to be a mere proposition expresses this independence of the result on the witness, since all possible inhabitants of such a type must be equal.

For the purposes of set-level mathematics in \autoref{cha:real-numbers,cha:set-math},
where we deal mostly with sets and mere propositions, it is convenient to use the
traditional logical notations to refer only to ``propositionally truncated logic''.

\begin{defn} \label{defn:logical-notation}
  We define \define{traditional logical notation}
  \indexdef{implication}%
  \indexdef{traditional logical notation}%
  \indexdef{logical notation, traditional}%
  \index{quantifier}%
  \indexsee{existential quantifier}{quantifier, existential}%
  \index{quantifier!existential}%
  \indexsee{universal!quantifier}{quantifier, universal}%
  \index{quantifier!universal}%
  \indexdef{conjunction}%
  \indexdef{disjunction}%
  \indexdef{true}%
  \indexdef{false}%
  using truncation as follows, where $P$ and $Q$ denote mere propositions (or families thereof):
  {\allowdisplaybreaks
  \begin{align*}
    \top            &\ \defeq \ \unit \\
    \bot            &\ \defeq \ \emptyt \\
    P \land Q       &\ \defeq \ P \times Q \\
    P \Rightarrow Q &\ \defeq \ P \to Q \\
    P \Leftrightarrow Q &\ \defeq \ P = Q \\
    \neg P          &\ \defeq \ P \to \emptyt \\
    P \lor Q        &\ \defeq \ \brck{P + Q} \\
    \fall{x : A} P(x) &\ \defeq \ \prd{x : A} P(x) \\
    \exis{x : A} P(x) &\ \defeq \ \Brck{\sm{x : A} P(x)}
  \end{align*}}
\end{defn}

The notations $\land$ and $\lor$ are also used in homotopy theory for the smash product and the wedge of pointed spaces, which we will introduce in \autoref{cha:hits}.
This technically creates a potential for conflict, but no confusion will generally arise.

Similarly, when discussing subsets as in \autoref{subsec:prop-subsets}, we may use the traditional notation for intersections, unions, and complements:
\indexdef{intersection!of subsets}%
\symlabel{intersection}%
\indexdef{union!of subsets}%
\symlabel{union}%
\indexdef{complement, of a subset}%
\symlabel{complement}%
\begin{align*}
  \setof{x:A | P(x)} \cap \setof{x:A | Q(x)}
  &\defeq \setof{x:A | P(x) \land Q(x)},\\
  \setof{x:A | P(x)} \cup \setof{x:A | Q(x)}
  &\defeq \setof{x:A | P(x) \lor Q(x)},\\
  A \setminus \setof{x:A | P(x)}
  &\defeq \setof{x:A | \neg P(x)}.
\end{align*}
Of course, in the absence of \LEM{}, the latter are not ``complements'' in the usual sense: we may not have $B \cup (A\setminus B) = A$.

\index{truncation!propositional|)}%
\index{mere proposition|)}%
\index{logic!of mere propositions|)}%

\section{The axiom of choice}
\label{sec:axiom-choice}

\index{axiom!of choice|(defstyle}%
\index{denial|(}%
We can now properly formulate the axiom of choice in homotopy type theory.
Assume a type $X$ and type families
\begin{equation*}
  A:X\to\type
  \qquad\text{and}\qquad 
  P:\prd{x:X} A(x)\to\type,
\end{equation*}
and moreover that
\begin{itemize}
\item $X$ is a set,
\item $A(x)$ is a set for all $x:X$, and
\item $P(x,a)$ is a mere proposition for all $x:X$ and $a:A(x)$.
\end{itemize}
The \define{axiom of choice}
$\choice{}$ asserts that under these assumptions,
\begin{equation}\label{eq:ac}
  \Parens{\prd{x:X} \Brck{\sm{a:A(x)} P(x,a)}}
  \to
  \Brck{\sm{g:\prd{x:X} A(x)} \prd{x:X} P(x,g(x))}.
\end{equation}
Of course, this is a direct translation of~\eqref{eq:english-ac} where we read ``there exists $x:A$ such that $B(x)$'' as $\brck{\sm{x:A}B(x)}$, so we could have written the statement in the familiar logical notation as
\begin{narrowmultline*}
  \textstyle
  \Big(\fall{x:X}\exis{a:A(x)} P(x,a)\Big)
  \Rightarrow \narrowbreak
  \Big(\exis{g : \prd{x:X} A(x)} \fall{x : X} P(x,g(x))\Big).
\end{narrowmultline*}
In particular, note that the propositional truncation appears twice.
The truncation in the domain means we assume that for every $x$ there exists some $a:A(x)$ such that $P(x,a)$, but that these values are not chosen or specified in any known way.
The truncation in the codomain means we conclude that there exists some function $g$, but this function is not determined or specified in any known way.

In fact, because of \autoref{thm:ttac}, this axiom can also be expressed in a simpler form.

\begin{lem}\label{thm:ac-epis-split}
  The axiom of choice~\eqref{eq:ac} is equivalent to the statement that for any set $X$ and any $Y:X\to\type$ such that each $Y(x)$ is a set, we have
  \begin{equation}
    \Parens{\prd{x:X} \Brck{Y(x)}}
    \to
    \Brck{\prd{x:X} Y(x)}.\label{eq:epis-split}
  \end{equation}
\end{lem}

This corresponds to a well-known equivalent form of the classical\index{mathematics!classical} axiom of choice, namely ``the cartesian product of a family of nonempty sets is nonempty.''

\begin{proof}
  By \autoref{thm:ttac}, the codomain of~\eqref{eq:ac} is equivalent to
  \[\Brck{\prd{x:X} \sm{a:A(x)} P(x,a)}.\]
  Thus,~\eqref{eq:ac} is equivalent to the instance of~\eqref{eq:epis-split} where \narrowequation{Y(x) \defeq \sm{a:A(x)} P(x,a).}
  Conversely,~\eqref{eq:epis-split} is equivalent to the instance of~\eqref{eq:ac} where $A(x)\defeq Y(x)$ and $P(x,a)\defeq\unit$.
  Thus, the two are logically equivalent.
  Since both are mere propositions, by \autoref{lem:equiv-iff-hprop} they are equivalent types.
\end{proof}

As with \LEM{}, the equivalent forms~\eqref{eq:ac} and~\eqref{eq:epis-split} are not a consequence of our basic type theory, but they may consistently be assumed as axioms.

\begin{rmk}
  It is easy to show that the right side of~\eqref{eq:epis-split} always implies the left.
  Since both are mere propositions, by \autoref{lem:equiv-iff-hprop} the axiom of choice is also equivalent to asking for an equivalence
  \[ \eqv{\Parens{\prd{x:X} \Brck{Y(x)}}}{\Brck{\prd{x:X} Y(x)}} \]
  This illustrates a common pitfall: although dependent function types preserve mere propositions (\autoref{thm:isprop-forall}), they do not commute with truncation: $\brck{\prd{x:A} P(x)}$ is not generally equivalent to $\prd{x:A} \brck{P(x)}$.
  The axiom of choice, if we assume it, says that this is true \emph{for sets}; as we will see below, it fails in general.
\end{rmk}

The restriction in the axiom of choice to types that are sets can be relaxed to a certain extent.
For instance, we may allow $A$ and $P$ in~\eqref{eq:ac}, or $Y$ in~\eqref{eq:epis-split}, to be arbitrary type families; this results in a seemingly stronger statement that is equally consistent.
We may also replace the propositional truncation by the more general $n$-truncations to be considered in \autoref{cha:hlevels}, obtaining a spectrum of axioms $\choice n$ interpolating between~\eqref{eq:ac}, which we call simply \choice{} (or $\choice{-1}$ for emphasis), and \autoref{thm:ttac}, which we shall call $\choice\infty$.
See also \autoref{ex:acnm,ex:acconn}.
However, observe that we cannot relax the requirement that $X$ be a set.  

\begin{lem}\label{thm:no-higher-ac}
  There exists a type $X$ and a family $Y:X\to \type$ such that each $Y(x)$ is a set, but such that~\eqref{eq:epis-split} is false.
\end{lem}
\begin{proof}
  Define $X\defeq \sm{A:\type} \brck{\bool = A}$, and let $x_0 \defeq (\bool, \bproj{\refl{\bool}}) : X$.
  Then by the identification of paths in $\Sigma$-types, the fact that $\brck{A=\bool}$ is a mere proposition, and univalence, for any $(A,p),(B,q):X$ we have $\eqv{(\id[X]{(A,p)}{(B,q)})}{(\eqv AB)}$.
  In particular, $\eqv{(\id[X]{x_0}{x_0})}{(\eqv \bool\bool)}$, so as in \autoref{thm:type-is-not-a-set}, $X$ is not a set.

  On the other hand, if $(A,p):X$, then $A$ is a set; this follows by induction on truncation for $p:\brck{\bool=A}$ and the fact that $\bool$ is a set.
  Since $\eqv A B$ is a set whenever $A$ and $B$ are, it follows that $\id[X]{x_1}{x_2}$ is a set for any $x_1,x_2:X$, i.e.\ $X$ is a 1-type.
  In particular, if we define $Y:X\to\UU$ by $Y(x) \defeq (x_0=x)$, then each $Y(x)$ is a set.

  Now by definition, for any $(A,p):X$ we have $\brck{\bool=A}$, and hence $\brck{x_0 = (A,p)}$.
  Thus, we have $\prd{x:X} \brck{Y(x)}$.
  If~\eqref{eq:epis-split} held for this $X$ and $Y$, then we would also have $\brck{\prd{x:X} Y(x)}$.
  Since we are trying to derive a contradiction ($\emptyt$), which is a mere proposition, we may assume $\prd{x:X} Y(x)$, i.e.\ that $\prd{x:X} (x_0=x)$.
  But this implies $X$ is a mere proposition, and hence a set, which is a contradiction.
\end{proof}

\index{denial|)}%
\index{axiom!of choice|)}%

\section{The principle of unique choice}
\label{sec:unique-choice}

\index{unique!choice|(defstyle}%
\indexsee{axiom!of choice!unique}{unique choice}%

The following observation is trivial, but very useful.

\begin{lem}
  If $P$ is a mere proposition, then $\eqv P {\brck P}$.
\end{lem}
\begin{proof}
  Of course, we have $P\to \brck{P}$ by definition.
  And since $P$ is a mere proposition, the universal property of $\brck P$ applied to $\idfunc[P] :P\to P$ yields $\brck P \to P$.
  These functions are quasi-inverses by \autoref{lem:equiv-iff-hprop}.
\end{proof}

Among its important consequences is the following.

\begin{cor}[The principle of unique choice]\label{cor:UC}
  Suppose a type family $P:A\to \type$ such that
  \begin{enumerate}
  \item For each $x$, the type $P(x)$ is a mere proposition, and
  \item For each $x$ we have $\brck {P(x)}$.
  \end{enumerate}
  Then we have $\prd{x:A} P(x)$.
\end{cor}
\begin{proof}
  Immediate from the two assumptions and the previous lemma.
\end{proof}

The corollary also encapsulates a very useful technique of reasoning.
Namely, suppose we know that $\brck A$, and we want to use this to construct an element of some other type $B$.
We would like to use an element of $A$ in our construction of an element of $B$, but this is allowed only if $B$ is a mere proposition, so that we can apply the induction principle for the propositional truncation $\brck A$; the most we could hope to do in general is to show $\brck B$.
Instead, we can extend $B$ with additional data which characterizes \emph{uniquely} the object we wish to construct.
Specifically, we define a predicate $Q:B\to\type$ such that $\sm{x:B} Q(x)$ is a mere proposition.
Then from an element of $A$ we construct an element $b:B$ such that $Q(b)$, hence from $\brck A$ we can construct $\brck{\sm{x:B} Q(x)}$, and because $\brck{\sm{x:B} Q(x)}$ is equivalent to $\sm{x:B} Q(x)$ an element of $B$ may be projected from it.
An example can be found in \autoref{ex:decidable-choice}.

A similar issue arises in set-theoretic mathematics, although it manifests slightly
differently. If we are trying to define a function $f: A \to B$, and depending on an
element $a : A$ we are able to prove mere existence of some $b : B$, we are not done yet
because we need to actually pinpoint an element of~$B$, not just prove its existence.
One option is of course to refine the argument to unique existence of $b : B$, as we did in type theory. But in set theory the problem can often be avoided more simply by an application of the axiom of choice, which picks the required elements for us.
In homotopy type theory, however, quite apart from any desire to avoid choice, the available forms of choice are simply less applicable, since they require that the domain of choice be a \emph{set}.
Thus, if $A$ is not a set (such as perhaps a universe $\UU$), there is no consistent form of choice that will allow us to simply pick an element of $B$ for each $a : A$ to use in defining $f(a)$.

\index{unique!choice|)}%

\section{When are propositions truncated?}
\label{subsec:when-trunc}

\index{logic!of mere propositions|(}%
\index{mere proposition|(}%
\index{logic!truncated}%

At first glance, it may seem that the truncated versions of $+$ and $\Sigma$ are actually closer to the informal mathematical meaning of ``or'' and ``there exists'' than the untruncated ones.
Certainly, they are closer to the \emph{precise} meaning of ``or'' and ``there exists'' in the first-order logic \index{first-order!logic} which underlies formal set theory, since the latter makes no attempt to remember any witnesses to the truth of propositions.
However, it may come as a surprise to realize that the practice of \emph{informal} mathematics is often more accurately described by the untruncated forms.

\index{prime number}%
For example, consider a statement like ``every prime number is either $2$ or odd.''
The working mathematician feels no compunction about using this fact not only to prove \emph{theorems} about prime numbers, but also to perform \emph{constructions} on prime numbers, perhaps doing one thing in the case of $2$ and another in the case of an odd prime.
The end result of the construction is not merely the truth of some statement, but a piece of data which may depend on the parity of the prime number.
Thus, from a type-theoretic perspective, such a construction is naturally phrased using the induction principle for the coproduct type ``$(p=2)+(p\text{ is odd})$'', not its propositional truncation.

Admittedly, this is not an ideal example, since ``$p=2$'' and ``$p$ is odd'' are mutually exclusive, so that $(p=2)+(p\text{ is odd})$ is in fact already a mere proposition and hence equivalent to its truncation (see \autoref{ex:disjoint-or}).
More compelling examples come from the existential quantifier.
It is not uncommon to prove a theorem of the form ``there exists an $x$ such that \dots'' and then refer later on to ``the $x$ constructed in Theorem Y'' (note the definite article).
Moreover, when deriving further properties of this $x$, one may use phrases such as ``by the construction of $x$ in the proof of Theorem Y''.

A very common example is ``$A$ is isomorphic to $B$'', which strictly speaking means only that there exists \emph{some} isomorphism between $A$ and $B$.
But almost invariably, when proving such a statement, one exhibits a specific isomorphism or proves that some previously known map is an isomorphism, and it often matters later on what particular isomorphism was given.

Set-theoretically trained mathematicians often feel a twinge of guilt at such ``abuses of language''.\index{abuse!of language}
We may attempt to apologize for them, expunge them from final drafts, or weasel out of them with vague words like ``canonical''.
The problem is exacerbated by the fact that in formalized set theory, there is technically no way to ``construct'' objects at all --- we can only prove that an object with certain properties exists.
Untruncated logic in type theory thus captures some common practices of informal mathematics that the set theoretic reconstruction obscures.
(This is similar to how the univalence axiom validates the common, but formally unjustified, practice of identifying isomorphic objects.)

On the other hand, sometimes truncated logic is essential.
We have seen this in the statements of \LEM{} and \choice{}; some other examples will appear later on in the book.
Thus, we are faced with the problem: when writing informal type theory, what should we mean by the words ``or'' and ``there exists'' (along with common synonyms such as ``there is'' and ``we have'')?

A universal consensus may not be possible.
Perhaps depending on the sort of mathematics being done, one convention or the other may be more useful --- or, perhaps, the choice of convention may be irrelevant.
In this case, a remark at the beginning of a mathematical paper may suffice to inform the reader of the linguistic conventions in use therein.
However, even after one overall convention is chosen, the other sort of logic will usually arise at least occasionally, so we need a way to refer to it.
More generally, one may consider replacing the propositional truncation with another operation on types that behaves similarly, such as the double negation operation $A\mapsto \neg\neg A$, or the $n$-truncations to be considered in \autoref{cha:hlevels}.
As an experiment in exposition,  in what follows we will occasionally use \emph{adverbs}\index{adverb} to denote the application of such ``modalities'' as propositional truncation.

For instance, if untruncated logic is the default convention, we may use the adverb \define{merely}
\indexdef{merely}%
to denote propositional truncation.
Thus the phrase
\begin{center}
  ``there merely exists an $x:A$ such that $P(x)$''
\end{center}
indicates the type $\brck{\sm{x:A} P(x)}$.
Similarly, we will say that a type $A$ is \define{merely inhabited}
\indexdef{merely!inhabited}%
\indexdef{inhabited type!merely}%
to mean that its propositional truncation $\brck A$ is inhabited (i.e.\ that we have an unnamed element of it).
Note that this is a \emph{definition}\index{definition!of adverbs} of the adverb ``merely'' as it is to be used in our informal mathematical English, in the same way that we define nouns\index{noun} like ``group'' and ``ring'', and adjectives\index{adjective} like ``regular'' and ``normal'', to have precise mathematical meanings.
We are not claiming that the dictionary definition of ``merely'' refers to propositional truncation; the choice of word is meant only to remind the mathematician reader that a mere proposition contains ``merely'' the information of a truth value and nothing more.

On the other hand, if truncated logic is the current default convention, we may use an adverb such as \define{purely}
\indexdef{purely}%
or \define{constructively} to indicate its absence, so that
\begin{center}
``there purely exists an $x:A$ such that $P(x)$''
\end{center}
would denote the type $\sm{x:A} P(x)$.
We may also use ``purely'' or ``actually'' just to emphasize the absence of truncation, even when that is the default convention.

In this book we will continue using untruncated logic as the default convention, for a number of reasons.
\begin{enumerate}[label=(\arabic*)]
\item We want to encourage the newcomer to experiment with it, rather than sticking to truncated logic simply because it is more familiar.
\item Using truncated logic as the default in type theory suffers from the same sort of ``abuse of language''\index{abuse!of language} problems as set-theoretic foundations, which untruncated logic avoids.
  For instance, our definition of ``$\eqv A B$'' as the type of equivalences between $A$ and $B$, rather than its propositional truncation, means that to prove a theorem of the form ``$\eqv A B$'' is literally to construct a particular such equivalence.
  This specific equivalence can then be referred to later on.
\item We want to emphasize that the notion of ``mere proposition'' is not a fundamental part of type theory.
  As we will see in \autoref{cha:hlevels}, mere propositions are just the second rung on an infinite ladder, and there are also many other modalities not lying on this ladder at all.
\item Many statements that classically are mere propositions are no longer so in homotopy type theory.
  Of course, foremost among these is equality.
\item On the other hand, one of the most interesting observations of homotopy type theory is that a surprising number of types are \emph{automatically} mere propositions, or can be slightly modified to become so, without the need for any truncation.
  (See \autoref{thm:isprop-isprop,cha:equivalences,cha:hlevels,cha:category-theory,cha:set-math}.)
  Thus, although these types contain no data beyond a truth value, we can nevertheless use them to construct untruncated objects, since there is no need to use the induction principle of propositional truncation.
  This useful fact is more clumsy to express if propositional truncation is applied to all statements by default.
\item Finally, truncations are not very useful for most of the mathematics we will be doing in this book, so it is simpler to notate them explicitly when they occur.
\end{enumerate}

\index{mere proposition|)}%
\index{logic!of mere propositions|)}%

\section{Contractibility}
\label{sec:contractibility}

\index{type!contractible|(defstyle}%
\index{contractible!type|(defstyle}%

In \autoref{thm:inhabprop-eqvunit} we observed that a mere proposition which is inhabited must be equivalent to $\unit$,
\index{type!unit}%
and it is not hard to see that the converse also holds.
A type with this property is called \emph{contractible}.
Another equivalent definition of contractibility, which is also sometimes convenient, is the following.

\begin{defn}\label{defn:contractible}
  A type $A$ is \define{contractible},
  or a \define{singleton},
  \indexdef{type!singleton}%
  \indexsee{singleton type}{type, singleton}%
  if there is $a:A$, called the \define{center of contraction},
  \indexdef{center!of contraction}%
  such that $a=x$ for all $x:A$.
  We denote the specified path $a=x$ by $\contr_x$.
\end{defn}

In other words, the type $\iscontr(A)$ is defined to be
\[ \iscontr(A) \defeq \sm{a:A} \prd{x:A}(a=x). \]
Note that under the usual propositions-as-types reading, we can pronounce $\iscontr(A)$ as ``$A$ contains exactly one element'', or more precisely ``$A$ contains an element, and every element of $A$ is equal to that element''.

\begin{rmk}
  We can also pronounce $\iscontr(A)$ more topologically as ``there is a point $a:A$ such that for all $x:A$ there exists a path from $a$ to $x$''.
  Note that to a classical ear, this sounds like a definition of \emph{connectedness} rather than contractibility.
  \index{continuity of functions in type theory@``continuity'' of functions in type theory}%
  \index{functoriality of functions in type theory@``functoriality'' of functions in type theory}%
  The point is that the meaning of ``there exists'' in this sentence is a continuous/natural one.
  A more correct way to express connectedness would be $\sm{a:A}\prd{x:A} \brck{a=x}$; see \autoref{thm:connected-pointed}.
\end{rmk}

\begin{lem}\label{thm:contr-unit}
  For a type $A$, the following are logically equivalent.
  \begin{enumerate}
  \item $A$ is contractible in the sense of \autoref{defn:contractible}.\label{item:contr}
  \item $A$ is a mere proposition, and there is a point $a:A$.\label{item:contr-inhabited-prop}
  \item $A$ is equivalent to \unit.\label{item:contr-eqv-unit}
  \end{enumerate}
\end{lem}
\begin{proof}
  If $A$ is contractible, then it certainly has a point $a:A$ (the center of contraction), while for any $x,y:A$ we have $x=a=y$; thus $A$ is a mere proposition.
  Conversely, if we have $a:A$ and $A$ is a mere proposition, then for any $x:A$ we have $x=a$; thus $A$ is contractible.
  And we showed~\ref{item:contr-inhabited-prop}$\Rightarrow$\ref{item:contr-eqv-unit} in \autoref{thm:inhabprop-eqvunit}, while the converse follows since \unit easily has property~\ref{item:contr-inhabited-prop}.
\end{proof}

\begin{lem}\label{thm:isprop-iscontr}
  For any type $A$, the type $\iscontr(A)$ is a mere proposition.
\end{lem}
\begin{proof}
  Suppose given $c,c':\iscontr(A)$.
  We may assume $c\jdeq(a,p)$ and $c'\jdeq(a',p')$ for $a,a':A$ and $p:\prd{x:A} (a=x)$ and $p':\prd{x:A} (a'=x)$.
  By the characterization of paths in $\Sigma$-types, to show $c=c'$ it suffices to exhibit $q:a=a'$ such that $\trans{q}{p}=p'$.

  We choose $q\defeq p(a')$.
  For the other equality, by function extensionality we must show that $(\trans q p)(x)=p'(x)$ for any $x:A$.
  For this, it will suffice to show that for any $x,y:A$ and $u:x=y$ we have $u= \opp{p(x)} \ct p(y)$, since then we would have $(\trans q p)(x) = \opp{p(a')} \ct p(x) = p'(x)$.
  But now we can invoke path induction to assume that $x\jdeq y$ and $u\jdeq \refl{x}$.
  In this case our goal is to show that $\refl x = \opp{p(x)} \ct p(x)$, which is just the inversion law for paths.
\end{proof}

\begin{cor}\label{thm:contr-contr}
  If $A$ is contractible, then so is $\iscontr(A)$.
\end{cor}
\begin{proof}
  By \autoref{thm:isprop-iscontr} and \autoref{thm:contr-unit}\ref{item:contr-inhabited-prop}.
\end{proof}

Like mere propositions, contractible types are preserved by many type constructors.
For instance, we have:

\begin{lem}\label{thm:contr-forall}
  If $P:A\to\type$ is a type family such that each $P(a)$ is contractible, then $\prd{x:A} P(x)$ is contractible.
\end{lem}
\begin{proof}
  By \autoref{thm:isprop-forall}, $\prd{x:A} P(x)$ is a mere proposition since each $P(x)$ is.
  But it also has an element, namely the function sending each $x:A$ to the center of contraction of $P(x)$.
  Thus by \autoref{thm:contr-unit}\ref{item:contr-inhabited-prop}, $\prd{x:A} P(x)$ is contractible.
\end{proof}

\index{function extensionality}%
(In fact, the statement of \autoref{thm:contr-forall} is equivalent to the function extensionality axiom.
See~\autoref{sec:univalence-implies-funext}.)

Of course, if $A$ is equivalent to $B$ and $A$ is contractible, then so is $B$.
More generally, it suffices for $B$ to be a \emph{retract} of $A$.
By definition, a \define{retraction}
\indexdef{retraction}%
\indexdef{function!retraction}%
is a function $r : A \to B$ such that there exists a function $s : B \to A$, called its \define{section},
\indexdef{section}%
\indexdef{function!section}%
and a homotopy $\epsilon:\prd{y:B} (r(s(y))=y)$; then we say that $B$ is a \define{retract}%
\indexdef{retract!of a type}
of $A$.

\begin{lem}\label{thm:retract-contr}
  If $B$ is a retract of $A$, and $A$ is contractible, then so is $B$.
\end{lem}
\begin{proof}
  Let $a_0 : A$ be the center of contraction.
  We claim that $b_0 \defeq r(a_0) : B$ is a center of contraction for $B$.
  Let $b : B$; we need a path $b = b_0$.
  But we have $\epsilon_b : r(s(b)) = b$ and $\contr_{s(b)} : s(b) = a_0$, so by composition
  \[ \opp{\epsilon_b} \ct \ap{r}{\contr_{s(b)}} : b = r(a_0) \jdeq b_0. \qedhere\]
\end{proof}

Contractible types may not seem very interesting, since they are all equivalent to \unit.
One reason the notion is useful is that sometimes a collection of individually nontrivial data will collectively form a contractible type.
An important example is the space of paths with one free endpoint.
As we will see in \autoref{sec:identity-systems}, this fact essentially
encapsulates the based path induction principle for identity types.

\begin{lem}\label{thm:contr-paths}
  For any $A$ and any $a:A$, the type $\sm{x:A} (a=x)$ is contractible.
\end{lem}
\begin{proof}
  We choose as center the point $(a,\refl a)$.
  Now suppose $(x,p):\sm{x:A}(a=x)$; we must show $(a,\refl a) = (x,p)$.
  By the characterization of paths in $\Sigma$-types, it suffices to exhibit $q:a=x$ such that $\trans{q}{\refl a} = p$.
  But we can take $q\defeq p$, in which case $\trans{q}{\refl a} = p$ follows from the characterization of transport in path types.
\end{proof}

When this happens, it can allow us to simplify a complicated construction up to equivalence, using the informal principle that contractible data can be freely ignored.
This principle consists of many lemmas, most of which we leave to the reader; the following is an example.

\begin{lem}\label{thm:omit-contr}
  Let $P:A\to\type$ be a type family.
  \begin{enumerate}
  \item If each $P(x)$ is contractible, then $\sm{x:A} P(x)$ is equivalent to $A$.\label{item:omitcontr1}
  \item If $A$ is contractible with center $a$, then $\sm{x:A} P(x)$ is equivalent to $P(a)$.\label{item:omitcontr2}
  \end{enumerate}
\end{lem}
\begin{proof}
  In the situation of~\ref{item:omitcontr1}, we show that $\proj1:\sm{x:A} P(x) \to A$ is an equivalence.
  For quasi-inverse we define $g(x)\defeq (x,c_x)$ where $c_x$ is the center of $P(x)$.
  The composite $\proj1 \circ g$ is obviously $\idfunc[A]$, whereas the opposite composite is homotopic to the identity by using the contractions of each $P(x)$.

  We leave the proof of~\ref{item:omitcontr2} to the reader (see \autoref{ex:omit-contr2}).
\end{proof}

Another reason contractible types are interesting is that they extend the ladder of $n$-types mentioned in \autoref{sec:basics-sets} downwards one more step.

\begin{lem}\label{thm:prop-minusonetype}
  A type $A$ is a mere proposition if and only if for all $x,y:A$, the type $\id[A]xy$ is contractible.
\end{lem}
\begin{proof}
  For ``if'', we simply observe that any contractible type is inhabited.
  For ``only if'', we observed in \autoref{subsec:hprops} that every mere proposition is a set, so that each type $\id[A]xy$ is a mere proposition.
  But it is also inhabited (since $A$ is a mere proposition), and hence by \autoref{thm:contr-unit}\ref{item:contr-inhabited-prop} it is contractible.
\end{proof}

Thus, contractible types may also be called \define{$(-2)$-types}.
They are the bottom rung of the ladder of $n$-types, and will be the base case of the recursive definition of $n$-types in \autoref{cha:hlevels}.

\index{type!contractible|)}%
\index{contractible!type|)}%

\sectionNotes

The fact that it is possible to define sets, mere propositions, and contractible types in type theory, with all higher homotopies automatically taken care of as in \autoref{sec:basics-sets,subsec:hprops,sec:contractibility}, was first observed by Voevodsky.
In fact, he defined the entire hierarchy of $n$-types by induction, as we will do in \autoref{cha:hlevels}.\index{n-type@$n$-type!definable in type theory}%

\autoref{thm:not-dneg,thm:not-lem} rely in essence on a classical theorem of Hedberg,
\index{Hedberg's theorem}%
\index{theorem!Hedberg's}%
which we will prove in \autoref{sec:hedberg}.
The implication that the propositions-as-types form of \LEM{} contradicts univalence was observed by Mart\'\i n Escard\'o on the \Agda mailing list.
The proof we have given of \autoref{thm:not-dneg} is due to Thierry Coquand.

The propositional truncation was introduced in the extensional type theory of
\index{proof!assistant!\NuPRL}%
\NuPRL in 1983 by Constable~\cite{Con85} as an
application of ``subset'' and ``quotient'' types.  What is here called the
``propositional truncation'' was called ``squashing'' in the \NuPRL type theory~\cite{constable+86nuprl-book}.
Rules characterizing the propositional truncation directly, still in extensional type theory, were given in~\cite{ab:bracket-types}.
The intensional version in homotopy type theory was constructed by Voevodsky using an impredicative\index{impredicative!truncation} quantification, and later by Lumsdaine using higher inductive types (see \autoref{sec:hittruncations}).

\index{propositional!resizing}%
Voevodsky~\cite{Universe-poly} has proposed resizing rules of the kind considered in \autoref{subsec:prop-subsets}.
\index{axiom!of reducibility}\index{resizing}%
These are clearly related to the notorious \emph{axiom of reducibility} proposed by Russell in his and Whitehead's \emph{Principia Mathematica}~\cite{PM2}.\index{Russell, Bertrand}

The adverb ``purely'' as used to refer to untruncated logic is a reference to the use of monadic modalities to model effects in programming languages; see \autoref{sec:modalities} and the Notes to \autoref{cha:hlevels}.

There are many different ways in which logic can be treated relative to type theory.
For instance, in addition to the plain propositions-as-types logic described in \autoref{sec:pat}, and the alternative which uses mere propositions only as described in \autoref{subsec:logic-hprop}, one may introduce a separate ``sort'' of propositions, which behave somewhat like types but are not identified with them.
This is the approach taken in logic enriched type theory~\cite{aczel2002collection} and in some presentations of the internal languages of toposes\index{topos} and related categories (e.g.~\cite{jacobs1999categorical,elephant}), as well as in the proof assistant \Coq.\index{proof!assistant!Coq@\textsc{Coq}}
Such an approach is more general, but less powerful.
For instance, the principle of unique choice (\autoref{sec:unique-choice}) fails in the category of so-called setoids in \Coq~\cite{Spiwack}, in logic enriched type theory~\cite{aczel2002collection}, and in minimal type theory~\cite{maietti2005toward}.\index{setoid}
Thus, the univalence axiom makes our type theory behave more like the internal logic of a topos;~see also \autoref{cha:set-math}.\index{topos}

Martin-L\"of~\cite{martin2006100} provides a discussion on the history of axioms of choice.
Of course, constructive and intuitionistic mathematics has a long and complicated history, which we will not delve into here; see for instance~\cite{TroelstraI,TroelstraII}.

\sectionExercises

\begin{ex}
  Prove that if $\eqv A B$ and $A$ is a set, then so is $B$.
\end{ex}

\begin{ex}\label{ex:isset-coprod}
  Prove that if $A$ and $B$ are sets, then so is $A+B$.
\end{ex}

\begin{ex}\label{ex:isset-sigma}
  Prove that if $A$ is a set and $B:A\to \type$ is a type family such that $B(x)$ is a set for all $x:A$, then $\sm{x:A} B(x)$ is a set.
\end{ex}

\begin{ex}\label{ex:prop-endocontr}
  Show that $A$ is a mere proposition if and only if $A\to A$ is contractible.
\end{ex}

\begin{ex}\label{ex:prop-inhabcontr}
  Show that $\eqv{\isprop(A)}{(A\to\iscontr(A))}$.
\end{ex}

\begin{ex}\label{ex:lem-mereprop}
  Show that if $A$ is a mere proposition, then so is $A+(\neg A)$.
  Thus, there is no need to insert a propositional truncation in~\eqref{eq:lem}.
\end{ex}

\begin{ex}\label{ex:disjoint-or}
  More generally, show that if $A$ and $B$ are mere propositions and $\neg(A\times B)$, then $A+B$ is also a mere proposition.
\end{ex}




\begin{ex}\label{ex:brck-qinv}
  Assuming that some type $\isequiv(f)$ satisfies conditions~\ref{item:be1}--\ref{item:be3} of \autoref{sec:basics-equivalences}, show that the type $\brck{\qinv(f)}$ satisfies the same conditions and is equivalent to $\isequiv(f)$.
\end{ex}

\begin{ex}
  Show that if \LEM{} holds, then the type $\prop \defeq \sm{A:\type} \isprop(A)$ is equivalent to \bool.
\end{ex}

\begin{ex}\label{ex:lem-impred}
  Show that if $\UU_{i+1}$ satisfies \LEM{}, then the canonical inclusion $\prop_{\UU_i} \to \prop_{\UU_{i+1}}$ is an equivalence.
\end{ex}

\begin{ex}
  Show that it is not the case that for all $A:\type$ we have $\brck{A} \to A$.
  (However, there can be particular types for which $\brck{A}\to A$.
  \autoref{ex:brck-qinv} implies that $\qinv(f)$ is such.)
\end{ex}

\begin{ex}
  \index{axiom!of choice}%
  Show that if \LEM{} holds, then for all $A:\type$ we have $\bbrck{(\brck A \to A)}$.
  (This property is a very simple form of the axiom of choice, which can fail in the absence of \LEM{}; see~\cite{krausgeneralizations}.)
\end{ex}

\begin{ex}
  We showed in \autoref{thm:not-lem} that the following naive form of \LEM{} is inconsistent with univalence:
  \[ \prd{A:\type} (A+(\neg A)) \]
  In the absence of univalence, this axiom is consistent.
  However, show that it implies the axiom of choice~\eqref{eq:ac}.
\end{ex}

\begin{ex}\label{ex:lem-brck}
  Show that assuming \LEM{}, the double negation $\neg \neg A$ has the same universal property as the propositional truncation $\brck A$, and is therefore equivalent to it.
  Thus, under \LEM{}, the propositional truncation can be defined rather than taken as a separate type former.
\end{ex}

\begin{ex}\label{ex:impred-brck}
  \index{propositional!resizing}%
  Show that if we assume propositional resizing as in \autoref{subsec:prop-subsets}, then the type
  \[\prd{P:\prop} \Parens{(A\to P)\to P}\]
  has the same universal property as $\brck A$.
  Thus, we can also define the propositional truncation in this case.
\end{ex}

\begin{ex}
  Assuming \LEM{}, show that double negation commutes with universal quantification of mere propositions over sets.
  That is, show that if $X$ is a set and each $Y(x)$ is a mere proposition, then \LEM{} implies
  \begin{equation}
    \eqv{\Parens{\prd{x:X} \neg\neg Y(x)}}{\Parens{\neg\neg \prd{x:X} Y(x)}}.\label{eq:dnshift}
  \end{equation}
  Observe that if we assume instead that each $Y(x)$ is a set, then~\eqref{eq:dnshift} becomes equivalent to the axiom of choice~\eqref{eq:epis-split}.
\end{ex}

\begin{ex}\label{ex:prop-trunc-ind}
  \index{induction principle!for truncation}%
  Show that the rules for the propositional truncation given in \autoref{subsec:prop-trunc} are sufficient to imply the following induction principle: for any type family $B:\brck A \to \type$ such that each $B(x)$ is a mere proposition, if for every $a:A$ we have $B(\bproj a)$, then for every $x:\brck A$ we have $B(x)$.
\end{ex}

\begin{ex}\label{ex:lem-ldn}
  Show that the law of excluded middle~\eqref{eq:lem} and the law of double negation~\eqref{eq:ldn} are logically equivalent.
\end{ex}

\begin{ex}\label{ex:decidable-choice}
  Suppose $P:\nat\to\type$ is a decidable family of mere propositions.
  Prove that
  \[ \Brck{\sm{n:\nat} P(n)} \;\to\; \sm{n:\nat}P(n).\]
\end{ex}

\begin{ex}\label{ex:omit-contr2}
  Prove \autoref{thm:omit-contr}\ref{item:omitcontr2}: if $A$ is contractible with center $a$, then $\sm{x:A} P(x)$ is equivalent to $P(a)$.
\end{ex}


\chapter{Equivalences}
\label{cha:equivalences}

We now study in more detail the notion of \emph{equivalence of types} that was introduced briefly in \autoref{sec:basics-equivalences}.
Specifically, we will give several different ways to define a type $\isequiv(f)$ having the properties mentioned there.
Recall that we wanted $\isequiv(f)$ to have the following properties, which we restate here:
\begin{enumerate}
\item $\qinv(f) \to \isequiv (f)$.\label{item:beb1}
\item $\isequiv (f) \to \qinv(f)$.\label{item:beb2}
\item $\isequiv(f)$ is a mere proposition.\label{item:beb3}
\end{enumerate}
Here $\qinv(f)$ denotes the type of quasi-inverses to $f$:
\begin{equation*}
  \sm{g:B\to A} \big((f \circ g \htpy \idfunc[B]) \times (g\circ f \htpy \idfunc[A])\big).
\end{equation*}
By function extensionality, it follows that $\qinv(f)$ is equivalent to the type
\begin{equation*}
  \sm{g:B\to A} \big((f \circ g = \idfunc[B]) \times (g\circ f = \idfunc[A])\big).
\end{equation*}
We will define three different types having properties~\ref{item:beb1}--\ref{item:beb3}, which we call
\begin{itemize}
\item half adjoint equivalences,
\item bi-invertible maps,
  \index{function!bi-invertible}
  and
\item contractible functions.
\end{itemize}
We will also show that all these types are equivalent.
These names are intentionally somewhat cumbersome, because after we know that they are all equivalent and have properties~\ref{item:beb1}--\ref{item:beb3}, we will revert to saying simply ``equivalence'' without needing to specify which particular definition we choose.
But for purposes of the comparisons in this chapter, we need different names for each definition.

Before we examine the different notions of equivalence, however, we give a little more explanation of why a different concept than quasi-invertibility is needed.

\section{Quasi-inverses}
\label{sec:quasi-inverses}

\index{quasi-inverse|(}%
We have said that $\qinv(f)$ is unsatisfactory because it is not a mere proposition, whereas we would rather that a given function can ``be an equivalence'' in at most one way.
However, we have given no evidence that $\qinv(f)$ is not a mere proposition.
In this section we exhibit a specific counterexample.

\begin{lem}\label{lem:qinv-autohtpy}
  If $f:A\to B$ is such that $\qinv (f)$ is inhabited, then
  \[\eqv{\qinv(f)}{\Parens{\prd{x:A}(x=x)}}.\]
\end{lem}
\begin{proof}
  By assumption, $f$ is an equivalence; that is, we have $e:\isequiv(f)$ and so $(f,e):\eqv A B$.
  By univalence, $\idtoeqv:(A=B) \to (\eqv A B)$ is an equivalence, so we may assume that $(f,e)$ is of the form $\idtoeqv(p)$ for some $p:A=B$.
  Then by path induction, we may assume $p$ is $\refl{A}$, in which case $\idtoeqv(p)$ is $\idfunc[A]$.
  Thus we are reduced to proving $\eqv{\qinv(\idfunc[A])}{(\prd{x:A}(x=x))}$.
  Now by definition we have
  \[ \qinv(\idfunc[A]) \jdeq
  \sm{g:A\to A} \big((g \htpy \idfunc[A]) \times (g \htpy \idfunc[A])\big).
  \]
  By function extensionality, this is equivalent to
  \[ \sm{g:A\to A} \big((g = \idfunc[A]) \times (g = \idfunc[A])\big).
  \]
  And by \autoref{ex:sigma-assoc}, this is equivalent to
  \[ \sm{h:\sm{g:A\to A} (g = \idfunc[A])} (\proj1(h) = \idfunc[A])
  \]
  However, by \autoref{thm:contr-paths}, $\sm{g:A\to A} (g = \idfunc[A])$ is contractible with center $\idfunc[A]$; therefore by \autoref{thm:omit-contr} this type is equivalent to $\idfunc[A] = \idfunc[A]$.
  And by function extensionality, $\idfunc[A] = \idfunc[A]$ is equivalent to $\prd{x:A} x=x$.
\end{proof}

\noindent
We remark that \autoref{ex:qinv-autohtpy-no-univalence} asks for a proof of the above lemma which avoids univalence.

Thus, what we need is some $A$ which admits a nontrivial element of $\prd{x:A}(x=x)$.
Thinking of $A$ as a higher groupoid, an inhabitant of $\prd{x:A}(x=x)$ is a natural transformation\index{natural!transformation} from the identity functor of $A$ to itself.
Such transformations are said to form the \define{center of a category},
\index{center!of a category}%
\index{category!center of}%
since the naturality axiom requires that they commute with all morphisms.
Classically, if $A$ is simply a group regarded as a one-object groupoid, then this yields precisely its center in the usual group-theoretic sense.
This provides some motivation for the following.

\begin{lem}\label{lem:autohtpy}
  Suppose we have a type $A$ with $a:A$ and $q:a=a$ such that
  \begin{enumerate}
  \item The type $a=a$ is a set.\label{item:autohtpy1}
  \item For all $x:A$ we have $\brck{a=x}$.\label{item:autohtpy2}
  \item For all $p:a=a$ we have $p\ct q = q \ct p$.\label{item:autohtpy3}
  \end{enumerate}
  Then there exists $f:\prd{x:A} (x=x)$ with $f(a)=q$.
\end{lem}
\begin{proof}
  Let $g:\prd{x:A} \brck{a=x}$ be as given by~\ref{item:autohtpy2}.  First we
  observe that each type $\id[A]xy$ is a set.  For since being a set is a mere
  proposition, we may apply the induction principle of propositional truncation, and assume that $g(x)=\bproj
  p$ and $g(y)=\bproj q$ for $p:a=x$ and $q:a=y$.  In this case, composing with
  $p$ and $\opp{q}$ yields an equivalence $\eqv{(x=y)}{(a=a)}$.  But $(a=a)$ is
  a set by~\ref{item:autohtpy1}, so $(x=y)$ is also a set.

  Now, we would like to define $f$ by assigning to each $x$ the path $\opp{g(x)}
  \ct q \ct g(x)$, but this does not work because $g(x)$ does not inhabit $a=x$
  but rather $\brck{a=x}$, and the type $(x=x)$ may not be a mere proposition,
  so we cannot use induction on propositional truncation.  Instead we can apply
  the technique mentioned in \autoref{sec:unique-choice}: we characterize
  uniquely the object we wish to construct.  Let us define, for each $x:A$, the
  type
  \[ B(x) \defeq \sm{r:x=x} \prd{s:a=x} (r = \opp s \ct q\ct s).\]
  We claim that $B(x)$ is a mere proposition for each $x:A$.
  Since this claim is itself a mere proposition, we may again apply induction on
  truncation and assume that $g(x) = \bproj p$ for some $p:a=x$.
  Now suppose given $(r,h)$ and $(r',h')$ in $B(x)$; then we have
  \[ h(p) \ct \opp{h'(p)} : r = r'. \]
  It remains to show that $h$ is identified with $h'$ when transported along this equality, which by transport in identity types and function types (\autoref{sec:compute-paths,sec:compute-pi}), reduces to showing
  \[ h(s) = h(p) \ct \opp{h'(p)} \ct h'(s) \]
  for any $s:a=x$.
  But each side of this is an equality between elements of $(x=x)$, so it follows from our above observation that $(x=x)$ is a set.

  Thus, each $B(x)$ is a mere proposition; we claim that $\prd{x:A} B(x)$.
  Given $x:A$, we may now invoke the induction principle of propositional truncation to assume that $g(x) = \bproj p$ for $p:a=x$.
  We define $r \defeq \opp p \ct q \ct p$; to inhabit $B(x)$ it remains to show that for any $s:a=x$ we have
  $r = \opp s \ct q \ct s$.
  Manipulating paths, this reduces to showing that $q\ct (p\ct \opp s) = (p\ct \opp s) \ct q$.
  But this is just an instance of~\ref{item:autohtpy3}.
\end{proof}

\begin{thm}
  There exist types $A$ and $B$ and a function $f:A\to B$ such that $\qinv(f)$ is not a mere proposition.
\end{thm}
\begin{proof}
  It suffices to exhibit a type $X$ such that $\prd{x:X} (x=x)$ is not a mere proposition.
  Define $X\defeq \sm{A:\type} \brck{\bool=A}$, as in the proof of \autoref{thm:no-higher-ac}.
  It will suffice to exhibit an $f:\prd{x:X} (x=x)$ which is unequal to $\lam{x} \refl{x}$.

  Let $a \defeq (\bool,\bproj{\refl{\bool}}) : X$, and let $q:a=a$ be the path corresponding to the nonidentity equivalence $e:\eqv\bool\bool$ defined by $e(\bfalse)\defeq\btrue$ and $e(\btrue)\defeq\bfalse$.
  We would like to apply \autoref{lem:autohtpy} to build an $f$.
  By definition of $X$, equalities in subset types (\autoref{subsec:prop-subsets}), and univalence, we have $\eqv{(a=a)}{(\eqv{\bool}{\bool})}$, which is a set, so~\ref{item:autohtpy1} holds.
  Similarly, by definition of $X$ and equalities in subset types we have~\ref{item:autohtpy2}.
  Finally, \autoref{ex:eqvboolbool} implies that every equivalence $\eqv\bool\bool$ is equal to either $\idfunc[\bool]$ or $e$, so we can show~\ref{item:autohtpy3} by a four-way case analysis.

  Thus, we have $f:\prd{x:X} (x=x)$ such that $f(a) = q$.
  Since $e$ is not equal to $\idfunc[\bool]$, $q$ is not equal to $\refl{a}$, and thus $f$ is not equal to $\lam{x} \refl{x}$.
  Therefore, $\prd{x:X} (x=x)$ is not a mere proposition.
\end{proof}

More generally, \autoref{lem:autohtpy} implies that any ``Eilenberg--Mac Lane space'' $K(G,1)$, where $G$ is a nontrivial abelian\index{group!abelian} group, will provide a counterexample; see \autoref{cha:homotopy}.
The type $X$ we used turns out to be equivalent to $K(\mathbb{Z}_2,1)$.
In \autoref{cha:hits} we will see that the circle $\Sn^1 = K(\mathbb{Z},1)$ is another easy-to-describe example.

We now move on to describing better notions of equivalence.

\index{quasi-inverse|)}%

\section{Half adjoint equivalences}
\label{sec:hae}

\index{equivalence!half adjoint|(defstyle}%
\index{half adjoint equivalence|(defstyle}%
\index{adjoint!equivalence!of types, half|(defstyle}%

In \autoref{sec:quasi-inverses} we concluded that $\qinv(f)$ is equivalent to $\prd{x:A} (x=x)$ by discarding a contractible type.
Roughly, the type $\qinv(f)$ contains three data $g$, $\eta$, and $\epsilon$, of which two ($g$ and $\eta$) could together be seen to be contractible when $f$ is an equivalence.
The problem is that removing these data left one remaining ($\epsilon$).
In order to solve this problem, the idea is to add one \emph{additional} datum which, together with $\epsilon$, forms a contractible type.

\begin{defn}\label{defn:ishae}
  A function $f:A\to B$ is a \define{half adjoint equivalence}
  if there are $g:B\to A$ and homotopies $\eta: g \circ f \htpy \idfunc[A]$ and $\epsilon:f \circ g \htpy \idfunc[B]$ such that there exists a homotopy
  \[\tau : \prd{x:A} \map{f}{\eta x} = \epsilon(fx).\]
\end{defn}

Thus we have a type $\ishae(f)$, defined to be
\begin{equation*}
  \sm{g:B\to A}{\eta: g \circ f \htpy \idfunc[A]}{\epsilon:f \circ g \htpy \idfunc[B]} \prd{x:A} \map{f}{\eta x} = \epsilon(fx).
\end{equation*}
Note that in the above definition, the coherence\index{coherence} condition relating $\eta$ and $\epsilon$ only involves $f$.
We might consider instead an analogous coherence condition involving $g$:
\[\upsilon : \prd{y:B} \map{g}{\epsilon y} = \eta(gy)\]
and a resulting analogous definition $\ishae'(f)$.

Fortunately, it turns out each of the conditions implies the other one:

\begin{lem}\label{lem:coh-equiv}
For functions $f : A \to B$ and $g:B\to A$ and homotopies $\eta: g \circ f \htpy \idfunc[A]$ and $\epsilon:f \circ g \htpy \idfunc[B]$, the following conditions are logically equivalent:
\begin{itemize}
\item $\prd{x:A} \map{f}{\eta x} = \epsilon(fx)$
\item $\prd{y:B} \map{g}{\epsilon y} = \eta(gy)$
\end{itemize}
\end{lem}
\begin{proof}
  It suffices to show one direction; the other one is obtained by replacing $A$, $f$, and $\eta$ by $B$, $g$, and $\epsilon$ respectively.
  Let $\tau : \prd{x:A}\;\map{f}{\eta x} = \epsilon(fx)$.
  Fix $y : B$.
  Using naturality of $\epsilon$ and applying $g$, we get the following commuting diagram of paths:
\[\uppercurveobject{{ }}\lowercurveobject{{ }}\twocellhead{{ }}
  \xymatrix@C=3pc{gfgfgy \ar@{=}^-{gfg(\epsilon y)}[r] \ar@{=}_{g(\epsilon (fgy))}[d] & gfgy \ar@{=}^{g(\epsilon y)}[d] \\ gfgy \ar@{=}_{g(\epsilon y)}[r] & gy
  }\]
Using $\tau(gy)$ on the left side of the diagram gives us
\[\uppercurveobject{{ }}\lowercurveobject{{ }}\twocellhead{{ }}
  \xymatrix@C=3pc{gfgfgy \ar@{=}^-{gfg(\epsilon y)}[r] \ar@{=}_{gf(\eta (gy))}[d] & gfgy \ar@{=}^{g(\epsilon y)}[d] \\ gfgy \ar@{=}_{g(\epsilon y)}[r] & gy
  }\]
Using the commutativity of $\eta$ with $g \circ f$ (\autoref{cor:hom-fg}), we have
\[\uppercurveobject{{ }}\lowercurveobject{{ }}\twocellhead{{ }}
  \xymatrix@C=3pc{gfgfgy \ar@{=}^-{gfg(\epsilon y)}[r] \ar@{=}_{\eta (gfgy)}[d] & gfgy \ar@{=}^{g(\epsilon y)}[d] \\ gfgy \ar@{=}_{g(\epsilon y)}[r] & gy
  }\]
However, by naturality of $\eta$ we also have
\[\uppercurveobject{{ }}\lowercurveobject{{ }}\twocellhead{{ }}
  \xymatrix@C=3pc{gfgfgy \ar@{=}^-{gfg(\epsilon y)}[r] \ar@{=}_{\eta (gfgy)}[d] & gfgy \ar@{=}^{\eta(gy)}[d] \\ gfgy \ar@{=}_{g(\epsilon y)}[r] & gy 
  }\]
Thus, canceling all but the right-hand homotopy, we have $g(\epsilon y) = \eta(g y)$ as desired.
\end{proof}

However, it is important that we do not include \emph{both} $\tau$ and $\upsilon$ in the definition of $\ishae (f)$ (whence the name ``\emph{half} adjoint equivalence'').
If we did, then after canceling contractible types we would still have one remaining datum --- unless we added another higher coherence condition.
In general, we expect to get a well-behaved type if we cut off after an odd number of coherences.

Of course, it is obvious that $\ishae(f) \to\qinv(f)$: simply forget the coherence datum.
The other direction is a version of a standard argument from homotopy theory and category theory.

\begin{thm}\label{thm:equiv-iso-adj}
  For any $f:A\to B$ we have $\qinv(f)\to\ishae(f)$.
\end{thm}
\begin{proof}
Suppose that $(g,\eta,\epsilon)$ is a quasi-inverse for $f$. We have to provide
a quadruple $(g',\eta',\epsilon',\tau)$ witnessing that $f$ is a half adjoint equivalence. To
define $g'$ and $\eta'$, we can just make the obvious choice by setting $g'
\defeq g$ and $\eta'\defeq \eta$. However, in the definition of $\epsilon'$ we
need start worrying about the construction of $\tau$, so we cannot just follow our nose
and take $\epsilon'$ to be $\epsilon$. Instead, we take
\begin{equation*}
\epsilon'(b) \defeq (\epsilon(b)\ct \ap{f}{\eta(g(b))})\ct \opp{\epsilon(f(g(b)))}.
\end{equation*}
Now we need to find
\begin{equation*}
\tau(a):(\epsilon(f(a))\ct \ap{f}{\eta(g(f(a)))})\ct \opp{\epsilon(f(g(f(a))))}=\ap{f}{\eta(a)}.
\end{equation*}
Note first that by \autoref{cor:hom-fg}, we have 
$\eta(g(f(a)))=\ap{g}{\ap{f}{\eta(a)}}$. Therefore, we can apply
\autoref{lem:htpy-natural} to compute
\begin{align*}
\epsilon(f(a))\ct \ap{f}{\eta(g(f(a)))}
& = \epsilon(f(a))\ct \ap{f}{\ap{g}{\ap{f}{\eta(a)}}}\\
& = \ap{f}{\eta(a)}\ct\epsilon(f(g(f(a))))
\end{align*}
from which we get the desired path $\tau(a)$.
\end{proof}

Combining this with \autoref{lem:coh-equiv} (or symmetrizing the proof), we also have $\qinv(f)\to\ishae'(f)$.

It remains to show that $\ishae(f)$ is a mere proposition.
For this, we will need to know that the fibers of an equivalence are contractible.

\begin{defn}\label{defn:homotopy-fiber}
  The \define{fiber}
  \indexdef{fiber}%
  \indexsee{function!fiber of}{fiber}%
  of a map $f:A\to B$ over a point $y:B$ is
  \[ \hfib f y \defeq \sm{x:A} (f(x) = y).\]
\end{defn}

In homotopy theory, this is what would be called the \emph{homotopy fiber} of $f$.
The path lemmas in \autoref{sec:computational} yield the following characterization of paths in fibers:

\begin{lem}\label{lem:hfib}
  For any $f : A \to B$, $y : B$, and $(x,p),(x',p') : \hfib{f}{y}$, we have
  \[ \big((x,p) = (x',p')\big) \eqvsym \Parens{\sm{\gamma : x = x'} f(\gamma) \ct p' = p} \qedhere\]
\end{lem}

\begin{thm}\label{thm:contr-hae}
  If $f:A\to B$ is a half adjoint equivalence, then for any $y:B$ the fiber $\hfib f y$ is contractible.
\end{thm}
\begin{proof}
  Let $(g,\eta,\epsilon,\tau) : \ishae(f)$, and fix $y : B$.
  As our center of contraction for $\hfib{f}{y}$ we choose $(gy, \epsilon y)$.
  Now take any $(x,p) : \hfib{f}{y}$; we want to construct a path from $(gy, \epsilon y)$ to $(x,p)$.
  By \autoref{lem:hfib}, it suffices to give a path $\gamma : \id{gy}{x}$ such that $\ap f\gamma \ct p = \epsilon y$.
  We put $\gamma \defeq \opp{g(p)} \ct \eta x$.
  Then we have 
  \begin{align*}
    f(\gamma) \ct p & = \opp{fg(p)} \ct f (\eta x) \ct p \\
    & = \opp{fg(p)} \ct \epsilon(fx) \ct p \\
    & = \epsilon y
  \end{align*}
  where the second equality follows by $\tau x$ and the third equality is naturality of $\epsilon$.
\end{proof}

We now define the types which encapsulate contractible pairs of data.
The following types put together the quasi-inverse $g$ with one of the homotopies.

\begin{defn}\label{defn:linv-rinv}
  Given a function $f:A\to B$, we define the types 
    \begin{align*}
      \linv(f) &\defeq \sm{g:B\to A} (g\circ f\htpy \idfunc[A])\\
      \rinv(f) &\defeq \sm{g:B\to A} (f\circ g\htpy \idfunc[B])
    \end{align*}
  of \define{left inverses}
  \indexdef{left!inverse}%
  \indexdef{inverse!left}%
  and \define{right inverses}
  \indexdef{right!inverse}%
  \indexdef{inverse!right}%
  to $f$, respectively.
  We call $f$ \define{left invertible}
  \indexdef{function!left invertible}%
  \indexdef{function!right invertible}%
  if $\linv(f)$ is inhabited, and similarly \define{right invertible}
  \indexdef{left!invertible function}%
  \indexdef{right!invertible function}%
  if $\rinv(f)$ is inhabited.
\end{defn}

\begin{lem}\label{thm:equiv-compose-equiv}
  If $f:A\to B$ has a quasi-inverse, then so do
  \begin{align*}
    (f\circ \blank) &: (C\to A) \to (C\to B)\\
    (\blank\circ f) &: (B\to C) \to (A\to C).
  \end{align*}
\end{lem}
\begin{proof}
  If $g$ is a quasi-inverse of $f$, then $(g\circ \blank)$ and $(\blank\circ g)$ are quasi-inverses of $(f\circ \blank)$ and $(\blank\circ f)$ respectively.
\end{proof}

\begin{lem}\label{lem:inv-hprop}
  If $f : A \to B$ has a quasi-inverse, then the types $\rinv(f)$ and $\linv(f)$ are contractible.
\end{lem}
\begin{proof}
  By function extensionality, we have
  \[\eqv{\linv(f)}{\sm{g:B\to A} (g\circ f = \idfunc[A])}.\]
  But this is the fiber of $(\blank\circ f)$ over $\idfunc[A]$, and so
  by \cref{thm:equiv-compose-equiv,thm:equiv-iso-adj,thm:contr-hae}, it is contractible.
  Similarly, $\rinv(f)$ is equivalent to the fiber of $(f\circ \blank)$ over $\idfunc[B]$ and hence contractible.
\end{proof}

Next we define the types which put together the other homotopy with the additional coherence datum.\index{coherence}%

\begin{defn}\label{defn:lcoh-rcoh}
For $f : A \to B$, a left inverse $(g,\eta) : \linv(f)$, and a right inverse $(g,\epsilon) : \rinv(f)$, we denote
\begin{align*}
\lcoh{f}{g}{\eta} & \defeq \sm{\epsilon : f\circ g \htpy \idfunc[B]} \prd{y:B} g(\epsilon y) = \eta (gy), \\
\rcoh{f}{g}{\epsilon} & \defeq \sm{\eta : g\circ f \htpy \idfunc[A]} \prd{x:A} f(\eta x) = \epsilon (fx).
\end{align*}
\end{defn}

\begin{lem}\label{lem:coh-hfib}
For any $f,g,\epsilon,\eta$, we have
\begin{align*}
\lcoh{f}{g}{\eta} & \eqvsym {\prd{y:B} \id[\hfib{g}{gy}]{(fgy,\eta(gy))}{(y,\refl{gy})}}, \\
\rcoh{f}{g}{\epsilon} & \eqvsym {\prd{x:A} \id[\hfib{f}{fx}]{(gfx,\epsilon(fx))}{(x,\refl{fx})}}.
\end{align*}
\end{lem}
\begin{proof}
Using \autoref{lem:hfib}.
\end{proof}

\begin{lem}\label{lem:coh-hprop}
  If $f$ is a half adjoint equivalence, then for any $(g,\epsilon) : \rinv(f)$, the type $\rcoh{f}{g}{\epsilon}$ is contractible.
\end{lem}
\begin{proof}
  By \autoref{lem:coh-hfib} and the fact that dependent function types preserve contractible spaces, it suffices to show that for each $x:A$, the type $\id[\hfib{f}{fx}]{(fgx,\epsilon(fx))}{(x,\refl{fx})}$ is contractible.
  But by \autoref{thm:contr-hae}, $\hfib{f}{fx}$ is contractible, and any path space of a contractible space is itself contractible.
\end{proof}

\begin{thm}\label{thm:hae-hprop}
  For any $f : A \to B$, the type $\ishae(f)$ is a mere proposition.
\end{thm}
\begin{proof}
  By \autoref{thm:contr-unit} it suffices to assume $f$ to be a half adjoint equivalence and show that $\ishae(f)$ is contractible.
  Now by associativity of $\Sigma$ (\autoref{ex:sigma-assoc}), the type $\ishae(f)$ is equivalent to
  \[\sm{u : \rinv(f)} \rcoh{f}{\proj{1}(u)}{\proj{2}(u)}.\]
  But by \cref{lem:inv-hprop,lem:coh-hprop} and the fact that $\Sigma$ preserves contractibility, the latter type is also contractible.
\end{proof}

Thus, we have shown that $\ishae(f)$ has all three desiderata for the type $\isequiv(f)$.
In the next two sections we consider a couple of other possibilities.

\index{equivalence!half adjoint|)}%
\index{half adjoint equivalence|)}%
\index{adjoint!equivalence!of types, half|)}%

\section{Bi-invertible maps}
\label{sec:biinv}

\index{function!bi-invertible|(defstyle}%
\index{bi-invertible function|(defstyle}%
\index{equivalence!as bi-invertible function|(defstyle}%

Using the language introduced in \autoref{sec:hae}, we can restate the definition proposed in \autoref{sec:basics-equivalences} as follows.

\begin{defn}\label{defn:biinv}
  We say $f:A\to B$ is \define{bi-invertible}
  if it has both a left inverse and a right inverse:
  \[ \biinv (f) \defeq \linv(f) \times \rinv(f). \]
\end{defn}

In \autoref{sec:basics-equivalences} we proved that $\qinv(f)\to\biinv(f)$ and $\biinv(f)\to\qinv(f)$.
What remains is the following.

\begin{thm}\label{thm:isprop-biinv}
  For any $f:A\to B$, the type $\biinv(f)$ is a mere proposition.
\end{thm}
\begin{proof}
  We may suppose $f$ to be bi-invertible and show that $\biinv(f)$ is contractible.
  But since $\biinv(f)\to\qinv(f)$, by \autoref{lem:inv-hprop} in this case both $\linv(f)$ and $\rinv(f)$ are contractible, and the product of contractible types is contractible.
\end{proof}

Note that this also fits the proposal made at the beginning of \autoref{sec:hae}: we combine $g$ and $\eta$ into a contractible type and add an additional datum which combines with $\epsilon$ into a contractible type.
The difference is that instead of adding a \emph{higher} datum (a 2-dimensional path) to combine with $\epsilon$, we add a \emph{lower} one (a right inverse that is separate from the left inverse).

\begin{cor}\label{thm:equiv-biinv-isequiv}
  For any $f:A\to B$ we have $\eqv{\biinv(f)}{\ishae(f)}$.
\end{cor}
\begin{proof}
  We have $\biinv(f) \to \qinv(f) \to \ishae(f)$ and $\ishae(f) \to \qinv(f) \to \biinv(f)$.
  Since both $\iscontr(f)$ and $\biinv(f)$ are mere propositions, the equivalence follows from \autoref{lem:equiv-iff-hprop}.
\end{proof}

\index{function!bi-invertible|)}%
\index{bi-invertible function|)}%
\index{equivalence!as bi-invertible function|)}%

\section{Contractible fibers}
\label{sec:contrf}

\index{function!contractible|(defstyle}%
\index{contractible!function|(defstyle}%
\index{equivalence!as contractible function|(defstyle}%

Note that our proofs about $\ishae(f)$ and $\biinv(f)$ made essential use of the fact that the fibers of an equivalence are contractible.
In fact, it turns out that this property is itself a sufficient definition of equivalence.

\begin{defn}[Contractible maps] \label{defn:equivalence}
  A map $f:A\to B$ is \define{contractible}
  if for all $y:B$, the fiber $\hfib f y$ is contractible.
\end{defn}

Thus, the type $\iscontr(f)$ is defined to be
\begin{align}
  \iscontr(f) &\defeq \prd{y:B} \iscontr(\hfib f y)\label{eq:iscontrf}
\end{align}
Note that in \autoref{sec:contractibility} we defined what it means for a \emph{type} to be contractible.
Here we are defining what it means for a \emph{map} to be contractible.
Our terminology follows the general homotopy-theoretic practice of saying that a map has a certain property if all of its (homotopy) fibers have that property.
Thus, a type $A$ is contractible just when the map $A\to\unit$ is contractible.
From \autoref{cha:hlevels} onwards we will also call contractible maps and types \emph{$(-2)$-truncated}.

We have already shown in \autoref{thm:contr-hae} that $\ishae(f) \to \iscontr(f)$.
Conversely:

\begin{thm}\label{thm:lequiv-contr-hae}
For any $f:A\to B$ we have ${\iscontr(f)} \to {\ishae(f)}$.
\end{thm}
\begin{proof}
Let $P : \iscontr(f)$. We define an inverse mapping $g : B \to A$ by sending each $y : B$ to the center of contraction of the fiber at $y$:
\[ g(y) \defeq \proj{1}(\proj{1}(Py)) \]
We can thus define the homotopy $\epsilon$ by mapping $y$ to the witness that $g(y)$ indeed belongs to the fiber at $y$:
\[ \epsilon(y) \defeq \proj{2}(\proj{1}(P y)) \]
It remains to define $\eta$ and $\tau$. This of course amounts to giving an element of $\lcoh{f}{g}{\epsilon}$. By \autoref{lem:coh-hfib}, this is the same as giving for each $x:A$ a path from $(x,\refl{fx})$ to $(gfx,\epsilon(fx))$ in the fiber at $fx$. But this is easy: for any $x : A$, the type $\hfib{f}{fx}$ 
is contractible by assumption, hence such a path must exist. We can construct it explicitly as
\[\opp{\big(P(fx) \; (x,\refl{fx})\big)} \ct \big(P(fx) \; (gfx,\epsilon(fx))\big). \qedhere \]
\end{proof}

It is also easy to see:

\begin{lem}\label{thm:contr-hprop}
  For any $f$, the type $\iscontr(f)$ is a mere proposition.
\end{lem}
\begin{proof}
  By \autoref{thm:isprop-iscontr}, each type $\iscontr (\hfib f y)$ is a mere proposition.
  Thus, by \autoref{thm:isprop-forall}, so is~\eqref{eq:iscontrf}.
\end{proof}

\begin{thm}\label{thm:equiv-contr-hae}
  For any $f:A\to B$ we have $\eqv{\iscontr(f)}{\ishae(f)}$.
\end{thm}
\begin{proof}
  We have already established a logical equivalence ${\iscontr(f)} \Leftrightarrow {\ishae(f)}$, and both are mere propositions (\cref{thm:contr-hprop,thm:hae-hprop}).
  Thus, \autoref{lem:equiv-iff-hprop} applies.
\end{proof}

Usually, we prove that a function is an equivalence by exhibiting a quasi-inverse, but sometimes this definition is more convenient.
For instance, it implies that when proving a function to be an equivalence, we are free to assume that its codomain is inhabited.

\begin{cor}\label{thm:equiv-inhabcod}
  If $f:A\to B$ is such that $B\to \isequiv(f)$, then $f$ is an equivalence.
\end{cor}
\begin{proof}
  To show $f$ is an equivalence, it suffices to show that $\hfib f y$ is contractible for any $y:B$.
  But if $e:B\to \isequiv(f)$, then given any such $y$ we have $e(y):\isequiv(f)$, so that $f$ is an equivalence and hence $\hfib f y$ is contractible, as desired.
\end{proof}

\index{function!contractible|)}%
\index{contractible!function|)}%
\index{equivalence!as contractible function|)}%

\section{On the definition of equivalences}
\label{sec:concluding-remarks}

\indexdef{equivalence}
We have shown that all three definitions of equivalence satisfy the three desirable properties and are pairwise equivalent:
\[ \iscontr(f) \eqvsym \ishae(f) \eqvsym \biinv(f). \]
(There are yet more possible definitions of equivalence, but we will stop with these three.
See \autoref{ex:brck-qinv} and the exercises in this chapter for some more.)
Thus, we may choose any one of them as ``the'' definition of $\isequiv (f)$.
For definiteness, we choose to define
\[ \isequiv(f) \defeq \ishae(f).\]
\index{mathematics!formalized}%
This choice is advantageous for formalization, since $\ishae(f)$ contains the most directly useful data.
On the other hand, for other purposes, $\biinv(f)$ is often easier to deal with, since it contains no 2-dimensional paths and its two symmetrical halves can be treated independently.
However, for purposes of this book, the specific choice will make little difference.

In the rest of this chapter, we study some other properties and characterizations of equivalences.
\index{equivalence!properties of}%

\section{Surjections and embeddings}
\label{sec:mono-surj}

\index{set}
When $A$ and $B$ are sets and $f:A\to B$ is an equivalence, we also call it as \define{isomorphism}
\indexdef{isomorphism!of sets}%
or a \define{bijection}.
\indexdef{bijection}%
\indexsee{function!bijective}{bijection}%
(We avoid these words for types that are not sets, since in homotopy theory and higher category theory they often denote a stricter notion of ``sameness'' than homotopy equivalence.)
In set theory, a function is a bijection just when it is both injective and surjective.
The same is true in type theory, if we formulate these conditions appropriately.
For clarity, when dealing with types that are not sets, we will speak of \emph{embeddings} instead of injections.

\begin{defn}
  Let $f:A\to B$.
  \begin{enumerate}
  \item We say $f$ is \define{surjective}
    \indexsee{surjective!function}{function, surjective}%
    \indexdef{function!surjective}%
    (or a \define{surjection})
    \indexsee{surjection}{function, surjective}%
    if for every $b:B$ we have $\brck{\hfib f b}$.
  \item We say $f$ is a \define{embedding}
    \indexdef{function!embedding}%
    \indexsee{embedding}{function, embedding}%
    if for every $x,y:A$ the function $\apfunc f : (\id[A]xy) \to (\id[B]{f(x)}{f(y)})$ is an equivalence.
  \end{enumerate}
\end{defn}

In other words, $f$ is surjective if every fiber of $f$ is merely inhabited, or equivalently if for all $b:B$ there merely exists an $a:A$ such that $f(a)=b$.
In traditional logical notation, $f$ is surjective if $\fall{b:B}\exis{a:A} (f(a)=b)$.
This must be distinguished from the stronger assertion that $\prd{b:B}\sm{a:A} (f(a)=b)$; if this holds we say that $f$ is a \define{split surjection}.
\indexsee{split!surjection}{function, split surjective}%
\indexsee{surjection!split}{function, split surjective}%
\indexsee{surjective!function!split}{function, split surjective}%
\indexdef{function!split surjective}%

If $A$ and $B$ are sets, then by \autoref{lem:equiv-iff-hprop}, $f$ is an embedding just when
\begin{equation}
  \prd{x,y:A} (\id[B]{f(x)}{f(y)}) \to (\id[A]xy).\label{eq:injective}
\end{equation}
In this case we say that $f$ is \define{injective},
\indexsee{injective function}{function, injective}%
\indexdef{function!injective}%
or an \define{injection}.
\indexsee{injection}{function, injective}%
We avoid these word for types that are not sets, because they might be interpreted as~\eqref{eq:injective}, which is an ill-behaved notion for non-sets.
It is also true that any function between sets is surjective if and only if it is an \emph{epimorphism} in a suitable sense, but this also fails for more general types, and surjectivity is generally the more important notion.

\begin{thm}\label{thm:mono-surj-equiv}
  A function $f:A\to B$ is an equivalence if and only if it is both surjective and an embedding.
\end{thm}
\begin{proof}
  If $f$ is an equivalence, then each $\hfib f b$ is contractible, hence so is $\brck{\hfib f b}$, so $f$ is surjective.
  And we showed in \autoref{thm:paths-respects-equiv} that any equivalence is an embedding.

  Conversely, suppose $f$ is a surjective embedding.
  Let $b:B$; we show that $\sm{x:A}(f(x)=b)$ is contractible.
  Since $f$ is surjective, there merely exists an $a:A$ such that $f(a)=b$.
  Thus, the fiber of $f$ over $b$ is inhabited; it remains to show it is a mere proposition.
  For this, suppose given $x,y:A$ with $p:f(x)=b$ and $q:f(y)=b$.
  Then since $\apfunc f$ is an equivalence, there exists $r:x=y$ with $\apfunc f (r) = p \ct \opp q$.
  However, using the characterization of paths in $\Sigma$-types, the latter equality rearranges to $\trans{r}{p} = q$.
  Thus, together with $r$ it exhibits $(x,p) = (y,q)$ in the fiber of $f$ over $b$.
\end{proof}

\begin{cor}
  For any $f:A\to B$ we have
  \[ \isequiv(f) \eqvsym (\mathsf{isEmbedding}(f) \times \mathsf{isSurjective}(f)).\]
\end{cor}
\begin{proof}
  Being a surjection and an embedding are both mere propositions; now apply \autoref{lem:equiv-iff-hprop}.
\end{proof}

Of course, this cannot be used as a definition of ``equivalence'', since the definition of embeddings refers to equivalences.
However, this characterization can still be useful; see \autoref{sec:whitehead}.
We will generalize it in \autoref{cha:hlevels}.

\section{Closure properties of equivalences}
\label{sec:equiv-closures}
\label{sec:fiberwise-equivalences}
\index{equivalence!properties of}%

We have already seen in \autoref{thm:equiv-eqrel} that equivalences are closed under composition.
Furthermore, we have:

\begin{thm}[The 2-out-of-3 property]\label{thm:two-out-of-three}
  \index{2-out-of-3 property}%
  Suppose $f:A\to B$ and $g:B\to C$.
  If any two of $f$, $g$, and $g\circ f$ are equivalences, so is the third.
\end{thm}
\begin{proof}[Sketch of proof]
  If $g\circ f$ and $g$ are equivalences, then $\opp{(g\circ f)} \circ g$ is a quasi-inverse to $f$.
  Similarly, if $g\circ f$ and $f$ are equivalences, then $f\circ \opp{(g\circ f)}$ is a quasi-inverse to $g$.
\end{proof}

This is a standard closure condition on equivalences from homotopy theory.
Also well-known is that they are closed under retracts, in the following sense.

\index{retract!of a function|(defstyle}%

\begin{defn}\label{defn:retract}
A function $g:A\to B$ is said to be a \define{retract}
of a function $f:X\to Y$ if there is a diagram
\begin{equation*}
  \xymatrix{
    {A} \ar[r]^{s} \ar[d]_{g}
    &
    {X} \ar[r]^{r} \ar[d]_{f}
    &
    {A} \ar[d]^{g}
    \\
    {B} \ar[r]_{s'}
    &
    {Y} \ar[r]_{r'}
    &
    {B}
  }
\end{equation*}
for which there are
\begin{enumerate}
\item a homotopy $R:r\circ s \htpy \idfunc[A]$.
\item a homotopy $R':r'\circ s' \htpy\idfunc[B]$.
\item a homotopy $L:f\circ s\htpy s'\circ g$.
\item a homotopy $K:g\circ r\htpy r'\circ f$.
\item for every $a:A$, a path $H(a)$ witnessing the commutativity of the square
\begin{equation*}
  \xymatrix@C=3pc{
    {g(r(s(a)))} \ar@{=}[r]^-{K(s(a))} \ar@{=}[d]_{\ap g{R(a)}}
    &
    {r'(f(s(a)))} \ar@{=}[d]^{\ap{r'}{L(a)}}
    \\
    {g(a)} \ar@{=}[r]_-{\opp{R'(g(a))}}
    &
    {r'(s'(g(a)))}
  }
\end{equation*}
\end{enumerate}
\end{defn}

Recall that in \autoref{sec:contractibility} we defined what it means for a type to be a retract of another.
This is a special case of the above definition where $B$ and $Y$ are $\unit$.
Conversely, just as with contractibility, retractions of maps induce retractions of their fibers.

\begin{lem}\label{lem:func_retract_to_fiber_retract}
If a function $g:A\to B$ is a retract of a function $f:X\to Y$, then $\hfib{g}b$ is a retract of $\hfib{f}{s'(b)}$
for every $b:B$, where $s':B\to Y$ is as in \autoref{defn:retract}.
\end{lem}

\begin{proof}
Suppose that $g:A\to B$ is a retract of $f:X\to Y$. Then for any $b:B$ we have the functions
\begin{align*}
\varphi_b &:\hfiber{g}b\to\hfib{f}{s'(b)}, &
\varphi_b(a,p) & \defeq \pairr{s(a),s'(p)\ct L(a)},\\
\psi_b &:\hfib{f}{s'(b)}\to\hfib{g}b, &
\psi_b(x,q) &\defeq \pairr{r(x),R'(b)\ct r'(q)\ct K(x)}.
\end{align*}
Then we have $\psi_b(\varphi_b({a,p}))\equiv\pairr{r(s(a)),R'(b)\ct r'(s'(p)\ct L(a))\ct K(s(a))}$.
We claim $\psi_b$ is a retraction with section $\varphi_b$ for all $b:B$, which is to say that for all $(a,p):\hfib g b$ we have $\psi_b(\varphi_b({a,p}))= \pairr{a,p}$.
In other words, we want to show
\begin{equation*}
\prd{b:B}{a:A}{p:f(a)=b} \psi_b(\varphi_b({a,p}))= \pairr{a,p}.
\end{equation*}
By reordering the first two $\Pi$s and applying a version of \autoref{thm:omit-contr}, this is equivalent to
\begin{equation*}
\prd{a:A}\psi_{f(a)}(\varphi_{f(a)}({a,\refl{g(a)}}))=\pairr{a,\refl{g(a)}}.
\end{equation*}
By assumption, we have $R(a):r(s(a))= a$. So it is left to show that there is a path
\begin{equation*}
\trans{R(a)}{R'(g(a))\ct r'(L(a))\ct K(s(a))} = \refl{g(a)}.
\end{equation*}
But this transportation computes as $R'(g(a))\ct r'(L(a))\ct K(s(a))\ct g(R(a))$, so the required path is given by $H(a)$.
\end{proof}

\begin{thm}\label{thm:retract-equiv}
  If $g$ is a retract of an equivalence $f$, then $g$ is also an equivalence.
\end{thm}
\begin{proof}
  By \autoref{lem:func_retract_to_fiber_retract}, every fiber of $g$ is a retract of a fiber of $f$.
  Thus, by \autoref{thm:retract-contr}, if the latter are all contractible, so are the former.
\end{proof}

\index{retract!of a function|)}%

\index{fibration}%
\index{total!space}%
Finally, we show that fiberwise equivalences can be characterized in terms of equivalences of total spaces.
To explain the terminology, recall from \autoref{sec:fibrations} that a type family $P:A\to\type$ can be viewed as a fibration over $A$ with total space $\sm{x:A} P(x)$, the fibration being is the projection $\proj1:\sm{x:A} P(x) \to A$.
From this point of view, given two type families $P,Q:A\to\type$, we may refer to a function $f:\prd{x:A} (P(x)\to Q(x))$ as a \define{fiberwise map} or a \define{fiberwise transformation}.
\indexsee{transformation!fiberwise}{fiberwise transformation}%
\indexsee{function!fiberwise}{fiberwise transformation}%
\index{fiberwise!transformation|(defstyle}%
\indexsee{fiberwise!map}{fiberwise transformation}%
\indexsee{map!fiberwise}{fiberwise transformation}
Such a map induces a function on total spaces:

\begin{defn}\label{defn:total-map}
  Given type families $P,Q:A\to\type$ and a map $f:\prd{x:A} P(x)\to Q(x)$, we define
  \begin{equation*}
    \total f  \defeq \lam{w}\pairr{\proj{1}w,f(\proj{1}w,\proj{2}w)} : \sm{x:A}P(x)\to\sm{x:A}Q(x).
  \end{equation*}
\end{defn}

\begin{thm}\label{fibwise-fiber-total-fiber-equiv}
Suppose that $f$ is a fiberwise transformation between families $P$ and
$Q$ over a type $A$ and let $x:A$ and $v:Q(x)$. Then we have an equivalence
\begin{equation*}
\eqv{\hfib{\total{f}}{\pairr{x,v}}}{\hfib{f(x)}{v}}.
\end{equation*}
\end{thm}
\begin{proof}
  We calculate:
\begin{align}
  \hfib{\total{f}}{\pairr{x,v}} 
  & \jdeq \sm{w:\sm{x:A}P(x)}\pairr{\proj{1}w,f(\proj{1}w,\proj{2}w)}=\pairr{x,v}
  \notag \\
  & \eqv{}{} \sm{a:A}{u:P(a)}\pairr{a,f(a,u)}=\pairr{x,v}
  \tag{by~\eqref{eq:sigma-lump}} \\
  & \eqv{}{} \sm{a:A}{u:P(a)}{p:a=x}\trans{p}{f(a,u)}=v
  \tag{by \autoref{thm:path-sigma}} \\
  & \eqv{}{} \sm{a:A}{p:a=x}{u:P(a)}\trans{p}{f(a,u)}=v
  \notag \\
  & \eqv{}{} \sm{u:P(x)}f(x,u)=v
  \tag{$*$}\label{eq:uses-sum-over-paths} \\
  & \jdeq \hfib{f(x)}{v}. \notag
\end{align}
The equivalence~\eqref{eq:uses-sum-over-paths} follows from \autoref{thm:omit-contr,thm:contr-paths,ex:sigma-assoc}.
\end{proof}

We say that a fiberwise transformation $f:\prd{x:A} P(x)\to Q(x)$ is a \define{fiberwise equivalence}%
\indexdef{fiberwise!equivalence}%
\indexdef{equivalence!fiberwise}
if each $f(x):P(x) \to Q(x)$ is an equivalence.

\begin{thm}\label{thm:total-fiber-equiv}
Suppose that $f$ is a fiberwise transformation between families
$P$ and $Q$ over a type $A$.
Then $f$ is a fiberwise equivalence if and only if $\total{f}$ is an equivalence.
\end{thm}

\begin{proof}
Let $f$, $P$, $Q$ and $A$ be as in the statement of the theorem.
By \autoref{fibwise-fiber-total-fiber-equiv} it follows for all
$x:A$ and $v:Q(x)$ that
$\hfib{\total{f}}{\pairr{x,v}}$ is contractible if and only if
$\hfib{f(x)}{v}$ is contractible.
Thus, $\hfib{\total{f}}{w}$ is contractible for all $w:\sm{x:A}Q(x)$ if and only if $\hfib{f(x)}{v}$ is contractible for all $x:A$ and $v:Q(x)$.
\end{proof}

\index{fiberwise!transformation|)}%

\section{The object classifier}
\label{sec:object-classification}

In type theory we have a basic notion of \emph{family of types}, namely a function $B:A\to\type$.
We have seen that such families behave somewhat like \emph{fibrations} in homotopy theory, with the fibration being the projection $\proj1:\sm{a:A} B(a) \to A$.
A basic fact in homotopy theory is that every map is equivalent to a fibration.
With univalence at our disposal, we can prove the same thing in type theory.

\begin{lem}\label{thm:fiber-of-a-fibration}
  For any type family $B:A\to\type$, the fiber of $\proj1:\sm{x:A} B(x) \to A$ over $a:A$ is equivalent to $B(a)$:
  \[ \eqv{\hfib{\proj1}{a}}{B(a)} \]
\end{lem}
\begin{proof}
  We have
  \begin{align*}
    \hfib{\proj1}{a} &\defeq \sm{u:\sm{x:A} B(x)} \proj1(u)=a\\
    &\eqvsym \sm{x:A}{b:B(x)} (x=a)\\
    &\eqvsym \sm{x:A}{p:x=a} B(x)\\
    &\eqvsym B(a)
  \end{align*}
  using the left universal property of identity types.
\end{proof}

\begin{lem}\label{thm:total-space-of-the-fibers}
  For any function $f:A\to B$, we have $\eqv{A}{\sm{b:B}\hfib{f}{b}}$.
\end{lem}
\begin{proof}
  We have
  \begin{align*}
    \sm{b:B}\hfib{f}{b} &\defeq \sm{b:B}{a:A} (f(a)=b)\\
    &\eqvsym \sm{a:A}{b:B} (f(a)=b)\\
    &\eqvsym A
  \end{align*}
  using the fact that $\sm{b:B} (f(a)=b)$ is contractible.
\end{proof}

\begin{thm}\label{thm:nobject-classifier-appetizer}
For any type $B$ there is an equivalence
\begin{equation*}
\chi:\Parens{\sm{A:\type} (A\to B)}\eqvsym (B\to\type).
\end{equation*}
\end{thm}
\begin{proof}
We have to construct quasi-inverses
\begin{align*}
\chi & : \Parens{\sm{A:\type} (A\to B)}\to B\to\type\\
\psi & : (B\to\type)\to\Parens{\sm{A:\type} (A\to B)}.
\end{align*}
We define $\chi$ by $\chi((A,f),b)\defeq\hfiber{f}b$, and $\psi$ by $\psi(P)\defeq\Pairr{(\sm{b:B} P(b)),\proj1}$.
Now we have to verify that $\chi\circ\psi\htpy\idfunc{}$ and that $\psi\circ\chi \htpy\idfunc{}$.
\begin{enumerate}
\item Let $P:B\to\type$.
  By \autoref{thm:fiber-of-a-fibration},
$\hfiber{\proj1}{b}\eqvsym P(b)$ for any $b:B$, so it follows immediately
that $P\htpy\chi(\psi(P))$.
\item Let $f:A\to B$ be a function. We have to find a path
\begin{equation*}
\Pairr{\tsm{b:B} \hfiber{f}b,\,\proj1}=\pairr{A,f}.
\end{equation*}
First note that by \autoref{thm:total-space-of-the-fibers}, we have
$e:\sm{b:B} \hfiber{f}b\eqvsym A$ with $e(b,a,p)\defeq a$ and $e^{-1}(a)
\defeq(f(a),a,\refl{f(a)})$.
By \autoref{thm:path-sigma}, it remains to show $\trans{(\ua(e))}{\proj1} = f$.
But by the computation rule for univalence and~\eqref{eq:transport-arrow}, we have $\trans{(\ua(e))}{\proj1} = \proj1\circ e^{-1}$, and the definition of $e^{-1}$ immediately yields $\proj1 \circ e^{-1} \jdeq f$.\qedhere
\end{enumerate}
\end{proof}

\noindent
\indexdef{object!classifier}%
\indexdef{classifier!object}%
\index{.infinity1-topos@$(\infty,1)$-topos}%
In particular, this implies that we have an \emph{object classifier} in the sense of higher topos theory.
Recall from \autoref{def:pointedtype} that $\pointed\type$ denotes the type $\sm{A:\type} A$ of pointed types.

\begin{thm}\label{thm:object-classifier}
Let $f:A\to B$ be a function. Then the diagram
\begin{equation*}
  \vcenter{\xymatrix{
      A\ar[r]^-{\vartheta_f} \ar[d]_{f} &
      \pointed{\type}\ar[d]^{\proj1}\\
      B\ar[r]_{\chi_f} &
      \type
      }}
\end{equation*}
is a pullback\index{pullback} square (see \autoref{ex:pullback}).
Here the function $\vartheta_f$ is defined by
\begin{equation*}
 \lam{a} \pairr{\hfiber{f}{f(a)},\pairr{a,\refl{f(a)}}}.
\end{equation*}
\end{thm}
\begin{proof}
Note that we have the equivalences
\begin{align*}
A & \eqvsym \sm{b:B} \hfiber{f}b\\
& \eqvsym \sm{b:B}{X:\type}{p:\hfiber{f}b= X} X\\
& \eqvsym \sm{b:B}{X:\type}{x:X} \hfiber{f}b= X\\
& \eqvsym \sm{b:B}{Y:\pointed{\type}} \hfiber{f}b = \proj1 Y\\
& \jdeq B\times_{\type}\pointed{\type}.
\end{align*}
which gives us a composite equivalence $e:A\eqvsym B\times_\type\pointed{\type}$. 
We may display the action of this composite equivalence step by step by
\begin{align*}
a & \mapsto \pairr{f(a),\; \pairr{a,\refl{f(a)}}}\\
& \mapsto \pairr{f(a), \; \hfiber{f}{f(a)}, \; \refl{\hfiber{f}{f(a)}}, \; \pairr{a,\refl{f(a)}}}\\
& \mapsto \pairr{f(a), \; \hfiber{f}{f(a)}, \; \pairr{a,\refl{f(a)}}, \; \refl{\hfiber{f}{f(a)}}}.
\end{align*}
Therefore, we get homotopies $f\htpy\proj1\circ e$ and $\vartheta_f\htpy \proj2\circ e$. 
\end{proof}

\section{Univalence implies function extensionality}
\label{sec:univalence-implies-funext}

\index{function extensionality!proof from univalence}%
In the last section of this chapter we include a proof that the univalence axiom implies function
extensionality. Thus, in this section we work \emph{without} the function extensionality axiom.
The proof consists of two steps. First we show
in \autoref{uatowfe} that the univalence
axiom implies a weak form of function extensionality, defined in \autoref{weakfunext} below. The
principle of weak function extensionality in turn implies the usual function extensionality,
and it does so without the univalence axiom (\autoref{wfetofe}).

\index{univalence axiom}%
Let $\type$ be a universe; we will explicitly indicate where we assume that it is univalent.

\begin{defn}\label{weakfunext}
The \define{weak function extensionality principle}
\indexdef{function extensionality!weak}%
asserts that there is a function
\begin{equation*}
\Parens{\prd{x:A}\iscontr(P(x))} \to\iscontr\Parens{\prd{x:A}P(x)}
\end{equation*}
for any family $P:A\to\type$ of types over any type $A$.
\end{defn}

The following lemma is easy to prove using function extensionality; the point here is that it also follows from univalence without assuming function extensionality separately.

\begin{lem} \label{UA-eqv-hom-eqv}
Assuming $\type$ is univalent, for any $A,B,X:\type$ and any $e:\eqv{A}{B}$, there is an equivalence
\begin{equation*}
\eqv{(X\to A)}{(X\to B)}
\end{equation*}
of which the underlying map is given by post-composition with the underlying function of $e$.
\end{lem}

\begin{proof}
  As in the proof of \autoref{lem:qinv-autohtpy}, we may assume that $e = \idtoeqv(p)$ for some $p:A=B$.
  Then by path induction, we may assume $p$ is $\refl{A}$, so that $e = \idfunc[A]$.
  But in this case, post-composition with $e$ is the identity, hence an equivalence.
\end{proof}

\begin{cor}\label{contrfamtotalpostcompequiv}
Let $P:A\to\type$ be a family of contractible types, i.e.\ \narrowequation{\prd{x:A}\iscontr(P(x)).}
Then the projection $\proj{1}:(\sm{x:A}P(x))\to A$ is an equivalence. Assuming $\type$ is univalent, it follows immediately that precomposition with $\proj{1}$ gives an equivalence
\begin{equation*}
\alpha : \eqv{\Parens{A\to\sm{x:A}P(x)}}{(A\to A)}.
\end{equation*}
\end{cor}

\begin{proof}
  By \autoref{thm:fiber-of-a-fibration}, for $\proj{1}:\sm{x:A}P(X)\to A$ and $x:A$ we have an equivalence
  \begin{equation*}
    \eqv{\hfiber{\proj{1}}{x}}{P(x)}.
  \end{equation*}
  Therefore $\proj{1}$ is an equivalence whenever each $P(x)$ is contractible. The assertion is now a consequence of  \autoref{UA-eqv-hom-eqv}.
\end{proof}

In particular, the homotopy fiber of the above equivalence at $\idfunc[A]$ is contractible. Therefore, we can show that univalence implies weak function extensionality by showing that the dependent function type $\prd{x:A}P(x)$ is a retract of $\hfiber{\alpha}{\idfunc[A]}$.

\begin{thm}\label{uatowfe}
In a univalent universe $\type$, suppose that $P:A\to\type$ is a family of contractible types
and let $\alpha$ be the function of \autoref{contrfamtotalpostcompequiv}. 
Then $\prd{x:A}P(x)$ is a retract of $\hfiber{\alpha}{\idfunc[A]}$. As a consequence, $\prd{x:A}P(x)$ is contractible. In other words, the univalence axiom implies the weak function extensionality principle.
\end{thm}

\begin{proof}
Define the functions
\begin{align*}
  \varphi &: \tprd{x:A}P(x)\to\hfiber{\alpha}{\idfunc[A]},\\
  \varphi(f) &\defeq (\lam{x} (x,f(x)),\refl{\idfunc[A]}),
\intertext{and}
  \psi &: \hfiber{\alpha}{\idfunc[A]}\to \tprd{x:A}P(x), \\
  \psi(g,p) &\defeq \lam{x} \trans p{\proj{2} (g(x))}.
\end{align*}
Then $\psi(\varphi(f))=\lam{x} f(x)$, which is $f$, by the uniqueness principle for dependent function types.
\end{proof}

We now show that weak function extensionality implies the usual function extensionality.
Recall from~\eqref{eq:happly} the function $\happly (f,g) : (f = g)\to(f\htpy g)$ which
converts equality of functions to homotopy. In the proof that follows, the univalence
axiom is not used.

\begin{thm}\label{wfetofe}
  \index{function extensionality}%
Weak function extensionality implies the function extensionality \autoref{axiom:funext}.
\end{thm}

\begin{proof}
We want to show that
\begin{equation*}
\prd{A:\type}{P:A\to\type}{f,g:\prd{x:A}P(x)}\isequiv(\happly (f,g)).
\end{equation*}
Since a fiberwise map induces an equivalence on total spaces if and only if it is fiberwise an equivalence by \autoref{thm:total-fiber-equiv}, it suffices to show that the function of type
\begin{equation*}
\Parens{\sm{g:\prd{x:A}P(x)}(f= g)} \to \sm{g:\prd{x:A}P(x)}(f\htpy g)
\end{equation*}
induced by $\lam{g:\prd{x:A}P(x)} \happly (f,g)$ is an equivalence.
Since the type on the left is contractible by \autoref{thm:contr-paths}, it suffices to show that the type on the right:
\begin{equation}\label{eq:uatofesp}
\sm{g:\prd{x:A}P(x)}\prd{x:A}f(x)= g(x)
\end{equation}
is contractible.
Now \autoref{thm:ttac} says that this is equivalent to
\begin{equation}\label{eq:uatofeps}
\prd{x:A}\sm{u:P(x)}f(x)= u.
\end{equation}
The proof of \autoref{thm:ttac} uses function extensionality, but only for one of the composites.
Thus, without assuming function extensionality, we can conclude that~\eqref{eq:uatofesp} is a retract\index{retract!of a type} of~\eqref{eq:uatofeps}.
And~\eqref{eq:uatofeps} is a product of contractible types, which is contractible by the weak function extensionality principle; hence~\eqref{eq:uatofesp} is also contractible.
\end{proof}

\sectionNotes

The fact that the space of continuous maps equipped with quasi-inverses has the wrong homotopy type to be the ``space of homotopy equivalences'' is well-known in algebraic topology.
In that context, the ``space of homotopy equivalences'' $(\eqv AB)$ is usually defined simply as the subspace of the function space $(A\to B)$ consisting of the functions that are homotopy equivalences.
In type theory, this would correspond most closely to $\sm{f:A\to B} \brck{\qinv(f)}$; see \autoref{ex:brck-qinv}.

The first definition of equivalence given in homotopy type theory was the one that we have called $\iscontr(f)$, which was due to Voevodsky.
The possibility of the other definitions was subsequently observed by various people.
The basic theorems about adjoint equivalences\index{adjoint!equivalence} such as \autoref{lem:coh-equiv,thm:equiv-iso-adj} are adaptations of standard facts in higher category theory and homotopy theory.
Using bi-invertibility as a definition of equivalences was suggested by Andr\'e Joyal.

The properties of equivalences discussed in \autoref{sec:mono-surj,sec:equiv-closures} are well-known in homotopy theory.
Most of them were first proven in type theory by Voevodsky.

The fact that every function is equivalent to a fibration is a standard fact in homotopy theory.
The notion of object classifier
\index{object!classifier}%
\index{classifier!object}%
in $(\infty,1)$-category
\index{.infinity1-category@$(\infty,1)$-category}%
theory (the categorical analogue of \autoref{thm:nobject-classifier-appetizer}) is due to Rezk (see~\cite{Rezk05,lurie:higher-topoi}).

Finally, the fact that univalence implies function extensionality (\autoref{sec:univalence-implies-funext}) is due to Voevodsky.
Our proof is a simplification of his.

\sectionExercises

\begin{ex}
  Consider the type of ``two-sided adjoint equivalence\index{adjoint!equivalence} data'' for $f:A\to B$,
  \begin{narrowmultline*}
    \sm{g:B\to A}{\eta: g \circ f \htpy \idfunc[A]}{\epsilon:f \circ g \htpy \idfunc[B]}
    \narrowbreak
    \Parens{\prd{x:A} \map{f}{\eta x} = \epsilon(fx)} \times
    \Parens{\prd{y:B} \map{g}{\epsilon y} = \eta(gy) }.
  \end{narrowmultline*}
  By \autoref{lem:coh-equiv}, we know that if $f$ is an equivalence, then this type is inhabited.
  Give a characterization of this type analogous to \autoref{lem:qinv-autohtpy}.

  Can you give an example showing that this type is not generally a mere proposition?
  (This will be easier after \autoref{cha:hits}.)
\end{ex}

\begin{ex}
  Show that for any $f:A\to B$, the following type also satisfies the three desiderata of $\isequiv(f)$:
  \begin{equation*}
    \sm{R:A\to B\to \type}
    \Parens{\prd{a:A} \iscontr\Parens{\sm{b:B} R(a,b)}} \times
    \Parens{\prd{b:B} \iscontr\Parens{\sm{a:A} R(a,b)}}.
  \end{equation*}
\end{ex}

\begin{ex} \label{ex:qinv-autohtpy-no-univalence}
  Reformulate the proof of \autoref{lem:qinv-autohtpy} without using univalence.
\end{ex}

\begin{ex}[The unstable octahedral axiom]\label{ex:unstable-octahedron}
  \index{axiom!unstable octahedral}%
  \index{octahedral axiom, unstable}%
  Suppose $f:A\to B$ and $g:B\to C$ and $b:B$.
  \begin{enumerate}
  \item Show that there is a natural map $\hfib{g\circ f}{g(b)} \to \hfib{g}{g(b)}$ whose fiber over $(b,\refl{g(b)})$ is equivalent to $\hfib f b$.
  \item Show that $\eqv{\hfib{g\circ f}{g(b)}}{\sm{w:\hfib{g}{g(b)}} \hfib f {\proj1 w}}$.
  \end{enumerate}
\end{ex}

\begin{ex}
  \index{2-out-of-6 property}%
  Prove that equivalences satisfy the \emph{2-out-of-6 property}: given $f:A\to B$ and $g:B\to C$ and $h:C\to D$, if $g\circ f$ and $h\circ g$ are equivalences, so are $f$, $g$, $h$, and $h\circ g\circ f$.
\end{ex}


\chapter{Induction}
\label{cha:induction}

In \autoref{cha:typetheory}, we introduced many ways to form new types from old ones.
Except for (dependent) function types and universes, all these rules are special cases of the general notion of \emph{inductive definition}.\index{definition!inductive}\index{inductive!definition}
In this chapter we study inductive definitions more generally.

\section{Introduction to inductive types}
\label{sec:bool-nat}

\index{type!inductive|(}%
\index{generation!of a type, inductive|(}
An \emph{inductive type} $X$ can be intuitively understood as a type ``freely generated'' by a certain finite collection of \emph{constructors}, each of which is a function (of some number of arguments) with codomain $X$.
This includes functions of zero arguments, which are simply elements of $X$.

When describing a particular inductive type, we list the constructors with bullets.
For instance, the type \bool from \autoref{sec:type-booleans} is inductively generated by the following constructors:
\begin{itemize}
\item $\bfalse:\bool$
\item $\btrue:\bool$
\end{itemize}
Similarly, $\unit$ is inductively generated by the constructor:
\begin{itemize}
\item $\ttt:\unit$
\end{itemize}
while $\emptyt$ is inductively generated by no constructors at all.
An example where the constructor functions take arguments is the coproduct $A+B$, which is generated by the two constructors
\begin{itemize}
\item $\inl:A\to A+B$
\item $\inr:B\to A+B$.
\end{itemize}
And an example with a constructor taking multiple arguments is the cartesian product $A\times B$, which is generated by one constructor
\begin{itemize}
\item $\pairr{\blank, \blank} : A\to B \to A\times B$.
\end{itemize}
Crucially, we also allow constructors of inductive types that take arguments from the inductive type being defined.
For instance, the type $\nat$ of natural numbers has constructors
\begin{itemize}
\item $0:\nat$
\item $\suc:\nat\to\nat$.
\end{itemize}
\symlabel{lst}
\index{type!of lists}%
\indexsee{list}{type of lists}%
Another useful example is the type $\lst A$ of finite lists of elements of some type $A$, which has constructors
\begin{itemize}
\item $\nil:\lst A$
\item $\cons:A\to \lst A \to \lst A$.
\end{itemize}

\index{free!generation of an inductive type}%
Intuitively, we should understand an inductive type as being \emph{freely generated by} its constructors.
That is, the elements of an inductive type are exactly what can be obtained by starting from nothing and applying the constructors repeatedly.
(We will see in \autoref{sec:identity-systems,cha:hits} that this conception has to be modified slightly for more general kinds of inductive definitions, but for now it is sufficient.)
For instance, in the case of \bool, we should expect that the only elements are $\bfalse$ and $\btrue$.
Similarly, in the case of \nat, we should expect that every element is either $0$ or obtained by applying $\suc$ to some ``previously constructed'' natural number.

\index{induction principle}%
Rather than assert properties such as this directly, however, we express them by means of an \emph{induction principle}, also called a \emph{(dependent) elimination rule}.
We have seen these principles already in \autoref{cha:typetheory}.
For instance, the induction principle for \bool is:
\begin{itemize}
\item When proving a statement $E : \bool \to \type$ about \emph{all} inhabitants of \bool, it suffices to prove it for $\bfalse$ and $\btrue$, i.e., to give proofs $e_0 : E(\bfalse)$ and $e_1 : E(\btrue)$.
\end{itemize}

Furthermore, the resulting proof $\ind\bool(E,e_0,e_1): \prd{b : \bool}E(b)$ behaves as expected when applied to the constructors $\bfalse$ and $\btrue$; this principle is expressed by the \emph{computation rules}\index{computation rule!for type of booleans}:
\begin{itemize}
\item We have $\ind\bool(E,e_0,e_1,\bfalse) \jdeq e_0$.
\item We have $\ind\bool(E,e_0,e_1,\btrue) \jdeq e_1$.
\end{itemize}

\index{case analysis}%
Thus, the induction principle for the type $\bool$ of booleans allow us to reason by \emph{case analysis}.
Since neither of the two constructors takes any arguments, this is all we need for booleans.

\index{natural numbers}%
For natural numbers, however, case analysis is generally not sufficient: in the case corresponding to the inductive step $\suc(n)$, we also want to presume that the statement being proven has already been shown for $n$.
This gives us the following induction principle:
\begin{itemize}
\item When proving a statement $E : \nat \to \type$ about \emph{all} natural numbers, it suffices to prove it for $0$ and for $\suc(n)$, assuming it holds
for $n$, i.e., we construct $e_z : E(0)$ and $e_s : \prd{n : \nat}{y : E(n)} E(\suc(n))$.
\end{itemize}
The variable
\index{variable}%
$y$ represents our inductive hypothesis.  As in the case of booleans, we also have the associated computation rules for the function $\ind\nat(E,e_z,e_s) : \prd{x:\nat} E(x)$:
\index{computation rule!for natural numbers}%
\begin{itemize}
\item $\ind\nat(E,e_z,e_s,0) \jdeq e_z$.
\item $\ind\nat(E,e_z,e_s,\suc(n)) \jdeq e_s(n,\ind\nat(E,e_z,e_s,n))$ for any $n : \nat$.
\end{itemize}
The dependent function $\ind\nat(E,e_z,e_s)$ can thus be understood as being defined recursively on the argument $x : \nat$, via the functions $e_z$ and $e_s$ which we call the \define{recurrences}\indexdef{recurrence}.
When $x$ is zero,\index{zero} the function simply returns $e_z$.
When $x$ is the successor\index{successor} of another natural number $n$, the result is obtained by taking the recurrence $e_s$ and substituting the specific predecessor\index{predecessor} $n$ and the recursive call value $\ind\nat(E,e_z,e_s,n)$.

The induction principles for all the examples mentioned above share this family resemblance.
In \autoref{sec:strictly-positive} we will discuss a general notion of ``inductive definition'' and how to derive an appropriate \emph{induction principle} for it, but first we investigate various commonalities between inductive definitions.

\index{recursion principle}%
For instance, we have remarked in every case in \autoref{cha:typetheory} that from the induction principle we can derive a \emph{recursion principle} in which the codomain is a simple type (rather than a family).
Both induction and recursion principles may seem odd, since they yield only the \emph{existence} of a function without seeming to characterize it uniquely.
However, in fact the induction principle is strong enough also to prove its own \emph{uniqueness principle}\index{uniqueness!principle, propositional!for functions on N@for functions on $\nat$}, as in the following theorem.

\begin{thm}\label{thm:nat-uniq}
Let $f,g : \prd{x:\nat} E(x)$ be two functions which satisfy the recurrences 
\begin{equation*}
  e_z : E(0)
  \qquad\text{and}\qquad
  e_s : \prd{n : \nat}{y : E(n)} E(\suc(n))
\end{equation*}
up to propositional equality, i.e., such that
\begin{equation*}
  \id{f(0)}{e_z}
  \qquad\text{and}\qquad
  \id{g(0)}{e_z}
\end{equation*}
as well as
\begin{gather*}
  \prd{n : \nat} \id{f(\suc(n))}{e_s(n, f(n))},\\
  \prd{n : \nat} \id{g(\suc(n))}{e_s(n, g(n))}.
\end{gather*}
Then $f$ and $g$ are equal.
\end{thm}

\begin{proof}
We use induction on the type family $D(x) \defeq \id{f(x)}{g(x)}$. For the base case, we have \[f(0) = e_z = g(0)\]
For the inductive case, assume $n : \nat$ such that $f(n) = g(n)$. Then
\[ f(\suc(n)) = e_s(n, f(n)) = e_s(n, g(n)) = g(\suc(n)) \]
The first and last equality follow from the assumptions on $f$ and $g$. The middle equality follows from the inductive hypothesis and the fact that application preserves equality. This gives us pointwise equality between $f$ and $g$; invoking function extensionality finishes the proof.
\end{proof}

Note that the uniqueness principle applies even to functions that only satisfy the recurrences \emph{up to propositional equality}, i.e.\ a path.
Of course, the particular function obtained from the induction principle satisfies these recurrences judgmentally\index{judgmental equality}; we will return to this point in \autoref{sec:htpy-inductive}.
On the other hand, the theorem itself only asserts a propositional equality between functions (see also \autoref{ex:same-recurrence-not-defeq}).
From a homotopical viewpoint it is natural to ask whether this path is \emph{coherent}, i.e.\ whether the equality $f=g$ is unique up to higher paths; in \autoref{sec:initial-alg} we will see that this is in fact the case.

Of course, similar uniqueness theorems for functions can generally be formulated and shown for other inductive types as well.
In the next section, we show how this uniqueness property, together with univalence, implies that an inductive type such as the natural numbers is completely characterized by its introduction, elimination, and computation rules.

\index{type!inductive|)}%
\index{generation!of a type, inductive|)}%


\section{Uniqueness of inductive types}
\label{sec:appetizer-univalence}

\index{uniqueness!of identity types}%
We have defined ``the'' natural numbers to be a particular type \nat with particular inductive generators $0$ and $\suc$.
However, by the general principle of inductive definitions in type theory described in the previous section, there is nothing preventing us from defining \emph{another} type in an identical way.
\index{natural numbers!isomorphic definition of}
That is, suppose we let $\natp$ be the inductive type generated by the constructors
\begin{itemize}
\item $\zerop:\natp$
\item $\sucp:\natp\to\natp$.
\end{itemize}
Then $\natp$ will have identical-looking induction and recursion principles to $\nat$.
When proving a statement $E : \natp \to \type$ for all of these ``new'' natural numbers, it suffices to give the proofs $e_z : E(\zerop)$ and \narrowequation{e_s : \prd{n : \natp}{x : E(n)} E(\sucp(n)).}
And the function $\rec\natp(E,e_z,e_s) : \prd{n:\natp} E(n)$ has the following computation rules:
\begin{itemize}
\item $\rec\natp(E,e_z,e_s,\zerop) \jdeq e_z$,
\item $\rec\natp(E,e_z,e_s,\sucp(n)) \jdeq e_s(n,\rec\natp(E,e_z,e_s,n))$ for any $n : \natp$.
\end{itemize}
But what is the relation between $\nat$ and $\natp$?

This is not just an academic question, since structures that ``look like'' the natural numbers can be found in many other places.
\index{type!of lists}%
For instance, we may identify natural numbers with lists over the type with one element (this is arguably the oldest appearance, found on walls of caves\index{caves, walls of}), with the non-negative integers, with subsets of the rationals and the reals, and so on.
And from a programming\index{programming} point of view, the ``unary'' representation of our natural numbers is very inefficient, so we might prefer sometimes to use a binary one instead.
We would like to be able to identify all of these versions of ``the natural numbers'' with each other, in order to transfer constructions and results from one to another.

Of course, if two versions of the natural numbers satisfy identical induction principles, then they have identical induced structure.
For instance, recall the example of the function $\dbl$ defined in \autoref{sec:inductive-types}. A similar function
for our new natural numbers is readily defined by duplication and adding primes:
\[ \dblp \defeq \rec\natp(\natp, \; \zerop, \;  \lamu{n:\natp}{m:\natp} \sucp(\sucp(m))). \]
Simple as this may seem, it has the obvious drawback of leading to a
proliferation of duplicates. Not only functions have to be
duplicated, but also all lemmas and their proofs. For example,
an easy result such as  $\prd{n : \nat} \dbl(\suc(n))=\suc(\suc(\dbl(n)))$, as well
as its proof by induction, also has to be ``primed''.

In traditional mathematics, one just proclaims that $\nat$ and $\natp$ are obviously ``the same'', and can be substituted for each other whenever the need arises.
This is usually unproblematic, but it sweeps a fair amount under the rug, widening the gap between informal mathematics and its precise description.
In homotopy type theory, we can do better.

First observe that we have the following definable maps:
\begin{itemize}
\item $f \defeq \rec\nat(\nat, \; \zerop, \;  \lamu{n:\nat} \sucp)
       : \nat \to\natp$,
\item $g \defeq \rec\natp(\natp, \; 0, \;  \lamu{n:\natp}, \suc)
       : \natp \to\nat$.
\end{itemize}
Since the composition of $g$ and $f$ satisfies the same recurrences as the identity function on $\nat$, \autoref{thm:nat-uniq} gives that $\prd{n : \nat} \id{g(f(n))}{n}$, and the ``primed'' version of the same theorem gives $\prd{n : \natp} \id{f(g(n))}{n}$.
Thus, $f$ and $g$ are quasi-inverses, so that $\eqv{\nat}{\natp}$.
We can now transfer functions on $\nat$ directly to functions on $\natp$ (and vice versa) along this equivalence, e.g.
\[ \dblp \defeq \lamu{n:\natp} f(\dbl(g(n))). \]
It is an easy exercise to show that this version of $\dblp$ is equal to the earlier one.

Of course, there is nothing surprising about this; such an isomorphism is exactly how a mathematician will envision ``identifying'' $\nat$ with $\natp$.
However, the mechanism of ``transfer'' across an isomorphism depends on the thing being transferred; it is not always as simple as pre- and post-composing a single function with $f$ and $g$.
Consider, for instance, a simple lemma such as
\[\prd{n : \natp} \dblp(\sucp(n))=\sucp(\sucp(\dblp(n))).\]
Inserting the correct $f$s and $g$s is only a little easier than re-proving it by induction on $n:\natp$ directly.
\index{isomorphism!transfer across}%

\index{univalence axiom}%
Here is where the univalence axiom steps in: since $\eqv{\nat}{\natp}$, we also have $\id[\type]{\nat}{\natp}$, i.e.\ $\nat$ and $\natp$ are
\emph{equal} as types.
Now the induction principle for identity guarantees that any construction or proof relating to $\nat$ can automatically be transferred to $\natp$ in the same way.
We simply consider the type of the function or theorem as a type-indexed family of types $P:\type\to\type$, with the given object being an element of $P(\nat)$, and transport along the path $\id \nat\natp$.
This involves considerably less overhead.

For simplicity, we have described this method in the case of two types \nat and \natp with \emph{identical}-looking definitions.
However, a more common situation in practice is when the definitions are not literally identical, but nevertheless one induction principle implies the other.
\index{type!unit}%
\index{natural numbers!encoded as list(unit)@encoded as $\lst\unit$}%
Consider, for instance, the type of lists from a one-element type, $\lst\unit$, which is generated by
\begin{itemize}
\item an element $\nil:\lst\unit$, and
\item a function $\cons:\unit \times \lst\unit \to\lst\unit$.
\end{itemize}
This is not identical to the definition of \nat, and it does not give rise to an identical induction principle.
The induction principle of $\lst\unit$ says that for any $E:\lst\unit\to\type$ together with recurrence data $e_\nil:E(\nil)$ and $e_\cons : \prd{u:\unit}{\ell:\lst\unit} E(\ell) \to E(\cons(u,\ell))$, there exists $f:\prd{\ell:\lst\unit} E(\ell)$ such that $f(\nil)\jdeq e_\nil$ and $f(\cons(u,\ell))\jdeq e_\cons(u,\ell,f(\ell))$.
(We will see how to derive the induction principle of an inductive definition in \autoref{sec:strictly-positive}.)

Now suppose we define $0'' \defeq \nil: \lst\unit$, and $\suc'':\lst\unit\to\lst\unit$ by $\suc''(\ell) \defeq \cons(\ttt,\ell)$.
Then for any $E:\lst\unit\to\type$ together with $e_0:E(0'')$ and $e_s:\prd{\ell:\lst\unit} E(\ell) \to E(\suc''(\ell))$, we can define
\begin{align*}
  e_\nil &\defeq e_0\\
  e_\cons(\ttt,\ell,x) &\defeq e_s(\ell,x).
\end{align*}
(In the definition of $e_\cons$ we use the induction principle of \unit to assume that $u$ is $\ttt$.)
Now we can apply the induction principle of $\lst\unit$, obtaining $f:\prd{\ell:\lst\unit} E(\ell)$ such that
\begin{gather*}
  f(0'') \jdeq f(\nil) \jdeq e_\nil \jdeq e_0\\
  f(\suc''(\ell)) \jdeq f(\cons(\ttt,\ell)) \jdeq e_\cons(\ttt,\ell,f(\ell)) \jdeq e_s(\ell,f(\ell)).
\end{gather*}
Thus, $\lst\unit$ satisfies the same induction principle as $\nat$, and hence (by the same arguments above) is equal to it.

Finally, these conclusions are not confined to the natural numbers: they apply to any inductive type.
If we have an inductively defined type $W$, say, and some other type $W'$ which satisfies the same induction principle as $W$, then it follows that $\eqv{W}{W'}$, and hence $W=W'$.
We use the derived recursion principles for $W$ and $W'$ to construct maps $W\to W'$ and $W'\to W$, respectively, and then the induction principles for each to prove that both composites are equal to identities.
For instance, in \autoref{cha:typetheory} we saw that the coproduct $A+B$ could also have been defined as $\sm{x:\bool} \rec{\bool}(\UU,A,B,x)$.
The latter type satisfies the same induction principle as the former; hence they are canonically equivalent.

This is, of course, very similar to the familiar fact in category theory that if two objects have the same \emph{universal property}, then they are equivalent.
In \autoref{sec:initial-alg} we will see that inductive types actually do have a universal property, so that this is a manifestation of that general principle.
\index{universal property}

\section{$\w$-types}
\label{sec:w-types}

Inductive types are very general, which is excellent for their usefulness and applicability, but makes them difficult to study as a whole.
Fortunately, they can all be formally reduced to a few special cases.
It is beyond the scope of this book to discuss this reduction --- which is anyway irrelevant to the mathematician using type theory in practice --- but we will take a little time to discuss the one of the basic special cases that we have not yet met.
These are Martin-L{\"o}f's \emph{$\w$-types}, also known as the types of \emph{well-founded trees}.
\index{tree, well-founded}%
$\w$-types are a generalization of such types as natural numbers, lists, and binary trees, which are sufficiently general to encapsulate the ``recursion'' aspect of \emph{any} inductive type.

A particular $\w$-type is specified by giving two parameters $A : \type$ and $B : A \to \type$, in which case the resulting $\w$-type is written $\wtype{a:A} B(a)$.
The type $A$ represents the type of \emph{labels} for $\wtype{a :A} B(a)$, which function as constructors (however, we reserve that word for the actual functions which arise in inductive definitions). For instance, when defining natural numbers as a $\w$-type,%
\index{natural numbers!encoded as a W-type@encoded as a $\w$-type}
the type $A$ would be the type $\bool$ inhabited by the two elements $\bfalse$ and $\btrue$, since there are precisely two ways to obtain a natural number --- either it will be zero or a successor of another natural number. 

The type family $B : A \to \type$ is used to record the arity\index{arity} of labels: a label $a : A$ will take a family of inductive arguments, indexed over $B(a)$. We can therefore think of the ``$B(a)$-many'' arguments of $a$. These arguments are represented by a function $f : B(a) \to \wtype{a :A} B(a)$, with the understanding that for any $b : B(a)$, $f(b)$ is the ``$b$-th'' argument to the label $a$. The $\w$-type $\wtype{a :A} B(a)$ can thus be thought of as the type of well-founded trees, where nodes are labeled by elements of $A$ and each node labeled by $a : A$ has $B(a)$-many branches.

In the case of natural numbers, the label $\bfalse$ has arity 0, since it constructs the constant zero\index{zero}; the label $\btrue$ has arity 1, since it constructs the successor\index{successor} of its argument. We can capture this by using simple elimination on $\bool$ to define a function $\rec\bool(\bbU,\emptyt,\unit)$ into a universe of types; this function returns the empty type $\emptyt$ for $\bfalse$ and the unit\index{type!unit} type $\unit$ for $\btrue$. We can thus define
\symlabel{natw}
\index{type!unit}%
\index{type!empty}%
\[ \natw \defeq \wtype{b:\bool} \rec\bool(\bbU,\emptyt,\unit) \]
where the superscript $\mathbf{w}$ serves to distinguish this version of natural numbers from the previously used one.
Similarly, we can define the type of lists\index{type!of lists} over $A$ as a $\w$-type with $\unit + A$ many labels: one nullary label for the empty list, plus one unary label for each $a : A$, corresponding to appending $a$ to the head of a list:
\[ \lst A \defeq \wtype{x: \unit + A} \rec{\unit + A}(\bbU, \; \emptyt, \; \lamu{a:A} \unit). \]
\index{W-type@$\w$-type}%
\indexsee{type!W-@$\w$-}{$\w$-type}%
In general, the $\w$-type $\wtype{x:A} B(x)$ specified by  $A : \type$ and $B : A \to \type$ is the inductive type generated by the following constructor:
\begin{itemize}
\item \label{defn:supp}
  $\supp : \prd{a:A} \Big(B(a) \to \wtype{x:A} B(x)\Big) \to \wtype{x:A} B(x)$.
\end{itemize}
The constructor $\supp$ (short for supremum\index{supremum!constructor of a W-type@constructor of a $\w$-type}) takes a label $a : A$ and a function $f : B(a) \to \wtype{x:A} B(x)$ representing the arguments to $a$, and constructs a new element of $\wtype{x:A} B(x)$. Using our previous encoding of natural numbers as $\w$-types, we can for instance define
\begin{equation*}
\zerow \defeq \supp(\bfalse, \; \lamu{x:\emptyt} \rec\emptyt(\natw,x)).
\end{equation*}
Put differently, we use the label $\bfalse$ to construct $\zerow$. Then, $\rec\bool(\bbU,\emptyt,\unit, \bfalse)$ evaluates to $\emptyt$, as it should since $\bfalse$ is a nullary label. Thus, we need to construct a function $f : \emptyt \to \natw$, which represents the (zero) arguments supplied to $\bfalse$. This is of course trivial, using simple elimination on $\emptyt$ as shown. Similarly, we can define
\begin{align*}
1^{\mathbf{w}} &\defeq \supp(\btrue, \; \lamu{x:\unit} 0^{\mathbf{w}}) \\
2^{\mathbf{w}} &\defeq \supp(\btrue, \; \lamu{x:\unit} 1^{\mathbf{w}})
\end{align*}
and so on.

\index{induction principle!for W-types@for $\w$-types}%
We have the following induction principle for $\w$-types:
\begin{itemize}
\item When proving a statement $E : \big(\wtype{x:A} B(x)\big) \to \type$ about \emph{all} elements of the $\w$-type $\wtype{x:A} B(x)$, it suffices to prove it for $\supp(a,f)$, assuming it holds for all $f(b)$ with $b : B(a)$. 
In other words, it suffices to give a proof 
\begin{equation*}
e : \prd{a:A}{f : B(a) \to \wtype{x:A} B(x)}{g : \prd{b : B(a)} E(f(b))} E(\supp(a,f))
\end{equation*}
\end{itemize}

\index{variable}%
The variable $g$ represents our inductive hypothesis, namely that all arguments of $a$ satisfy $E$. To state this, we quantify over all elements of type $B(a)$, since each $b : B(a)$ corresponds to one argument $f(b)$ of $a$.

How would we define the function $\dbl$ on natural numbers encoded as a $\w$-type? We would like to use the recursion principle of $\natw$ with a codomain of $\natw$ itself. We thus need to construct a suitable function
\[e : \prd{a : \bool}{f : B(a) \to \natw}{g : B(a) \to \natw} \natw\]
which will represent the recurrence for the $\dbl$ function; for simplicity we denote the type family $\rec\bool(\bbU,\emptyt,\unit)$ by $B$.

Clearly, $e$ will be a function taking $a : \bool$ as its first argument. The next step is to perform case analysis on $a$ and proceed based on whether it is $\bfalse$ or $\btrue$. This suggests the following form
\[ e \defeq \lamu{a:\bool} \rec\bool(C,e_0,e_1,a) \]
where \[C \defeq \prd{f : B(a) \to \natw}{g : B(a) \to \natw} \natw\]
If $a$ is $\bfalse$, the type $B(a)$ becomes $\emptyt$. Thus, given $f : \emptyt \to \natw$ and $g : \emptyt \to \natw$, we want to construct an element of $\natw$. Since the label $\bfalse$ represents $\emptyt$, it needs zero inductive arguments and the variables $f$ and $g$ are irrelevant. We return $\zerow$\index{zero} as a result:
\[ e_0 \defeq \lamu{f:\emptyt \to \natw}{g:\emptyt \to \natw} \zerow \]
Analogously, if $a$ is $\btrue$, the type $B(a)$ becomes $\unit$.
Since the label $\btrue$ represents the successor\index{successor} operator, it needs one inductive argument --- the predecessor\index{predecessor} --- which is represented by the variable $f : \unit \to \natw$.
The value of the recursive call on the predecessor is represented by the variable $g : \unit \to \natw$.
Thus, taking this value (namely $g(\ttt)$) and applying the successor operator twice thus yields the desired result:
\begin{equation*}
e_1 \defeq \; \lamu{f:\unit \to \natw}{g:\unit \to \natw}
  \supp(\btrue, (\lamu{x:\unit} \supp(\btrue, (\lamu{y : \unit} g(\ttt))))).
\end{equation*}
Putting this together, we thus have
\[ \dbl \defeq \rec\natw(\natw, e) \]
with $e$ as defined above.

\symlabel{defn:recursor-wtype}
The associated computation rule for the function $\rec{\wtype{x:A} B(x)}(E,e) : \prd{w : \wtype{x:A} B(x)} E(w)$ is as follows.
\index{computation rule!for W-types@for $\w$-types}%
\begin{itemize}
\item
  For any $a : A$ and $f : B(a) \to \wtype{x:A} B(x)$ we have 
  \begin{equation*}
    \rec{\wtype{x:A} B(x)}(E,e,\supp(a,f)) \jdeq
    e(a,f,\big(\lamu{b:B(a)} \rec{\wtype{x:A} B(x)}(E,f(b))\big)).
  \end{equation*}
\end{itemize}
In other words, the function $\rec{\wtype{x:A} B(x)}(E,e)$ satisfies the recurrence $e$.

By the above computation rule, the function $\dbl$ behaves as expected:
\begin{align*}
\dbl(\zerow) & \jdeq \rec\natw(\natw, e, \supp(\bfalse, \; \lamu{x:\emptyt} \rec\emptyt(\natw,x))) \\
& \jdeq e(\bfalse, \big(\lamu{x:\emptyt} \rec\emptyt(\natw,x)\big), 
   \big(\lamu{x:\emptyt} \dbl(\rec\emptyt(\natw,x))\big)) \\
 & \jdeq e_t(\big(\lamu{x:\emptyt} \rec\emptyt(\natw,x)\big), \big(\lamu{x:\emptyt} \dbl(\rec\emptyt(\natw,x))\big))\\
 & \jdeq \zerow \\
 \intertext{and}
\dbl(1^{\mathbf{w}}) & \jdeq \rec\natw(\natw, e, \supp(\btrue, \; \lamu{x:\unit} \zerow)) \\
& \jdeq e(\btrue, \big(\lamu{x:\unit} \zerow\big), \big(\lamu{x:\unit} \dbl(\zerow)\big)) \\
 & \jdeq e_f(\big(\lamu{x:\unit} \zerow\big), \big(\lamu{x:\unit} \dbl(\zerow)\big)) \\
 & \jdeq \supp(\btrue, (\lamu{x:\unit} \supp(\btrue,(\lamu{y:\unit} \dbl(\zerow))))) \\
 & \jdeq \supp(\btrue, (\lamu{x:\unit} \supp(\btrue,(\lamu{y:\unit} \zerow)))) \\
 & \jdeq 2^{\mathbf{w}}
\end{align*}
and so on.

Just as for natural numbers, we can prove a uniqueness theorem for 
$\w$-types:
\begin{thm}\label{thm:w-uniq}
  \index{uniqueness!principle, propositional!for functions on W-types@for functions on $\w$-types}%
Let $g,h : \prd{w:\wtype{x:A}B(x)} E(w)$ be two functions which satisfy the recurrence
\begin{equation*}
  e : \prd{a,f} \Parens{\prd{b : B(a)} E(f(b))} \to  E(\supp(a,f)),
\end{equation*}
i.e., such that
\begin{gather*}
 \prd{a,f} \id{g(\supp(a,f))} {e(a,f,\lamu{b:B(a)} g(f(b)))}, \\
 \prd{a,f} \id{h(\supp(a,f))}{e(a,f,\lamu{b:B(a)} h(f(b)))}.
\end{gather*}
Then $g$ and $h$ are equal. 
\end{thm}

\section{Inductive types are initial algebras}
\label{sec:initial-alg}

\indexsee{initial!algebra characterization of inductive types}{homotopy-initial}%

As suggested earlier, inductive types also have a category-theoretic universal property.
They are \emph{homotopy-initial algebras}: initial objects (up to coherent homotopy) in a category of ``algebras'' determined by the specified constructors.
As a simple example, consider the natural numbers.
The appropriate sort of ``algebra'' here is a type equipped with the same structure that the constructors of $\nat$ give to it.

\index{natural numbers!as homotopy-initial algebra}
\begin{defn}\label{defn:nalg}
  A \define{$\nat$-algebra}
  \indexdef{N-algebra@$\nat$-algebra}%
  \indexdef{algebra!N-@$\nat$-}%
  is a type $C$ with two elements $c_0 : C$, $c_s : C \to C$. The type of such algebras is
\begin{equation*}
\nalg \defeq \sm {C : \type} C \times (C \to C).
\end{equation*}
\end{defn}

\begin{defn}\label{defn:nhom}
  A \define{$\nat$-homomorphism}
  \indexdef{N-homomorphism@$\nat$-homomorphism}%
  \indexdef{homomorphism!N-@$\nat$-}%
  between $\nat$-algebras $(C,c_0,c_s)$ and $(D,d_0,d_s)$ is a function $h : C \to D$ such that $h(c_0) = d_0$ and $h(c_s(c)) = d_s(h(c))$ for all $c : C$. The type of such homomorphisms is
\begin{narrowmultline*}
\nhom((C,c_0,c_s),(D,d_0,d_s)) \defeq \narrowbreak
 \dsm {h : C \to D} (\id{h(c_0)}{d_0}) \times \tprd{c:C} (\id{h(c_s(c))}{d_s(h(c))}).
\end{narrowmultline*}
\end{defn}

We thus have a category of $\nat$-algebras and $\nat$-homomorphisms, and the claim is that $\nat$ is the initial object of this category.
A category theorist will immediately recognize this as the definition of a \emph{natural numbers object} in a category.

Of course, since our types behave like $\infty$-groupoids\index{.infinity-groupoid@$\infty$-groupoid}, we actually have an $(\infty,1)$-category\index{.infinity1-category@$(\infty,1)$-category} of $\nat$-algebras, and we should ask $\nat$ to be initial in the appropriate $(\infty,1)$-categorical sense.
Fortunately, we can formulate this without needing to define $(\infty,1)$-categories.

\begin{defn}
  \index{universal!property!of natural numbers}%
  A $\nat$-algebra $I$ is called \define{homotopy-initial},
  \indexdef{homotopy-initial!N-algebra@$\nat$-algebra}%
  \indexdef{N-algebra@$\nat$-algebra!homotopy-initial (h-initial)}%
  or \define{h-initial}
  \indexsee{h-initial}{homotopy-initial}%
  for short, if for any other $\nat$-algebra $C$, the type of $\nat$-homomorphisms from $I$ to $C$ is contractible. Thus,
\begin{equation*}
\ishinitn(I) \defeq \prd{C : \nalg} \iscontr(\nhom(I,C)).
\end{equation*}
\end{defn}

When they exist, h-initial algebras are unique --- not just up to isomorphism, as usual in category theory, but up to equality, by the univalence axiom.

\begin{thm}
  Any two h-initial $\nat$-algebras are equal.
  Thus, the type of h-initial $\nat$-algebras is a mere proposition.
\end{thm}
\begin{proof}
  Suppose $I$ and $J$ are h-initial $\nat$-algebras.
  Then $\nhom(I,J)$ is contractible, hence inhabited by some $\nat$-homomorphism $f:I\to J$, and likewise we have an $\nat$-homomorphism $g:J\to I$.
  Now the composite $g\circ f$ is a $\nat$-homomorphism from $I$ to $I$, as is $\idfunc[I]$; but $\nhom(I,I)$ is contractible, so $g\circ f = \idfunc[I]$.
  Similarly, $f\circ g = \idfunc[J]$.
  Hence $\eqv IJ$, and so $I=J$. Since being contractible is a mere proposition and dependent products preserve mere propositions, it follows that being h-initial is itself a mere proposition. Thus any two proofs that $I$ (or $J$) is h-initial are necessarily equal, which finishes the proof. 
\end{proof}

We now have the following theorem.

\begin{thm}\label{thm:nat-hinitial}
The $\nat$-algebra $(\nat, \emptyt, \suc)$ is homotopy initial.
\end{thm}
\begin{proof}[Sketch of proof]
  Fix an arbitrary $\nat$-algebra $(C,c_0,c_s)$.
  The recursion principle of $\nat$ yields a function $f:\nat\to C$ defined by
  \begin{align*}
    f(0) &\defeq c_0\\
    f(\suc(n)) &\defeq c_s(f(n)).
  \end{align*}
  These two equalities make $f$ an $\nat$-homomorphism, which we can take as the center of contraction for $\nhom(\nat,C)$.
  The uniqueness theorem (\autoref{thm:nat-uniq}) then implies that any other $\nat$-homomorphism is equal to $f$.
\end{proof}

To place this in a more general context, it is useful to consider the notion of \emph{algebra for an endofunctor}. \index{algebra!for an endofunctor}
Note that to make a type $C$ into a $\nat$-algebra is the same as to give a function $c:C+\unit\to C$, and a function $f:C\to D$ is a $\nat$-homomorphism just when $f \circ c \htpy d \circ (f+\unit)$.
In categorical language, this means the \nat-algebras are the algebras for the endofunctor $F(X)\defeq X+1$ of the category of types.

\indexsee{functor!polynomial}{endofunctor, polynomial}%
\indexsee{polynomial functor}{endofunctor, polynomial}%
\indexdef{endofunctor!polynomial}%
\index{W-type@$\w$-type!as homotopy-initial algebra}
For a more generic case, consider the $\w$-type associated to $A : \type$ and $B : A \to \type$.
In this case we have an associated \define{polynomial functor}:
\begin{equation}
\label{eq:polyfunc}
P(X) = \sm{x : A} (B(x) \rightarrow X).
\end{equation}
Actually, this assignment is functorial only up to homotopy, but this makes no difference in what follows.
By definition, a \define{$P$-algebra}
\indexdef{algebra!for a polynomial functor}%
\indexdef{algebra!W-@$\w$-}%
is then a type $C$ equipped a function $s_C :  PC \rightarrow C$.
By the universal property of $\Sigma$-types, this is equivalent to giving a function $\prd{a:A} (B(a) \to C) \to C$.
We will also call such objects \define{$\w$-algebras}
\indexdef{W-algebra@$\w$-algebra}%
for $A$ and $B$, and we write
\symlabel{walg}
\begin{equation*}
\walg(A,B) \defeq \sm {C : \type} \prd{a:A} (B(a) \to C) \to C.
\end{equation*}

Similarly, for $P$-algebras $(C,s_C)$ and $(D,s_D)$, a \define{homomorphism}
\indexdef{homomorphism!of algebras for a functor}%
between them $(f, s_f) : (C, s_C) \rightarrow (D, s_D)$ consists of a function $f : C \rightarrow D$ and a homotopy between maps $PC \rightarrow D$
\[
s_f :  f \circ s_C \, = s_{D} \circ Pf,
\]
where $Pf : PC\rightarrow PD$ is the result of the easily-definable action of $P$ on $f: C \rightarrow D$. Such an algebra homomorphism can be represented suggestively in the form:
\[
\xymatrix{
 PC \ar[d]_{s_C} \ar[r]^{Pf}  \ar@{}[dr]|{s_f} &  PD \ar[d]^{s_D}\\
C \ar[r]_{f}   & D }
\]
In terms of elements, $f$ is a $P$-homomorphism (or \define{$\w$-homomorphism}) if
\indexdef{W-homomorphism@$\w$-homomorphism}%
\indexdef{homomorphism!W-@$\w$-}%
\[f(s_C(a,h)) = s_D(a,\lam{b} f(h(b))).\]
We have the type of $\w$-homomorphisms:
\symlabel{whom}
\begin{equation*}
  \whom_{A,B}((C, c),(D,d)) \defeq \sm{h : C \to D} \prd{a:A}{f:B(a)\to C} \id{h(c(a,f))}{\lamu{b:B(a)} h(f(b))}
\end{equation*}

\index{universal!property!of $\w$-type}%
Finally, a $P$-algebra $(C, s_C)$ is said to be \define{homotopy-initial}
\indexdef{homotopy-initial!algebra for a functor}%
\indexdef{homotopy-initial!W-algebra@$\w$-algebra}%
if for every $P$-algebra $(D, s_D)$, the type of all algebra homomorphisms $(C, s_C) \rightarrow (D, s_D)$ is contractible.
That is,
\begin{equation*}
\ishinitw(A,B,I) \defeq \prd{C : \walg(A,B)} \iscontr(\whom_{A,B}(I,C)).
\end{equation*}
Now the analogous theorem to \autoref{thm:nat-hinitial} is:

\begin{thm}\label{thm:w-hinit}
For any type $A : \type$ and type family $B : A \to \type$, the $\w$-algebra $(\wtype{x:A}B(x), \supp)$ is h-initial.
\end{thm}

\begin{proof}[Sketch of proof]
Suppose we have $A : \type$ and $B : A \to \type$,
and consider the associated polynomial functor $P(X)\defeq\sm{x:A}(B(x)\to X)$.
Let $W \defeq \wtype{x:A}B(x)$.  Then using
the $\w$-introduction rule from \autoref{sec:w-types}, we have a structure map $s_W\defeq\supp: PW \rightarrow W$. 
We want to show that the algebra $(W, s_W)$ is
h-initial. So, let us consider another algebra $(C,s_C)$ and show that the type $T\defeq \whom_{A,B}((W, s_W),(C,s_C)) $ 
of  $\w$-homomorphisms from $(W, s_W)$ to $(C, s_C)$ is contractible. To do
so, observe that the $\w$-elimination rule and the $\w$-computation
rule allow us to define a $\w$-homomorphism $(f, s_f) : (W, s_W) \rightarrow (C, s_C)$, 
thus showing that $T$ is inhabited. It is furthermore necessary to show that for every $\w$-ho\-mo\-mor\-phism $(g, s_g) : (W, s_W) \rightarrow (C, s_C)$, there is an identity proof 
\begin{equation}
\label{equ:prequired}
p :  (f,s_f) = (g,s_g).
\end{equation}
This uses the fact that, in general, a type of the form $(f,s_f) = (g,s_g) $
is  equivalent to the type of what we call \define{algebra $2$-cells},
\indexdef{algebra!2-cell}%
whose canonical
elements are pairs of the form $(e, s_e)$, where $e : f=g$ and $s_e$ is a higher identity proof between the identity proofs represented by the following pasting diagrams:
\[
\xymatrix{
PW \ar@/^1pc/[r]^{Pg}   \ar[d]_{s_W} \ar@{}[r]_(.52){s_g}  & PD \ar[d]^{s_D}  \\
W \ar@/^1pc/[r]^g  \ar@/_1pc/[r]_f  \ar@{}[r]|{e} & D } \qquad
\xymatrix{
PW \ar@/^1pc/[r]^{Pg}   \ar[d]_{s_W}   \ar@/_1pc/[r]_{Pf} \ar@{}[r]|{Pe}
& PD \ar[d]^{s_D}  \\
W  \ar@/_1pc/[r]_f  \ar@{}[r]^{s_f} & D }
\]
In light of this fact, to prove that there exists an element as in~\eqref{equ:prequired}, it is 
sufficient to show that there is an algebra 2-cell 
\[
(e,s_e) : (f,s_f) = (g,s_g).
\]
The identity proof $e : f=g$ is now constructed by function extensionality and
$\w$-elimination so as to guarantee the existence of the required identity
proof $s_e$. 
\end{proof}


\section{Homotopy-inductive types}
\label{sec:htpy-inductive}

In \autoref{sec:w-types} we showed how to encode natural numbers as $\w$-types, with 
\begin{align*}
\natw & \defeq \wtype{b:\bool} \rec\bool(\bbU,\emptyt,\unit), \\
\zerow & \defeq \supp(\bfalse, (\lamu{x:\emptyt} \rec\emptyt(\natw,x))), \\
\sucw & \defeq \lamu{n:\natw} \supp(\btrue, (\lamu{x:\unit} n)).
\end{align*}
We also showed how one can define a $\dbl$ function on $\natw$ using the recursion principle.
When it comes to the induction principle, however, this encoding is no longer satisfactory: given $E : \natw \to \type$ and recurrences $e_z : E(\zerow)$ and $e_s : \prd{n : \natw}{y : E(n)} E(\sucw(n))$, we can only construct a dependent function $r(E,e_z,e_s) : \prd{n : \natw} E(n)$ satisfying the given recurrences \emph{propositionally}, i.e.\ up to a path.
This means that the computation rules for natural numbers, which give judgmental equalities, cannot be derived from the rules for $\w$-types in any obvious way.

\index{type!homotopy-inductive}%
\index{homotopy-inductive type}%
This problem goes away if instead of the conventional inductive types we consider \emph{homotopy-inductive types}, where all computation rules are stated up to a path, i.e.\ the symbol $\jdeq$ is replaced by $=$. For instance, the computation rule for the homotopy version of $\w$-types $\mathsf{W^h}$ becomes:
\index{computation rule!propositional}%
\begin{itemize}
\item For any $a : A$ and $f : B(a) \to \wtypeh{x:A} B(x)$ we have 
\begin{equation*}
  \rec{\wtypeh{x:A} B(x)}(E,\supp(a,f)) = e\Big(a,f,\big(\lamu{b:B(a)} \rec{\wtypeh{x:A} B(x)}(E,f(b))\big)\Big)
\end{equation*}
\end{itemize}

Homotopy-inductive types have an obvious disadvantage when it comes to computational properties --- the behavior of any function constructed using the induction principle can now only be characterized propositionally.
But numerous other considerations drive us to consider homotopy-inductive types as well.
For instance, while we showed in \autoref{sec:initial-alg} that inductive types are homotopy-initial algebras, not every homotopy-initial algebra is an inductive type (i.e.\ satisfies the corresponding induction principle) --- but every homotopy-initial algebra \emph{is} a homotopy-inductive type.
Similarly, we might want to apply the uniqueness argument from \autoref{sec:appetizer-univalence} when one (or both) of the types involved is only a homotopy-inductive type --- for instance, to show that the $\w$-type encoding of $\nat$ is equivalent to the usual $\nat$.

Additionally, the notion of a homotopy-inductive type is now internal to the type theory.
For example, this means we can form a type of all natural numbers objects and make assertions about it.
In the case of $\w$-types, we can characterize a homotopy $\w$-type $\wtype{x:A} B(x)$ as any type endowed with a supremum function and an induction principle satisfying the appropriate (propositional) computation rule:
\begin{multline*}
\w_d(A,B) \defeq \sm{W : \type} 
                 \sm{\supp : \prd {a} (B(a) \to W) \to W} 
                 \prd{E : W \to \type} \\
                 \prd{e : \prd{a,f} (\prd{b : B(a)} E(f(b))) \to E(\supp(a,f))}
                 \sm{\ind{} : \prd{w : W} E(w)}
                 \prd{a,f} \\
                 \ind{}(\supp(a,f)) = e(a,\lamu{b:B(a)} \ind{}(f(b))).
\end{multline*}
In \autoref{cha:hits} we will see some other reasons why propositional computation rules are worth considering.

In this section, we will state some basic facts about homotopy-inductive types.
We omit most of the proofs, which are somewhat technical.

\begin{thm}
  For any $A : \type$ and $B : A \to \type$, the type $\w_d(A,B)$ is a mere proposition.
\end{thm}

It turns out that there is an equivalent characterization of $\w$-types using a recursion principle, plus certain \emph{uniqueness} and \emph{coherence} laws. First we give the recursion principle:
\begin{itemize}
\item When constructing a function from the the $\w$-type $\wtypeh{x:A} B(x)$ into the type $C$, it suffices to give its value for $\supp(a,f)$, assuming we are given the values of all $f(b)$ with $b : B(a)$.
In other words, it suffices to construct a function
\begin{equation*}
  c : \prd{a:A} (B(a) \to C) \to C.
\end{equation*}
\end{itemize}
\index{computation rule!propositional}%
The associated computation rule for $\rec{\wtypeh{x:A} B(x)}(C,c) : (\wtype{x:A} B(x)) \to C$ is as follows:
\begin{itemize}
\item For any $a : A$ and $f : B(a) \to \wtypeh{x:A} B(x)$ we have
a witness $\beta(C,c,a,f)$ for equality
\begin{equation*}
  \rec{\wtypeh{x:A} B(x)}(C,c,\supp(a,f)) = 
  c(a,\lamu{b:B(a)} \rec{\wtypeh{x:A} B(x)}(C,c,f(b))).
\end{equation*}
\end{itemize}

Furthermore, we assert the following uniqueness principle, saying that any two functions defined by the same recurrence are equal:
\index{uniqueness!principle, propositional!for homotopy W-types@for homotopy $\w$-types}%
\begin{itemize}
\item Let $C : \type$ and $c : \prd{a:A} (B(a) \to C) \to C$ be given. Let $g,h : (\wtypeh{x:A} B(x)) \to C$ be two functions which satisfy the recurrence $c$ up to propositional equality, i.e., such that we have
\begin{align*}
  \beta_g &: \prd{a,f} \id{g(\supp(a,f))}{c(a,\lamu{b: B(a)} g(f(b)))}, \\
  \beta_h &: \prd{a,f} \id{h(\supp(a,f))}{c(a,\lamu{b: B(a)} h(f(b)))}.
\end{align*}
Then $g$ and $h$ are equal, i.e.\ there is $\alpha(C,c,f,g,\beta_g,\beta_h)$ of type $g = h$.
\end{itemize}

\index{coherence}%
Recall that when we have an induction principle rather than only a recursion principle, this propositional uniqueness principle is derivable (\autoref{thm:w-uniq}).
But with only recursion, the uniqueness principle is no longer derivable --- and in fact, the statement is not even true (exercise).  Hence, we postulate it as an axiom.
We also postulate the following coherence\index{coherence} law, which tells us how the proof of uniqueness behaves on canonical elements:
\begin{itemize}
\item
For any $a : A$ and $f : B(a) \to C$, the following diagram commutes propositionally:
\[\xymatrix{
  g(\supp(x,f)) \ar_{\alpha(\supp(x,f))}[d] \ar^-{\beta_g}[r] & c(a,\lamu{b:B(a)} g(f(b)))
  \ar^{c(a,\blank)(\funext (\lam{b} \alpha(f(b))))}[d] \\
  h(\supp(x,f)) \ar_-{\beta_h}[r] & c(a,\lamu{b: B(a)} h(f(b))) \\
}\]
where $\alpha$ abbreviates the path $\alpha(C,c,f,g,\beta_g,\beta_h) : g = h$.
\end{itemize}

Putting all of this data together yields another characterization of $\wtype{x:A} B(x)$, as a type with a supremum function, satisfying simple elimination, computation, uniqueness, and coherence rules:
\begin{multline*}
\w_s(A,B) \defeq \sm{W : \type}
                      \; \sm{\supp : \prd {a} (B(a) \to W) \to W}
                      \; \prd{C : \type}
                      \; \prd{c : \prd{a} (B(a) \to C) \to C}
                      \; \sm{\rec{} : W \to C} \\
                      \; \sm{\beta : \prd{a,f} \rec{}(\supp(a,f)) = c(a,\lamu{b: B(a)} \rec{}(f(b)))} \narrowbreak
                      \; \prd{g : W \to C}
                      \; \prd{h : W \to C}
                      \; \prd{\beta_g : \prd{a,f} g(\supp(a,f)) = c(a,\lamu{b: B(a)} g(f(b)))} \\
                      \; \prd{\beta_h : \prd{a,f} h(\supp(a,f)) = c(a,\lamu{b: B(a)} h(f(b)))}
                      \; \sm{\alpha : \prd {w : W} g(w) = h(w)} \\
                      \; \alpha(\supp(x,f)) \ct \beta_h = \beta_g \ct c(a,-)(\funext \; \lam{b} \alpha(f(b)))
\end{multline*}

\begin{thm}
For any $A : \type$ and $B : A \to \type$, the type $\w_s (A,B)$ is a mere proposition.
\end{thm}

Finally, we have a third, very concise characterization of $\wtype{x:A} B(x)$ as an h-initial $\w$-algebra:
\begin{equation*}
\w_h(A,B) \defeq \sm{I : \walg(A,B)} \ishinitw(A,B,I).
\end{equation*}

\begin{thm}
For any $A : \type$ and $B : A \to \type$, the type $\w_h (A,B)$ is a mere proposition.
\end{thm}

It turns out all three characterizations of $\w$-types are in fact equivalent:
\begin{lem}\label{lem:homotopy-induction-times-3}
For any $A : \type$ and $B : A \to \type$, we have
\[ \w_d(A,B) \eqvsym \w_s(A,B) \eqvsym \w_h(A,B) \]
\end{lem}

Indeed, we have the following theorem, which is an improvement over \autoref{thm:w-hinit}:

\begin{thm}
The types satisfying the formation, introduction, elimination, and propositional computation rules for $\w$-types are precisely the homotopy-initial $\w$-algebras.
\end{thm}

\begin{proof}[Sketch of proof]
Inspecting the proof of \autoref{thm:w-hinit}, we see that only the \emph{propositional} computation rule was required to establish the h-initiality of $\wtype{x:A}B(x)$. 
For the converse implication, let us assume that the polynomial functor associated
to $A : \type$ and $B : A \to \UU$, has an h-initial algebra $(W,s_W)$; we show that $W$ satisfies the propositional rules of $\w$-types.
The $\w$-introduction rule is simple; namely, for $a : A$ and $t : B(a) \rightarrow W$,  we define $\supp(a,t) : W$ to be the 
result of applying the structure map $s_W : PW \rightarrow W$ to $(a,t) : PW$.
For the $\w$-elimination rule, let us assume its premisses and in particular that $C' : W \to \type$.
Using the other premisses, one shows that the type $C \defeq \sm{ w : W} C'(w)$
can be equipped with a structure map $s_C : PC \rightarrow C$. By the h-initiality of $W$,
we obtain an algebra homomorphism $(f, s_f) : (W, s_W) \rightarrow (C, s_C)$. Furthermore,
the first projection $\proj1 : C \rightarrow W$ can be equipped with the structure of a
homomorphism, so that we obtain a diagram of the form
\[
\xymatrix{
PW \ar[r]^{Pf} \ar[d]_{s_W}  & PC \ar[d]^{s_C}  \ar[r]^{P \proj1}  & PW  \ar[d]^{s_W}  \\
W \ar[r]_f  & C \ar[r]_{\proj1}  & W.}
\]
But the identity function $1_W : W \rightarrow W$ has a canonical structure of an
algebra homomorphism and so, by the contractibility of the type of homomorphisms
from $(W,s_W)$ to itself, there must be an identity proof between the composite
of $(f,s_f)$ with $(\proj1, s_{\proj1})$ and $(1_W, s_{1_W})$. This implies, in particular,
that there is an identity proof $p :  \proj1 \circ f = 1_W$. 

Since $(\proj2 \circ f) w : C( (\proj1 \circ f) w)$, we can define
\[
\rec{}(w,c) \defeq
p_{\, * \,}( ( \proj2 \circ  f)   w )   : C(w) 
\]
where the transport $p_{\, * \,}$ is with respect to the family
\[
\lamu{u}C\circ u : (W\to W)\to W\to \UU.
\]
The verification of the propositional $\w$-computation rule is a calculation,
involving the naturality properties of operations of the form $p_{\, * \,}$.
\end{proof}

\index{natural numbers!encoded as a W-type@encoded as a $\w$-type}%
Finally, as desired, we can encode homotopy-natural-numbers as homo\-topy-$\w$-types:

\begin{thm}
The rules for natural numbers with propositional computation rules can be derived from the rules for $\w$-types with propositional computation rules.
\end{thm}


\section{The general syntax of inductive definitions}
\label{sec:strictly-positive}

\index{type!inductive|(}%
\indexsee{inductive!type}{type, inductive}%

So far, we have been discussing only particular inductive types: $\emptyt$, $\unit$, $\bool$, $\nat$, coproducts, products, $\Sigma$-types, $\w$-types, etc.
However, an important aspect of type theory is the ability to define \emph{new} inductive types, rather than being restricted only to some particular fixed list of them.
In order to be able to do this, however, we need to know what sorts of ``inductive definitions'' are valid or reasonable.

To see that not everything which ``looks like an inductive definition'' makes sense, consider the following ``constructor'' of a type $C$:
\begin{itemize}
\item $g:(C\to \nat) \to C$.
\end{itemize}
The recursion principle for such a type $C$ ought to say that given a type $P$, in order to construct a function $f:C\to P$, it suffices to consider the case when the input $c:C$ is of the form $g(\alpha)$ for some $\alpha:C\to\nat$.
Moreover, we would expect to be able to use the ``recursive data'' of $f$ applied to $\alpha$ in some way.
However, it is not at all clear how to ``apply $f$ to $\alpha$'', since both are functions with domain $C$.

We could write down a ``recursion principle'' for $C$ by just supposing (unjustifiably) that there is some way to apply $f$ to $\alpha$ and obtain a function $P\to\nat$.
Then the input to the recursion rule would ask for a type $P$ together with a function
\begin{equation}
  h:(C\to\nat) \to (P\to\nat) \to P\label{eq:fake-recursor}
\end{equation}
where the two arguments of $h$ are $\alpha$ and ``the result of applying $f$ to $\alpha$''.
However, what would the computation rule for the resulting function $f:C\to P$ be?
Looking at other computation rules, we would expect something like ``$f(g(\alpha)) \jdeq h(\alpha,f(\alpha))$'' for $\alpha:C\to\nat$, but as we have seen, ``$f(\alpha)$'' does not make sense.
The induction principle of $C$ is even more problematic; it's not even clear how to write down the hypotheses.
(See also \autoref{ex:loop,ex:loop2}.)

This example suggests one restriction on inductive definitions: the domains of all the constructors must be \emph{covariant functors}\index{functor!covariant}\index{covariant functor} of the type being defined, so that we can ``apply $f$ to them'' to get the result of the ``recursive call''.
In other words, if we replace all occurrences of the type being defined with a variable
\index{variable!type}%
$X:\type$, then each domain of a constructor
\index{domain!of a constructor}%
must be an expression that can be made into a covariant functor of $X$.
This is the case for all the examples we have considered so far.
For instance, with the constructor $\inl:A\to A+B$, the relevant functor is constant at $A$ (i.e.\ $X\mapsto A$), while for the constructor $\suc:\nat\to\nat$, the functor is the identity functor ($X\mapsto X$).

However, this necessary condition is also not sufficient.
Covariance prevents the inductive type from occurring on the left of a single function type, as in the argument $C\to\nat$ of the ``constructor'' $g$ considered above, since this yields a contravariant\index{functor!contravariant}\index{contravariant functor} functor rather than a covariant one.
However, since the composite of two contravariant functors is covariant, \emph{double} function types such as $((X\to \nat)\to \nat)$ are once again covariant.
This enables us to reproduce Cantorian-style paradoxes\index{paradox}.

For instance, consider an ``inductive type'' $D$ with the following constructor:
\begin{itemize}
\item $k:((D\to\prop)\to\prop)\to D$.
\end{itemize}
Assuming such a type exists, we define functions
\begin{align*}
  r&:D\to (D\to\prop)\to\prop,\\
  f&:(D\to\prop) \to D,\\
  p&:(D\to \prop) \to (D\to\prop)\to \prop,\\
  \intertext{by}
  r(k(\theta)) &\defeq \theta,\\
  f(\delta) &\defeq k(\lam{x} (x=\delta)),\\
  p(\delta) &\defeq \lam{x} \delta(f(x)).
\end{align*}
Here $r$ is defined by the recursion principle of $D$, while $f$ and $p$ are defined explicitly.
Then for any $\delta:D\to\prop$, we have $r(f(\delta)) = \lam{x}(x=\delta)$.

In particular, therefore, if $f(\delta)=f(\delta')$, then we have a path $s:(\lam{x}(x=\delta)) = (\lam{x}(x=\delta'))$.
Thus, $\happly(s,\delta) : (\delta=\delta) = (\delta=\delta')$, and so in particular $\delta=\delta'$ holds.
Hence, $f$ is ``injective'' (although \emph{a priori} $D$ may not be a set).
This already sounds suspicious --- we have an ``injection'' of the ``power set''\index{power set} of $D$ into $D$ --- and with a little more work we can massage it into a contradiction.

Suppose given $\theta:(D\to\prop)\to\prop$, and define $\delta:D\to\prop$ by
\begin{equation}
  \delta(d) \defeq \exis{\gamma:D\to\prop} (f(\gamma) = d) \times \theta(\gamma).\label{eq:Pinj}
\end{equation}
We claim that $p(\delta)=\theta$.
By function extensionality, it suffices to show $p(\delta)(\gamma) =_\prop \theta(\gamma)$ for any $\gamma:D\to\prop$.
And by univalence, for this it suffices to show that each implies the other.
Now by definition of $p$, we have
\begin{align*}
  p(\delta)(\gamma) &\jdeq \delta(f(\gamma))\\
  &\jdeq \exis{\gamma':D\to\prop} (f(\gamma') = f(\gamma)) \times \theta(\gamma').
\end{align*}
Clearly this holds if $\theta(\gamma)$, since we may take $\gamma'\defeq \gamma$.
On the other hand, if we have $\gamma'$ with $f(\gamma') = f(\gamma)$ and $\theta(\gamma')$, then $\gamma'=\gamma$ since $f$ is injective, hence also $\theta(\gamma)$.

This completes the proof that $p(\delta)=\theta$.
Thus, every element $\theta:(D\to\prop)\to\prop$ is the image under $p$ of some element $\delta:D\to\prop$.
However, if we define $\theta$ by a classic diagonalization:
\[ \theta(\gamma) \defeq \neg p(\gamma)(\gamma) \quad\text{for all $\gamma:D\to\prop$} \]
then from $\theta = p(\delta)$ we deduce $p(\delta)(\delta) = \neg p(\delta)(\delta)$.
This is a contradiction: no proposition can be equivalent to its negation.
(Supposing $P\Leftrightarrow \neg P$, if $P$, then $\neg P$, and so $\emptyt$; hence $\neg P$, but then $P$, and so $\emptyt$.)

\begin{rmk}
  There is a question of universe size to be addressed.
  In general, an inductive type must live in a universe that already contains all the types going into its definition.
  Thus if in the definition of $D$, the ambiguous notation \prop means $\prop_{\UU}$, then we do not have $D:\UU$ but only $D:\UU'$ for some larger universe $\UU'$ with $\UU:\UU'$.
  \index{mathematics!predicative}%
  \indexsee{impredicativity}{mathematics, predicative}%
  \indexsee{predicative mathematics}{mathematics, predicative}%
  In a predicative theory, therefore, the right-hand side of~\eqref{eq:Pinj} lives in $\prop_{\UU'}$, not $\prop_\UU$.
  So this contradiction does require the propositional resizing axiom
  \index{propositional!resizing}%
  mentioned in \autoref{subsec:prop-subsets}.
\end{rmk}

\index{consistency}%
This counterexample suggests that we should ban an inductive type from ever appearing on the left of an arrow in the domain of its constructors, even if that appearance is nested in other arrows so as to eventually become covariant.
(Similarly, we also forbid it from appearing in the domain of a dependent function type.)
This restriction is called \define{strict positivity}
\indexdef{strict!positivity}%
\indexsee{positivity, strict}{strict positivity}%
(ordinary ``positivity'' being essentially covariance), and it turns out to suffice.

\index{constructor}%
In conclusion, therefore, a valid inductive definition of a type $W$ consists of a list of \emph{constructors}.
Each constructor is assigned a type that is a function type taking some number (possibly zero) of inputs (possibly dependent on one another) and returning an element of $W$.
Finally, we allow $W$ itself to occur in the input types of its constructors, but only strictly positively.
This essentially means that each argument of a constructor is either a type not involving $W$, or some iterated function type with codomain $W$.
For instance, the following is a valid constructor type:
\begin{equation}
  c:(A\to W) \to (B\to C \to W) \to D \to W \to W.\label{eq:example-constructor}
\end{equation}
All of these function types can also be dependent functions ($\Pi$-types).%
\footnote{In the language of \autoref{sec:initial-alg}, the condition of strict positivity ensures that the relevant endofunctor is polynomial.\indexfoot{endofunctor!polynomial}\indexfoot{algebra!initial}\indexsee{algebra!initial}{homotopy-initial} It is well-known in category theory that not \emph{all} endofunctors can have initial algebras; restricting to polynomial functors ensures consistency.
One can consider various relaxations of this condition, but in this book we will restrict ourselves to strict positivity as defined here.}

Note we require that an inductive definition is given by a \emph{finite} list of constructors.
This is simply because we have to write it down on the page.
If we want an inductive type which behaves as if it has an infinite number of constructors, we can simply parametrize one constructor by some infinite type.
For instance, a constructor such as $\nat \to W \to W$ can be thought of as equivalent to countably many constructors of the form $W\to W$.
(Of course, the infinity is now \emph{internal} to the type theory, but this is as it should be for any foundational system.)
Similarly, if we want a constructor that takes ``infinitely many arguments'', we can allow it to take a family of arguments parametrized by some infinite type, such as $(\nat\to W) \to W$ which takes an infinite sequence\index{sequence} of elements of $W$.

\index{recursion principle!for an inductive type}%
Now, once we have such an inductive definition, what can we do with it?
Firstly, there is a \define{recursion principle} stating that in order to define a function $f:W\to P$, it suffices to consider the case when the input $w:W$ arises from one of the constructors, allowing ourselves to recursively call $f$ on the inputs to that constructor.
For the example constructor~\eqref{eq:example-constructor}, we would require $P$ to be equipped with a function of type
\begin{narrowmultline}\label{eq:example-rechyp}
  d : (A\to W) \to (A\to P) \to (B\to C\to W) \to 
  \narrowbreak
  (B\to C \to P) \to D \to W \to P \to P.
\end{narrowmultline}
Under these hypotheses, the recursion principle yields $f:W\to P$, which moreover ``preserves the constructor data'' in the evident way --- this is the computation rule, where we use covariance of the inputs.
\index{computation rule!for inductive types}%
For instance, in the example~\eqref{eq:example-constructor}, the computation rule says that for any $\alpha:A\to W$, $\beta:B\to C\to W$, $\delta:d$, and $\omega:W$, we have
\begin{equation}
  f(c(\alpha,\beta,\delta,\omega)) \jdeq d(\alpha,f\circ \alpha,\beta, f\circ \beta, \delta, \omega,f(\omega)).\label{eq:example-comp}
\end{equation}

\index{induction principle!for an inductive type}%
The \define{induction principle} for a general inductive type $W$ is only a little more complicated.
Of course, we start with a type family $P:W\to\type$, which we require to be equipped with constructor data ``lying over'' the constructor data of $W$.
That means the ``recursive call'' arguments such as $A\to P$ above must be replaced by dependent functions with types such as $\prd{a:A} P(\alpha(a))$.
In the full example of~\eqref{eq:example-constructor}, the corresponding hypothesis for the induction principle would require
\begin{multline}\label{eq:example-indhyp}
d : \prd{\alpha:A\to W}\Parens{\prd{a:A} P(\alpha(a))} \to \narrowbreak
\prd{\beta:B\to C\to W} \Parens{\prd{b:B}{c:C} P(\beta(b,c))} \to\\
\prd{\delta:D}
\prd{\omega:W} P(\omega) \to
P(c(\alpha,\beta,\delta,\omega)).
\end{multline}
The corresponding computation rule looks identical to~\eqref{eq:example-comp}.
Of course, the recursion principle is the special case of the induction principle where $P$ is a constant family.
As we have mentioned before, the induction principle is also called the \define{eliminator}, and the recursion principle the \define{non-dependent eliminator}.

As discussed in \autoref{sec:pattern-matching}, we also allow ourselves to invoke the induction and recursion principles implicitly, writing a definitional equation with $\defeq$ for each expression that would be the hypotheses of the induction principle.
This is called giving a definition by (dependent) \define{pattern matching}.
\index{pattern matching}%
\index{definition!by pattern matching}%
In our running example, this means we could define $f:\prd{w:W} P(w) $ by
\[ f(c(\alpha,\beta,\delta,\omega)) \defeq \cdots \]
where $\alpha:A\to W$ and $\beta:B\to C\to W$ and $\delta:D$ and $\omega:W$ are variables
\index{variable}%
that are bound in the right-hand side.
Moreover, the right-hand side may involve recursive calls to $f$ of the form $f(\alpha(a))$, $f(\beta(b,c))$, and $f(\omega)$.
When this definition is repackaged in terms of the induction principle, we replace such recursive calls by $\bar\alpha(a)$, $\bar\beta(b,c)$, and $\bar\omega$, respectively, for new variables
\begin{align*}
  \bar\alpha &: \prd{a:A} P(\alpha(a))\\
  \bar\beta &: \prd{b:B}{c:C} P(\beta(b,c))\\
  \bar\omega &: P(\omega).
\end{align*}
\symlabel{defn:induction-wtype}%
Then we could write
\[ f \defeq \ind{W}(P,\, \lam{\alpha}{\bar\alpha}{\beta}{\bar\beta}{\delta}{\omega}{\bar\omega} \cdots ) \]
where the second argument to $\ind{W}$ has the type of~\eqref{eq:example-indhyp}.

We will not attempt to give a formal presentation of the grammar of a valid inductive definition and its resulting induction and recursion principles and pattern matching rules.
This is possible to do (indeed, it is necessary to do if implementing a computer proof assistant), but provides no additional insight.
With practice, one learns to automatically deduce the induction and recursion principles for any inductive definition, and to use them without having to think twice.

\section{Generalizations of inductive types}
\label{sec:generalizations}

\index{type!inductive!generalizations}%
The notion of inductive type has been studied in type theory for many years, and admits of many, many generalizations: inductive type families, mutual inductive types, inductive-inductive types, inductive-recursive types, etc.
In this section we give an overview of some of these, a few of which will be used later in the book.
(In \autoref{cha:hits} we will study in more depth a very different generalization of inductive types, which is particular to \emph{homotopy} type theory.)

Most of these generalizations involve allowing ourselves to define more than one type by induction at the same time.
One very simple example of this, which we have already seen, is the coproduct $A+B$.
It would be tedious indeed if we had to write down separate inductive definitions for $\nat+\nat$, for $\nat+\bool$, for $\bool+\bool$, and so on every time we wanted to consider the coproduct of two types.
Instead, we make one definition in which $A$ and $B$ are variables standing for types;
\index{variable!type}%
in type theory they are called \define{parameters}.%
\indexdef{parameter!of an inductive definition}
Thus technically speaking, what results from the definition is not a single type, but a family of types $+ : \type\to\type\to\type$, taking two types as input and producing their coproduct.
Similarly, the type $\lst A$ of lists\index{type!of lists} is a family $\lst{\nameless}:\type\to\type$ in which the type $A$ is a parameter.

In mathematics, this sort of thing is so obvious as to not be worth mentioning, but we bring it up in order to contrast it with the next example.
Note that each type $A+B$ is \emph{independently} defined inductively, as is each type $\lst A$.
\index{type!family of!inductive}%
\index{inductive!type family}%
By contrast, we might also consider defining a whole type family $B:A\to\type$ by induction \emph{together}.
The difference is that now the constructors may change the index $a:A$, and as a consequence we cannot say that the individual types $B(a)$ are inductively defined, only that the entire family is inductively defined.

\index{type!of vectors}%
\index{vector}%
The standard example is the type of \emph{lists of specified length}, traditionally called \define{vectors}.
We fix a parameter type $A$, and define a type family $\vect n A$, for $n:\nat$, generated by the following constructors:
\begin{itemize}
\item a vector $\nil:\vect 0 A$ of length zero,
\item a function $\cons:\prd{n:\nat} A\to \vect n A \to \vect{\suc (n)} A$.
\end{itemize}
In contrast to lists, vectors (with elements from a fixed type $A$) form a family of types indexed by their length.
While $A$ is a parameter, we say that $n:\nat$ is an \define{index}
\indexdef{index of an inductive definition}%
of the inductive family.
An individual type such as $\vect3A$ is not inductively defined: the constructors which build elements of $\vect3A$ take input from a different type in the family, such as $\cons:A \to \vect2A \to \vect3A$.

\index{induction principle!for type of vectors}
\index{vector!induction principle for}
In particular, the induction principle must refer to the entire type family as well; thus the hypotheses and the conclusion must quantify over the indices appropriately.
In the case of vectors, the induction principle states that given a type family $C:\prd{n:\nat} \vect n A \to \type$, together with
\begin{itemize}
\item an element $c_\nil : C(0,\nil)$, and
\item a function \narrowequation{c_\cons : \prd{n:\nat}{a:A}{\ell:\vect n A} C(n,\ell) \to C(\suc(n),\cons(a,\ell))}
\end{itemize}
there exists a function $f:\prd{n:\nat}{\ell:\vect n A} C(n,\ell)$ such that
\begin{align*}
  f(0,\nil) &\jdeq c_\nil\\
  f(\suc(n),\cons(a,\ell)) &\jdeq c_\cons(n,a,\ell,f(\ell)).
\end{align*}

\index{predicate!inductive}%
\index{inductive!predicate}%
One use of inductive families is to define \emph{predicates} inductively.
For instance, we might define the predicate $\mathsf{iseven}:\nat\to\type$ as an inductive family indexed by $\nat$, with the following constructors:
\begin{itemize}
\item an element $\mathsf{even}_0 : \mathsf{iseven}(0)$,
\item a function $\mathsf{even}_{ss} : \prd{n:\nat} \mathsf{iseven}(n) \to \mathsf{iseven}(\suc(\suc(n)))$.
\end{itemize}
In other words, we stipulate that $0$ is even, and that if $n$ is even then so is $\suc(\suc(n))$.
These constructors ``obviously'' give no way to construct an element of, say, $\mathsf{iseven}(1)$, and since $\mathsf{iseven}$ is supposed to be freely generated by these constructors, there must be no such element.
(Actually proving that $\neg \mathsf{iseven}(1)$ is not entirely trivial, however).
The induction principle for $\mathsf{iseven}$ says that to prove something about all even natural numbers, it suffices to prove it for $0$ and verify that it is preserved by adding two.

\index{mathematics!formalized}%
Inductively defined predicates are much used in computer formalization of mathematics and software verification.
But we will not have much use for them, with a couple of exceptions in \autoref{sec:ordinals,sec:compactness-interval}.

\index{type!mutual inductive}%
\index{mutual inductive type}%
Another important special case is when the indexing type of an inductive family is finite.
In this case, we can equivalently express the inductive definition as a finite collection of types defined by \emph{mutual induction}.
For instance, we might define the types $\mathsf{even}$ and $\mathsf{odd}$ of even and odd natural numbers by mutual induction, where $\mathsf{even}$ is generated by constructors
\begin{itemize}
\item $0:\mathsf{even}$ and
\item $\mathsf{esucc} : \mathsf{odd}\to\mathsf{even}$,
\end{itemize}
while $\mathsf{odd}$ is generated by the one constructor
\begin{itemize}
\item $\mathsf{osucc} : \mathsf{even}\to \mathsf{odd}$.
\end{itemize}
Note that $\mathsf{even}$ and $\mathsf{odd}$ are simple types (not type families), but their constructors can refer to each other.
If we expressed this definition as an inductive type family $\mathsf{paritynat} : \bool \to \type$, with $\mathsf{paritynat}(\bfalse)$ and $\mathsf{paritynat}(\btrue)$ representing $\mathsf{even}$ and $\mathsf{odd}$ respectively, it would instead have constructors:
\begin{itemize}
\item $0 : \mathsf{paritynat}(\bfalse)$,
\item $\mathsf{esucc} : \mathsf{paritynat}(\bfalse) \to \mathsf{paritynat}(\btrue)$,
\item $\mathsf{oesucc} : \mathsf{paritynat}(\btrue) \to \mathsf{paritynat}(\bfalse)$.
\end{itemize}
When expressed explicitly as a mutual inductive definition, the induction principle for $\mathsf{even}$ and $\mathsf{odd}$ says that given $C:\mathsf{even}\to\type$ and $D:\mathsf{odd}\to\type$, along with
\begin{itemize}
\item $c_0 : C(0)$,
\item $c_s : \prd{n:\mathsf{odd}} D(n) \to C(\mathsf{esucc}(n))$,
\item $d_s : \prd{n:\mathsf{even}} C(n) \to D(\mathsf{osucc}(n))$,
\end{itemize}
there exist $f:\prd{n:\mathsf{even}} C(n)$ and $g:\prd{n:\mathsf{odd}}D(n)$ such that
\begin{align*}
  f(0) &\jdeq c_0\\
  f(\mathsf{esucc}(n)) &\jdeq c_s(g(n))\\
  g(\mathsf{osucc}(n)) &\jdeq d_s(f(n)).
\end{align*}
In particular, just as we can only induct over an inductive family ``all at once'', we have to induct on $\mathsf{even}$ and $\mathsf{odd}$ simultaneously.
We will not have much use for mutual inductive definitions in this book either.

\index{type!inductive-inductive}%
\index{inductive-inductive type}%
A further, more radical, generalization is to allow definition of a type family $B:A\to \type$ in which not only the types $B(a)$, but the type $A$ itself, is defined as part of one big induction.
In other words, not only do we specify constructors for the $B(a)$s which can take inputs from other $B(a')$s, as with inductive families, we also at the same time specify constructors for $A$ itself, which can take inputs from the $B(a)$s.
This can be regarded as an inductive family in which the indices are inductively defined simultaneously with the indexed types, or as a mutual inductive definition in which one of the types can depend on the other.
More complicated dependency structures are also possible.
In general, these are called \define{inductive-inductive definitions}.
For the most part, we will not use them in this book, but their higher variant (see \autoref{cha:hits}) will appear in a couple of experimental examples in \autoref{cha:real-numbers}.

\index{type!inductive-recursive}%
\index{inductive-recursive type}%
The last generalization we wish to mention is \define{inductive-recursive definitions}, in which a type is defined inductively at the same time as a \emph{recursive} function on it.
That is, we fix a known type $P$, and give constructors for an inductive type $A$ and at the same time define a function $f:A\to P$ using the recursion principle for $A$ resulting from its constructors --- with the twist that the constructors of $A$ are allowed to refer also to the values of $f$.
We do not yet know how to justify such definitions from a homotopical perspective, and we will not use any of them in this book.

\index{type!inductive|)}%

\section{Identity types and identity systems}
\label{sec:identity-systems}

\index{type!identity!as inductive}%
We now wish to point out that the \emph{identity types}, which play so central a role in homotopy type theory, may also be considered to be defined inductively.
Specifically, they are an ``inductive family'' with indices, in the sense of \autoref{sec:generalizations}.
In fact, there are \emph{two} ways to describe identity types as an
inductive family, resulting in the two induction principles described in
\cref{cha:typetheory}, path induction and based path induction.  

In both definitions, the type $A$ is a parameter.
For the first definition, we inductively define a family $=_A : A\to A\to \type$, with two indices belonging to $A$, by the following constructor:
\begin{itemize}
\item for any $a:A$, an element $\refl A : a=_A a$.
\end{itemize}
By analogy with the other inductive families, we may extract the induction principle from this definition.
It states that given any \narrowequation{C:\prd{a,b:A} (a=_A b) \to \type,} along with $d:\prd{a:A} C(a,a,\refl{a})$, there exists \narrowequation{f:\prd{a,b:A}{p:a=_A b} C(a,b,p)} such that $f(a,a,\refl a)\jdeq d(a)$.
This is exactly the path induction principle for identity types.

For the second definition, we consider one element $a_0:A$ to be a parameter along with $A:\type$, and we inductively define a family $(a_0 =_A \blank):A\to \type$, with \emph{one} index belonging to $A$, by the following constructor:
\begin{itemize}
\item an element $\refl{a_0} : a_0 =_A a_0$.
\end{itemize}
Note that because $a_0:A$ was fixed as a parameter, the constructor $\refl{a_0}$ does not appear inside the inductive definition as a function, but only an element.
The induction principle for this definition says that given $C:\prd{b:A} (a_0 =_A b) \to \type$ along with an element $d:C(a_0,\refl{a_0})$, there exists $f:\prd{b:A}{p:a_0 =_A b} C(b,p)$ with $f(a_0,\refl{a_0})\jdeq d$.
This is exactly the based path induction principle for identity types.

The view of identity types as inductive types has historically caused some confusion, because of the intuition mentioned in \autoref{sec:bool-nat} that all the elements of an inductive type should be obtained by repeatedly applying its constructors.
For ordinary inductive types such as \bool and \nat, this is the case: we saw in \autoref{thm:allbool-trueorfalse} that indeed every element of \bool is either $\bfalse$ or $\btrue$, and similarly one can prove that every element of \nat is either $0$ or a successor.

However, this is \emph{not} true for identity types: there is only one constructor $\refl{}$, but not every path is equal to the constant path.
More precisely, we cannot prove, using only the induction principle for identity types (either one), that every inhabitant of $a=_A a$ is equal to $\refl a$.
In order to actually exhibit a counterexample, we need some additional principle such as the univalence axiom --- recall that in \autoref{thm:type-is-not-a-set} we used univalence to exhibit a particular path $\bool=_\type\bool$ which is not equal to $\refl{\bool}$.

\index{free!generation of an inductive type}%
\index{generation!of a type, inductive|(}%
The point is that, as validated by the study of homotopy-initial algebras, an inductive definition should be regarded as \emph{freely generated} by its constructors.
Of course, a freely generated structure may contain elements other than its generators: for instance, the free group on two symbols $x$ and $y$ contains not only $x$ and $y$ but also words such as $xy$, $yx^{-1}y$, and $x^3y^2x^{-2}yx$.
In general, the elements of a free structure are obtained by applying not only the generators, but also the operations of the ambient structure, such as the group operations if we are talking about free groups.

In the case of inductive types, we are talking about freely generated \emph{types} --- so what are the ``operations'' of the structure of a type?
If types are viewed as like \emph{sets}, as was traditionally the case in type theory, then there are no such operations, and hence we expect there to be no elements in an inductive type other than those resulting from its constructors.
In homotopy type theory, we view types as like \emph{spaces} or $\infty$-groupoids,%
\index{.infinity-groupoid@$\infty$-groupoid}
in which case there are many operations on the \emph{paths} (concatenation, inversion, etc.) --- this will be important in \autoref{cha:hits} --- but there are still no operations on the \emph{objects} (elements).
Thus, it is still true for us that, e.g., every element of \bool is either $\bfalse$ or $\btrue$, and every element of $\nat$ is either $0$ or a successor.

However, as we saw in \autoref{cha:basics}, viewing types as $\infty$-groupoids entails also viewing functions as functors, and this includes type families $B:A\to\type$.
Thus, the identity type $(a_0 =_A \blank)$, viewed as an inductive type family, is actually a \emph{freely generated functor} $A\to\type$.
Specifically, it is the functor $F:A\to\type$ freely generated by one element $\refl{a_0}: F(a_0)$.
And a functor does have operations on objects, namely the action of the morphisms (paths) of $A$.

In category theory, the \emph{Yoneda lemma}\index{Yoneda!lemma} tells us that for any category $A$ and object $a_0$, the functor freely generated by an element of $F(a_0)$ is the representable functor $\hom_A(a_0,\blank)$.
Thus, we should expect the identity type $(a_0 =_A \blank)$ to be this representable functor, and this is indeed exactly how we view it: $(a_0 =_A b)$ is the space of morphisms (paths) in $A$ from $a_0$ to $b$.

\index{generation!of a type, inductive|)}

\mentalpause

One reason for viewing identity types as inductive families is to apply the uniqueness principles of \autoref{sec:appetizer-univalence,sec:htpy-inductive}.
Specifically, we can characterize the family of identity types of a type $A$, up to equivalence, by giving another family of types over $A\times A$ satisfying the same induction principle.
This suggests the following definitions and theorem.

\indexsee{system, identity}{identity system}%
\index{identity!system!at a point|(defstyle}%

\begin{defn}\label{defn:identity-systems}
  Let $A$ be a type and $a_0:A$ an element.
  \begin{itemize}
  \item A \define{pointed predicate}
    \indexdef{predicate!pointed}%
    \indexdef{pointed!predicate}%
    over $(A,a_0)$ is a family $R:A\to\type$ equipped with an element $r_0:R(a_0)$.
  \item For pointed predicates $(R,r_0)$ and $(S,s_0)$, a family of maps $g:\prd{b:A} R(b) \to S(b)$ is \define{pointed} if $g(a_0, r_0)=s_0$.
    We have
    \[ \mathsf{ppmap}(R,S) \defeq \sm{g:\prd{b:A} R(b) \to S(b)} (g(a_0, r_0)=s_0).\]
  \item An \define{identity system at $a_0$}
    is a pointed predicate $(R,r_0)$ such that for any type family $D:\prd{b:A} R(b) \to \type$ and $d:D(a_0,r_0)$, there exists a function $f:\prd{b:A}{r:R(b)} D(b,r)$ such that $f(a_0,r_0)=d$.
\end{itemize}
\end{defn}

\begin{thm}\label{thm:identity-systems}
  For a pointed predicate $(R,r_0)$, the following are logically equivalent.
  \begin{enumerate}
  \item $(R,r_0)$ is an identity system at $a_0$.\label{item:identity-systems1}
  \item For any pointed predicate $(S,s_0)$, the type $\mathsf{ppmap}(R,S)$ is contractible.\label{item:identity-systems2}
  \item For any $b:A$, the function $\transfib{R}{\nameless}{r_0} : (a_0 =_A b) \to R(b)$ is an equivalence.\label{item:identity-systems3}
  \item The type $\sm{b:A} R(b)$ is contractible.\label{item:identity-systems4}
  \end{enumerate}
\end{thm}

Note that the equivalences~\ref{item:identity-systems1}$\Leftrightarrow$\ref{item:identity-systems2}$\Leftrightarrow$\ref{item:identity-systems3} are a version of \autoref{lem:homotopy-induction-times-3} for identity types $a_0 =_A \blank$, regarded as inductive families varying over one element of $A$.
Of course,~\ref{item:identity-systems2}--\ref{item:identity-systems4} are mere propositions, so that logical equivalence implies actual equivalence.
(Condition~\ref{item:identity-systems1} is also a mere proposition, but we will not prove this.)

\begin{proof}
  First, assume~\ref{item:identity-systems1} and let $(S,s_0)$ be a pointed predicate.
  Define $D(b,r) \defeq S(b)$ and $d\defeq s_0: S(a_0) \jdeq D(a_0,r_0)$.
  Since $R$ is an identity system, we have $f:\prd{b:A} R(b) \to S(b)$ with $f(a_0,r_0) = s_0$; hence $\mathsf{ppmap}(R,S)$ is inhabited.
  Now suppose $(f,f_r),(g,g_r) : \mathsf{ppmap}(R,S)$, and define $D(b,r) \defeq (f(b,r) = g(b,r))$, and let $d \defeq f_r \ct \opp{g_r} : f(a_0,r_0) = s_0 = g(a_0,r_0)$.
  Then again since $R$ is an identity system, we have $h:\prd{b:A}{r:R(b)} D(b,r)$ such that $h(a_0,r_0) = f_r \ct \opp{g_r}$.
  By the characterization of paths in $\Sigma$-types and path types, these data yield an equality $(f,f_r) = (g,g_r)$.
  Hence $\mathsf{ppmap}(R,S)$ is an inhabited mere proposition, and thus contractible; so~\ref{item:identity-systems2} holds.

  Now suppose~\ref{item:identity-systems2}, and define $S(b) \defeq (a_0=b)$ with $s_0 \defeq \refl{a_0}:S(a_0)$.
  Then $(S,s_0)$ is a pointed predicate, and $\lamu{b:B}{p:a_0=b} \transfib{R}{p}{r} : \prd{b:A} S(b) \to R(b)$ is a pointed family of maps from $S$ to $R$.
  By assumption, $\mathsf{ppmap}(R,S)$ is contractible, hence inhabited, so there also exists a pointed family of maps from $R$ to $S$.
  And the composites in either direction are pointed families of maps from $R$ to $R$ and from $S$ to $S$, respectively, hence equal to identities since $\mathsf{ppmap}(R,R)$ and $\mathsf{ppmap}(S,S)$ are contractible.
  Thus~\ref{item:identity-systems3} holds.

  Now supposing~\ref{item:identity-systems3}, condition~\ref{item:identity-systems4} follows from \autoref{thm:contr-paths}, using the fact that $\Sigma$-types respect equivalences (the ``if'' direction of \autoref{thm:total-fiber-equiv}).

  Finally, assume~\ref{item:identity-systems4}, and let $D:\prd{b:A} R(b)\to  \type$ and $d:D(a_0,r_0)$.
  We can equivalently express $D$ as a family $D':(\sm{b:A} R(b)) \to \type$.
  Now since $\sm{b:A} R(b)$ is contractible, we have
  \[p:\prd{u:\sm{b:A} R(b)} (a_0,r_0) = u. \]
  Moreover, since the path types of a contractible type are again contractible, we have $p((a_0,r_0)) = \refl{(a_0,r_0)}$.
  Define $f(u) \defeq \transfib{D'}{p(u)}{d}$, yielding $f:\prd{u:\sm{b:A} R(b)} D'(u)$, or equivalently $f:\prd{b:A}{r:R(b)} D(b,r)$.
  Finally, we have
  \[f(a_0,r_0) \jdeq \transfib{D'}{p((a_0,r_0))}{d} = \transfib{D'}{\refl{(a_0,r_0)}}{d} = d.\]
  Thus,~\ref{item:identity-systems1} holds.
\end{proof}

\index{identity!system!at a point|)}%

We can deduce a similar result for identity types $=_A$, regarded as a family varying over two elements of $A$.

\index{identity!system|(defstyle}%

\begin{defn}
  An \define{identity system}
  over a type $A$ is a family $R:A\to A\to \type$ equipped with a function $r_0:\prd{a:A} R(a,a)$ such that for any type family $D:\prd{a,b:A} R(a,b) \to \type$ and $d:\prd{a:A} D(a,a,r_0(a))$, there exists a function $f:\prd{a,b:A}{r:R(b)} D(a,b,r)$ such that $f(a,a,r_0(a))=d(a)$ for all $a:A$.
\end{defn}

\begin{thm}\label{thm:ML-identity-systems}
  For $R:A\to A\to\type$ equipped with $r_0:\prd{a:A} R(a,a)$, the following are logically equivalent.
  \begin{enumerate}
  \item $(R,r_0)$ is an identity system over $A$.\label{item:MLis1}
  \item For all $a_0:A$, the pointed predicate $(R(a_0),r_0(a_0))$ is an identity system at $a_0$.\label{item:MLis2}
  \item For any $S:A\to A\to\type$ and $s_0:\prd{a:A} S(a,a)$, the type
    \[ \sm{g:\prd{a,b:A} R(a,b) \to S(a,b)} \prd{a:A} g(a,a,r_0(a)) = s_0(a) \]
    is contractible.\label{item:MLis3}
  \item For any $a,b:A$, the map $\transfib{R(a)}{\nameless}{r_0(a)} : (a =_A b) \to R(a,b)$ is an equivalence.\label{item:MLis4}
  \item For any $a:A$, the type $\sm{b:A} R(a,b)$ is contractible.\label{item:MLis5}
  \end{enumerate}
\end{thm}
\begin{proof}
  The equivalence~\ref{item:MLis1}$\Leftrightarrow$\ref{item:MLis2} follows exactly the proof of equivalence between the path induction and based path induction principles for identity types; see \autoref{sec:identity-types}.
  The equivalence with~\ref{item:MLis4} and~\ref{item:MLis5} then follows from \autoref{thm:identity-systems}, while~\ref{item:MLis3} is straightforward.
\end{proof}

\index{identity!system|)}%

One reason this characterization is interesting is that it provides an alternative way to state univalence and function extensionality.
\index{univalence axiom}%
The univalence axiom for a universe \UU says exactly that the type family
\[ (\eqv{\blank}{\blank}) : \UU\to\UU\to\UU \]
together with $\idfunc : \prd{A:\UU} (\eqv AA)$ satisfies \autoref{thm:ML-identity-systems}\ref{item:MLis4}.
Therefore, it is equivalent to the corresponding version of~\ref{item:MLis1}, which we can state as follows.

\begin{cor}[Equivalence induction]\label{thm:equiv-induction}
  \index{induction principle!for equivalences}%
  \index{equivalence!induction}%
  Given any type family \narrowequation{D:\prd{A,B:\UU} (\eqv AB) \to \type} and function $d:\prd{A:\UU} D(A,A,\idfunc[A])$, there exists \narrowequation{f : \prd{A,B:\UU}{e:\eqv AB} D(A,B,e)} such that $f(A,A,\idfunc[A]) = d(A)$ for all $A:\UU$.
\end{cor}

In other words, to prove something about all equivalences, it suffices to prove it about identity maps.
We have already used this principle (without stating it in generality) in \autoref{lem:qinv-autohtpy}.

Similarly, function extensionality says that for any $B:A\to\type$, the type family
\[ (\blank\htpy\blank) : \Parens{\prd{a:A} B(a)} \to \Parens{\prd{a:A} B(a)} \to \type
\]
together with $\lamu{f:\prd{a:A} B(a)}{a:A} \refl{f(a)}$ satisfies \autoref{thm:ML-identity-systems}\ref{item:MLis4}.
Thus, it is also equivalent to the corresponding version of~\ref{item:MLis1}.

\begin{cor}[Homotopy induction]\label{thm:htpy-induction}
  \index{induction principle!for homotopies}%
  \index{homotopy!induction}%
  Given any \narrowequation{D:\prd{f,g:\prd{a:A} B(a)} (f\htpy g) \to \type} and $d:\prd{f:\prd{a:A} B(a)} D(f,f,\lam{x}\refl{f(x)})$, there exists
  \begin{equation*}
    k:\prd{f,g:\prd{a:A} B(a)}{h:f\htpy g} D(f,g,h)    
  \end{equation*}
  such that $k(f,f,\lam{x}\refl{f(x)}) = d(f)$ for all $f$.
\end{cor}

\sectionNotes

Inductive definitions have a long pedigree in mathematics, arguably going back at least to Frege and Peano's axioms for the natural numbers.\index{Frege}\index{Peano} %
More general ``inductive predicates'' are not uncommon, but in set theoretic foundations they are usually constructed explicitly, either as an intersection of an appropriate class of subsets or using transfinite iteration along the ordinals, rather than regarded as a basic notion.

In type theory, particular cases of inductive definitions date back to Martin-L\"of's original papers: \cite{martin-lof-hauptsatz} presents a general notion of inductively defined predicates and relations; the notion of inductive type was present (but only with instances, not as a general notion) in Martin-L\"of's first papers in type theory \cite{Martin-Lof-1973};
and then as a general notion with $\w$-types in \cite{Martin-Lof-1979}.\index{Martin-L\"of}%

A general notion of inductive type was introduced in 1985 by Constable and Mendler~\cite{DBLP:conf/lop/ConstableM85}.  A general schema for inductive types in intensional type theory was suggested in
\cite{PfenningPaulinMohring}.  Further developments included \cite{CoquandPaulin, Dybjer:1991}.

The notion of inductive-recursive definition appears in \cite{Dybjer:2000}. An important  type-theoretic notion is the notion of tree types (a general expression of the notion of Post system in type theory) which appears in \cite{PeterssonSynek}.

The universal property of the natural numbers as an initial object of the category of $\nat$-algebras is due to Lawvere \cite{lawvere:adjinfound}\index{Lawvere}.
This was later generalized to a description of $\w$-types as initial algebras for polynomial endofunctors by~\cite{mp:wftrees}.\index{endofunctor!algebra for}
The coherently homotopy-theoretic equivalence between such universal properties and the corresponding induction principles (\autoref{sec:initial-alg,sec:htpy-inductive}) is due to~\cite{ags:it-hott}.

For actual constructions of inductive types in homotopy-theoretic semantics of type theory, see~\cite{klv:ssetmodel,mvdb:wtypes,ls:hits}.

\sectionExercises

\begin{ex}
  Derive the induction principle for the type $\lst{A}$ of lists from its definition as an inductive type in \autoref{sec:bool-nat}.\index{type!of lists}
\end{ex}

\begin{ex}\label{ex:same-recurrence-not-defeq}
  Construct two functions on natural numbers which satisfy the same recurrence\index{recurrence} $(e_z, e_s)$ but are not definitionally equal.
\end{ex}

\begin{ex}\label{ex:one-function-two-recurrences}
  Construct two different recurrences $(e_z,e_s)$ on the same type $E$ which are both satisfied by the same function $f:\nat\to E$.
\end{ex}

\begin{ex}\label{ex:bool}
  Show that for any type family $E : \bool \to \type$, the induction operator
  \[ \ind{\bool}(E) : \big(E(\bfalse) \times E(\btrue)\big) \to \prd{b : \bool} E(b) \]
  is an equivalence.
\end{ex}

\begin{ex}
  Show that the analogous statement to \autoref{ex:bool} for $\nat$ fails.
\end{ex}

\begin{ex}
  Show that if we assume simple instead of dependent elimination for $\w$-types, the uniqueness property (analogue of \autoref{thm:w-uniq}) fails to hold.
  That is, exhibit a type satisfying the recursion principle of a $\w$-type, but for which functions are not determined uniquely by their recurrence\index{recurrence}.
\end{ex}

\begin{ex}\label{ex:loop}
  Suppose that in the ``inductive definition'' of the type $C$ at the beginning of \autoref{sec:strictly-positive}, we replace the type \nat by \emptyt.
  Using only a ``recursion principle'' for such a definition with hypotheses analogous to~\eqref{eq:fake-recursor}, construct an element of \emptyt.
\end{ex}

\begin{ex}\label{ex:loop2}
  Similarly to the previous exercise, derive a contradiction from an ``inductive type'' $D$ with one constructor $(D\to D) \to D$.
\end{ex}


\chapter{Higher inductive types}
\label{cha:hits}

\index{type!higher inductive|(}%
\indexsee{inductive!type!higher}{type, higher inductive}%
\indexsee{higher inductive type}{type, higher inductive}%

\section{Introduction}
\label{sec:intro-hits}

\index{generation!of a type, inductive|(}

Like the general inductive types we discussed in \autoref{cha:induction}, \emph{higher inductive types} are a general schema for defining new types generated by some constructors.
But unlike ordinary inductive types, in defining a higher inductive type we may have ``constructors'' which generate not only \emph{points} of that type, but also \emph{paths} and higher paths in that type.
\index{type!circle}%
\indexsee{circle type}{type,circle}%
For instance, we can consider the higher inductive type $\Sn^1$ generated by
\begin{itemize}
\item A point $\base:\Sn^1$, and
\item A path $\lloop : {\id[\Sn^1]\base\base}$.
\end{itemize}
This should be regarded as entirely analogous to the definition of, for instance, $\bool$, as being generated by
\begin{itemize}
\item A point $\bfalse:\bool$ and
\item A point $\btrue:\bool$,
\end{itemize}
or the definition of $\nat$ as generated by
\begin{itemize}
\item A point $0:\nat$ and
\item A function $\suc:\nat\to\nat$.
\end{itemize}
When we think of types as higher groupoids, the more general notion of ``generation'' is very natural:
since a higher groupoid is a ``multi-sorted object'' with paths and higher paths as well as points, we should allow ``generators'' in all dimensions.

We will refer to the ordinary sort of constructors (such as $\base$) as \define{point constructors}
\indexdef{constructor!point}%
\indexdef{point!constructor}%
or \emph{ordinary constructors}, and to the others (such as $\lloop$) as \define{path constructors}
\indexdef{constructor!path}%
\indexdef{path!constructor}%
or \emph{higher constructors}.
Each path constructor must specify the starting and ending point of the path, which we call its \define{source}
\indexdef{source!of a path constructor}%
and \define{target};
\indexdef{target!of a path constructor}%
for $\lloop$, both source and target are $\base$.

Note that a path constructor such as $\lloop$ generates a \emph{new} inhabitant of an identity type, which is not (at least, not \emph{a priori}) equal to any previously existing such inhabitant.
In particular, $\lloop$ is not \emph{a priori} equal to $\refl{\base}$ (although proving that they are definitely unequal takes a little thought; see \autoref{thm:loop-nontrivial}).
This is what distinguishes $\Sn^1$ from the ordinary inductive type \unit.

There are some important points to be made regarding this generalization.

\index{free!generation of an inductive type}%
First of all, the word ``generation'' should be taken seriously, in the same sense that a group can be freely generated by some set.
In particular, because a higher groupoid comes with \emph{operations} on paths and higher paths, when such an object is ``generated'' by certain constructors, the operations create more paths that do not come directly from the constructors themselves.
For instance, in the higher inductive type $\Sn^1$, the constructor $\lloop$ is not the only nontrivial path from $\base$ to $\base$; we have also ``$\lloop\ct\lloop$'' and ``$\lloop\ct\lloop\ct\lloop$'' and so on, as well as $\opp{\lloop}$, etc., all of which are different.
This may seem so obvious as to be not worth mentioning, but it is a departure from the behavior of ``ordinary'' inductive types, where one can expect to see nothing in the inductive type except what was ``put in'' directly by the constructors.

Secondly, this generation is really \emph{free} generation: higher inductive types do not technically allow us to impose ``axioms'', such as forcing ``$\lloop\ct\lloop$'' to equal $\refl{\base}$.
However, in the world of $\infty$-groupoids,%
\index{.infinity-groupoid@$\infty$-groupoid}
there is little difference between ``free generation'' and ``presentation'',
\index{presentation!of an infinity-groupoid@of an $\infty$-groupoid}%
\index{generation!of an infinity-groupoid@of an $\infty$-groupoid}%
since we can make two paths equal \emph{up to homotopy} by adding a new 2-di\-men\-sion\-al generator relating them (e.g.\ a path $\lloop\ct\lloop = \refl{\base}$ in $\base=\base$).
We do then, of course, have to worry about whether this new generator should satisfy its own ``axioms'', and so on, but in principle any ``presentation'' can be transformed into a ``free'' one by making axioms into constructors.
As we will see, by adding ``truncation constructors'' we can use higher inductive types to express classical notions such as group presentations as well.

Thirdly, even though a higher inductive type contains ``constructors'' which generate \emph{paths in} that type, it is still an inductive definition of a \emph{single} type.
In particular, as we will see, it is the higher inductive type itself which is given a universal property (expressed, as usual, by an induction principle), and \emph{not} its identity types.
The identity type of a higher inductive type retains the usual induction principle of any identity type (i.e.\ path induction), and does not acquire any new induction principle.

Thus, it may be nontrivial to identify the identity types of a higher inductive type in a concrete way, in contrast to how in \autoref{cha:basics} we were able to give explicit descriptions of the behavior of identity types under all the traditional type forming operations.
For instance, are there any paths from $\base$ to $\base$ in $\Sn^1$ which are not simply composites of copies of $\lloop$ and its inverse?
Intuitively, it seems that the answer should be no (and it is), but proving this is not trivial.
Indeed, such questions bring us rapidly to problems such as calculating the homotopy groups of spheres, a long-standing problem in algebraic topology for which no simple formula is known.
Homotopy type theory brings a new and powerful viewpoint to bear on such questions, but it also requires type theory to become as complex as the answers to these questions.

\index{dimension!of path constructors}%
Fourthly, the ``dimension'' of the constructors (i.e.\ whether they output points, paths, paths between paths, etc.)\ does not have a direct connection to which dimensions the resulting type has nontrivial homotopy in.
As a simple example, if an inductive type $B$ has a constructor of type $A\to B$, then any paths and higher paths in $A$ result in paths and higher paths in $B$, even though the constructor is not a ``higher'' constructor at all.
The same thing happens with higher constructors too: having a constructor of type $A\to (\id[B]xy)$ means not only that points of $A$ yield paths from $x$ to $y$ in $B$, but that paths in $A$ yield paths between these paths, and so on.
As we will see, this possibility is responsible for much of the power of higher inductive types.

On the other hand, it is even possible for constructors \emph{without} higher types in their inputs to generate ``unexpected'' higher paths.
For instance, in the 2-dimensional sphere $\Sn^2$ generated by
\symlabel{s2a}
\index{type!2-sphere}%
\begin{itemize}
\item A point $\base:\Sn^2$, and
\item A 2-dimensional path $\surf:\refl{\base} = \refl{\base}$ in ${\base=\base}$,
\end{itemize}
there is a nontrivial \emph{3-dimensional path} from $\refl{\refl{\base}}$ to itself.
Topologists will recognize this path as an incarnation of the \emph{Hopf fibration}.
From a category-theoretic point of view, this is the same sort of phenomenon as the fact mentioned above that $\Sn^1$ contains not only $\lloop$ but also $\lloop\ct\lloop$ and so on: it's just that in a \emph{higher} groupoid, there are \emph{operations} which raise dimension.
Indeed, we saw many of these operations back in \autoref{sec:equality}: the associativity and unit laws are not just properties, but operations, whose inputs are 1-paths and whose outputs are 2-paths.

\index{generation!of a type, inductive|)}%

\vspace*{0pt plus 20ex}

\section{Induction principles and dependent paths}
\label{sec:dependent-paths}

When we describe a higher inductive type such as the circle as being generated by certain constructors, we have to explain what this means by giving rules analogous to those for the basic type constructors from \autoref{cha:typetheory}.
The constructors themselves give the \emph{introduction} rules, but it requires a bit more thought to explain the \emph{elimination} rules, i.e.\ the induction and recursion principles.
In this book we do not attempt to give a general formulation of what constitutes a ``higher inductive definition'' and how to extract the elimination rule from such a definition --- indeed, this is a subtle question and the subject of current research.
Instead we will rely on some general informal discussion and numerous examples.

\index{type!circle}%
\index{recursion principle!for S1@for $\Sn^1$}%
The recursion principle is usually easy to describe: given any type equipped with the same structure with which the constructors equip the higher inductive type in question, there is a function which maps the constructors to that structure.
For instance, in the case of $\Sn^1$, the recursion principle says that given any type $B$ equipped with a point $b:B$ and a path $\ell:b=b$, there is a function $f:\Sn^1\to B$ such that $f(\base)=b$ and $\apfunc f (\lloop) = \ell$.

\index{computation rule!for S1@for $\Sn^1$}%
\index{equality!definitional}%
The latter two equalities are the \emph{computation rules}.
\index{computation rule!for higher inductive types|(}%
\index{computation rule!propositional|(}%
There is, however, a question of whether these computation rules are judgmental\index{judgmental equality} equalities or propositional equalities (paths).
For ordinary inductive types, we had no qualms about making them judgmental, although we saw in \autoref{cha:induction} that making them propositional would still yield the same type up to equivalence.
In the ordinary case, one may argue that the computation rules are really \emph{definitional} equalities, in the intuitive sense described in the Introduction.

\index{equality!judgmental}%
For higher inductive types, this is less clear. 
Moreover, since the operation $\apfunc f$ is not really a fundamental part of the type theory, but something that we \emph{defined} using the induction principle of identity types (and which we might have defined in some other, equivalent, way), it seems inappropriate to refer to it explicitly in a \emph{judgmental} equality.
Judgmental equalities are part of the deductive system, which should not depend on particular choices of definitions that we may make \emph{within} that system.
There are also semantic and implementation issues to consider; see the Notes.

It does seem unproblematic to make the computational rules for the \emph{point} constructors of a higher inductive type judgmental.
In the example above, this means we have $f(\base)\jdeq b$, judgmentally.
This choice facilitates a computational view of higher inductive types.
Moreover, it also greatly simplifies our lives, since otherwise the second computation rule $\apfunc f (\lloop) = \ell$ would not even be well-typed as a propositional equality; we would have to compose one side or the other with the specified identification of $f(\base)$ with $b$.
(Such problems do arise eventually, of course, when we come to talk about paths of higher dimension, but that will not be of great concern to us here.
See also \autoref{sec:hubs-spokes}.)
Thus, we take the computation rules for point constructors to be judgmental, and those for paths and higher paths to be propositional.%
\footnote{In particular, in the language of \autoref{sec:types-vs-sets}, this means that our higher inductive types are a mix of \emph{rules} (specifying how we can introduce such types and their elements, their induction principle, and their computation rules for point constructors) and \emph{axioms} (the computation rules for path constructors, which assert that certain identity types are inhabited by otherwise unspecified terms).
We may hope that eventually, there will be a better type theory in which higher inductive types, like univalence, will be presented using only rules and no axioms.%
\indexfoot{axiom!versus rules}%
\indexfoot{rule!versus axioms}%
}

\begin{rmk}\label{rmk:defid}
Recall that for ordinary inductive types, we regard the computation rules for a recursively defined function as not merely judgmental equalities, but \emph{definitional} ones, and thus we may use the notation $\defeq$ for them.
For instance, the truncated predecessor\index{predecessor!function, truncated} function $p:\nat\to\nat$ is defined by $p(0)\defeq 0$ and $p(\suc(n))\defeq n$.
In the case of higher inductive types, this sort of notation is reasonable for the point constructors (e.g.\ $f(\base)\defeq b$), but for the path constructors it could be misleading, since equalities such as $\ap f \lloop = \ell$ are not judgmental.
Thus, we hybridize the notations, writing instead $\ap f \lloop \defid \ell$ for this sort of ``propositional equality by definition''.
\end{rmk}
\index{computation rule!for higher inductive types|)}%
\index{computation rule!propositional|)}%

\index{type!circle|(}%
\index{induction principle!for S1@for $\Sn^1$}%
Now, what about the the induction principle (the dependent eliminator)?
Recall that for an ordinary inductive type $W$, to prove by induction that $\prd{x:W} P(x)$, we must specify, for each constructor of $W$, an operation on $P$ which acts on the ``fibers'' above that constructor in $W$.
For instance, if $W$ is the natural numbers \nat, then to prove by induction that $\prd{x:\nat} P(x)$, we must specify
\begin{itemize}
\item An element $b:P(0)$ in the fiber over the constructor $0:\nat$, and
\item For each $n:\nat$, a function $P(n) \to P(\suc(n))$.
\end{itemize}
The second can be viewed as a function ``$P\to P$'' lying \emph{over} the constructor $\suc:\nat\to\nat$, generalizing how $b:P(0)$ lies over the constructor $0:\nat$.

By analogy, therefore, to prove that $\prd{x:\Sn^1} P(x)$, we should specify
\begin{itemize}
\item An element $b:P(\base)$ in the fiber over the constructor $\base:\Sn^1$, and
\item A path from $b$ to $b$ ``lying over the constructor $\lloop:\base=\base$''.
\end{itemize}
Note that even though $\Sn^1$ contains paths other than $\lloop$ (such as $\refl{\base}$ and $\lloop\ct\lloop$), we only need to specify a path lying over the constructor \emph{itself}.
This expresses the intuition that $\Sn^1$ is ``freely generated'' by its constructors.

The question, however, is what it means to have a path ``lying over'' another path.
It definitely does \emph{not} mean simply a path $b=b$, since that would be a path in the fiber $P(\base)$ (topologically, a path lying over the \emph{constant} path at $\base$).
Actually, however, we have already answered this question in \autoref{cha:basics}: in the discussion preceding \autoref{lem:mapdep} we concluded that a path from $u:P(x)$ to $v:P(y)$ lying over $p:x=y$ can be represented by a path $\trans p u = v$ in the fiber $P(y)$.
Since we will have a lot of use for such \define{dependent paths}
\index{path!dependent}%
in this chapter, we introduce a special notation for them:
\begin{equation}
  (\dpath P p u v) \defeq (\transfib{P} p u = v).\label{eq:dpath}
\end{equation}

\begin{rmk}
There are other possible ways to define dependent paths.
For instance, instead of $\trans p u = v$ we could consider $u = \trans{(\opp p)}{v}$.
We could also obtain it as a special case of a more general ``heterogeneous equality'',
\index{heterogeneous equality}%
\index{equality!heterogeneous}%
or with a direct definition as an inductive type family.
All these definitions result in equivalent types, so in that sense it doesn't much matter which we pick.
However, choosing $\trans p u = v$ as the definition makes it easiest to conclude other things about dependent paths, such as the fact that $\apdfunc{f}$ produces them, or that we can compute them in particular type families using the transport lemmas in \autoref{sec:computational}.
\end{rmk}

With the notion of dependent paths in hand, we can now state more precisely the induction principle for $\Sn^1$: given $P:\Sn^1\to\type$ and
\begin{itemize}
\item An element $b:P(\base)$, and
\item A path $\ell : \dpath P \lloop b b$,
\end{itemize}
there is a function $f:\prd{x:\Sn^1} P(x)$ such that $f(\base)\jdeq b$ and $\apd f \lloop = \ell$.
As in the non-dependent case, we speak of defining $f$ by $f(\base)\defeq b$ and $\apd f \lloop \defid \ell$.

\begin{rmk}\label{rmk:varies-along}
  When describing an application of this induction principle informally, we regard it as a splitting of the goal ``$P(x)$ for all $x:\Sn^1$'' into two cases, which we will sometimes introduce with phrases such as ``when $x$ is $\base$'' and ``when $x$ varies along $\lloop$'', respectively.
  \index{vary along a path constructor}%
  There is no specific mathematical meaning assigned to ``varying along a path'': it is just a convenient way to indicate the beginning of the corresponding section of a proof; see \autoref{thm:S1-autohtpy} for an example.
\end{rmk}

Topologically, the induction principle for $\Sn^1$ can be visualized as shown in \autoref{fig:topS1ind}.
Given a fibration over the circle (which in the picture is a torus), to define a section of this fibration is the same as to give a point $b$ in the fiber over $\base$ along with a path from $b$ to $b$ lying over $\lloop$.
The way we interpret this type-theoretically, using our definition of dependent paths, is shown in \autoref{fig:ttS1ind}: the path from $b$ to $b$ over $\lloop$ is represented by a path from $\trans \lloop b$ to $b$ in the fiber over $\base$.

\begin{figure}
  \centering
  \begin{tikzpicture}
    \draw (0,0) ellipse (3 and .5);
    \draw (0,3) ellipse (3.5 and 1.5);
    \begin{scope}[yshift=4]
      \clip (-3,3) -- (-1.8,3) -- (-1.8,3.7) -- (1.8,3.7) -- (1.8,3) -- (3,3) -- (3,0) -- (-3,0) -- cycle;
      \draw[clip] (0,3.5) ellipse (2.25 and 1);
      \draw (0,2.5) ellipse (1.7 and .7);
    \end{scope}
    \node (P) at (4.5,3) {$P$};
    \node (S1) at (4.5,0) {$\Sn^1$};
    \draw[->>,thick] (P) -- (S1);
    \node[fill,circle,inner sep=1pt,label={below right:$\base$}] at (0,-.5) {};
    \node at (-2.6,.6) {$\lloop$};
    \node[fill,circle,\OPTblue,inner sep=1pt] (b) at (0,2.3) {};
    \node[\OPTblue] at (-.2,2.1) {$b$};
      \begin{scope}
        \draw[\OPTblue] (b) to[out=180,in=-150] (-2.7,3.5) to[out=30,in=180] (0,3.35);
        \draw[\OPTblue,dotted] (0,3.35) to[out=0,in=175] (1.4,4.35);
        \draw[\OPTblue] (1.4,4.35) to[out=-5,in=90] (2.5,3) to[out=-90,in=0,looseness=.8] (b);
      \end{scope}
      \node[\OPTblue] at (-2.2, 3.3) {$\ell$};
  \end{tikzpicture}
  \caption{The topological induction principle for $\Sn^1$}
  \label{fig:topS1ind}
\end{figure}

\begin{figure}
  \centering
  \begin{tikzpicture}
    \draw (0,0) ellipse (3 and .5);
    \draw (0,3) ellipse (3.5 and 1.5);
    \begin{scope}[yshift=4]
      \clip (-3,3) -- (-1.8,3) -- (-1.8,3.7) -- (1.8,3.7) -- (1.8,3) -- (3,3) -- (3,0) -- (-3,0) -- cycle;
      \draw[clip] (0,3.5) ellipse (2.25 and 1);
      \draw (0,2.5) ellipse (1.7 and .7);
    \end{scope}
    \node (P) at (4.5,3) {$P$};
    \node (S1) at (4.5,0) {$\Sn^1$};
    \draw[->>,thick] (P) -- (S1);
    \node[fill,circle,inner sep=1pt,label={below right:$\base$}] at (0,-.5) {};
    \node at (-2.6,.6) {$\lloop$};
    \node[fill,circle,\OPTblue,inner sep=1pt] (b) at (0,2.3) {};
      \node[\OPTblue] at (-.3,2.3) {$b$};
      \node[fill,circle,\OPTpurple,inner sep=1pt] (tb) at (0,1.8) {};
      \draw[\OPTpurple,dashed] (b) arc (-90:90:2.9 and 0.85) arc (90:270:2.8 and 1.1);
      \begin{scope}
        \clip (b) -- ++(.1,0) -- (.1,1.8) -- ++(-.2,0) -- ++(0,-1) -- ++(3,2) -- ++(-3,0) -- (-.1,2.3) -- cycle;
        \draw[\OPTred,dotted,thick] (.2,2.07) ellipse (.2 and .57);
        \begin{scope}
          \clip (.2,0) rectangle (-2,3);
          \draw[\OPTred,thick] (.2,2.07) ellipse (.2 and .57);
        \end{scope}
      \end{scope}
      \node[\OPTred] at (1,1.2) {$\ell: \trans \lloop b=b$};
  \end{tikzpicture}
  \caption{The type-theoretic induction principle for $\Sn^1$}
  \label{fig:ttS1ind}
\end{figure}

Of course, we expect to be able to prove the recursion principle from the induction principle, by taking $P$ to be a constant type family.
This is in fact the case, although deriving the non-dependent computation rule for $\lloop$ (which refers to $\apfunc f$) from the dependent one (which refers to $\apdfunc f$) is surprisingly a little tricky.

\begin{lem}\label{thm:S1rec}
  \index{recursion principle!for S1@for $\Sn^1$}%
  \index{computation rule!for S1@for $\Sn^1$}%
  If $A$ is a type together with $a:A$ and $p:\id[A]aa$, then there is a
  function $f:\Sn^1\to{}A$ with
  \begin{align*}
    f(\base)&\defeq a \\
    \apfunc f(\lloop)&\defid p.
  \end{align*}
\end{lem}
\begin{proof}
  We would like to apply the induction principle of $\Sn^1$ to the constant type family, $(\lam{x} A): \Sn^1\to \UU$.
  The required hypotheses for this are a point of $(\lam{x} A)(\base) \jdeq A$, which we have (namely $a:A$), and a dependent path in $\dpath {x \mapsto A}{\lloop} a a$, or equivalently $\transfib{x \mapsto A}{\lloop} a = a$.
  This latter type is not the same as the type $\id[A]aa$ where $p$ lives, but it is equivalent to it, because by \autoref{thm:trans-trivial} we have $\transconst{A}{\lloop}{a} : \transfib{x \mapsto A}{\lloop} a= a$.
  Thus, given $a:A$ and $p:a=a$, we can consider the composite
  \[\transconst{A}{\lloop}{a} \ct p:(\dpath {x \mapsto A}\lloop aa).\]
  Applying the induction principle, we obtain $f:\Sn^1\to A$ such that
  \begin{align}
    f(\base) &\jdeq a \qquad\text{and}\label{eq:S1recindbase}\\
    \apdfunc f(\lloop) &= \transconst{A}{\lloop}{a} \ct p.\label{eq:S1recindloop}
  \end{align}
  It remains to derive the equality $\apfunc f(\lloop)=p$.
  However, by \autoref{thm:apd-const}, we have
  \[\apdfunc f(\lloop) = \transconst{A}{\lloop}{f(\base)} \ct \apfunc f(\lloop).\]
  Combining this with~\eqref{eq:S1recindloop} and canceling the occurrences of $\transconstf$ (which are the same by~\eqref{eq:S1recindbase}), we obtain $\apfunc f(\lloop)=p$.
\end{proof}

We also have a corresponding uniqueness principle.

\begin{lem}
  \index{uniqueness!principle, propositional!for functions on the circle}%
  If $A$ is a type and $f,g:\Sn^1\to{}A$ are two maps together with two
  equalities $p,q$:
  \begin{align*}
    p:f(\base)&=_Ag(\base),\\
    q:\map{f}\lloop&=^{\lam{x} x=_Ax}_p\map{g}\lloop.
  \end{align*}
  Then for all $x:\Sn^1$ we have $f(x)=g(x)$.
\end{lem}
\begin{proof}
  This is just the induction principle for the type family $P(x)\defeq(f(x)=g(x))$.
\end{proof}

\index{universal!property!of S1@of $\Sn^1$}%
These two lemmas imply the expected universal property of the circle:

\begin{lem}\label{thm:S1ump}
  For any type $A$ we have a natural equivalence
  \[ (\Sn^1 \to A) \;\eqvsym\;
  \sm{x:A} (x=x).
  \]
\end{lem}
\begin{proof}
  We have a canonical function $f:(\Sn^1 \to A) \to \sm{x:A} (x=x)$ defined by $f(g) \defeq (g(\base),\ap g \lloop)$.
  The induction principle shows that the fibers of $f$ are inhabited, while the uniqueness principle shows that they are mere propositions.
  Hence they are contractible, so $f$ is an equivalence.
\end{proof}

\index{type!circle|)}%

As in \autoref{sec:htpy-inductive}, we can show that the conclusion of \autoref{thm:S1ump} is equivalent to having an induction principle with propositional computation rules.
Other higher inductive types also satisfy lemmas analogous to \autoref{thm:S1rec,thm:S1ump}; we will generally leave their proofs to the reader.
We now proceed to consider many examples.

\section{The interval}
\label{sec:interval}

\index{type!interval|(defstyle}%
\indexsee{interval!type}{type, interval}%
The \define{interval}, which we denote $\interval$, is perhaps an even simpler higher inductive type than the circle.
It is generated by:
\begin{itemize}
\item a point $\izero:\interval$,
\item a point $\ione:\interval$, and
\item a path $\seg : \id[\interval]\izero\ione$.
\end{itemize}
\index{recursion principle!for interval type}%
The recursion principle for the interval says that given a type $B$ along with
\begin{itemize}
\item a point $b_0:B$,
\item a point $b_1:B$, and
\item a path $s:b_0=b_1$,
\end{itemize}
there is a function $f:\interval\to B$ such that $f(\izero)\jdeq b_0$, $f(\ione)\jdeq b_1$, and $\ap f \seg = s$.
\index{induction principle!for interval type}%
The induction principle says that given $P:\interval\to\type$ along with
\begin{itemize}
\item a point $b_0:P(\izero)$,
\item a point $b_1:P(\ione)$, and
\item a path $s:\dpath{P}{\seg}{b_0}{b_1}$,
\end{itemize}
there is a function $f:\prd{x:\interval} P(x)$ such that $f(\izero)\jdeq b_0$, $f(\ione)\jdeq b_1$, and $\apd f \seg = s$.

Regarded purely up to homotopy, the interval is not really interesting:

\begin{lem}
  The type $\interval$ is contractible.
\end{lem}

\begin{proof}
  We prove that for all $x:\interval$ we have $x=_\interval\ione$. In other words we want a
  function $f$ of type $\prd{x:\interval}(x=_\interval\ione)$. We begin to define $f$ in the following way:
  \begin{alignat*}{2}
    f(\izero)&\defeq \seg  &:\izero&=_\interval\ione,\\
    f(\ione)&\defeq \refl\ione &:\ione &=_\interval\ione.
  \end{alignat*}
  It remains to define $\apd{f}\seg$, which must have type $\seg =_\seg^{\lam{x} x=_\interval\ione}\refl \ione$.
  By definition this type is $\trans\seg\seg=_{\ione=_\interval\ione}\refl\ione$, which in turn is equivalent to $\rev\seg\ct\seg=\refl\ione$.
  But there is a canonical element of that type, namely the proof that path inverses are in fact inverses.
\end{proof}

However, type-theoretically the interval does still have some interesting features, just like the topological interval in classical homotopy theory.
For instance, it enables us to give an easy proof of function extensionality.
(Of course, as in \autoref{sec:univalence-implies-funext}, for the duration of the following proof we suspend our overall assumption of the function extensionality axiom.)

\begin{lem}\label{thm:interval-funext}
  \index{function extensionality!proof from interval type}%
  If $f,g:A\to{}B$ are two functions such that $f(x)=g(x)$ for every $x:A$, then
  $f=g$ in the type $A\to{}B$.
\end{lem}

\begin{proof}
  Let's call the proof we have $p:\prd{x:A}(f(x)=g(x))$. For all $x:A$ we define
  a function $\widetilde{p}_x:\interval\to{}B$ by
  \begin{align*}
    \widetilde{p}_x(\izero) &\defeq f(x), \\
    \widetilde{p}_x(\ione) &\defeq g(x), \\
    \map{\widetilde{p}_x}\seg &\defid p(x).
  \end{align*}
  We now define $q:\interval\to(A\to{}B)$ by
  \[q(i)\defeq(\lam{x} \widetilde{p}_x(i))\]
  Then $q(\izero)$ is the function $\lam{x} \widetilde{p}_x(\izero)$, which is equal to $f$ because $\widetilde{p}_x(\izero)$ is defined by $f(x)$.
  Similarly, we have $q(\ione)=g$, and hence
  \[\map{q}\seg:f=_{(A\to{}B)}g \qedhere\]
\end{proof}

\index{type!interval|)}%

\section{Circles and spheres}
\label{sec:circle}

\index{type!circle|(}%
We have already discussed the circle $\Sn^1$ as the higher inductive type generated by
\begin{itemize}
\item A point $\base:\Sn^1$, and
\item A path $\lloop : {\id[\Sn^1]\base\base}$.
\end{itemize}
\index{induction principle!for S1@for $\Sn^1$}%
Its induction principle says that given $P:\Sn^1\to\type$ along with $b:P(\base)$ and $\ell :\dpath P \lloop b b$, we have $f:\prd{x:\Sn^1} P(x)$ with $f(\base)\jdeq b$ and $\apd f \lloop = \ell$.
Its non-dependent recursion principle says that given $B$ with $b:B$ and $\ell:b=b$, we have $f:\Sn^1\to B$ with $f(\base)\jdeq b$ and $\ap f \lloop = \ell$.

We observe that the circle is nontrivial.

\begin{lem}\label{thm:loop-nontrivial}
  $\lloop\neq\refl{\base}$.
\end{lem}
\begin{proof}
  Suppose that $\lloop=\refl{\base}$.
  Then since for any type $A$ with $x:A$ and $p:x=x$, there is a function $f:\Sn^1\to A$ defined by $f(\base)\defeq x$ and $\ap f \lloop \defid p$, we have
  \[p = f(\lloop) = f(\refl{\base}) = \refl{x}.\]
  But this implies that every type is a set, which as we have seen is not the case (see \autoref{thm:type-is-not-a-set}).
\end{proof}

The circle also has the following interesting property, which is useful as a source of counterexamples.

\begin{lem}\label{thm:S1-autohtpy}
  There exists an element of $\prd{x:\Sn^1} (x=x)$ which is not equal to $x\mapsto \refl{x}$.
\end{lem}
\begin{proof}
  We define $f:\prd{x:\Sn^1} (x=x)$ by $\Sn^1$-induction.
  When $x$ is $\base$, we let $f(\base)\defeq \lloop$.
  Now when $x$ varies along $\lloop$ (see \autoref{rmk:varies-along}), we must show that $\transfib{x\mapsto x=x}{\lloop}{\lloop} = \lloop$.
  However, in \autoref{sec:compute-paths} we observed that $\transfib{x\mapsto x=x}{p}{q} = \opp{p} \ct q \ct p$, so what we have to show is that $\opp{\lloop} \ct \lloop \ct \lloop = \lloop$.
  But this is clear by canceling an inverse.

  To show that $f\neq (x\mapsto \refl{x})$, it suffices by function extensionality to show that $f(\base) \neq \refl{\base}$.
  But $f(\base)=\lloop$, so this is just the previous lemma.
\end{proof}

For instance, this enables us to extend \autoref{thm:type-is-not-a-set} by showing that any universe which contains the circle cannot be a 1-type.

\begin{cor}
  If the type $\Sn^1$ belongs to some universe \type, then \type is not a 1-type.
\end{cor}
\begin{proof}
  The type $\Sn^1=\Sn^1$ in \type is, by univalence, equivalent to the type $\eqv{\Sn^1}{\Sn^1}$ of auto\-equivalences of $\Sn^1$, so it suffices to show that $\eqv{\Sn^1}{\Sn^1}$ is not a set.
  \index{automorphism!of S1@of $\Sn^1$}%
  For this, it suffices to show that its equality type $\id[(\eqv{\Sn^1}{\Sn^1})]{\idfunc[\Sn^1]}{\idfunc[\Sn^1]}$ is not a mere proposition.
  Since being an equivalence is a mere proposition, this type is equivalent to $\id[(\Sn^1\to\Sn^1)]{\idfunc[\Sn^1]}{\idfunc[\Sn^1]}$.
  But by function extensionality, this is equivalent to $\prd{x:\Sn^1} (x=x)$, which as we have seen in \autoref{thm:S1-autohtpy} contains two unequal elements.
\end{proof}

\index{type!circle|)}%

\index{type!2-sphere|(}%
\indexsee{sphere type}{type, sphere}%
We have also mentioned that the 2-sphere $\Sn^2$ should be the higher inductive type generated by
\symlabel{s2b}
\begin{itemize}
\item A point $\base:\Sn^2$, and
\item A 2-dimensional path $\surf:\refl{\base} = \refl{\base}$ in ${\base=\base}$.
\end{itemize}
\index{recursion principle!for S2@for $\Sn^2$}%
The recursion principle for $\Sn^2$ is not hard: it says that given $B$ with $b:B$ and $s:\refl b = \refl b$, we have $f:\Sn^2\to B$ with $f(\base)\jdeq b$ and $\aptwo f \surf = s$.
Here by ``$\aptwo f \surf$'' we mean an extension of the functorial action of $f$ to two-dimensional paths, which can be stated precisely as follows.

\begin{lem}\label{thm:ap2}
  Given $f:A\to B$ and $x,y:A$ and $p,q:x=y$, and $r:p=q$, we have a path $\aptwo f r : \ap f p = \ap f q$.
\end{lem}
\begin{proof}
  By path induction, we may assume $p\jdeq q$ and $r$ is reflexivity.
  But then we may define $\aptwo f {\refl p} \defeq \refl{\ap f p}$.
\end{proof}

In order to state the general induction principle, we need a version of this lemma for dependent functions, which in turn requires a notion of dependent two-dimensional paths.
As before, there are many ways to define such a thing; one is by way of a two-dimensional version of transport.

\begin{lem}\label{thm:transport2}
  Given $P:A\to\type$ and $x,y:A$ and $p,q:x=y$ and $r:p=q$, for any $u:P(x)$ we have $\transtwo r u : \trans p u = \trans q u$.
\end{lem}
\begin{proof}
  By path induction.
\end{proof}

Now suppose given $x,y:A$ and $p,q:x=y$ and $r:p=q$ and also points $u:P(x)$ and $v:P(y)$ and dependent paths $h:\dpath P p u v$ and $k:\dpath P q u v$.
By our definition of dependent paths, this means $h:\trans p u = v$ and $k:\trans q u = v$.
Thus, it is reasonable to define the type of dependent 2-paths over $r$ to be
\[ (\dpath P r h k )\defeq (h = \transtwo r u \ct k). \]
We can now state the dependent version of \autoref{thm:ap2}.

\begin{lem}\label{thm:apd2}
  Given $P:A\to\type$ and $x,y:A$ and $p,q:x=y$ and $r:p=q$ and a function $f:\prd{x:A} P(x)$, we have
  $\apdtwo f r : \dpath P r {\apd f p}{\apd f q}$.
\end{lem}
\begin{proof}
  Path induction.
\end{proof}

\index{induction principle!for S2@for $\Sn^2$}%
Now we can state the induction principle for $\Sn^2$: given $P:\Sn^2\to P$ with $b:P(\base)$ and $s:\dpath P \surf {\refl b}{\refl b}$, there is a function $f:\prd{x:\Sn^2} P(x)$ such that $f(\base)\jdeq b$ and $\apdtwo f \surf = s$.

\index{type!2-sphere|)}%

Of course, this explicit approach gets more and more complicated as we go up in dimension.
Thus, if we want to define $n$-spheres for all $n$, we need some more systematic idea.
One approach is to work with $n$-dimensional loops\index{loop!n-@$n$-} directly, rather than general $n$-dimensional paths.\index{path!n-@$n$-}

\index{type!pointed}%
Recall from \autoref{sec:equality} the definitions of \emph{pointed types} $\type_*$, and the $n$-fold loop space\index{loop space!iterated} $\Omega^n : \type_* \to \type_*$
(\cref{def:pointedtype,def:loopspace}).  Now we can define the
$n$-sphere $\Sn^n$ to be the higher inductive type generated by
\index{type!n-sphere@$n$-sphere}%
\begin{itemize}
\item A point $\base:\Sn^n$, and
\item An $n$-loop $\lloop_n : \Omega^n(\Sn^n,\base)$.
\end{itemize}
In order to write down the induction principle for this presentation, we would need to define a notion of ``dependent $n$-loop\indexdef{loop!dependent n-@dependent $n$-}'', along with the action of dependent functions on $n$-loops.
We leave this to the reader (see \autoref{ex:nspheres}); in the next section we will discuss a different way to define the spheres that is sometimes more tractable.

\section{Suspensions}
\label{sec:suspension}

\indexsee{type!suspension of}{suspension}%
\index{suspension|(defstyle}%
The \define{suspension} of a type $A$ is the universal way of making the points of $A$ into paths (and hence the paths in $A$ into 2-paths, and so on).
It is a type $\susp A$ defined by the following generators:\footnote{There is an unfortunate clash of notation with dependent pair types, which of course are also written with a $\Sigma$.
  However, context usually disambiguates.}
\begin{itemize}
\item a point $\north:\susp A$,
\item a point $\south:\susp A$, and
\item a function $\merid:A \to (\id[\susp A]\north\south)$.
\end{itemize}
The names are intended to suggest a ``globe'' of sorts, with a north pole, a south pole, and an $A$'s worth of meridians
\indexdef{pole}%
\indexdef{meridian}%
from one to the other.
Indeed, as we will see, if $A=\Sn^1$, then its suspension is equivalent to the surface of an ordinary sphere, $\Sn^2$.

\index{recursion principle!for suspension}%
The recursion principle for $\susp A$ says that given a type $B$ together with
\begin{itemize}
\item points $n,s:B$ and
\item a function $m:A \to (n=s)$,
\end{itemize}
we have a function $f:\susp A \to B$ such that $f(\north)\jdeq n$ and $f(\south)\jdeq s$, and for all $a:A$ we have $\ap f {\merid(a)} = m(a)$.
\index{induction principle!for suspension}%
Similarly, the induction principle says that given $P:\susp A \to \type$ together with
\begin{itemize}
\item a point $n:P(\north)$,
\item a point $s:P(\south)$, and
\item for each $a:A$, a path $m(a):\dpath P{\merid(a)}ns$,
\end{itemize}
there exists a function $f:\prd{x:\susp A} P(x)$ such that $f(\north)\jdeq n$ and $f(\south)\jdeq s$ and for each $a:A$ we have $\apd f {\merid(a)} = m(a)$.

Our first observation about suspension is that it gives another way to define the circle.

\begin{lem}\label{thm:suspbool}
  \index{type!circle}%
  $\eqv{\susp\bool}{\Sn^1}$.
\end{lem}
\begin{proof}
  Define $f:\susp\bool\to\Sn^1$ by recursion such that $f(\north)\defeq \base$ and $f(\south)\defeq\base$, while $\ap f{\merid(\bfalse)}\defid\lloop$ but $\ap f{\merid(\btrue)} \defid \refl{\base}$.
  Define $g:\Sn^1\to\susp\bool$ by recursion such that $g(\base)\defeq \north$ and $\ap g \lloop \defid \merid(\bfalse) \ct \opp{\merid(\btrue)}$.
  We now show that $f$ and $g$ are quasi-inverses.

  First we show by induction that $g(f(x))=x$ for all $x:\susp \bool$.
  If $x\jdeq\north$, then $g(f(\north)) \jdeq g(\base)\jdeq \north$, so we have $\refl{\north} : g(f(\north))=\north$.
  If $x\jdeq\south$, then $g(f(\south)) \jdeq g(\base)\jdeq \north$, and we choose the equality $\merid(\btrue) : g(f(\south)) = \south$.
  It remains to show that for any $y:\bool$, these equalities are preserved as $x$ varies along $\merid(y)$, which is to say that when $\refl{\north}$ is transported along $\merid(y)$ it yields $\merid(\btrue)$.
  By transport in path spaces and pulled back fibrations, this means we are to show that
  \[ \opp{\ap g {\ap f {\merid(y)}}} \ct \refl{\north} \ct \merid(y) = \merid(\btrue). \]
  Of course, we may cancel $\refl{\north}$.
  Now by \bool-induction, we may assume either $y\jdeq \bfalse$ or $y\jdeq \btrue$.
  If $y\jdeq \bfalse$, then we have
  \begin{align*}
    \opp{\ap g {\ap f {\merid(\bfalse)}}} \ct \merid(\bfalse)
    &= \opp{\ap g {\lloop}} \ct \merid(\bfalse)\\
    &= \opp{(\merid(\bfalse) \ct \opp{\merid(\btrue)})} \ct \merid(\bfalse)\\
    &= \merid(\btrue) \ct \opp{\merid(\bfalse)} \ct \merid(\bfalse)\\
    &= \merid(\btrue)
  \end{align*}
  while if $y\jdeq \btrue$, then we have
  \begin{align*}
    \opp{\ap g {\ap f {\merid(\btrue)}}} \ct \merid(\btrue)
    &= \opp{\ap g {\refl{\base}}} \ct \merid(\btrue)\\
    &= \opp{\refl{\north}} \ct \merid(\btrue)\\
    &= \merid(\btrue).
  \end{align*}
  Thus, for all $x:\susp \bool$, we have $g(f(x))=x$.

  Now we show by induction that $f(g(x))=x$ for all $x:\Sn^1$.
  If $x\jdeq \base$, then $f(g(\base))\jdeq f(\north)\jdeq\base$, so we have $\refl{\base} : f(g(\base))=\base$.
  It remains to show that this equality is preserved as $x$ varies along $\lloop$, which is to say that it is transported along $\lloop$ to itself.
  Again, by transport in path spaces and pulled back fibrations, this means to show that
  \[ \opp{\ap f {\ap g {\lloop}}} \ct \refl{\base} \ct \lloop = \refl{\base}.\]
  However, we have
  \begin{align*}
    \ap f {\ap g {\lloop}} &= \ap f {\merid(\bfalse) \ct \opp{\merid(\btrue)}}\\
    &= \ap f {\merid(\bfalse)} \ct \opp{\ap f {\merid(\btrue)}}\\
    &= \lloop \ct \refl{\base}
  \end{align*}
  so this follows easily.
\end{proof}

Topologically, the two-point space \bool is also known as the \emph{0-dimensional sphere}, $\Sn^0$.
(For instance, it is the space of points at distance $1$ from the origin in $\mathbb{R}^1$, just as the topological 1-sphere is the space of points at distance $1$ from the origin in $\mathbb{R}^2$.)
Thus, \autoref{thm:suspbool} can be phrased suggestively as $\eqv{\susp\Sn^0}{\Sn^1}$.
\index{type!n-sphere@$n$-sphere|defstyle}%
\indexsee{n-sphere@$n$-sphere}{type, $n$-sphere}%
In fact, this pattern continues: we can define all the spheres inductively by
\begin{equation}\label{eq:Snsusp}
  \Sn^0 \defeq \bool
  \qquad\text{and}\qquad
  \Sn^{n+1} \defeq \susp \Sn^n.
\end{equation}
We can even start one dimension lower by defining $\Sn^{-1}\defeq \emptyt$, and observe that $\eqv{\susp\emptyt}{\bool}$.

To prove carefully that this agrees with the definition of $\Sn^n$ from the previous section would require making the latter more explicit.
However, we can show that the recursive definition has the same universal property that we would expect the other one to have.
If $(A,a_0)$ and $(B,b_0)$ are pointed types (with basepoints often left implicit), let $\Map_*(A,B)$ denote the type of based maps:
\index{based map}
\symlabel{based-maps}
\[ \Map_*(A,B) \defeq \sm{f:A\to B} (f(a_0)=b_0). \]
Note that any type $A$ gives rise to a pointed type $A_+ \defeq A+\unit$ with basepoint $\inr(\ttt)$; this is called \emph{adjoining a disjoint basepoint}.
\indexdef{basepoint!adjoining a disjoint}%
\index{disjoint!basepoint}%
\index{adjoining a disjoint basepoint}%

\begin{lem}
  For a type $A$ and a pointed type $(B,b_0)$, we have
  \[ \eqv{\Map_*(A_+,B)}{(A\to B)} \]
\end{lem}
Note that on the right we have the ordinary type of \emph{unbased} functions from $A$ to $B$.
\begin{proof}
  From left to right, given $f:A_+ \to B$ with $p:f(\inr(\ttt)) = b_0$, we have $f\circ \inl : A \to B$.
  And from right to left, given $g:A\to B$ we define $g':A_+ \to B$ by $g'(\inl(a))\defeq g(a)$ and $g'(\inr(u)) \defeq b_0$.
  We leave it to the reader to show that these are quasi-inverse operations.
\end{proof}

In particular, note that $\eqv{\bool}{\unit_+}$.
Thus, for any pointed type $B$ we have
\[{\Map_*(\bool,B)} \eqvsym {(\unit \to B)}\eqvsym B.\]
Now recall that the loop space\index{loop space} operation $\Omega$ acts on pointed types, with definition $\Omega(A,a_0) \defeq (\id[A]{a_0}{a_0},\refl{a_0})$.
We can also make the suspension $\susp$ act on pointed types, by $\susp(A,a_0)\defeq (\susp A,\north)$.

\begin{lem}\label{lem:susp-loop-adj}
  \index{universal!property!of suspension}%
  For pointed types $(A,a_0)$ and $(B,b_0)$ we have
  \[ \eqv{\Map_*(\susp A, B)}{\Map_*(A,\Omega B)}.\]
\end{lem}
\begin{proof}
  From left to right, given $f:\susp A \to B$ with $p:f(\north) = b_0$, we define $g:A \to \Omega B$ by
  \[g(a) \defeq \opp p \ct \ap f{\merid(a) \ct \opp{\merid(a_0)}} \ct p.\]
  Then we have
  \begin{align*}
    g(a_0) &\jdeq \opp p \ct \ap f{\merid(a_0) \ct \opp{\merid(a_0)}} \ct p\\
    &= \opp p \ct \ap f{\refl{\north}} \ct p\\
    &= \opp p \ct p\\
    &= \refl{b_0}.
  \end{align*}
  Thus, denoting this path by $q:g(a_0)=\refl{b_0}$, we have $(g,q):\Map_*(A,\Omega B)$.

  On the other hand, from right to left, given $g:A\to \Omega B$ and $q:g(a_0)=\refl{b_0}$, we define $f:\susp A \to B$ by $\susp$-recursion, such that $f(\north)\defeq b_0$ and $f(\south)\defeq b_0$ and
  \[ \ap f {\merid(a)} \defid g(a). \]
  Then we can simply take $p$ to be $\refl{b_0} : f(\north)= b_0$.

  Now given $(f,p)$, by passing back and forth we obtain $(f',p')$ where $f'$ is defined by $f'(\north)\jdeq b_0$ and $f'(\south)\jdeq b_0$ and
  \[ \ap {f'} {\merid(a)} = \opp p \ct \ap f{\merid(a) \ct \opp{\merid(a_0)}} \ct p, \]
  while $p' \jdeq \refl{b_0}$.
  To show $f=f'$, by function extensionality it suffices to show $f(x)=f'(x)$ for all $x:\susp A$, so we can use the induction principle of suspension.
  First, we have
  \begin{equation}
    f(\north) \overset{p}{=} b_0 \jdeq f'(\north). \label{eq:ffprime-north}
  \end{equation}
  Second, we have
  \[\xymatrix@C=4pc{ f(\south) \ar@{=}[r]^-{\opp{\ap f {\merid(a_0)}}} & f(\north) \overset{\smash p}{=} b_0 \jdeq f'(\south).}\]
  And thirdly, as $x$ varies along $\merid(a)$ we must show that the following diagram of paths commutes (invoking the definition of $\ap{f'}{\merid(a)}$):
  \[ \xymatrix{
    f(\north) \ar@{=}[rrr]^-{p} \ar@{=}[ddd]_{f(\merid(a))} &&&
    b_0 \ar@{=}[r]^-{\refl{}} &
    f'(\north) \ar@{=}[d]^{\opp p}\\
    &&&& f(\north) \ar@{=}[d]^{\ap f{\merid(a) \ct \opp{\merid(a_0)}}}\\
    &&&& f(\north) \ar@{=}[d]^p\\
    f(\south) \ar@{=}[rr]_-{\opp{\ap f {\merid(a_0)}}} &&
    f(\north) \ar@{=}[r]_-p &
    b_0 \ar@{=}[r]_-{\refl{}} &
    f'(\south) }
  \]
  This is clear.
  Thus, to show that $(f,p)=(f',p')$, it remains only to show that $p$ is identified with $p'$ when transported along this equality $f=f'$.
  Since the type of $p$ is $f(\north)=b_0$, this means essentially that when $p$ is composed on the left with the inverse of the equality~\eqref{eq:ffprime-north}, it becomes $p'$.
  But this is obvious, since~\eqref{eq:ffprime-north} is just $p$ itself, while $p'$ is reflexivity.

  On the other side, suppose given $(g,q)$.
  By passing back and forth we obtain $(g',q')$ with
  \begin{align*}
    g'(a) &= \opp{\refl{b_0}} \ct g(a) \ct \opp{g(a_0)} \ct \refl{b_0}\\
    &= g(a) \ct \opp{g(a_0)}\\
    &= g(a)
  \end{align*}
  using $q:g(a_0) = \refl{b_0}$ in the last equality.
  Thus, $g'=g$ by function extensionality, so it remains to show that when transported along this equality $q$ is identified with $q'$.
  At $a_0$, the induced equality $g(a_0)=g'(a_0)$ consists essentially of $q$ itself, while the definition of $q'$ involves only canceling inverses and reflexivities.
  Thus, some tedious manipulations of naturality finish the proof.
\end{proof}

\index{type!n-sphere@$n$-sphere|defstyle}%
In particular, for the spheres defined as in~\eqref{eq:Snsusp} we have
\index{universal!property!of Sn@of $\Sn^n$}%
\[ \Map_*(\Sn^n,B) \eqvsym \Map_*(\Sn^{n-1}, \Omega B) \eqvsym \cdots \eqvsym \Map_*(\bool,\Omega^n B) \eqvsym \Omega^n B. \]
Thus, these spheres $\Sn^n$ have the universal property that we would expect from the spheres defined directly in terms of $n$-fold loop spaces\index{loop space!iterated} as in \autoref{sec:circle}.

\index{suspension|)}%

\section{Cell complexes}
\label{sec:cell-complexes}

\index{cell complex|(defstyle}%
\index{CW complex|(defstyle}%
In classical topology, a \emph{cell complex} is a space obtained by successively attaching discs along their boundaries.
It is called a \emph{CW complex} if the boundary of an $n$-dimensional disc\index{disc} is constrained to lie in the discs of dimension strictly less than $n$ (the $(n-1)$-skeleton).\index{skeleton!of a CW-complex}

Any finite CW complex can be presented as a higher inductive type, by turning $n$-dimensional discs into $n$-dimensional paths and partitioning the image of the attaching\index{attaching map} map into a source\index{source!of a path constructor} and a target\index{target!of a path constructor}, with each written as a composite of lower dimensional paths.
Our explicit definitions of $\Sn^1$ and $\Sn^2$ in \autoref{sec:circle} had this form.

\index{torus}%
Another example is the torus $T^2$, which is generated by:
\begin{itemize}
\item a point $b:T^2$,
\item a path $p:b=b$,
\item another path $q:b=b$, and
\item a 2-path $t: p\ct q = q \ct p$.
\end{itemize}
Perhaps the easiest way to see that this is a torus is to start with a rectangle, having four corners $a,b,c,d$, four edges $p,q,r,s$, and an interior which is manifestly a 2-path $t$ from $p\ct q$ to $r\ct s$:
\begin{equation*}
  \xymatrix{
      a\ar@{=}[r]^p\ar@{=}[d]_r \ar@{}[dr]|{\Downarrow t} &
      b\ar@{=}[d]^q\\
      c\ar@{=}[r]_s &
      d
      }
\end{equation*}
Now identify the edge $r$ with $q$ and the edge $s$ with $p$, resulting in also identifying all four corners.
Topologically, this identification can be seen to produce a torus.

\index{induction principle!for torus}%
\index{torus!induction principle for}%
The induction principle for the torus is the trickiest of any we've written out so far.
Given $P:T^2\to\type$, for a section $\prd{x:T^2} P(x)$ we require
\begin{itemize}
\item a point $b':P(b)$,
\item a path $p' : \dpath P p b b$,
\item a path $q' : \dpath P q b b$, and
\item a 2-path $t'$ between the ``composites'' $p'\ct q'$ and $q'\ct p'$, lying over $t$.
\end{itemize}
In order to make sense of this last datum, we need a composition operation for dependent paths, but this is not hard to define.
Then the induction principle gives a function $f:\prd{x:T^2} P(x)$ such that $f(b)\jdeq b'$ and $\apd f {p} = p'$ and $\apd f {q} = q'$ and something like ``$\apdtwo f t = t'$''.
However, this is not well-typed as it stands, firstly because the equalities $\apd f {p} = p'$ and $\apd f {q} = q'$ are not judgmental, and secondly because $\apdfunc f$ only preserves path concatenation up to homotopy.
We leave the details to the reader (see \autoref{ex:torus}).

Of course, another definition of the torus is $T^2 \defeq \Sn^1 \times \Sn^1$ (in \autoref{ex:torus-s1-times-s1} we ask the reader to verify the equivalence of the two).
\index{Klein bottle}%
\index{projective plane}%
The cell-complex definition, however, generalizes easily to other spaces without such descriptions, such as the Klein bottle, the projective plane, etc.
But it does get increasingly difficult to write down the induction principles, requiring us to define notions of dependent $n$-paths and of $\apdfunc{}$ acting on $n$-paths.
Fortunately, once we have the spheres in hand, there is a way around this.

\section{Hubs and spokes}
\label{sec:hubs-spokes}

\indexsee{spoke}{hub and spoke}%
\index{hub and spoke|(defstyle}%

In topology, one usually speaks of building CW complexes by attaching $n$-dimensional discs along their $(n-1)$-dimensional boundary spheres.
\index{attaching map}%
However, another way to express this is by gluing in the \emph{cone}\index{cone!of a sphere} on an $(n-1)$-dimensional sphere.
That is, we regard a disc\index{disc} as consisting of a cone point (or ``hub''), with meridians
\index{meridian}%
(or ``spokes'') connecting that point to every point on the boundary, continuously, as shown in \autoref{fig:hub-and-spokes}.

\begin{figure}
  \centering
  \begin{tikzpicture}
    \draw (0,0) circle (2cm);
    \foreach \x in {0,20,...,350}
      \draw[\OPTblue] (0,0) -- (\x:2cm);
    \node[\OPTblue,circle,fill,inner sep=2pt] (hub) at (0,0) {};
  \end{tikzpicture}
  \caption{A 2-disc made out of a hub and spokes}
  \label{fig:hub-and-spokes}
\end{figure}

We can use this idea to express higher inductive types containing $n$-dimensional path-con\-struc\-tors for $n>1$ in terms of ones containing only 1-di\-men\-sion\-al path-con\-struc\-tors.
The point is that we can obtain an $n$-dimensional path as a continuous family of 1-dimensional paths parametrized by an $(n-1)$-di\-men\-sion\-al object.
The simplest $(n-1)$-dimensional object to use is the $(n-1)$-sphere, although in some cases a different one may be preferable.
(Recall that we were able to define the spheres in \autoref{sec:suspension} inductively using suspensions, which involve only 1-dimensional path constructors.
Indeed, suspension can also be regarded as an instance of this idea, since it involves a family of 1-dimensional paths parametrized by the type being suspended.)

\index{torus}
For instance, the torus $T^2$ from the previous section could be defined instead to be generated by:
\begin{itemize}
\item a point $b:T^2$,
\item a path $p:b=b$,
\item another path $q:b=b$,
\item a point $h:T^2$, and
\item for each $x:\Sn^1$, a path $s(x) : f(x)=h$, where $f:\Sn^1\to T^2$ is defined by $f(\base)\defeq b$ and $\ap f \lloop \defid p \ct q \ct \opp p \ct \opp q$.
\end{itemize}
The induction principle for this version of the torus says that given $P:T^2\to\type$, for a section $\prd{x:T^2} P(x)$ we require
\begin{itemize}
\item a point $b':P(b)$,
\item a path $p' : \dpath P p b b$,
\item a path $q' : \dpath P q b b$,
\item a point $h':P(h)$, and
\item for each $x:\Sn^1$, a path $\dpath {P}{s(x)}{g(x)}{h'}$, where $g:\prd{x:\Sn^1} P(f(x))$ is defined by $g(\base)\defeq b'$ and $\apd g \lloop \defid p' \ct q' \ct \opp{(p')} \ct \opp{(q')}$.
\end{itemize}
Note that there is no need for dependent 2-paths or $\apdtwofunc{}$.
We leave it to the reader to write out the computation rules.

\begin{rmk}\label{rmk:spokes-no-hub}
One might question the need for introducing the hub point $h$; why couldn't we instead simply add paths continuously relating the boundary of the disc to a point \emph{on} that boundary, as shown in \autoref{fig:spokes-no-hub}?
This does work, but not as well.
For if, given some $f:\Sn^1 \to X$, we give a path constructor connecting each $f(x)$ to $f(\base)$, then what we end up with is more like the picture in \autoref{fig:spokes-no-hub-ii} of a cone whose vertex is twisted around and glued to some point on its base.
The problem is that the specified path from $f(\base)$ to itself may not be reflexivity.
We could add a 2-dimensional path constructor ensuring this, but using a separate hub avoids the need for any path constructors of dimension above~$1$.
\end{rmk}

\begin{figure}
  \centering
  \begin{minipage}{2in}
    \begin{center}
      \begin{tikzpicture}
        \draw (0,0) circle (2cm);
        \clip (0,0) circle (2cm);
        \foreach \x in {0,15,...,165}
        \draw[\OPTblue] (0,-2cm) -- (\x:4cm);
      \end{tikzpicture}
    \end{center}
    \caption{Hubless spokes}
    \label{fig:spokes-no-hub}
  \end{minipage}
  \qquad
  \begin{minipage}{2in}
    \begin{center}
      \begin{tikzpicture}[xscale=1.3]
        \draw (0,0) arc (-90:90:.7cm and 2cm) ;
        \draw[dashed] (0,4cm) arc (90:270:.7cm and 2cm) ;
        \draw[\OPTblue] (0,0) to[out=90,in=0] (-1,1) to[out=180,in=180] (0,0);
        \draw[\OPTblue] (0,4cm) to[out=180,in=180,looseness=2] (0,0);
        \path (0,0) arc (-90:-60:.7cm and 2cm) node (a) {};
        \draw[\OPTblue] (a.center) to[out=120,in=10] (-1.2,1.2) to[out=190,in=180] (0,0);
        \path (0,0) arc (-90:-30:.7cm and 2cm) node (b) {};
        \draw[\OPTblue] (b.center) to[out=150,in=20] (-1.4,1.4) to[out=200,in=180] (0,0);
        \path (0,0) arc (-90:0:.7cm and 2cm) node (c) {};
        \draw[\OPTblue] (c.center) to[out=180,in=30] (-1.5,1.5) to[out=210,in=180] (0,0);
        \path (0,0) arc (-90:30:.7cm and 2cm) node (d) {};
        \draw[\OPTblue] (d.center) to[out=190,in=50] (-1.7,1.7) to[out=230,in=180] (0,0);
        \path (0,0) arc (-90:60:.7cm and 2cm) node (e) {};
        \draw[\OPTblue] (e.center) to[out=200,in=70] (-2,2) to[out=250,in=180] (0,0);
        \clip (0,0) to[out=90,in=0] (-1,1) to[out=180,in=180] (0,0);
        \draw (0,4cm) arc (90:270:.7cm and 2cm) ;
      \end{tikzpicture}
    \end{center}
    \caption{Hubless spokes, II}
    \label{fig:spokes-no-hub-ii}
  \end{minipage}
\end{figure}

\begin{rmk}
  \index{computation rule!propositional}%
  Note also that this ``translation'' of higher paths into 1-paths does not preserve judgmental computation rules for these paths, though it does preserve propositional ones.
\end{rmk}

\index{cell complex|)}%
\index{CW complex|)}%

\index{hub and spoke|)}%

\section{Pushouts}
\label{sec:colimits}

\index{type!limit}%
\index{type!colimit}%
\index{limit!of types}%
\index{colimit!of types}%
From a category-theoretic point of view, one of the important aspects of any foundational system is the ability to construct limits and colimits.
In set-theoretic foundations, these are limits and colimits of sets, whereas in our case they are limits and colimits of \emph{types}.
We have seen in \autoref{sec:universal-properties} that cartesian product types have the correct universal property of a categorical product of types, and in \autoref{ex:coprod-ump} that coproduct types likewise have their expected universal property.

As remarked in \autoref{sec:universal-properties}, more general limits can be constructed using identity types and $\Sigma$-types, e.g.\ the pullback\index{pullback} of $f:A\to C$ and $g:B\to C$ is $\sm{a:A}{b:B} (f(a)=g(b))$ (see \autoref{ex:pullback}).
However, more general \emph{colimits} require identifying elements coming from different types, for which higher inductives are well-adapted.
Since all our constructions are homotopy-invariant, all our colimits are necessarily \emph{homotopy colimits}, but we drop the ubiquitous adjective in the interests of concision.

In this section we discuss \emph{pushouts}, as perhaps the simplest and one of the most useful colimits.
Indeed, one expects all finite colimits (for a suitable homotopical definition of ``finite'') to be constructible from pushouts and finite coproducts.
It is also possible to give a direct construction of more general colimits using higher inductive types, but this is somewhat technical, and also not completely satisfactory since we do not yet have a good fully general notion of homotopy coherent diagrams.

\indexsee{type!pushout of}{pushout}%
\index{pushout|(defstyle}%
\index{span}%
Suppose given a span of types and functions:
\[\Ddiag=\;\vcenter{\xymatrix{C \ar^g[r] \ar_f[d] & B \\ A & }}\]
The \define{pushout} of this span is the higher inductive type $A\sqcup^CB$ presented by
\begin{itemize}
\item a function $\inl:A\to A\sqcup^CB$,
\item a function $\inr:B \to A\sqcup^CB$, and
\item for each $c:C$ a path $\glue(c):(\inl(f(c))=\inr(g(c)))$.
\end{itemize}
In other words, $A\sqcup^CB$ is the disjoint union of $A$ and $B$, together with for every $c:C$ a witness that $f(c)$ and $g(c)$ are equal.
The recursion principle says that if $D$ is another type, we can define a map $s:A\sqcup^CB\to{}D$ by defining
\begin{itemize}
\item for each $a:A$, the value of $s(\inl(a)):D$,
\item for each $b:B$, the value of $s(\inr(b)):D$, and
\item for each $c:C$, the value of $\mapfunc{s}(\glue(c)):s(\inl(f(c)))=s(\inr(g(c)))$.
\end{itemize}
We leave it to the reader to formulate the induction principle.
It also implies the uniqueness principle that if $s,s':A\sqcup^CB\to{}D$ are two maps such that
\index{uniqueness!principle, propositional!for functions on a pushout}%
\begin{align*}
  s(\inl(a))&=s'(\inl(a))\\
  s(\inr(b))&=s'(\inr(b))\\
  \mapfunc{s}(\glue(c))&=\mapfunc{s'}(\glue(c))
  \qquad\text{(modulo the previous two equalities)}
\end{align*}
for every $a,b,c$, then $s=s'$.

To formulate the universal property of a pushout, we introduce the following.

\begin{defn}\label{defn:cocone}
  Given a span $\Ddiag= (A \xleftarrow{f} C \xrightarrow{g} B)$ and a type $D$, a \define{cocone under $\Ddiag$ with vertex $D$}
  \indexdef{cocone}%
  \index{vertex of a cocone}%
  consists of functions $i:A\to{}D$ and $j:B\to{}D$ and a homotopy $h : \prd{c:C} (i(f(c))=j(g(c)))$:
  \[\uppercurveobject{{ }}\lowercurveobject{{ }}\twocellhead{{ }}
  \xymatrix{C \ar^g[r] \ar_f[d] \drtwocell{^h} & B \ar^j[d] \\ A \ar_i[r] & D
  }\]
  We denote by $\cocone{\Ddiag}{D}$ the type of all such cocones, i.e.
  \[ \cocone{\Ddiag}{D} \defeq
  \sm{i:A\to D}{j:B\to D} \prd{c:C} (i(f(c))=j(g(c))).
  \]
\end{defn}

Of course, there is a canonical cocone under $\Ddiag$ with vertex $A\sqcup^C B$ consisting of $\inl$, $\inr$, and $\glue$.
\[\uppercurveobject{{ }}\lowercurveobject{{ }}\twocellhead{{ }}
\xymatrix{C \ar^g[r] \ar_f[d] \drtwocell{^\glue\ \ } & B \ar^\inr[d] \\
  A \ar_-\inl[r] & A\sqcup^CB }\]
The following lemma says that this is the universal such cocone.

\begin{lem}\label{thm:pushout-ump}
  \index{universal!property!of pushout}%
  For any type $E$, there is an equivalence
  \[ (A\sqcup^C B \to E) \;\eqvsym\; \cocone{\Ddiag}{E}. \]
\end{lem}
\begin{proof}
  Let's consider an arbitrary type $E:\type$.
  There is a canonical function
  \[\function{(A\sqcup^CB\to{}E)}{\cocone{\Ddiag}{E}}
  {t}{\composecocone{t}c_\sqcup}\]
  defined by sending $(i,j,h)$ to $(t\circ{}i,t\circ{}j,\mapfunc{t}\circ{}h)$.
  We show that this is an equivalence.

  Firstly, given a $c=(i,j,h):\cocone{\mathscr{D}}{E}$, we need to construct a
  map $\mathsf{s}(c)$ from $A\sqcup^CB$ to $E$.
  \[\uppercurveobject{{ }}\lowercurveobject{{ }}\twocellhead{{ }}
  \xymatrix{C \ar^g[r] \ar_f[d] \drtwocell{^h} & B \ar^{j}[d] \\
    A \ar_-{i}[r] & E }\]
 The map $\mathsf{s}(c)$ is defined in the following way
  \begin{align*}
    \mathsf{s}(c)(\inl(a))&\defeq i(a),\\
    \mathsf{s}(c)(\inr(b))&\defeq j(b),\\
    \mapfunc{\mathsf{s}(c)}(\glue(x))&\defid h(x).
  \end{align*}
We have defined a map
\[\function{\cocone{\Ddiag}{E}}{(A\sqcup^BC\to{}E)}{c}{\mathsf{s}(c)}\]
and we need to prove that this map is an inverse to
$t\mapsto{}\composecocone{t}c_\sqcup$.
On the one hand, if $c=(i,j,h):\cocone{\Ddiag}{E}$, we have
\begin{align*}
  \composecocone{\mathsf{s}(c)}c_\sqcup & =
  (\mathsf{s}(c)\circ\inl,\mathsf{s}(c)\circ\inr,
  \mapfunc{\mathsf{s}(c)}\circ\glue) \\
  & = (\lamu{a:A} \mathsf{s}(c)(\inl(a)),\;
  \lamu{b:B} \mathsf{s}(c)(\inr(b)),\;
  \lamu{x:C} \mapfunc{\mathsf{s}(c)}(\glue(x))) \\
  & = (\lamu{a:A} i(a),\;
  \lamu{b:B} j(b),\;
  \lamu{x:C} h(x)) \\
  & \jdeq (i, j, h) \\
  & = c.
\end{align*}
On the other hand, if $t:A\sqcup^BC\to{}E$, we want to prove that
$\mathsf{s}(\composecocone{t}c_\sqcup)=t$.
For $a:A$, we have
\[\mathsf{s}(\composecocone{t}c_\sqcup)(\inl(a))=t(\inl(a))\]
because the first component of $\composecocone{t}c_\sqcup$ is $t\circ\inl$. In
the same way, for $b:B$ we have
\[\mathsf{s}(\composecocone{t}c_\sqcup)(\inr(b))=t(\inr(b))\]
and for $x:C$ we have
\[\mapfunc{\mathsf{s}(\composecocone{t}c_\sqcup)}(\glue(x))
=\mapfunc{t}(\glue(x))\]
hence $\mathsf{s}(\composecocone{t}c_\sqcup)=t$.

This proves that $c\mapsto\mathsf{s}(c)$ is a quasi-inverse to $t\mapsto{}\composecocone{t}c_\sqcup$, as desired.
\end{proof}

A number of standard homotopy-theoretic constructions can be expressed as (homotopy) pushouts.
\begin{itemize}
\item The pushout of the span $\unit \leftarrow A \to \unit$ is the \define{suspension} $\susp A$ (see \autoref{sec:suspension}).%
  \index{suspension}
\symlabel{join}
\item The pushout of $A \xleftarrow{\proj1} A\times B \xrightarrow{\proj2} B$ is called the \define{join} of $A$ and $B$, written $A*B$.%
  \indexdef{join!of types}
\item The pushout of $\unit \leftarrow A \xrightarrow{f} B$ is the \define{cone} or \define{cofiber} of $f$.%
  \indexdef{cone!of a function}%
  \indexsee{mapping cone}{cone of a function}%
  \indexdef{cofiber of a function}%
\symlabel{wedge}
\item If $A$ and $B$ are equipped with basepoints $a_0:A$ and $b_0:B$, then the pushout of $A \xleftarrow{a_0} \unit \xrightarrow{b_0} B$ is the \define{wedge} $A\vee B$.%
  \indexdef{wedge}
\symlabel{smash}
\item If $A$ and $B$ are pointed as before, define $f:A\vee B \to A\times B$ by $f(\inl(a))\defeq (a,b_0)$ and $f(\inr(b))\defeq (a_0,b)$, with $\ap f \glue \defid \refl{(a_0,b_0)}$.
  Then the cone of $f$ is called the \define{smash product} $A\wedge B$.%
  \indexdef{smash product}
\end{itemize}
We will discuss pushouts further in \autoref{cha:hlevels,cha:homotopy}.

\begin{rmk}
  As remarked in \autoref{subsec:prop-trunc}, the notations $\wedge$ and $\vee$ for the smash product and wedge of pointed spaces are also used in logic for ``and'' and ``or'', respectively.
  Since types in homotopy type theory can behave either like spaces or like propositions, there is technically a potential for conflict --- but since they rarely do both at once, context generally disambiguates.
  Furthermore, the smash product and wedge only apply to \emph{pointed} spaces, while the only pointed mere proposition is $\top\jdeq\unit$ --- and we have $\unit\wedge \unit = \unit$ and $\unit\vee\unit=\unit$ for either meaning of $\wedge$ and $\vee$.
\end{rmk}

\index{pushout|)}%

\begin{rmk}
  Note that colimits do not in general preserve truncatedness.
  For instance, $\Sn^0$ and \unit are both sets, but the pushout of $\unit \leftarrow \Sn^0 \to \unit$ is $\Sn^1$, which is not a set.
  If we are interested in colimits in the category of $n$-types, therefore (and, in particular, in the category of sets), we need to ``truncate'' the colimit somehow.
  We will return to this point in \autoref{sec:hittruncations,cha:hlevels,cha:set-math}.
\end{rmk}

\section{Truncations}
\label{sec:hittruncations}

\index{truncation!propositional|(}%
In \autoref{subsec:prop-trunc} we introduced the propositional truncation as a new type forming operation;
we now observe that it can be obtained as a special case of higher inductive types.
This reduces the problem of understanding truncations to the problem of understanding higher inductives, which at least are amenable to a systematic treatment.
It is also interesting because it provides our first example of a higher inductive type which is truly \emph{recursive}, in that its constructors take inputs from the type being defined (as does the successor $\suc:\nat\to\nat$).

Let $A$ be a type; we define its propositional truncation $\brck A$ to be the higher inductive type generated by:
\begin{itemize}
\item A function $\bprojf : A \to \brck A$, and
\item For each $x,y:\brck A$, a path $x=y$.
\end{itemize}
Note that the second constructor is by definition the assertion that $\brck A$ is a mere proposition.
Thus, the definition of $\brck A$ can be interpreted as saying that $\brck A$ is freely generated by a function $A\to\brck A$ and the fact that it is a mere proposition.

The recursion principle for this higher inductive definition is easy to write down: it says that given any type $B$ together with
\begin{itemize}
\item A function $g:A\to B$, and
\item For any $x,y:B$, a path $x=_B y$,
\end{itemize}
there exists a function $f:\brck A \to B$ such that
\begin{itemize}
\item $f(\bproj a) \jdeq g(a)$ for all $a:A$, and
\item for any $x,y:\brck A$, the function $\apfunc f$ takes the specified path $x=y$ in $\brck A$ to the specified path $f(x) = f(y)$ in $B$ (propositionally).
\end{itemize}
\index{recursion principle!for truncation}%
These are exactly the hypotheses that we stated in \autoref{subsec:prop-trunc} for the recursion principle of propositional truncation --- a function $A\to B$ such that $B$ is a mere proposition --- and the first part of the conclusion is exactly what we stated there as well.
The second part (the action of $\apfunc f$) was not mentioned previously, but it turns out to be vacuous in this case, because $B$ is a mere proposition, so \emph{any} two paths in it are automatically equal.

\index{induction principle!for truncation}%
There is also an induction principle for $\brck A$, which says that given any $B:\brck A \to \type$ together with
\begin{itemize}
\item a function $g:\prd{a:A} B(\bproj a)$, and
\item for any $x,y:\brck A$ and $u:B(x)$ and $v:B(y)$, a dependent path $q:\dpath{B}{p(x,y)}{u}{v}$, where $p(x,y)$ is the path coming from the second constructor of $\brck A$,
\end{itemize}
there exists $f:\prd{x:\brck A} B(x)$ such that $f(\bproj a)\jdeq a$ for $a:A$, and also another computation rule.
However, because there can be at most one function between any two mere propositions (up to homotopy), this induction principle is not really useful (see also \autoref{ex:prop-trunc-ind}).

\index{truncation!propositional|)}%
\index{truncation!set|(}%

\index{set|(}%
We can, however, extend this idea to construct similar truncations landing in $n$-types, for any $n$.
For instance, we might define the \emph{0-trun\-ca\-tion} $\trunc0A$ to be generated by
\begin{itemize}
\item A function $\tprojf0 : A \to \trunc0 A$, and
\item For each $x,y:\trunc0A$ and each $p,q:x=y$, a path $p=q$.
\end{itemize}
Then $\trunc0A$ would be freely generated by a function $A\to \trunc0A$ together with the assertion that $\trunc0A$ is a set.
A natural induction principle for it would say that given $B:\trunc0 A \to \type$ together with
\begin{itemize}
\item a function $g:\prd{a:A} B(\tproj0a)$, and
\item for any $x,y:\trunc0A$ with $z:B(x)$ and $w:B(y)$, and each $p,q:x=y$ with $r:\dpath{B}{p}{z}{w}$ and $s:\dpath{B}{q}{z}{w}$, a 2-path $v:\dpath{B}{u(x,y,p,q)}{p}{q}$, where $u(x,y,p,q):p=q$ is obtained from the second constructor of $\trunc0A$,
\end{itemize}
there exists $f:\prd{x:\trunc0A} B(x)$ such that $f(\tproj0a)\jdeq g(a)$ for all $a:A$, and also $\apdtwo{f}{u(x,y,p,q)}$ is the 2-path specified above.
(As in the propositional case, the latter condition turns out to be uninteresting.)
From this, however, we can prove a more useful induction principle.

\begin{lem}\label{thm:trunc0-ind}
  Suppose given $B:\trunc0 A \to \type$ together with $g:\prd{a:A} B(\tproj0a)$, and assume that each $B(x)$ is a set.
  Then there exists $f:\prd{x:\trunc0A} B(x)$ such that $f(\tproj0a)\jdeq g(a)$ for all $a:A$.
\end{lem}
\begin{proof}
  It suffices to construct, for any $x,y,z,w,p,q,r,s$ as above, a 2-path $v:\dpath{B}{u(x,y,p,q)}{p}{q}$.
  However, by the definition of dependent 2-paths, this is an ordinary 2-path in the fiber $B(y)$.
  Since $B(y)$ is a set, a 2-path exists between any two parallel paths.
\end{proof}

This implies the expected universal property.

\begin{lem}\label{thm:trunc0-lump}
  \index{universal!property!of truncation}%
  For any set $B$ and any type $A$, composition with $\tprojf0:A\to \trunc0A$ determines an equivalence
  \[ \eqvspaced{(\trunc0A\to B)}{(A\to B)}. \]
\end{lem}
\begin{proof}
  The special case of \autoref{thm:trunc0-ind} when $B$ is the constant family gives a map from right to left, which is a right inverse to the ``compose with $\tprojf0$'' function from left to right.
  To show that it is also a left inverse, let $h:\trunc0A\to B$, and define $h':\trunc0A\to B$ by applying \autoref{thm:trunc0-ind} to the composite $a\mapsto h(\tproj0a)$.
  Thus, $h'(\tproj0a)=h(\tproj0a)$.

  However, since $B$ is a set, for any $x:\trunc0A$ the type $h(x)=h'(x)$ is a mere proposition, and hence also a set.
  Therefore, by \autoref{thm:trunc0-ind}, the observation that $h'(\tproj0a)=h(\tproj0a)$ for any $a:A$ implies $h(x)=h'(x)$ for any $x:\trunc0A$, and hence $h=h'$.
\end{proof}

\index{limit!of sets}%
\index{colimit!of sets}%
For instance, this enables us to construct colimits of sets.
We have seen that if $A \xleftarrow{f} C \xrightarrow{g} B$ is a span of sets, then the pushout $A\sqcup^C B$ may no longer be a set.
(For instance, if $A$ and $B$ are \unit and $C$ is \bool, then the pushout is $\Sn^1$.)
However, we can construct a pushout that is a set, and has the expected universal property with respect to other sets, by truncating.

\begin{lem}\label{thm:set-pushout}
  \index{universal!property!of pushout}%
  Let $A \xleftarrow{f} C \xrightarrow{g} B$ be a span\index{span} of sets.
  Then for any set $E$, there is a canonical equivalence
  \[ \Parens{\trunc0{A\sqcup^C B} \to E} \;\eqvsym\; \cocone{\Ddiag}{E}. \]
\end{lem}
\begin{proof}
  Compose the equivalences in \autoref{thm:pushout-ump,thm:trunc0-lump}.
\end{proof}

We refer to $\trunc0{A\sqcup^C B}$ as the \define{set-pushout}
\indexdef{set-pushout}%
\index{pushout!of sets}
of $f$ and $g$, to distinguish it from the (homotopy) pushout $A\sqcup^C B$.
Alternatively, we could modify the definition of the pushout in \autoref{sec:colimits} to include the $0$-truncation constructor directly, avoiding the need to truncate afterwards.
Similar remarks apply to any sort of colimit of sets; we will explore this further in \autoref{cha:set-math}.

However, while the above definition of the 0-truncation works --- it gives what we want, and is consistent --- it has a couple of issues.
Firstly, it doesn't fit so nicely into the general theory of higher inductive types.
In general, it is tricky to deal directly with constructors such as the second one we have given for $\trunc0A$, whose \emph{inputs} involve not only elements of the type being defined, but paths in it.

This can be gotten round fairly easily, however.
Recall in \autoref{sec:bool-nat} we mentioned that we can allow a constructor of an inductive type $W$ to take ``infinitely many arguments'' of type $W$ by having it take a single argument of type $\nat\to W$.
There is a general principle behind this: to model a constructor with funny-looking inputs, use an auxiliary inductive type (such as \nat) to parametrize them, reducing the input to a simple function with inductive domain.

For the 0-truncation, we can consider the auxiliary \emph{higher} inductive type $S$ generated by two points $a,b:S$ and two paths $p,q:a=b$.
Then the fishy-looking constructor of $\trunc 0A$ can be replaced by the unobjectionable
\begin{itemize}
\item For every $f:S\to A$, a path $\apfunc{f}(p) = \apfunc{f}(q)$.
\end{itemize}
Since to give a map out of $S$ is the same as to give two points and two parallel paths between them, this yields the same induction principle.

\index{set|)}%

\index{truncation!set|)}%
\index{truncation!n-truncation@$n$-truncation}%
A more serious problem with our current definition of $0$-truncation, however, is that it doesn't generalize very well.
If we want to describe a notion of definition of ``$n$-truncation'' into $n$-types uniformly for all $n:\nat$, then this approach is unfeasible, since the second constructor would need a number of arguments that increases with $n$.
In \autoref{sec:truncations}, therefore, we will use a different idea to construct these, based on the observation that the type $S$ introduced above is equivalent to the circle $\Sn^1$.
This includes the 0-truncation as a special case, and satisfies generalized versions of \autoref{thm:trunc0-ind,thm:trunc0-lump}.

\section{Quotients}
\label{sec:set-quotients}

A particularly important sort of colimit of sets is the \emph{quotient} by a relation.
That is, let $A$ be a set and $R:A\times A \to \prop$ a family of mere propositions (a \define{mere relation}).
\indexdef{relation!mere}%
\indexdef{mere relation}%
Its quotient should be the set-coequalizer of the two projections
\[ \tsm{a,b:A} R(a,b) \rightrightarrows A. \]
We can also describe this directly, as the higher inductive type $A/R$ generated by
\index{set-quotient|(defstyle}%
\indexsee{quotient of sets}{set-quotient}%
\indexsee{type!quotient}{set-quotient}%
\begin{itemize}
\item A function $q:A\to A/R$;
\item For each $a,b:A$ such that $R(a,b)$, an equality $q(a)=q(b)$; and
\item The $0$-truncation constructor: for all $x,y:A/R$ and $r,s:x=y$, we have $r=s$.
\end{itemize}
We may sometimes refer to $A/R$ as the \define{set-quotient} of $A$ by $R$, to emphasize that it produces a set by definition.
(There are more general notions of ``quotient'' in homotopy theory, but they are mostly beyond the scope of this book.
However, in \autoref{sec:rezk} we will consider the ``quotient'' of a type by a 1-groupoid, which is the next level up from set-quotients.)

\begin{rmk}
  It is not actually necessary for the definition of set-quotients, and most of their properties, that $A$ be a set.
  However, this is generally the case of most interest.
\end{rmk}

\begin{lem}\label{thm:quotient-surjective}
  The function $q:A\to A/R$ is surjective.
\end{lem}
\begin{proof}
  We must show that for any $x:A/R$ there merely exists an $a:A$ with $q(a)=x$.
  We use the induction principle of $A/R$.
  The first case is trivial: if $x$ is $q(a)$, then of course there merely exists an $a$ such that $q(a)=q(a)$.
  And since the goal is a mere proposition, it automatically respects all path constructors, so we are done.
\end{proof}

\begin{lem}\label{thm:quotient-ump}
  For any set $B$, precomposing with $q$ yields an equivalence
  \[ \eqvspaced{(A/R \to B)}{\Parens{\sm{f:A\to B} \prd{a,b:A} R(a,b) \to (f(a)=f(b))}}.\]
\end{lem}
\begin{proof}
  The quasi-inverse of $\blank\circ q$, going from right to left, is just the recursion principle for $A/R$.
  That is, given $f:A\to B$ such that
  \narrowequation{\prd{a,b:A} R(a,b) \to (f(a)=f(b)),} we define $\bar f:A/R\to B$ by $\bar f(q(a))\defeq f(a)$.
  This defining equation says precisely that $(f\mapsto \bar f)$ is a right inverse to $(\blank\circ q)$.

  For it to also be a left inverse, we must show that for any $g:A/R\to B$ and $x:A/R$ we have $g(x) = \overline{g\circ q}$.
  However, by \autoref{thm:quotient-surjective} there merely exists $a$ such that $q(a)=x$.
  Since our desired equality is a mere proposition, we may assume there purely exists such an $a$, in which case $g(x) = g(q(a)) = \overline{g\circ q}(q(a)) = \overline{g\circ q}(x)$.  
\end{proof}

Of course, classically the usual case to consider is when $R$ is an \define{equivalence relation}, i.e.\ we have
\indexdef{relation!equivalence}%
\indexsee{equivalence!relation}{relation, equivalence}%
\begin{itemize}
\item \define{reflexivity}: $\prd{a:A} R(a,a)$,
  \indexdef{reflexivity!of a relation}%
  \indexdef{relation!reflexive}%
\item \define{symmetry}: $\prd{a,b:A} R(a,b) \to R(b,a)$, and
  \indexdef{symmetry!of a relation}%
  \indexdef{relation!symmetric}%
\item \define{transitivity}: $\prd{a,b,c:C} R(a,b) \times R(b,c) \to R(a,c)$.
  \indexdef{transitivity!of a relation}%
  \indexdef{relation!transitive}%
\end{itemize}
In this case, the set-quotient $A/R$ has additional good properties, as we will see in \autoref{sec:piw-pretopos}: for instance, we have $R(a,b) \eqvsym (\id[A/R]{q(a)}{q(b)})$.
\symlabel{equivalencerelation}
We often write an equivalence relation $R(a,b)$ infix as $a\eqr b$.

The quotient by an equivalence relation can also be constructed in other ways.
The set theoretic approach is to consider the set of equivalence classes, as a subset of the power set\index{power set} of $A$.
We can mimic this ``impredicative'' construction in type theory as well.
\index{impredicative!quotient}

\begin{defn}
  A predicate $P:A\to\prop$ is an \define{equivalence class}
  \indexdef{equivalence!class}%
  of a relation $R : A \times A \to \prop$ if there merely exists an $a:A$ such that for all $b:A$ we have $\eqv{R(a,b)}{P(b)}$.
\end{defn}

As $R$ and $P$ are mere propositions, the equivalence $\eqv{R(a,b)}{P(b)}$ is the same thing as implications $R(a,b) \to P(b)$ and $P(b) \to R(a,b)$.
And of course, for any $a:A$ we have the canonical equivalence class $P_a(b) \defeq R(a,b)$.

\begin{defn}\label{def:VVquotient}
  We define
  \begin{equation*}
    A\sslash R \defeq \setof{ P:A\to\prop | P \text{ is an equivalence class of } R}.
  \end{equation*}
  The function $q':A\to A\sslash R$ is defined by $q'(a) \defeq P_a$.
\end{defn}

\begin{thm}
  For any equivalence relation $R$ on $A$, the two set-quotients $A/R$ and $A\sslash R$ are equivalent.
\end{thm}
\begin{proof}
  First, note that if $R(a,b)$, then since $R$ is an equivalence relation we have $R(a,c) \Leftrightarrow R(b,c)$ for any $c:A$.
  Thus, $R(a,c) = R(b,c)$ by univalence, hence $P_a=P_b$ by function extensionality, i.e.\ $q'(a)=q'(b)$.
  Therefore, by \autoref{thm:quotient-ump} we have an induced map $f:A/R \to A\sslash R$ such that $f\circ q = q'$.

  We show that $f$ is injective and surjective, hence an equivalence.
  Surjectivity follows immediately from the fact that $q'$ is surjective, which in turn is true essentially by definition of $A\sslash R$.
  For injectivity, if $f(x)=f(y)$, then to show the mere proposition $x=y$, by surjectivity of $q$ we may assume $x=q(a)$ and $y=q(b)$ for some $a,b:A$.
  Then $R(a,c) = f(q(a))(c) = f(q(b))(c) = R(b,c)$ for any $c:A$, and in particular $R(a,b) = R(b,b)$.
  But $R(b,b)$ is inhabited, since $R$ is an equivalence relation, hence so is $R(a,b)$.
  Thus $q(a)=q(b)$ and so $x=y$.
\end{proof}

In \autoref{subsec:quotients} we will give an alternative proof of this theorem.
Note that unlike $A/R$, the construction $A\sslash R$ raises universe level: if $A:\UU_i$ and $R:A\to A\to \prop_{\UU_i}$, then in the definition of $A\sslash R$ we must also use $\prop_{\UU_i}$ to include all the equivalence classes, so that $A\sslash R : \UU_{i+1}$.
Of course, we can avoid this if we assume the propositional resizing axiom from \autoref{subsec:prop-subsets}.

\begin{rmk}\label{defn-Z}
The previous two constructions provide quotients in generality, but in particular cases there may be easier constructions.
For instance, we may define the integers \Z as a set-quotient
\indexdef{integers}%
\indexdef{number!integers}%
\[ \Z \defeq (\N \times \N)/{\eqr} \]
where $\eqr$ is the equivalence relation defined by
\[ (a,b) \eqr (c,d) \defeq (a + d = b + c). \]
In other words, a pair $(a,b)$ represents the integer $a - b$.
In this case, however, there are \emph{canonical representatives} of the equivalence classes: those of the form $(n,0)$ or $(0,n)$.
\end{rmk}

The following lemma says that when this sort of thing happens, we don't need either general construction of quotients.
(A function $r:A\to A$ is called \define{idempotent}
\indexdef{function!idempotent}%
\indexdef{idempotent!function}%
if $r\circ r = r$.)

\begin{lem}\label{lem:quotient-when-canonical-representatives}
  Suppose $\eqr$ is an equivalence relation on a set $A$, and there exists an idempotent $r
  : A \to A$ such that, for all $x, y \in A$, $\eqv{(r(x) = r(y))}{(x \eqr y)}$. Then the
  type
  \begin{equation*}
    (A/{\eqr}) \defeq \sm{x : A} r(x) = x
  \end{equation*}
  is the set-quotient of $A$ by~$\eqr$.
  In other words, there is a map $q : A \to (A/{\eqr})$ such that for every set $B$, the type $(A/{\eqr}) \to B$ is equivalent to
  \begin{equation}
    \label{eq:quotient-when-canonical}
    \sm{g : A \to B} \prd{x, y : A} (x \eqr y) \to (g(x) = g(y))
  \end{equation}
  with the map being induced by precomposition with $q$.
\end{lem}

\begin{proof}
  Let $i : \prd{x : A} r(r(x)) = r(x)$ witness idempotence of~$r$.
  The map $q : A \to A/{\eqr}$ is defined by $q(x) \defeq (r(x), i(x))$. An equivalence $e$
  from $A/{\eqr} \to B$ to~\eqref{eq:quotient-when-canonical} is defined by
  \[ e(f) \defeq (f \circ q, \nameless), \]
  where the underscore $\nameless$ denotes the following proof: if $x, y : A$ and $x \eqr y$ then by assumption
  $r(x) = r(y)$, hence $(r(x), i(x)) = (r(y), i(y))$ as $A$ is a set, therefore $f(q(x)) =
  f(q(y))$. To see that $e$ is an equivalence, consider the map $e'$ in the opposite
  direction,
  \[ e'(g, p) (x, q) \jdeq g(x). \]
  Given any $f : A/{\eqr} \to B$,
  \[ e'(e(f))(x, p) \jdeq f(q(x)) \jdeq f(r(x), i(x)) = f(x, p) \]
  where the last equality holds because $p : r(x) = x$ and so $(x,p) = (r(x), i(x))$
  because $A$ is a set. Similarly we compute
  \[ e(e'(g, p)) \jdeq e(g \circ \proj{1}) \jdeq (f \circ \proj{1} \circ q, {\nameless}). \]
  Because $B$ is a set we need not worry about the $\nameless$ part, while for the first
  component we have
  \[ f(\proj{1}(q(x))) \defeq f(r(x)) = f(x), \]
  where the last equation holds because $r(x) \eqr x$ and $f$ respects $\eqr$ by
  assumption.
\end{proof}

\begin{cor}\label{thm:retraction-quotient}
  Suppose $p:A\to B$ is a retraction between sets.
  Then $B$ is the quotient of $A$ by the equivalence relation $\eqr$ defined by
  \[ (a_1 \eqr a_2) \defeq (p(a_1) = p(a_2)). \]
\end{cor}
\begin{proof}
  Suppose $s:B\to A$ is a section of $p$.
  Then $s\circ p : A\to A$ is an idempotent which satisfies the condition of \autoref{lem:quotient-when-canonical-representatives} for this $\eqr$, and $s$ induces an isomorphism from $B$ to its set of fixed points.
\end{proof}

\begin{rmk}\label{Z-quotient-by-canonical-representatives}
\autoref{lem:quotient-when-canonical-representatives} applies to $\Z$ with the idempotent $r : \N \times \N \to \N \times \N$
defined by
\begin{equation*}
  r(a, b) =
  \begin{cases}
    (a - b, 0) & \text{if $a \geq b$,} \\
    (0, b - a) & \text{otherwise.}
  \end{cases}  
\end{equation*}
(This is a valid definition even constructively, since the relation $\geq$ on $\N$ is decidable.)
Thus a non-negative integer is canonically represented as $(k, 0)$ and a non-positive one by $(0, m)$, for $k,m:\N$.
This division into cases implies the following induction principle for integers, which will be useful in \autoref{cha:homotopy}.
\index{natural numbers}%
(As usual, we identify natural numbers with the corresponding non-negative integers.)
\end{rmk}

\begin{lem}\label{thm:sign-induction}
  \index{integers!induction principle for}%
  \index{induction principle!for integers}%
  Suppose $P:\Z\to\type$ is a type family and that we have
  \begin{itemize}
  \item $d_0: P(0)$,
  \item $d_+: \prd{n:\N} P(n) \to P(\suc(n))$, and
  \item $d_- : \prd{n:\N} P(-n) \to P(-\suc(n))$.
  \end{itemize}
  Then we have $f:\prd{z:\Z} P(z)$ such that $f(0)\jdeq d_0$ and $f(\suc(n))\jdeq d_+(f(n))$, and $f(-\suc(n))\jdeq d_-(f(-n))$ for all $n:\N$.
\end{lem}
\begin{proof}
  We identify $\Z$ with $\sm{x:\N\times\N}(r(x)=x)$, where $r$ is the above idempotent.
  Now define $Q\defeq P\circ r:\N\times \N \to \type$.
  We can construct $g:\prd{x:\N\times \N} Q(x)$ by double induction on $n$:
  \begin{align*}
    g(0,0) &\defeq d_0,\\
    g(\suc(n),0) &\defeq d_+(g(n,0)),\\
    g(0,\suc(m)) &\defeq d_-(g(0,m)),\\
    g(\suc(n),\suc(m)) &\defeq g(n,m).
  \end{align*}
  Let $f$ be the restriction of $g$ to $\Z$.
\end{proof}

For example, we can define the $n$-fold concatenation of a loop for any integer $n$.

\begin{cor}\label{thm:looptothe}
  \indexdef{path!concatenation!n-fold@$n$-fold}%
  Let $A$ be a type with $a:A$ and $p:a=a$.
  There is a function $\prd{n:\Z} (a=a)$, denoted $n\mapsto p^n$, defined by
  \begin{align*}
    p^0 &\defeq \refl{\base}\\
    p^{n+1} &\defeq p^n \ct p
    & &\text{for $n\ge 0$}\\
    p^{n-1} &\defeq p^n \ct \opp p
    & &\text{for $n\le 0$.}
  \end{align*}
\end{cor}

We will discuss the integers further in \autoref{sec:free-algebras,sec:field-rati-numb}.

\index{set-quotient|)}%

\section{Algebra}
\label{sec:free-algebras}

In addition to constructing higher-dimensional objects such as spheres and cell complexes, higher inductive types are also very useful even when working only with sets.
We have seen one example already in \autoref{thm:set-pushout}: they allow us to construct the colimit of any diagram of sets, which is not possible in the base type theory of \autoref{cha:typetheory}.
Higher inductive types are also very useful when we study sets with algebraic structure.

As a running example in this section, we consider \emph{groups}, which are familiar to most mathematicians and exhibit the essential phenomena (and will be needed in later chapters).
However, most of what we say applies equally well to any sort of algebraic structure.

\index{monoid|(}%

\begin{defn}
  A \define{monoid}
  \indexdef{monoid}%
  is a set $G$ together with
  \begin{itemize}
  \item a \emph{multiplication}
    \indexdef{multiplication!in a monoid}%
    \indexdef{multiplication!in a group}%
    function $G\times G\to G$, written infix as $(x,y) \mapsto x\cdot y$; and
  \item a \emph{unit}
    \indexdef{unit!of a monoid}%
    \indexdef{unit!of a group}%
    element $e:G$; such that
  \item for any $x:G$, we have $x\cdot e = x$ and $e\cdot x = x$; and
  \item for any $x,y,z:G$, we have $x\cdot (y\cdot z) = (x\cdot y)\cdot z$.
    \index{associativity!in a monoid}%
    \index{associativity!in a group}%
  \end{itemize}
  A \define{group}
  \indexdef{group}%
  is a monoid $G$ together with
  \begin{itemize}
  \item an \emph{inversion} function $i:G\to G$, written $x\mapsto \opp x$; such that
    \index{inverse!in a group}%
  \item for any $x:G$ we have $x\cdot \opp x = e$ and $\opp x \cdot x = e$.
  \end{itemize}
\end{defn}

\begin{rmk}\label{rmk:infty-group}
Note that we require a group to be a set.
We could consider a more general notion of ``$\infty$-group''%
\index{.infinity-group@$\infty$-group}
which is not a set, but this would take us further afield than is appropriate at the moment.
With our current definition, we may expect the resulting ``group theory'' to behave similarly to the way it does in set-theoretic mathematics (with the caveat that, unless we assume \LEM{}, it will be ``constructive'' group theory).\index{mathematics!constructive}
\end{rmk}

\begin{eg}
  The natural numbers \N are a monoid under addition, with unit $0$, and also under multiplication, with unit $1$.
  If we define the arithmetical operations on the integers \Z in the obvious way, then as usual they are a group under addition and a monoid under multiplication (and, of course, a ring).
  For instance, if $u, v \in \Z$ are represented by $(a,b)$ and $(c,d)$, respectively, then $u + v$ is represented by $(a + c, b + d)$, $-u$ is represented by $(b, a)$, and $u v$ is represented by $(a c + b d, a d + b c)$.
\end{eg}

\begin{eg}\label{thm:homotopy-groups}
  We essentially observed in \autoref{sec:equality} that if $(A,a)$ is a pointed type, then its loop space\index{loop space} $\Omega(A,a)\defeq (\id[A]aa)$ has all the structure of a group, except that it is not in general a set.
  It should be an ``$\infty$-group'' in the sense mentioned in \autoref{rmk:infty-group}, but we can also make it a group by truncation.
  Specifically, we define the \define{fundamental group}
  \indexsee{group!fundamental}{fundamental group}%
  \indexdef{fundamental!group}%
  of $A$ based at $a:A$ to be
  \[\pi_1(A,a)\defeq \trunc0{\Omega(A,a)}.\]
  This inherits a group structure; for instance, the multiplication $\pi_1(A,a) \times \pi_1(A,a) \to \pi_1(A,a)$ is defined by double induction on truncation from the concatenation of paths.

  More generally, the \define{$n^{\mathrm{th}}$ homotopy group}
  \index{homotopy!group}%
  \indexsee{group!homotopy}{homotopy group}%
  of $(A,a)$ is $\pi_n(A,a)\defeq \trunc0{\Omega^n(A,a)}$.
  \index{loop space!iterated}%
  Then $\pi_n(A,a) = \pi_1(\Omega^{n-1}(A,a))$ for $n\ge 1$, so it is also a group.
  (When $n=0$, we have $\pi_0(A) \jdeq \trunc0 A$, which is not a group.)
  Moreover, the Eckmann--Hilton argument \index{Eckmann--Hilton argument} (\autoref{thm:EckmannHilton}) implies that if $n\ge 2$, then $\pi_n(A,a)$ is an \emph{abelian}\index{group!abelian} group, i.e.\ we have $x\cdot y = y\cdot x$ for all $x,y$.
  \autoref{cha:homotopy} will be largely the study of these groups.
\end{eg}

\index{algebra!free}%
\index{free!algebraic structure}%
One important notion in group theory is that of the \emph{free group} generated by a set, or more generally of a group \emph{presented} by generators\index{generator!of a group} and relations.
It is well-known in type theory that \emph{some} free algebraic objects can be defined using \emph{ordinary} inductive types.
\symlabel{lst-freemonoid}%
\indexdef{type!of lists}%
\indexsee{list type}{type, of lists}%
\index{monoid!free|(}%
For instance, the free monoid on a set $A$ can be identified with the type $\lst A$ of \emph{finite lists} \index{finite!lists, type of} of elements of $A$, which is inductively generated by
\begin{itemize}
\item a constructor $\nil:\lst A$, and
\item for each $\ell:\lst A$ and $a:A$, an element $\cons(a,\ell):\lst A$.
\end{itemize}
We have an obvious inclusion $\eta : A\to \lst A$ defined by $a\mapsto \cons(a,\nil)$.
The monoid operation on $\lst A$ is concatenation, defined recursively by
\begin{align*}
  \nil \cdot \ell &\defeq \ell\\
  \cons (a,\ell_1) \cdot \ell_2 &\defeq \cons(a, \ell_1\cdot\ell_2).
\end{align*}
It is straightforward to prove, using the induction principle for $\lst A$, that $\lst A$ is a set and that concatenation of lists is associative
\index{associativity!of list concatenation}%
and has $\nil$ as a unit.
Thus, $\lst A$ is a monoid.

\begin{lem}\label{thm:free-monoid}
  \indexsee{free!monoid}{monoid, free}%
  For any set $A$, the type $\lst A$ is the free monoid on $A$.
  In other words, for any monoid $G$, composition with $\eta$ is an equivalence
  \[ \eqv{\hom_{\mathrm{Monoid}}(\lst A,G)}{(A\to G)}, \]
  where $\hom_{\mathrm{Monoid}}(\blank,\blank)$ denotes the set of monoid homomorphisms (functions which preserve the multiplication and unit).
  \indexdef{homomorphism!monoid}%
  \indexdef{monoid!homomorphism}%
\end{lem}
\begin{proof}
  Given $f:A\to G$, we define $\bar{f}:\lst A \to G$ by recursion:
  \begin{align*}
    \bar{f}(\nil) &\defeq e\\
    \bar{f}(\cons(a,\ell)) &\defeq f(a) \cdot \bar{f}(\ell).
  \end{align*}
  It is straightforward to prove by induction that $\bar{f}$ is a monoid homomorphism, and that $f\mapsto \bar f$ is a quasi-inverse of $(\blank\circ \eta)$; see \autoref{ex:free-monoid}.
\end{proof}

\index{monoid!free|)}%

This construction of the free monoid is possible essentially because elements of the free monoid have computable canonical forms (namely, finite lists).
However, elements of other free (and presented) algebraic structures --- such as groups --- do not in general have \emph{computable} canonical forms.
For instance, equality of words in group presentations is algorithmically\index{algorithm} undecidable.
However, we can still describe free algebraic objects as \emph{higher} inductive types, by simply asserting all the axiomatic equations as path constructors.

\indexsee{free!group}{group, free}%
\index{group!free|(}%
For example, let $A$ be a set, and define a higher inductive type $\freegroup{A}$ with the following generators.
\begin{itemize}
\item A function $\eta:A\to \freegroup{A}$.
\item A function $m: \freegroup{A} \times \freegroup{A} \to \freegroup{A}$.
\item An element $e:\freegroup{A}$.
\item A function $i:\freegroup{A} \to \freegroup{A}$.
\item For each $x,y,z:\freegroup{A}$, an equality $m(x,m(y,z)) = m(m(x,y),z)$.
\item For each $x:\freegroup{A}$, equalities $m(x,e) = x$ and $m(e,x) = x$.
\item For each $x:\freegroup{A}$, equalities $m(x,i(x)) = e$ and $m(i(x),x) = e$.
\item The $0$-truncation constructor: for any $x,y:\freegroup{A}$ and $p,q:x=y$, we have $p=q$.
\end{itemize}
The first constructor says that $A$ maps to $\freegroup{A}$.
The next three give $\freegroup{A}$ the operations of a group: multiplication, an identity element, and inversion.
The three constructors after that assert the axioms of a group: associativity\index{associativity}, unitality, and inverses.
Finally, the last constructor asserts that $\freegroup{A}$ is a set.

Therefore, $\freegroup{A}$ is a group.
It is also straightforward to prove:

\begin{thm}
  \index{universal!property!of free group}%
  $\freegroup{A}$ is the free group on $A$.
  In other words, for any (set) group $G$, composition with $\eta:A\to \freegroup{A}$ determines an equivalence
  \[ \hom_{\mathrm{Group}}(\freegroup{A},G) \eqvsym (A\to G) \]
  where $\hom_{\mathrm{Group}}(\blank,\blank)$ denotes the set of group homomorphisms between two groups.
  \indexdef{group!homomorphism}%
  \indexdef{homomorphism!group}%
\end{thm}
\begin{proof}
  The recursion principle of the higher inductive type $\freegroup{A}$ says \emph{precisely} that if $G$ is a group and we have $f:A\to G$, then we have $\bar{f}:\freegroup{A} \to G$.
  Its computation rules say that $\bar{f}\circ \eta \jdeq f$, and that $\bar f$ is a group homomorphism.
  Thus, $(\blank\circ \eta) :  \hom_{\mathrm{Group}}(\freegroup{A},G) \to (A\to G)$ has a right inverse.
  It is straightforward to use the induction principle of $\freegroup{A}$ to show that this is also a left inverse.
\end{proof}

\index{acceptance}
It is worth taking a step back to consider what we have just done.
We have proven that the free group on any set exists \emph{without} giving an explicit construction of it.
Essentially all we had to do was write down the universal property that it should satisfy.
In set theory, we could achieve a similar result by appealing to black boxes such as the adjoint functor theorem\index{adjoint!functor theorem}; type theory builds such constructions into the foundations of mathematics.

Of course, it is sometimes also useful to have a concrete description of free algebraic structures.
In the case of free groups, we can provide one, using quotients.
Consider $\lst{A+A}$, where in $A+A$ we write $\inl(a)$ as $a$, and $\inr(a)$ as $\hat{a}$ (intended to stand for the formal inverse of $a$).
The elements of $\lst{A+A}$ are \emph{words} for the free group on $A$.

\begin{thm}
  Let $A$ be a set, and let $\freegroupx{A}$ be the set-quotient of $\lst{A+A}$ by the following relations.
  \begin{align*}
    (\dots,a_1,a_2,\widehat{a_2},a_3,\dots) &=
    (\dots,a_1,a_3,\dots)\\
    (\dots,a_1,\widehat{a_2},a_2,a_3,\dots) &=
    (\dots,a_1,a_3,\dots).
  \end{align*}
  Then $\freegroupx{A}$ is also the free group on the set $A$.
\end{thm}
\begin{proof}
  First we show that $\freegroupx{A}$ is a group.
  We have seen that $\lst{A+A}$ is a monoid; we claim that the monoid structure descends to the quotient.
  We define $\freegroupx{A} \times \freegroupx{A} \to \freegroupx{A}$ by double quotient recursion; it suffices to check that the equivalence relation generated by the given relations is preserved by concatenation of lists.
  Similarly, we prove the associativity and unit laws by quotient induction.

  In order to define inverses in $\freegroupx{A}$, we first define $\mathsf{reverse}:\lst B\to\lst B$ by recursion on lists:
  \begin{align*}
    \mathsf{reverse}(\nil) &\defeq \nil,\\
    \mathsf{reverse}(\cons(b,\ell))&\defeq \mathsf{reverse}(\ell)\cdot \cons(b,\nil).
  \end{align*}
  Now we define $i:\freegroupx{A}\to \freegroupx{A}$ by quotient recursion, acting on a list $\ell:\lst{A+A}$ by switching the two copies of $A$ and reversing the list.
  This preserves the relations, hence descends to the quotient.
  And we can prove that $i(x) \cdot x = e$ for $x:\freegroupx{A}$ by induction.
  First, quotient induction allows us to assume $x$ comes from $\ell:\lst{A+A}$, and then we can do list induction:
  \begin{align*}
    i(\nil) \ct \nil &= \nil \ct \nil\\
    &= \nil\\
    i(\cons(a,\ell)) \ct \cons(a,\ell) &= i(\ell) \ct \cons(\hat{a},\nil) \ct \cons(a,\ell)\\
    &= i(\ell) \ct \cons(\hat{a},\cons(a,\ell))\\
    &= i(\ell) \ct \ell\\
    &= \nil. \tag{by the inductive hypothesis}
  \end{align*}
  (We have omitted a number of fairly evident lemmas about the behavior of concatenation of lists, etc.)

  This completes the proof that $\freegroupx{A}$ is a group.
  Now if $G$ is any group with a function $f:A\to G$, we can define $A+A\to G$ to be $f$ on the first copy of $A$ and $f$ composed with the inversion map of $G$ on the second copy.
  Now the fact that $G$ is a monoid yields a monoid homomorphism $\lst{A+A} \to G$.
  And since $G$ is a group, this map respects the relations, hence descends to a map $\freegroupx{A}\to G$.
  It is straightforward to prove that this is a group homomorphism, and the unique one which restricts to $f$ on $A$.
\end{proof}

\index{monoid|)}%

If $A$ has decidable equality\index{decidable!equality} (such as if we assume excluded middle), then the quotient defining $\freegroupx{A}$ can be obtained from an idempotent as in \autoref{lem:quotient-when-canonical-representatives}.
We define a word, which we recall is just an element of $\lst{A+A}$, to be \define{reduced}
\indexdef{reduced word in a free group}
if it contains no adjacent pairs of the form $(a,\hat a)$ or $(\hat a,a)$.
When $A$ has decidable equality, it is straightforward to define the \define{reduction}
\index{reduction!of a word in a free group}%
of a word, which is an idempotent generating the appropriate quotient; we leave the details to the reader.

If $A\defeq \unit$, which has decidable equality, a reduced word must consist either entirely of $\ttt$'s or entirely of $\hat{\ttt}$'s.
Thus, the free group on $\unit$ is equivalent to the integers \Z, with $0$ corresponding to $\nil$, the positive integer $n$ corresponding to a reduced word of $n$ $\ttt$'s, and the negative integer $(-n)$ corresponding to a reduced word of $n$ $\hat{\ttt}$'s.
One could also, of course, show directly that \Z has the universal property of $\freegroup{\unit}$.

\begin{rmk}\label{thm:freegroup-nonset}
  Nowhere in the construction of $\freegroup{A}$ and $\freegroupx{A}$, and the proof of their universal properties, did we use the assumption that $A$ is a set.
  Thus, we can actually construct the free group on an arbitrary type.
  Comparing universal properties, we conclude that $\eqv{\freegroup{A}}{\freegroup{\trunc0A}}$.
\end{rmk}

\index{group!free|)}%

\index{algebra!colimits of}%
We can also use higher inductive types to construct colimits of algebraic objects.
For instance, suppose $f:G\to H$ and $g:G\to K$ are group homomorphisms.
Their pushout in the category of groups, called the \define{amalgamated free product}
\indexdef{amalgamated free product}%
\indexdef{free!product!amalgamated}%
$H *_G K$, can be constructed as the higher inductive type generated by
\begin{itemize}
\item Functions $h:H\to H *_G K$ and $k:K\to H *_G K$.
\item The operations and axioms of a group, as in the definition of $\freegroup{A}$.
\item Axioms asserting that $h$ and $k$ are group homomorphisms.
\item For $x:G$, we have $h(f(x)) = k(g(x))$.
\item The $0$-truncation constructor.
\end{itemize}
On the other hand, it can also be constructed explicitly, as the set-quotient of $\lst{H+K}$ by the following relations:
\begin{align*}
  (\dots, x_1, x_2, \dots) &= (\dots, x_1\cdot x_2, \dots)
  & &\text{for $x_1,x_2:H$}\\
  (\dots, y_1, y_2, \dots) &= (\dots, y_1\cdot y_2, \dots)
  & &\text{for $y_1,y_2:K$}\\
  (\dots, 1_G, \dots) &= (\dots, \dots) &&  \\
  (\dots, 1_H, \dots) &= (\dots, \dots) &&  \\
  (\dots, f(x), \dots) &= (\dots, g(x), \dots)
  & &\text{for $x:G$.}
\end{align*}
We leave the proofs to the reader.
In the special case that $G$ is the trivial group, the last relation is unnecessary, and we obtain the \define{free product}
\indexdef{free!product}%
$H*K$, the coproduct in the category of groups.
(This notation unfortunately clashes with that for the \emph{join} of types, as in \autoref{sec:colimits}, but context generally disambiguates.)

\index{presentation!of a group}%
Note that groups defined by \emph{presentations} can be regarded as a special case of colimits.
Suppose given a set (or more generally a type) $A$, and a pair of functions $R\rightrightarrows \freegroup{A}$.
We regard $R$ as the type of ``relations'', with the two functions assigning to each relation the two words that it sets equal.
For instance, in the presentation $\langle a \mid a^2 = e \rangle$ we would have $A\defeq \unit$ and $R\defeq \unit$, with the two morphisms $R\rightrightarrows \freegroup{A}$ picking out the list $(a,a)$ and the empty list $\nil$, respectively.
Then by the universal property of free groups, we obtain a pair of group homomorphisms $\freegroup{R} \rightrightarrows \freegroup{A}$.
Their coequalizer in the category of groups, which can be built just like the pushout, is the group \emph{presented} by this presentation.

\mentalpause

Note that all these sorts of construction only apply to \emph{algebraic} theories,\index{theory!algebraic} which are theories whose axioms are (universally quantified) equations referring to variables, constants, and operations from a given signature\index{signature!of an algebraic theory}.
They can be modified to apply also to what are called \emph{essentially algebraic theories}:\index{theory!essentially algebraic} those whose operations are partially defined on a domain specified by equalities between previous operations.
They do not apply, for instance, to the theory of fields, in which the ``inversion'' operation is partially defined on a domain $\setof{x | x \mathrel{\#} 0}$ specified by an \emph{apartness} $\#$ between previous operations, see \autoref{RD-inverse-apart-0}.
And indeed, it is well-known that the category of fields has no initial object.
\index{initial!field}%

On the other hand, these constructions do apply just as well to \emph{infinitary}\index{infinitary!algebraic theory} algebraic theories, whose ``operations'' can take infinitely many inputs.
In such cases, there may not be any presentation of free algebras or colimits of algebras as a simple quotient, unless we assume the axiom of choice.
This means that higher inductive types represent a significant strengthening of constructive type theory (not necessarily in terms of proof-theoretic strength, but in terms of practical power), and indeed are stronger in some ways than Zermelo--Fraenkel\index{set theory!Zermelo--Fraenkel} set theory (without choice).

\section{The flattening lemma}
\label{sec:flattening}

As we will see in \autoref{cha:homotopy}, amazing things happen when we combine higher inductive types with univalence.
The principal way this comes about is that if $W$ is a higher inductive type and \UU is a type universe, then we can define a type family $P:W\to \UU$ by using the recursion principle for $W$.
When we come to the clauses of the recursion principle dealing with the path constructors of $W$, we will need to supply paths in \UU, and this is where univalence comes in.

For example, suppose we have a type $X$ and a self-equivalence $e:\eqv X X$.
Then we can define a type family $P:\Sn^1 \to \UU$ by using $\Sn^1$-recursion:
\begin{equation*}
  P(\base) \defeq X
  \qquad\text{and}\qquad
  \ap P\lloop \defid \ua(e).
\end{equation*}
The type $X$ thus appears as the fiber $P(\base)$ of $P$ over the basepoint.
The self-equivalence $e$ is a little more hidden in $P$, but the following lemma says that it can be extracted by transporting along \lloop.

\begin{lem}\label{thm:transport-is-given}
  Given $B:A\to\type$ and $x,y:A$, with a path $p:x=y$ and an equivalence $e:\eqv{P(x)}{P(y)}$ such that $\ap{B}p = \ua(e)$, then for any $u:P(x)$ we have
  \begin{align*}
    \transfib{B}{p}{u} &= e(u).
  \end{align*}
\end{lem}
\begin{proof}
  Applying \autoref{thm:transport-is-ap}, we have
  \begin{align*}
    \transfib{B}{p}{u} &= \idtoeqv(\ap{B}p)(u)\\
    &= \idtoeqv(\ua(e))(u)\\
    &= e(u).\qedhere
  \end{align*}
\end{proof}

We have seen type families defined by recursion before: in \autoref{sec:compute-coprod,sec:compute-nat} we used them to characterize the identity types of (ordinary) inductive types.
In \autoref{cha:homotopy}, we will use similar ideas to calculate homotopy groups of higher inductive types.

In this section, we describe a general lemma about type families of this sort which will be useful later on.
We call it the \define{flattening lemma}:
\indexdef{flattening lemma}%
\indexdef{lemma!flattening}%
it says that if $P:W\to\UU$ is defined recursively as above, then its total space $\sm{x:W} P(x)$ is equivalent to a ``flattened'' higher inductive type, whose constructors may be deduced from those of $W$ and the definition of $P$.
From a category-theoretic point of view, $\sm{x:W} P(x)$ is the ``Grothendieck\index{Grothendieck construction} construction'' of $P$, and this expresses its universal property as a ``lax\index{lax colimit} colimit''.

We prove here one general case of the flattening lemma, which directly implies many particular cases and suggests the method to prove others.
Suppose we have $A,B:\type$ and $f,g:B\to{}A$, and that the higher inductive type $W$ is generated by
\begin{itemize}
\item $\cc:A\to{}W$ and
\item $\pp:\prd{b:B} (\cc(f(b))=_W\cc(g(b)))$.
\end{itemize}
Thus, $W$ is the \define{(homotopy) coequalizer}
\indexdef{coequalizer}%
\indexdef{type!coequalizer}%
of $f$ and $g$.
Using binary sums (coproducts) and dependent sums ($\Sigma$-types), a lot of interesting nonrecursive higher
inductive types can be represented in this form. All point constructors have to
be bundled in the type $A$ and all path constructors in the type $B$.
For instance:
\begin{itemize}
\item The circle $\Sn^1$ can be represented by taking $A\defeq \unit$ and $B\defeq \unit$, with $f$ and $g$ the identity.
\item The pushout of $j:X\to Y$ and $k:X\to Z$ can be represented by taking $A\defeq Y+Z$ and $B\defeq X$, with $f\defeq \inl \circ j$ and $g\defeq \inr\circ k$.
\end{itemize}
Now suppose in addition that
\begin{itemize}
\item $C:A\to\type$ is a family of types over $A$, and
\item $D:\prd{b:B}\eqv{C(f(b))}{C(g(b))}$ is a family of equivalences over $B$.
\end{itemize}
Define a type family $P : W\to\type$ inductively by
\begin{align*}
  P(\cc(a)) &\defeq C(a)\\
  \map{P}{\pp(b)} &\defid \ua(D(b)).
\end{align*}
Let \Wtil be the higher inductive type generated by
\begin{itemize}
\item $\cct:\prd{a:A} C(a) \to \Wtil$ and
\item $\ppt:\prd{b:B}{y:C(f(b))} (\cct(f(b),y)=_{\Wtil}\cct(g(b),D(b)(y)))$.
\end{itemize}

The flattening lemma is:

\begin{lem}[Flattening lemma]\label{thm:flattening}
  In the above situation, we have
  \[ \eqvspaced{\Parens{\sm{x:W} P(x)}}{\widetilde{W}}. \]
\end{lem}

\index{universal!property!of dependent pair type}%
As remarked above, this equivalence can be seen as expressing the universal property of $\sm{x:W} P(x)$ as a ``lax\index{lax colimit} colimit'' of $P$ over $W$.
It can also be seen as part of the \emph{stability and descent} property of colimits, which characterizes higher toposes.%
\index{.infinity1-topos@$(\infty,1)$-topos}%
\index{stability!and descent}%

The proof of \autoref{thm:flattening} occupies the rest of this section.
It is somewhat technical and can be skipped on a first reading.
But it is also a good example of ``proof-relevant mathematics'',
\index{mathematics!proof-relevant}%
so we recommend it on a second reading.

The idea is to show that $\sm{x:W} P(x)$ has the same universal property as \Wtil.
We begin by showing that it comes with analogues of the constructors $\cct$ and $\ppt$.

\begin{lem}
  There are functions
  \begin{itemize}
  \item $\cct':\prd{a:A} C(a) \to \sm{x:W} P(x)$ and
  \item $\ppt':\prd{b:B}{y:C(f(b))} \Big(\cct'(f(b),y)=_{\sm{w:W}P(w)}\cct'(g(b),D(b)(y))\Big)$.
  \end{itemize}
\end{lem}
\begin{proof}
  The first is easy; define $\cct'(a,x) \defeq (\cc(a),x)$ and note that by definition $P(\cc(a))\jdeq C(a)$.
  For the second, suppose given $b:B$ and $y:C(f(b))$; we must give an equality
  \[ (\cc(f(b)),y) = (\cc(g(b),D(b)(y))). \]
  Since we have $\pp(b):f(b)=g(b)$, by equalities in $\Sigma$-types it suffices to give an equality $\trans{\pp(b)}{y} = D(b)(y)$.
  But this follows from \autoref{thm:transport-is-given}, using the definition of $P$.
\end{proof}

Now the following lemma says to define a section of a type family over $\sm{w:W} P(w)$, it suffices to give analogous data as in the case of \Wtil.

\begin{lem}\label{thm:flattening-rect}
  Suppose $Q:\big(\sm{x:W} P(x)\big) \to \type$ is a type family and that we have
  \begin{itemize}
  \item $c : \prd{a:A}{x:C(a)} Q(\cct'(a,x))$ and
  \item $p : \prd{b:B}{y:C(f(b))} \Big(\trans{\ppt'(b,y)}{c(f(b),y)} = c(g(b),D(b)(y))\Big)$. 
  \end{itemize}
  Then there exists $f:\prd{z:\sm{w:W} P(w)} Q(z)$ such that $f(\cct'(a,x)) \jdeq c(a,x)$.
\end{lem}
\begin{proof}
  Suppose given $w:W$ and $x:P(w)$; we must produce an element $f(w,x):Q(w,x)$.
  By induction on $w$, it suffices to consider two cases.
  When $w\jdeq \cc(a)$, then we have $x:C(a)$, and so $c(a,x):Q(\cc(a),x)$ as desired.
  (This part of the definition also ensures that the stated computational rule holds.)

  Now we must show that this definition is preserved by transporting along $\pp(b)$ for any $b:B$.
  Since what we are defining, for all $w:W$, is a function of type $\prd{x:P(w)} Q(w,x)$, by \autoref{thm:dpath-forall} it suffices to show that for any $y:C(f(b))$, we have
  \[ \transfib{Q}{\pairpath(\pp(b),\refl{\trans{\pp(b)}{y}})}{c(f(b),y)} = c(g(b),\trans{\pp(b)}{y}). \]
  Let $q:\trans{\pp(b)}{y} = D(b)(y)$ be the path obtained from \autoref{thm:transport-is-given}.
  Then we have
  \begin{align}
    c(g(b),\trans{\pp(b)}{y})
    &= \transfib{x\mapsto Q(c(g(b),x))}{\opp{q}}{c(g(b),D(b)(y))}
    \tag{by $\apdfunc{x\mapsto c(g(b),x)}(\opp q)$} \\
    &= \transfib{Q}{\apfunc{x\mapsto c(g(b),x)}(\opp q)}{c(g(b),D(b)(y))}
    \tag{by \autoref{thm:transport-compose}}.
  \end{align}
  Thus, it suffices to show
  \begin{multline*}
    \Transfib{Q}{\pairpath(\pp(b),\refl{\trans{\pp(b)}{y}})}{c(f(b),y)} = {}\\
    \Transfib{Q}{\apfunc{x\mapsto c(g(b),x)}(\opp q)}{c(g(b),D(b)(y))}.
  \end{multline*}
  Moving the right-hand transport to the other side, and combining two transports, this is equivalent to
  \begin{narrowmultline*}
    \Transfib{Q}{\apfunc{x\mapsto c(g(b),x)}(q) \ct
      \pairpath(\pp(b),\refl{\trans{\pp(b)}{y}})}{c(f(b),y)} =
    \narrowbreak
    c(g(b),D(b)(y)).
  \end{narrowmultline*}
  However, we have
  \begin{multline*}
    \apfunc{x\mapsto c(g(b),x)}(q) \ct \pairpath(\pp(b),\refl{\trans{\pp(b)}{y}})
    = {} \\
    \pairpath(\refl{g(b)},q) \ct \pairpath(\pp(b),\refl{\trans{\pp(b)}{y}})
    = \pairpath(\pp(b),q)
    = \ppt'(b,y)
  \end{multline*}
  so the construction is completed by the assumption $p(b,y)$ of type
  \[ \transfib{Q}{\ppt'(b,y)}{c(f(b),y)} = c(g(b),D(b)(y)). \qedhere \]
\end{proof}

\autoref{thm:flattening-rect} \emph{almost} gives $\sm{w:W}P(w)$ the same induction principle as \Wtil.
The missing bit is the equality $\apdfunc{f}(\ppt'(b,y)) = p(b,y)$.
In order to prove this, we would need to analyze the proof of \autoref{thm:flattening-rect}, which of course is the definition of $f$.

It should be possible to do this, but it turns out that we only need the computation rule for the non-dependent recursion principle.
Thus, we now give a somewhat simpler direct construction of the recursor, and a proof of its computation rule.

\begin{lem}\label{thm:flattening-rectnd}
  Suppose $Q$ is a type and that we have
  \begin{itemize}
  \item $c : \prd{a:A} C(a) \to Q$ and
  \item $p : \prd{b:B}{y:C(f(b))} \Big(c(f(b),y) =_Q c(g(b),D(b)(y))\Big)$.
  \end{itemize}
  Then there exists $f:\big(\sm{w:W} P(w)\big) \to Q$ such that $f(\cct'(a,x)) \jdeq c(a,x)$.
\end{lem}
\begin{proof}
  As in \autoref{thm:flattening-rect}, we define $f(w,x)$ by induction on $w:W$.
  When $w\jdeq \cc(a)$, we define $f(\cc(a),x)\defeq c(a,x)$.
  Now by \autoref{thm:dpath-arrow}, it suffices to consider, for $b:B$ and $y:C(f(b))$, the composite path
  \begin{equation}\label{eq:flattening-rectnd}
    \transfib{x\mapsto Q}{\pp(b)}{c(f(b),y)}
    = c(g(b),\transfib{P}{\pp(b)}{y})
  \end{equation}
  defined as the composition
  \begin{align}
    \transfib{x\mapsto Q}{\pp(b)}{c(f(b),y)}
    &= c(f(b),y) \tag{by \autoref{thm:trans-trivial}}\\
    &= c(g(b),D(b)(y)) \tag{by $p(b,y)$}\\
    &= c(g(b),\transfib{P}{\pp(b)}{y}). \tag{by \autoref{thm:transport-is-given}}
  \end{align}
  The computation rule $f(\cct'(a,x)) \jdeq c(a,x)$ follows by definition, as before.
\end{proof}

For the second computation rule, we need the following lemma.

\begin{lem}\label{thm:ap-sigma-rect-path-pair}
  Let $Y:X\to\type$ be a type family and let $f:(\sm{x:X}Y(x)) \to Z$ be defined componentwise by $f(x,y) \defeq d(x)(y)$ for a curried function $d:\prd{x:X} Y(x)\to Z$.
  Then for any $s:\id[X]{x_1}{x_2}$ and any $y_1:P(x_1)$ and $y_2:P(x_2)$ with a path $r:\trans{s}{y_1}=y_2$, the path
  \[\apfunc f (\pairpath(s,r)) :f(x_1,y_1) = f(x_2,y_2)\]
  is equal to the composite
  \begin{align}
    f(x_1,y_1)
    &\jdeq d(x_1)(y_1) \notag\\
    &= \transfib{x\mapsto Q}{s}{d(x_1)(y_1)}
    \tag{by $\opp{\text{(\autoref{thm:trans-trivial})}}$}\\
    &= \transfib{x\mapsto Q}{s}{d(x_1)(\trans{\opp s}{\trans{s}{y_1}})}
    \notag\\
    &= \big(\transfib{x\mapsto (Y(x)\to Z)}{s}{d(x_1)}\big)(\trans{s}{y_1})
    \tag{by~\eqref{eq:transport-arrow}}\\
    &= d(x_2)(\trans{s}{y_1})
    \tag{by $\happly(\apdfunc{d}(s))(\trans{s}{y_1}$}\\
    &= d(x_2)(y_2)
    \tag{by $\apfunc{d(x_2)}(r)$}\\
    &\jdeq f(x_2,y_2).
    \notag
  \end{align}
\end{lem}
\begin{proof}
  After path induction on $s$ and $r$, both equalities reduce to reflexivities.
\end{proof}

At first it may seem surprising that \autoref{thm:ap-sigma-rect-path-pair} has such a complicated statement, while it can be proven so simply.
The reason for the complication is to ensure that the statement is well-typed: $\apfunc f (\pairpath(s,r))$ and the composite path it is claimed to be equal to must both have the same start and end points.
Once we have managed this, the proof is easy by path induction.

\begin{lem}\label{thm:flattening-rectnd-beta-ppt}
  In the situation of \autoref{thm:flattening-rectnd}, we have $\apfunc{f}(\ppt'(b,y)) = p(b,y)$.
\end{lem}
\begin{proof}
  Recall that $\ppt'(b,y) \defeq \pairpath(\pp(b),q)$ where $q:\trans{\pp(b)}{y} = D(b)(y)$ comes from \autoref{thm:transport-is-given}.
  Thus, since $f$ is defined componentwise, we may compute $\apfunc{f}(\ppt'(b,y))$ by \autoref{thm:ap-sigma-rect-path-pair}, with
  \begin{align*}
    x_1 &\defeq \cc(f(b)) & y_1 &\defeq y\\
    x_2 &\defeq \cc(g(b)) & y_2 &\defeq D(b)(y)\\
    s &\defeq \pp(b)      &   r &\defeq q.
  \end{align*}
  The curried function $d:\prd{w:W} P(w) \to Q$ was defined by induction on $w:W$;
  to apply \autoref{thm:ap-sigma-rect-path-pair} we need to understand $\apfunc{d(x_2)}(r)$ and $\happly(\apdfunc{d}(s),\trans s {y_1})$.

  For the first, since $d(\cc(a),x)\jdeq c(a,x)$, we have
  \[ \apfunc{d(x_2)}(r) \jdeq \apfunc{c(g(b),-)}(q). \]
  For the second, the computation rule for the induction principle of $W$ tells us that $\apdfunc{d}(\pp(b))$ is equal to the composite~\eqref{eq:flattening-rectnd}, passed across the equivalence of \autoref{thm:dpath-arrow}.
  Thus, the computation rule given in \autoref{thm:dpath-arrow} implies that $\happly(\apdfunc{d}(\pp(b)),\trans {\pp(b)}{y})$ is equal to the composite
  \begin{align}
    \big(\trans{\pp(b)}{c(f(b),-)}\big)(\trans {\pp(b)}{y})
    &= \trans{\pp(b)}{c(f(b),\trans{\opp {\pp(b)}}{\trans {\pp(b)}{y}})}
    \tag{by~\eqref{eq:transport-arrow}}\\
    &= \trans{\pp(b)}{c(f(b),y)}
    \notag \\
    &= c(f(b),y)
    \tag{by \autoref{thm:trans-trivial}}\\
    &= c(f(b),D(b)(y))
   \tag{by $p(b,y)$}\\
    &= c(f(b),\trans{\pp(b)}{y}).
    \tag{by $\opp{\apfunc{c(g(b),-)}(q)}$}
  \end{align}
  Finally, substituting these values of $\apfunc{d(x_2)}(r)$ and $\happly(\apdfunc{d}(s),\trans s {y_1})$ into \autoref{thm:ap-sigma-rect-path-pair}, we see that all the paths cancel out in pairs, leaving only $p(b,y)$.
\end{proof}

Now we are finally ready to prove the flattening lemma.

\begin{proof}[Proof of \autoref{thm:flattening}]
  We define $h:\Wtil \to \sm{w:W}P(w)$ by using the recursion principle for \Wtil, with $\cct'$ and $\ppt'$ as input data.
  Similarly, we define $k:(\sm{w:W}P(w)) \to \Wtil$ by using the recursion principle of \autoref{thm:flattening-rectnd}, with $\cct$ and $\ppt$ as input data.

  On the one hand, we must show that for any $z:\Wtil$, we have $k(h(z))=z$.
  By induction on $z$, it suffices to consider the two constructors of \Wtil.
  But we have
  \[k(h(\cct(a,x))) \jdeq k(\cct'(a,x)) \jdeq \cct(a,x)\]
  by definition, while similarly
  \[\ap k{\ap h{\ppt(b,y)}} = \ap k{\ppt'(b,y)} = \ppt(b,y) \]
  using the propositional computation rule for $\Wtil$ and \autoref{thm:flattening-rectnd-beta-ppt}.

  On the other hand, we must show that for any $z:\sm{w:W}P(w)$, we have $h(k(z))=z$.
  But this is essentially identical, using \autoref{thm:flattening-rect} for ``induction on $\sm{w:W}P(w)$'' and the same computation rules.
\end{proof}

\section{The general syntax of higher inductive definitions}
\label{sec:naturality}

In \autoref{sec:strictly-positive}, we discussed the conditions on a putative ``inductive definition'' which make it acceptable, namely that all inductive occurrences of the type in its constructors are ``strictly positive''.\index{strict!positivity}
In this section, we say something about the additional conditions required for \emph{higher} inductive definitions.
Finding a general syntactic description of valid higher inductive definitions is an area of current research, and all of the solutions proposed to date are somewhat technical in nature; thus we only give a general description and not a precise definition.
Fortunately, the corner cases never seem to arise in practice.

Like an ordinary inductive definition, a higher inductive definition is specified by a list of \emph{constructors}, each of which is a (dependent) function.
For simplicity, we may require the inputs of each constructor to satisfy the same condition as the inputs for constructors of ordinary inductive types.
In particular, they may contain the type being defined only strictly positively.
Note that this excludes definitions such as the $0$-truncation as presented in \autoref{sec:hittruncations}, where the input of a constructor contains not only the inductive type being defined, but its identity type as well.
It may be possible to extend the syntax to allow such definitions; but also, in \autoref{sec:truncations} we will give a different construction of the $0$-truncation whose constructors do satisfy the more restrictive condition.

The only difference between an ordinary inductive definition and a higher one, then, is that the \emph{output} type of a constructor may be, not the type being defined ($W$, say), but some identity type of it, such as $\id[W]uv$, or more generally an iterated identity type such as $\id[({\id[W]uv})]pq$.
Thus, when we give a higher inductive definition, we have to specify not only the inputs of each constructor, but the expressions $u$ and $v$ (or $u$, $v$, $p$, and $q$, etc.)\ which determine the source\index{source!of a path constructor} and target\index{target!of a path constructor} of the path being constructed.

Importantly, these expressions may refer to \emph{other} constructors of $W$.
For instance, in the definition of $\Sn^1$, the constructor $\lloop$ has both $u$ and $v$ being $\base$, the previous constructor.
To make sense of this, we require the constructors of a higher inductive type to be specified \emph{in order}, and we allow the source and target expressions $u$ and $v$ of each constructor to refer to previous constructors, but not later ones.
(Of course, in practice the constructors of any inductive definition are written down in some order, but for ordinary inductive types that order is irrelevant.)

Note that this order is not necessarily the order of ``dimension'': in principle, a 1-dimensional path constructor could refer to a 2-dimensional one and hence need to come after it.
However, we have not given the 0-dimensional constructors (point constructors) any way to refer to previous constructors, so they might as well all come first.
And if we use the hub-and-spoke construction (\autoref{sec:hubs-spokes}) to reduce all constructors to points and 1-paths, then we might assume that all point constructors come first, followed by all 1-path constructors --- but the order among the 1-path constructors continues to matter.

The remaining question is, what sort of expressions can $u$ and $v$ be?
We might hope that they could be any expression at all involving the previous constructors.
However, the following example shows that a naive approach to this idea does not work.

\begin{eg}\label{eg:unnatural-hit}
  Consider a family of functions $f:\prd{X:\type} (X\to X)$.
  Of course, $f_X$ might be just $\idfunc[X]$ for all $X$, but other such $f$s may also exist.
  For instance, nothing prevents $f_{\bool}:\bool\to\bool$ from being the nonidentity automorphism\index{automorphism!of 2, nonidentity@of $\bool$, nonidentity} (see \autoref{ex:unnatural-endomorphisms}).

  Now suppose that we attempt to define a higher inductive type $K$ generated by:
  \begin{itemize}
  \item two elements $a,b:K$, and
  \item a path $\sigma:f_K(a)=f_K(b)$.
  \end{itemize}
  What would the induction principle for $K$ say?
  We would assume a type family $P:K\to\type$, and of course we would need $x:P(a)$ and $y:P(b)$.
  The remaining datum should be a dependent path in $P$ living over $\sigma$, which must therefore connect some element of $P(f_K(a))$ to some element of $P(f_K(b))$.
  But what could these elements possibly be?
  We know that $P(a)$ and $P(b)$ are inhabited by $x$ and $y$, respectively, but this tells us nothing about $P(f_K(a))$ and $P(f_K(b))$.
\end{eg}

Clearly some condition on $u$ and $v$ is required in order for the definition to be sensible.
It seems that, just as the domain of each constructor is required to be (among other things) a \emph{covariant functor}, the appropriate condition on the expressions $u$ and $v$ is that they define \emph{natural transformations}.
Making precise sense of this requirement is beyond the scope of this book, but informally it means that $u$ and $v$ must only involve operations which are preserved by all functions between types.

For instance, it is permissible for $u$ and $v$ to refer to concatenation of paths, as in the case of the final constructor of the torus in \autoref{sec:cell-complexes}, since all functions in type theory preserve path concatenation (up to homotopy).
However, it is not permissible for them to refer to an operation like the function $f$ in \autoref{eg:unnatural-hit}, which is not necessarily natural: there might be some function $g:X\to Y$ such that $f_Y \circ g \neq g\circ f_X$.
(Univalence implies that $f_X$ must be natural with respect to all \emph{equivalences}, but not necessarily with respect to functions that are not equivalences.)

The intuition of naturality supplies only a rough guide for when a higher inductive definition is permissible.
Even if it were possible to give a precise specification of permissible forms of such definitions in this book, such a specification would probably be out of date quickly, as new extensions to the theory are constantly being explored.
For instance, the presentation of $n$-spheres in terms of ``dependent $n$-loops\index{loop!dependent n-@dependent $n$-}'' referred to in \autoref{sec:circle}, and the ``higher inductive-recursive definitions'' used in \autoref{cha:real-numbers}, were innovations introduced while this book was being written.
We encourage the reader to experiment --- with caution.

\sectionNotes

The general idea of higher inductive types was conceived in discussions between Andrej Bauer, Peter Lumsdaine, Mike Shulman, and Michael Warren at the Oberwolfach meeting in 2011, although there are some suggestions of some special cases in earlier work.  Subsequently, Guillaume Brunerie and Dan Licata contributed substantially to the general theory, especially by finding convenient ways to represent them in computer proof assistants
\index{proof!assistant}
and do homotopy theory with them (see \autoref{cha:homotopy}).

A general discussion of the syntax of higher inductive types, and their semantics in higher-categorical models, appears in~\cite{ls:hits}.
As with ordinary inductive types, models of higher inductive types can be constructed by transfinite iterative processes; a slogan is that ordinary inductive types describe \emph{free} monads while higher inductive types describe \emph{presentations} of monads.\index{monad}%
The introduction of path constructors also involves the model-category-theoretic equivalence between ``right homotopies'' (defined using path spaces) and ``left homotopies'' (defined using cylinders) --- the fact that this equivalence is generally only up to homotopy provides a semantic reason to prefer propositional computation rules for path constructors.

Another (temporary) reason for this preference comes from the limitations of existing computer implementations.
Proof assistants\index{proof!assistant} like \Coq and \Agda have ordinary inductive types built in, but not yet higher inductive types.
We can of course introduce them by assuming lots of axioms, but this results in only propositional computation rules.
However, there is a trick due to Dan Licata which implements higher inductive types using private data types; this yields judgmental rules for point constructors but not path constructors.

The type-theoretic description of higher spheres using loop spaces and suspensions in \autoref{sec:circle,sec:suspension} is largely due to  Brunerie and  Licata; Favonia has given a type-theoretic version of the alternative description that uses $n$-dimensional paths\index{path!n-@$n$-}.
The reduction of higher paths to 1-dimensional paths with hubs and spokes (\autoref{sec:hubs-spokes}) is due to  Lumsdaine and  Shulman.
The description of truncation as a higher inductive type is due to  Lumsdaine; the $(-1)$-truncation is closely related to the ``bracket types'' of~\cite{ab:bracket-types}.
The flattening lemma was first formulated in generality by  Brunerie.

\index{set-quotient}
Quotient types are unproblematic in extensional type theory, such as \NuPRL~\cite{constable+86nuprl-book}.
They are often added by passing to an extended system of setoids.\index{setoid}
However, quotients are a trickier issue in intensional type theory (the starting point for homotopy type theory), because one cannot simply add new propositional equalities without specifying how they are to behave. Some solutions to this problem have been studied~\cite{hofmann:thesis,Altenkirch1999,altenkirch+07ott}, and several different notions of quotient types have been considered.  The construction of set-quotients using higher-inductives provides an argument for our particular approach (which is similar to some that have previously been considered), because it arises as an instance of a general mechanism.  Our construction does not yet provide a new solution to all the computational problems related to quotients, since we still lack a good computational understanding of higher inductive types in general---but it does mean that ongoing work on the computational interpretation of higher inductives applies to the quotients as well.  The construction of quotients in terms of equivalence classes is, of
course, a standard set-theoretic idea, and a well-known aspect of elementary topos theory; its use in type theory (which depends on the univalence axiom, at least for mere propositions) was proposed by Voevodsky.  The fact that quotient types in intensional type theory imply function extensionality was proved by~\cite{hofmann:thesis}, inspired by the work of~\cite{carboni} on exact completions; \autoref{thm:interval-funext} is an adaptation of such arguments.

\sectionExercises

\begin{ex}\label{ex:torus}
  Define concatenation of dependent paths, prove that application of dependent functions preserves concatenation, and write out the precise induction principle for the torus $T^2$ with its computation rules.\index{torus}
\end{ex}

\begin{ex}\label{ex:suspS1}
  Prove that $\eqv{\susp \Sn^1}{\Sn^2}$, using the explicit definition of $\Sn^2$ in terms of $\base$ and $\surf$ given in \autoref{sec:circle}.
\end{ex}

\begin{ex}\label{ex:torus-s1-times-s1}
  Prove that the torus $T^2$ as defined in \autoref{sec:cell-complexes} is equivalent to $\Sn^1\times \Sn^1$.
  (Warning: the path algebra for this is rather difficult.)
\end{ex}

\begin{ex}\label{ex:nspheres}
  Define dependent $n$-loops\index{loop!dependent n-@dependent $n$-} and the action of dependent functions on $n$-loops, and write down the induction principle for the $n$-spheres as defined at the end of \autoref{sec:circle}.
\end{ex}

\begin{ex}
  Prove that $\eqv{\susp \Sn^n}{\Sn^{n+1}}$, using the definition of $\Sn^n$ in terms of $\Omega^n$ from \autoref{sec:circle}.
\end{ex}

\begin{ex}
  Prove that if the type $\Sn^2$ belongs to some universe \type, then \type is not a 2-type.
\end{ex}

\begin{ex}
  Prove that if $G$ is a monoid and $x:G$, then $\sm{y:G}((x\cdot y = e) \times (y\cdot x =e))$ is a mere proposition.
  Conclude, using the principle of unique choice (\autoref{cor:UC}), that it would be equivalent to define a group to be a monoid such that for every $x:G$, there merely exists a $y:G$ such that $x\cdot y = e$ and $y\cdot x=e$.
\end{ex}

\begin{ex}\label{ex:free-monoid}
  Prove that if $A$ is a set, then $\lst A$ is a monoid.
  Then complete the proof of \autoref{thm:free-monoid}.\index{monoid!free}
\end{ex}

\begin{ex}\label{ex:unnatural-endomorphisms}
  Assuming \LEM{}, construct a family $f:\prd{X:\type}(X\to X)$ such that $f_\bool:\bool\to\bool$ is the nonidentity automorphism.\index{automorphism!of 2, nonidentity@of $\bool$, nonidentity}
\end{ex}

\index{type!higher inductive|)}%


\chapter{Homotopy \texorpdfstring{$n$}{n}-types}
\label{cha:hlevels}

\index{n-type@$n$-type|(}%
\indexsee{h-level}{$n$-type}

One of the basic notions of homotopy theory is that of a \emph{homotopy $n$-type}: a space containing no interesting homotopy above dimension $n$.
For instance, a homotopy $0$-type is essentially a set, containing no nontrivial paths, while a homotopy $1$-type may contain nontrivial paths, but no nontrivial paths between paths.
Homotopy $n$-types are also called \emph{$n$-truncated spaces}.
We have mentioned this notion already in \autoref{sec:basics-sets}; our first goal in this chapter is to give it a precise definition in homotopy type theory.

A dual notion to truncatedness is connectedness: a space is \emph{$n$-connected} if it has no interesting homotopy in dimensions $n$ and \emph{below}.
For instance, a space is $0$-connected (also called just ``connected'') if it has only one connected component, and $1$-connected (also called ``simply connected'') if it also has no nontrivial loops (though it may have nontrivial higher loops between loops\index{loop!n-@$n$-}).

The duality between truncatedness and connectedness is most easily seen by extending both notions to maps.
We call a map \emph{$n$-truncated} or \emph{$n$-connected} if all its fibers are so.
Then $n$-connected and $n$-truncated maps form the two classes of maps in an \emph{orthogonal factorization system},
\index{orthogonal factorization system}
\indexsee{factorization!system, orthogonal}{orthogonal factorization system}
i.e.\ every map factors uniquely as an $n$-connected map followed by an $n$-truncated one.

In the case $n={-1}$, the $n$-truncated maps are the embeddings and the $n$-connected maps are the surjections, as defined in \autoref{sec:mono-surj}.
Thus, the $n$-connected factorization system is a massive generalization of the standard image factorization of a function between sets into a surjection followed by an injection.
At the end of this chapter, we sketch briefly an even more general theory: any type-theoretic \emph{modality} gives rise to an analogous factorization system.

\section{Definition of \texorpdfstring{$n$}{n}-types}
\label{sec:n-types}

As mentioned in \autoref{sec:basics-sets,sec:contractibility}, it turns out to be convenient to define $n$-types starting two levels below zero, with the $(-1)$-types being the mere propositions and the $(-2)$-types the contractible ones.

\begin{defn}\label{def:hlevel}
  Define the predicate $\istype{n} : \type \to \type$ for $n \geq -2$ by recursion as follows:
  \[ \istype{n}(X) \defeq
  \begin{cases}
    \iscontr(X) & \text{ if } n = -2, \\
    \prd{x,y : X} \istype{n'}(\id[X]{x}{y}) & \text{ if } n = n'+1.
  \end{cases}
  \]
  We say that $X$ is an \define{$n$-type}, or sometimes that it is \emph{$n$-truncated},
  \indexdef{n-type@$n$-type}%
  \indexsee{n-truncated@$n$-truncated!type}{$n$-type}%
  \indexsee{type!n-type@$n$-type}{$n$-type}%
  \indexsee{type!n-truncated@$n$-truncated}{$n$-type}%
 if $\istype{n}(X)$ is inhabited.
\end{defn}

\begin{rmk}
  The number $n$ in \autoref{def:hlevel} ranges over all integers greater than or equal to $-2$.
  We could make sense of this formally by defining a type $\Z_{{\geq}-2}$ of such integers (a type whose induction principle is identical to that of $\nat$), or instead defining a predicate $\istype{(k-2)}$ for $k : \nat$.
  Either way, we can prove theorems about $n$-types by induction on $n$, with $n = -2$ as the base case.
\end{rmk}

\begin{eg}
  \index{set}
  We saw in \autoref{thm:prop-minusonetype} that $X$ is a $(-1)$-type if and only if it is a mere proposition.
  Therefore, $X$ is a $0$-type if and only if it is a set.
\end{eg}

We have also seen that there are types which are not sets (\autoref{thm:type-is-not-a-set}).
So far, however, we have not shown for any $n>0$ that there exist types which are not $n$-types.
In \autoref{cha:homotopy}, however, we will show that the $(n+1)$-sphere $\Sn^{n+1}$ is not an $n$-type.
(Kraus has also shown that the $n^{\mathrm{th}}$ nested univalent universe is also not an $n$-type, without using any higher inductive types.)
Moreover, in \autoref{sec:whitehead} will give an example of a type that is not an $n$-type for \emph{any} (finite) number $n$.

We begin the general theory of $n$-types by showing they are closed under certain operations and constructors.

\begin{thm}\label{thm:h-level-retracts}
  \index{retract!of a type}%
  \index{retraction}%
 Let $p : X \to Y$ be a retraction and suppose that $X$ is an $n$-type, for any $n\geq -2$.
 Then $Y$ is also an $n$-type.
\end{thm}

\begin{proof}
 We proceed by induction on $n$.
 The base case $n=-2$ is handled by \autoref{thm:retract-contr}.

 For the inductive step, assume that any retract of an $n$-type is an $n$-type, and that $X$ is an $\nplusone$-type.
 Let $y, y' : Y$; we must show that $\id{y}{y'}$ is an $n$-type.
 Let $s$ be a section of $p$, and let $\epsilon$ be a homotopy $\epsilon : p \circ s \htpy 1$.
 Since $X$ is an $\nplusone$-type, $\id[X]{s(y)}{s(y')}$ is an $n$-type.
 We claim that $\id{y}{y'}$ is a retract of $\id[X]{s(y)}{s(y')}$.
 For the section, we take
 \[ \apfunc s : (y=y') \to (s(y)=s(y')). \]
 For the retraction, we define $t:(s(y)=s(y'))\to(y=y')$ by
 \[ t(q) \defeq  \opp{\epsilon_y} \ct \ap p q \ct \epsilon_{y'}.\]
 To show that $t$ is a retraction of $\apfunc s$, we must show that
 \[ \opp{\epsilon_y} \ct \ap p {\ap sr} \ct \epsilon_{y'} = r \]
 for any $r:y=y'$.
 But this follows from \autoref{lem:htpy-natural}.
\end{proof}

As an immediate corollary we obtain the stability of $n$-types under equivalence (which is also immediate from univalence):

\begin{cor}\label{cor:preservation-hlevels-weq}
 If $\eqv{X}{Y}$ and $X$ is an $n$-type, then so is $Y$.
\end{cor}

Recall also the notion of embedding from \autoref{sec:mono-surj}.

\begin{thm}\label{thm:isntype-mono}
  \index{function!embedding}
  If $f:X\to Y$ is an embedding and $Y$ is an $n$-type for some $n\ge -1$, then so is $X$.
\end{thm}
\begin{proof}
  Let $x,x':X$; we must show that $\id[X]{x}{x'}$ is an $\nminusone$-type.
  But since $f$ is an embedding, we have $(\id[X]{x}{x'}) \eqvsym (\id[Y]{f(x)}{f(x')})$, and the latter is an $\nminusone$-type by assumption.
\end{proof}

Note that this theorem fails when $n=-2$: the map $\emptyt \to \unit$ is an embedding, but $\unit$ is a $(-2)$-type while $\emptyt$ is not.

\begin{thm}\label{thm:hlevel-cumulative}
 The hierarchy of $n$-types is cumulative in the following sense:
   given a number $n \geq -2$, if $X$ is an $n$-type, then it is also an $\nplusone$-type.
\end{thm}

\begin{proof}
 We proceed by induction on $n$.

 For $n = -2$, we need to show that a contractible type, say, $A$, has contractible path spaces.
       Let $a_0: A$ be the center of contraction of $A$, and let $x, y : A$. We show that $\id[A]{x}{y}$
       is contractible.
       By contractibility of $A$ we have a path $\contr_x \ct \opp{\contr_y} : x = y$, which we choose as
       the center of contraction for $\id{x}{y}$.
       Given any $p : x = y$, we need to show $p = \contr_x \ct \opp{\contr_y}$.
           By path induction, it suffices to show that
        $\refl{x} = \contr_x \ct \opp{\contr_x}$, which is trivial.

 For the inductive step, we need to show that $\id[X]{x}{y}$ is an $\nplusone$-type, provided
          that $X$ is an $\nplusone$-type. Applying the inductive hypothesis to $\id[X]{x}{y}$
         yields the desired result.
\end{proof}


We now show that  $n$-types are preserved by most of the type forming operations.

\begin{thm}\label{thm:ntypes-sigma}
 Let $n \geq -2$, and let $A : \type$ and $B : A \to \type$.
 If $A$ is an $n$-type and for all $a : A$, $B(a)$ is an $n$-type, then so is $\sm{x : A} B(x)$.
\end{thm}

\begin{proof}
 We proceed by induction on $n$.

 For $n = -2$, we choose the center of contraction for $\sm{x : A} B(x)$ to be the pair
       $(a_0, b_0)$, where $a_0 : A$ is the center of contraction of $A$ and $b_0 : B(a_0)$ is the center of contraction of $B(a_0)$.
       Given any other element $(a,b)$ of $\sm{x : A} B(x)$, we provide a path $\id{(a, b)}{(a_0,b_0)}$
       by contractibility of $A$ and $B(a_0)$, respectively.

 For the inductive step, suppose that $A$ is an $\nplusone$-type and
         for any $a : A$, $B(a)$ is an $\nplusone$-type. We show that $\sm{x : A} B(x)$ is an $\nplusone$-type:
      fix $(a_1, b_1)$ and $(a_2,b_2)$ in $\sm{x : A} B(x)$,
     we show that $\id{(a_1, b_1)}{(a_2,b_2)}$ is an $n$-type.
      By \autoref{thm:path-sigma} we have
      \[ \eqvspaced{(\id{(a_1, b_1)}{(a_2,b_2)})}{\sm{p : \id{a_1}{a_2}} (\id[B(a_2)]{\trans{p}{b_1}}{b_2})} \]
   and by preservation of $n$-types under equivalences (\autoref{cor:preservation-hlevels-weq})
   it suffices to prove that the latter is an $n$-type. This follows from the
   inductive hypothesis.
\end{proof}

As a special case, if $A$ and $B$ are $n$-types, so is $A\times B$.
Note also that \autoref{thm:hlevel-cumulative} implies that if $A$ is an $n$-type, then so is $\id[A]xy$ for any $x,y:A$.
Combining this with \autoref{thm:ntypes-sigma}, we see that for any functions $f:A\to C$ and $g:B\to C$ between $n$-types, their pullback\index{pullback}
\[ A\times_C B \defeq \sm{x:A}{y:B} (f(x)=g(y)) \]
(see \autoref{ex:pullback}) is also an $n$-type.
More generally, $n$-types are closed under all \emph{limits}.

\begin{thm}\label{thm:hlevel-prod}
 Let $n\geq -2$, and let $A : \type$ and $B : A \to \type$.
 If for all $a : A$, $B(a)$ is an $n$-type, then so is $\prd{x : A} B(x)$.
\end{thm}

\begin{proof}
  We proceed by induction on $n$.
  For $n = -2$, the result is simply \autoref{thm:contr-forall}.

  For the inductive step, assume the result is true for $n$-types, and that each $B(a)$ is an $\nplusone$-type.
  Let $f, g : \prd{a:A}B(a)$.
  We need to show that $\id{f}{g}$ is an $n$-type.
  By function extensionality and closure of $n$-types under equivalence, it suffices to show that $\prd{a : A} (\id[B(a)]{f(a)}{g(a)})$ is an $n$-type.
  This follows from the inductive hypothesis.
\end{proof}

As a special case of the above theorem, the function space $A \to B$ is an $n$-type provided that $B$ is an $n$-type.
We can now generalize our observations in \autoref{cha:basics} that $\isset(A)$ and $\isprop(A)$ are mere propositions.

\begin{thm}\label{thm:isaprop-isofhlevel}
 For any $n \geq -2$ and any type $X$, the type $\istype{n}(X)$ is a mere proposition.
\end{thm}
\begin{proof}
  We proceed by induction with respect to $n$.

 For the base case, we need to show that for any $X$, the type $\iscontr(X)$ is a mere proposition.
 This is \autoref{thm:isprop-iscontr}.

For the inductive step we need to show
\[\prd{X : \type} \isprop (\istype{n}(X)) \to \prd{X : \type} \isprop (\istype{\nplusone}(X)) \]
To show the conclusion of this implication, we need to show that for any type $X$, the type
\[\prd{x, x' : X}\istype{n}(x = x')\]
is a mere proposition. By \autoref{thm:isprop-forall} or \autoref{thm:hlevel-prod}, it suffices to show that for any $x, x' : X$, the type $\istype{n}(x =_X x')$ is a mere
proposition.
But this follows from the inductive hypothesis applied to the type $(x =_X x')$.
\end{proof}

Finally, we show that the type of $n$-types is itself an $\nplusone$-type.
We define this to be:
\symlabel{universe-of-ntypes}
\[\ntype{n} \defeq \sm{X : \type} \istype{n}(X) \]
If necessary, we may specify the universe $\UU$ by writing $\ntypeU{n}$.
In particular, we have $\prop \defeq \ntype{(-1)}$ and $\set \defeq \ntype{0}$, as defined in \autoref{cha:basics}.
Note that just as for \prop and \set, because $\istype{n}(X)$ is a mere proposition, by \autoref{thm:path-subset} for any $(X,p), (X',p'):\ntype{n}$ we have
\begin{align*}
  \Big(\id[\ntype{n}]{(X, p)}{(X', p')}\Big) &\eqvsym (\id[\type] X X')\\
  &\eqvsym (\eqv{X}{X'}).
\end{align*}

\begin{thm}\label{thm:hleveln-of-hlevelSn}
 For any $n \geq -2$, the type $\ntype{n}$ is an $\nplusone$-type.
\end{thm}
\begin{proof}
  Let $(X, p), (X', p') : \ntype{n}$; we need to show that $\id{(X, p)}{(X', p')}$ is an $n$-type.
  By the above observation, this type is equivalent to $\eqv{X}{X'}$.
  Next, we observe that the projection
  \[(\eqv{X}{X'}) \to (X \rightarrow X').\]
  is an embedding, so that if $n\geq -1$, then by \autoref{thm:isntype-mono} it suffices to show that $X \rightarrow X'$ is an $n$-type.
  But since $n$-types are preserved under the arrow type, this reduces to an assumption that $X'$ is an $n$-type.

  In the case $n=-2$, this argument shows that $\eqv{X}{X'}$ is a $(-1)$-type --- but it is also inhabited, since any two contractible types
are equivalent to \unit, and hence to each other.
  Thus, $\eqv{X}{X'}$ is also a $(-2)$-type.
\end{proof}

\section{Uniqueness of identity proofs and Hedberg's theorem}
\label{sec:hedberg}

\index{set|(}%

In \autoref{sec:basics-sets} we defined a type $X$ to be a \emph{set} if for all $x, y : X$ and $p, q : x =_X y$ we have $p = q$.
In conventional type theory, this property goes by the name of \define{uniqueness of identity proofs (UIP)}.
\indexdef{uniqueness!of identity proofs}%
We have seen also that it is equivalent to being a $0$-type in the sense of the previous section.
Here is another equivalent characterization, involving Streicher's ``Axiom K'' \cite{Streicher93}:

\begin{thm}\label{thm:h-set-uip-K}
 A type $X$ is a set if and only if it satisfies \define{Axiom K}:
 \indexdef{axiom!Streicher's Axiom K}%
 for all $x : X$ and $p : (x =_A x)$ we have $p = \refl{x}$.
\end{thm}

\begin{proof}
  Clearly Axiom K is a special case of UIP.
  Conversely, if $X$ satisfies Axiom K, let $x, y : X$ and $p, q : (\id{x}{y})$; we want to show $p=q$.
  But induction on $q$ reduces this goal precisely to Axiom K.
\end{proof}

We stress that \emph{we} are not assuming UIP or the K principle as axioms!
They are simply properties which a particular type may or may not satisfy (which are equivalent to being a set).
Recall from \autoref{thm:type-is-not-a-set} that \emph{not} all types are sets.

The following theorem is another useful way to show that types are sets.

\begin{thm}\label{thm:h-set-refrel-in-paths-sets}
  \index{relation!reflexive}%
  Suppose $R$ is a reflexive\index{reflexivity!of a relation} mere relation on a type $X$ implying identity.
  Then $X$ is a set, and $R(x,y)$ is equivalent to $\id[X]{x}{y}$ for all $x,y:X$.
\end{thm}

\begin{proof}
  Let $\rho : \prd{x:X} R(x,x)$ witness reflexivity of $R$, and let \narrowequation{f : \prd{x,y:X} R(x,y) \to (\id[X]{x}{y})} be a witness that $R$
implies identity.
  Note first that the two statements in the theorem are equivalent.
  For on one hand, if $X$ is a set, then $\id[X]xy$ is a mere proposition, and since it is logically equivalent to the mere proposition
$R(x,y)$ by hypothesis, it must also be equivalent to it.
  On the other hand, if $\id[X]xy$ is equivalent to $R(x,y)$, then like the latter it is a mere proposition for all $x,y:X$, and hence $X$
is a set.

  We give two proofs of this theorem.
  The first shows directly that $X$ is a set; the second shows directly that $R(x,y)\eqvsym (x=y)$.

  \emph{First proof:} we show that $X$ is a set.
  The idea is the same as that of \autoref{thm:prop-set}: the function $f$ must be continuous in its arguments $x$ and $y$.
  However, it is slightly more notationally complicated because we have to deal with the additional argument of type $R(x,y)$.

  Firstly, for any $x:X$ and $p:\id[X]xx$, consider $\apdfunc{f(x)}(p)$.
  This is a dependent path from $f(x,x)$ to itself.
  Since $f(x,x)$ is still a function $R(x,x) \to (\id[X]xy)$, by \autoref{thm:dpath-arrow} this yields for any $r:R(x,x)$ a path
  \[\trans{p}{f(x,x,r)} = f(x,x,\trans{p}r).
  \]
  On the left-hand side, we have transport in an identity type, which is concatenation.
  And on the right-hand side, we have $\trans{p}r = r$, since both lie in the mere proposition $R(x,x)$.
  Thus, substituting $r\defeq \rho(x)$, we obtain
  \[ f(x,x,\rho(x)) \ct p = f(x,x,\rho(x)). \]
  By cancellation, $p=\refl{x}$.
  So $X$ satisfies Axiom K, and hence is a set.

  \emph{Second proof:} we show that each $f(x,y) : R(x,y) \to \id[X]{x}{y}$ is an equivalence.
  By \autoref{thm:total-fiber-equiv}, it suffices to show that $f$ induces an equivalence of total spaces:
  \begin{equation*}
    \eqv{\Parens{\sm{y:X}R(x,y)}}{\Parens{\sm{y:X}\id[X]{x}{y}}}.
  \end{equation*}
  By \autoref{thm:contr-paths}, the type on the right is contractible, so it
  suffices to show that the type on the left is contractible. As the center of
  contraction we take the pair $\pairr{x,\rho(x)}$.  It remains to show, for
  every ${y:X}$ and every ${H:R(x,y)}$ that
  \begin{equation*}
    \id{\pairr{x,\rho(x)}}{\pairr{y,H}}.
  \end{equation*}
  But since $R(x,y)$ is a mere proposition, by \autoref{thm:path-sigma} it suffices to show that
  $\id[X]{x}{y}$, which we get from $f(H)$.
\end{proof}

\begin{cor}\label{notnotstable-equality-to-set}
  If a type $X$ has the property that $\neg\neg(x=y)\to(x=y)$ for any $x,y:X$, then $X$ is a set.
\end{cor}

Another convenient way to show that a type is a set is the following.
Recall from \autoref{sec:intuitionism} that a type $X$ is said to have \emph{decidable equality}
\index{decidable!equality|(}%
if for all $x, y : X$ we have
\[(x =_X y) + \neg (x =_X y).\]
\index{continuity of functions in type theory@``continuity'' of functions in type theory}%
\index{functoriality of functions in type theory@``functoriality'' of functions in type theory}%
This is a very strong condition: it says that a path $x=y$ can be chosen, when it exists, continuously (or computably, or functorially) in $x$ and $y$.
This turns out to imply that $X$ is a set, by way of \autoref{thm:h-set-refrel-in-paths-sets} and the following lemma.

\begin{lem}
For any type $A$ we have $(A+\neg A)\to(\neg\neg A\to A)$.
\end{lem}

\begin{proof}
Suppose $x:A+\neg A$. We have two cases to consider.
If $x$ is $\inl(a)$ for some $a:A$, then we have the constant function $\neg\neg A
\to A$ which maps everything to $a$. If $x$ is $\inr(f)$ for some $t:\neg A$,
we have $g(t):\emptyt$ for every $g:\neg\neg A$. Hence we may use
\emph{ex falso quodlibet}, that is $\rec{\emptyt}$, to obtain an element of $A$ for any $g:\neg\neg A$.
\end{proof}

\index{anger}
\begin{thm}[Hedberg]\label{thm:hedberg}
  \index{Hedberg's theorem}%
  \index{theorem!Hedberg's}%
  If $X$ has decidable equality, then $X$ is a set.
\end{thm}

\begin{proof}
If $X$ has decidable equality, it follows that $\neg\neg(x=y)\to(x=y)$ for any
$x,y:X$. Therefore, Hedberg's theorem follows from 
\autoref{notnotstable-equality-to-set}.
\end{proof}

There is, of course, a strong connection between this theorem and \autoref{thm:not-lem}.
The statement \LEM{\infty} that is denied by \autoref{thm:not-lem} clearly implies that every type has decidable equality, and hence is a set, which we know is not the case.
\index{excluded middle}%
Note that the consistent axiom \LEM{} from \autoref{sec:intuitionism} implies only that every type has \emph{merely decidable equality}, i.e.\ that for any $A$ we have
\indexdef{equality!merely decidable}%
\indexdef{merely!decidable equality}%
\[ \prd{a,b:A} (\brck{a=b} + \neg\brck{a=b}). \]

\index{decidable!equality|)}%

As an example application of \autoref{thm:hedberg}, recall that in \autoref{thm:nat-set} we observed that $\nat$ is a set, using our characterization of its equality types in
\autoref{sec:compute-nat}.
A more traditional proof of this theorem uses only~\eqref{eq:zero-not-succ} and~\eqref{eq:suc-injective}, rather than the full
characterization of \autoref{thm:path-nat}, with \autoref{thm:hedberg} to fill in the blanks.

\begin{thm}\label{prop:nat-is-set}
 The type $\nat$ of natural numbers has decidable equality, and hence is a set.
\end{thm}

\begin{proof}
  Let $x, y : \nat$ be given; we proceed by induction on $x$ and case analysis on $y$ to prove $(x=y)+\neg(x=y)$.
  If $x \jdeq 0$ and $y \jdeq 0$, we take $\inl(\refl{}(0))$.
  If $x \jdeq 0$ and $y \jdeq \suc(n)$, then by~\eqref{eq:zero-not-succ} we get $\neg (0 = \suc (n))$.

  For the inductive step, let $x \jdeq \suc (n)$.
  If $y \jdeq 0$, we use~\eqref{eq:zero-not-succ} again.
  Finally, if $y \jdeq \suc (m)$, the inductive hypothesis gives $(m = n)+\neg(m = n)$.
  In the first case, if $p:m=n$, then $\ap \suc p:\suc(m)=\suc(n)$.
  And in the second case,~\eqref{eq:suc-injective} yields $\neg(\suc(m)=\suc(n))$.
\end{proof}

\index{set|)}%

\index{axiom!Streicher's Axiom K!generalization to n-types@generalization to $n$-types}%
Although Hedberg's theorem appears rather special to sets ($0$-types), ``Axiom K'' generalizes naturally to $n$-types.
Note that the ordinary Axiom K (as a property of a type $X$) states that for all $x:X$, the loop space\index{loop space} $\Omega(X,x)$ (see \cref{def:loopspace}) is contractible.
Since $\Omega(X,x)$ is always inhabited (by $\refl{x}$), this is equivalent to its being a mere proposition (a $(-1)$-type).
Since $0 = (-1)+1$, this suggests the following generalization.

\begin{thm}\label{thm:hlevel-loops}
  For any $n\geq -1$, a type $X$ is an $\nplusone$-type if and only if for all $x : X$, the type $\Omega(X, x)$ is an $n$-type.
\end{thm}

Before proving this, we prove an auxiliary lemma:

\begin{lem}\label{lem:hlevel-if-inhab-hlevel}
  Given $n \geq -1$ and $X : \type$.
  If, given any inhabitant of $X$ it follows that $X$ is an $n$-type, then $X$ is an $n$-type.
\end{lem}
\begin{proof}
  Let $f : X \to \istype{n}(X)$ be the given map.
  We need to show that for any $x, x' : X$, the type $\id{x}{x'}$ is an $\nminusone$-type.
  But then $f(x)$ shows that $X$ is an $n$-type, hence all its path spaces are $\nminusone$-types.
\end{proof}

\begin{proof}[Proof of \autoref{thm:hlevel-loops}]
  The ``only if'' direction is obvious, since $\Omega(X,x)\defeq (\id[X]xx)$.
  Conversely, in order to show that $X$ is an $\nplusone$-type, we need to show that for any $x, x' : X$, the type $\id{x}{x'}$ is an
$n$-type.
  Following \autoref{lem:hlevel-if-inhab-hlevel} it suffices to give a map
  \[ (\id{x}{x'}) \to \istype{n}(\id{x}{x'}). \]
  By path induction, it suffices to do this when $x\jdeq x'$, in which case it follows from the assumption that $\Omega(X, x)$ is an
$n$-type.
\end{proof}

\index{whiskering}
By induction and some slightly clever whiskering, we can obtain a generalization of the K property to $n>0$.

\begin{thm}\label{thm:ntype-nloop}
  \index{loop space!iterated}%
  For every $n\ge -1$, a type $A$ is an $n$-type if and only if $\Omega^{n+1}(A,a)$ is contractible for all $a:A$.
\end{thm}
\begin{proof}
  Recalling that $\Omega^0(A,a) = (A,a)$, the case $n=-1$ is \autoref{ex:prop-inhabcontr}.
  The case $n=0$ is \autoref{thm:h-set-uip-K}.
  Now we use induction; suppose the statement holds for $n:\N$.
  By \autoref{thm:hlevel-loops}, $A$ is an $(n+1)$-type iff $\Omega(A,a)$ is an $n$-type for all $a:A$.
  By the inductive hypothesis, the latter is equivalent to saying that $\Omega^{n+1}(\Omega(A,a),p)$ is contractible for all $p:\Omega(A,a)$.

  Since $\Omega^{n+2}(A,a) \defeq \Omega^{n+1}(\Omega(A,a),\refl{a})$, and $\Omega^{n+1} = \Omega^n \circ \Omega$, it will suffice to show that $\Omega(\Omega(A,a),p)$ is equal to $\Omega(\Omega(A,a),\refl{a})$, in the type $\pointed\type$ of pointed types.
  For this, it suffices to give an equivalence
  \[ g : \Omega(\Omega(A,a),p) \eqvsym \Omega(\Omega(A,a),\refl{a}) \]
  which carries the basepoint $\refl{p}$ to the basepoint $\refl{\refl{a}}$.
  For $q:p=p$, define $g(q):\refl{a} = \refl{a}$ to be the following composite:
  \[ \refl{a} = p\ct \opp p \overset{q}{=} p\ct\opp p = \refl{a}, \]
  where the path labeled ``$q$'' is actually $\apfunc{\lam{r} r\ct\opp p} (q)$.
  Then $g$ is an equivalence because it is a composite of equivalences
  \[ (p=p) \xrightarrow{\apfunc{\lam{r} r\ct\opp p}} (p\ct \opp p = p\ct \opp p) \xrightarrow{i\ct - \ct \opp i} (\refl{a} = \refl{a}). \]
  using \autoref{eg:concatequiv,thm:paths-respects-equiv}, where $i:\refl{a} = p\ct \opp p$ is the canonical equality.
  And it is evident that $g(\refl{p}) = \refl{\refl{a}}$.
\end{proof}

\section{Truncations}
\label{sec:truncations}

\indexsee{n-truncation@$n$-truncation}{truncation}%
\index{truncation!n-truncation@$n$-truncation|(defstyle}%

In \autoref{subsec:prop-trunc} we introduced the propositional truncation, which makes the ``best approximation'' of a type that is a mere
proposition, i.e.\ a $(-1)$-type.
In \autoref{sec:hittruncations} we constructed this truncation as a higher inductive type, and gave one way to generalize it to a
0-truncation.
We now explain a better generalization of this, which truncates any type into an $n$-type for any $n\geq -2$; in classical homotopy theory this would be called its \define{$n^{\mathrm{th}}$ Postnikov section}.\index{Postnikov tower}

The idea is to make use of \autoref{thm:ntype-nloop}, which states that $A$ is an $n$-type just when $\Omega^{n+1}(A,a)$
\index{loop space!iterated}%
is contractible for
all $a:A$, and \autoref{lem:susp-loop-adj}, which implies that
\narrowequation{\Omega^{n+1}(A,a) \eqvsym \Map_{*}(\Sn^{n+1},(A,a)),} where $\Sn^{n+1}$ is equipp\-ed with some basepoint which we may as well call \base.
However, contractibility of $\Map_*(\Sn^{n+1},(A,a))$ is something that we can ensure directly by giving path constructors.

We might first of all try to define $\trunc nA$ to be generated by a function $\tprojf n : A \to \trunc n A$, together with for each
$r:\Sn^{n+1} \to \trunc n A$ and each $x:\Sn^{n+1}$, a path $s_r(x):r(x) = r(\base)$.
But this does not quite work, for the same reason that \autoref{rmk:spokes-no-hub} fails.
\index{hub and spoke}%
Instead, we use the full ``hub and spoke'' construction as in \autoref{sec:hubs-spokes}.

Thus, for $n\ge -1$, we take $\trunc nA$ to be the higher inductive type generated by:
\begin{itemize}
\item a function $\tprojf n : A \to \trunc n A$,
\item for each $r:\Sn^{n+1} \to \trunc n A$, a \emph{hub} point $h(r):\trunc n A$, and
\item for each $r:\Sn^{n+1} \to \trunc n A$ and each $x:\Sn^{n+1}$, a \emph{spoke} path $s_r(x):r(x) = h(r)$.
\end{itemize}

\noindent
The existence of these constructors is now enough to show:

\begin{lem}
  $\trunc n A$ is an $n$-type.
\end{lem}
\begin{proof}
  By \autoref{thm:ntype-nloop}, it suffices to show that $\Omega ^{n+1}(\trunc nA,b)$ is contractible for all $b:\trunc nA$, which by
\autoref{lem:susp-loop-adj} is equivalent to \narrowequation{\Map_*(\Sn^{n+1},(\trunc nA,b)).}
  As center of contraction for the latter, we choose the function $c_b:\Sn^{n+1} \to \trunc nA$ which is constant at $b$, together with
$\refl b : c_b(\base) = b$.

  Now, an arbitrary element of $\Map_*(\Sn^{n+1},(\trunc nA,b))$ consists of a map $r:\Sn^{n+1} \to \trunc n A$ together with a path
$p:r(\base)=b$.
  By function extensionality, to show $r = c_b$ it suffices to give, for each $x:\Sn^{n+1}$, a path $r(x)=c_b(x) \jdeq b$.
  We choose this to be the composite $s_r(x) \ct \opp{s_r(\base)} \ct p$, where $s_r(x)$ is the spoke at $x$.

  Finally, we must show that when transported along this equality $r=c_b$, the path $p$ becomes $\refl b$.
  By transport in path types, this means we need
  \[\opp{(s_r(\base) \ct \opp{s_r(\base)} \ct p)} \ct p = \refl b.\]
  But this is immediate from path operations.
\end{proof}

(This construction fails for $n=-2$, but in that case we can simply define $\trunc{-2}{A}\defeq \unit$ for all $A$.
From now on we assume $n\ge -1$.)

\index{induction principle!for truncation}%
To show the desired universal property of the $n$-truncation, we need the induction principle.
We extract this from the constructors in the usual way; it says that given $P:\trunc nA\to\type$ together with
\begin{itemize}
\item For each $a:A$, an element $g(a) : P(\tproj na)$,
\item For each $r:\Sn^{n+1} \to \trunc n A$ and $r':\prd{x:\Sn^{n+1}} P(r(x))$, an element $h'(r,r'):P(h(r))$,
\item For each $r:\Sn^{n+1} \to \trunc n A$ and $r':\prd{x:\Sn^{n+1}} P(r(x))$, and each $x:\Sn^{n+1}$, a dependent path
$\dpath{P}{s_r(x)}{r'(x)}{h'(r,r')}$,
\end{itemize}
there exists a section $f:\prd{x:\trunc n A} P(x)$ with $f(\tproj n a) \jdeq g(a)$ for all $a:A$.
To make this more useful, we reformulate it as follows.

\begin{thm}\label{thm:truncn-ind}
  For any type family $P:\trunc n A \to \type$ such that each $P(x)$ is an $n$-type, and any function $g : \prd{a:A} P(\tproj n a)$, there
exists a section $f:\prd{x:\trunc n A} P(x)$ such that $f(\tproj n a)\defeq g(a)$ for all $a:A$.
\end{thm}
\begin{proof}
  It will suffice to construct the second and third data listed above, since $g$ has exactly the type of the first datum.
  Given $r:\Sn^{n+1} \to \trunc n A$ and $r':\prd{x:\Sn^{n+1}} P(r(x))$, we have $h(r):\trunc n A$ and $s_r :\prd{x:\Sn^{n+1}} (r(x) =
h(r))$.
  Define $t:\Sn^{n+1} \to P(h(r))$ by $t(x) \defeq \trans{s_r(x)}{r'(x)}$.
  Then since $P(h(r))$ is $n$-truncated, there exists a point $u:P(h(r))$ and a contraction $v:\prd{x:\Sn^{n+1}} (t(x) = u)$.
  Define $h'(r,r') \defeq u$, giving the second datum.
  Then (recalling the definition of dependent paths), $v$ has exactly the type required of the third datum.
\end{proof}

In particular, if $E$ is some $n$-type, we can consider the constant family of types equal to $E$ for every point of $A$.
\symlabel{extend}
\index{recursion principle!for truncation}%
Thus, every map $f:A\to{}E$ can be extended to a map $\extend{f}:\trunc nA\to{}E$ defined by $\extend{f}(\tproj na)\defeq f(a)$; this is the \emph{recursion principle} for $\trunc n A$.

The induction principle also implies a uniqueness principle for functions of this form.
\index{uniqueness!principle, propositional!for functions on a truncation}%
Namely, if $E$ is an $n$-type and $g,g':\trunc nA\to{}E$ are such
that $g(\tproj na)=g'(\tproj na)$ for every $a:A$, then $g(x)=g'(x)$ for all $x:\trunc nA$, since the type $g(x)=g'(x)$ is an $n$-type.
Thus, $g=g'$.
This yields the following universal property.

\begin{lem}[Universal property of truncations]\label{thm:trunc-reflective}
  \index{universal!property!of truncation}%
  Let $n\ge-2$, $A:\type$ and $B:\typele{n}$. The following map is an
  equivalence:
  \[\function{(\trunc nA\to{}B)}{(A\to{}B)}{g}{g\circ\tprojf n}\]
\end{lem}

\begin{proof}
  Given that $B$ is $n$-truncated, any $f:A\to{}B$ can be extended to a map $\extend{f}:\trunc nA\to{}B$.
  The map $\extend{f}\circ\tprojf n$ is equal to $f$, because for every $a:A$ we have $\extend{f}(\tproj na)=f(a)$ by definition.
  And the map $\extend{g\circ\tprojf n}$ is equal to $g$, because they both send $\tproj na$ to $g(\tproj na)$.
\end{proof}

In categorical language, this says that the $n$-types form a \emph{reflective subcategory} of the category of types.
\index{reflective!subcategory}%
(To state this fully precisely, one ought to use the language of $(\infty,1)$-categories.)
\index{.infinity1-category@$(\infty,1)$-category}%
In particular, this implies that the $n$-truncation is functorial:
given $f:A\to B$, applying the recursion principle to the composite $A\xrightarrow{f} B \to \trunc n B$ yields a map $\trunc n f: \trunc n A \to \trunc n B$.
By definition, we have a homotopy
\begin{equation}
  \mathsf{nat}^f_n : \prd{a:A} \trunc n f(\tproj n a) = \tproj n {f(a)},\label{eq:trunc-nat}
\end{equation}
expressing \emph{naturality} of the maps $\tprojf n$.

Uniqueness implies functoriality laws such as $\trunc n {g\circ f} = \trunc n g \circ \trunc n f$ and $\trunc n{\idfunc[A]} = \idfunc[\trunc n A]$, with attendant coherence laws.
We also have higher functoriality, for instance:

\begin{lem}\label{thm:trunc-htpy}
  Given $f,g:A\to B$ and a homotopy $h:f\htpy g$, there is an induced homotopy $\trunc n h : \trunc n f \htpy \trunc n g$ such that the composite
  \begin{equation}
    \xymatrix@C=3.6pc{\tproj n{f(a)} \ar@{=}[r]^-{\opp{\mathsf{nat}^f_n(a)}} &
      \trunc n f(\tproj n a) \ar@{=}[r]^-{\trunc n h(\tproj na)} &
      \trunc n g(\tproj n a) \ar@{=}[r]^-{\mathsf{nat}^g_n(a)} &
      \tproj n{g(a)}}\label{eq:trunc-htpy}
  \end{equation}
  is equal to $\apfunc{\tprojf n}(h(a))$.
\end{lem}
\begin{proof}
  First, we indeed have a homotopy with components $\apfunc{\tprojf n}(h(a)) : \tproj n{f(a)} = \tproj n{g(a)}$.
  Composing on either sides with the paths $\tproj n{f(a)} = \trunc n f(\tproj n a)$ and $\tproj n{g(a)} = \trunc n g(\tproj n a)$, which arise from the definitions of $\trunc n f$ and $\trunc ng$, we obtain a homotopy $(\trunc n f \circ \tprojf n) \htpy (\trunc n g \circ \tprojf n)$, and hence an equality by function extensionality.
  But since $(\blank\circ \tprojf n)$ is an equivalence, there must be a path $\trunc nf = \trunc ng$ inducing it, and the coherence laws for function extensionality imply~\eqref{eq:trunc-htpy}.
\end{proof}

The following observation about reflective subcategories is also standard.

\begin{cor}
  A type $A$ is an $n$-type if and only if $\tprojf n : A \to \trunc n A$ is an equivalence.
\end{cor}
\begin{proof}
  ``If'' follows from closure of $n$-types under equivalence.
  On the other hand, if $A$ is an $n$-type, we can define $\ext(\idfunc[A]):\trunc n A\to{}A$.
  Then we have $\ext(\idfunc[A])\circ\tprojf n=\idfunc[A]:A\to{}A$ by
  definition.  In order to prove that
  $\tprojf n\circ\ext(\idfunc[A])=\idfunc[\trunc nA]$, we only need to prove
  that $\tprojf n\circ\ext(\idfunc[A])\circ\tprojf n=
  \idfunc[\trunc nA]\circ\tprojf n$.
  This is again true:
  \[\raisebox{\depth-\height+1em}{\xymatrix{
    A \ar^-{\tprojf n}[r] \ar_{\idfunc[A]}[rd] &
    \trunc nA \ar^>>>{\ext(\idfunc[A])}[d] \ar@/^40pt/^{\idfunc[\trunc nA]}[dd] \\
    & A \ar_{\tprojf n}[d] \\
    & \trunc nA}}
  \qedhere\]
\end{proof}

The category of $n$-types also has some special properties not possessed by all reflective subcategories.
For instance, the reflector $\trunc n-$ preserves finite products.

\begin{thm}\label{cor:trunc-prod}
  For any types $A$ and $B$, the induced map $\trunc n{A\times B} \to \trunc nA \times \trunc nB$ is an equivalence.
\end{thm}
\begin{proof}
  It suffices to show that $\trunc nA \times \trunc nB$ has the same universal property as $\trunc n{A\times B}$.
  Thus, let $C$ be an $n$-type; we have
  \begin{align*}
    (\trunc nA \times \trunc nB \to C)
    &= (\trunc nA \to (\trunc nB \to C))\\
    &= (\trunc nA \to (B \to C))\\
    &= (A \to (B \to C))\\
    &= (A \times B \to C)
  \end{align*}
  using the universal properties of $\trunc nB$ and $\trunc nA$, along with the fact that $B\to C$ is an $n$-type since $C$ is.
  It is straightforward to verify that this equivalence is given by composing with $\tprojf n \times \tprojf n$, as needed.
\end{proof}

The following related fact about dependent sums is often useful.

\begin{thm}\label{thm:trunc-in-truncated-sigma}
Let $P:A\to\type$ be a family of types. Then there is an equivalence
\begin{equation*}
\eqv{\Trunc n{\sm{x:A}\trunc n{P(x)}}}{\Trunc n{\sm{x:A}P(x)}}.
\end{equation*}
\end{thm}

\begin{proof}
We use the induction principle of $n$-truncation several times to construct
functions
\begin{align*}
\varphi & : \Trunc n{\sm{x:A}\trunc n{P(x)}}\to\Trunc n{\sm{x:A}P(x)}\\
\psi & : \Trunc n{\sm{x:A}P(x)}\to \Trunc n{\sm{x:A} \trunc n{P(x)}}
\end{align*}
and homotopies $H:\varphi\circ\psi\htpy \idfunc$ and $K:\psi\circ\varphi\htpy
\idfunc$ exhibiting them as quasi-inverses.
We define $\varphi$ by setting $\varphi(\tproj n{\pairr{x,\tproj nu}})\defeq\tproj n{\pairr{x,u}}$.
We define $\psi$ by setting $\psi(\tproj n{\pairr{x,u}})\defeq\tproj n{\pairr{x,\tproj nu}}$.
Then we define $H(\tproj n{\pairr{x,u}})\defeq \refl{\tproj n{\pairr{x,u}}}$ and
$K(\tproj n{\pairr{x,\tproj nu}})\defeq \refl{\tproj n{\pairr{x,\tproj nu}}}$.
\end{proof}

\begin{cor}\label{thm:refl-over-ntype-base}
  If $A$ is an $n$-type and $P:A\to\type$ is any type family, then
  \[ \eqv{\sm{a:A} \trunc n{P(a)}}{\Trunc n{\sm{a:A}P(a)}} \]
\end{cor}
\begin{proof}
  If $A$ is an $n$-type, then the left-hand type above is already an $n$-type, hence equivalent to its $n$-truncation; thus this follows from \autoref{thm:trunc-in-truncated-sigma}.
\end{proof}

We can characterize the path spaces of a truncation using the same method that we used in \autoref{sec:compute-coprod,sec:compute-nat} for
coproducts and natural numbers (and which we will use in \autoref{cha:homotopy} to calculate homotopy groups).
Unsurprisingly, the path spaces in the $(n+1)$-truncation of $A$ are the $n$-truncations of the path spaces of $A$.
Indeed, for any $x,y:A$ there is a canonical map
\begin{equation}
  f:\ttrunc n{x=_Ay}\to \Big(\tproj {n+1}x=_{\trunc{n+1}A}\tproj {n+1}y\Big)\label{eq:path-trunc-map}
\end{equation}
defined by
\[f(\tproj n{p})\defeq \apfunc{\tprojf {n+1}}(p). \]
This definition uses the recursion principle for $\truncf n$, which is correct because $\trunc {n+1}A$ is $(n+1)$-truncated, so that the
codomain of $f$ is $n$-truncated.

\begin{thm} \label{thm:path-truncation}
  For any $A$ and $x,y:A$ and $n\ge -2$, the map~\eqref{eq:path-trunc-map} is an equivalence; thus we have
  \[ \eqv{\ttrunc n{x=_Ay}}{\Big(\tproj {n+1}x=_{\trunc{n+1}A}\tproj {n+1}y\Big)}. \]
\end{thm}

\begin{proof}
  The proof is a simple application of the encode-decode method: 
  As in previous situations, we cannot directly define a quasi-inverse to the map~\eqref{eq:path-trunc-map} because there is no way to induct on an
equality between $\tproj {n+1}x$ and $\tproj {n+1}y$.
  Thus, instead we generalize its type, in order to have general elements of the type $\trunc{n+1}A$ instead of $\tproj {n+1}x$ and $\tproj
{n+1}y$.
  Define $P:\trunc {n+1}A\to\trunc {n+1}A\to\typele{n}$ by
  \[P(\tproj {n+1}x,\tproj {n+1}y)\defeq \trunc n{x=_Ay}\]
  This definition is correct because $\trunc n{x=_Ay}$ is $n$-truncated, and $\typele{n}$ is $(n+1)$-truncated by
\autoref{thm:hleveln-of-hlevelSn}.
  Now for every $u,v:\trunc{n+1}A$, there is a map
  \[\encode:P(u,v) \to \big(u=_{\trunc{n+1}A}v\big)\]
  defined for $u=\tproj {n+1}x$ and $v=\tproj {n+1}y$ and $p:x=y$ by
  \[\encode(\tproj n{p})\defeq \apfunc{\tproj{n+1}\nameless} (p).\]
  Since the codomain of $\encode$ is $n$-truncated, it suffices to define it only for $u$ and $v$ of this form, and then it's just the same
definition as before.
  We also define a function
  \[ r : \prd{u:\trunc{n+1} A} P(u,u) \]
  by induction on $u$, where $r(\tproj{n+1} x) \defeq \tproj n {\refl x}$.

  Now we can define an inverse map
  \[\decode: (u=_{\trunc{n+1}A}v) \to P(u,v)\]
  by
  \[\decode(p) \defeq \transfib{v\mapsto P(u,v)}{p}{r(u)}. \]
  To show that the composite
  \[ (u=_{\trunc{n+1}A}v) \xrightarrow{\decode} P(u,v) \xrightarrow{\encode} (u=_{\trunc{n+1}A}v) \]
  is the identity function, by path induction it suffices to check it for $\refl u : u=u$, in which case what we need to know is that
$\decode(r(u)) = \refl{u}$.
  But since this is an $n$-type, hence also an $(n+1)$-type, we may assume $u\jdeq \tproj {n+1} x$, in which case it follows by definition
of $r$ and $\decode$.
  Finally, to show that 
  \[ P(u,v) \xrightarrow{\encode} (u=_{\trunc{n+1}A}v) \xrightarrow{\decode} P(u,v) \]
  is the identity function, since this goal is again an $n$-type, we may assume that $u=\tproj {n+1}x$ and $v=\tproj {n+1}y$ and that we are
considering $\tproj n p:P(\tproj{n+1}x,\tproj{n+1}y)$ for some $p:x=y$.
  Then we have
  \begin{align*}
    \decode(\encode(\tproj n p)) &= \decode(\apfunc{\tproj{n+1}\nameless}(p))\\
    &= \transfib{v\mapsto P(\tproj{n+1}x,v)}{\apfunc{\tproj{n+1}\nameless}(p)}{\tproj n {\refl x}}\\
    &= \transfib{v\mapsto \trunc n{u=v}}{p}{\tproj n {\refl x}}\\
    &= \tproj n {\transfib{v \mapsto (u=v)}{p}{\refl x}}\\
    &= \tproj n p.
  \end{align*}
  This completes the proof that \encode and \decode are quasi-inverses.
  The stated result is then the special case where $u=\tproj {n+1}x$ and $v=\tproj {n+1}y$.
\end{proof}

\begin{cor}
  Let $n\ge-2$ and $(A,a)$ be a pointed type. Then
  \[\ttrunc n{\Omega(A,a)}=\Omega\mathopen{}\left(\trunc{n+1}{(A,a)}\right)\]
\end{cor}
\begin{proof}
  This is a special case of the previous lemma where $x=y=a$.
\end{proof}

\begin{cor}
  Let $n\ge -2$ and $k\ge 0$ and $(A,a)$ a pointed type.  Then
  \[\ttrunc n{\Omega^k(A,a)} = \Omega^k\mathopen{}\left(\trunc{n+k}{(A,a)}\right). \]
\end{cor}
\begin{proof}
  By induction on $k$, using the recursive definition of $\Omega^k$.
\end{proof}

We also observe that ``truncations are cumulative'': if we truncate to an $n$-type and then to a $k$-type with $k\le n$, then we might as
well have truncated directly to a $k$-type.

\begin{lem} \label{lem:truncation-le}
  Let $k,n\ge-2$ with $k\le{}n$ and $A:\type$. Then
  $\trunc k{\trunc nA}=\trunc kA$.
\end{lem}
\begin{proof}
  We define two maps $f:\trunc k{\trunc nA}\to\trunc kA$ and
  $g:\trunc kA\to\trunc k{\trunc nA}$ by
  \[
   f(\tproj k{\tproj na}) \defeq \tproj ka
   \qquad\text{and}\qquad
   g(\tproj ka) \defeq \tproj k{\tproj na}.
  \]
  The map $f$ is well-defined because $\trunc kA$ is $k$-truncated and also
  $n$-truncated (because $k\le{}n$), and the map $g$ is well-defined because
  $\trunc k{\trunc nA}$ is $k$-truncated.

  The composition $f\circ{}g:\trunc kA\to\trunc kA$ satisfies
  $(f\circ{}g)(\tproj ka)=\tproj ka$, hence $f\circ{}g=\idfunc[\trunc kA]$.
  Similarly, we have $(g\circ{}f)(\tproj k{\tproj na})=\tproj k{\tproj na}$ and hence $g\circ{}f=\idfunc[\trunc k{\trunc nA}]$.
\end{proof}


\index{truncation!n-truncation@$n$-truncation|)}%

\section{Colimits of \texorpdfstring{$n$}{n}-types}
\label{sec:pushouts}

Recall that in \autoref{sec:colimits}, we used higher inductive types to define pushouts of types, and proved their universal property.
In general, a (homotopy) colimit of $n$-types may no longer be an $n$-type (for an extreme counterexample, see \autoref{ex:s2-colim-unit}).
However, if we $n$-truncate it, we obtain an $n$-type which satisfies the correct universal property with respect to other $n$-types.

In this section we prove this for pushouts, which are the most important and nontrivial case of colimits.
Recall the following definitions from \autoref{sec:colimits}.

\begin{defn}
  A \define{span} 
  \indexdef{span} %
  is a 5-tuple $\Ddiag=(A,B,C,f,g)$ with 
  $f:C\to{}A$ and $g:C\to{}B$.
  \[\Ddiag=\quad\vcenter{\xymatrix{C \ar^g[r] \ar_f[d] & B \\ A & }}\]
\end{defn}

\begin{defn}
  Given a span $\Ddiag=(A,B,C,f,g)$ and a type $D$, a 
  \define{cocone under $\Ddiag$ with base $D$} is a triple $(i, j, h)$
  \index{cocone} %
  with $i:A\to{}D$, $j:B\to{}D$ and $h : \prd{c:C}i(f(c))=j(g(c))$:
  \[\uppercurveobject{{ }}\lowercurveobject{{ }}\twocellhead{{ }}
  \xymatrix{C \ar^g[r] \ar_f[d] \drtwocell{^h} & B \ar^j[d] \\ A \ar_i[r] & D
  }\]
  We denote by $\cocone{\Ddiag}{D}$ the type of all such cocones.
\end{defn}

The type of cocones is (covariantly) functorial.
For instance, given $D,E$ 
and a map $t:D\to{}E$, there is a map
  \[\function{\cocone{\Ddiag}{D}}{\cocone{\Ddiag}{E}}{c}{\composecocone{t}c}\]
  defined by:
  \[\composecocone{t}(i,j,h)=(t\circ{}i,t\circ{}j,\mapfunc{t}\circ{}h)\]
And given $D,E,F$, 
functions $t:D\to{}E$, $u:E\to{}F$ and $c:\cocone{\Ddiag}{D}$, we have
\begin{align}
  \composecocone{\idfunc[D]}c &= c \label{eq:composeconeid}\\
  \composecocone{(u\circ{}t)}c&=\composecocone{u}(\composecocone{t}c). \label{eq:composeconefunc}
\end{align}

\begin{defn}
  Given a span $\Ddiag$ of $n$-types, an $n$-type $D$, and a cocone
  $c:\cocone{\Ddiag}{D}$, the pair $(D,c)$ is said to be a \define{pushout
  of $\Ddiag$ in $n$-types}
  \indexdef{pushout!in ntypes@in $n$-types}%
  if for every $n$-type $E$, the map
  \[\function{(D\to{}E)}{\cocone{\Ddiag}{E}}{t}{\composecocone{t}c}\]
  is an equivalence.
\end{defn}

In order to construct pushouts of $n$-types, we need to explain how to reflect spans and cocones.

\bgroup
\def\reflect(#1){\trunc n{#1}}

\begin{defn}
  Let
  \[\Ddiag=\quad\vcenter{\xymatrix{C \ar^g[r] \ar_f[d] & B \\ A & }}\]
  be a span. We denote by $\reflect(\Ddiag)$ the following
  span of $n$-types:
  \[\reflect(\Ddiag)\defeq\quad \vcenter{\xymatrix{\reflect(C) \ar^{\reflect(g)}[r]
      \ar_{\reflect(f)}[d] & \reflect(B) \\ \reflect(A) & }}\]
\end{defn}

\begin{defn}
  Let $D:\type$ and $c=(i,j,h):\cocone{\Ddiag}{D}$.
  We define
  \[\reflect(c)=(\reflect(i),\reflect(j),\reflect(h)):
  \cocone{\reflect(\Ddiag)}{\reflect(D)}\]
  where $\reflect(h): \reflect(i) \circ \reflect(f) \htpy \reflect(j) \circ \reflect(g)$ is defined as in \autoref{thm:trunc-htpy}.
\end{defn}

\egroup

We now observe that the maps from each type to its $n$-truncation assemble into a map of spans, in the following sense.

\begin{defn}
  Let 
  \[\Ddiag=\quad\vcenter{\xymatrix{C \ar^g[r] \ar_f[d] & B \\ A & }}
  \qquad\text{and}\qquad
  \Ddiag'=\quad\vcenter{\xymatrix{C' \ar^{g'}[r] \ar_{f'}[d] & B' \\ A' & }}
  \]
  be spans.
  A \define{map of spans}
  \indexdef{map!of spans}%
  $\Ddiag \to \Ddiag'$ consists of functions $\alpha:A\to A'$, $\beta:B\to B'$, and $\gamma:C\to C'$ and homotopies $\phi: \alpha\circ f \htpy f'\circ \gamma$ and $\psi:\beta\circ g \htpy g' \circ \gamma$.
\end{defn}

Thus, for any span $\Ddiag$, we have a map of spans $\tprojf[\Ddiag] n : \Ddiag \to \trunc n\Ddiag$ consisting of $\tprojf[A]n$, $\tprojf[B]n$, $\tprojf[C]n$, and the naturality homotopies $\mathsf{nat}^f_n$ and $\mathsf{nat}^g_n$ from~\eqref{eq:trunc-nat}.

We also need to know that maps of spans behave functorially.
Namely, if $(\alpha,\beta,\gamma,\phi,\psi):\Ddiag \to \Ddiag'$ is a map of spans and $D$ any type, then we have
\[ \function{\cocone{\Ddiag'}{D}}{\cocone{\Ddiag}{D}}{(i,j,h)}{(i\circ \alpha,j\circ\beta, k)} \]
where $k: \prd{z:C} i(\alpha(f(z))) = j(\beta(g(z)))$ is the composite
\begin{equation}\label{eq:mapofspans-htpy}
\xymatrix{
  i(\alpha(f(z))) \ar@{=}[r]^{\apfunc{i}(\phi)} &
  i(f'(\gamma(z))) \ar@{=}[r]^{h(\gamma(z))} &
  j(g'(\gamma(z))) \ar@{=}[r]^{\apfunc{j}(\psi)} &
  j(\beta(g(z))). }
\end{equation}
We denote this cocone by $(i,j,h) \circ (\alpha,\beta,\gamma,\phi,\psi)$.
Moreover, this functorial action commutes with the other functoriality of cocones:

\begin{lem}\label{thm:conemap-funct}
  Given $(\alpha,\beta,\gamma,\phi,\psi):\Ddiag \to \Ddiag'$ and $t:D\to E$, the following diagram commutes:
  \begin{equation*}
    \vcenter{\xymatrix{
        \cocone{\Ddiag'}{D}\ar[r]^{t \circ {\blank}}\ar[d] &
        \cocone{\Ddiag'}{E}\ar[d]\\
        \cocone{\Ddiag}{D}\ar[r]_{t \circ {\blank}} &
        \cocone{\Ddiag}{E}
      }}
  \end{equation*}
\end{lem}
\begin{proof}
  Given $(i,j,h):\cocone{\Ddiag'}{D}$, note that both composites yield a cocone whose first two components are $t\circ i\circ \alpha$ and $t\circ j\circ\beta$.
  Thus, it remains to verify that the homotopies agree.
  For the top-right composite, the homotopy is~\eqref{eq:mapofspans-htpy} with $(i,j,h)$ replaced by $(t\circ i, t\circ j, \apfunc{t}\circ h)$:
  \begin{equation*}
    \xymatrix@+2.8em{
      {t \, i \, \alpha \, f \, z} \ar@{=}[r]^{\apfunc{t\circ i}(\phi)} &
      {t \, i \, f' \, \gamma \, z} \ar@{=}[r]^{\apfunc{t}(h(\gamma(z)))} &
      {t \, j \, g' \, \gamma \, z} \ar@{=}[r]^{\apfunc{t\circ j}(\psi)} &
      {t \, j \, \beta \, g \, z}
    }
  \end{equation*}
  (For brevity, we are omitting the parentheses around the arguments of functions.)
  On the other hand, for the left-bottom composite, the homotopy is $\apfunc{t}$ applied to~\eqref{eq:mapofspans-htpy}.
  Since $\apfunc{}$ respects path-concatenation, this is equal to
  \begin{equation*}
    \xymatrix@+2.8em{
      {t \, i \, \alpha \, f \, z} \ar@{=}[r]^{\apfunc{t}(\apfunc{i}(\phi))} &
      {t \, i \, f' \, \gamma \, z} \ar@{=}[r]^{\apfunc{t}(h(\gamma(z)))} &
      {t \, j \, g' \, \gamma \, z} \ar@{=}[r]^{\apfunc{t}(\apfunc{j}(\psi))} &
      {t \, j \, \beta \, g \, z}. }
  \end{equation*}
  But $\apfunc{t}\circ \apfunc{i} = \apfunc{t\circ i}$ and similarly for $j$, so these two homotopies are equal.
\end{proof}

Finally, note that since we defined $\trunc nc : \cocone{\trunc n \Ddiag}{\trunc n D}$ using \autoref{thm:trunc-htpy}, the additional condition~\eqref{eq:trunc-htpy} implies
\begin{equation}
  \tprojf[D] n \circ c = \trunc n c \circ \tprojf[\Ddiag]n. \label{eq:conetrunc}
\end{equation}
for any $c:\cocone{\Ddiag}{D}$.
Now we can prove our desired theorem.

\begin{thm}
  \label{reflectcommutespushout}
  \index{universal!property!of pushout}%
  Let $\Ddiag$ be a span and $(D,c)$ its pushout.
  Then $(\trunc nD,\trunc n c)$ is a pushout of $\trunc n\Ddiag$ in $n$-types.
\end{thm}
\begin{proof}
  Let $E$ be an $n$-type, and consider the following diagram:
\bgroup
\def\reflect(#1){\trunc n{#1}}
  \begin{equation*}
  \vcenter{\xymatrix{
      (\trunc nD \to E)\ar[r]^-{\blank\circ \tprojf[D] n}\ar[d]_{\blank\circ \trunc nc} &
      (D\to E)\ar[d]^{\blank\circ c}\\
      \cocone{\trunc n \Ddiag}{E}\ar[r]^-{\blank\circ \tprojf[\Ddiag]n}\ar@{<-}[d]_{\ell_1} &
      \cocone{\Ddiag}{E}\ar@{<-}[d]^{\ell_2}\\
      (\reflect(A)\to{}E)\times_{(\reflect(C)\to{}E)}(\reflect(B)\to{}E)\ar[r] &
      (A\to{}E)\times_{(C\to{}E)}(B\to{}E)
      }}
  \end{equation*}
\egroup
  The upper horizontal arrow is an equivalence since $E$ is an $n$-type, while $\blank\circ c$ is an equivalence since $c$ is a pushout cocone.
  Thus, by the 2-out-of-3 property, to show that $\blank\circ \trunc nc$ is an equivalence, it will suffice to show that the upper square commutes and that the middle horizontal arrow is an equivalence.
  To see that the upper square commutes, let $t:\trunc nD \to E$; then
  \begin{align}
    \big(t \circ \trunc n c\big) \circ \tprojf[\Ddiag] n
    &= t \circ \big(\trunc n c \circ \tprojf[\Ddiag] n\big)
    \tag{by \autoref{thm:conemap-funct}}\\
    &= t\circ \big(\tprojf[D]n \circ c\big)
    \tag{by~\eqref{eq:conetrunc}}\\
    &= \big(t\circ \tprojf[D]n\big) \circ c
    \tag{by~\eqref{eq:composeconefunc}}.
  \end{align}
  To show that the middle horizontal arrow is an equivalence, consider the lower square.
  The two lower vertical arrows are simply applications of $\happly$:
  \begin{align*}
    \ell_1(i,j,p) &\defeq (i,j,\happly(p))\\
    \ell_2(i,j,p) &\defeq (i,j,\happly(p))
  \end{align*}
  and hence are equivalences by function extensionality.
  The lowest horizontal arrow is defined by
  \[ (i,j,p) \mapsto \big( i\circ \tprojf[A]n,\;\; j \circ \tprojf[B] n,\;\; q\big) \]
  where $q$ is the composite
  \begin{align}
    i\circ \tprojf[A]n \circ f
    &= i\circ \trunc nf \circ \tprojf[C]n
    \tag{by $\funext(\lam{z} \apfunc{i}(\mathsf{nat}^f_n(z)))$}\\
    &= j\circ \trunc ng \circ \tprojf[C]n
    \tag{by $\apfunc{\blank\circ \tprojf[C] n}(p)$}\\
    &= j\circ \tprojf[B]n \circ g.
    \tag{by $\funext(\lam{z} \apfunc{j}(\mathsf{nat}^g_n(z)))$}
  \end{align}
  This is an equivalence, because it is induced by an equivalence of cospans.
  Thus, by 2-out-of-3, it will suffice to show that the lower square commutes.
  But the two composites around the lower square agree definitionally on the first two components, so it suffices to show that for $(i,j,p)$ in the lower left corner and $z:C$, the path
  \[ \happly(q,z) : i(\tproj n{f(z)}) = j(\tproj n{g(z)}) \]
  (with $q$ as above)
  is equal to the composite
  \begin{align}
    i(\tproj n{f(z)})
    &= i(\trunc nf(\tproj nz))
    \tag{by $\apfunc{i}(\mathsf{nat}^f_n(z))$}\\
    &= j(\trunc ng(\tproj nz))
    \tag{by $\happly(p,\tproj nz)$}\\
    &= j(\tproj n{g(z)}).
    \tag{by $\apfunc{j}(\mathsf{nat}^g_n(z))$}
  \end{align}
  However, since $\happly$ is functorial, it suffices to check equality for the three component paths:
  \begin{align*}
    \happly({\funext(\lam{z} \apfunc{i}(\mathsf{nat}^f_n(z)))},z)
    &= {\apfunc{i}(\mathsf{nat}^f_n(z))}\\
    \happly(\apfunc{\blank\circ \tprojf[C] n}(p), z)
    &= {\happly(p,\tproj nz)}\\
    \happly({\funext(\lam{z} \apfunc{j}(\mathsf{nat}^g_n(z)))},z)
    &= {\apfunc{j}(\mathsf{nat}^g_n(z))}.
  \end{align*}
  The first and third of these are just the fact that $\happly$ is quasi-inverse to $\funext$, while
  the second is an easy general lemma about $\happly$ and precomposition.
\end{proof}

\section{Connectedness}
\label{sec:connectivity}

An $n$-type is one that has no interesting information above dimension $n$.
By contrast, an \emph{$n$-connected type} is one that has no interesting information \emph{below} dimension $n$.
It turns out to be natural to study a more general notion for functions as well.

\begin{defn}
A function $f:A\to B$ is said to be \define{$n$-connected}
\indexdef{function!n-connected@$n$-connected}%
\indexsee{n-connected@$n$-connected!function}{function, $n$-connected}%
if for all $b:B$, the type $\trunc n{\hfiber f b}$ is contractible:
\begin{equation*}
  \mathsf{conn}_n(f)\defeq \prd{b:B}\iscontr(\trunc n{\hfiber{f}b}). 
\end{equation*}
A type $A$ is said to be \define{$n$-connected}
\indexsee{n-connected@$n$-connected!type}{type, $n$-connected}%
\indexdef{type!n-connected@$n$-connected}%
 if the unique function $A\to\unit$ is $n$-connected, i.e.\ if $\trunc nA$ is contractible.
\end{defn}
\indexsee{connected!function}{function, $n$-connected}

Thus, a function $f:A\to B$ is $n$-connected if and only if $\hfib{f}b$ is $n$-connected for every $b:B$.
Of course, every function is $(-2)$-connected.
At the next level, we have:

\begin{lem}\label{thm:minusoneconn-surjective}
  \index{function!surjective}%
  A function $f$ is $(-1)$-connected if and only if it is surjective in the sense of \autoref{sec:mono-surj}.
\end{lem}
\begin{proof}
  We defined $f$ to be surjective if $\trunc{-1}{\hfiber f b}$ is inhabited for all $b$.
  But since it is a mere proposition, inhabitation is equivalent to contractibility.
\end{proof}

Thus, $n$-connectedness of a function for $n\ge 0$ can be thought of as a strong form of surjectivity.
Category-theoretically, $(-1)$-connectedness corresponds to essential surjectivity on objects, while $n$-connectedness corresponds to essential surjectivity on $k$-morphisms for $k\le n+1$.

\autoref{thm:minusoneconn-surjective} also implies that a type $A$ is $(-1)$-connected if and only if it is merely inhabited.
When a type is $0$-connected we may simply say that it is \define{connected},
\indexdef{connected!type}%
\indexdef{type!connected}%
and when it is $1$-connected we say it is \define{simply connected}.
\indexdef{simply connected type}%
\indexdef{type!simply connected}%

\begin{rmk}\label{rmk:connectedness-indexing}
  While our notion of $n$-connectedness for types agrees with the standard notion in homotopy theory, our notion of $n$-connectedness for \emph{functions} is off by one from a common indexing in classical homotopy theory.
  Whereas we say a function $f$ is $n$-connected if all its fibers are $n$-connected, some classical homotopy theorists would call such a function $(n+1)$-connected.
  (This is due to a historical focus on \emph{cofibers} rather than fibers.)
\end{rmk}

We now observe a few closure properties of connected maps.
\index{function!n-connected@$n$-connected}

\begin{lem}
\index{retract!of a function}%
Suppose that $g$ is a retract of a $n$-connected function $f$.  Then $g$ is
$n$-connected.
\end{lem}
\begin{proof}
This is a direct consequence of \autoref{lem:func_retract_to_fiber_retract}.
\end{proof}

\begin{cor}
If $g$ is homotopic to a $n$-connected function $f$, then $g$ is $n$-connected.
\end{cor}

\begin{lem}\label{lem:nconnected_postcomp}
Suppose that $f:A\to B$ is $n$-connected. Then $g:B\to C$ is $n$-connected if and only if $g\circ f$ is
$n$-connected.
\end{lem}

\begin{proof}
For any $c:C$, we have
\begin{align*}
  \trunc n{\hfib{g\circ f}c}
  & \eqvsym \Trunc n{ \sm{w:\hfib{g}c}\hfib{f}{\proj1 w}}
  \tag{by \autoref{ex:unstable-octahedron}}\\
  & \eqvsym \Trunc n{\sm{w:\hfib{g}c} \trunc n{\hfib{f}{\proj1 w}}}
  \tag{by \autoref{thm:trunc-in-truncated-sigma}}\\
  & \eqvsym \trunc n{\hfib{g}c}.
  \tag{since $\trunc n{\hfib{f}{\proj1 w}}$ is contractible}
\end{align*}
It follows that $\trunc n{\hfib{g}c}$ is contractible if and only if $\trunc n{\hfib{g\circ f}c}$ is
contractible.
\end{proof}

Importantly, $n$-connected functions can be equivalently characterized as those which satisfy an ``induction principle'' with respect to $n$-types.\index{induction principle!for connected maps} 
This idea will lead directly into our proof of the Freudenthal suspension theorem in \autoref{sec:freudenthal}.

\begin{lem}\label{prop:nconnected_tested_by_lv_n_dependent types}
For $f:A\to B$ and $P:B\to\type$, consider the following function:
\begin{equation*}
\lam{s} s\circ f :\Parens{\prd{b:B} P(b)}\to\Parens{\prd{a:A}P(f(a))}.
\end{equation*}
For a fixed $f$ and $n\ge -2$, the following are equivalent.
\begin{enumerate}
\item $f$ is $n$-connected.\label{item:conntest1}
\item For every $P:B\to\ntype{n}$, the map $\lam{s} s\circ f$ is an equivalence.\label{item:conntest2}
\item For every $P:B\to\ntype{n}$, the map $\lam{s} s\circ f$ has a section.\label{item:conntest3}
\end{enumerate}
\end{lem}

\begin{proof}
Suppose that $f$ is $n$-connected and let $P:B\to\ntype{n}$. Then we have the equivalences
\begin{align}
  \prd{b:B} P(b) & \eqvsym \prd{b:B} \Parens{\trunc n{\hfib{f}b} \to P(b)}
  \tag{since $\trunc n{\hfib{f}b}$ is contractible}\\
  & \eqvsym \prd{b:B} \Parens{\hfib{f}b\to P(b)}
  \tag{since $P(b)$ is an $n$-type}\\
  & \eqvsym \prd{b:B}{a:A}{p:f(a)= b} P(b)
  \tag{by the left universal property of $\Sigma$-types}\\
  & \eqvsym \prd{a:A} P(f(a)).
  \tag{by the left universal property of path types}
\end{align}
We omit the proof that this equivalence is indeed given by $\lam{s} s\circ f$.
Thus,~\ref{item:conntest1}$\Rightarrow$\ref{item:conntest2}, and clearly~\ref{item:conntest2}$\Rightarrow$\ref{item:conntest3}.
To show~\ref{item:conntest3}$\Rightarrow$\ref{item:conntest1}, consider the type family
\begin{equation*}
P(b)\defeq \trunc n{\hfib{f}b}.
\end{equation*}
Then~\ref{item:conntest3} yields a map $c:\prd{b:B} \trunc n{\hfib{f}b}$ with
$c(f(a))=\tproj n{\pairr{a,\refl{f(a)}}}$. To show that each $\trunc n{\hfib{f}b}$ is contractible,
we will find a function of type
\begin{equation*}
\prd{b:B}{w:\trunc n{\hfib{f}b}} w= c(b).
\end{equation*}
By \autoref{thm:truncn-ind}, for this it suffices to find a function of type
\begin{equation*}
\prd{b:B}{a:A}{p:f(a)= b} \tproj n{\pairr{a,p}}= c(b).
\end{equation*}
But by rearranging variables and path induction, this is equivalent to the type
\begin{equation*}
\prd{a:A} \tproj n{\pairr{a,\refl{f(a)}}}= c(f(a)).
\end{equation*}
This property holds by our choice of $c(f(a))$. 
\end{proof}

\begin{cor}\label{cor:totrunc-is-connected}
For any $A$, the canonical function $\tprojf n:A\to\trunc n A$ is $n$-connected.
\end{cor}
\begin{proof}
By \autoref{thm:truncn-ind} and the associated uniqueness principle, the condition of \autoref{prop:nconnected_tested_by_lv_n_dependent types} holds.
\end{proof}

For instance, when $n=-1$, \autoref{cor:totrunc-is-connected} says that the map $A\to \brck A$ from a type to its propositional truncation is surjective.

\begin{cor}\label{thm:nconn-to-ntype-const}\label{connectedtotruncated}
A type $A$ is $n$-connected if and only if the map
\begin{equation*}
  \lam{b}{a} b: B \to (A\to B)
\end{equation*}
is an equivalence for every $n$-type $B$.
In other words, ``every map from $A$ to an $n$-type is constant''.
\end{cor}
\begin{proof}
  By \autoref{prop:nconnected_tested_by_lv_n_dependent types} applied to a function with codomain $\unit$.
\end{proof}

\begin{lem}\label{lem:nconnected_to_leveln_to_equiv}
Let $B$ be an $n$-type and let $f:A\to B$ be a function. Then the induced function $g:\trunc n A\to B$ is an
equivalence if and only if $f$ is $n$-connected.
\end{lem}

\begin{proof}
By \autoref{cor:totrunc-is-connected}, $\tprojf n$ is $n$-connected.
Thus, since $f = g\circ \tprojf n$, by
\autoref{lem:nconnected_postcomp} $f$ is $n$-connected if and only if $g$ is $n$-connected.
But since $g$ is a function between $n$-types, its fibers are also $n$-types.
Thus, $g$ is $n$-connected if and only if it is an equivalence.
\end{proof}

We can also characterize connected pointed types in terms of connectivity of the inclusion of their basepoint.

\begin{lem}\label{thm:connected-pointed}
  \index{basepoint}%
  Let $A$ be a type and $a_0:\unit\to A$ a basepoint, with $n\ge -1$.
  Then $A$ is $n$-connected if and only if the map $a_0$ is $(n-1)$-connected.
\end{lem}
\begin{proof}
  First suppose $a_0:\unit\to A$ is $(n-1)$-connected and let $B$ be an $n$-type; we will use \autoref{thm:nconn-to-ntype-const}.
  The map $\lam{b}{a} b: B \to (A\to B)$ has a retraction given by $f\mapsto f(a_0)$, so it suffices to show it also has a section, i.e.\ that for any $f:A\to B$ there is $b:B$ such that $f = \lam{a}b$.
  We choose $b\defeq f(a_0)$.
  Define $P:A\to\type$ by $P(a) \defeq (f(a)=f(a_0))$.
  Then $P$ is a family of $(n-1)$-types and we have $P(a_0)$; hence we have $\prd{a:A} P(a)$ since $a_0:\unit\to A$ is $(n-1)$-connected.
  Thus, $f = \lam{a} f(a_0)$ as desired.

  Now suppose $A$ is $n$-connected, and let $P:A\to\ntype{(n-1)}$ and $u:P(a_0)$ be given.
  By \autoref{prop:nconnected_tested_by_lv_n_dependent types}, it will suffice to construct $f:\prd{a:A} P(a)$ such that $f(a_0)=u$.
  Now $\ntype{(n-1)}$ is an $n$-type and $A$ is $n$-connected, so by \autoref{thm:nconn-to-ntype-const}, there is an $n$-type $B$ such that $P = \lam{a} B$.
  Hence, we have a family of equivalences $g:\prd{a:A} (\eqv{P(a)}{B})$.
  Define $f(a) \defeq \opp{g_a}(g_{a_0}(u))$; then $f:\prd{a:A} P(a)$ and $f(a_0) = u$ as desired.
\end{proof}

In particular, a pointed type $(A,a_0)$ is 0-connected if and only if $a_0:\unit\to A$ is surjective, which is to say $\prd{x:A} \brck{x=a_0}$.

A useful variation on \autoref{lem:nconnected_postcomp} is:

\begin{lem}\label{lem:nconnected_postcomp_variation}
Let $f:A\to B$ be a function and $P:A\to\type$ and $Q:B\to\type$ be type families. Suppose that $g:\prd{a:A} P(a)\to Q(f(a))$
is a fiberwise $n$-connected%
\index{fiberwise!n-connected family of functions@$n$-connected family of functions}
family of functions, i.e.\ each function $g_a : P(a) \to Q(f(a))$ is $n$-connected. Then the function
\begin{align*}
\varphi &:\Parens{\sm{a:A} P(a)}\to\Parens{\sm{b:B} Q(b)}\\
\varphi(a,u) &\defeq \pairr{f(a),g_a(u)}
\end{align*}
is $n$-connected if and only if $f$ is $n$-connected.
\end{lem}

\begin{proof}
For $b:B$ and $v:Q(b)$ we have
{\allowdisplaybreaks
\begin{align*}
\trunc n{\hfib{\varphi}{\pairr{b,v}}} & \eqvsym \Trunc n{\sm{a:A}{u:P(a)}{p:f(a)= b} \trans{\ap f p }{g_a(u)}= v}\\
& \eqvsym \Trunc n{\sm{w:\hfib{f}b}{u:P(\proj1(w))} g_{\proj 1 w}(u)= \trans{\opp{\ap f {\proj2 w}}}{v}}\\
& \eqvsym \Trunc n{\sm{w:\hfib{f}b} \hfib{g(\proj1 w)}{\trans{\opp{\ap f {\proj 2 w}}}{v}}}\\
& \eqvsym \Trunc n{\sm{w:\hfib{f}b} \trunc n{\hfib{g(\proj1 w)}{\trans{\opp{\ap f {\proj 2 w}}}{v}}}}\\
& \eqvsym \trunc n{\hfib{f}b}
\end{align*}}
where the transportations along $f(p)$ and $f(p)^{-1}$ are with respect to $Q$.
Therefore, if either is contractible, so is the other.
\end{proof}

In the other direction, we have

\begin{lem}\label{prop:nconn_fiber_to_total}
Let $P,Q:A\to\type$ be type families and consider a fiberwise transformation\index{fiberwise!transformation}
\begin{equation*}
f:\prd{a:A} \Parens{P(a)\to Q(a)}
\end{equation*}
from $P$ to $Q$. Then the induced map $\total f: \sm{a:A}P(a) \to \sm{a:A} Q(a)$ is $n$-connected if and only if each $f(a)$ is $n$-connected. 
\end{lem}

\begin{proof}
By \autoref{fibwise-fiber-total-fiber-equiv}, we have
$\hfib{\total f}{\pairr{x,v}}\eqvsym\hfib{f(x)}v$
for each $x:A$ and $v:Q(x)$. Hence $\trunc n{\hfib{\total f}{\pairr{x,v}}}$ is contractible if and only if
$\trunc n{\hfib{f(x)}v}$ is contractible.
\end{proof}

Another useful fact about connected maps is that they induce an
equivalence on $n$-truncations:

\begin{lem} \label{lem:connected-map-equiv-truncation}
If $f : A \to B$ is $n$-connected, then it induces an equivalence
$\eqv{\trunc{n}{A}}{\trunc{n}{B}}$.
\end{lem}
\begin{proof}
Let $c$ be the proof that $f$ is $n$-connected.  From left to right, we
use the map $\trunc{n}{f} : \trunc{n}{A} \to \trunc{n}{B}$.
To define the map from right to left, by the universal property of
truncations, it suffices to give a map $\mathsf{back} : B \to {\trunc{n}{A}}$.  We can
define this map as follows:
\[
\mathsf{back}(y) \defeq \trunc{n}{\proj{1}}{(\proj{1}{(c(y))})}
\]
By definition, $c(y)$ has type $\iscontr(\trunc n {\hfiber{f}y})$, so its
first component has type $\trunc n{\hfiber{f}y}$, and we can obtain an
element of $\trunc n A$ from this by projection.

Next, we show that the composites are the identity.  In both directions,
because the goal is a path in an $n$-truncated type, it suffices to
cover the case of the constructor $\tprojf{n}$.

In one direction, we must show that for all $x:A$, 
\[
\trunc{n}{\proj{1}}{(\proj{1}{(c(f(x)))})} = \tproj{n}{x}
\]
But $\tproj{n}{(x, \refl{})} : \trunc n{\hfiber{f}y}$, and
$c(y)$ says that this type is contractible, so 
\[
\proj{1}{(c(f(x)))} = \tproj{n}{(x, \refl{})}
\]
Applying $\trunc{n}{\proj{1}}$ to both sides of this equation gives the
result.  

In the other direction, we must show that for all $y:B$, 
\[
\trunc{n}{f}(\trunc{n}{\proj{1}} (\proj{1}{(c(y))})) = \tproj{n}{y}
\]
$\proj{1}{(c(y))}$ has type $\trunc n {\hfiber{f}y}$, and the path we
want is essentially the second component of the $\hfiber{f}y$, but we
need to make sure the truncations work out.  

In general, suppose we are given $p:\trunc{n}{\sm{x:A} B(x)}$ and wish to prove
$P(\trunc{n}{\proj{1}{}}(p))$. By truncation induction, it suffices to
prove $P(\tproj{n}{a})$ for all $a:A$ and $b:B(a)$.  Applying this
principle in this case, it suffices to prove
\[
\trunc{n}{f}(\tproj{n}{a}) = \tproj{n}{y}
\]
given $a:A$ and $b:f (a) = y$.  But the left-hand side equals $\tproj{n}{f (a)}$,
so applying $\tprojf{n}$ to both sides of $b$ gives the result.
\end{proof}

One might guess that this fact characterizes the $n$-connected maps, but in fact being $n$-connected is a bit stronger than this.
For instance, the inclusion $\bfalse:\unit \to\bool$ induces an equivalence on $(-1)$-truncations, but is not surjective (i.e.\ $(-1)$-connected).
In \autoref{sec:long-exact-sequence-homotopy-groups} we will see that the difference in general is an analogous extra bit of surjectivity.

\section{Orthogonal factorization}
\label{sec:image-factorization}

\index{unique!factorization system|(}%
\index{orthogonal factorization system|(}%
In set theory, the surjections and the injections form a unique factorization system: every function factors essentially uniquely as a surjection followed by an injection.
We have seen that surjections generalize naturally to $n$-connected maps, so it is natural to inquire whether these also participate in a factorization system.
Here is the corresponding generalization of injections.

\begin{defn}
  A function $f:A\to B$ is \define{$n$-truncated}
  \indexdef{n-truncated@$n$-truncated!function}%
  \indexdef{function!n-truncated@$n$-truncated}%
if the fiber $\hfib f b$ is an $n$-type for all $b:B$.
\end{defn}

In particular, $f$ is $(-2)$-truncated if and only if it is an equivalence.
And of course, $A$ is an $n$-type if and only if $A\to\unit$ is $n$-truncated.
Moreover, $n$-truncated maps could equivalently be defined recursively, like $n$-types.

\begin{lem}\label{thm:modal-mono}
  For any $n\ge -2$, a function $f:A\to B$ is $(n+1)$-truncated if and only if for all $x,y:A$, the map $\apfunc{f}:(x=y) \to (f(x)=f(y))$ is $n$-truncated.
  \index{function!embedding}%
  \index{function!injective}%
  In particular, $f$ is $(-1)$-truncated if and only if it is an embedding in the sense of \autoref{sec:mono-surj}.
\end{lem}
\begin{proof}
  Note that for any $(x,p),(y,q):\hfib f b$, we have
  \begin{align*}
    \big((x,p) = (y,q)\big)
    &= \sm{r:x=y} (p = \apfunc f(r)\ct q)\\
    &= \sm{r:x=y} (\apfunc f (r) = p\ct \opp q)\\
    &= \hfib{\apfunc{f}}{p\ct \opp q}.
  \end{align*}
  Thus, any path space in any fiber of $f$ is a fiber of $\apfunc{f}$.
  On the other hand, choosing $b\defeq f(y)$ and $q\defeq \refl{f(y)}$ we see that any fiber of $\apfunc f$ is a path space in a fiber of $f$.
  The result follows, since $f$ is $\nplusone$-truncated if all path spaces of its fibers are $n$-types.
\end{proof}

We can now construct the factorization, in a fairly obvious way.

\begin{defn}\label{defn:modal-image}
Let $f:A\to B$ be a function. The \define{$n$-image}
\indexdef{image}%
\indexdef{image!n-image@$n$-image}%
\indexdef{n-image@$n$-image}%
\indexdef{function!n-image of@$n$-image of}%
of $f$ is defined as
\begin{equation*}
\im_n(f)\defeq \sm{b:B} \trunc n{\hfib{f}b}.
\end{equation*}
When $n=-1$, we write simply $\im(f)$ and call it the \define{image} of $f$.
\end{defn}

\begin{lem}\label{prop:to_image_is_connected}
For any function $f:A\to B$, the canonical function $\tilde{f}:A\to\im_n(f)$ is $n$-connected. 
Consequently, any function factors as an $n$-connected function followed by an $n$-truncated function.
\end{lem}

\begin{proof}
Note that $A\eqvsym\sm{b:B}\hfib{f}b$. The function $\tilde{f}$ is the function on total spaces induced by the canonical fiberwise transformation
\begin{equation*}
\prd{b:B} \Parens{\hfib{f}b\to\trunc n{\hfib{f}b}}.
\end{equation*}
Since each map $\hfib{f}b\to\trunc n{\hfib{f}b}$ is $n$-connected by \autoref{cor:totrunc-is-connected}, $\tilde{f}$ is $n$-connected by \autoref{prop:nconn_fiber_to_total}.
Finally, the projection $\proj1:\im_n(f) \to B$ is $n$-truncated, since its fibers are equivalent to the $n$-truncations of the fibers of $f$.
\end{proof}

In the following lemma we set up some machinery to prove the unique factorization theorem.

\begin{lem}\label{prop:factor_equiv_fiber}
Suppose we have a commutative diagram of functions
\begin{equation*}
  \xymatrix{
    {A} \ar[r]^{g_1} \ar[d]_{g_2} &
    {X_1} \ar[d]^{h_1} &
    \\
    {X_2} \ar[r]_{h_2}
    &
    {B}
  }
\end{equation*}
with $H:h_1\circ g_1\htpy h_2\circ g_2$, where $g_1$ and $g_2$ are $n$-connected and where $h_1$ and $h_2$ are $n$-truncated.
Then there is an equivalence
\begin{equation*}
E(H,b):\hfib{h_1}b\eqvsym\hfib{h_2}b
\end{equation*}
for any $b:B$, such that for any $a:A$ we have an identification
\[\overline{E}(H,a) :  E(H,h_1(g_1(a)))({g_1(a),\refl{h_1(g_1(a))}}) = \pairr{g_2(a),\opp{H(a)}}.\]
\end{lem}

\begin{proof}
Let $b:B$. Then we have the following equivalences:
\begin{align}
\hfib{h_1}b
& \eqvsym \sm{w:\hfib{h_1}b} \trunc n{ \hfib{g_1}{\proj1 w}}
\tag{since $g_1$ is $n$-connected}\\
& \eqvsym \Trunc n{\sm{w:\hfib{h_1}b}\hfib{g_1}{\proj1 w}}
\tag{by \autoref{thm:refl-over-ntype-base}, since $h_1$ is $n$-truncated}\\
& \eqvsym \trunc n{\hfib{h_1\circ g_1}b}
\tag{by \autoref{ex:unstable-octahedron}}
\end{align}
and likewise for $h_2$ and $g_2$.
Also, since we have a homotopy $H:h_1\circ g_1\htpy h_2\circ g_2$, there is an obvious equivalence $\hfib{h_1\circ g_1}b\eqvsym\hfib{h_2\circ g_2}b$.
Hence we obtain
\begin{equation*}
\hfib{h_1}b\eqvsym\hfib{h_2}b
\end{equation*}
for any $b:B$. By analyzing the underlying functions, we get the following representation of what happens to the element
$\pairr{g_1(a),\refl{h_1(g_1(a))}}$ after applying each of the equivalences of which $E$ is composed.
Some of the identifications are definitional, but others (marked with a $=$ below) are only propositional; putting them together we obtain $\overline E(H,a)$.
{\allowdisplaybreaks
\begin{align*}
\pairr{g_1(a),\refl{h_1(g_1(a))}} & 
    \overset{=}{\mapsto} \Pairr{\pairr{g_1(a),\refl{h_1(g_1(a))}}, \tproj n{ \pairr{a,\refl{g_1(a)}}}}\\
  & \mapsto \tproj n { \pairr{\pairr{g_1(a),\refl{h_1(g_1(a))}}, \pairr{a,\refl{g_1(a)}} }}\\
  & \mapsto \tproj n { \pairr{a,\refl{h_1(g_1(a))}}}\\
  & \overset{=}{\mapsto} \tproj n { \pairr{a,\opp{H(a)}}}\\
  & \mapsto \tproj n { \pairr{\pairr{g_2(a),\opp{H(a)}},\pairr{a,\refl{g_2(a)}}} }\\
  & \mapsto \Pairr{\pairr{g_2(a),\opp{H(a)}}, \tproj n {\pairr{a,\refl{g_2(a)}}} }\\
  & \mapsto \pairr{g_2(a),\opp{H(a)}}
\end{align*}}
The first equality is because for general $b$, the map
\narrowequation{ \hfib{h_1}b \to \sm{w:\hfib{h_1}b} \trunc n{ \hfib{g_1}{\proj1 w}} }
inserts the center of contraction for $\trunc n{ \hfib{g_1}{\proj1 w}}$ supplied by the assumption that $g_1$ is $n$-truncated; whereas in the case in question this type has the obvious inhabitant $\tproj n{ \pairr{a,\refl{g_1(a)}}}$, which by contractibility must be equal to the center.
The second propositional equality is because the equivalence $\hfib{h_1\circ g_1}b\eqvsym\hfib{h_2\circ g_2}b$ concatenates the second components with $\opp{H(a)}$, and we have $\opp{H(a)} \ct \refl{} = \opp{H(a)}$.
The reader may check that the other equalities are definitional (assuming a reasonable solution to \autoref{ex:unstable-octahedron}).
\end{proof}


Combining \autoref{prop:to_image_is_connected,prop:factor_equiv_fiber}, we have the following unique factorization result:

\begin{thm}\label{thm:orth-fact}
For each $f:A\to B$, the space $\fact_n(f)$ defined by
\begin{equation*}
\sm{X:\type}{g:A\to X}{h:X\to B} (h\circ g\htpy f)\times\mathsf{conn}_n(g)\times\mathsf{trunc}_n(h).
\end{equation*}
is contractible.
Its center of contraction is the element
\begin{equation*}
\pairr{\im_n(f),\tilde{f},\proj1,\theta,\varphi,\psi}:\fact_n(f)
\end{equation*}
arising from \autoref{prop:to_image_is_connected},
where $\theta:\proj1\circ\tilde{f}\htpy f$ is the canonical homotopy, where $\varphi$ is the proof of
\autoref{prop:to_image_is_connected}, and where $\psi$ is the obvious proof that $\proj1:\im_n(f)\to B$ has $n$-truncated fibers.
\end{thm}

\begin{proof}
By \autoref{prop:to_image_is_connected} we know that there is an element of $\fact_n(f)$, hence it is enough to
show that $\fact_n(f)$ is a mere proposition. Suppose we have two $n$-factorizations
\begin{equation*}
\pairr{X_1,g_1,h_1,H_1,\varphi_1,\psi_1}\qquad\text{and}\qquad\pairr{X_2,g_2,h_2,H_2,\varphi_2,\psi_2}
\end{equation*}
of $f$. Then we have the pointwise-concatenated homotopy
\[ H\defeq (\lam{a} H_1(a) \ct H_2^{-1}(a)) \,:\, (h_1\circ g_1\htpy h_2\circ g_2).\]
By univalence and the characterization of paths and transport in $\Sigma$-types, function types, and path types, it suffices to show that
\begin{enumerate}
\item there is an equivalence $e:X_1\eqvsym X_2$,
\item there is a homotopy $\zeta:e\circ g_1\htpy g_2$,
\item there is a homotopy $\eta:h_2\circ e\htpy h_1$,
\item for any $a:A$ we have $\opp{\apfunc{h_2}(\zeta(a))} \ct \eta(g_1(a)) \ct H_1(a) = H_2(a)$.
\end{enumerate}
We prove these four assertions in that order.
\begin{enumerate}
\item By \autoref{prop:factor_equiv_fiber}, we have a fiberwise equivalence
\begin{equation*}
E(H) : \prd{b:B} \eqv{\hfib{h_1}b}{\hfib{h_2}b}.
\end{equation*}
This induces an equivalence of total spaces, i.e.\ we have
\begin{equation*}
\eqvspaced{\Parens{\sm{b:B} \hfib{h_1}b}}{\Parens{\sm{b:B}\hfib{h_2}b}}.
\end{equation*}
Of course, we also have the equivalences $X_1\eqvsym\sm{b:B}\hfib{h_1}b$ and $X_2\eqvsym\sm{b:B}
\hfib{h_2}b$ from \autoref{thm:total-space-of-the-fibers}.
This gives us our equivalence $e:X_1\eqvsym X_2$; the reader may verify that the underlying function of $e$ is given by
\begin{equation*}
e(x) \jdeq \proj1(E(H,h_1(x))(x,\refl{h_1(x)})).
\end{equation*}
\item By \autoref{prop:factor_equiv_fiber}, we may choose
  $\zeta(a) \defeq \apfunc{\proj1}(\overline E(H,a)) : e(g_1(a)) = g_2(a)$.
  \label{item:orth-fact-2}
\item For every $x:X_1$, we have
\begin{equation*}
\proj2(E(H,h_1(x))({x,\refl{h_1(x)}})) :h_2(e(x))= h_1(x),
\end{equation*}
giving us a homotopy $\eta:h_2\circ e\htpy h_1$.
\item By the characterization of paths in fibers (\autoref{lem:hfib}), the path $\overline E(H,a)$ from \autoref{prop:factor_equiv_fiber} gives us
  $\eta(g_1(a)) = \apfunc{h_2}(\zeta(a)) \ct \opp{H(a)}$.
  The desired equality follows by substituting the definition of $H$ and rearranging paths.\qedhere
\end{enumerate}
\end{proof}

%

By standard arguments, this yields the following orthogonality principle.

\begin{thm}
  Let $e:A\to B$ be $n$-connected and $m:C\to D$ be $n$-truncated.
  Then the map
  \[ \varphi: (B\to C) \;\to\; \sm{h:A\to C}{k:B\to D} (m\circ h \htpy k \circ e) \]
  is an equivalence.
\end{thm}
\begin{proof}[Sketch of proof]
  For any $(h,k,H)$ in the codomain, let $h = h_2 \circ h_1$ and $k = k_2 \circ k_1$, where $h_1$ and $k_1$ are $n$-connected and $h_2$ and $k_2$ are $n$-truncated.
  Then $f = (m\circ h_2) \circ h_1$ and $f = k_2 \circ (k_1\circ e)$ are both $n$-factorizations of $m \circ h = k\circ e$.
  Thus, there is a unique equivalence between them.
  It is straightforward (if a bit tedious) to extract from this that $\hfib\varphi{(h,k,H)}$ is contractible.
\end{proof}

\index{orthogonal factorization system|)}%
\index{unique!factorization system|)}%

We end by showing that images are stable under pullback.
\index{image!stability under pullback}
\index{factorization!stability under pullback}

\begin{lem}\label{lem:hfiber_wrt_pullback}
Suppose that the square
\begin{equation*}
  \vcenter{\xymatrix{
      A\ar[r]\ar[d]_f &
      C\ar[d]^g\\
      B\ar[r]_-h &
      D
      }}
\end{equation*}
is a pullback square and let $b:B$. Then $\hfib{f}b\eqvsym\hfib{g}{h(b)}$.
\end{lem}

\begin{proof}
This follows from pasting of pullbacks (\autoref{ex:pullback-pasting}), since the type $X$ in the diagram
\begin{equation*}
  \vcenter{\xymatrix{
      X\ar[r]\ar[d] &
      A\ar[r]\ar[d]_f &
      C\ar[d]^g\\
      \unit\ar[r]_b &
      B\ar[r]_h &
      D
      }}
\end{equation*}
is the pullback of the left square if and only if it is the pullback of the outer rectangle, while $\hfib{f}b$ is the pullback of the square on the left and $\hfib{g}{h(b)}$ is the pullback of the outer rectangle.
\end{proof}

\begin{thm}\label{thm:stable-images}
\index{stability!of images under pullback}%
Consider functions $f:A\to B$, $g:C\to D$ and the diagram
\begin{equation*}
  \vcenter{\xymatrix{
      A\ar[r]\ar[d]_{\tilde{f}_n} &
      C\ar[d]^{\tilde{g}_n}\\
      \im_n(f)\ar[r]\ar[d]_{\proj1} &
      \im_n(g)\ar[d]^{\proj1}\\
      B\ar[r]_h &
      D
      }}
\end{equation*}
If the outer rectangle is a pullback, then so is the bottom square (and hence so is the top square, by \autoref{ex:pullback-pasting}). Consequently, images are stable under pullbacks.
\end{thm}

\begin{proof}
Assuming the outer square is a pullback, we have equivalences
\begin{align*}
B\times_D\im_n(g) & \jdeq \sm{b:B}{w:\im_n(g)} h(b)=\proj1 w\\
& \eqvsym \sm{b:B}{d:D}{w:\trunc n{\hfib{g}d}} h(b)= d\\
& \eqvsym \sm{b:B} \trunc n{\hfib{g}{h(b)}}\\
& \eqvsym \sm{b:B} \trunc n{\hfib{f}b} &&
\text{(by \autoref{lem:hfiber_wrt_pullback})}\\
& \equiv \im_n(f). && \qedhere
\end{align*}
\end{proof}

\index{n-type@$n$-type|)}%

\section{Modalities}
\label{sec:modalities}

\index{modality|(}

Nearly all of the theory of $n$-types and connectedness can be done in much greater generality.
This section will not be used in the rest of the book.

Our first thought regarding generalizing the theory of $n$-types might be to take \autoref{thm:trunc-reflective} as a definition.

\begin{defn}\label{defn:reflective-subuniverse}
  A \define{reflective subuniverse}
  \indexdef{reflective!subuniverse}%
  \indexdef{subuniverse, reflective}%
  is a predicate $P:\type\to\prop$ such that
  for every $A:\type$ we have a type $\reflect A$ such that $P(\reflect A)$ and a map
  $\project_A:A\to\reflect A$, with the property that for every $B:\type$ with $P(B)$, the following map is an equivalence:
  \[\function{(\reflect A\to{}B)}{(A\to{}B)}{f}{f\circ\project_A}.\]
\end{defn}

We write $\P \defeq \setof{A:\type | P(A)}$, so $A:\P$ means that $A:\type$ and we have $P(A)$.
We also write $\rec{\modal}$ for the quasi-inverse of the above map.
The notation $\reflect$ may seem slightly odd, but it will make more sense soon.

For any reflective subuniverse, we can prove all the familiar facts about reflective subcategories from category theory, in the usual way.
For instance, we have:
\begin{itemize}
\item A type $A$ lies in $\P$ if and only if $\project_A:A\to\reflect A$ is an equivalence.
\item $\P$ is closed under retracts.
  In particular, $A$ lies in $\P$ as soon as $\project_A$ admits a retraction.
\item The operation $\reflect$ is a functor in a suitable up-to-coherent-homotopy sense, which we can make precise at as high levels as necessary.
\item The types in $\P$ are closed under all limits such as products and pullbacks.
  In particular, for any $A:\P$ and $x,y:A$, the identity type $(x=_A y)$ is also in $\P$, since it is a pullback of two functions $\unit\to A$.
\item Colimits in $\P$ can be constructed by applying $\reflect$ to ordinary colimits of types.
\end{itemize}

Importantly, closure under products extends also to ``infinite products'', i.e.\ dependent function types.

\begin{thm}\label{thm:reflsubunv-forall}
  If $B:A\to\P$ is any family of types in a reflective subuniverse \P, then $\prd{x:A} B(x)$ is also in \P.
\end{thm}
\begin{proof}
  For any $x:A$, consider the function $\mathsf{ev}_x : (\prd{x:A} B(x)) \to B(x)$ defined by $\mathsf{ev}_x(f) \defeq f(x)$.
  Since $B(x)$ lies in $P$, this extends to a function
  \[ \rec{\modal}(\mathsf{ev}_x) : \reflect\Parens{\prd{x:A} B(x)} \to B(x). \]
  Thus we can define $h:\reflect(\prd{x:A} B(x)) \to \prd{x:A} B(x)$ by $h(z)(x) \defeq \rec{\modal}(\mathsf{ev}_x)(z)$.
  Then $h$ is a retraction of $\project_{\prd{x:A} B(x)}$, so that ${\prd{x:A} B(x)}$ is in $\P$.
\end{proof}

In particular, if $B:\P$ and $A$ is any type, then $(A\to B)$ is in \P.
In categorical language, this means that any reflective subuniverse is an \define{exponential ideal}.
\indexdef{exponential ideal}%
This, in turn, implies by a standard argument that the reflector preserves finite products.

\begin{cor}\label{cor:trunc_prod}
  For any types $A$ and $B$ and any reflective subuniverse, the induced map $\reflect(A\times B) \to \reflect(A) \times \reflect(B)$ is an equivalence.
\end{cor}
\begin{proof}
  It suffices to show that $\reflect(A) \times \reflect(B)$ has the same universal property as $\reflect(A\times B)$.
  Thus, let $C:\P$; we have
  \begin{align*}
    (\reflect(A) \times \reflect(B) \to C)
    &= (\reflect(A) \to (\reflect(B) \to C))\\
    &= (\reflect(A) \to (B \to C))\\
    &= (A \to (B \to C))\\
    &= (A \times B \to C)
  \end{align*}
  using the universal properties  of $\reflect(B)$ and $\reflect(A)$, along with the fact that $B\to C$ is in \P since $C$ is.
  It is straightforward to verify that this equivalence is given by composing with $\mreturn_A \times \mreturn_B$, as needed.
\end{proof}

It may seem odd that every reflective subcategory of types is automatically an exponential ideal, with a product-preserving reflector.
However, this is also the case classically in the category of \emph{sets}, for the same reasons.
It's just that this fact is not usually remarked on, since the classical category of sets---in contrast to the category of homotopy
types---does not have many interesting reflective subcategories.

Two basic properties of $n$-types are \emph{not} shared by general reflective subuniverses: \autoref{thm:ntypes-sigma} (closure under $\Sigma$-types) and \autoref{thm:truncn-ind} (truncation induction).
However, the analogues of these two properties are equivalent to each other.

\begin{thm}\label{thm:modal-char}
  For a reflective subuniverse \P, the following are logically equivalent.
  \begin{enumerate}
  \item If $A:\P$ and $B:A\to \P$, then $\sm{x:A} B(x)$ is in \P.\label{item:mchr1}
  \item for every $A:\type$, type family $B:\reflect A\to\P$, and map $g:\prd{a:A} B(\project(a))$, there exists $f:\prd{z:\reflect A} B(z)$ such that $f(\project(a)) = g(a)$ for all $a:A$.\label{item:mchr2}
  \end{enumerate}
\end{thm}
\begin{proof}
  Suppose~\ref{item:mchr1}.
  Then in the situation of~\ref{item:mchr2}, the type $\sm{z:\reflect A} B(z)$ lies in $\P$, and we have $g':A\to \sm{z:\reflect A} B(z)$ defined by $g'(a)\defeq (\project(a),g(a))$.
  Thus, we have $\rec{\modal}(g'):\reflect A \to \sm{z:\reflect A} B(z)$ such that $\rec{\modal}(g')(\project(a)) = (\project(a),g(a))$.

  Now consider the functions $\proj2 \circ \rec{\modal}(g') : \reflect A \to \reflect A$ and $\idfunc[\reflect A]$.
  By assumption, these become equal when precomposed with $\project$.
  Thus, by the universal property of $\reflect$, they are equal already, i.e.\ we have $p_z:\proj2(\rec{\modal}(g')(z)) = z$ for all $z$.
  Now we can define
  \narrowequation{f(z) \defeq \trans{p_z}{\proj2(\rec{\modal}(g')(z))},}
  and the second component of
  \narrowequation{\rec{\modal}(g')(\project(a)) = (\project(a),g(a))}
  yields $f(\project(a)) = g(a)$.

  Conversely, suppose~\ref{item:mchr2}, and that $A:\P$ and $B:A\to\P$.
  Let $h$ be the composite
  \[ \reflect\Parens{\sm{x:A} B(x)} \xrightarrow{\reflect(\proj1)} \reflect A \xrightarrow{\opp{(\project_A)}} A. \]
  Then for $z:\sm{x:A} B(x)$ we have
  \begin{align*}
    h(\project(z)) &= \opp\project(\reflect(\proj1)(\project(z)))\\
    &= \opp\project(\project(\proj1(z)))\\
    &= \proj1(z).
  \end{align*}
  Denote this path by $p_z$.
  Now if we define $C:\reflect(\sm{x:A} B(x)) \to \type$ by $C(w) \defeq B(h(w))$, we have
  \[ g \defeq \lam{z} \trans{p_z}{\proj2(z)} \;:\; \prd{z:\sm{x:A} B(x)} C(\project(z)). \]
  Thus, the assumption yields
  \narrowequation{f:\prd{w:\reflect(\sm{x:A}B(x))} C(w)}
  such that $f(\project(z)) = g(z)$.
  Together, $h$ and $f$ give a function
  \narrowequation{k:\reflect(\sm{x:A}B(x)) \to \sm{x:A}B(x)}
  defined by $k(w) \defeq (h(w),f(w))$, while $p_z$ and the equality $f(\project(z)) = g(z)$ show that $k$ is a retraction of $\project_{\sm{x:A}B(x)}$.
  Therefore, $\sm{x:A}B(x)$ is in \P.
\end{proof}

Note the similarity to the discussion in \autoref{sec:htpy-inductive}.
\index{recursion principle!for a modality}%
\index{induction principle!for a modality}%
\index{uniqueness!principle, propositional!for a modality}%
The universal property of the reflector of a reflective subuniverse is like a recursion principle with its uniqueness property, while \autoref{thm:modal-char}\ref{item:mchr2} is like the corresponding induction principle.
Unlike in \autoref{sec:htpy-inductive}, the two are not equivalent here, because of the restriction that we can only eliminate into types that lie in $\P$.
Condition~\ref{item:mchr1} of \autoref{thm:modal-char} is what fixes the disconnect.

Unsurprisingly, of course, if we have the induction principle, then we can derive the recursion principle.
We can also derive its uniqueness property, as long as we allow ourselves to eliminate into path types.
This suggests the following definition.
Note that any reflective subuniverse can be characterized by the operation $\reflect:\type\to\type$ and the functions $\project_A:A\to \reflect A$, since we have $P(A) = \isequiv(\project_A)$.

\begin{defn}\label{defn:modality}
A \define{modality}
\indexdef{modality}
is an operation $\modal:\type\to\type$ for which there are
\begin{enumerate}
\item functions $\mreturn^\modal_A:A\to\modal(A)$ for every type $A$.\label{item:modal1}
\item for every $A:\type$ and every type family $B:\modal(A)\to\type$, a function\label{item:modal2}
\begin{equation*}
\ind{\modal}:\Parens{\prd{a:A}\modal(B(\mreturn^\modal_A(a)))}\to\prd{z:\modal(A)}\modal(B(z)).
\end{equation*}
\item A path $\ind\modal(f)(\mreturn^\modal_A(a)) = f(a)$ for each $f:\prd{a:A}\modal(B(\mreturn^\modal_A(a)))$.\label{item:modal3}
\item For any $z,z':\modal(A)$, the function $\mreturn^\modal_{z=z'} : (z=z') \to \modal(z=z')$ is an equivalence.\label{item:modal4}
\end{enumerate}
We say that $A$ is \define{modal}
\indexdef{modal!type}%
\indexdef{type!modal}%
for $\modal$ if $\mreturn^\modal_A:A\to\modal(A)$ is an equivalence, and we write
\begin{equation}
  \modaltype\defeq\setof{X:\type | X \text{ is $\modal$-modal} }\label{eq:modaltype}
\end{equation}
for the type of modal types.
\end{defn}

Conditions~\ref{item:modal2} and~\ref{item:modal3} are very similar to \autoref{thm:modal-char}\ref{item:mchr2}, but phrased using $\modal B(z)$ rather than assuming $B$ to be valued in $\P$.
This allows us to state the condition purely in terms of the operation $\modal$, rather than requiring the predicate $P:\type\to\prop$ to be given in advance.
(It is not entirely satisfactory, since we still have to refer to $P$ not-so-subtly in clause~\ref{item:modal4}.
We do not know whether~\ref{item:modal4} follows from~\ref{item:modal1}--\ref{item:modal3}.)
However, the stronger-looking property of \autoref{thm:modal-char}\ref{item:mchr2} follows from \autoref{defn:modality}\ref{item:modal2} and~\ref{item:modal3}, since for any $C:\modal A \to \modaltype$ we have $C(z) \eqvsym \modal C(z)$, and we can pass back across this equivalence.

\index{universal!property!of a modality}%
As with other induction principles, this implies a universal property.

\begin{thm}\label{prop:lv_n_deptype_sec_equiv_by_precomp}
Let $A$ be a type and let $B:\modal(A)\to\modaltype$. Then the function
\begin{equation*}
(\blank\circ \mreturn^\modal_A) : \Parens{\prd{z:\modal(A)}B(z)} \to \Parens{\prd{a:A}B(\mreturn^\modal_A(a))}
\end{equation*}
is an equivalence.
\end{thm}
\begin{proof}
By definition, the operation $\ind{\modal}$ is a right inverse to $(\blank\circ \mreturn^\modal_A)$.
Thus, we only need to find a homotopy
\begin{equation*}
\prd{z:\modal(A)}s(z)= \ind{\modal}(s\circ \mreturn^\modal_A)(z)
\end{equation*}
for each $s:\prd{z:\modal(A)}B(z)$, exhibiting it as a left inverse as well.
By assumption, each $B(z)$ is modal, and hence each type $s(z)= R^\modal_X(s\circ \mreturn^\modal_A)(z)$
is also modal.
Thus, it suffices to find a function of type
\begin{equation*}
\prd{a:A}s(\mreturn^\modal_A(a))= \ind{\modal}(s\circ \mreturn^\modal_A)(\mreturn^\modal_A(a)).
\end{equation*}
which follows from \autoref{defn:modality}\ref{item:modal3}.
\end{proof}

In particular, for every type $A$ and every modal type $B$, we have an equivalence $(\modal A\to B)\eqvsym (A\to B)$.

\begin{cor}
  For any modality $\modal$, the $\modal$-modal types form a reflective subuniverse satisfying the equivalent conditions of \autoref{thm:modal-char}.
\end{cor}

Thus, modalities can be identified with reflective subuniverses closed under $\Sigma$-types.
The name \emph{modality} comes, of course, from \emph{modal logic}\index{modal!logic}, which studies logic where we can form statements such as ``possibly $A$'' (usually written $\diamond A$) or ``necessarily $A$'' (usually written $\Box A$).
The symbol $\modal$ is somewhat common for an arbitrary modal operator\index{modal!operator}. 
Under the propositions-as-types principle, a modality in the sense of modal logic corresponds to an operation on \emph{types}, and \autoref{defn:modality} seems a reasonable candidate for how such an operation should be defined.
(More precisely, we should perhaps call these \emph{idempotent, monadic} modalities; see the Notes.)
\index{idempotent!modality}%
As mentioned in \autoref{subsec:when-trunc}, we may in general use adverbs\index{adverb} to speak informally about such modalities, such as ``merely''\index{merely} for the propositional truncation and ``purely''\index{purely} for the identity modality
\index{identity!modality}%
\index{modality!identity}%
(i.e.\ the one defined by $\modal A \defeq A$).

For any modality $\modal$, we define a map $f:A\to B$ to be \define{$\modal$-connected}
\indexdef{function!.circle-connected@$\modal$-connected}%
\indexdef{.circle-connected function@$\modal$-connected function}%
if $\modal(\hfib f b)$ is contractible for all $b:B$, and to be \define{$\modal$-truncated}
\indexdef{function!.circle-truncated@$\modal$-truncated}%
\indexdef{.circle-truncated function@$\modal$-truncated function}%
if $\hfib f b$ is modal for all $b:B$.
All of the theory of \autoref{sec:connectivity,sec:image-factorization} which doesn't involve relating $n$-types for different values of $n$ applies verbatim in this generality.
\index{orthogonal factorization system}%
\index{unique!factorization system}%
In particular, we have an orthogonal factorization system.

An important class of modalities which does \emph{not} include the $n$-trun\-ca\-tions is the \emph{left exact} modalities: those for which the functor $\modal$ preserves pullbacks as well as finite products.
\index{topology!Lawvere-Tierney}%
These are a categorification of ``Lawvere-Tierney\index{Lawvere}\index{Tierney} topologies'' in elementary topos\index{topos} theory,
and correspond in higher-categorical semantics to sub-$(\infty,1)$-toposes.
\index{.infinity1-topos@$(\infty,1)$-topos}%
However, this is beyond the scope of this book.

Some particular examples of modalities other than $n$-truncation can be found in the exercises.

\index{modality|)}

\sectionNotes

The notion of homotopy $n$-type in classical homotopy theory is quite old.
It was Voevodsky who realized that the notion can be defined recursively in homotopy type theory, starting from contractibility.

\index{axiom!Streicher's Axiom K}%
The property ``Axiom K'' was so named by Thomas Streicher, as a property of identity types which comes after J, the latter being the traditional name for the eliminator of identity types.
\autoref{thm:hedberg} is due to Hedberg~\cite{hedberg1998coherence}; \cite{krausgeneralizations} contains more information and generalizations.

The notions of $n$-connected spaces and functions are also classical in homotopy theory, although as mentioned before, our indexing for connectedness of functions is off by one from the classical indexing.
The importance of the resulting factorization system has been emphasized by recent work in higher topos theory by Rezk, Lurie, and others.%
\index{.infinity1-topos@$(\infty,1)$-topos}
In particular, the results of this chapter should be compared with~\cite[\S6.5.1]{lurie:higher-topoi}.
In \autoref{sec:freudenthal}, the theory of $n$-connected maps will be crucial to our proof of the Freudenthal suspension theorem.

Modal operators\index{modal!operator} in \emph{simple} type theory have been studied extensively; see e.g.~\cite{modalTT}.  In the setting of dependent type theory, \cite{ab:bracket-types} treats the special case of propositional truncation ($(-1)$-truncation) as a modal operator\index{modal!operator}.  The development presented here greatly extends and generalizes this work, while drawing also on ideas from topos theory.\index{topos}

Generally, modal operators\index{modal!operator} come in (at least) two flavors: those such as $\diamond$ (``possibly'') for which $A\Rightarrow \diamond A$, and those such as $\Box$ (``necessarily'') for which $\Box A \Rightarrow A$.
When they are also \emph{idempotent} (i.e.\ $\diamond A = \diamond{\diamond A}$ or $\Box A = \Box{\Box A}$), the former may be identified with reflective subcategories (or equivalently, idempotent monads), and the latter with coreflective subcategories (or idempotent comonads).
\index{monad}
\index{comonad}
However, in dependent type theory it is trickier to deal with the comonadic sort, since they are more rarely stable under pullback, and thus cannot be interpreted as operations on the universe \UU.
Sometimes there are ways around this (see e.g.~\cite{QGFTinCHoTT12}), but for simplicity, here we stick to the monadic sort.

On the computational side, monads (and hence modalities\index{modality}) are used to model computational effects in functional programming~\cite{Moggi89}.%
\index{programming}%
\index{computational effect}%
A computation is said to be \emph{pure} if its execution results in no side effects (such as printing a message to the screen, playing music, or sending data over the Internet).
There exist ``purely functional'' programming languages, such as Haskell\index{Haskell}, in which it is technically only possible to write pure functions: side effects are represented by applying ``monads'' to output types.
For instance, a function of type $\mathsf{Int}\to\mathsf{Int}$ is pure, while a function of type $\mathsf{Int}\to \mathsf{IO}(\mathsf{Int})$ may perform input and output along the way to computing its result; the operation $\mathsf{IO}$ is a monad.
\index{purely}%
(This is the origin of our use of the adverb ``purely'' for the identity monad, since it corresponds computationally to pure functions with no side-effects.)
The modalities we have considered in this chapter are all idempotent, whereas those used in functional programming rarely are, but the ideas are still closely related.

\sectionExercises

\begin{ex}\ 
  \begin{enumerate}
    \item Use \autoref{thm:h-set-refrel-in-paths-sets} to show 
    that if $\brck{A}\to A$ for every type $A$, 
    then every type is a set.
    \item Show that if every surjective function (purely) splits, 
    i.e.~if
    \narrowequation{\prd{b:B}\brck{\hfib{f}{b}}\to\prd{b:B}\hfib{f}{b}}
    for every $f:A\to B$, then every type is a set.
  \end{enumerate}
\end{ex}

\begin{ex}\label{ex:s2-colim-unit}
  Express $\Sn^2$ as a colimit of a diagram consisting entirely of copies of $\unit$.
  Note that $\unit$ is a $(-2)$-type, while $\Sn^2$ is not expected to be an $n$-type for any finite $n$.
\end{ex}

\begin{ex}\label{ex:ntypes-closed-under-wtypes}
  Show that if $A$ is an $n$-type and $B:A\to \ntype{n}$ is a family of $n$-types, where $n\ge -1$, then the $W$-type $\wtype{a:A} B(a)$ (see \autoref{sec:w-types}) is also an $n$-type.
\end{ex}

\begin{ex}
  Use \autoref{prop:nconn_fiber_to_total} to extend \autoref{thm:connected-pointed} to any section-retraction pair.
\end{ex}

\begin{ex}\label{ex:ntype-from-nconn-const}
  Show that \autoref{thm:nconn-to-ntype-const} also works as a characterization in the other direction: $B$ is an $n$-type if and only if every map into $B$ from an $n$-con\-nect\-ed type is constant.
  Ideally, your proof should work for any modality as in \autoref{sec:modalities}.
\end{ex}

\begin{ex}\label{ex:connectivity-inductively}
  Prove that for $n\ge 0$, a type $A$ is $n$-connected if and only if it is $(-1)$-connected (i.e.\ merely inhabited) and for all $a,b:A$ the type $\id[A]ab$ is $(n-1)$-connected.
\end{ex}

\begin{ex}\label{ex:lemnm}
  \indexdef{excluded middle!LEMnm@$\LEM{n,m}$}%
  For $-1\le n,m \le\infty$, let $\LEM{n,m}$ denote the statement
  \[ \prd{A:\ntype{n}} \trunc m{A + \neg A},\]
  where $\ntype{\infty} \defeq \type$ and $\trunc{\infty}{X}\defeq X$.
  Show that:
  \begin{enumerate}
  \item If $n=-1$ or $m=-1$, then $\LEM{n,m}$ is equivalent to $\LEM{}$ from \autoref{sec:intuitionism}.
  \item If $n\ge 0$ and $m\ge 0$, then $\LEM{n,m}$ is inconsistent with univalence.
  \end{enumerate}
\end{ex}

\begin{ex}\label{ex:acnm}
  \indexdef{axiom!of choice!ACnm@$\choice{n,m}$}%
  For $-1\le n,m\le\infty$, let $\choice{n,m}$ denote the statement
  \[ \prd{X:\set}{Y:X\to\ntype{n}}
  \Parens{\prd{x:X} \trunc m{Y(x)}}
  \to
  \Trunc m{\prd{x:X} Y(x)},
  \]
  with conventions as in \autoref{ex:lemnm}.
  Thus $\choice{0,-1}$ is the axiom of choice from \autoref{sec:axiom-choice}, while $\choice{\infty,\infty}$ is \autoref{thm:ttac}.
  It is known that $\choice{\infty,-1}$ is consistent with univalence, since it holds in Voevodsky's simplicial model.
  \begin{enumerate}
  \item Without using univalence, show that $\LEM{n,\infty}$ implies $\choice{n,m}$ for all $m$.
    (On the other hand, in \autoref{subsec:emacinsets} we will show that $\choice{}=\choice{0,-1}$ implies $\LEM{}=\LEM{-1,-1}$.)
  \item Of course, $\choice{n,m}\Rightarrow \choice{k,m}$ if $k\le n$.
    Are there any other implications between the principles $\choice{n,m}$?
    Is $\choice{n,m}$ consistent with univalence for any $m\ge -1$?
    (These are open questions.)\index{open!problem}
  \end{enumerate}
\end{ex}

\begin{ex}\label{ex:acnm-surjset}
  Show that $\choice{n,-1}$ implies that for any $n$-type $A$, there merely exists a set $B$ and a surjection $B\to A$.
\end{ex}

\begin{ex}\label{ex:acconn}
  Define the \define{$n$-connected axiom of choice}
  \indexdef{n-connected@$n$-connected!axiom of choice}%
  \indexdef{axiom!of choice!n-connected@$n$-connected}%
  to be the statement
  \begin{quote}
    If $X$ is a set and $Y:X\to \type$ is a family of types such that each $Y(x)$ is $n$-connected, then $\prd{x:X} Y(x)$ is $n$-connected.
  \end{quote}
  Note that the $(-1)$-connected axiom of choice is $\choice{\infty,-1}$ from \autoref{ex:acnm}.
  \begin{enumerate}
  \item Prove that the $(-1)$-connected axiom of choice implies the $n$-con\-nect\-ed axiom of choice for all $n\ge -1$.
  \item Are there any other implications between the $n$-connected axioms of choice and the principles $\choice{n,m}$?
    (This is an open question.)\index{open!problem}
  \end{enumerate}
\end{ex}

\begin{ex}
  Show that the $n$-truncation modality is not left exact for any $n\ge -1$.
  That is, exhibit a pullback which it fails to preserve.
\end{ex}

\begin{ex}
  Show that $X\mapsto (\neg\neg X)$ is a modality.\index{modal!operator}%
\end{ex}

\begin{ex}
  Let $P$ be a mere proposition.
  \begin{enumerate}
  \item Show that $X\mapsto (P\to X)$ is a left exact modality.
    This is called the \define{open modality}
    \indexdef{open!modality}%
    \indexdef{modality!open}%
    associated to $P$.
  \item Show that $X\mapsto P*X$ is a left exact modality, where $*$ denotes the join (see \autoref{sec:colimits}).
    This is called the \define{closed modality}
    \indexdef{closed!modality}%
    \indexdef{modality!closed}%
    associated to $P$.
  \end{enumerate}
\end{ex}

\begin{ex}
  Let $f:A\to B$ be a map; a type $Z$ is \define{$f$-local}
  \indexdef{f-local type@$f$-local type}%
  \indexdef{type!f-local@$f$-local}%
  if $(\blank\circ f):(B\to Z) \to (A\to Z)$ is an equivalence.
  \begin{enumerate}
  \item Prove that the $f$-local types form a reflective subuniverse.
    You will want to use a higher inductive type to define the reflector (localization).
  \item Prove that if $B=\unit$, then this subuniverse is a modality.
  \end{enumerate}
\end{ex}


\cleartooddpage[\thispagestyle{empty}] 
\part{Mathematics}
\label{part:mathematics}

\chapter{Homotopy theory}
\label{cha:homotopy}

\index{acceptance|(}

In this chapter, we develop some homotopy theory within type theory.  We
use the \emph{synthetic approach} to homotopy theory introduced in
\cref{cha:basics}: Spaces, points, paths, and homotopies are basic
notions, which are represented by types and elements of types, particularly
the identity type.  The algebraic structure of paths and homotopies is
represented by the natural $\infty$-groupoid
\index{.infinity-groupoid@$\infty$-groupoid}%
structure on types, which is generated
by the rules for the identity type.  Using higher inductive types, as
introduced in \cref{cha:hits}, we can describe spaces directly by their
universal properties.  

\index{synthetic mathematics}%
There are several interesting aspects of this synthetic approach.
First, it combines advantages of concrete models (such as topological
spaces\index{topological!space}
or simplicial sets)\index{simplicial!sets}
with advantages of abstract categorical frameworks
for homotopy theory (such as Quillen model categories).\index{Quillen model category} 
 On the one hand,
our proofs feel elementary, and refer concretely to
points, paths, and homotopies in types.  On the other hand, our approach nevertheless abstracts away from
any concrete presentation of these objects --- for example,
associativity of path concatenation is proved by path induction, rather
than by reparametrization of maps $[0,1] \to X$ or by horn-filling conditions.
Type theory seems to be a very convenient way to study the abstract homotopy theory
of $\infty$-groupoids: by using the rules for the identity type, we
can avoid the complicated combinatorics involved in many definitions of
$\infty$-groupoids, and explicate only as much of the
structure as is needed in any particular proof.  

The abstract nature of type theory means that our proofs apply automatically in a variety of settings.
In particular, as mentioned previously, homotopy type theory has one interpretation in
Kan\index{Kan complex} simplicial sets\index{simplicial!sets},
which is one model for the homotopy theory of $\infty$-groupoids.  Thus,
our proofs apply to this model, and transferring them along the geometric
realization\index{geometric realization} functor from simplicial sets to topological spaces gives
proofs of corresponding theorems in classical homotopy theory.
However, though the details are work in progress, we can also
interpret type theory in a wide variety of other categories
that look like the category of $\infty$-groupoids, such as
$(\infty,1)$-toposes\index{.infinity1-topos@$(\infty,1)$-topos}.
Thus, proving a result in type theory will show
that it holds in these settings as well.
This sort of extra generality is well-known as a property of ordinary categorical logic:
univalent foundations extends it to homotopy theory as well.

Second, our synthetic approach has suggested new type-theoretic
methods and proofs.  Some of our proofs are fairly
direct transcriptions of classical proofs.  Others have a more
type-theoretic feel, and consist mainly of calculations with
$\infty$-groupoid operations, in a style that is very similar to how
computer scientists use type theory to reason about computer programs.
One thing that seems to have permitted these new proofs is the fact that type theory
emphasizes different aspects of homotopy theory than other approaches:
while tools like path induction and the universal properties of higher
inductives are available in a setting like Kan\index{Kan complex} simplicial sets, type
theory elevates their importance, because they are the \emph{only}
primitive tools available for working with these types.  Focusing on
these tools had led to new descriptions of familiar constructions such
as the universal cover of the circle and the Hopf fibration, using just the
recursion principles for higher inductive types.  These descriptions are
very direct, and many of the proofs in this chapter involve
computational calculations with such fibrations.
\index{mathematics!constructive}%
Another new aspect of our proofs is that they are constructive (assuming
univalence and higher inductives types are constructive); we describe an
application of this to homotopy groups of spheres in
\cref{sec:moreresults}.    

\index{mathematics!formalized}%
\indexsee{formalization of mathematics}{mathematics, formalized}%
Third, our synthetic approach is very amenable to computer-checked
proofs in proof assistants\index{proof!assistant} such as \Coq and \Agda.
Almost all of the proofs
described in this chapter have been computer-checked, and many of these
proofs were first given in a proof assistant, and then ``unformalized''
for this book.  The computer-checked proofs are comparable in length and
effort to the informal proofs presented here, and in some cases they are
even shorter and easier to do.  

\mentalpause

Before turning to the presentation of our results, we briefly review some
 basic concepts and theorems from homotopy theory for the benefit of the reader who is not familiar with them.
 We also give an overview of the
results proved in this chapter.  

\index{classical!homotopy theory|(}%
Homotopy theory is a branch of algebraic topology, and uses tools from abstract algebra,
such as group theory, to investigate properties of spaces.  One question
homotopy theorists investigate is how to tell whether two spaces are the
same, where ``the same'' means \emph{homotopy equivalence}
\index{homotopy!equivalence!topological}%
(continuous
maps back and forth that compose to the identity up to homotopy---this
gives the opportunity to ``correct'' maps that don't exactly compose to
the identity).  One common way to tell whether two spaces are the same
is to calculate \emph{algebraic invariants} associated with a space,
which include its \emph{homotopy groups} and \emph{homology} and
\emph{cohomology groups}.
\index{homology}%
\index{cohomology}%
Equivalent spaces have isomorphic
homotopy/(co)homology groups, so if two spaces have different groups,
then they are not equivalent.  Thus, these algebraic invariants provide
global information about a space, which can be used to tell spaces
apart, and complements the local information provided by notions such as
continuity.  For example, the torus locally looks like the $2$-sphere, but it
has a global difference, because it as a hole in it, and this difference
is visible in the homotopy groups of these two spaces.

The simplest example of a homotopy group is the \emph{fundamental group}
\index{fundamental!group}%
of a space, which is written $\pi_1(X,x_0$): Given a space $X$ and a
point $x_0$ in it, one can make a group whose elements are loops at
$x_0$ (continuous paths from $x_0$ to $x_0$), considered up to homotopy, with the
group operations given by the identity path (standing still), path
concatenation, and path reversal.  For example, the fundamental group of
the $2$-sphere is trivial, but the fundamental group of the torus is not,
which shows that the sphere and the torus are not homotopy equivalent.
The intuition is that every loop on the sphere is homotopic to the
identity, because its inside can be filled in.  In contrast, a loop on
the torus that goes through the donut's hole is not homotopic to the
identity, so there are non-trivial elements in the fundamental group.

\index{homotopy!group}
The \emph{higher homotopy groups} provide additional information about a
space.  Fix a point $x_0$ in $X$, and consider the constant path
$\refl{x_0}$.  Then the homotopy classes of homotopies between $\refl{x_0}$
and itself form a group $\pi_2(X,x_0)$, which tells us something about
the two-dimensional structure of the space.  Then $\pi_3(X,x_0$) is the
group of homotopy classes of homotopies between homotopies, and so on.
One of the basic problems of algebraic topology  is
\emph{calculating the homotopy groups of a space $X$}, which means
giving a group isomorphism between $\pi_k(X,x_0)$ and some more direct
description of a group (e.g., by a multiplication table or
presentation).  Somewhat surprisingly, this is a very difficult
question, even for spaces as simple as the spheres.  As can be seen from
\cref{tab:homotopy-groups-of-spheres}, some patterns emerge in the
higher homotopy groups of spheres, but there is no general formula, and
many homotopy groups of spheres are currently still unknown.

\bgroup
\definecolor{xA}{\OPTcolormodel}{\OPTcolxA}
\definecolor{xB}{\OPTcolormodel}{\OPTcolxB}
\definecolor{xC}{\OPTcolormodel}{\OPTcolxC}
\definecolor{xD}{\OPTcolormodel}{\OPTcolxD}
\definecolor{xE}{\OPTcolormodel}{\OPTcolxE}
\definecolor{xF}{\OPTcolormodel}{\OPTcolxF}
\definecolor{xG}{\OPTcolormodel}{\OPTcolxG}
\definecolor{xH}{\OPTcolormodel}{\OPTcolxH}
\definecolor{xI}{\OPTcolormodel}{\OPTcolxI}
\definecolor{xJ}{\OPTcolormodel}{\OPTcolxJ}
\definecolor{xK}{\OPTcolormodel}{\OPTcolxK}
\definecolor{xL}{\OPTcolormodel}{\OPTcolxL}
\definecolor{xM}{\OPTcolormodel}{\OPTcolxM}

\newcommand{\cA}{\colorbox{xG}{\hbox to 20pt {\hfil$0$\hfil}}}
\newcommand{\cB}{\colorbox{xH}{\hbox to 20pt {\hfil$\Z$\hfil}}}
\newcommand{\cC}{\colorbox{xI}{\hbox to 20pt {\hfil$\Z_{2}$\hfil}}}
\newcommand{\cD}{\colorbox{xJ}{\hbox to 20pt {\hfil$\Z_{2}$\hfil}}}
\newcommand{\cE}{\colorbox{xK}{\hbox to 20pt {\hfil$\Z_{24}$\hfil}}}
\newcommand{\cF}{\colorbox{xL}{\hbox to 20pt {\hfil$0$\hfil}}}
\newcommand{\cG}{\colorbox{xM}{\hbox to 20pt {\hfil$0$\hfil}}}
\newcommand{\cH}{\colorbox{xF}{\hbox to 20pt {\hfil$0$\hfil}}}
\newcommand{\cI}{\colorbox{xE}{\hbox to 20pt {\hfil$0$\hfil}}}
\newcommand{\cJ}{\colorbox{xD}{\hbox to 20pt {\hfil$0$\hfil}}}
\newcommand{\cK}{\colorbox{xC}{\hbox to 20pt {\hfil$0$\hfil}}}
\newcommand{\cL}{\colorbox{xB}{\hbox to 20pt {\hfil$0$\hfil}}}
\newcommand{\cM}{\colorbox{xA}{\hbox to 20pt {\hfil$0$\hfil}}}

\begin{table}[htb]
\centering\small
\begin{tabular}{p{15pt}>{\centering\arraybackslash}p{\OPTspherescolwidth}>{\centering\arraybackslash}p{\OPTspherescolwidth}>{\centering\arraybackslash}p{\OPTspherescolwidth}>{\centering\arraybackslash}p{\OPTspherescolwidth}>{\centering\arraybackslash}p{\OPTspherescolwidth}>{\centering\arraybackslash}p{\OPTspherescolwidth}>{\centering\arraybackslash}p{\OPTspherescolwidth}>{\centering\arraybackslash}p{\OPTspherescolwidth}>{\centering\arraybackslash}p{\OPTspherescolwidth}}
\toprule
           & $\Sn^0$ & $\Sn^1$ & $\Sn^{2}$ & $\Sn^{3}$ & $\Sn^{4}$ & $\Sn^{5}$ & $\Sn^{6}$ & $\Sn^{7}$ & $\Sn^{8}$ \\ \addlinespace[3pt] \midrule
$\pi_{1}$  & $0$     & $\Z$    & \cA       & \cH       & \cI       & \cJ       & \cK       & \cL       & \cM       \\ \addlinespace[3pt]
$\pi_{2}$  & $0$     & $0$     & \cB       & \cA       & \cH       & \cI       & \cJ       & \cK       & \cL       \\ \addlinespace[3pt]
$\pi_{3}$  & $0$     & $0$     & $\Z$      & \cB       & \cA       & \cH       & \cI       & \cJ       & \cK       \\ \addlinespace[3pt]
$\pi_{4}$  & $0$     & $0$     & $\Z_{2}$  & \cC       & \cB       & \cA       & \cH       & \cI       & \cJ       \\ \addlinespace[3pt]
$\pi_{5}$  & $0$     & $0$     & $\Z_{2}$  & $\Z_{2}$  & \cC       & \cB       & \cA       & \cH       & \cI       \\ \addlinespace[3pt]
$\pi_{6}$  & $0$     & $0$     & $\Z_{12}$ & $\Z_{12}$ & \cD       & \cC       & \cB       & \cA       & \cH       \\ \addlinespace[3pt]
$\pi_{7}$  & $0$     & $0$     & $\Z_{2}$  & $\Z_{2}$  & {\footnotesize $\Z {\times} \Z_{12}$} & \cD & \cC & \cB     & \cA    \\ \addlinespace[3pt]
$\pi_{8}$  & $0$     & $0$     & $\Z_{2}$  & $\Z_{2}$  & $\Z_{2}^{2}$ & \cE & \cD & \cC & \cB \\ \addlinespace[3pt]
$\pi_{9}$  & $0$     & $0$     & $\Z_{3}$  & $\Z_{3}$  & $\Z_{2}^{2}$ & $\Z_{2}$ & \cE & \cD & \cC \\ \addlinespace[3pt]
$\pi_{10}$ & $0$     & $0$     & $\Z_{15}$ & $\Z_{15}$ & \footnotesize{$\Z_{24} {\times} \Z_{3}$} & $\Z_{2}$ & \cF & \cE & \cD \\ \addlinespace[3pt]
$\pi_{11}$ & $0$     & $0$     & $\Z_{2}$  & $\Z_{2}$  & $\Z_{15}$ & $\Z_{2}$ & $\Z$ & \cF & \cE \\ \addlinespace[3pt]
$\pi_{12}$ & $0$     & $0$     & $\Z_{2}^{2}$ & $\Z_{2}^{2}$ & $\Z_{2}$ & $\Z_{30}$ & $\Z_{2}$ & \cG & \cF \\ \addlinespace[3pt]
$\pi_{13}$ & $0$     & $0$     & {\footnotesize $\Z_{12} {\times} \Z_{2}$} & {\footnotesize $\Z_{12} {\times} \Z_{2}$} & $\Z_{2}^{3}$ & $\Z_{2}$ & $\Z_{60}$ & $\Z_{2}$ & \cG \\ \addlinespace[3pt]
\bottomrule
\end{tabular}

\caption{Homotopy groups of spheres~\cite{wikipedia-groups}.\index{homotopy!group!of sphere}
The $k^{\textrm{th}}$ homotopy group $\pi_k$ of the $n$-dimensional sphere
$\Sn^n$ is isomorphic to the group listed in each entry, where $\Z$ is
the additive group of integers,\index{integers} and $\Z_{m}$ is the cyclic group\index{group!cyclic}\index{cyclic group} of order~$m$.
}
\label{tab:homotopy-groups-of-spheres}
\end{table}
\egroup

One way of understanding this complexity is through the correspondence
between spaces and $\infty$-groupoids
\index{.infinity-groupoid@$\infty$-groupoid}%
introduced in \cref{cha:basics}.
As discussed in \cref{sec:circle}, the 2-sphere is presented by a higher
inductive type with one point and one 2-dimensional loop.  Thus, one
might wonder why $\pi_3(\Sn ^2)$ is $\Z$, when the type $\Sn ^2$ has no
generators creating 3-dimensional cells.  It turns out that the
generating element of $\pi_3(\Sn ^2)$ is constructed using the
interchange law described in the proof of \cref{thm:EckmannHilton}: the
algebraic structure of an $\infty$-groupoid includes non-trivial
interactions between levels, and these interactions create elements of
higher homotopy groups.

\index{classical!homotopy theory|)}%

Type theory provides a natural setting for investigating this structure, as
we can easily define the higher homotopy groups.  Recall from \cref{def:loopspace} that for $n:\N$, the
$n$-fold iterated loop space\index{loop space!iterated} of a pointed type $(A,a)$ is defined
recursively by:
\begin{align*}
  \Omega^0(A,a)&=(A,a)\\
  \Omega^{n+1}(A,a)&=\Omega^n(\Omega(A,a)).
\end{align*}
This gives a \emph{space} (i.e.\ a type) of $n$-dimensional loops\index{loop!n-@$n$-}, which itself has
higher homotopies.  We obtain the set of $n$-dimensional loops by
truncation (this was also defined as an example in
\autoref{sec:free-algebras}):

\begin{defn}[Homotopy Groups]\label{def-of-homotopy-groups}
  Given $n\ge 1$ and $(A,a)$ a pointed type, we define the \define{homotopy groups} of $A$
  \indexdef{homotopy!group}%
  at $a$ by
  \[\pi_n(A,a)\defeq \Trunc0{\Omega^n(A,a)}\]
\end{defn}

\noindent
Since $n\ge 1$, the path concatenation and inversion operations on
$\Omega^n(A)$ induce operations on $\pi_n(A)$ making it into a group in
a straightforward way.  If $n\ge 2$, then the group $\pi_n(A)$ is
abelian\index{group!abelian}, by the Eckmann--Hilton argument \index{Eckmann--Hilton argument} (\autoref{thm:EckmannHilton}).
It is convenient to also write $\pi_0(A) \defeq \trunc0 A$,
but this case behaves somewhat differently: not only is it not a group,
it is defined without reference to any basepoint in $A$.

\index{presentation!of an infinity-groupoid@of an $\infty$-groupoid}%
This definition is a suitable one for investigating homotopy groups
because the (higher) inductive definition of a type $X$ presents $X$ as
a free type, analogous to a free $\infty$-groupoid,
\index{.infinity-groupoid@$\infty$-groupoid}%
and this
presentation \emph{determines} but does not \emph{explicitly describe}
the higher identity types of $X$.  The identity types are populated by
both the generators ($\lloop$, for the circle) and the results of applying to them all of the groupoid
operations (identity, composition, inverses, associativity, interchange,
\ldots).  Thus, the higher-inductive presentation of a space allows us
to pose the question ``what does the identity type of $X$ really turn out
to be?'' though it can take some significant mathematics to answer it.
This is a higher-dimensional generalization of a familiar fact in type
theory: characterizing the identity type of $X$ can take some work,
even if $X$ is an ordinary inductive type, such as the natural numbers
or booleans.  For example, the theorem that $\bfalse$ is different
from $\btrue$ does not follow immediately from the definition;
see \autoref{sec:compute-coprod}.

\index{univalence axiom}%
The univalence axiom plays an essential role in calculating homotopy
groups (without univalence, type theory is compatible with an
interpretation where all paths, including, for example, the loop on the
circle, are reflexivity).  We will see this in the calculation of the
fundamental group of the circle below: the map from $\Omega(\Sn^1)$ to $\Z$ is defined by mapping a loop on the circle to an
automorphism\index{automorphism!of Z, successor@of $\Z$, successor} of the set $\Z$, so that, for example, $\lloop \ct \opp
\lloop$ is sent to $\mathsf{successor} \ct \mathsf{predecessor}$ (where
$\mathsf{successor}$ and $\mathsf{predecessor}$ are automorphisms of
$\Z$ viewed, by univalence, as paths in the universe), and then applying
the automorphism to 0. Univalence produces non-trivial paths in the
universe, and this is used to extract information from paths in higher
inductive types.

In this chapter, we first calculate some homotopy groups of spheres,
including $\pi_k(\Sn ^1)$ (\autoref{sec:pi1-s1-intro}), $\pi_{k}(\Sn
^n)$ for $k<n$ (\autoref{sec:conn-susp,sec:pik-le-n}), $\pi_2(\Sn ^2)$ and $\pi_3(\Sn ^2)$ by
way of the Hopf fibration (\autoref{sec:hopf}) and a long-exact-sequence
argument (\autoref{sec:long-exact-sequence-homotopy-groups}), and
$\pi_n(\Sn ^n)$ by way of the Freudenthal suspension theorem
(\autoref{sec:freudenthal}).  Next, we discuss the van Kampen theorem
(\autoref{sec:van-kampen}), which characterizes the fundamental group of
a pushout, and the status of Whitehead's principle (when is a map that
induces an equivalence on all homotopy groups an equivalence?)
(\autoref{sec:whitehead}).  Finally, we include brief summaries of
additional results that are not included in the book, such as
$\pi_{n+1}(\Sn ^n)$ for $n\ge 3$, the Blakers--Massey theorem, and a
construction of Eilenberg--Mac Lane spaces (\autoref{sec:moreresults}).
Prerequisites for this chapter include \autoref{cha:typetheory,cha:basics,cha:hits,cha:hlevels} as well as parts of \autoref{cha:logic}.

\index{fundamental!group!of circle|(}
\section{\texorpdfstring{$\pi_1(S^1)$}{π₁(S¹)}}
\label{sec:pi1-s1-intro}

In this section, our goal is to show that $\id {\pi_1(\Sn ^1)} {\Z}$.
In fact, we will show that the loop space\index{loop space} ${\Omega(\Sn ^1)}$ is equivalent to $\Z$.
This is a stronger statement, because $\id{\pi_1(\Sn ^1)} {\trunc 0 {\Omega(\Sn ^1)}}$ by
definition; so if $\id {\Omega(\Sn ^1)} {\Z}$, then $\id {\trunc
  0 {\Omega(\Sn ^1)}} {\trunc 0 {\Z}}$ by congruence, and
$\Z$ is a set by definition (being a set-quotient; see \autoref{defn-Z,Z-quotient-by-canonical-representatives}), so $\id {\trunc 0 {\Z}} {\Z}$.
Moreover, knowing that ${\Omega(\Sn ^1)}$ is a set will imply that $\pi_n(\Sn^1)$ is trivial for $n>1$, so we will actually have calculated \emph{all} the homotopy groups of $\Sn^1$.

\subsection{Getting started}
\label{sec:pi1s1-initial-thoughts}

It is not too hard to define functions in both directions between $\Omega(\Sn^1)$ and \Z.
By specializing \autoref{thm:looptothe} to $\lloop:\base=\base$, we have a function $\lloop^{\blank} : \Z \rightarrow (\id{\base}{\base})$ defined (loosely speaking) by
\[
  \lloop^n =
  \begin{cases}
    \underbrace{\lloop \ct \lloop \ct \cdots \ct \lloop}_{n}  & \text{if $n > 0$,} \\
    \underbrace{\opp \lloop \ct \opp \lloop \ct \cdots \ct \opp \lloop}_{-n} & \text{if $n < 0$,} \\
    \refl{\base} & \text{if $n = 0$.}
\end{cases}
\]
Defining a function $g:\Omega(\Sn^1)\to\Z$ in the other direction is a bit trickier.
Note that the successor function $\Zsuc:\Z\to\Z$ is an equivalence,
\index{successor!isomorphism on Z@isomorphism on $\Z$}%
and hence induces a path $\ua(\Zsuc):\Z=\Z$ in the universe \type.
Thus, the recursion principle of $\Sn^1$ induces a map $c:\Sn^1\to\type$ by $c(\base)\defeq \Z$ and $\apfunc c (\lloop) \defid \ua(\Zsuc)$.
Then we have $\apfunc{c} : (\base=\base) \to (\Z=\Z)$, and we can define $g(p)\defeq \transfib{X\mapsto X}{\apfunc{c}(p)}{0}$.

With these definitions, we can even prove that $g(\lloop^n)=n$ for any $n:\Z$, using the induction principle \autoref{thm:sign-induction} for $n$.
(We will prove something more general a little later on.)
However, the other equality $\lloop^{g(p)}=p$ is significantly harder.
The obvious thing to try is path induction, but path induction does not apply to loops such as $p:(\base=\base)$ that have \emph{both} endpoints fixed!
A new idea is required, one which can be explained both in terms of classical homotopy theory and in terms of type theory.
We begin with the former.

\subsection{The classical proof}
\label{sec:pi1s1-classical-proof}

\index{classical!homotopy theory|(}%
In classical homotopy theory, there is a standard proof of $\pi_1(\Sn^1)=\Z$ using universal covering spaces.
Our proof can be regarded as a type-theoretic version of this proof, with covering spaces appearing here as fibrations whose fibers are sets.
\index{fibration}%
\index{total!space}%
Recall that \emph{fibrations} over a space $B$ in homotopy theory correspond to type families $B\to\type$ in type theory.
\index{path!fibration}%
\index{fibration!of paths}%
In particular, for a point $x_0:B$, the type family $(x\mapsto (x_0=x))$ corresponds to the \emph{path fibration} $P_{x_0} B \to B$, in which the points of $P_{x_0} B$ are paths in $B$ starting at $x_0$, and the map to $B$ selects the other endpoint of such a path.
This total space $P_{x_0} B$ is contractible, since we can ``retract'' any path to its initial endpoint $x_0$ --- we have seen the type-theoretic version of this as \autoref{thm:contr-paths}.
Moreover, the fiber over $x_0$ is the loop space\index{loop space} $\Omega(B,x_0)$ --- in type theory this is obvious by definition of the loop space.

\begin{figure}\centering
  \begin{tikzpicture}[xscale=1.4,yscale=.6]
    \node (R) at (2,1) {$\mathbb{R}$};
    \node (S1) at (2,-2) {$S^1$};
    \draw[->] (R) -- node[auto] {$w$} (S1);
    \draw (0,-2) ellipse (1 and .4);
    \draw[dotted] (1,0) arc (0:-30:1 and .8);
    \draw (1,0) arc (0:90:1 and .8) arc (90:270:1 and .3) coordinate (t1);
    \draw[white,line width=4pt] (t1) arc (-90:90:1 and .8);
    \draw (t1) arc (-90:90:1 and .8) arc (90:270:1 and .3) coordinate (t2);
    \draw[white,line width=4pt] (t2) arc (-90:90:1 and .8);
    \draw (t2) arc (-90:90:1 and .8) arc (90:270:1 and .3) coordinate (t3);
    \draw[white,line width=4pt] (t3) arc (-90:90:1 and .8);
    \draw (t3) arc (-90:-30:1 and .8) coordinate (t4);
    \draw[dotted] (t4) arc (-30:0:1 and .8);
    \node[fill,circle,inner sep=1pt,label={below:\scriptsize \base}] at (0,-2.4) {};
    \node[fill,circle,inner sep=1pt,label={above left:\scriptsize 0}] at (0,.2) {};
    \node[fill,circle,inner sep=1pt,label={above left:\scriptsize 1}] at (0,1.2) {};
    \node[fill,circle,inner sep=1pt,label={above left:\scriptsize 2}] at (0,2.2) {};
  \end{tikzpicture}
  \caption{The winding map in classical topology}\label{fig:winding}
\end{figure}

Now in classical homotopy theory, where $\Sn^1$ is regarded as a topological space, we may proceed as follows.
\index{winding!map}%
Consider the ``winding'' map $w:\mathbb{R}\to \Sn ^1$, which looks like a helix projecting down onto the circle (see \autoref{fig:winding}).
This map $w$ sends each point on the helix\index{helix} to the point on the circle that it is ``sitting above''.
It is a fibration, and the fiber over each point is isomorphic to the integers.\index{integers}
If we lift the path that goes counterclockwise around the loop on the bottom, we go up one level in the helix, incrementing the integer in the fiber.
Similarly, going clockwise around the loop on the bottom corresponds to going down one level in the helix, decrementing this count.
This fibration is called the \emph{universal cover} of the circle.
\index{universal!cover}%
\index{cover!universal}%
\index{covering space!universal}%

Now a basic fact in classical homotopy theory is that a map $E_1\to E_2$ of fibrations over $B$ which is a homotopy equivalence between $E_1$ and $E_2$ induces a homotopy equivalence on all fibers.
(We have already seen the type-theoretic version of this as well in \autoref{thm:total-fiber-equiv}.)
Since $\mathbb{R}$ and $P_{\base} S^1$ are both contractible topological spaces, they are homotopy equivalent, and thus their fibers $\Z$ and $\Omega(\Sn ^1)$ over the basepoint are also homotopy equivalent.

\index{classical!homotopy theory|)}%

\subsection{The universal cover in type theory}
\label{sec:pi1s1-universal-cover}

\index{universal!cover|(}%
\index{cover!universal|(}%
\index{covering space!universal|(}%

Let us consider how we might express the preceding proof in type theory.
We have already remarked that the path fibration of $\Sn^1$ is represented by the type family $(x\mapsto (\base=x))$.
We have also already seen a good candidate for the universal cover of $\Sn^1$: it's none other than the type family $c:\Sn^1\to\type$ which we defined in \autoref{sec:pi1s1-initial-thoughts}!
By definition, the fiber of this family over $\base$ is $\Z$, while the effect of transporting around $\lloop$ is to add one --- thus it behaves just as we would expect from \autoref{fig:winding}.

However, since we don't know yet that this family behaves like a universal cover is supposed to (for instance, that its total space is simply connected), we use a different name for it.
For reference, therefore, we repeat the definition.
\index{integers}

\begin{defn}[Universal Cover of $\Sn^1$] \label{S1-universal-cover}
  Define $\code : \Sn ^1 \to \type$ by circle-recursion, with 
  \begin{align*}
    \code(\base) &\defeq \Z \\
    \apfunc{\code}({\lloop}) &\defid \ua(\Zsuc).
  \end{align*}
\end{defn}

We emphasize briefly the definition of this family, since it is so different from how one usually defines covering spaces in classical homotopy theory.
To define a function by circle recursion, we need to find a point and a
loop in the codomain.  In this case, the codomain is $\type$, and the point
we choose is $\Z$, corresponding to our expectation that the
fiber of the universal cover should be the integers.  The loop we choose
is the successor/predecessor
\index{successor!isomorphism on Z@isomorphism on $\Z$}%
\index{predecessor!isomorphism on Z@isomorphism on $\Z$}%
isomorphism on $\Z$, which
corresponds to the fact that going around the loop in the base goes up
one level on the helix.  Univalence is necessary for this part of the
proof, because we need to convert a \emph{non-trivial} equivalence on $\Z$ into an identity.  

We call this the fibration of ``codes'', because its elements are combinatorial data that act as codes for paths on the circle: the integer $n$ codes for the path which loops around the circle $n$ times.

From this definition, it is simple to calculate that transporting with
$\code$ takes $\lloop$ to the successor function, and 
$\opp{\lloop}$ to the predecessor function:
\begin{lem} \label{lem:transport-s1-code}
\id{\transfib \code \lloop x} {x + 1} and 
\id{\transfib \code {\opp \lloop} x} {x - 1}.
\end{lem}
\begin{proof}
For the first equation, we calculate as follows:
\begin{align}
{\transfib \code \lloop x} 
&= \transfib {A \mapsto A} {(\ap{\code}{\lloop})} x \tag{by \autoref{thm:transport-compose}}\\
&= \transfib {A \mapsto A} {\ua (\Zsuc)} x \tag{by computation for $\rec{\Sn^1}$}\\
&= x + 1 \tag{by computation for \ua}.
\end{align}
The second equation follows from the first, because $\transfib{B}{p}{\blank}$ and $\transfib{B}{\opp p}{\blank}$ are always inverses, so $\transfib\code {\opp \lloop}{\blank}$ must be the inverse of $\Zsuc$.
\end{proof}

We can now see what was wrong with our first approach: we defined $f$ and $g$ only on the fibers $\Omega(\Sn^1)$ and \Z, when we should have defined a whole morphism \emph{of fibrations} over $\Sn^1$.
In type theory, this means we should have defined functions having types
\begin{align}
  \prd{x:\Sn^1} &((\base=x) \to \code(x)) \qquad\text{and/or} \label{eq:pi1s1-encode}\\
  \prd{x:\Sn^1} &(\code(x) \to (\base=x))\label{eq:pi1s1-decode}
\end{align}
instead of only the special cases of these when $x$ is \base.
This is also an instance of a common observation in type theory: when attempting to prove something about particular inhabitants of some inductive type, it is often easier to generalize the statement so that it refers to \emph{all} inhabitants of that type, which we can then prove by induction.
Looked at in this way, the proof of $\Omega(\Sn^1)=\Z$ fits into the same pattern as the characterization of the identity types of coproducts and natural numbers in \autoref{sec:compute-coprod,sec:compute-nat}.

At this point, there are two ways to finish the proof.
We can continue mimicking the classical argument by constructing~\eqref{eq:pi1s1-encode} or~\eqref{eq:pi1s1-decode} (it doesn't matter which), proving that a homotopy equivalence between total spaces induces an equivalence on fibers, and then that the total space of the universal cover is contractible.
The first type-theoretic proof of $\Omega(\Sn^1)=\Z$ followed this pattern; we call it the \emph{homotopy-theoretic} proof.

Later, however, we discovered that there is an alternative proof, which has a more type-theoretic feel and more closely follows the proofs in \autoref{sec:compute-coprod,sec:compute-nat}.
In this proof, we directly construct both~\eqref{eq:pi1s1-encode} and~\eqref{eq:pi1s1-decode}, and prove that they are mutually inverse by calculation.
We will call this the \emph{encode-decode} proof, because we call the functions~\eqref{eq:pi1s1-encode} and~\eqref{eq:pi1s1-decode} \emph{encode} and \emph{decode} respectively.
Both proofs use the same construction of the cover given above.
Where the classical proof induces an equivalence on fibers from an equivalence between total spaces, the encode-decode proof constructs the inverse map (\emph{decode}) explicitly as a map between fibers.
And where the classical proof uses contractibility, the encode-decode proof uses path induction, circle induction, and integer induction.
These are the same tools used to prove contractibility---indeed, path induction \emph{is} essentially contractibility of the path fibration composed with $\mathsf{transport}$---but they are applied in a different way.

Since this is a book about homotopy type theory, we present the encode-decode proof first.
A homotopy theorist who gets lost is encouraged to skip to the homotopy-theoretic proof (\autoref{subsec:pi1s1-homotopy-theory}).

\index{universal!cover|)}%
\index{cover!universal|)}%
\index{covering space!universal|)}%

\subsection{The encode-decode proof}
\label{subsec:pi1s1-encode-decode}

\index{encode-decode method|(}%
\indexsee{encode}{encode-decode method}%
\indexsee{decode}{encode-decode method}%
We begin with the function~\eqref{eq:pi1s1-encode} that maps paths to codes:
\begin{defn}
Define $\encode : \prd{x : \Sn ^1} (\base=x) \rightarrow  \code(x)$ by 
\[
\encode \: p \defeq \transfib{\code} p 0
\]
(we leave the argument $x$ implicit).  
\end{defn}
Encode is defined by lifting a path into the universal cover, which
determines an equivalence, and then applying the resulting equivalence
to $0$.  
The interesting thing about this function is that it computes a concrete
number from a loop on the circle, when this loop is represented using
the abstract groupoidal framework of homotopy type theory.  To gain an
intuition for how it does this, observe that by the above lemmas,
$\transfib \code \lloop x$ is the successor map and $\transfib \code {\opp
  \lloop} x$ is the predecessor map.
Further, $\mathsf{transport}$ is functorial (\autoref{cha:basics}), so
$\transfib{\code} {\lloop \ct \lloop}{\blank}$ is
\[(\transfib \code \lloop-) \circ (\transfib \code \lloop-)\]
and so on.
Thus, when $p$ is a composition like 
\[
\lloop \ct \opp \lloop \ct \lloop \ct \cdots
\]
$\transfib{\code}{p}{\blank}$ will compute a composition of functions like
\[
\Zsuc \circ \Zpred \circ \Zsuc \circ \cdots 
\]
Applying this composition of functions to 0 will compute the
\index{winding!number}%
\emph{winding number} of the path---how many times it goes around the
circle, with orientation marked by whether it is positive or negative,
after inverses have been canceled.  Thus, the computational behavior of
$\encode$ follows from the reduction rules for higher-inductive types and
univalence, and the action of $\mathsf{transport}$ on compositions and inverses.

Note that the instance $\encode' \defeq \encode_{\base}$ has type 
$(\id \base \base) \rightarrow \Z$.
This will be one half of our desired equivalence; indeed, it is exactly the function $g$ defined in \autoref{sec:pi1s1-initial-thoughts}.

Similarly, the function~\eqref{eq:pi1s1-decode} is a generalization of the function $\lloop^{\blank}$ from \autoref{sec:pi1s1-initial-thoughts}.

\begin{defn}\label{thm:pi1s1-decode}
Define $\decode : \prd{x : \Sn ^1}\code(x) \rightarrow (\base=x)$ by 
circle induction on $x$.  It suffices to give a function 
${\code(\base) \rightarrow (\base=\base)}$, for which we use $\lloop^{\blank}$, and 
to show that $\lloop^{\blank}$ respects the loop.  
\end{defn}

\begin{proof}
To show that $\lloop^{\blank}$ respects the loop, it suffices to give a path
from $\lloop^{\blank}$ to itself that lies over $\lloop$. 
By the definition of dependent paths, this means a path from
\[\transfib {(x' \mapsto \code(x') \rightarrow (\base=x'))} {\lloop} {\lloop^{\blank}}\]
to $\lloop^{\blank}$.  We define such a
path as follows:
\begin{align*}
  \MoveEqLeft \transfib {(x' \mapsto \code(x') \rightarrow (\base=x'))} \lloop {\lloop^{\blank}}\\
&= \transfibf {x'\mapsto (\base=x')}(\lloop) \circ {\lloop^{\blank}} \circ \transfibf \code ({\opp \lloop}) \\
&= (- \ct \lloop) \circ (\lloop^{\blank}) \circ \transfibf \code ({\opp \lloop}) \\
&= (- \ct \lloop) \circ (\lloop^{\blank}) \circ \Zpred \\
&= (n \mapsto \lloop^{n - 1} \ct \lloop).
\end{align*}
On the first line, we apply the characterization of $\mathsf{transport}$
when the outer connective of the fibration is $\rightarrow$, which
reduces the $\mathsf{transport}$ to pre- and post-composition with
$\mathsf{transport}$ at the domain and codomain types.  On the second line,
we apply the characterization of $\mathsf{transport}$ when the type family
is $x\mapsto \id{\base}{x}$, which is post-composition of paths.  On the third line,
we use the action of $\code$ on $\opp \lloop$ from
\autoref{lem:transport-s1-code}.  And on the fourth line, we simply
reduce the function composition.  Thus, it suffices to show that for all
$n$, \id{\lloop^{n - 1} \ct \lloop}{\lloop^{n}}, which is an easy
induction, using the groupoid laws.  
\end{proof}

We can now show that $\encode$ and $\decode$ are quasi-inverses.
What used to be the difficult direction is now easy!

\begin{lem} \label{lem:s1-decode-encode}  For all 
for all $x: \Sn ^1$ and $p : \id \base x$, $\id
{\decode_x({{\encode_x(p)}})} p$.  
\end{lem}

\begin{proof}
By path induction, it suffices to show that 
\narrowequation{\id {\decode_{\base}({{\encode_{\base}(\refl{\base})}})} {\refl{\base}}.}
But
\narrowequation{\encode_{\base}(\refl{\base}) \jdeq \transfib{\code}{\refl{\base}} 0 \jdeq 0,}
and $\decode_{\base}(0) \jdeq \lloop^ 0 \jdeq \refl{\base}$.  
\end{proof}

The other direction is not much harder.

\begin{lem} \label{lem:s1-encode-decode} For all 
$x: \Sn ^1$ and $c : \code(x)$, we have $\id
{\encode_x({{\decode_x(c)}})} c$.  
\end{lem}

\begin{proof}
The proof is by circle induction.  It suffices to show the case for
\base, because the case for \lloop is a path between paths in
$\Z$, which is immediate because $\Z$ is a set.  

Thus, it suffices to show, for all $n : \Z$, that
\[
\id {\encode'(\lloop^n)} {n}
\]
The proof is by induction, using \cref{thm:looptothe}.
\begin{itemize}

\item In the case for $0$, the result is true by definition.

\item In the case for $n+1$, 
\begin{align}
 {\encode'(\lloop^{n+1})}
&= {\encode'(\lloop^{n} \ct \lloop)} \tag{by definition of $\lloop^{\blank}$} \\
&= \transfib{\code}{(\lloop^{n} \ct \lloop)}{0} \tag{by definition of $\encode$}\\
&= \transfib{\code}{\lloop}{(\transfib{\code}{\lloop^n}{0})} \tag{by functoriality}\\
&= {(\transfib{\code}{\lloop^n}{0})} + 1 \tag{by \autoref{lem:transport-s1-code}}\\
&= n + 1. \tag{by the inductive hypothesis}
\end{align}

\item The case for negatives is analogous.  \qedhere
\end{itemize}
\end{proof}

Finally, we conclude the theorem.

\begin{thm} 
There is a family of equivalences $\prd{x : \Sn ^1} (\eqv {(\base=x)} {\code(x)})$.
\end{thm}
\begin{proof}
The maps $\encode$ and $\decode$ are quasi-inverses by
\autoref{lem:s1-decode-encode,lem:s1-decode-encode}.
\end{proof}

Instantiating at {\base} gives
\begin{cor}\label{cor:omega-s1}
$\eqv {\Omega(\Sn^1,\base)} {\Z}$.
\end{cor}

A simple induction shows that this equivalence takes addition to
composition, so that $\Omega(\Sn ^1) = \Z$ as groups.

\begin{cor} \label{cor:pi1s1}
$\id{\pi_1(\Sn ^1)} {\Z}$, while $\id{\pi_n(\Sn ^1)}0$ for $n>1$.
\end{cor}
\begin{proof}
For $n=1$, we sketched the proof from \autoref{cor:omega-s1} above.
For $n > 1$, we have $\trunc 0 {\Omega^{n}(\Sn ^1)} = \trunc 0 {\Omega^{n-1}(\Omega{\Sn ^1})} = \trunc 0 {\Omega^{n-1}(\Z)}$.
And since $\Z$ is a set, $\Omega^{n-1}(\Z)$ is contractible, so this is trivial.
\end{proof}

\index{encode-decode method|)}%

\subsection{The homotopy-theoretic proof}
\label{subsec:pi1s1-homotopy-theory}

In \autoref{sec:pi1s1-universal-cover}, we defined the putative universal cover $\code:\Sn^1\to\type$ in type theory, and in \autoref{subsec:pi1s1-homotopy-theory} we defined a map $\encode : \prd{x:\Sn^1} (\base=x) \to \code(x)$ from the path fibration to the universal cover.
What remains for the classical proof is to show that this map induces an equivalence on total spaces because both are contractible, and to deduce from this that it must be an equivalence on each fiber.

\index{total!space}%
In \autoref{thm:contr-paths} we saw that the total space $\sm{x:\Sn^1} (\base=x)$ is contractible.
For the other, we have:

\begin{lem}
  The type $\sm{x:\Sn^1}\code(x)$ is contractible.
\end{lem}
\begin{proof}
  We apply the flattening lemma (\autoref{thm:flattening}) with the following values:
  \begin{itemize}
  \item $A\defeq\unit$ and $B\defeq\unit$, with $f$ and $g$ the obvious functions.
    Thus, the base higher inductive type $W$ in the flattening lemma is equivalent to $\Sn^1$.
  \item $C:A\to\type$ is constant at \Z.
  \item $D:\prd{b:B} (\eqv{\Z}{\Z})$ is constant at $\Zsuc$.
  \end{itemize}
  Then the type family $P:\Sn^1\to\type$ defined in the flattening lemma is equivalent to $\code:\Sn^1\to\type$.
  Thus, the flattening lemma tells us that $\sm{x:\Sn^1}\code(x)$ is equivalent to a higher inductive type with the following generators, which we denote $R$:
  \begin{itemize}
  \item A function $\mathsf{c}: \Z \to R$.
  \item For each $z:\Z$, a path $\mathsf{p}_z:\mathsf{c}(z) = \mathsf{c}(\Zsuc(z))$.
  \end{itemize}
  We might call this type the \define{homotopical reals};
  \indexdef{real numbers!homotopical}%
  it plays the same role as the topological space $\mathbb{R}$ in the classical proof.

  Thus, it remains to show that $R$ is contractible.
  As center of contraction we choose $\mathsf{c}(0)$; we must now show that $x=\mathsf{c}(0)$ for all $x:R$.
  We do this by induction on $R$.
  Firstly, when $x$ is $\mathsf{c}(z)$, we must give a path $q_z:\mathsf{c}(0) = \mathsf{c}(z)$, which we can do by induction on $z:\Z$:
  \begin{align*}
    q_0 &\defeq \refl{\mathsf{c}(0)}\\
    q_{n+1} &\defeq q_n \ct \mathsf{p}_n & &\text{for $n\ge 0$}\\
    q_{n-1} &\defeq q_n \ct \opp{\mathsf{p}_{n-1}} & &\text{for $n\le 0$.}
  \end{align*}
  Secondly, we must show that for any $z:\Z$, the path $q_z$ is transported along $\mathsf{p}_z$ to $q_{z+1}$.
  By transport of paths, this means we want $q_z \ct \mathsf{p}_z = q_{z+1}$.
  This is easy by induction on $z$, using the definition of $q_z$.
  This completes the proof that $R$ is contractible, and thus so is $\sm{x:\Sn^1}\code(x)$.
\end{proof}

\begin{cor}\label{thm:encode-total-equiv}
  The map induced by \encode:
  \[ \tsm{x:\Sn^1} (\base=x) \to \tsm{x:\Sn^1}\code(x) \]
  is an equivalence.
\end{cor}
\begin{proof}
  Both types are contractible.
\end{proof}

\begin{thm}
  $\eqv {\Omega(\Sn^1,\base)} {\Z}$.
\end{thm}
\begin{proof}
  Apply \autoref{thm:total-fiber-equiv} to $\encode$, using \autoref{thm:encode-total-equiv}.
\end{proof}

In essence, the two proofs are not very different: the encode-decode one may be seen as a ``reduction'' or ``unpackaging'' of the homotopy-theoretic one.
Each has its advantages; the interplay between the two points of view is part of the interest of the subject.
\index{fundamental!group!of circle|)}

\subsection{The universal cover as an identity system}
\label{sec:pi1s1-idsys}

Note that the fibration $\code:\Sn^1\to\type$ together with $0:\code(\base)$ is a \emph{pointed predicate} in the sense of \autoref{defn:identity-systems}.
From this point of view, we can see that the encode-decode proof in \autoref{subsec:pi1s1-encode-decode} consists of proving that \code satisfies \autoref{thm:identity-systems}\ref{item:identity-systems3}, while the homotopy-theoretic proof in \autoref{subsec:pi1s1-homotopy-theory} consists of proving that it satisfies \autoref{thm:identity-systems}\ref{item:identity-systems4}.
This suggests a third approach.

\begin{thm}
  The pair $(\code,0)$ is an identity system at $\base:\Sn^1$ in the sense of \autoref{defn:identity-systems}.
\end{thm}
\begin{proof}
  Let $D:\prd{x:\Sn^1} \code(x) \to \type$ and $d:D(\base,0)$ be given; we want to define a function $f:\prd{x:\Sn^1}{c:\code(x)} D(x,c)$.
  By circle induction, it suffices to specify $f(\base):\prd{c:\code(\base)} D(\base,c)$ and verify that $\trans{\lloop}{f(\base)} = f(\base)$.

  Of course, $\code(\base)\jdeq \Z$.
  By \autoref{lem:transport-s1-code} and induction on $n$, we may obtain a path $p_n : \transfib{\code}{\lloop^n}{0} = n$ for any integer $n$.
  Therefore, by paths in $\Sigma$-types, we have a path $\pairpath(\lloop^n,p_n) : (\base,0) = (\base,n)$ in $\sm{x:\Sn^1} \code(x)$.
  Transporting $d$ along this path in the fibration $\widehat{D}:(\sm{x:\Sn^1} \code(x)) \to\type$ associated to $D$, we obtain an element of $D(\base,n)$ for any $n:\Z$.
  We define this element to be $f(\base)(n)$:
  \[ f(\base)(n) \defeq \transfib{\widehat{D}}{\pairpath(\lloop^n,p_n)}{d}. \]
  Now we need $\transfib{\lam{x} \prd{c:\code(x)} D(x,c)}{\lloop}{f(\base)} = f(\base)$.
  By \autoref{thm:dpath-forall}, this means we need to show that for any $n:\Z$,
  \begin{narrowmultline*}
    \transfib{\widehat D}{\pairpath(\lloop,\refl{\trans\lloop n})}{f(\base)(n)} 
    =_{D(\base,\trans\lloop n)} \narrowbreak
    f(\base)(\trans\lloop n).
  \end{narrowmultline*}
  Now we have a path $q:\trans\lloop n = n+1$, so transporting along this, it suffices to show
  \begin{multline*}
    \transfib{D(\base)}{q}{\transfib{\widehat D}{\pairpath(\lloop,\refl{\trans\lloop n})}{f(\base)(n)}}\\
    =_{D(\base,n+1)} \transfib{D(\base)}{q}{f(\base)(\trans\lloop n)}.
  \end{multline*}
  By a couple of lemmas about transport and dependent application, this is equivalent to
  \[ \transfib{\widehat D}{\pairpath(\lloop,q)}{f(\base)(n)} =_{D(\base,n+1)} f(\base)(n+1). \]
  However, expanding out the definition of $f(\base)$, we have
  \begin{narrowmultline*}
    \transfib{\widehat D}{\pairpath(\lloop,q)}{f(\base)(n)}
    \narrowbreak
    \begin{aligned}[t]
      &= \transfib{\widehat D}{\pairpath(\lloop,q)}{\transfib{\widehat{D}}{\pairpath(\lloop^n,p_n)}{d}}\\
      &= \transfib{\widehat D}{\pairpath(\lloop^n,p_n) \ct \pairpath(\lloop,q)}{d}\\
      &= \transfib{\widehat D}{\pairpath(\lloop^{n+1},p_{n+1})}{d}\\
      &= f(\base)(n+1).
    \end{aligned}
  \end{narrowmultline*}
  We have used the functoriality of transport, the characterization of composition in $\Sigma$-types (which was an exercise for the reader), and a lemma relating $p_n$ and $q$ to $p_{n+1}$ which we leave it to the reader to state and prove.
  
  This completes the construction of $f:\prd{x:\Sn^1}{c:\code(x)} D(x,c)$.
  Since
  \[f(\base,0) \jdeq \trans{\pairpath(\lloop^0,p_0)}{d} = \trans{\refl{\base}}{d} = d,\]
  we have shown that $(\code,0)$ is an identity system.
\end{proof}

\begin{cor}
  For any $x:\Sn^1$, we have $\eqv{(\base=x)}{\code(x)}$.
\end{cor}
\begin{proof}
  By \autoref{thm:identity-systems}.
\end{proof}

Of course, this proof also contains essentially the same elements as the previous two.
Roughly, we can say that it unifies the proofs of \autoref{thm:pi1s1-decode,lem:s1-encode-decode}, performing the requisite inductive argument only once in a generic case.

\section{Connectedness of suspensions}
\label{sec:conn-susp}

Recall from \autoref{sec:connectivity} that a type $A$ is called \define{$n$-connected} if $\trunc nA$ is contractible.
The aim of this section is to prove that the operation of suspension from \autoref{sec:suspension} increases connectedness.

\begin{thm} \label{thm:suspension-increases-connectedness}
  If $A$ is $n$-connected then the suspension of $A$ is $(n+1)$-connected.
\end{thm}

\begin{proof}
  We remarked in \autoref{sec:colimits} that the suspension of $A$ is the pushout $\unit\sqcup^A\unit$, so we need to
  prove that the following type is contractible:
  \[\trunc{n+1}{\unit\sqcup^A\unit}.\]
  By \autoref{reflectcommutespushout} we know that
  $\trunc{n+1}{\unit\sqcup^A\unit}$ is a pushout in $\typelep{n+1}$ of the diagram
  \[\xymatrix{\trunc{n+1}A \ar[d] \ar[r] & \trunc{n+1}{\unit} \\
    \trunc{n+1}{\unit} & }.\]
  Given that $\trunc{n+1}{\unit}=\unit$, the type
  $\trunc{n+1}{\unit\sqcup^A\unit}$ is also a pushout of the following diagram in
  $\typelep{n+1}$ (because both diagrams are equal)
  \[\Ddiag=\vcenter{\xymatrix{\trunc{n+1}A \ar[d] \ar[r] & \unit \\
    \unit & }}.\]
  We will now prove that $\unit$ is also a pushout of $\Ddiag$ in
  $\typelep{n+1}$.
  Let $E$ be an $(n+1)$-truncated type; we need to prove that the following map
  is an equivalence
  \[\function{(\unit \to E)}{\cocone{\Ddiag}{E}}{y}
  {(y,y,\lamu{u:{\trunc{n+1}A}} \refl{y(\ttt)})}.\]
  where we recall that $\cocone{\Ddiag}{E}$ is the type
  \[\sm{f:\unit \to E}{g:\unit \to E}(\trunc{n+1}A\to
  (f(\ttt)=_E{}g(\ttt))).\]
  The map $\function{(\unit\to E)}{E}{f}{f(\ttt)}$ is an equivalence, hence
  we also have
  \[\cocone{\Ddiag}{E}=\sm{x:E}{y:E}(\trunc{n+1}A\to(x=_Ey)).\]
  Now $A$ is $n$-connected hence so is $\trunc{n+1}A$ because
  $\trunc n{\trunc{n+1}A}=\trunc nA=\unit$, and $(x=_Ey)$ is $n$-truncated because
  $E$ is $(n+1)$-connected. Hence by \autoref{connectedtotruncated} the
  following map is an equivalence
  \[\function{(x=_Ey)}{(\trunc{n+1}A\to(x=_Ey))}{p}{\lam{z} p}\]
  Hence we have
  \[\cocone{\Ddiag}{E}=\sm{x:E}{y:E}(x=_Ey).\]
  But the following map is an equivalence
  \[\function{E}{\sm{x:E}{y:E}(x=_Ey)}{x}{(x,x,\refl{x})}.\]
  Hence
  \[\cocone{\Ddiag}{E}=E.\]
  Finally we get an equivalence
  \[(\unit \to E)\eqvsym\cocone{\Ddiag}{E}\]
  We can now unfold the definitions in order to get the explicit expression of
  this map, and we see easily that this is exactly the map we had at the
  beginning.

  Hence we proved that $\unit$ is a pushout of $\Ddiag$ in $\typelep{n+1}$. Using
  uniqueness of pushouts we get that $\trunc{n+1}{\unit\sqcup^A\unit}=\unit$
  which proves that the suspension of $A$ is $(n+1)$-connected.
\end{proof}

\begin{cor} \label{cor:sn-connected}
  For all $n:\N$, the sphere $\Sn^n$ is $(n-1)$-connected.
\end{cor}

\begin{proof}
  We prove this by induction on $n$.
  For $n=0$ we have to prove that $\Sn^0$ is merely inhabited, which is clear.
  Let $n:\N$ be such that $\Sn^n$ is $(n-1)$-connected. By definition $\Sn^{n+1}$
  is the suspension of $\Sn^n$, hence by the previous lemma $\Sn^{n+1}$ is
  $n$-connected.
\end{proof}

\section{\texorpdfstring{$\pi_{k \le n}$}{π\_(k≤n)} of an \texorpdfstring{$n$}{n}-connected space and \texorpdfstring{$\pi_{k < n}(\Sn ^n)$}{π\_(k<n)(Sⁿ)}}
\label{sec:pik-le-n}

Let $(A,a)$ be a pointed type and $n:\N$.  Recall from
\autoref{thm:homotopy-groups} that if $n>0$ the set $\pi_n(A,a)$ has a group
structure, and if $n>1$ the group is abelian\index{group!abelian}.

We can now say something about homotopy groups of $n$-truncated and
$n$-connected types.

\begin{lem}
  If $A$ is $n$-truncated and $a:A$, then $\pi_k(A,a)=\unit$ for all $k>n$.
\end{lem}

\begin{proof}
  The loop space of an $n$-type  is an
  $(n-1)$-type, hence $\Omega^k(A,a)$ is an $(n-k)$-type, and we have
  $(n-k)\le-1$ so $\Omega^k(A,a)$ is a mere proposition. But $\Omega^k(A,a)$ is inhabited,
  so it is actually contractible and
  $\pi_k(A,a)=\trunc0{\Omega^k(A,a)}=\trunc0{\unit}=\unit$.
\end{proof}

\begin{lem} \label{lem:pik-nconnected}
  If $A$ is $n$-connected and $a:A$, then $\pi_k(A,a)=\unit$ for all $k\le{}n$.
\end{lem}

\begin{proof}
  We have the following sequence of equalities:
  \begin{narrowmultline*}
    \pi_k(A,a) = \trunc0{\Omega^k(A,a)}
    = \Omega^k(\trunc k{(A,a)})
    = \Omega^k(\trunc k{\trunc n{(A,a)}})
    = \narrowbreak
      \Omega^k(\trunc k{\unit})
    = \Omega^k(\unit)
    = \unit.
  \end{narrowmultline*}
  The third equality uses the fact that $k\le{}n$ in order to use that
  $\truncf k\circ\truncf n=\truncf k$ and the fourth equality uses the fact that $A$ is
  $n$-connected.
\end{proof}

\begin{cor}
  $\pi_k(\Sn ^n) = \unit$ for $k < n$.
\end{cor}
\begin{proof}
  The sphere $\Sn^n$ is $(n-1)$-connected by \autoref{cor:sn-connected}, so
  we can apply \autoref{lem:pik-nconnected}.
\end{proof}

\section{Fiber sequences and the long exact sequence}
\label{sec:long-exact-sequence-homotopy-groups}

\index{fiber sequence|(}%
\index{sequence!fiber|(}%

If the codomain of a function $f:X\to Y$ is equipped with a basepoint $y_0:Y$, then we refer to the fiber $F\defeq \hfib f {y_0}$ of $f$ over $y_0$ as \define{the fiber of $f$}\index{fiber}.
(If $Y$ is connected, then $F$ is determined up to mere equivalence; see \autoref{ex:unique-fiber}.)
We now show that if $X$ is also pointed and $f$ preserves basepoints, then there is a relation between the homotopy groups of $F$, $X$, and $Y$ in the form of a \emph{long exact sequence}.
We derive this by way of the \emph{fiber sequence} associated to such an $f$.

\begin{defn}\label{def:pointedmap}
  A \define{pointed map}
  \indexdef{pointed!map}%
  \indexsee{function!pointed}{pointed map}%
  between pointed types $(X,x_0)$ and $(Y,y_0)$ is a
  map $f:X\to Y$ together with a path $f_0:f(x_0)=y_0$.
\end{defn}

For any pointed types $(X,x_0)$ and $(Y,y_0)$, there is a pointed map $(\lam{x} y_0) : X\to Y$ which is constant at the basepoint.
We call this the \define{zero map}\indexdef{zero!map}\indexdef{function!zero} and sometimes write it as $0:X\to Y$.

Recall that every pointed type $(X,x_0)$ has a loop space\index{loop space} $\Omega (X,x_0)$.
We now note that this operation is functorial on pointed maps.\index{loop space!functoriality of}\index{functor!loop space}

\begin{defn}
  Given a pointed map between pointed types $f:X \to Y$, we define a pointed
  map $\Omega f:\Omega X
  \to \Omega Y$ by
  \[(\Omega f)(p) \defeq \rev{f_0}\ct\ap{f}{p}\ct f_0.\]
  The path $(\Omega f)_0 : (\Omega f) (\refl{x_0}) = \refl{y_0}$, which exhibits $\Omega f$ as a pointed map, is the obvious path of type
  \[\rev{f_0}\ct\ap{f}{\refl{x_0}}\ct f_0=\refl{y_0}.\]
\end{defn}

There is another functor on pointed maps, which takes $f:X\to Y$ to $\proj1 : \hfib f {y_0} \to X$.
When $f$ is pointed, we always consider $\hfib f {y_0}$ to be pointed with basepoint $(x_0,f_0)$, in which case $\proj1$ is also a pointed map, with witness $(\proj1)_0 \defeq \refl{x_0}$.
Thus, this operation can be iterated.

\begin{defn}
  The \define{fiber sequence}
  \indexdef{fiber sequence}%
  \indexdef{sequence!fiber}%
  of a pointed map $f:X\to Y$ is the infinite sequence of pointed types and pointed maps
  \[\xymatrix{\dots \ar[r]^{f^{(n+1)}} & X^{(n+1)} \ar[r]^{f^{(n)}} & X^{(n)} \ar^-{f^{(n-1)}}[r] & \dots \ar[r] & X^{(2)} \ar^-{f^{(1)}}[r] & X^{(1)} \ar[r]^{f^{(0)}} & X^{(0)}}\]
  defined recursively by
  \[ X^{(0)} \defeq Y\qquad
    X^{(1)} \defeq X\qquad
    f^{(0)} \defeq f\qquad
    \]    
  and
  \begin{alignat*}{2}
    X^{(n+1)} &\defeq \hfib {f^{(n-1)}}{x^{(n-1)}_0}\\
    f^{(n)} &\defeq \proj1 &: X^{(n+1)} \to X^{(n)}.
  \end{alignat*}
  where $x^{(n)}_0$ denotes the basepoint of $X^{(n)}$, chosen recursively as above.
\end{defn}

Thus, any adjacent pair of maps in this fiber sequence is of the form
\[ \xymatrix{ X^{(n+1)} \jdeq \hfib{f^{(n-1)}}{x^{(n-1)}_0} \ar[rr]^-{f^{(n)}\jdeq \proj1} && X^{(n)} \ar[r]^{f^{(n-1)}} & X^{(n-1)}. } \]
In particular, we have $f^{(n-1)} \circ f^{(n)} = 0$.
We now observe that the types occurring in this sequence are the iterated loop spaces\index{loop space!iterated} of the base
space $Y$, the total space\index{total!space} $X$, and the fiber $F\defeq \hfib f {y_0}$, and similarly for the maps.

\begin{lem}\label{thm:fiber-of-the-fiber}
  Let $f:X\to Y$ be a pointed map of pointed spaces.  Then:
  \begin{enumerate}
  \item The fiber of $f^{(1)}\defeq \proj1 : \hfib f {y_0} \to X$ is equivalent to $\Omega Y$.\label{item:fibseq1}
  \item Similarly, the fiber of $f^{(2)} : \Omega Y \to \hfib f {y_0}$ is equivalent to $\Omega X$.\label{item:fibseq2}
  \item Under these equivalences, the map $f^{(3)} : \Omega X\to \Omega Y$ is identified with $\Omega f\circ\rev{(\blank)}$.\label{item:fibseq3}
  \end{enumerate}
\end{lem}
\begin{proof}
  For~\ref{item:fibseq1}, we have
  \begin{align*}
    \hfib{f^{(1)}}{x_0}
    &\defeq \sm{z:\hfib{f}{y_0}} (\proj{1}(z) = x_0)\\
    &\eqvsym \sm{x:A}{p:f(x)=y_0} (x = x_0) &\text{(by \autoref{ex:sigma-assoc})}\\
    &\eqvsym (f(x_0) = y_0) &\text{(as $\tsm{x:A} (x=x_0)$ is contractible)}\\
    &\eqvsym (y_0 = y_0) &\text{(by $(f_0 \ct \blank)$)}\\
    &\jdeq \Omega Y.
  \end{align*}
  Tracing through, we see that this equivalence sends $((x,p),q)$ to $\opp{f_0} \ct \ap{f}{\opp q} \ct p$, while its inverse sends $r:y_0=y_0$ to $((x_0, f_0 \ct r), \refl{x_0})$.
  In particular, the basepoint $((x_0,f_0),\refl{x_0})$ of $\hfib{f^{(1)}}{x_0}$ is sent to $\opp{f_0} \ct \ap{f}{\opp {\refl{x_0}}} \ct f_0$, which equals $\refl{y_0}$.
  Hence this equivalence is a pointed map (see \autoref{ex:pointed-equivalences}).
  Moreover, under this equivalence, $f^{(2)}$ is identified with $\lam{r} (x_0, f_0 \ct r): \Omega Y \to \hfib f {y_0}$.

  \cref{item:fibseq2} follows immediately by applying~\ref{item:fibseq1} to $f^{(1)}$ in place of $f$.
  Since $(f^{(1)})_0 \defeq \refl{x_0}$, under this equivalence $f^{(3)}$ is identified with the map $\Omega X \to \hfib {f^{(1)}}{x_0}$ defined by $s \mapsto ((x_0,f_0),s)$.
  Thus, when we compose with the previous equivalence $\hfib {f^{(1)}}{x_0} \eqvsym \Omega Y$, we see that $s$ maps to $\opp{f_0} \ct \ap{f}{\opp s} \ct f_0$, which is by definition $(\Omega f)(\opp s)$, giving~\ref{item:fibseq3}.
\end{proof}

Thus, the fiber sequence of $f:X\to Y$ can be pictured as:
\[\xymatrix@C=1.5pc{
  \dots \ar[r] &
  \Omega^2 X \ar[r]^-{\Omega^2 f} &
  \Omega^2 Y \ar[r]^-{-\Omega \partial} &
  \Omega F \ar[r]^-{-\Omega i} &
  \Omega X \ar[r]^-{-\Omega f} &
  \Omega Y \ar[r]^-{\partial} &
  F \ar[r]^-{i} &
  X \ar[r]^{f} & Y.}\]
where the minus signs denote composition with path inversion $\rev{(\blank)}$.
Note that by \autoref{ex:ap-path-inversion}, we have
\[ \Omega\left(\Omega f\circ \opp{(\blank)}\right) \circ \opp{(\blank)}
= \Omega^2 f \circ \opp{(\blank)} \circ \opp{(\blank)}
= \Omega^2 f.
\]
Thus, there are minus signs on the $k$-fold loop maps whenever $k$ is odd.

From this fiber sequence we will deduce an \emph{exact sequence of pointed sets}.
\index{image}%
Let $A$ and $B$ be sets and $f:A\to B$ a function, and recall from \autoref{defn:modal-image} the definition of the \emph{image} $\im(f)$, which can be regarded as a subset of $B$:
\[\im(f) \defeq \setof{b:B | \exis{a:A} f(a)=b}. \]
If $A$ and $B$ are moreover pointed with basepoints $a_0$ and $b_0$, and $f$ is a pointed map, we define the \define{kernel}
\indexdef{kernel}%
\indexdef{pointed!map!kernel of}%
of $f$ to be the following subset of $A$:
\symlabel{kernel}
\[\ker(f) \defeq \setof{x:A | f(x) = b_0}. \]
Of course, this is just the fiber of $f$ over the basepoint $b_0$; it a subset of $A$ because $B$ is a set.

Note that any group is a pointed set, with its unit element as basepoint, and any group homomorphism is a pointed map.
In this case, the kernel and image agree with the usual notions from group theory.

\begin{defn}
  An \define{exact sequence of pointed sets}
  \indexdef{exact sequence}%
  \indexdef{sequence!exact}%
  is a (possibly bounded) sequence of pointed sets and pointed maps:
  \[\xymatrix{\dots \ar[r] & A^{(n+1)} \ar[r]^-{f^{(n)}} & A^{(n)} \ar[r]^{f^{(n-1)}} & A^{(n-1)} \ar[r] &
    \dots}\]
  such that for every $n$, the image of $f^{(n)}$ is equal, as a subset of $A^{(n)}$, to the kernel of $f^{(n-1)}$.
  In other words, for all $a:A^{(n)}$ we have
  \[ (f^{(n-1)}(a) = a^{(n-1)}_0) \iff \exis{b:A^{(n+1)}} (f^{(n)}(b)=a). \]
  where $a^{(n)}_0$ denotes the basepoint of $A^{(n)}$.
\end{defn}

Usually, most or all of the pointed sets in an exact sequence are groups, and often abelian groups.
When we speak of an \define{exact sequence of groups}, it is assumed moreover that the maps are group homomorphisms and not just pointed maps.

\index{exact sequence}
\begin{thm}\label{thm:les}
  Let $f:X \to Y$ be a pointed map between pointed spaces with fiber $F\defeq \hfib f {y_0}$.
  Then we have the following long exact sequence, which consists of groups except for the last three terms, and abelian groups except for the last six.

  \[
  \xymatrix@R=1.2pc@C=3pc{
    \vdots & \vdots & \vdots \ar[lld] \\
    \pi_k(F) \ar[r] & \pi_k(X) \ar[r] & \pi_k(Y) \ar[lld] \\
    \vdots & \vdots & \vdots \ar[lld] \\
    \pi_2(F) \ar[r] & \pi_2(X) \ar[r] & \pi_2(Y) \ar[lld] \\
    \pi_1(F) \ar[r] & \pi_1(X) \ar[r] & \pi_1(Y) \ar[lld]\\
    \pi_0(F) \ar[r] & \pi_0(X) \ar[r] & \pi_0(Y)}
  \]
\end{thm}

\vspace*{0pt plus 10ex}
\goodbreak

\begin{proof}
  We begin by showing that the 0-truncation of a fiber sequence is an exact sequence of pointed sets.
  Thus, we need to show that for any adjacent pair of maps in a fiber sequence:
  \[\xymatrix{\hfib{f}{z_0} \ar^-g[r] & W \ar^-f[r] & Z}\]
  with $g\defeq \proj1$, the sequence
  \[\xymatrix{\trunc{0}{\hfib{f}{z_0}} \ar^-{\trunc0g}[r] & \trunc0{W}
    \ar^-{\trunc0{f}}[r] & \trunc0{Z}}\]
  is exact, i.e.\ that $\im(\trunc0g)\subseteq\ker(\trunc0f)$ and $\ker(\trunc0f)\subseteq\im(\trunc0g)$.

  The first inclusion is equivalent to $\trunc0g\circ\trunc0f=0$, which holds by functoriality of $\truncf0$ and the fact that $g\circ f=0$.
  For the second, we assume $w':\trunc0W$ and $p':\trunc0f(w')=\tproj0{z_0}$ and show there merely exists ${t:\hfib{f}{z_0}}$ such that $g(t)=w'$.
  Since our goal is a mere proposition, we can assume that $w'$ is of the
  form $\tproj0w$ for some $w:W$.
  Now by \autoref{thm:path-truncation}, $p' : \tproj0{f(w)}=\tproj0{z_0}$ yields $p'': \trunc{-1}{f(w)=z_0}$, so by a further truncation induction we may assume some $p:f(w)=z_0$.
  But now we have $\tproj0{(w,p)}:\tproj0{\hfib{f}{z_0}}$ whose image under $\trunc0g$ is $\tproj0w\jdeq w'$, as desired.

  Thus, applying $\truncf0$ to the fiber sequence of $f$, we obtain a long exact sequence involving the pointed sets $\pi_k(F)$, $\pi_k(X)$, and $\pi_k(Y)$ in the desired order.
  And of course, $\pi_k$ is a group for $k\ge1$, being the 0-truncation of a loop space, and an abelian group for $k\ge 2$ by the Eckmann--Hilton argument
  \index{Eckmann--Hilton argument} (\autoref{thm:EckmannHilton}).
  Moreover, \autoref{thm:fiber-of-the-fiber} allows us to identify the maps $\pi_k(F) \to \pi_k(X)$ and $\pi_k(X) \to \pi_k(Y)$ in this exact sequence as $(-1)^k \pi_k(i)$ and $(-1)^k \pi_k(f)$ respectively.

  More generally, every map in this long exact sequence except the last three is of the form $\trunc0{\Omega h}$ or $\trunc0{-\Omega h}$ for some $h$.
  In the former case it is a group homomorphism, while in the latter case it is a homomorphism if the groups are abelian; otherwise it is an ``anti-homomorphism''.
  However, the kernel and image of a group homomorphism are unchanged when we replace it by its negative, and hence so is the exactness of any sequence involving it.
  Thus, we can modify our long exact sequence to obtain one involving $\pi_k(i)$ and $\pi_k(f)$ directly and in which all the maps are group homomorphisms (except the last three).
\end{proof}

The usual properties of exact sequences of abelian groups\index{group!abelian!exact sequence of} can be proved as
usual. In particular we have:
\begin{lem}\label{thm:ses}
  Suppose given an exact sequence of abelian groups:
  \[\xymatrix{K \ar[r]& G \ar^f[r] & H \ar[r] & Q.}\]
  \begin{enumerate}
  \item If $K=0$, then $f$ is injective.\label{item:sesinj}
  \item If $Q=0$, then $f$ is surjective.\label{item:sessurj}
  \item If $K=Q=0$, then $f$ is an isomorphism.\label{item:sesiso}
  \end{enumerate}
\end{lem}
\begin{proof}
  Since the kernel of $f$ is the image of $K\to G$, if $K=0$ then the kernel of $f$ is $\{0\}$;
  hence $f$ is injective because it's a group morphism.
  Similarly, since the image of $f$ is the kernel of $H\to Q$, if $Q=0$ then the image of $f$ is all of $H$, so $f$ is surjective.
  Finally,~\ref{item:sesiso} follows from~\ref{item:sesinj} and~\ref{item:sessurj} by \autoref{thm:mono-surj-equiv}.
\end{proof}

As an immediate application, we can now quantify in what way $n$-connectedness of a map is stronger than inducing an equivalence on $n$-truncations.

\begin{cor}\label{thm:conn-pik}
  Let $f:A\to B$ be $n$-connected and $a:A$, and define $b\defeq f(a)$.  Then:
  \begin{enumerate}
  \item If $k\le n$, then $\pi_k(f):\pi_k(A,a) \to \pi_k(B,b)$ is an isomorphism.
  \item If $k=n+1$, then $\pi_k(f):\pi_k(A,a) \to \pi_k(B,b)$ is surjective.
  \end{enumerate}
\end{cor}
\begin{proof}
  As part of the long exact sequence, for each $k$ we have an exact sequence
  \[\xymatrix{\pi_k(\hfib f b) \ar[r]& \pi_k(A,a) \ar^f[r] & \pi_k(B,b) \ar[r] & \pi_{k-1}(\hfib f b).}\]
  Now since $f$ is $n$-connected, $\trunc n{\hfib f b}$ is contractible.
  Therefore, if $k\le n$, then $\pi_k(\hfib f b) = \trunc0{\Omega^k(\hfib f b)} = \Omega^k(\trunc k{\hfib f b})$ is also contractible.
  Thus, $\pi_k(f)$ is an isomorphism for $k\le n$ by \autoref{thm:ses}\ref{item:sesiso}, while for $k=n+1$ it is surjective by \autoref{thm:ses}\ref{item:sessurj}.
\end{proof}

In \autoref{sec:whitehead} we will see that the converse of \autoref{thm:conn-pik} also holds.

\index{fiber sequence|)}%
\index{sequence!fiber|)}%

\section{The Hopf fibration}
\label{sec:hopf}

In this section we will define the \define{Hopf fibration}.
\indexdef{Hopf!fibration}%

\begin{thm}[Hopf Fibration]\label{thm:hopf-fibration}
There is a fibration $H$ over $\Sn ^2$ whose fiber over the basepoint is $\Sn ^1$ and
whose total space is $\Sn ^3$.
\end{thm}

The Hopf fibration will allow us to compute several homotopy groups of
spheres.
Indeed, it yields the following long exact sequence
\index{homotopy!group!of sphere}
\index{sequence!exact}
of homotopy groups
(see
\autoref{sec:long-exact-sequence-homotopy-groups}):
\[
\xymatrix@R=1.2pc{
  \pi_k(\Sn^1) \ar[r] & \pi_k(\Sn^3) \ar[r] & \pi_k(\Sn^2) \ar[lld] \\
  \vdots & \vdots & \vdots \ar[lld] \\
  \pi_2(\Sn^1) \ar[r] & \pi_2(\Sn^3) \ar[r] & \pi_2(\Sn^2) \ar[lld] \\
  \pi_1(\Sn^1) \ar[r] & \pi_1(\Sn^3) \ar[r] & \pi_1(\Sn^2)}
\]
We've already computed all $\pi_n(\Sn^1)$, and $\pi_k(\Sn^n)$ for $k<n$, so this
becomes the following:
\[
\xymatrix@R=1.2pc{
  0 \ar[r] & \pi_k(\Sn^3) \ar[r] & \pi_k(\Sn^2) \ar[lld] \\
  \vdots & \vdots & \vdots \ar[lld] \\
  0 \ar[r] & \pi_3(\Sn^3) \ar[r] & \pi_3(\Sn^2) \ar[lld] \\
  0 \ar[r] & 0 \ar[r] & \pi_2(\Sn^2) \ar[lld] \\
  \Z \ar[r] & 0 \ar[r] & 0}
\]
In particular we get the following result:

\begin{cor} \label{cor:pis2-hopf}
  We have $\eqv{\pi_2(\Sn^2)}{\Z}$ and $\eqv{\pi_k(\Sn^3)}{\pi_k(\Sn^2)}$ for
  every $k\ge3$ (where the map is induced by the Hopf fibration, seen as a map
  from the total space $\Sn^3$ to the base space $\Sn^2$).
\end{cor}

In fact, we can say more: the fiber sequence of the Hopf fibration will show that $\Omega^3(\Sn^3)$ is the fiber of a map from $\Omega^3(\Sn^2)$ to $\Omega^2(\Sn^1)$.
Since $\Omega^2(\Sn^1)$ is contractible, we have $\eqv{\Omega^3(\Sn^3)}{\Omega^3(\Sn^2)}$.
In classical homotopy theory, this fact would be a consequence of \autoref{cor:pis2-hopf} and Whitehead's theorem, but Whitehead's theorem is not necessarily valid in homotopy type theory (see \autoref{sec:whitehead}).
We will not use the more precise version here though.

\subsection{Fibrations over pushouts}
\label{sec:fib-over-pushout}

We first start with a lemma explaining how to construct fibrations over
pushouts.
\index{pushout}%

\begin{lem}\label{lem:fibration-over-pushout}
  Let $\Ddiag=(Y\xleftarrow{j}X\xrightarrow{k}Z)$ be a span\index{span} and assume
  that we have
  \begin{itemize}
  \item Two fibrations $E_Y:Y\to\type$ and $E_Z:Z\to\type$.
  \item An equivalence $e_X$ between $E_Y\circ j:X\to\type$ and $E_Z\circ
    k:X\to\type$, i.e.
    \[e_X:\prd{x:X}\eqv{E_Y(j(x))}{E_Z(k(x))}.\]
  \end{itemize}

  Then we can construct a fibration $E:Y\sqcup^XZ\to\type$ such that
  \begin{itemize}
  \item For all $y:Y$, $E(\inl(y))\judgeq E_Y(y)$.
  \item For all $z:Z$, $E(\inr(z))\judgeq E_Z(z)$.
  \item For all $x:X$, $\map E{\glue(x)}=\ua(e_X(x))$ (note that both sides of
    the equation are paths in $\type$ from $E_Y(j(x))$ to $E_Z(k(x))$).
  \end{itemize}
  Moreover, the total space of this fibration fits in the following pushout
  square:
  \[\xymatrix{ \sm{x:X}E_Y(j(x)) \ar[r]_\sim^{\idfunc\times e_X}
    \ar[d]_{j\times\idfunc} &
    \sm{x:X}E_Z(k(x)) \ar[r]^-{k\times\idfunc}
    & \sm{z:Z}E_Z(z) \ar[d]^\inr \\
    \sm{y:Y}E_Y(y) \ar[rr]_\inl & & \sm{t:Y\sqcup^XZ}E(t) }\]
\end{lem}

\begin{proof}
  We define $E$ by the recursion principle of the pushout $Y\sqcup^XZ$. For
  that, we need to specify the value of $E$ on elements of the form $\inl(y)$, $\inr(z)$
  and the action of $E$ on paths $\glue(x)$, so we can just choose the following
  values:
  \begin{align*}
    E(\inl(y)) & \defeq E_Y(y),\\
    E(\inr(z)) & \defeq E_Z(z),\\
    \map E{\glue(x)} & \defid \ua(e_X(x)).
  \end{align*}
  To see that the total space of this fibration is a pushout, we apply the
  flattening lemma (\autoref{thm:flattening}) with the following values:\index{flattening lemma}
  \begin{itemize}
  \item $A\defeq Y+Z$, $B\defeq X$ and $f,g:B\to A$ are defined by
    $f(x)\defeq\inl(j(x))$, $g(x)\defeq\inr(k(x))$,
  \item the type family $C:A\to\type$ is defined by
    \begin{equation*}
      C(\inl(y)) \defeq E_Y(y)
      \qquad\text{and}\qquad
      C(\inr(z)) \defeq E_Z(z),
    \end{equation*}
  \item the family of equivalences $D:\prd{b:B}C(f(b))\eqvsym C(g(b))$ is defined
    to be $e_X$.
  \end{itemize}
  The base higher inductive type $W$ in the flattening lemma is equivalent to
  the pushout $Y\sqcup^XZ$ and the type family $P:Y\sqcup^XZ\to\type$ is
  equivalent to the $E$ defined above.

  Thus the flattening lemma tells us that $\sm{t:Y\sqcup^XZ}E(t)$ is equivalent
  the higher inductive type ${E^{\mathrm{tot}}}'$ with the following generators:
  \begin{itemize}
  \item a function $\mathsf{z}:\sm{a:Y+Z}C(a)\to {E^{\mathrm{tot}}}'$,
  \item for each $x:X$ and $t:E_Y(j(x))$, a path
    \narrowequation{\mathsf{z}(\inl(j(x)),t)=\mathsf{z}(\inr(k(x)),e_C(t)).}
  \end{itemize}
  Using the flattening lemma again or a direct computation, it is easy to see
  that $\eqv{\sm{a:Y+Z}C(a)}{\sm{y:Y}E_Y(y)+\sm{z:Z}E_Z(z)}$, hence
  ${E^{\mathrm{tot}}}'$ is equivalent to the higher inductive type
  $E^{\mathrm{tot}}$ with the following generators:
  \begin{itemize}
  \item a function $\inl:\sm{y:Y}E_Y(y)\to E^{\mathrm{tot}}$,
  \item a function $\inr:\sm{z:Z}E_Z(z)\to E^{\mathrm{tot}}$,
  \item for each $(x,t):\sm{x:X}E_Y(j(x))$ a path
    \narrowequation{\glue(x,t):\inl(j(x),t) = \inr(k(x),e_X(t)).}
  \end{itemize}
  Thus the total space of $E$ is the pushout of the total spaces of
  $E_Y$ and $E_Z$, as required.
\end{proof}

\subsection{The Hopf construction}

\begin{defn}
  An \define{H-space}
  \indexdef{H-space}%
  consists of
  \begin{itemize}
  \item a type $A$,
  \item a base point $e:A$,
  \item a binary operation $\mu:A\times A\to A$, and
  \item for every $a:A$, equalities $\mu(e,a)=a$ and $\mu(a,e)=a$.
  \end{itemize}
\end{defn}

\begin{lem}
  Let $A$ be a connected H-space. Then for every $a:A$, the maps $\mu(a,\blank):A\to
  A$ and $\mu(\blank,a):A\to A$ are equivalences.
\end{lem}

\begin{proof}
  Let us prove that for every $a:A$ the map $\mu(a,\blank)$ is an equivalence. The
  other statement is symmetric.
  The statement that $\mu(a,\blank)$ is an equivalence corresponds to a type family
  $P:A\to\prop$ and proving it corresponds to finding a section of this type
  family.

  The type $\prop$ is a set (\autoref{thm:hleveln-of-hlevelSn}) hence we can
  define a new type family $P':\trunc0A\to\prop$ by $P'(\tproj0a)\defeq
  P(a)$. But $A$ is connected by assumption, hence $\trunc0A$ is
  contractible. This implies that in order to find a section of $P'$, it is
  enough to find a point in the fiber of $P'$ over $\tproj0e$. But we have
  $P'(\tproj0e)=P(e)$ which is inhabited because $\mu(e,\blank)$ is equal to the
  identity map by definition of an H-space, hence is an equivalence.

  We have proved that for every $x:\trunc0A$ the proposition $P'(x)$ is true,
  hence in particular for every $a:A$ the proposition $P(a)$ is true because
  $P(a)$ is $P'(\tproj0a)$.
\end{proof}

\begin{defn}
  Let $A$ be a connected H-space. We define a fibration over $\susp A$ using
  \autoref{lem:fibration-over-pushout}.

  Given that $\susp A$ is the pushout $\unit\sqcup^A\unit$, we can define a
  fibration over $\susp A$ by specifying
  \begin{itemize}
  \item two fibrations over $\unit$ (i.e. two types $F_1$ and $F_2$), and
  \item a family $e:A\to(\eqv{F_1}{F_2})$ of equivalences between
    $F_1$ and $F_2$, one for every element of $A$.
  \end{itemize}
  We take $A$ for $F_1$ and $F_2$, and for $a:A$ we take the equivalence
  $\mu(a,\blank)$ for $e(a)$.
\end{defn}

According to \autoref{lem:fibration-over-pushout}, we have the following
diagram:
\[\xymatrix{A \ar@{->>}[d] & A \times A \ar[l]_-{\proj2} \ar@{->>}_{\proj1}[d]
  \ar[r]^-{\mu} & A \ar@{->>}[d] \\
  1 & A \ar[r] \ar[l] & 1}\]
and the fibration we just constructed is a fibration over $\susp A$ whose total
space is the pushout of the top line.

Moreover, with $f(x,y)\defeq(\mu(x,y),y)$ we have the following diagram:
\[\xymatrix{A \ar_\idfunc[d] & A \times A \ar[l]_-{\proj2} \ar^f[d]
  \ar[r]^-{\mu} & A \ar^\idfunc[d] \\
  A & A\times A \ar^-{\proj2}[l] \ar_-{\proj1}[r] & A}\]
The diagram commutes and the three vertical maps are equivalences, the inverse
of $f$ being the function $g$ defined by
\[g(u,v)\defeq(\opp{\mu(\blank,v)}(u),v).\]
This shows that the two lines are equivalent (hence equal) spans, so the total
space of the fibration we constructed is equivalent to the pushout of the bottom
line.
And by definition, this latter pushout is the \emph{join} of $A$ with itself (see \autoref{sec:colimits}).
We have proven:

\begin{lem}\label{lem:hopf-construction}
  Given a connected H-space $A$, there is a fibration, called the
  \define{Hopf construction},
  \indexdef{Hopf!construction}%
  over $\susp A$ with fiber $A$ and total space $A*A$.
\end{lem}

\subsection{The Hopf fibration}

\index{Hopf fibration|(}
\indexsee{fibration!Hopf}{Hopf fibration}
We will first construct a structure of H-space on the circle $\Sn^1$, hence by
\autoref{lem:hopf-construction} we will get a fibration over $\Sn^2$ with fiber
$\Sn^1$ and total space $\Sn^1*\Sn^1$. We will then prove that this join is
equivalent to $\Sn^3$.

\begin{lem}\label{lem:hspace-S1}
  There is an H-space structure on the circle $\Sn^1$.
\end{lem}
\begin{proof}
  For the base point of the H-space structure we choose $\base$.
  Now we need to define the multiplication operation
  $\mu:\Sn^1\times\Sn^1\to\Sn^1$.
  We will define the curried form $\widetilde\mu:\Sn^1\to(\Sn^1\to\Sn^1)$ of $\mu$
  by recursion on $\Sn^1$:
  \begin{equation*}
    \widetilde\mu(\base) \defeq\idfunc[\Sn^1],
    \qquad\text{and}\qquad
    \ap{\widetilde\mu}{\lloop} \defid\funext(h).
  \end{equation*}
  where $h:\prd{x:\Sn^1}(x=x)$ is the function defined in \autoref{thm:S1-autohtpy},
  which has the property that $h(\base) \defeq\lloop$.

  Now we just have to prove that $\mu(x,\base)=\mu(\base,x)=x$ for every
  $x:\Sn^1$.
  By definition, if $x:\Sn^1$ we have
  $\mu(\base,x)=\widetilde\mu(\base)(x)=\idfunc[\Sn^1](x)=x$. For the equality
  $\mu(x,\base)=x$ we do it by induction on $x:\Sn^1$:
  \begin{itemize}
  \item If $x$ is $\base$ then $\mu(\base,\base)=\base$ by definition, so we
    have $\refl\base:\mu(\base,\base)=\base$.
  \item When $x$ varies along $\lloop$, we need to prove that
    \[\refl\base\ct\apfunc{\lam{x}x}({\lloop})=
    \apfunc{\lam{x}\mu(x,\base)}({\lloop})\ct\refl\base.\]
    The left-hand side is equal to $\lloop$, and for the right-hand side we have:
    \begin{align*}
      \apfunc{\lam{x}\mu(x,\base)}({\lloop})\ct\refl\base &=
      \apfunc{\lam{x}(\widetilde\mu(x))(\base)}({\lloop})\\
      &=\happly(\apfunc{\lam{x}(\widetilde\mu(x))}({\lloop}),\base)\\
      &=\happly(\funext(h),\base)\\
      &=h(\base)\\
      &=\lloop. \qedhere
    \end{align*}
  \end{itemize}
\end{proof}

Now recall from \autoref{sec:colimits} that the \emph{join} $A*B$ of types $A$ and $B$ is the pushout of the diagram
\index{join!of types}%
\[A \xleftarrow{\proj1}A\times B \xrightarrow{\proj2} B. \]

\begin{lem}
  \index{associativity!of join}%
  The operation of join is associative: if $A$, $B$ and $C$ are three types
  then we have an equivalence $\eqv{(A*B)*C}{A*(B*C)}$.
\end{lem}

\begin{proof}
  We define a map $f:(A*B)*C\to A*(B*C)$ by induction. We first need to define
  $f\circ\inl:A*B\to A*(B*C)$ which will be done by induction, then
  $f\circ\inr:C\to A*(B*C)$, and then $\apfunc{f}\circ\glue:\prd{t:(A*B)\times
    C}f(\inl(\fst(t)))=f(\inr(\snd(t)))$ which will be done by induction on the
  first component of $t$:
  \begin{align*}
    (f\circ\inl)(\inl(a))) &\defeq \inl(a), \\
    (f\circ\inl)(\inr(b))) &\defeq \inr(\inl(b)), \\
    \apfunc{f\circ\inl}(\glue(a,b)) &\defid \glue(a,\inl(b)), \\
    f(\inr(c)) &\defeq \inr(\inr(c)),\\
    \apfunc{f}(\glue(\inl(a),c)) &\defid\glue(a,\inr(c)),\\
    \apfunc{f}(\glue(\inr(b),c)) &\defid\apfunc{\inr}(\glue(b,c)),\\
    \apdfunc{\lam{x}\apfunc{f}(\glue(x,c))}(\glue(a,b)) &\defid
    ``\apdfunc{\lam{x}\glue(a,x)}(\glue(b,c))''.
  \end{align*}
  For the last equation, note that the right-hand side is of type
  \[\transfib{\lam{x}\inl(a)=\inr(x)}{\glue(b,c)}{\glue(a,\inl(b))}=
  \glue(a,\inr(c))\]
  whereas it is supposed to be of type
  \begin{narrowmultline*}
    \transfib{\lam{x}f(\inl(x))=f(\inr(c))}{\glue(a,b)}{\apfunc{f}(\glue(\inl(a),c))}
    = \narrowbreak
    \apfunc{f}(\glue(\inr(b),c)).
  \end{narrowmultline*}
  But by the previous clauses in the definition, both of these types are equivalent to the following type:
  \[\glue(a,\inr(c))=\glue(a,\inl(b))\ct\apfunc\inr(\glue(b,c)),\]
  and so we can coerce by an equivalence to obtain the necessary element.
  Similarly, we can define a map $g:A*(B*C)\to(A*B)*C$, and checking that $f$ and
  $g$ are inverse to each other is a long and tedious but essentially
  straightforward computation.
\end{proof}

A more conceptual proof sketch is as follows.

\begin{proof}
  Let us consider the following diagram where the maps are the obvious
  projections:
  \[\xymatrix{
    A & A\times C \ar[l] \ar[r] & A\times C\\
    A\times B \ar[u] \ar[d] & A\times B\times C \ar[l]\ar[u]\ar[r]\ar[d] &
    A\times C \ar[u] \ar[d] \\
    B & B\times C \ar[l] \ar[r] & C}\]
  Taking the colimit of the columns gives the following
  diagram, whose colimit is $(A*B)*C$:
  \[\xymatrix{A*B & (A*B)\times C \ar[l]\ar[r] & C}\]
  On the other hand, taking the colimit of the lines gives a diagram whose
  colimit is $A*(B*C)$.

  Hence using a Fubini-like theorem for colimits (that we haven’t proved) we
  have an equivalence $\eqv{(A*B)*C}{A*(B*C)}$. The proof of this Fubini theorem
  \index{Fubini theorem for colimits}
  for colimits still requires the long and tedious computation, though.
\end{proof}

\begin{lem}
  For any type $A$, there is an equivalence $\eqv{\susp A}{\bool*A}$.
\end{lem}

\begin{proof}
  It is easy to define the two maps back and forth and to prove that they are
  inverse to each other. The details are left as an exercise to the reader.
\end{proof}

We can now construct the Hopf fibration:

\begin{thm}
  There is a fibration over $\Sn^2$ of fiber $\Sn^1$ and total space $\Sn^3$.
  \index{total!space}%
\end{thm}
\begin{proof}
  We proved that $\Sn^1$ has a structure of H-space (cf \autoref{lem:hspace-S1})
  hence by \autoref{lem:hopf-construction} there is a fibration over $\Sn^2$ of
  fiber $\Sn^1$ and total space $\Sn^1*\Sn^1$. But by the two previous results
  and \autoref{thm:suspbool} we have:
  \begin{equation*}
    \Sn^1*\Sn^1 = (\susp\bool)*\Sn^1
    =(\bool*\bool)*\Sn^1
    =\bool*(\bool*\Sn^1)
    =\susp(\susp\Sn^1)
    =\Sn^3. \qedhere
  \end{equation*}
\end{proof}
\index{Hopf fibration|)}

\section{The Freudenthal suspension theorem}
\label{sec:freudenthal}

\index{Freudenthal suspension theorem|(}%
\index{theorem!Freudenthal suspension|(}%

Before proving the Freudenthal suspension theorem, we need some auxiliary lemmas about connectedness.
In \autoref{cha:hlevels} we proved a number of facts about $n$-connected maps and $n$-types for fixed $n$; here we are now interested in what happens when we vary $n$.
For instance, in \autoref{prop:nconnected_tested_by_lv_n_dependent types} we showed that $n$-connected maps are characterized by an ``induction principle'' relative to families of $n$-types.
If we want to ``induct along'' an $n$-connected map into a family of $k$-types for $k> n$, we don't immediately know that there is a function by such an induction principle, but the following lemma says that at least our ignorance can be quantified.

\begin{lem}\label{thm:conn-trunc-variable-ind}
  If $f:A\to B$ is $n$-connected and $P:B\to \ntype{k}$ is a family of $k$-types for $k\ge n$, then the induced function
  \[ (\blank\circ f) : \Parens{\prd{b:B} P(b)} \to \Parens{\prd{a:A} P(f(a)) } \]
  is $(k-n-2)$-truncated.
\end{lem}
\begin{proof}
  We induct on the natural number $k-n$.
  When $k=n$, this is \autoref{prop:nconnected_tested_by_lv_n_dependent types}.
  For the inductive step, suppose $f$ is $n$-connected and $P$ is a family of $k+1$-types.
  To show that $(\blank\circ f)$ is $(k-n-1)$-truncated, let $k:\prd{a:A} P(a)$; then we have
  \[ \hfib{(\blank\circ f)}{k} \eqvsym \sm{g:\prd{b:B} P(b)} \prd{a:A} g(f(a)) = k(a).\]
  Let $(g,p)$ and $(h,q)$ lie in this type, so $p:g\circ f \htpy k$ and $q:h\circ f \htpy k$; then we also have
  \[ \big((g,p) = (h,q)\big) \eqvsym
  \Parens{\sm{r:g\htpy h} r\circ f = p \ct \opp{q}}.
  \]
  However, here the right-hand side is a fiber of the map
  \[ (\blank\circ f) : \Parens{\prd{b:B} Q(b)} \to \Parens{\prd{a:A} Q(f(a)) } \]
  where $Q(b) \defeq (g(b)=h(b))$.
  Since $P$ is a family of $(k+1)$-types, $Q$ is a family of $k$-types, so the the inductive hypothesis implies that this fiber is a $(k-n-2)$-type.
  Thus, all path spaces of $\hfib{(\blank\circ f)}{k}$ are $(k-n-2)$-types, so it is a $(k-n-1)$-type.
\end{proof}

Recall that if $\pairr{A,a_0}$ and $\pairr{B,b_0}$ are pointed types, then
their \define{wedge}
\index{wedge}%
$A\vee B$ is defined to be the pushout of $A\xleftarrow{a_0}
\unit\xrightarrow{b_0} B$.
There is a canonical map $i:A\vee B \to A\times B$ defined by the two maps $\lam{a} (a,b_0)$ and $\lam{b} (a_0,b)$; the following lemma essentially says that this map is highly connected if $A$ and $B$ are so.
It is a bit more convenient both to prove and use, however, if we use the characterization of connectedness from \autoref{prop:nconnected_tested_by_lv_n_dependent types} and substitute in the universal property of the wedge (generalized to type families).

\begin{lem}[Wedge connectivity lemma]\label{thm:wedge-connectivity}
  Suppose that $\pairr{A,a_0}$ and $\pairr{B,b_0}$ are $n$- and $m$-connected pointed types, respectively, with $n,m\geq0$, and let 
\narrowequation{P:A\to B\to \ntype{(n+m)}.}
Then for any ${f:\prd{a:A} P(a,b_0)}$ and ${g:\prd{b:B} P(a_0,b)}$ with $p:f(a_0) = g(b_0)$, there exists $h:\prd{a:A}{b:B} P(a,b)$ with homotopies
\begin{equation*}
  q:\prd{a:A} h(a,b_0)=f(a)
  \qquad\text{and}\qquad
  r:\prd{b:B} h(a_0,b)=g(b)
 \end{equation*}
such that $p = \opp{q(a_0)} \ct r(b_0)$.
\end{lem}
\begin{proof}
  Define $P:A\to\type$ by
  \[ P(a) \defeq \sm{k:\prd{b:B} P(a,b)} (f(a) = k(b_0)). \]
  Then we have $(g,p):P(a_0)$.
  Since $a_0:\unit\to A$ is $(n-1)$-connected, if $P$ is a family of $(n-1)$-types then we will have $\ell:\prd{a:A} P(a)$ such that $\ell(a_0) = (g,p)$, in which case we can define $h(a,b) \defeq \proj1(\ell(a))(b)$.
  However, for fixed $a$, the type $P(a)$ is the fiber over $f(a)$ of the map
  \[ \Parens{\prd{b:B} P(a,b) } \to P(a,b_0) \]
  given by precomposition with $b_0:\unit\to B$.
  Since $b_0:\unit\to B$ is $(m-1)$-connected, for this fiber to be $(n-1)$-connected, by \autoref{thm:conn-trunc-variable-ind} it suffices for each type $P(a,b)$ to be an $(n+m)$-type, which we have assumed.
\end{proof}

Let $(X,x_0)$ be a pointed type, and recall the definition of the suspension $\susp X$ from \autoref{sec:suspension}, with constructors $\north,\south:\susp X$ and $\merid:X \to (\north=\south)$.
We regard $\susp X$ as a pointed space with basepoint $\north$, so that we have $\Omega\susp X \defeq (\id[\susp X]\north\north)$.
Then there is a canonical map
\begin{align*}
  \sigma &: X \to \Omega\susp X\\
  \sigma(x) &\defeq \merid(x) \ct \opp{\merid(x_0)}.
\end{align*}

\begin{rmk}
  In classical algebraic topology, one considers the \emph{reduced suspension}, in which the path $\merid(x_0)$ is collapsed down to a point, identifying $\north$ and $\south$.
  The reduced and unreduced suspensions are homotopy equivalent, so the distinction is invisible to our purely homotopy-theoretic eyes --- and higher inductive types only allow us to ``identify'' points up to a higher path anyway, there is no purpose to considering reduced suspensions in homotopy type theory.
  However, the ``unreducedness'' of our suspension is the reason for the (possibly unexpected) appearance of $\opp{\merid(x_0)}$ in the definition of $\sigma$.
\end{rmk}

Our goal is now to prove the following.

\begin{thm}[The Freudenthal suspension theorem]\label{thm:freudenthal}
  Suppose that $X$ is $n$-connected and pointed, with $n\geq 0$.
  Then the map $\sigma:X\to \Omega\susp(X)$ is $2n$-connected.
\end{thm}

\index{encode-decode method|(}%

We will use the encode-decode method, but applied in a slightly different way.
In most cases so far, we have used it to characterize the loop space $\Omega (A,a_0)$ of some type as equivalent to some other type $B$, by constructing a family $\code:A\to \type$ with $\code(a_0)\defeq B$ and a family of equivalences $\decode:\prd{x:A}\code(x) \eqvsym (a_0=x)$.

In this case, however, we want to show that $\sigma:X\to \Omega \susp X$ is $2n$-connected.
We could use a truncated version of the previous method, such as we will see in \autoref{sec:van-kampen}, to prove that $\trunc{2n}X\to \trunc{2n}{\Omega \susp X}$ is an equivalence---but this is a slightly weaker statement than the map being $2n$-connected (see \autoref{thm:conn-pik,thm:pik-conn}).
However, note that in the general case, to prove that $\decode(x)$ is an equivalence, we could equivalently be proving that its fibers are contractible, and we would still be able to use induction over the base type.
This we can generalize to prove connectedness of a map into a loop space, i.e.\ that the \emph{truncations} of its fibers are contractible.
Moreover, instead of constructing $\code$ and $\decode$ separately, we can construct directly a family of \emph{codes for the truncations of the fibers}.

\begin{defn}\label{thm:freudcode}
  If $X$ is $n$-connected and pointed with $n\geq 0$, then there is a family
  \begin{equation}
    \code:\prd{y:\susp X} (\north=y) \to \type\label{eq:freudcode}
  \end{equation}
  such that
  \begin{align}
    \code(\north,p) &\defeq \trunc{2n}{\hfib{\sigma}{p}}
    \jdeq \trunc{2n}{\tsm{x:X} (\merid(x) \ct \opp{\merid(x_0)} = p)}\label{eq:freudcodeN}\\
    \code(\south,q) &\defeq \trunc{2n}{\hfib{\merid}{q}}
    \jdeq \trunc{2n}{\tsm{x:X} (\merid(x) = q)}.\label{eq:freudcodeS}
  \end{align}
\end{defn}

Our eventual goal will be to prove that $\code(y,p)$ is contractible for all $y:\susp X$ and $p:\north=y$.
Applying this with $y\defeq \north$ will show that all fibers of $\sigma$ are $2n$-connected, and thus $\sigma$ is $2n$-connected.

\begin{proof}[Proof of \autoref{thm:freudcode}]
  We define $\code(y,p)$ by induction on $y:\susp X$, where the first two cases are~\eqref{eq:freudcodeN} and~\eqref{eq:freudcodeS}.
  It remains to construct, for each $x_1:X$, a dependent path
  \[ \dpath{\lam{y}(\north=y)\to\type}{\merid(x_1)}{\code(\north)}{\code(\south)}. \]
  By \autoref{thm:dpath-arrow}, this is equivalent to giving a family of paths
  \[ \prd{q:\north=\south} \code(\north)(\transfib{\lam{y}(\north=y)}{\opp{\merid(x_1)}}{q}) = \code(\south)(q). \]
  And by univalence and transport in path types, this is equivalent to a family of equivalences
  \[ \prd{q:\north=\south} \code(\north,q \ct \opp{\merid(x_1)}) \eqvsym \code(\south,q). \]
  We will define a family of maps
  \begin{equation}\label{eq:freudmap}
    \prd{q:\north=\south} \code(\north,q \ct \opp{\merid(x_1)}) \to \code(\south,q).
  \end{equation}
  and then show that they are all equivalences.
  Thus, let $q:\north=\south$; by the universal property of truncation and the definitions of $\code(\north,\blank)$ and $\code(\south,\blank)$, it will suffice to define for each $x_2:X$, a map
  \begin{equation*}
    \big(\merid(x_2)\ct \opp{\merid(x_0)} = q \ct \opp{\merid(x_1)}\big)
    \to \trunc{2n}{\tsm{x:X} (\merid(x) = q)}.
  \end{equation*}
  Now for each $x_1,x_2:X$, this type is $2n$-truncated, while $X$ is $n$-connected.
  Thus, by \autoref{thm:wedge-connectivity}, it suffices to define this map when $x_1$ is $x_0$, when $x_2$ is $x_0$, and check that they agree when both are $x_0$.

  When $x_1$ is $x_0$, the hypothesis is $r:\merid(x_2)\ct \opp{\merid(x_0)} = q \ct \opp{\merid(x_0)}$.
  Thus, by canceling $\opp{\merid(x_0)}$ from $r$ to get $r':\merid(x_2)=q$, so we can define the image to be $\tproj{2n}{(x_2,r')}$.

  When $x_2$ is $x_0$, the hypothesis is $r:\merid(x_0)\ct \opp{\merid(x_0)} = q \ct \opp{\merid(x_1)}$.
  Rearranging this, we obtain $r'':\merid(x_1)=q$, and we can define the image to be $\tproj{2n}{(x_1,r'')}$.

  Finally, when both $x_1$ and $x_2$ are $x_0$, it suffices to show the resulting $r'$ and $r''$ agree; this is an easy lemma about path composition.
  This completes the definition of~\eqref{eq:freudmap}.
  To show that it is a family of equivalences, since being an equivalence is a mere proposition and $x_0:\unit\to X$ is (at least) $(-1)$-connected, it suffices to assume $x_1$ is $x_0$.
  In this case, inspecting the above construction we see that it is essentially the $2n$-truncation of the function that cancels $\opp{\merid(x_0)}$, which is an equivalence.
\end{proof}

In addition to~\eqref{eq:freudcodeN} and~\eqref{eq:freudcodeS}, we will need to extract from the construction of $\code$ some information about how it acts on paths.
For this we use the following lemma.

\begin{lem}\label{thm:freudlemma}
  Let $A:\UU$, $B:A\to \UU$, and $C:\prd{a:A} B(a)\to\UU$, and also $a_1,a_2:A$ with $m:a_1=a_2$ and $b:B(a_2)$.
  Then the function
  \[\transfib{\widehat{C}}{\pairpath(m,t)}{\blank} : C(a_1,\transfib{B}{\opp m}{b}) \to C(a_2,b),\]
  where $t:\transfib{B}{m}{\transfib{B}{\opp m}{b}} = b$ is the obvious coherence path and $\widehat{C}:(\sm{a:A} B(a)) \to\type$ is the uncurried form of $C$, is equal to the equivalence obtained by univalence from the composite
  \begin{align}
    C(a_1,\transfib{B}{\opp m}{b})
    &= \transfib{\lam{a} B(a)\to \UU}{m}{C(a_1)}(b)
    \tag{by~\eqref{eq:transport-arrow}}\\
    &= C(a_2,b). \tag{by $\happly(\apd{C}{m},b)$}
  \end{align}
\end{lem}
\begin{proof}
  By path induction, we may assume $a_2$ is $a_1$ and $m$ is $\refl{a_1}$, in which case both functions are the identity.
\end{proof}

We apply this lemma with $A\defeq\susp X$ and $B\defeq \lam{y}(\north=y)$ and $C\defeq\code$, while $a_1\defeq\north$ and $a_2\defeq\south$ and $m\defeq \merid(x_1)$ for some $x_1:X$, and finally $b\defeq q$ is some path $\north=\south$.
The computation rule for induction over $\susp X$ identifies $\apd{C}{m}$ with a path constructed in a certain way out of univalence and function extensionality.
The second function described in \autoref{thm:freudlemma} essentially consists of undoing these applications of univalence and function extensionality, reducing back to the particular functions~\eqref{eq:freudmap} that we defined using \autoref{thm:wedge-connectivity}.
Therefore, \autoref{thm:freudlemma} says that transporting along $\pairpath(q,t)$ essentially recovers these functions.

Finally, by construction, when $x_1$ or $x_2$ coincides with $x_0$ and the input is in the image of $\tproj{2n}{\blank}$, we know more explicitly what these functions are.
Thus, for any $x_2:X$, we have
\begin{equation}
  \transfib{\hat{\code}}{\pairpath(\merid(x_0),t)}{\tproj{2n}{(x_2,r)}}
  =\tproj{2n}{(x_1,r')}\label{eq:freudcompute1}
\end{equation}
where $r:\merid(x_2) \ct \opp{\merid(x_0)} = \transfib{B}{\opp{\merid(x_0)}}{q}$ is arbitrary as before, and $r':\merid(x_2)=q$ is obtained from $r$ by identifying its end point with $q \ct \opp{\merid(x_0)}$ and canceling $\opp{\merid(x_0)}$.
Similarly, for any $x_1:X$, we have
\begin{equation}
  \transfib{\hat{\code}}{\pairpath(\merid(x_1),t)}{\tproj{2n}{(x_0,r)}}
  = \tproj{2n}{(x_1,r'')}\label{eq:freudcompute2}
\end{equation}
where $r:\merid(x_0) \ct \opp{\merid(x_0)} = \transfib{B}{\opp{\merid(x_1)}}{q}$, and $r'':\merid(x_1)=q$ is obtained by identifying its end point and rearranging paths.

\begin{proof}[Proof of \autoref{thm:freudenthal}]
  It remains to show that $\code(y,p)$ is contractible for each $y:\susp X$ and $p:\north=y$.
  First we must choose a center of contraction, say $c(y,p):\code(y,p)$.
  This corresponds to the definition of the function $\encode$ in our previous proofs, so we define it by transport.
  Note that in the special case when $y$ is $\north$ and $p$ is $\refl{\north}$, we have
  \[\code(\north,\refl{\north}) \jdeq \trunc{2n}{\tsm{x:X} (\merid(x) \ct \opp{\merid(x_0)} = \refl{\north})}.\]
  Thus, we can choose $c(\north,\refl{\north})\defeq \tproj{2n}{(x_0,\mathsf{rinv}_{\merid(x_0)})}$, where $\mathrm{rinv}_q$ is the obvious path $q\ct\opp q = \refl{}$ for any $q$.
  We can now obtain $c:\prd{y:\susp X}{p:\north=y} \code(y,p)$ by path induction on $p$, but it will be important below that we can also give a concrete definition in terms of transport:
  \[ c(y,p) \defeq \transfib{\hat{\code}}{\pairpath(p,\mathsf{tid}_p)}{c(\north,\refl{\north})}
  \]
  where $\hat{\code}: \big(\sm{y:\susp X} (\north=y)\big) \to \type$ is the uncurried version of \code, and $\mathsf{tid}_p:\trans{p}{\refl{}} = p$ is a standard lemma.

  Next, we must show that every element of $\code(y,p)$ is equal to $c(y,p)$.
  Again, by path induction, it suffices to assume $y$ is $\north$ and $p$ is $\refl{\north}$.
  In fact, we will prove it more generally when $y$ is $\north$ and $p$ is arbitrary.
  That is, we will show that for any $p:\north=\north$ and $d:\code(\north,p)$ we have $d = c(\north,p)$.
  Since this equality is a $(2n-1)$-type, we may assume $d$ is of the form $\tproj{2n}{(x_1,r)}$ for some $x_1:X$ and $r:\merid(x_1) \ct \opp{\merid(x_0)} = p$.

  Now by a further path induction, we may assume that $r$ is reflexivity, and $p$ is $\merid(x_1) \ct \opp{\merid(x_0)}$.
  (This is why we generalized to arbitrary $p$ above.)
  Thus, we have to prove that
  \begin{equation}
    \tproj{2n}{(x_1, \refl{\merid(x_1) \ct \opp{\merid(x_0)}})}
    \;=\;
    c\left(\north,\refl{\merid(x_1) \ct \opp{\merid(x_0)}}\right).\label{eq:freudgoal}
  \end{equation}
  By definition, the right-hand side of this equality is
  \begin{multline*}
    \Transfib{\hat{\code}}{\pairpath(\merid(x_1) \ct \opp{\merid(x_0)}, \nameless)}{\tproj{2n}{(x_0,\nameless)}} \\
    = \transfibf{\hat{\code}}
    \begin{aligned}[t]
      \Big(
      &{\pairpath(\opp{\merid(x_0)}, \nameless)},\\
      &{\Transfib{\hat{\code}}{\pairpath(\merid(x_1), \nameless)}{\tproj{2n}{(x_0,\nameless)}}}
      \Big)
    \end{aligned}
    \\
    = \Transfib{\hat{\code}}{\pairpath(\opp{\merid(x_0)}, \nameless)}{\tproj{2n}{(x_1,\nameless)}}
    = \tproj{2n}{(x_1,\nameless)}
  \end{multline*}
  where the underscore $\nameless$ ought to be filled in with suitable coherence paths.
  Here the first step is functoriality of transport, the second invokes~\eqref{eq:freudcompute2}, and the third invokes~\eqref{eq:freudcompute1} (with transport moved to the other side).
  Thus we have the same first component as the left-hand side of~\eqref{eq:freudgoal}.
  We leave it to the reader to verify that the coherence paths all cancel, giving reflexivity in the second component.
\end{proof}


\begin{cor}[Freudenthal Equivalence] \label{cor:freudenthal-equiv}
Suppose that $X$ is $n$-connected and pointed, with $n\geq 0$.
Then $\eqv{\trunc{2n}{X}}{\trunc{2n}{\Omega\susp(X)}}$.
\end{cor}
\begin{proof}
By \cref{thm:freudenthal}, $\sigma$ is $2n$-connected.  By
\cref{lem:connected-map-equiv-truncation}, it is therefore an
equivalence on $2n$-truncations.  
\end{proof}

\index{encode-decode method|)}%

\index{Freudenthal suspension theorem|)}%
\index{theorem!Freudenthal suspension|)}%

\index{homotopy!group!of sphere}%
\index{stability!of homotopy groups of spheres}%
\index{type!n-sphere@$n$-sphere}%
One important corollary of the Freudenthal suspension theorem is that the homotopy groups of
spheres are stable in a certain range (these are the northeast-to-southwest diagonals
in \autoref{tab:homotopy-groups-of-spheres}):

\begin{cor}[Stability for Spheres] \label{cor:stability-spheres}
If $k \le 2n-2$, then $\pi_{k+1}(S^{n+1}) = \pi_{k}(S^{n})$.
\end{cor}
\begin{proof}
Assume $k \le 2n-2$.  
By \cref{cor:sn-connected}, $\Sn ^{n}$ is $\nminusone$-connected.  Therefore,
by \cref{cor:freudenthal-equiv}, 
\[
\trunc{2(n-1)}{\Omega(\susp(\Sn^{n}))} = \trunc{2(n-1)}{\Sn^{n}}.
\]
By \cref{lem:truncation-le}, because $k \le 2(n-1)$, applying $\trunc{k}{\blank}$
to both sides shows that this equation holds for $k$:
\begin{equation}\label{eq:freudenthal-for-spheres}
\trunc{k}{\Omega(\susp(\Sn^{n}))} = \trunc{k}{\Sn^{n}}.
\end{equation}
Then, the main idea of the proof is as follows; we omit checking that these
equivalences act appropriately on the base points of these spaces:
\begin{align*}
\pi_{k+1}(\Sn^{n+1}) &\jdeq \trunc{0}{\Omega^{k+1}(\Sn^{n+1})} \\
                     &\jdeq \trunc{0}{\Omega^k(\Omega(\Sn^{n+1}))} \\
                     &\jdeq \trunc{0}{\Omega^k(\Omega(\susp(\Sn^{n})))} \\
                     &= \Omega^k(\trunc{k}{(\Omega(\susp(\Sn^{n})))})
                     \tag{by \autoref{thm:path-truncation}}\\
                     &= \Omega^k(\trunc{k}{\Sn^{n}})
                     \tag{by \eqref{eq:freudenthal-for-spheres}}\\
                     &= \trunc{0}{\Omega^k(\Sn^{n})}
                     \tag{by \autoref{thm:path-truncation}}\\
                     &\jdeq \pi_k(\Sn^{n}). \qedhere
\end{align*}
\end{proof}

This means that once we have calculated one entry in one of these stable
diagonals, we know all of them.  For example:
\begin{thm}
$\pi_n(\Sn^n)=\Z$ for every $n\geq 1$. 
\end{thm}

\begin{proof}
The proof is by induction on $n$.  We already have $\pi_1(\Sn ^1) = \Z$
(\autoref{cor:pi1s1}) and $\pi_2(\Sn ^2) = \Z$ (\autoref{cor:pis2-hopf}).
When $n \ge 2$, $n \le (2n - 2)$. Therefore, by
\cref{cor:stability-spheres}, $\pi_{n+1}(S^{n+1}) = \pi_{n}(S^{n})$, and
this equivalence, combined with the inductive hypothesis, gives the result.  
\end{proof}

\begin{cor}
  $\Sn^{n+1}$ is not an $n$-type for any $n\ge -1$.
\end{cor}

\section{The van Kampen theorem}
\label{sec:van-kampen}

\index{van Kampen theorem|(}%
\index{theorem!van Kampen|(}%

\index{fundamental!group}%
The van Kampen theorem calculates the fundamental group $\pi_1$ of a (homotopy) pushout of spaces.
It is traditionally stated for a topological space $X$ which is the union of two open subspaces $U$ and $V$, but in homotopy-theoretic terms this is just a convenient way of ensuring that $X$ is the pushout of $U$ and $V$ over their intersection.
Thus, we will prove a version of the van Kampen theorem for arbitrary pushouts.

In this section we will describe a proof of the van Kampen theorem which uses the same encode-decode method that we used for $\pi_1(\Sn^1)$ in \autoref{sec:pi1-s1-intro}.
There is also a more homotopy-theoretic approach; see \autoref{ex:rezk-vankampen}.

We need a more refined version of the encode-decode method.
In \autoref{sec:pi1-s1-intro} (as well as in \autoref{sec:compute-coprod,sec:compute-nat}) we used it to characterize the path space of a (higher) inductive type $W$ --- deriving as a consequence a characterization of the loop space $\Omega(W)$, and thereby also of its 0-truncation $\pi_1(W)$.
In the van Kampen theorem, our goal is only to characterize the fundamental group $\pi_1(W)$, and we do not have any explicit description of the loop spaces or the path spaces to use.

It turns out that we can use the same technique directly for a truncated version of the path fibration, thereby characterizing not only the fundamental \emph{group} $\pi_1(W)$, but also the whole fundamental \emph{groupoid}.
\index{fundamental!pregroupoid}%
Spe\-cif\-ical\-ly, for a type $X$, write $\Pi_1 X: X\to X\to \type$ for the $0$-truncation of its identity type, i.e.\ $\Pi_1 X(x,y) \defeq \trunc0{x=y}$.
Note that we have induced groupoid operations
\begin{align*}
  (\blank\ct\blank) &\;:\; \Pi_1X(x,y) \to \Pi_1X(y,z) \to \Pi_1X(x,z)\\
  \opp{(\blank)} &\;:\; \Pi_1X(x,y) \to \Pi_1X(y,x)\\
  \refl{x} &\;:\; \Pi_1X(x,x)\\
  \apfunc{f} &\;:\; \Pi_1X(x,y) \to \Pi_1Y(fx,fy)
\end{align*}
for which we use the same notation as the corresponding operations on paths.

\subsection{Naive van Kampen}
\label{sec:naive-vankampen}

We begin with a ``naive'' version of the van Kampen theorem, which is useful but not quite as useful as the classical version.
In \autoref{sec:better-vankampen} we will improve it to a more useful version.

\index{encode-decode method|(}%

Given types $A,B,C$ and functions $f:A\to B$ and $g:A\to C$, let $P$ be their pushout $B\sqcup^A C$.
As we saw in \autoref{sec:colimits}, $P$ is the higher inductive type generated by
\begin{itemize}
\item $i:B\to P$,
\item $j:C\to P$, and
\item for all $x:A$, a path $k x:ifx = jgx$.
\end{itemize}
Define $\code:P\to P\to \type$ by double induction on $P$ as follows.
\begin{itemize}
\item $\code(ib,ib')$ is a set-quotient (see \autoref{sec:set-quotients}) of the type of sequences 
  \[ (b, p_0, x_1, q_1, y_1, p_1, x_2, q_2, y_2, p_2, \dots, y_n, p_n, b') \]
  where
  \begin{itemize}
  \item $n:\mathbb{N}$
  \item $x_k:A$ and $y_k:A$ for $0<k \le n$
  \item $p_0:\Pi_1B(b,f x_1)$ and $p_n:\Pi_1B(f y_n, b')$ for $n>0$, and $p_0:\Pi_1B(b,b')$ for $n=0$
  \item $p_k:\Pi_1B(f y_k, fx_{k+1})$ for $1\le k < n$
  \item $q_k:\Pi_1C(gx_k, gy_k)$ for $1\le k\le n$
  \end{itemize}
  The quotient is generated by the following equalities:
  \begin{align*}
    (\dots, q_k, y_k, \refl{fy_k}, y_k, q_{k+1},\dots)
    &= (\dots, q_k\ct q_{k+1},\dots)\\
    (\dots, p_k, x_k, \refl{gx_k}, x_k, p_{k+1},\dots)
    &= (\dots, p_k\ct p_{k+1},\dots)
  \end{align*}
  (see \autoref{rmk:naive} below).
  We leave it to the reader to define this type of sequences precisely as an inductive type.
\item $\code(jc,jc')$ is identical, with the roles of $B$ and $C$ reversed.
  We likewise notationally reverse the roles of $x$ and $y$, and of $p$ and $q$.
\item $\code(ib,jc)$ and $\code(jc,ib)$ are similar, with the parity changed so that they start in one type and end in the other.
\item For $a:A$ and $b:B$, we require an equivalence
  \begin{equation}
    \code(ib, ifa) \eqvsym \code(ib,jga).\label{eq:bfa-bga}
  \end{equation}
  We define this to consist of the two functions defined on sequences by
  \begin{align*}
    (\dots, y_n, p_n,fa) &\mapsto (\dots,y_n,p_n,a,\refl{ga},ga),\\
    (\dots, x_n, p_n, a, \refl{fa}, fa) &\mapsfrom (\dots, x_n, p_n, ga).
  \end{align*}
  Both of these functions are easily seen to respect the equivalence relations, and hence to define functions on the types of codes.
  The left-to-right-to-left composite is
  \[ (\dots, y_n, p_n,fa) \mapsto
  (\dots,y_n,p_n,a,\refl{ga},a,\refl{fa},fa)
  \]
  which is equal to the identity by a generating equality of the quotient.
  The other composite is analogous.
  Thus we have defined an equivalence~\eqref{eq:bfa-bga}.
\item Similarly, we require equivalences
  \begin{align*}
    \code(jc,ifa) &\eqvsym \code(jc,jga)\\
    \code(ifa,ib)&\eqvsym (jga,ib)\\
    \code(ifa,jc)&\eqvsym (jga,jc)
  \end{align*}
  all of which are defined in exactly the same way (the second two by adding reflexivity terms on the beginning rather than the end).
\item Finally, we need to know that for $a,a':A$, the following diagram commutes:
  \begin{equation}\label{eq:bfa-bga-comm}
  \vcenter{\xymatrix{
      \code(ifa,ifa') \ar[r]\ar[d] &
      \code(ifa,jga')\ar[d]\\
      \code(jga,ifa')\ar[r] &
      \code(jga,jga')
      }}
  \end{equation}
  This amounts to saying that if we add something to the beginning and then something to the end of a sequence, we might as well have done it in the other order.
\end{itemize}

\begin{rmk}\label{rmk:naive}
  One might expect to see in the definition of \code some additional generating equations for the set-quotient, such as
  \begin{align*}
    (\dots, p_{k-1} \ct fw, x_{k}', q_{k}, \dots) &=
    (\dots, p_{k-1}, x_{k}, gw \ct q_{k}, \dots)
    \tag{for $w:\Pi_1A(x_{k},x_{k}')$}\\
    (\dots, q_k \ct gw, y_k', p_k, \dots) &=
    (\dots, q_k, y_k, fw \ct p_k, \dots).
    \tag{for $w:\Pi_1A(y_k, y_k')$}
  \end{align*}
  However, these are not necessary!
  In fact, they follow automatically by path induction on $w$.
  This is the main difference between the ``naive'' van Kampen theorem and the more refined one we will consider in the next subsection.
\end{rmk}

Continuing on, we can characterize transporting in the fibration $\code$:
\begin{itemize}
\item For $p:b=_B b'$ and $u:P$, we have
  \[ \mathsf{transport}^{b\mapsto \code(u,ib)}(p, (\dots, y_n,p_n,b))
  = (\dots,y_n,p_n\ct p,b').
  \]
\item For $q:c=_C c'$ and $u:P$, we have
  \[ \mathsf{transport}^{c\mapsto \code(u,jc)}(q, (\dots, x_n,q_n,c))
  = (\dots,x_n,q_n\ct q,c').
  \]
\end{itemize}
Here we are abusing notation by using the same name for a path in $X$ and its image in $\Pi_1X$.
Note that transport in $\Pi_1X$ is also given by concatenation with (the image of) a path.
From this we can prove the above statements by induction on $u$.
We also have:
\begin{itemize}
\item For $a:A$ and $u:P$,
  \[ \mathsf{transport}^{v\mapsto \code(u,v)}(ha, (\dots, y_n,p_n,fa))
  = (\dots,y_n,p_n,a,\refl{ga},ga).
  \]
\end{itemize}
This follows essentially from the definition of $\code$.

We also construct a function
\[ r : \prd{u:P} \code(u,u) \]
by induction on $u$ as follows:
\begin{align*}
  rib &\defeq (b,\refl{b},b)\\
  rjc &\defeq (c,\refl{c},c)
\end{align*}
and for $rka$ we take the composite equality
\begin{align*}
  (ka,ka)_* (fa,\refl{fa},fa)
  &= (ga,\refl{ga},a,\refl{fa},a,\refl{ga},ga) \\
  &= (ga,\refl{ga},ga)
\end{align*}
where the first equality is by the observation above about transporting in $\code$, and the second is an instance of the set quotient relation used to define $\code$.

We will now prove:
\begin{thm}[Naive van Kampen theorem]\label{thm:naive-van-kampen}
  For all $u,v:P$ there is an equivalence
  \[ \Pi_1P(u,v) \eqvsym \code(u,v). \]
\end{thm}
\begin{proof}

To define a function
\[ \encode : \Pi_1P(u,v) \to \code(u,v) \]
it suffices to define a function $(u=_P v) \to \code(u,v)$,
since $\code(u,v)$ is a set.
We do this by transport:
\[\encode(p) \defeq \mathsf{transport}^{v\mapsto \code(u,v)}(p,r(u)).\]
Now to define
\[ \decode: \code(u,v) \to \Pi_1P(u,v) \]
we proceed as usual by induction on $u,v:P$.
In each case for $u$ and $v$, we apply $i$ or $j$ to all the equalities $p_k$ and $q_k$ as appropriate and concatenate the results in $P$, using $h$ to identify the endpoints.
For instance, when $u\jdeq ib$ and $v\jdeq ib'$, we define
\begin{narrowmultline}\label{eq:decode}
 \decode(b, p_0, x_1, q_1, y_1, p_1, \dots, y_n, p_n, b') \defeq\narrowbreak
 (p_0)\ct h(x_1) \ct j(q_1) \ct \opp{h(y_1)} \ct i(p_1) \ct \cdots \ct \opp{h(y_n)}\ct i(p_n).
\end{narrowmultline}
This respects the set-quotient equivalence relation and the equivalences such as~\eqref{eq:bfa-bga}, since $h: fi \htpy gj$ is natural and $f$ and $g$ are functorial.

As usual, to show that the composite
\[ \Pi_1P(u,v) \xrightarrow{\encode} \code(u,v) \xrightarrow{\decode} \Pi_1P(u,v) \]
is the identity, we first peel off the 0-truncation (since the codomain is a set) and then apply path induction.
The input $\refl{u}$ goes to $ru$, which then goes back to $\refl u$ (applying a further induction on $u$ to decompose $\decode(ru)$).

Finally, consider the composite
\[  \code(u,v) \xrightarrow{\decode} \Pi_1P(u,v) \xrightarrow{\encode} \code(u,v). \]
We proceed by induction on $u,v:P$.
When $u\jdeq ib$ and $v\jdeq ib'$, this composite is
\begin{narrowmultline*}
(b, p_0, x_1, q_1, y_1, p_1, \dots, y_n, p_n, b')
\narrowbreak
\begin{aligned}[t]
  &\mapsto \Big(ip_0\ct hx_1 \ct jq_1 \ct \opp{hy_1} \ct ip_1 \ct \cdots \ct \opp{hy_n}\ct ip_n\Big)_*(rib)\\
  &= (ip_n)_* \cdots(jq_1)_* (hx_1)_*(ip_0)_*(b,\refl{b},b)\\
  &= (ip_n)_* \cdots(jq_1)_* (hx_1)_*(b,p_0,ifx_1)\\
  &= (ip_n)_* \cdots(jq_1)_* (b,p_0,x_1,\refl{gx_1},jgx_1)\\
  &= (ip_n)_* \cdots (b,p_0,x_1,q_1,jgy_1)\\
  &= \quad\vdots\\
  &= (b, p_0, x_1, q_1, y_1, p_1, \dots, y_n, p_n, b').
\end{aligned}  
\end{narrowmultline*}
i.e., the identity function.
(To be precise, there is an implicit inductive argument needed here.)
The other three point cases are analogous, and the path cases are trivial since all the types are sets.
\end{proof}

\index{encode-decode method|)}%

\autoref{thm:naive-van-kampen} allows us to calculate the fundamental groups of many types, provided $A$ is a set,
for in that case, each $\code(u,v)$ is, by definition, a set-quotient of a \emph{set} by a relation.

\begin{eg}\label{eg:circle}
  Let $A\defeq \bool$, $B\defeq\unit$, and $C\defeq \unit$.
  Then $P \eqvsym S^1$.
  Inspecting the definition of, say, $\code(i(\ttt),i(\ttt))$, we see that the paths all may as well be trivial, so the only information is in the sequence of elements $x_1,y_1,\dots,x_n,y_n: \bool$.
  Moreover, if we have $x_k=y_k$ or $y_k=x_{k+1}$ for any $k$, then the set-quotient relations allow us to excise both of those elements.
  Thus, every such sequence is equal to a canonical \emph{reduced} one in which no two adjacent elements are equal.
  Clearly such a reduced sequence is uniquely determined by its length (a natural number $n$) together with, if $n>1$, the information of whether $x_1$ is $\bfalse$ or $\btrue$, since that determines the rest of the sequence uniquely.
  And these data can, of course, be identified with an integer, where $n$ is the absolute value and $x_1$ encodes the sign.
  Thus we recover $\pi_1(S^1)\cong \Z$.
\end{eg}

Since \autoref{thm:naive-van-kampen} asserts only a bijection of families of sets, this isomorphism $\pi_1(S^1)\cong \Z$ is likewise only a bijection of sets.
We could, however, define a concatenation operation on $\code$ (by concatenating sequences) and show that $\encode$ and $\decode$ form an isomorphism respecting this structure.
(In the language of \autoref{cha:category-theory}, these would be ``pregroupoids''.)
We leave the details to the reader.

\index{fundamental!group|(}%

\begin{eg}\label{eg:suspension}
  More generally, let $B\defeq\unit$ and $C\defeq \unit$ but $A$ be arbitrary, so that $P$ is the suspension of $A$.
  Then once again the paths $p_k$ and $q_k$ are trivial, so that the only information in a path code is a sequence of elements $x_1,y_1,\dots,x_n,y_n: A$.
  The first two generating equalities say that adjacent equal elements can be canceled, so it makes sense to think of this sequence as a word of the form
  \[ x_1 y_1^{-1} x_2 y_2^{-1} \cdots x_n y_n^{-1} \]
  in a group.
  Indeed, it looks similar to the free group on $A$ (or equivalently on $\trunc0A$; see \autoref{thm:freegroup-nonset}), but we are considering only words that start with a non-inverted element, alternate between inverted and non-inverted elements, and end with an inverted one.
  This effectively reduces the size of the generating set by one.
  For instance, if $A$ has a point $a:A$, then we can identify $\pi_1(\susp A)$ with the group presented by $\trunc0A$ as generators with the relation $\tproj0a = e$; see \autoref{ex:vksusppt,ex:vksuspnopt} for details.
\end{eg}

\begin{eg}\label{eg:wedge}
  Let $A\defeq\unit$ and $B$ and $C$ be arbitrary, so that $f$ and $g$ simply equip $B$ and $C$ with basepoints $b$ and $c$, say.
  Then $P$ is the \emph{wedge} $B\vee C$ of $B$ and $C$ (the coproduct in the category of based spaces).
  In this case, it is the elements $x_k$ and $y_k$ which are trivial, so that the only information is a sequence of loops $(p_0,q_1,p_1,\dots,p_n)$ with $p_k:\pi_1(B,b)$ and $q_k:\pi_1(C,c)$.
  Such sequences, modulo the equivalence relation we have imposed, are easily identified with the explicit description of the \emph{free product} of the groups $\pi_1(B,b)$ and $\pi_1(C,c)$, as constructed in \autoref{sec:free-algebras}.
  Thus, we have $\pi_1(B\vee C) \cong \pi_1(B) * \pi_1(C)$.
\end{eg}

\index{fundamental!group|)}%

However, \autoref{thm:naive-van-kampen} stops just short of being the full classical van Kampen theorem, which handles
the case where $A$ is not necessarily a set, 
and states that $\pi_1(B\sqcup^A C) \cong \pi_1(B) *_{\pi_1(A)} \pi_1(C)$ (with base point coming from $A$).
Indeed, the conclusion of \autoref{thm:naive-van-kampen} says nothing at all about $\pi_1(A)$; the paths in $A$ are ``built into the quotienting'' in a type-theoretic way that makes it hard to extract explicit information, in that $\code(u,v)$ is a set-quotient of a non-set by a relation.
For this reason, in the next subsection we consider a better version of the van Kampen theorem.

\subsection{The van Kampen theorem with a set of basepoints}
\label{sec:better-vankampen}

\index{basepoint!set of}%
The improvement of van Kampen we present now is closely analogous to a similar improvement in classical algebraic topology, where $A$ is equip\-ped with a \emph{set $S$ of base points}.
In fact, it turns out to be unnecessary for our proof to assume that the ``set of basepoints'' is a \emph{set} --- it might just as well be an arbitrary type; the utility of assuming $S$ is a set arises later, when applying the theorem to obtain computations.
What is important is that $S$ contains at least one point in each connected component of $A$.
We state this in type theory by saying that we have a type $S$ and a function $k:S \to A$ which is surjective, i.e.\ $(-1)$-connected.
If $S\jdeq A$ and $k$ is the identity function, then we will recover the naive van Kampen theorem.
Another example to keep in mind is when $A$ is pointed and (0-)connected, with $k:\unit\to A$ the point: by \autoref{thm:minusoneconn-surjective,thm:connected-pointed} this map is surjective just when $A$ is 0-connected.

Let $A,B,C,f,g,P,i,j,h$ be as in the previous section.
We now define, given our surjective map $k:S\to A$, an auxiliary type which improves the connectedness of $k$.
Let $T$ be the higher inductive type generated by
\begin{itemize}
\item A function $\ell:S\to T$, and
\item For each $s,s':S$, a function $m:(\id[A]{ks}{ks'}) \to (\id[T]{\ell s}{\ell s'})$.
\end{itemize}
There is an obvious induced function $\kbar:T\to A$ such that $\kbar \ell = k$, and any $p:ks=ks'$ is equal to the composite $ks = \kbar \ell s \overset{\kbar m p}{=} \kbar \ell s' = k s'$.

\begin{lem}\label{thm:kbar}
  $\kbar$ is 0-connected.
\end{lem}
\begin{proof}
  We must show that for all $a:A$, the 0-truncation of the type $\sm{t:T}(\kbar t = a)$ is contractible.
  Since contractibility is a mere proposition and $k$ is $(-1)$-connected, we may assume that $a=ks$ for some $s:S$.
  Now we can take the center of contraction to be $\tproj0{(\ell s,q)}$ where $q$ is the equality $\kbar\ell s = k s$.

  It remains to show that for any $\phi:\trunc0{\sm{t:T} (\kbar t = ks)}$ we have $\phi = \tproj0{(\ell s,q)}$.
  Since the latter is a mere proposition, and in particular a set, we may assume that $\phi=\tproj0{(t,p)}$ for $t:T$ and $p:\kbar t = ks$.

  Now we can do induction on $t:T$.
  If $t\jdeq\ell s'$, then $ks' = \kbar \ell s' \overset{p}{=} ks$ yields via $m$ an equality $\ell s = \ell s'$.
  Hence by definition of $\kbar$ and of equality in homotopy fibers, we obtain an equality $(ks',p) = (ks,q)$, and thus $\tproj0{(ks',p)} = \tproj0{(ks,q)}$.
  Next we must show that as $t$ varies along $m$ these equalities agree.
  But they are equalities in a set (namely $\trunc0{\sm{t:T} (\kbar t = ks)}$), and hence this is automatic.
\end{proof}

\begin{rmk}
  \index{kernel!pair}%
  $T$ can be regarded as the (homotopy) coequalizer of the ``kernel pair'' of $k$.
  If $S$ and $A$ were sets, then the $(-1)$-connectivity of $k$ would imply that $A$ is the $0$-truncation of this coequalizer (see \autoref{cha:set-math}).
  For general types, higher topos theory suggests that $(-1)$-con\-nec\-tiv\-i\-ty of $k$ will imply instead that $A$ is the colimit (a.k.a.\ ``geometric realization'') of the ``simplicial kernel'' of $k$.
  \index{.infinity1-topos@$(\infty,1)$-topos}%
  \index{geometric realization}%
  \index{simplicial!kernel}%
  \index{kernel!simplicial}%
  The type $T$ is the colimit of the ``1-skeleton'' of this simplicial kernel, so it makes sense that it improves the connectivity of $k$ by $1$.
  More generally, we might expect the colimit of the $n$-skeleton\index{skeleton!of a CW-complex} to improve connectivity by $n$.
\end{rmk}

\index{encode-decode method|(}%

Now we define $\code:P\to P\to \type$ by double induction as follows
\begin{itemize}
\item $\code(ib,ib')$ is now a set-quotient of the type of sequences
  \[ (b, p_0, x_1, q_1, y_1, p_1, x_2, q_2, y_2, p_2, \dots, y_n, p_n, b') \]
  where
  \begin{itemize}
  \item $n:\mathbb{N}$,
  \item $x_k:S$ and $y_k:S$ for $0<k \le n$,
  \item $p_0:\Pi_1B(b,f k x_1)$ and $p_n:\Pi_1B(f k y_n, b')$ for $n>0$, and $p_0:\Pi_1B(b,b')$ for $n=0$,
  \item $p_k:\Pi_1B(fk y_k, fkx_{k+1})$ for $1\le k < n$,
  \item $q_k:\Pi_1C(gkx_k, gky_k)$ for $1\le k\le n$.
  \end{itemize}
  The quotient is generated by the following equalities (see \autoref{rmk:naive}):
  \begin{align*}
    (\dots, q_k, y_k, \refl{fy_k}, y_k, q_{k+1},\dots)
    &= (\dots, q_k\ct q_{k+1},\dots)\\
    (\dots, p_k, x_k, \refl{gx_k}, x_k, p_{k+1},\dots)
    &= (\dots, p_k\ct p_{k+1},\dots)\\
    (\dots, p_{k-1} \ct fw, x_{k}', q_{k}, \dots) &=
    (\dots, p_{k-1}, x_{k}, gw \ct q_{k}, \dots)
    \tag{for $w:\Pi_1A(kx_{k},kx_{k}')$}\\
    (\dots, q_k \ct gw, y_k', p_k, \dots) &=
    (\dots, q_k, y_k, fw \ct p_k, \dots).
    \tag{for $w:\Pi_1A(ky_k, ky_k')$}
  \end{align*}
  We will need below the definition of the case of $\decode$ on such a sequence, which as before concatenates all the paths $p_k$ and $q_k$ together with instances of $h$ to give an element of $\Pi_1P(ifb,ifb')$, cf.~\eqref{eq:decode}.
  As before, the other three point cases are nearly identical.
\item For $a:A$ and $b:B$, we require an equivalence
  \begin{equation}
    \code(ib, ifa) \eqvsym \code(ib,jga).\label{eq:bfa-bga2}
  \end{equation}
  Since $\code$ is set-valued, by \autoref{thm:kbar} we may assume that $a=\kbar t$ for some $t:T$.
  Next, we can do induction on $t$.
  If $t\jdeq \ell s$ for $s:S$, then we define~\eqref{eq:bfa-bga2} as in \autoref{sec:naive-vankampen}:
  \begin{align*}
    (\dots, y_n, p_n,fks) &\mapsto (\dots,y_n,p_n,s,\refl{gks},gks),\\
    (\dots, x_n, p_n, s, \refl{fks}, fks) &\mapsfrom (\dots, x_n, p_n, gks).
  \end{align*}
  These respect the equivalence relations, and define quasi-inverses just as before.
  Now suppose $t$ varies along $m_{s,s'}(w)$ for some $w:ks=ks'$; we must show that~\eqref{eq:bfa-bga2} respects transporting along $\kbar mw$.
  By definition of $\kbar$, this essentially boils down to transporting along $w$ itself.
  By the characterization of transport in path types, what we need to show is that
  \[ w_*(\dots, y_n, p_n,fks) = (\dots,y_n, p_n \ct fw, fks') \]
  is mapped by~\eqref{eq:bfa-bga2} to
  \[ w_*(\dots,y_n,p_n,s,\refl{gks},gks) = (\dots, y_n, p_n, s, \refl{gks} \ct gw, gks') \]
  But this follows directly from the new generators we have imposed on the set-quotient relation defining \code.
\item The other three requisite equivalences are defined similarly.
\item Finally, since the commutativity~\eqref{eq:bfa-bga-comm} is a mere proposition, by $(-1)$-connectedness of $k$ we may assume that $a=ks$ and $a'=ks'$, in which case it follows exactly as before.
\end{itemize}

\begin{thm}[van Kampen with a set of basepoints]\label{thm:van-Kampen}
  For all $u,v:P$ there is an equivalence
  \[ \Pi_1P(u,v) \eqvsym \code(u,v). \]
  with \code defined as in this section.
\end{thm}

\begin{proof}
  Basically just like before.
  To show that $\decode$ respects the new generators of the quotient relation, we use the naturality of $h$.
  And to show that $\decode$ respects the equivalences such as~\eqref{eq:bfa-bga2}, we need to induct on $\kbar$ and on $T$ in order to decompose those equivalences into their definitions, but then it becomes again simply functoriality of $f$ and $g$.
  The rest is easy.
  In particular, no additional argument is required for $\encode\circ\decode$, since the goal is to prove an equality in a set, and so the case of $h$ is trivial.
\end{proof}

\index{encode-decode method|)}%

\index{fundamental!group|(}%
\autoref{thm:van-Kampen} allows us to calculate the fundamental group of a space~$A$,
even when $A$ is not a set, provided $S$ is a set, for in that case,
each $\code(u,v)$ is, by definition, a set-quotient of a \emph{set} by a
relation.  In that respect, it is an improvement over
\autoref{thm:naive-van-kampen}.

\begin{eg}\label{eg:clvk}
  Suppose $S\defeq \unit$, so that $A$ has a basepoint $a \defeq k(\ttt)$ and is connected.
  Then code for loops in the pushout can be identified with alternating sequences of loops in $\pi_1(B,f(a))$ and $\pi_1(C,g(a))$, modulo an equivalence relation which allows us to slide elements of $\pi_1(A,a)$ between them (after applying $f$ and $g$ respectively).
  Thus, $\pi_1(P)$ can be identified with the \emph{amalgamated free product}
  \index{amalgamated free product}%
  \index{free!product!amalgamated}%
  $\pi_1(B) *_{\pi_1(A)} \pi_1(C)$ (the pushout in the category of groups), as constructed in \autoref{sec:free-algebras}.
  This (in the case when $B$ and $C$ are open subspaces of $P$ and $A$ their intersection) is probably the most classical version of the van Kampen theorem.
\end{eg}

\begin{eg}\label{eg:cofiber}
  \index{cofiber}
  As a special case of \autoref{eg:clvk}, suppose additionally that $C\defeq\unit$, so that $P$ is the cofiber $B/A$.
  Then every loop in $C$ is equal to reflexivity, so the relations on path codes allow us to collapse all sequences to a single loop in $B$.
  The additional relations require that multiplying on the left, right, or in the middle by an element in the image of $\pi_1(A)$ is the identity.
  We can thus identify $\pi_1(B/A)$ with the quotient of the group $\pi_1(B)$ by the normal subgroup generated by the image of $\pi_1(A)$.
\end{eg}

\begin{eg}\label{eg:torus}
  \index{torus}
  As a further special case of \autoref{eg:cofiber}, let $B\defeq S^1 \vee S^1$, let $A\defeq S^1$, and let $f:A\to B$ pick out the composite loop $p \ct q \ct \opp p \ct \opp q$, where $p$ and $q$ are the generating loops in the two copies of $S^1$ comprising $B$.
  Then $P$ is a presentation of the torus $T^2$.
  Indeed, it is not hard to identify $P$ with the presentation of $T^2$ as described in \autoref{sec:hubs-spokes}, using the cone on a particular loop.
  Thus, $\pi_1(T^2)$ is the quotient of the free group on two generators\index{generator!of a group} (i.e., $\pi_1(B)$) by the relation $p \ct q \ct \opp p \ct \opp q = 1$.
  This clearly yields the free \emph{abelian}\index{group!abelian} group on two generators, which is $\Z\times\Z$.
\end{eg}

\begin{eg}
  \index{CW complex}
  \index{hub and spoke}
  More generally, any CW complex can be obtained by repeatedly ``coning off'' spheres, as described in \autoref{sec:hubs-spokes}.
  That is, we start with a set $X_0$ of points (``0-cells''), which is the ``0-skeleton'' of the CW complex.
  We take the pushout
  \begin{equation*}
    \vcenter{\xymatrix{
        S_1 \times \Sn^0\ar[r]^-{f_1}\ar[d] &
        X_0\ar[d]\\
        \unit \ar[r] &
        X_1
      }}
  \end{equation*}
  for some set $S_1$ of 1-cells and some family $f_1$ of ``attaching maps'', obtaining the ``1-skeleton''\index{skeleton!of a CW-complex} $X_1$.
  \index{attaching map}%
  Then we take the pushout
  \begin{equation*}
    \vcenter{\xymatrix{
        S_2 \times \Sn^1\ar[r]^{f_2}\ar[d] &
        X_1\ar[d]\\
        \unit \ar[r] &
        X_2
      }}
  \end{equation*}
  for some set $S_2$ of 2-cells and some family $f_2$ of attaching maps, obtaining the 2-skeleton $X_2$, and so on.
  The fundamental group of each pushout can be calculated from the van Kampen theorem: we obtain the group presented by generators derived from the 1-skeleton, and relations derived from $S_2$ and $f_2$.
  The pushouts after this stage do not alter the fundamental group, since $\pi_1(\Sn^n)$ is trivial for $n>1$ (see \autoref{sec:pik-le-n}).
\end{eg}

\begin{eg}\label{eg:kg1}
  In particular, suppose given any presentation\index{presentation!of a group} of a (set-)group $G = \langle X \mid R \rangle$, with $X$ a set of generators and $R$ a set of words in these generators\index{generator!of a group}.
  Let $B\defeq \bigvee_X S^1$ and $A\defeq \bigvee_R S^1$, with $f:A\to B$ sending each copy of $S^1$ to the corresponding word in the generating loops of $B$.
  It follows that $\pi_1(P) \cong G$; thus we have constructed a connected type whose fundamental group is $G$.
  Since any group has a presentation, any group is the fundamental group of some type.
  If we 1-truncate such a type, we obtain a type whose only nontrivial homotopy group is $G$; this is called an \define{Eilenberg--Mac Lane space} $K(G,1)$.%
  \indexdef{Eilenberg--Mac Lane space}%
\end{eg}

\index{fundamental!group|)}%

\index{van Kampen theorem|)}%
\index{theorem!van Kampen|)}%

\section{Whitehead's theorem and Whitehead's principle}
\label{sec:whitehead}

In classical homotopy theory, a map $f:A\to B$ which induces an isomorphism $\pi_n(A,a) \cong \pi_n(B,f(a))$ for all points $a$ in $A$ (and also an isomorphism $\pi_0(A)\cong\pi_0(B)$) is necessarily a homotopy equivalence, as long as the spaces $A$ and $B$ are well-behaved (e.g.\ have the homotopy types of CW-complexes).
\index{theorem!Whitehead's}%
\index{Whitehead's!theorem}%
This is known as \emph{Whitehead's theorem}.
In fact, the ``ill-behaved'' spaces for which Whitehead's theorem fails are invisible to type theory.
Roughly, the well-behaved topological spaces suffice to present $\infty$-groupoids,%
\index{.infinity-groupoid@$\infty$-groupoid}
and homotopy type theory deals with $\infty$-groupoids directly rather than actual topological spaces.
Thus, one might expect that Whitehead's theorem would be true in univalent foundations.

However, this is \emph{not} the case: Whitehead's theorem is not provable.
In fact, there are known models of type theory in which it fails to be true, although for entirely different reasons than its failure for ill-behaved topological spaces.
These models are ``non-hypercomplete $\infty$-toposes''
\index{.infinity1-topos@$(\infty,1)$-topos}%
\index{.infinity1-topos@$(\infty,1)$-topos!non-hypercomplete}%
(see~\cite{lurie:higher-topoi}); roughly speaking, they consist of sheaves of $\infty$-groupoids over $\infty$-dimensional base spaces.

\index{axiom!Whitehead's principle|(}%
\index{Whitehead's!principle|(}%

From a foundational point of view, therefore, we may speak of \emph{Whitehead's principle} as a ``classicality axiom'', akin to \LEM{} and \choice{}.
It may consistently be assumed, but it is not part of the computationally motivated type theory, nor does it hold in all natural models.
But when working from set-theoretic foundations, this principle is invisible: it cannot fail to be true in a world where $\infty$-groupoids are built up out of sets (using topological spaces, simplicial sets, or any other such model).

This may seem odd, but actually it should not be surprising.
Homotopy type theory is the \emph{abstract} theory of homotopy types, whereas the homotopy theory of topological spaces or simplicial sets in set theory is a \emph{concrete} model of this theory, in the same way that the integers are a concrete model of the abstract theory of rings.
It is to be expected that any concrete model will have special properties which are not intrinsic to the corresponding abstract theory, but which we might sometimes want to assume as additional axioms (e.g.\ the integers are a Principal Ideal Domain, but not all rings are).

It is beyond the scope of this book to describe any models of type theory, so we will not explain how Whitehead's principle might fail in some of them.
However, we can prove that it holds whenever the types involved are $n$-truncated for some finite $n$, by ``downward'' induction on $n$.
In addition to being of interest in its own right (for instance, it implies the essential uniqueness of Eilenberg--Mac Lane spaces), the proof of this result will hopefully provide some intuitive explanation for why we cannot hope to prove an analogous theorem without truncation hypotheses.

We begin with the following modification of \autoref{thm:mono-surj-equiv}, which will eventually supply the induction step in the proof of the truncated Whitehead's principle.
It may be regarded as a type-theoretic, $\infty$-group\-oid\-al version of the classical statement that a fully faithful and essentially surjective functor is an equivalence of categories.

\begin{thm}\label{thm:whitehead0}
  Suppose $f:A\to B$ is a function such that
  \begin{enumerate}
  \item $\trunc0 f : \trunc0 A \to \trunc0 B$ is surjective, and\label{item:whitehead01}
  \item for any $x,y:A$, the function $\apfunc f : (\id[A]xy) \to (\id[B]{f(x)}{f(y)})$ is an equivalence.\label{item:whitehead02}
  \end{enumerate}
  Then $f$ is an equivalence.
\end{thm}
\begin{proof}
  Note that~\ref{item:whitehead02} is precisely the statement that $f$ is an embedding, c.f.~\autoref{sec:mono-surj}.
  Thus, by \autoref{thm:mono-surj-equiv}, it suffices to show that $f$ is surjective, i.e.\ that for any $b:B$ we have $\trunc{-1}{\hfib f b}$.
  Suppose given $b$; then since $\trunc0 f$ is surjective, there merely exists an $a:A$ such that $\trunc 0 f(\tproj0a) = \tproj0b$.
  And since our goal is a mere proposition, we may assume given such an $a$.
  Then we have $\tproj0{f(a)} = \trunc 0 f(\tproj0a) =\tproj0b$, hence $\trunc{-1}{f(a)=b}$.
  Again, since our goal is still a mere proposition, we may assume $f(a)=b$.
  Hence $\hfib f b$ is inhabited, and thus merely inhabited.
\end{proof}

Since homotopy groups are truncations of loop spaces\index{loop space}, rather than path spaces, we need to modify this theorem to speak about these instead.

\begin{cor}\label{thm:whitehead1}
  Suppose $f:A\to B$ is a function such that
  \begin{enumerate}
  \item $\trunc0 f : \trunc0 A \to \trunc0 B$ is a bijection, and
  \item for any $x:A$, the function $\apfunc f : \Omega(A,x) \to \Omega(B,f(x))$ is an equivalence.
  \end{enumerate}
  Then $f$ is an equivalence.
\end{cor}
\begin{proof}
  By \autoref{thm:whitehead0}, it suffices to show that $\apfunc f : (\id[A]xy) \to (\id[B]{f(x)}{f(y)})$ is an equivalence for any $x,y:A$.
  And by \autoref{thm:equiv-inhabcod}, we may assume $\id[B]{f(x)}{f(y)}$.
  In particular, $\tproj0{f(x)} = \tproj0{f(y)}$, so since $\trunc0 f$ is an equivalence, we have $\tproj0 x = \tproj0y$, hence $\tproj{-1}{x=y}$.
  Since we are trying to prove a mere proposition ($f$ being an equivalence), we may assume given $p:x=y$.
  But now the following square commutes up to homotopy:
  \begin{equation*}
  \vcenter{\xymatrix@C=3pc{
      \Omega(A,x)\ar[r]^-{\blank\ct p}\ar[d]_{\apfunc f} &
      (\id[A]xy) \ar[d]^{\apfunc f}\\
      \Omega(B,f(x))\ar[r]_-{\blank\ct f(p)} &
      (\id[B]{f(x)}{f(y)}).
      }}
  \end{equation*}
  The top and bottom maps are equivalences, and the left-hand map is so by assumption.
  Hence, by the 2-out-of-3 property, so is the right-hand map.
\end{proof}

Now we can prove the truncated Whitehead's principle.

\begin{thm}\label{thm:whiteheadn}
  Suppose $A$ and $B$ are $n$-types and $f:A\to B$ is such that
  \begin{enumerate}
  \item $\trunc0f:\trunc0A \to \trunc0B$ is an isomorphism, and\label{item:wh0}
  \item $\pi_k(f):\pi_k(A,x) \to \pi_k(B,f(x))$ is a bijection for all $k\ge 1$ and all $x:A$.\label{item:whk}
  \end{enumerate}
  Then $f$ is an equivalence.
\end{thm}

\noindent
Condition~\ref{item:wh0} is almost the case of~\ref{item:whk} when $k=0$, except that it makes no reference to any basepoint $x:A$.

\begin{proof}
  We proceed by induction on $n$.
  When $n=-2$, the statement is trivial.
  Thus, suppose it to be true for all functions between $n$-types, and let $A$ and $B$ be $(n+1)$-types and $f:A\to B$ as above.
  The first condition in \autoref{thm:whitehead1} holds by assumption, so it will suffice to show that for any $x:A$, the function $\apfunc f: \Omega(A,x) \to \Omega(B,f(x))$ is an equivalence.
  However, $\Omega(A,x)$ and $\Omega(B,f(x))$ are $n$-types, and $\pi_k(\apfunc f) = \pi_{k+1}(f)$, so this follows from the inductive hypothesis.
\end{proof}

Note that if $A$ and $B$ are not $n$-types for any finite $n$, then there is no way for the induction to get started.

\begin{cor}\label{thm:whitehead-contr}
  If $A$ is a $0$-connected $n$-type and $\pi_k(A,a)=0$ for all $k$ and $a:A$, then $A$ is contractible.
\end{cor}
\begin{proof}
  Apply \autoref{thm:whiteheadn} to the map $A\to\unit$.
\end{proof}

As an application, we can deduce the converse of \autoref{thm:conn-pik}.

\begin{cor}\label{thm:pik-conn}
  For $n\ge 0$, a map $f:A\to B$ is $n$-connected if and only if the following all hold:
  \begin{enumerate}
  \item $\trunc0f:\trunc0A \to \trunc0B$ is an isomorphism.
  \item For any $a:A$ and $k\le n$, the map $\pi_k(f):\pi_k(A,a) \to \pi_k(B,f(a))$ is an isomorphism.
  \item For any $a:A$, the map $\pi_{n+1}(f):\pi_{n+1}(A,a) \to \pi_{n+1}(B,f(a))$ is surjective.
  \end{enumerate}
\end{cor}
\begin{proof}
  The ``only if'' direction is \autoref{thm:conn-pik}.
  Conversely, by the long exact sequence of a fibration (\autoref{thm:les}),
  \index{sequence!exact}%
  the hypotheses imply that $\pi_k(\hfib f {f(a)})=0$ for all $k\le n$ and $a:A$, and that $\trunc 0{\hfib f{f(a)}}$ is contractible.
  Since $\pi_k(\hfib f {f(a)}) = \pi_k(\trunc n{\hfib f{f(a)}})$ for $k\le n$, and $\trunc n{\hfib f{f(a)}}$ is $n$-connected, by \autoref{thm:whitehead-contr} it is contractible for any $a$.

  It remains to show that $\trunc n{\hfib f{b}}$ is contractible for $b:B$ not necessarily of the form $f(a)$.
  However, by assumption, there is $x:\trunc0A$ with $\tproj 0b = \trunc0f(x)$.
  Since contractibility is a mere proposition, we may assume $x$ is of the form $\tproj0a$ for $a:A$, in which case $\tproj 0 b = \trunc0f(\tproj0a) = \tproj0{f(a)}$, and therefore $\trunc{-1}{b=f(a)}$.
  Again since contractibility is a mere proposition, we may assume $b=f(a)$, and the result follows.
\end{proof}

A map $f$ such that $\trunc0f$ is a bijection and $\pi_k(f)$ is a bijection for all $k$ is called \define{$\infty$-connected}%
\indexdef{function!.infinity-connected@$\infty$-connected}%
\indexdef{.infinity-connected function@$\infty$-connected function}
or a \define{weak equivalence}.%
\indexdef{equivalence!weak}%
\indexdef{weak equivalence!of types}
This is equivalent to asking that $f$ be $n$-connected for all $n$.
A type $Z$ is called \define{$\infty$-truncated}%
\indexdef{type!.infinity-truncated@$\infty$-truncated}%
\indexdef{.infinity-truncated type@$\infty$-truncated type}
or \define{hypercomplete}%
\indexdef{type!hypercomplete}%
\indexdef{hypercomplete type}
if the induced map
\[(\blank\circ f):(B\to Z) \to (A\to Z)\]
is an equivalence whenever $f$ is $\infty$-connected --- that is, if $Z$ thinks every $\infty$-connected map is an equivalence.
\indexdef{axiom!Whitehead's principle}%
Then if we want to assume Whitehead's principle as an axiom, we may use either of the following equivalent forms.
\begin{itemize}
\item Every $\infty$-connected function is an equivalence.
\item Every type is $\infty$-truncated.
\end{itemize}
In higher topos models,
\index{.infinity1-topos@$(\infty,1)$-topos}%
the $\infty$-truncated types form a reflective subuniverse in the sense of 
\autoref{sec:modalities} (the ``hypercompletion'' of an $(\infty,1)$-topos), but we do not know whether this is 
true in general.

\index{axiom!Whitehead's principle|)}%
\index{Whitehead's!principle|)}%

It may not be obvious that there \emph{are} any types which are not $n$-types for any $n$, but in fact there are.
Indeed, in classical homotopy theory, $\Sn^n$ has this property for any $n\ge 2$.
We have not proven this fact in homotopy type theory yet, but there are other types which we can prove to have ``infinite truncation level''.

\begin{eg}
  Suppose we have $B:\nat\to\type$ such that for each $n$, the type $B(n)$ contains an $n$-loop\index{loop!n-@$n$-} which is not equal to $n$-fold reflexivity, say $p_n:\Omega^n(B(n),b_n)$ with $p_n \neq \refl{b_n}^n$.
  (For instance, we could define $B(n)\defeq \Sn^n$, with $p_n$ the image of $1:\Z$ under the isomorphism $\pi_n(\Sn^n)\cong \Z$.)
  Consider $C\defeq \prd{n:\nat} B(n)$, with the point $c:C$ defined by $c(n)\defeq b_n$.
  Since loop spaces commute with products, for any $m$ we have
  \[\eqvspaced{\Omega^m (C,c)}{\prd{n:\nat}\Omega^m(B(n),b_n)}.\]
  Under this equivalence, $\refl{c}^m$ corresponds to the function $(n\mapsto \refl{b_n}^m)$.
  Now define $q_m$ in the right-hand type by
  \[ q_m(n) \defeq
  \begin{cases}
    p_n &\quad m=n\\
    \refl{b_n}^m &\quad m\neq n.
  \end{cases}
  \]
  If we had $q_m = (n\mapsto \refl{b_n}^m)$, then we would have $p_n = \refl{b_n}^n$, which is not the case.
  Thus, $q_m \neq (n\mapsto \refl{b_n}^m)$, and so there is a point of $\Omega^m(C,c)$ which is unequal to $\refl{c}^m$.
  Hence $C$ is not an $m$-type, for any $m:\nat$.
\end{eg}

We expect it should also be possible to show that a universe $\UU$ itself is not an $n$-type for any $n$, using the fact that it contains higher inductive types such as $\Sn^n$ for all $n$.
However, this has not yet been done.

\section{A general statement of the encode-decode method}
\label{sec:general-encode-decode}

\indexdef{encode-decode method}

We have used the encode-decode method to characterize the path spaces
of various types, including coproducts (\cref{thm:path-coprod}), natural
numbers (\cref{thm:path-nat}), truncations (\cref{thm:path-truncation}),
the circle (\cref{{cor:omega-s1}}), suspensions (\autoref{thm:freudenthal}), and pushouts
(\cref{thm:van-Kampen}).  Variants of this technique are used in the
proofs of many of the other theorems mentioned in the introduction to
this chapter, such as a direct proof of $\pi_n(\Sn^n)$, the Blakers--Massey theorem, and the construction of Eilenberg--Mac Lane spaces.
While it is tempting to try to
abstract the method into a lemma, this is difficult because
slightly different variants are needed for different problems.  For
example, different variations on the same method  can be used to
characterize a loop space (as in \cref{thm:path-coprod,cor:omega-s1}) or
a whole path space (as in \cref{thm:path-nat}), to give a complete
characterization of a loop space (e.g.\ $\Omega^1(\Sn ^1)$) or only to
characterize some truncation of it (e.g.\ van Kampen), and to calculate
homotopy groups or to prove that a map is $n$-connected (e.g.\ Freudenthal and
Blakers--Massey).

However, we can state lemmas for specific variants of the method.
The proofs of these lemmas are almost trivial; the main point is to
clarify the method by stating them in generality.  The simplest
case is using an encode-decode method to characterize the loop space of a
type, as in \cref{thm:path-coprod} and \cref{{cor:omega-s1}}.

\begin{lem}[Encode-decode for Loop Spaces]
  \index{loop space}%
Given a pointed type $(A,a_0)$ and a fibration
$\code : A \to \type$, if 
\begin{enumerate}
\item $c_0 : \code(a_0)$,\label{item:ed1}
\item $\decode : \prd{x:A} \code(x) \to (\id{a_0}{x})$,\label{item:ed2}
\item for all $c : \code(a_0)$, \id{\transfib{\code}{\decode(c)}{c_0}}{c}, and\label{item:ed3}
\item $\id{\decode(c_0)}{\refl{}}$,\label{item:ed4}
\end{enumerate}
then $(\id{a_0}{a_0})$ is equivalent to $\code(a_0)$.
\end{lem}

\begin{proof}
Define
$\encode : \prd{x:A} (\id{a_0}{x}) \to \code(x)$ by
\[
\encode_x(\alpha) = \transfib{\code}{\alpha}{c_0}.
\]
We show that $\encode_{a_0}$ and $\decode_{a_0}$ are quasi-inverses.  
The composition $\encode_{a_0} \circ \decode_{a_0}$ is immediate by
assumption~\ref{item:ed3}.  For the other composition, we show
\[
\prd{x:A}{p : \id{a_0}{x}} \id{\decode_{x} (\encode_{x} p)}{p}.
\] 
By path induction, it suffices to show 
$\id{{\decode_{{a_0}} (\encode_{{a_o}} \refl{})}}{\refl{}}$.
After reducing the $\mathsf{transport}$, it suffices to show 
$\id{{\decode_{{a_0}} (c_0)}}{\refl{}}$, which is assumption~\ref{item:ed4}.
\end{proof}

If a fiberwise equivalence between $(\id{a_0}{\blank})$ and $\code$ is desired,
it suffices to strengthen condition (iii) to
\[
\prd{x:A}{c : \code(x)} \id{\encode_{x}(\decode_{x}(c))}{c}.a
\]
However, to calculate a loop space (e.g. $\Omega(\Sn ^1)$), this
stronger assumption is not necessary.  

Another variation, which comes up often when calculating homotopy
groups, characterizes the truncation of a loop space:

\begin{lem}[Encode-decode for Truncations of Loop Spaces]
Assume a pointed type $(A,a_0)$ and a fibration
$\code : A \to \type$, where for every $x$, $\code(x)$ is a $k$-type.
Define 
\[
\encode : \prd{x:A} \trunc{k}{\id{a_0}{x}} \to \code(x).
\]
by truncation recursion (using the fact
that $\code(x)$ is a $k$-type), mapping $\alpha : \id{a_0}{x}$ to 
\transfib{\code}{\alpha}{c_0}. Suppose:
\begin{enumerate}
\item $c_0 : \code(a_0)$,
\item $\decode : \prd{x:A} \code(x) \to \trunc{k}{\id{a_0}{x}}$,
\item \label{item:decode-encode-loop-iii}
  $\id{\encode_{a_0}(\decode_{a0}(c))}{c}$ for all $c : \code(a_0)$, and
\item \label{item:decode-encode-loop-iv}
  $\id{\decode(c_0)}{\tproj{}{\refl{}}}$.
\end{enumerate}
Then $\trunc{k}{\id{a_0}{a_0}}$ is equivalent to $\code(a_0)$.
\end{lem}

\begin{proof}
That $\decode \circ \encode$ is identity is immediate by \ref{item:decode-encode-loop-iii}.
To prove $\encode \circ \decode$, we first do a truncation induction, by
which it suffices to show
\[
\prd{x:A}{p : \id{a_0}{x}} \id{\decode_{x}(\encode_{x}(\tproj{k}{p}))}{\tproj{k}{p}}.
\] 
The truncation induction is allowed because paths in a $k$-type are a
$k$-type.  To show this type, we do a path induction, and after reducing
the \encode, use assumption~\ref{item:decode-encode-loop-iv}.
\end{proof}

\section{Additional Results}
\label{sec:moreresults}

Though we do not present the proofs in this chapter, following results have also been established in homotopy type theory.  

\begin{thm}
\index{homotopy!group!of sphere}
There exists a $k$ such that for all $n \ge 3$, $\pi_{n+1}(\Sn ^n) =
\Z_k$.  
\end{thm}

\begin{proof}[Notes on the proof.]
The proof consists of a calculation of $\pi_4(\Sn ^3)$, together with an
appeal to stability (\cref{cor:stability-spheres}).  In the classical
statement of this result, $k$ is $2$.  While we have not yet checked that
$k$ is in fact $2$, our calcluation of $\pi_4(\Sn ^3)$ is constructive,\index{mathematics!constructive}
like all the rest of the proofs in this chapter.
(More precisely, it doesn't use any additional axioms such as \LEM{} or \choice{}, making it as constructive as
univalence and higher inductive types are.)  Thus, given a
computational interpretation of homotopy type theory, we could run the
proof on a computer to verify that $k$ is $2$.  This example is quite
intriguing, because it is the first calculation of a homotopy group
for which we have not needed to know the answer in advance.
\end{proof}

\index{pushout}%

\begin{thm}[Blakers--Massey theorem]\label{Blakers-Massey}
  \indexdef{theorem!Blakers--Massey}%
  \indexsee{Blakers--Massey theorem}{theorem, Blakers--Massey}%
  Suppose we are given maps $f : C  \rightarrow X$, and $g : C \rightarrow Y$. Taking first the pushout $X \sqcup^C Y $ of $f$ and $g$ and then the pullback of its inclusions $\inl : X \rightarrow X \sqcup^C Y \leftarrow Y : \inr$, we have an induced map $C \to X \times_{(X \sqcup^C Y)} Y$.

  If $f$ is $i$-connected and $g$ is $j$-connected, then this induced map is $(i+j)$-connected. In other words, for any points $x:X$, $y:Y$, the corresponding fiber $C_{x,y}$ of $(f,g) : C \to X \times Y $ gives an approximation to the path space $\id[X \sqcup^C Y]{\inl(x)}{\inr(y)}$ in the pushout.
\end{thm}

It should be noted that in classical algebraic topology, the Blakers--Massey theorem is often stated in a somewhat different form, where the maps $f$ and $g$ are replaced by inclusions of subcomplexes of CW complexes, and the homotopy pushout and homotopy pullback by a union and intersection, respectively.
In order to express the theorem in homotopy type theory, we have to replace notions of this sort with ones that are homotopy-invariant.
We have seen another example of this in the van Kampen theorem (\autoref{sec:van-kampen}), where we had to replace a union of open subsets by a homotopy pushout.

\begin{thm}[Eilenberg--Mac Lane Spaces]\label{Eilenberg-Mac-Lane-Spaces}
\index{Eilenberg--Mac Lane space}
For any abelian\index{group!abelian} group $G$ and positive integer $n$, there is an $n$-type
$K(G,n)$ such that $\pi_n(K(G,n)) = G$, and  $\pi_k(K(G,n)) = 0$
for $k\neq n$.
\end{thm}

\begin{thm}[Covering spaces]\label{thm:covering-spaces}
  \index{covering space}%
  For a connected space $A$, there is an equivalence between covering spaces over $A$ and sets with an action of $\pi_1(A)$.
\end{thm}

\sectionNotes


{
\newcommand{\humancheck}{\ding{52}}
\newcommand{\computercheck}{\ding{52}\kern-0.5em\ding{52}}
\begin{table}[htb]
  \centering
\begin{tabular}{lcc}
\toprule
Theorem         & Status \\
\midrule
$\pi_1(\Sn ^1)$                     & \computercheck \\
$\pi_{k<n}(\Sn ^n)$                  & \computercheck \\
long-exact-sequence of homotopy groups & \computercheck    \\
total space of Hopf fibration is $\Sn ^3$ & \humancheck    \\
$\pi_2(\Sn ^2)$                     & \computercheck \\
$\pi_3(\Sn ^2)$                     & \humancheck    \\
$\pi_n(\Sn ^n)$                     & \computercheck \\
$\pi_4(\Sn ^3)$                     & \humancheck    \\
Freudenthal suspension theorem      & \computercheck \\
Blakers--Massey theorem              & \computercheck \\
Eilenberg--Mac Lane spaces $K(G,n)$ & \computercheck \\
van Kampen theorem                & \computercheck \\
covering spaces                     & \computercheck \\
Whitehead's principle for $n$-types & \computercheck \\
\bottomrule
\end{tabular}
\caption{Theorems from homotopy theory proved by
  hand (\humancheck) and by computer (\computercheck).}
  \label{tab:theorems}
\end{table}
}

The theorems described in this chapter are standard results in classical
homotopy theory; many are described by \cite{hatcher02topology}.  In these
notes, we review the development of the new synthetic proofs of them in homotopy
type theory.  \autoref{tab:theorems} lists the homotopy-theoretic
theorems that have been proven in homotopy type theory, and whether they
have been computer-checked.
Almost all of these results were developed during the spring term at IAS
in 2013, as part of a significant collaborative effort.  Many people
contributed to these results, for example by being the principal author
of a proof, by suggesting problems to work on, by participating in many
discussions and seminars about these problems, or by giving feedback on
results.  The following people were the principal authors of the first
homotopy type theory proofs of the above theorems. Unless indicated otherwise, for the
theorems that have been computer-checked, the principal authors were
also the first ones to formalize\index{mathematics!formalized} the proof using a computer proof
assistant.
\begin{itemize}
\item 
Shulman gave the homotopy-theoretic calculation of $\pi_1(\Sn^1)$.  Licata later discovered the
encode-decode proof and the encode-decode method.

\item 
Brunerie calculated $\pi_{k<n}(\Sn ^ n)$.
Licata later gave an encode-decode version.

\item Voevodsky constructed the long exact sequence of homotopy groups. 

\item Lumsdaine constructed the Hopf fibration.
  Brunerie proved that its total space is $\Sn ^3$, thereby calculating $\pi_2(\Sn ^2)$ and
  $\pi_3(\Sn ^3)$.

\item Licata and Brunerie gave a direct calculation of
$\pi_n(\Sn ^n)$.  

\item 
  Lumsdaine proved the Freudenthal suspension theorem; Licata and
  Lumsdaine formalized this proof.
\item Lumsdaine, Finster, and Licata proved the Blakers--Massey theorem;
  Lumsdaine, Brunerie, Licata, and Hou formalized it.

\item 
Licata gave an encode-decode calculation of $\pi_2(\Sn ^2)$, and a
calculation of $\pi_n(\Sn ^n)$ using the Freudenthal suspension theorem; using similar
techniques, he constructed $K(G,n)$.

\item 
Shulman proved the van Kampen theorem; Hou formalized this proof.

\item 
Licata proved Whitehead's theorem for $n$-types.

\item Brunerie calculated $\pi_4(\Sn ^3)$.

\item 
Hou established the theory of covering spaces and formalized it.
\end{itemize}

The interplay between homotopy theory and type theory was crucial to the
development of these results.  For example, the first proof that
$\pi_1(\Sn ^1)=\mathbb{Z}$ was the one given in \autoref{subsec:pi1s1-homotopy-theory}, which follows a classical homotopy theoretic one.  A
type-theoretic analysis of this proof resulted in the development of the
encode-decode method.  The first calculation of $\pi_2(\Sn ^2)$ also followed
classical methods, but this led quickly to an encode-decode proof of the
result.  The encode-decode calculation generalized to $\pi_n(\Sn ^n)$, which
in turn led to the proof of the Freudenthal suspension theorem, by
combining an encode-decode argument with classical homotopy-theoretic
reasoning about connectedness, which in turn led to the Blakers--Massey
theorem and Eilenberg--Mac Lane spaces.  The rapid development of this
series of results illustrates the promise of our new understanding of
the connections between these two subjects.

\sectionExercises

\begin{ex}
  Prove that homotopy groups respect products: $\eqv{\pi_n(A\times B)}{\pi_n(A)\times \pi_n(B)}$.
\end{ex}

\begin{ex}
  \index{decidable!equality}%
  Prove that if $A$ is a set with decidable equality (see \autoref{defn:decidable-equality}), then its suspension $\susp A$ is a 1-type.
  (It is an open question\index{open!problem} whether this is provable without the assumption of decidable equality.)
\end{ex}

\begin{ex}
  Define $\Sn^\infty$ to be the colimit of the sequence $\Sn^0 \to \Sn^1 \to \Sn^2 \to\cdots$.
  Prove that $\Sn^\infty$ is contractible.
\end{ex}

\begin{ex}
  Define $\Sn^\infty$ to be the higher inductive type generated by
  \begin{itemize}
  \item Two points $\north:\Sn^\infty$ and $\south:\Sn^\infty$, and
  \item For each $x:\Sn^\infty$, a path $\merid(x):\north=\south$.
  \end{itemize}
  In other words, $\Sn^\infty$ is its own suspension.
  Prove that $\Sn^\infty$ is contractible.
\end{ex}

\begin{ex}\label{ex:unique-fiber}
  Suppose $f:X\to Y$ is a function and $Y$ is connected.
  Show that for any $y_1,y_2:Y$ we have $\brck{\eqv{\hfib f {y_1}}{\hfib f {y_2}}}$.
\end{ex}

\begin{ex}\label{ex:ap-path-inversion}
  For any pointed type $A$, let $i_A : \Omega A \to \Omega A$ denote inversion of loops, $i_A \defeq \lam{p} \rev{p}$.
  Show that $i_{\Omega A} : \Omega^2 A \to \Omega^2 A$ is equal to $\Omega(i_A)$.
\end{ex}

\begin{ex}\label{ex:pointed-equivalences}
  Define a \textbf{pointed equivalence} to be a pointed map whose underlying function is an equivalence.
  \begin{enumerate}
  \item Show that the type of pointed equivalences between pointed types $(X,x_0)$ and $(Y,y_0)$ is equivalent to $\id[\pointed\type]{(X,x_0)}{(Y,y_0)}$.
  \item Reformulate the notion of pointed equivalence in terms of a pointed quasi-inverse and pointed homotopies, in one of the coherent styles from \autoref{cha:equivalences}.
  \end{enumerate}
\end{ex}

\begin{ex}\label{ex:HopfJr}
  Following the example of the Hopf fibration in \autoref{sec:hopf}, define the \define{junior Hopf fibration}
  \indexdef{Hopf!fibration!junior}%
as a fibration (that is, a type family) over $\Sn ^1$ whose fiber over the basepoint is $\Sn ^0$ and whose total space is $\Sn ^1$.  This is also called the ``twisted double cover'' of the circle $\Sn ^1$.
\end{ex}

\begin{ex}\label{ex:SuperHopf}
Again following the example of the Hopf fibration in \autoref{sec:hopf}, define an analogous fibration over $\Sn ^4$ whose fiber over the basepoint is $\Sn ^3$ and whose total space is $\Sn ^7$.  This is an open problem\index{open!problem} in homotopy type theory (such a fibration is known to exist in classical homotopy theory).
\end{ex}

\begin{ex}\label{ex:vksusppt}
  Continuing from \autoref{eg:suspension}, prove that if $A$ has a point $a:A$, then we can identify $\pi_1(\susp A)$ with the group presented by $\trunc0A$ as generators with the relation $\tproj0a = e$.
  Then show that if we assume excluded middle, this is also the free group on $\trunc0 A \setminus \{\tproj0 a \}$.
\end{ex}

\begin{ex}\label{ex:vksuspnopt}
  Again continuing from \autoref{eg:suspension}, but this time without assuming $A$ to be pointed, show that we can identify $\pi_1(\susp A)$ with the group presented by generators $\trunc0A \times \trunc0A$ and relations
  \begin{equation*}
    (a,b) = \opp{(b,a)},
    \qquad
    (a,c) = (a,b)\cdot (b,c),
    \qquad\text{and}\qquad
    (a,a) = e.
  \end{equation*}
\end{ex}

\index{acceptance|)}


\chapter{Category theory}
\label{cha:category-theory}

Of the branches of mathematics, category theory is one which perhaps fits the least comfortably in set theoretic foundations.
One problem is that most of category theory is invariant under weaker notions of ``sameness'' than equality, such as isomorphism in a category or equivalence of categories, in a way which set theory fails to capture.
But this is the same sort of problem that the univalence axiom solves for types, by identifying equality with equivalence.
Thus, in univalent foundations it makes sense to consider a notion of ``category'' in which equality of objects is identified with isomorphism in a similar way.

Ignoring size issues, in set-based mathematics a category consists of a \emph{set} $A_0$ of objects and, for each $x,y\in A_0$, a \emph{set} $\hom_A(x,y)$ of morphisms.
Under univalent foundations, a ``naive'' definition of category would simply mimic this with a \emph{type} of objects and \emph{types} of morphisms.
If we allowed these types to contain arbitrary higher homotopy, then we ought to impose higher coherence conditions, leading to some notion of $(\infty,1)$-category,
\index{.infinity1-category@$(\infty,1)$-category}%
but at present our goal is more modest.
We consider only 1-categories, and therefore we restrict the types $\hom_A(x,y)$ to be sets, i.e.\ 0-types.
If we impose no further conditions, we will call this notion a \emph{precategory}.

If we add the requirement that the type $A_0$ of objects is a set, then we end up with a definition that behaves much like the traditional set-theoretic one.
Following Toby Bartels, we call this notion a \emph{strict category}.
\index{strict!category}%
But we can also require a generalized version of the univalence axiom, identifying $(x=_{A_0} y)$ with the type $\mathsf{iso}(x,y)$ of isomorphisms from $x$ to $y$.
Since we regard this as usually the ``correct'' definition, we will call it simply a \emph{category}.

A good example of the difference between the three notions of category is provided by the statement ``every fully faithful and essentially surjective functor is an equivalence of categories'', which in classical set-based category theory is equivalent to the axiom of choice.
\index{mathematics!classical}%
\index{axiom!of choice}%
\index{classical!category theory}%
\begin{enumerate}
\item For strict categories, this is still equivalent to to the axiom of choice.
\item For precategories, there is no consistent axiom of choice which can make it true.
\item For categories, it is provable \emph{without} any axiom of choice.
\end{enumerate}
We will prove the latter statement in this chapter, as well as other pleasant properties of categories, e.g.\ that equivalent categories are equal (as elements of the type of categories).
We will also describe a universal way of ``saturating'' a precategory $A$ into a category $\widehat A$, which we call its \emph{Rezk completion},
\index{completion!Rezk}%
although it could also reasonably be called the \emph{stack completion} (see the Notes).

The Rezk completion also sheds further light on the notion of equivalence of categories.
For instance, the functor $A \to \widehat{A}$ is always fully faithful and essentially surjective, hence a ``weak equivalence''.
It follows that a precategory is a category exactly when it ``sees'' all fully faithful and essentially surjective functors as equivalences; thus our notion of ``category'' is already inherent in the notion of ``fully faithful and essentially surjective functor''.

We assume the reader has some basic familiarity with classical category theory.\index{classical!category theory}
Recall that whenever we write \type it denotes some universe of types, but perhaps a different one at different times; everything we say remains true for any consistent choice of universe levels\index{universe level}.
We will use the basic notions of homotopy type theory from \autoref{cha:typetheory,cha:basics} and the propositional truncation from \autoref{cha:logic}, but not much else from \autoref{part:foundations}, except that our second construction of the Rezk completion will use a higher inductive type.

\section{Categories and precategories}
\label{sec:cats}

In classical mathematics, there are many equivalent definitions of a category.
In our case, since we have dependent types, it is natural to choose the arrows to be a type family indexed by the objects.
This matches the way hom-types are always used in category theory: we never even consider comparing two arrows unless we know their domains and codomains agree.
Furthermore, it seems clear that for a theory of 1-categories, the hom-types should all be sets.
This leads us to the following definition.

\begin{defn}\label{ct:precategory}
  A \define{precategory}
  \indexdef{precategory}
  $A$ consists of the following.
  \begin{enumerate}
  \item A type $A_0$ of \define{objects}.%
    \indexdef{object!in a (pre)category}
    We write $a:A$ for $a:A_0$.
  \item For each $a,b:A$, a set $\hom_A(a,b)$ of \define{arrows} or \define{morphisms}.%
    \indexsee{arrow}{morphism}%
    \indexdef{morphism!in a (pre)category}%
    \indexdef{hom-set}%
  \item For each $a:A$, a morphism $1_a:\hom_A(a,a)$.%
    \indexdef{identity!morphism in a (pre)category}
  \item For each $a,b,c:A$, a function%
    \indexdef{composition!of morphisms in a (pre)category}
    \[  \hom_A(b,c) \to \hom_A(a,b) \to \hom_A(a,c) \]
    denoted infix by $g\mapsto f\mapsto g\circ f$, or sometimes simply by $gf$.
  \item For each $a,b:A$ and $f:\hom_A(a,b)$, we have $\id f {1_b\circ f}$ and $\id f {f\circ 1_a}$.
  \item For each $a,b,c,d:A$ and 
    \begin{equation*}
      f:\hom_A(a,b), \qquad
      g:\hom_A(b,c), \qquad
      h:\hom_A(c,d),
    \end{equation*}
    we have $\id {h\circ (g\circ f)}{(h\circ g)\circ f}$.
  \end{enumerate}
\end{defn}

The problem with the notion of precategory is that for objects $a,b:A$, we have two possibly-different notions of ``sameness''.
On the one hand, we have the type $(\id[A_0]{a}{b})$.
But on the other hand, there is the standard categorical notion of \emph{isomorphism}.

\begin{defn}\label{ct:isomorphism}
  A morphism $f:\hom_A(a,b)$ is an \define{isomorphism}
  \indexdef{isomorphism!in a (pre)category}%
  if there is a morphism $g:\hom_A(b,a)$ such that $\id{g\circ f}{1_a}$ and $\id{f\circ g}{1_b}$.
  We write $a\cong b$ for the type of such isomorphisms.
\end{defn}

\begin{lem}\label{ct:isoprop}
  For any $f:\hom_A(a,b)$, the type ``$f$ is an isomorphism'' is a mere proposition.
  Therefore, for any $a,b:A$ the type $a\cong b$ is a set.
\end{lem}
\begin{proof}
  Suppose given $g:\hom_A(b,a)$ and $\eta:(\id{1_a}{g\circ f})$ and $\epsilon:(\id{f\circ g}{1_b})$, and similarly $g'$, $\eta'$, and $\epsilon'$.
We must show $\id{(g,\eta,\epsilon)}{(g',\eta',\epsilon')}$.
  But since all hom-sets are sets, their identity types are mere propositions, so it suffices to show $\id g {g'}$.
  For this we have
  \[g' = 1_a\circ g' = (g\circ f)\circ g' = g\circ (f\circ g') = g\circ 1_b = g\]
  using $\eta$ and $\epsilon'$.
\end{proof}

\symlabel{ct:inv}
\index{inverse!in a (pre)category}%
If $f:a\cong b$, then we write $\inv f$ for its inverse, which by \autoref{ct:isoprop} is uniquely determined.

The only relationship between these two notions of sameness that we have in a precategory is the following.

\begin{lem}[\textsf{idtoiso}]\label{ct:idtoiso}
  If $A$ is a precategory and $a,b:A$, then
  \[(\id a b)\to (a \cong b).\]
\end{lem}
\begin{proof}
  By induction on identity, we may assume $a$ and $b$ are the same.
  But then we have $1_a:\hom_A(a,a)$, which is clearly an isomorphism.
\end{proof}

Evidently, this situation is analogous to the issue that motivated us to introduce the univalence axiom.
In fact, we have the following:

\begin{eg}\label{ct:precatset}
  \index{set}%
  There is a precategory \uset, whose type of objects is \set, and with $\hom_{\uset}(A,B) \defeq (A\to B)$.
  The identity morphisms are identity functions and the composition is function composition.
  For this precategory, \autoref{ct:idtoiso} is equal to (the restriction to sets of) the map $\idtoeqv$ from \autoref{sec:compute-universe}.

  Of course, to be more precise we should call this category $\uset_\UU$, since its objects are only the \emph{small sets}
  \index{small!set}%
  relative to a universe \UU.
\end{eg}

Thus, it is natural to make the following definition.

\begin{defn}\label{ct:category}
  A \define{category}
  \indexdef{category}
  is a precategory such that for all $a,b:A$, the function $\idtoiso_{a,b}$ from \autoref{ct:idtoiso} is an equivalence.
\end{defn}

In particular, in a category, if $a\cong b$, then $a=b$.

\begin{eg}\label{ct:eg:set}
  \index{univalence axiom}%
  The univalence axiom implies immediately that \uset is a category.
  One can also show, using univalence, that any precategory of set-level structures such as groups, rings, topological spaces, etc.\ is a category; see \autoref{sec:sip}.
\end{eg}

We also note the following.

\begin{lem}\label{ct:obj-1type}
  In a category, the type of objects is a 1-type.
\end{lem}
\begin{proof}
  It suffices to show that for any $a,b:A$, the type $\id a b$ is a set.
  But $\id a b$ is equivalent to $a \cong b$, which is a set.
\end{proof}

\symlabel{isotoid}
We write $\isotoid$ for the inverse $(a\cong b) \to (\id a b)$ of the map $\idtoiso$ from \autoref{ct:idtoiso}.
The following relationship between the two is important.

\begin{lem}\label{ct:idtoiso-trans}
  For $p:\id a a'$ and $q:\id b b'$ and $f:\hom_A(a,b)$, we have
  \begin{equation}\label{ct:idtoisocompute}
    \id{\trans{(p,q)}{f}}
    {\idtoiso(q)\circ f \circ \inv{\idtoiso(p)}}.
  \end{equation}
\end{lem}
\begin{proof}
  By induction, we may assume $p$ and $q$ are $\refl a$ and $\refl b$ respectively.
Then the left-hand side of~\eqref{ct:idtoisocompute} is simply $f$.
  But by definition, $\idtoiso(\refl a)$ is $1_a$, and $\idtoiso(\refl b)$ is $1_b$, so the right-hand side of~\eqref{ct:idtoisocompute} is $1_b\circ f\circ 1_a$, which is equal to $f$.
\end{proof}

Similarly, we can show
\begin{gather}
  \id{\idtoiso(\rev p)}{\inv {(\idtoiso(p))}}\\
  \id{\idtoiso(p\ct q)}{\idtoiso(q)\circ \idtoiso(p)}\\
  \id{\isotoid(f\circ e)}{\isotoid(e)\ct \isotoid(f)}
\end{gather}
and so on.

\begin{eg}\label{ct:orders}
  A precategory in which each set $\hom_A(a,b)$ is a mere proposition is equivalently a type $A_0$ equipped with a mere relation ``$\le$'' that is reflexive ($a\le a$) and transitive (if $a\le b$ and $b\le c$, then $a\le c$).
  We call this a \define{preorder}.
  \indexdef{preorder}

  In a preorder, a witness $f: a\le b$ is an isomorphism just when there exists some witness $g: b\le a$.
  Thus, $a\cong b$ is the mere proposition that $a\le b$ and $b\le a$.
  Therefore, a preorder $A$ is a category just when (1) each type $a=b$ is a mere proposition, and (2) for any $a,b:A_0$ there exists a function $(a\cong b) \to (a=b)$.
  In other words, $A_0$ must be a set, and $\le$ must be antisymmetric\index{relation!antisymmetric} (if $a\le b$ and $b\le a$, then $a=b$).
  We call this a \define{(partial) order} or a \define{poset}.
  \indexdef{partial order}%
  \indexdef{poset}%
\end{eg}

\begin{eg}\label{ct:gaunt}
  If $A$ is a category, then $A_0$ is a set if and only if for any $a,b:A_0$, the type $a\cong b$ is a mere proposition.
  This is equivalent to saying that every isomorphism in $A$ is an identity; thus it is rather stronger than the classical\index{mathematics!classical} notion of ``skeletal'' category.
  Categories of this sort are sometimes called \define{gaunt}~\cite{bsp12infncats}.
  \indexdef{category!gaunt}%
  \indexdef{gaunt category}%
  \index{skeletal category}%
  \index{category!skeletal}%
  There is not really any notion of ``skeletality'' for our categories, unless one considers \autoref{ct:category} itself to be such.
\end{eg}

\begin{eg}\label{ct:discrete}
  For any 1-type $X$, there is a category with $X$ as its type of objects and with $\hom(x,y) \defeq (x=y)$.
  If $X$ is a set, we call this the \define{discrete}
  \indexdef{category!discrete}%
  \indexdef{discrete!category}%
  category on $X$.
  In general, we call it a \define{groupoid}
  \indexdef{groupoid}
  (see \autoref{ct:groupoids}).
\end{eg}

\begin{eg}\label{ct:fundgpd}
  For \emph{any} type $X$, there is a precategory with $X$ as its type of objects and with $\hom(x,y) \defeq \pizero{x=y}$.
  The composition operation
  \[ \pizero{y=z} \to \pizero{x=y} \to \pizero{x=z} \]
  is defined by induction on truncation from concatenation $(y=z)\to(x=y)\to(x=z)$.
  We call this the \define{fundamental pregroupoid}
  \indexdef{fundamental!pregroupoid}%
  \indexsee{pregroupoid, fundamental}{fundamental pregroupoid}%
  of $X$.
  (In fact, we have met it already in \autoref{sec:van-kampen}; see also \autoref{ex:rezk-vankampen}.)
\end{eg}

\begin{eg}\label{ct:hoprecat}
  There is a precategory whose type of objects is \type and with $\hom(X,Y) \defeq \pizero{X\to Y}$, and composition defined by induction on truncation from ordinary composition $(Y\to Z) \to (X\to Y) \to (X\to Z)$.
  We call this the \define{homotopy precategory of types}.
  \indexdef{precategory!of types}%
  \index{homotopy!category of types@(pre)category of types}%
\end{eg}

\begin{eg}\label{ct:rel}
  Let \urel be the following precategory:
  \begin{itemize}
  \item Its objects are sets.
  \item $\hom_{\urel}(X,Y) = X\to Y\to \prop$.
  \item For a set $X$, we have $1_X(x,x') \defeq (x=x')$.
  \item For $R:\hom_{\urel}(X,Y)$ and $S:\hom_{\urel}(Y,Z)$, their composite is defined by
    \[ (S\circ R)(x,z) \defeq \Brck{\sm{y:Y} R(x,y) \times S(y,z)}.\]
  \end{itemize}
  Suppose $R:\hom_{\urel}(X,Y)$ is an isomorphism, with inverse $S$.
  We observe the following.
  \begin{enumerate}
  \item If $R(x,y)$ and $S(y',x)$, then $(R\circ S)(y',y)$, and hence $y'=y$.
    Similarly, if $R(x,y)$ and $S(y,x')$, then $x=x'$.\label{item:rel1}
  \item For any $x$, we have $x=x$, hence $(S\circ R)(x,x)$.
    Thus, there merely exists a $y:Y$ such that $R(x,y)$ and $S(y,x)$.\label{item:rel2}
  \item Suppose $R(x,y)$.
    By~\ref{item:rel2}, there merely exists a $y'$ with $R(x,y')$ and $S(y',x)$.
    But then by~\ref{item:rel1}, merely $y'=y$, and hence $y'=y$ since $Y$ is a set.
    Therefore, by transporting $S(y',x)$ along this equality, we have $S(y,x)$.
    In conclusion, $R(x,y)\to S(y,x)$.
    Similarly, $S(y,x) \to R(x,y)$.\label{item:rel3}
  \item If $R(x,y)$ and $R(x,y')$, then by~\ref{item:rel3}, $S(y',x)$, so that by~\ref{item:rel1}, $y=y'$.
    Thus, for any $x$ there is at most one $y$ such that $R(x,y)$.
    And by~\ref{item:rel2}, there merely exists such a $y$, hence there exists such a $y$.
  \end{enumerate}
  In conclusion, if $R:\hom_{\urel}(X,Y)$ is an isomorphism, then for each $x:X$ there is exactly one $y:Y$ such that $R(x,y)$, and dually.
  Thus, there is a function $f:X\to Y$ sending each $x$ to this $y$, which is an equivalence; hence $X=Y$.
  With a little more work, we conclude that \urel is a category.
\end{eg}

We might now restrict ourselves to considering categories rather than precategories.
Instead, we will develop many concepts for precategories as well as categories, in order to emphasize how much better-behaved categories are, as compared both to precategories and to ordinary categories in classical\index{mathematics!classical} mathematics.

We will also see in \crefrange{sec:strict-categories}{sec:dagger-categories} that in slightly more exotic contexts, there are uses for certain kinds of precategories other than categories, each of which ``fixes'' the equality of objects in different ways.
This emphasizes the ``pre''-ness of precategories: they are the raw material out of which multiple important categorical structures can be defined.

\section{Functors and transformations}
\label{sec:transfors}

The following definitions are fairly obvious, and need no modification.

\begin{defn}\label{ct:functor}
  Let $A$ and $B$ be precategories.
  A \define{functor}
  \indexdef{functor}%
  $F:A\to B$ consists of
  \begin{enumerate}
  \item A function $F_0:A_0\to B_0$, generally also denoted $F$.
  \item For each $a,b:A$, a function $F_{a,b}:\hom_A(a,b) \to \hom_B(Fa,Fb)$, generally also denoted $F$.
  \item For each $a:A$, we have $\id{F(1_a)}{1_{Fa}}$.
  \item For each $a,b,c:A$ and $f:\hom_A(a,b)$ and $g:\hom_B(b,c)$, we have
    \[\id{F(g\circ f)}{Fg\circ Ff}.\]
  \end{enumerate}
\end{defn}

Note that by induction on identity, a functor also preserves \idtoiso.

\begin{defn}\label{ct:nattrans}
  For functors $F,G:A\to B$, a \define{natural transformation}
  \indexdef{natural!transformation}%
  \indexsee{transformation!natural}{natural transformation}%
  $\gamma:F\to G$ consists of
  \begin{enumerate}
  \item For each $a:A$, a morphism $\gamma_a:\hom_B(Fa,Ga)$ (the ``components'').
  \item For each $a,b:A$ and $f:\hom_A(a,b)$, we have $\id{Gf\circ \gamma_a}{\gamma_b\circ Ff}$ (the ``naturality axiom'').
  \end{enumerate}
\end{defn}

Since each type $\hom_B(Fa,Gb)$ is a set, its identity type is a mere proposition.
Thus, the naturality axiom is a mere proposition, so identity of natural transformations is determined by identity of their components.
In particular, for any $F$ and $G$, the type of natural transformations from $F$ to $G$ is again a set.

Similarly, identity of functors is determined by identity of the functions $A_0\to B_0$ and (transported along this) of the corresponding functions on hom-sets.

\begin{defn}\label{ct:functor-precat}
  \indexdef{precategory!of functors}%
  For precategories $A,B$, there is a precategory $B^A$ defined by
  \begin{itemize}
  \item $(B^A)_0$ is the type of functors from $A$ to $B$.
  \item $\hom_{B^A}(F,G)$ is the type of natural transformations from $F$ to $G$.
  \end{itemize}
\end{defn}
\begin{proof}
  We define $(1_F)_a\defeq 1_{Fa}$.
  Naturality follows by the unit axioms of a precategory.
  For $\gamma:F\to G$ and $\delta:G\to H$, we define $(\delta\circ\gamma)_a\defeq \delta_a\circ \gamma_a$.
  Naturality follows by associativity.
  Similarly, the unit and associativity laws for $B^A$ follow from those for $B$.
\end{proof}

\begin{lem}\label{ct:natiso}
  \index{natural!isomorphism}%
  \index{isomorphism!natural}%
  A natural transformation $\gamma:F\to G$ is an isomorphism in $B^A$ if and only if each $\gamma_a$ is an isomorphism in $B$.
\end{lem}
\begin{proof}
  If $\gamma$ is an isomorphism, then we have $\delta:G\to F$ that is its inverse.
  By definition of composition in $B^A$, $(\delta\gamma)_a\jdeq \delta_a\gamma_a$ and similarly.
  Thus, $\id{\delta\gamma}{1_F}$ and $\id{\gamma\delta}{1_G}$ imply $\id{\delta_a\gamma_a}{1_{Fa}}$ and $\id{\gamma_a\delta_a}{1_{Ga}}$, so $\gamma_a$ is an isomorphism.

  Conversely, suppose each $\gamma_a$ is an isomorphism, with inverse called $\delta_a$, say.
We define a natural transformation $\delta:G\to F$ with components $\delta_a$; for the naturality axiom we have
  \[ Ff\circ \delta_a = \delta_b\circ \gamma_b\circ Ff \circ \delta_a = \delta_b\circ Gf\circ \gamma_a\circ \delta_a = \delta_b\circ Gf. \]
  Now since composition and identity of natural transformations is determined on their components, we have $\id{\gamma\delta}{1_G}$ and $\id{\delta\gamma}{1_F}$.
\end{proof}

The following result is fundamental.

\begin{thm}\label{ct:functor-cat}
  \indexdef{category!of functors}%
  \indexdef{functor!category of}%
  If $A$ is a precategory and $B$ is a category, then $B^A$ is a category.
\end{thm}
\begin{proof}
  Let $F,G:A\to B$; we must show that $\idtoiso:(\id{F}{G}) \to (F\cong G)$ is an equivalence.

  To give an inverse to it, suppose $\gamma:F\cong G$ is a natural isomorphism.
  Then for any $a:A$, we have an isomorphism $\gamma_a:Fa \cong Ga$, hence an identity $\isotoid(\gamma_a):\id{Fa}{Ga}$.
  By function extensionality, we have an identity $\bar{\gamma}:\id[(A_0\to B_0)]{F_0}{G_0}$.

  Now since the last two axioms of a functor are mere propositions, to show that $\id{F}{G}$ it will suffice to show that for any $a,b:A$, the functions
  \begin{align*}
    F_{a,b}&:\hom_A(a,b) \to \hom_B(Fa,Fb)\mathrlap{\qquad\text{and}}\\
    G_{a,b}&:\hom_A(a,b) \to \hom_B(Ga,Gb)
  \end{align*}
  become equal when transported along $\bar\gamma$.
  By computation for function extensionality, when applied to $a$, $\bar\gamma$ becomes equal to $\isotoid(\gamma_a)$.
  But by \autoref{ct:idtoiso-trans}, transporting $Ff:\hom_B(Fa,Fb)$ along $\isotoid(\gamma_a)$ and $\isotoid(\gamma_b)$ is equal to the composite $\gamma_b\circ Ff\circ \inv{(\gamma_a)}$, which by naturality of $\gamma$ is equal to $Gf$.

  This completes the definition of a function $(F\cong G) \to (\id F G)$.
  Now consider the composite
  \[ (\id F G) \to (F\cong G) \to (\id F G). \]
  Since hom-sets are sets, their identity types are mere propositions, so to show that two identities $p,q:\id F G$ are equal, it suffices to show that $\id[\id{F_0}{G_0}]{p}{q}$.
  But in the definition of $\bar\gamma$, if $\gamma$ were of the form $\idtoiso(p)$, then $\gamma_a$ would be equal to $\idtoiso(p_a)$ (this can easily be proved by induction on $p$).
  Thus, $\isotoid(\gamma_a)$ would be equal to $p_a$, and so by function extensionality we would have $\id{\bar\gamma}{p}$, which is what we need.

  Finally, consider the composite
  \[(F\cong G)\to  (\id F G) \to (F\cong G). \]
  Since identity of natural transformations can be tested componentwise, it suffices to show that for each $a$ we have $\id{\idtoiso(\bar\gamma)_a}{\gamma_a}$.
  But as observed above, we have $\id{\idtoiso(\bar\gamma)_a}{\idtoiso((\bar\gamma)_a)}$, while $\id{(\bar\gamma)_a}{\isotoid(\gamma_a)}$ by computation for function extensionality.
  Since $\isotoid$ and $\idtoiso$ are inverses, we have $\id{\idtoiso(\bar\gamma)_a}{\gamma_a}$ as desired.
\end{proof}

In particular, naturally isomorphic functors between categories (as opposed to precategories) are equal.

\mentalpause

We now define all the usual ways to compose functors and natural transformations.

\begin{defn}
  For functors $F:A\to B$ and $G:B\to C$, their composite $G\circ F:A\to C$ is given by
  \begin{itemize}
  \item The composite $(G_0\circ F_0) : A_0 \to C_0$
  \item For each $a,b:A$, the composite
    \[(G_{Fa,Fb}\circ F_{a,b}):\hom_A(a,b) \to \hom_C(GFa,GFb).\]
  \end{itemize}
  It is easy to check the axioms.
\end{defn}

\begin{defn}
  For functors $F:A\to B$ and $G,H:B\to C$ and a natural transformation $\gamma:G\to H$, the composite $(\gamma F):GF\to HF$ is given by
  \begin{itemize}
  \item For each $a:A$, the component $\gamma_{Fa}$.
  \end{itemize}
  Naturality is easy to check.
  Similarly, for $\gamma$ as above and $K:C\to D$, the composite $(K\gamma):KG\to KH$ is given by
  \begin{itemize}
  \item For each $b:B$, the component $K(\gamma_b)$.
  \end{itemize}
\end{defn}

\begin{lem}\label{ct:interchange}
  \index{interchange law}%
  For functors $F,G:A\to B$ and $H,K:B\to C$ and natural transformations $\gamma:F\to G$ and $\delta:H\to K$, we have
  \[\id{(\delta G)(H\gamma)}{(K\gamma)(\delta F)}.\]
\end{lem}
\begin{proof}
  It suffices to check componentwise: at $a:A$ we have
  \begin{align*}
    ((\delta G)(H\gamma))_a
    &\jdeq (\delta G)_{a}(H\gamma)_a\\
    &\jdeq \delta_{Ga}\circ H(\gamma_a)\\
    &= K(\gamma_a) \circ \delta_{Fa} \tag{by naturality of $\delta$}\\
    &\jdeq (K \gamma)_a\circ (\delta F)_a\\
    &\jdeq ((K \gamma)(\delta F))_a.\qedhere
  \end{align*}
\end{proof}

\index{horizontal composition!of natural transformations}%
\index{classical!category theory}%
Classically, one defines the ``horizontal composite'' of $\gamma:F\to G$ and $\delta:H\to K$ to be the common value of ${(\delta G)(H\gamma)}$ and ${(K\gamma)(\delta F)}$.
We will refrain from doing this, because while equal, these two transformations are not \emph{definitionally} equal.
This also has the consequence that we can use the symbol $\circ$ (or juxtaposition) for all kinds of composition unambiguously: there is only one way to compose two natural transformations (as opposed to composing a natural transformation with a functor on either side).

\begin{lem}\label{ct:functor-assoc}
  \index{associativity!of functor composition}
  Composition of functors is associative: $\id{H(GF)}{(HG)F}$.
\end{lem}
\begin{proof}
  Since composition of functions is associative, this follows immediately for the actions on objects and on homs.
  And since hom-sets are sets, the rest of the data is automatic.
\end{proof}

The equality in \autoref{ct:functor-assoc} is likewise not definitional.
(Composition of functions is definitionally associative, but the axioms that go into a functor must also be composed, and this breaks definitional associativity.)  For this reason, we need also to know about \emph{coherence}\index{coherence} for associativity.

\begin{lem}\label{ct:pentagon}
  \index{associativity!of functor composition!coherence of}%
  \autoref{ct:functor-assoc} is coherent, i.e.\ the following pentagon\index{pentagon, Mac Lane} of equalities commutes:
  \[ \xymatrix{ & K(H(GF)) \ar@{=}[dl] \ar@{=}[dr]\\
    (KH)(GF) \ar@{=}[d] && K((HG)F) \ar@{=}[d]\\
    ((KH)G)F && (K(HG))F \ar@{=}[ll] }
  \]
\end{lem}
\begin{proof}
  As in \autoref{ct:functor-assoc}, this is evident for the actions on objects, and the rest is automatic.
\end{proof}

We will henceforth abuse notation by writing $H\circ G\circ F$ or $HGF$ for either $H(GF)$ or $(HG)F$, transporting along \autoref{ct:functor-assoc} whenever necessary.
We have a similar coherence result for units.

\begin{lem}\label{ct:units}
  For a functor $F:A\to B$, we have equalities $\id{(1_B\circ F)}{F}$ and $\id{(F\circ 1_A)}{F}$, such that given also $G:B\to C$, the following triangle of equalities commutes.
  \[ \xymatrix{
    G\circ (1_B \circ F) \ar@{=}[rr] \ar@{=}[dr] &&
    (G\circ 1_B)\circ F \ar@{=}[dl] \\
    & G \circ F.}
  \]
\end{lem}

See \autoref{ct:pre2cat,ct:2cat} for further development of these ideas.

\section{Adjunctions}
\label{sec:adjunctions}

The definition of adjoint functors is straightforward; the main interesting aspect arises from proof-relevance.

\begin{defn}
  A functor $F:A\to B$ is a \define{left adjoint}
  \indexdef{left!adjoint}%
  \indexdef{adjoint!functor}%
  \indexdef{right!adjoint}%
  \indexdef{adjoint!functor}%
  \index{functor!adjoint}%
  if there exists
  \begin{itemize}
  \item A functor $G:B\to A$.
  \item A natural transformation $\eta:1_A \to GF$ (the \define{unit}\indexdef{unit!of an adjunction}).
  \item A natural transformation $\epsilon:FG\to 1_B$ (the \define{counit}\indexdef{counit of an adjunction}).
  \item $\id{(\epsilon F)(F\eta)}{1_F}$.
  \item $\id{(G\epsilon)(\eta G)}{1_G}$.
  \end{itemize}
\end{defn}

The last two equations are called the \define{triangle identities}\indexdef{triangle!identity} or \define{zigzag identities}\indexdef{zigzag identity}.
\indexdef{identity!triangle}\indexdef{identity!zigzag}
We leave it to the reader to define right adjoints analogously.

\begin{lem}\label{ct:adjprop}
  If $A$ is a category (but $B$ may be only a precategory), then the type ``$F$ is a left adjoint'' is a mere proposition.
\end{lem}
\begin{proof}
  Suppose we are given $(G,\eta,\epsilon)$ with the triangle identities and also $(G',\eta',\epsilon')$.
  Define $\gamma:G\to G'$ to be $(G'\epsilon)(\eta' G)$, and $\delta:G'\to G$ to be $(G\epsilon')(\eta G')$.
  Then
  \begin{align*}
    \delta\gamma &=
    (G\epsilon')(\eta G')(G'\epsilon)(\eta'G)\\
    &= (G\epsilon')(G F G'\epsilon)(\eta G' F G)(\eta'G)\\
    &= (G\epsilon)(G\epsilon'FG)(G F \eta' G)(\eta G)\\
    &= (G\epsilon)(\eta G)\\
    &= 1_G
  \end{align*}
  using \autoref{ct:interchange} and the triangle identities.
  Similarly, we show $\id{\gamma\delta}{1_{G'}}$, so $\gamma$ is a natural isomorphism $G\cong G'$.
  By \autoref{ct:functor-cat}, we have an identity $\id G {G'}$.

  Now we need to know that when $\eta$ and $\epsilon$ are transported along this identity, they become equal to $\eta'$ and $\epsilon'$.
  By \autoref{ct:idtoiso-trans}, this transport is given by composing with $\gamma$ or $\delta$ as appropriate.
  For $\eta$, this yields
  \begin{equation*}
    (G'\epsilon F)(\eta'GF)\eta
    = (G'\epsilon F)(G'F\eta)\eta'
    = \eta'
  \end{equation*}
  using \autoref{ct:interchange} and the triangle identity.
  The case of $\epsilon$ is similar.
  Finally, the triangle identities transport correctly automatically, since hom-sets are sets.
\end{proof}

In \autoref{sec:yoneda} we will give another proof of \autoref{ct:adjprop}.

\section{Equivalences}
\label{sec:equivalences}

It is usual in category theory to define an \emph{equivalence of categories} to be a functor $F:A\to B$ such that there exists a functor $G:B\to A$ and natural isomorphisms $F G \cong 1_B$ and $G F \cong 1_A$.
Unlike the property of being an adjunction, however, this would not be a mere proposition without truncating it, for the same reasons that the type of quasi-inverses is ill-behaved (see \autoref{sec:quasi-inverses}).
And as in \autoref{sec:hae}, we can avoid this by using the usual notion of \emph{adjoint} equivalence.
\indexdef{adjoint!equivalence!of (pre)categories}

\begin{defn}\label{ct:equiv}
  A functor $F:A\to B$ is an \define{equivalence of (pre)categories}
  \indexdef{equivalence!of (pre)categories}%
  \indexdef{category!equivalence of}%
  \indexdef{precategory!equivalence of}%
  \index{functor!equivalence}%
  if it is a left adjoint for which $\eta$ and $\epsilon$ are isomorphisms.
  We write $\cteqv A B$ for the type of equivalences of categories from $A$ to $B$.
\end{defn}

By \autoref{ct:adjprop,ct:isoprop}, if $A$ is a category, then the type ``$F$ is an equivalence of precategories'' is a mere proposition.

\begin{lem}\label{ct:adjointification}
  If for $F:A\to B$ there exists $G:B\to A$ and isomorphisms $GF\cong 1_A$ and $FG\cong 1_B$, then $F$ is an equivalence of precategories.
\end{lem}
\begin{proof}
  Just like the proof of \autoref{thm:equiv-iso-adj} for equivalences of types.
\end{proof}

\begin{defn}
  We say a functor $F:A\to B$ is \define{faithful}
  \indexdef{functor!faithful}%
  \index{faithful functor}
  if for all $a,b:A$, the function
  \[F_{a,b}:\hom_A(a,b) \to \hom_B(Fa,Fb)\]
  is injective, and \define{full}
  \indexdef{functor!full}%
  \indexdef{full functor}%
  if for all $a,b:A$ this function is surjective.
  If it is both (hence each $F_{a,b}$ is an equivalence) we say $F$ is \define{fully faithful}.
  \indexdef{functor!fully faithful}%
  \indexdef{fully faithful functor}%
\end{defn}

\begin{defn}
  We say a functor $F:A\to B$ is \define{split essentially surjective}
  \indexdef{functor!split essentially surjective}%
  \indexdef{split!essentially surjective functor}%
  if for all $b:B$ there exists an $a:A$ such that $Fa\cong b$.
\end{defn}

\begin{lem}\label{ct:ffeso}
  For any precategories $A$ and $B$ and functor $F:A\to B$, the following types are equivalent.
  \begin{enumerate}
  \item $F$ is an equivalence of precategories.\label{item:ct:ffeso1}
  \item $F$ is fully faithful and split essentially surjective.\label{item:ct:ffeso2}
  \end{enumerate}
\end{lem}
\begin{proof}
  Suppose $F$ is an equivalence of precategories, with $G,\eta,\epsilon$ specified.
  Then we have the function
  \begin{align*}
      \hom_B(Fa,Fb) &\to \hom_A(a,b),\\
      g &\mapsto \inv{\eta_b}\circ G(g)\circ \eta_a.
  \end{align*}
  For $f:\hom_A(a,b)$, we have
  \[ \inv{\eta_{b}}\circ G(F(f))\circ \eta_{a}  =
  \inv{\eta_{b}} \circ \eta_{b} \circ f=
  f
  \]
  while for $g:\hom_B(Fa,Fb)$ we have
  \begin{align*}
    F(\inv{\eta_b} \circ G(g)\circ\eta_a)
    &= F(\inv{\eta_b})\circ F(G(g))\circ F(\eta_a)\\
    &= \epsilon_{Fb}\circ F(G(g))\circ F(\eta_a)\\
    &= g\circ\epsilon_{Fa}\circ F(\eta_a)\\
    &= g
  \end{align*}
  using naturality of $\epsilon$, and the triangle identities twice.
  Thus, $F_{a,b}$ is an equivalence, so $F$ is fully faithful.
  Finally, for any $b:B$, we have $Gb:A$ and $\epsilon_b:FGb\cong b$.

  On the other hand, suppose $F$ is fully faithful and split essentially surjective.
  Define $G_0:B_0\to A_0$ by sending $b:B$ to the $a:A$ given by the specified essential splitting, and write $\epsilon_b$ for the likewise specified isomorphism $FGb\cong b$.

  Now for any $g:\hom_B(b,b')$, define $G(g):\hom_A(Gb,Gb')$ to be the unique morphism such that $\id{F(G(g))}{\inv{(\epsilon_{b'})}\circ g \circ \epsilon_b }$ (which exists since $F$ is fully faithful).
  Finally, for $a:A$ define $\eta_a:\hom_A(a,GFa)$ to be the unique morphism such that $\id{F\eta_a}{\inv{\epsilon_{Fa}}}$.
  It is easy to verify that $G$ is a functor and that $(G,\eta,\epsilon)$ exhibit $F$ as an equivalence of precategories.

  Now consider the composite~\ref{item:ct:ffeso1}$\to$\ref{item:ct:ffeso2}$\to$\ref{item:ct:ffeso1}.
  We clearly recover the same function $G_0:B_0 \to A_0$.
  For the action of $G$ on hom-sets, we must show that for $g:\hom_B(b,b')$, $G(g)$ is the (necessarily unique) morphism such that $F(G(g)) = \inv{(\epsilon_{b'})}\circ g \circ \epsilon_b$.
  But this equation holds by the assumed naturality of $\epsilon$.
  We also clearly recover $\epsilon$, while $\eta$ is uniquely characterized by $\id{F\eta_a}{\inv{\epsilon_{Fa}}}$ (which is one of the triangle identities assumed to hold in the structure of an equivalence of precategories).
  Thus, this composite is equal to the identity.

  Finally, consider the other composite~\ref{item:ct:ffeso2}$\to$\ref{item:ct:ffeso1}$\to$\ref{item:ct:ffeso2}.
  Since being fully faithful is a mere proposition, it suffices to observe that we recover, for each $b:B$, the same $a:A$ and isomorphism $F a \cong b$.
  But this is clear, since we used this function and isomorphism to define $G_0$ and $\epsilon$ in~\ref{item:ct:ffeso1}, which in turn are precisely what we used to recover~\ref{item:ct:ffeso2} again.
  Thus, the composites in both directions are equal to identities, hence we have an equivalence \eqv{\text{\ref{item:ct:ffeso1}}}{\text{\ref{item:ct:ffeso2}}}.
\end{proof}

However, if $B$ is not a category, then neither type in \autoref{ct:ffeso} may necessarily be a mere proposition.
This suggests considering as well the following notions.

\begin{defn}
  A functor $F:A\to B$ is \define{essentially surjective}
  \indexdef{functor!essentially surjective}%
  \indexdef{essentially surjective functor}%
  if for all $b:B$, there \emph{merely} exists an $a:A$ such that $Fa\cong b$.
  We say $F$ is a \define{weak equivalence}
  \indexsee{equivalence!of (pre)categories!weak}{weak equivalence}%
  \indexdef{weak equivalence!of precategories}%
  \indexsee{functor!weak equivalence}{weak equivalence}%
  if it is fully faithful and essentially surjective.
\end{defn}

Being a weak equivalence is \emph{always} a mere proposition.
For categories, however, there is no difference between equivalences and weak ones.

\index{acceptance}
\begin{lem}\label{ct:catweq}
  If $F:A\to B$ is fully faithful and $A$ is a category, then for any $b:B$ the type $\sm{a:A} (Fa\cong b)$ is a mere proposition.
  Hence a functor between categories is an equivalence if and only if it is a weak equivalence.
\end{lem}
\begin{proof}
  Suppose given $(a,f)$ and $(a',f')$ in $\sm{a:A} (Fa\cong b)$.
  Then $\inv{f'}\circ f$ is an isomorphism $Fa \cong Fa'$.
  Since $F$ is fully faithful, we have $g:a\cong a'$ with $Fg = \inv{f'}\circ f$.
  And since $A$ is a category, we have $p:a=a'$ with $\idtoiso(p)=g$.
  Now $Fg = \inv{f'}\circ f$ implies $\trans{(\map{(F_0)}{p})}{f} = f'$, hence (by the characterization of equalities in dependent pair types) $(a,f)=(a',f')$.

  Thus, for fully faithful functors whose domain is a category, essential surjectivity is equivalent to split essential surjectivity, and so being a weak equivalence is equivalent to being an equivalence.
\end{proof}

This is an important advantage of our category theory over set-based approaches.
With a purely set-based definition of category, the statement ``every fully faithful and essentially surjective functor is an equivalence of categories'' is equivalent to the axiom of choice \choice{}.
Here we have it for free, as a category-theoretic version of the principle of unique choice (\autoref{sec:unique-choice}).
(In fact, this property characterizes categories among precategories; see \autoref{sec:rezk}.)

On the other hand, the following characterization of equivalences of categories is perhaps even more useful.

\begin{defn}\label{ct:isocat}
  A functor $F:A\to B$ is an \define{isomorphism of (pre)cat\-ego\-ries}
  \indexdef{isomorphism!of (pre)categories}%
  \indexdef{category!isomorphism of}%
  \indexdef{precategory!isomorphism of}%
  if $F$ is fully faithful and $F_0:A_0\to B_0$ is an equivalence of types.
\end{defn}

This definition is an exception to our general rule (see \autoref{sec:basics-equivalences}) of only using the word ``isomorphism'' for sets and set-like objects.
However, it does carry an appropriate connotation here, because for general precategories, isomorphism is stronger than equivalence.

Note that being an isomorphism of precategories is always a mere property.
Let $A\cong B$ denote the type of isomorphisms of (pre)categories from $A$ to $B$.

\begin{lem}\label{ct:isoprecat}
  For precategories $A$ and $B$ and $F:A\to B$, the following are equivalent.
  \begin{enumerate}
  \item $F$ is an isomorphism of precategories.\label{item:ct:ipc1}
  \item There exist $G:B\to A$ and $\eta:1_A = GF$ and $\epsilon:FG=1_B$ such that\label{item:ct:ipc2}
    \begin{equation}
      \apfunc{(\lam{H} F H)}({\eta}) = \apfunc{(\lam{K} K F)}({\opp\epsilon}).\label{eq:ct:isoprecattri}
    \end{equation}
  \item There merely exist $G:B\to A$ and $\eta:1_A = GF$ and $\epsilon:FG=1_B$.\label{item:ct:ipc3}
  \end{enumerate}
\end{lem}

Note that if $B_0$ is not a 1-type, then~\eqref{eq:ct:isoprecattri} may not be a mere proposition.

\begin{proof}
  First note that since hom-sets are sets, equalities between equalities of functors are uniquely determined by their object-parts.
  Thus, by function extensionality,~\eqref{eq:ct:isoprecattri} is equivalent to
  \begin{equation}
    \map{(F_0)}{\eta_0}_a = \opp{(\epsilon_0)}_{F_0 a}.\label{eq:ct:ipctri}
  \end{equation}
  for all $a:A_0$.
  Note that this is precisely the triangle identity for $G_0$, $\eta_0$, and $\epsilon_0$ to be a proof that $F_0$ is a half adjoint equivalence of types.

  Now suppose~\ref{item:ct:ipc1}.
  Let $G_0:B_0 \to A_0$ be the inverse of $F_0$, with $\eta_0: \idfunc[A_0] = G_0 F_0$ and $\epsilon_0:F_0G_0 = \idfunc[B_0]$ satisfying the triangle identity, which is precisely~\eqref{eq:ct:ipctri}.
  Now define $G_{b,b'}:\hom_B(b,b') \to \hom_A(G_0b,G_0b')$ by
  \[ G_{b,b'}(g) \defeq
  \inv{(F_{G_0b,G_0b'})}\Big(\idtoiso(\opp{(\epsilon_0)}_{b'}) \circ g \circ \idtoiso((\epsilon_0)_b)\Big)
  \]
  (using the assumption that $F$ is fully faithful).
  Since \idtoiso takes inverses to inverses and concatenation to composition, and $F$ is a functor, it follows that $G$ is a functor.

  By definition, we have $(GF)_0 \jdeq G_0 F_0$, which is equal to $\idfunc[A_0]$ by $\eta_0$.
  To obtain $1_A = GF$, we need to show that when transported along $\eta_0$, the identity function of $\hom_A(a,a')$ becomes equal to the composite $G_{Fa,Fa'} \circ F_{a,a'}$.
  In other words, for any $f:\hom_A(a,a')$ we must have
  \begin{multline*}
    \idtoiso((\eta_0)_{a'}) \circ f \circ \idtoiso(\opp{(\eta_0)}_a)\\
    = \inv{(F_{GFa,GFa'})}\Big(\idtoiso(\opp{(\epsilon_0)}_{Fa'})
    \circ F_{a,a'}(f) \circ \idtoiso((\epsilon_0)_{Fa})\Big).
  \end{multline*}
  But this is equivalent to
  \begin{multline*}
    (F_{GFa,GFa'})\Big(\idtoiso((\eta_0)_{a'}) \circ f \circ \idtoiso(\opp{(\eta_0)}_a)\Big)\\
    = \idtoiso(\opp{(\epsilon_0)}_{Fa'})
    \circ F_{a,a'}(f) \circ \idtoiso((\epsilon_0)_{Fa}).
  \end{multline*}
  which follows from functoriality of $F$, the fact that $F$ preserves \idtoiso, and~\eqref{eq:ct:ipctri}.
  Thus we have $\eta:1_A = GF$.

  On the other side, we have $(FG)_0\jdeq F_0 G_0$, which is equal to $\idfunc[B_0]$ by $\epsilon_0$.
  To obtain $FG=1_B$, we need to show that when transported along $\epsilon_0$, the identity function of $\hom_B(b,b')$ becomes equal to the composite $F_{Gb,Gb'} \circ G_{b,b'}$.
  That is, for any $g:\hom_B(b,b')$ we must have
  \begin{multline*}
    F_{Gb,Gb'}\Big(\inv{(F_{Gb,Gb'})}\Big(\idtoiso(\opp{(\epsilon_0)}_{b'}) \circ g \circ \idtoiso((\epsilon_0)_b)\Big)\Big)\\
    = \idtoiso((\opp{\epsilon_0})_{b'}) \circ g \circ \idtoiso((\epsilon_0)_b).
  \end{multline*}
  But this is just the fact that $\inv{(F_{Gb,Gb'})}$ is the inverse of $F_{Gb,Gb'}$.
  And we have remarked that~\eqref{eq:ct:isoprecattri} is equivalent to~\eqref{eq:ct:ipctri}, so~\ref{item:ct:ipc2} holds.

  Conversely, suppose given~\ref{item:ct:ipc2}; then the object-parts of $G$, $\eta$, and $\epsilon$ together with~\eqref{eq:ct:ipctri} show that $F_0$ is an equivalence of types.
  And for $a,a':A_0$, we define $\overline{G}_{a,a'}: \hom_B(Fa,Fa') \to \hom_A(a,a')$ by
  \begin{equation}
    \overline{G}_{a,a'}(g) \defeq \idtoiso(\opp{\eta})_{a'} \circ G(g) \circ \idtoiso(\eta)_a.\label{eq:ct:gbar}
  \end{equation}
  By naturality of $\idtoiso(\eta)$, for any $f:\hom_A(a,a')$ we have
  \begin{align*}
    \overline{G}_{a,a'}(F_{a,a'}(f))
    &= \idtoiso(\opp{\eta})_{a'} \circ G(F(f)) \circ \idtoiso(\eta)_a\\
    &= \idtoiso(\opp{\eta})_{a'} \circ \idtoiso(\eta)_{a'} \circ f \\
    &= f.
  \end{align*}
  On the other hand, for $g:\hom_B(Fa,Fa')$ we have
  \begin{align*}
    F_{a,a'}(\overline{G}_{a,a'}(g))
    &= F(\idtoiso(\opp{\eta})_{a'}) \circ F(G(g)) \circ F(\idtoiso(\eta)_a)\\
    &= \idtoiso(\epsilon)_{Fa'}
    \circ F(G(g))
    \circ \idtoiso(\opp{\epsilon})_{Fa}\\
    &= \idtoiso(\epsilon)_{Fa'}
    \circ \idtoiso(\opp{\epsilon})_{Fa'}
    \circ g\\
    &= g.
  \end{align*}
  (There are lemmas needed here regarding the compatibility of \idtoiso and whiskering, which we leave it to the reader to state and prove.)
  Thus, $F_{a,a'}$ is an equivalence, so $F$ is fully faithful; i.e.~\ref{item:ct:ipc1} holds.

  Now the composite~\ref{item:ct:ipc1}$\to$\ref{item:ct:ipc2}$\to$\ref{item:ct:ipc1} is equal to the identity since~\ref{item:ct:ipc1} is a mere proposition.
  On the other side, tracing through the above constructions we see that the composite~\ref{item:ct:ipc2}$\to$\ref{item:ct:ipc1}$\to$\ref{item:ct:ipc2} essentially preserves the object-parts $G_0$, $\eta_0$, $\epsilon_0$, and the object-part of~\eqref{eq:ct:isoprecattri}.
  And in the latter three cases, the object-part is all there is, since hom-sets are sets.

  Thus, it suffices to show that we recover the action of $G$ on hom-sets.
  In other words, we must show that if $g:\hom_B(b,b')$, then
  \[ G_{b,b'}(g) =
  \overline{G}_{G_0b,G_0b'}\Big(\idtoiso(\opp{(\epsilon_0)}_{b'}) \circ g \circ \idtoiso((\epsilon_0)_b)\Big)
  \]
  where $\overline{G}$ is defined by~\eqref{eq:ct:gbar}.
  However, this follows from functoriality of $G$ and the \emph{other} triangle identity, which we have seen in \autoref{cha:equivalences} is equivalent to~\eqref{eq:ct:ipctri}.

  Now since~\ref{item:ct:ipc1} is a mere proposition, so is~\ref{item:ct:ipc2}, so it suffices to show they are logically equivalent to~\ref{item:ct:ipc3}.
  Of course,~\ref{item:ct:ipc2}$\to$\ref{item:ct:ipc3}, so let us assume~\ref{item:ct:ipc3}.
  Since~\ref{item:ct:ipc1} is a mere proposition, we may assume given $G$, $\eta$, and $\epsilon$.
  Then $G_0$ along with $\eta$ and $\epsilon$ imply that $F_0$ is an equivalence.
  Moreover, we also have natural isomorphisms $\idtoiso(\eta):1_A\cong GF$ and $\idtoiso(\epsilon):FG\cong 1_B$, so by \autoref{ct:adjointification}, $F$ is an equivalence of precategories, and in particular fully faithful.
\end{proof}

From \autoref{ct:isoprecat}\ref{item:ct:ipc2} and $\idtoiso$ in functor categories, we conclude immediately that any isomorphism of precategories is an equivalence.
For precategories, the converse can fail.

\begin{eg}\label{ct:chaotic}
  Let $X$ be a type and $x_0:X$ an element, and let $X_{\mathrm{ch}}$ denote the \emph{chaotic}\indexdef{chaotic precategory} or \emph{indiscrete}\indexdef{indiscrete precategory} precategory on $X$.
  By definition, we have $(X_{\mathrm{ch}})_0\defeq X$, and $\hom_{X_{\mathrm{ch}}}(x,x') \defeq \unit$ for all $x,x'$.
  Then the unique functor $X_{\mathrm{ch}}\to \unit$ is an equivalence of precategories, but not an isomorphism unless $X$ is contractible.

  This example also shows that a precategory can be equivalent to a category without itself being a category.
  Of course, if a precategory is \emph{isomorphic} to a category, then it must itself be a category.
\end{eg}

However, for categories, the two notions coincide.

\begin{lem}\label{ct:eqv-levelwise}
  For categories $A$ and $B$, a functor $F:A\to B$ is an equivalence of categories if and only if it is an isomorphism of categories.
\end{lem}
\begin{proof}
  Since both are mere properties, it suffices to show they are logically equivalent.
  So first suppose $F$ is an equivalence of categories, with $(G,\eta,\epsilon)$ given.
  We have already seen that $F$ is fully faithful.
  By \autoref{ct:functor-cat}, the natural isomorphisms $\eta$ and $\epsilon$ yield identities $\id{1_A}{GF}$ and $\id{FG}{1_B}$, hence in particular identities $\id{\idfunc[A]}{G_0\circ F_0}$ and $\id{F_0\circ G_0}{\idfunc[B]}$.
Thus, $F_0$ is an equivalence of types.

  Conversely, suppose $F$ is fully faithful and $F_0$ is an equivalence of types, with inverse $G_0$, say.
  Then for each $b:B$ we have $G_0 b:A$ and an identity $\id{FGb}{b}$, hence an isomorphism $FGb\cong b$.
  Thus, by \autoref{ct:ffeso}, $F$ is an equivalence of categories.
\end{proof}

Of course, there is yet a third notion of sameness for (pre)categories: equality.
However, the univalence axiom implies that it coincides with isomorphism.

\begin{lem}\label{ct:cat-eq-iso}
  If $A$ and $B$ are precategories, then the function
  \[(\id A B) \to (A\cong B)\]
  (defined by induction from the identity functor) is an equivalence of types.
\end{lem}
\begin{proof}
  As usual for dependent sum types, to give an element of $\id A B$ is equivalent to giving
  \begin{itemize}
  \item an identity $P_0:\id{A_0}{B_0}$,
  \item for each $a,b:A_0$, an identity
    \[P_{a,b}:\id{\hom_A(a,b)}{\hom_B(\trans {P_0} a,\trans {P_0} b)},\]
  \item identities $\id{\trans {(P_{a,a})} {1_a}}{1_{\trans {P_0} a}}$ and
    \narrowequation{\id{\trans {(P_{a,c})} {gf}}{\trans {(P_{b,c})} g \circ \trans {(P_{a,b})} f}.}
  \end{itemize}
  (Again, we use the fact that the identity types of hom-sets are mere propositions.)
  However, by univalence, this is equivalent to giving
  \begin{itemize}
  \item an equivalence of types $F_0:\eqv{A_0}{B_0}$,
  \item for each $a,b:A_0$, an equivalence of types
    \[F_{a,b}:\eqv{\hom_A(a,b)}{\hom_B(F_0 (a),F_0 (b))},\]
  \item and identities $\id{F_{a,a}(1_a)}{1_{F_0 (a)}}$ and $\id{F_{a,c}(gf)}{F_{b,c} (g)\circ F_{a,b} (f)}$.
  \end{itemize}
  But this consists exactly of a functor $F:A\to B$ that is an isomorphism of categories.
  And by induction on identity, this equivalence $\eqv{(\id A B)}{(A\cong B)}$ is equal to the one obtained by induction.
\end{proof}

Thus, for categories, equality also coincides with equivalence.
We can interpret this as saying that categories, functors, and natural transformations form, not just a pre-2-category, but a 2-category.

\begin{thm}\label{ct:cat-2cat}
  If $A$ and $B$ are categories, then the function
  \[(\id A B) \to (\cteqv A B)\]
  (defined by induction from the identity functor) is an equivalence of types.
\end{thm}
\begin{proof}
  By \autoref{ct:cat-eq-iso,ct:eqv-levelwise}.
\end{proof}

As a consequence, the type of categories is a 2-type.
For since $\cteqv A B$ is a subtype of the type of functors from $A$ to $B$, which are the objects of a category, it is a 1-type; hence the identity types $\id A B$ are also 1-types.

\section{The Yoneda lemma}
\label{sec:yoneda}
\index{Yoneda!lemma|(}
 
Recall that we have a category \uset whose objects are sets and whose morphisms are functions.
We now show that every precategory has a \uset-valued hom-functor.
First we need to define opposites and products of (pre)categories.

\begin{defn}\label{ct:opposite-category}
  For a precategory $A$, its \define{opposite}
  \indexdef{opposite of a (pre)category}%
  \indexdef{precategory!opposite}%
  \indexdef{category!opposite}%
  $A\op$ is a precategory with the same type of objects, with $\hom_{A\op}(a,b) \defeq \hom_A(b,a)$, and with identities and composition inherited from $A$.
\end{defn}

\begin{defn}
  For precategories $A$ and $B$, their \define{product}
  \index{precategory!product of}%
  \index{category!product of}%
  \index{product!of (pre)categories}%
  $A\times B$ is a precategory with $(A\times B)_0 \defeq A_0 \times B_0$ and
  \[\hom_{A\times B}((a,b),(a',b')) \defeq \hom_A(a,a') \times \hom_B(b,b').\]
  Identities are defined by $1_{(a,b)}\defeq (1_a,1_b)$ and composition by
  \narrowequation{(g,g')(f,f') \defeq ((gf),(g'f')).}
\end{defn}

\begin{lem}\label{ct:functorexpadj}
  For precategories $A,B,C$, the following types are equivalent.
  \begin{enumerate}
  \item Functors $A\times B\to C$.
  \item Functors $A\to C^B$.
  \end{enumerate}
\end{lem}
\begin{proof}
  Given $F:A\times B\to C$, for any $a:A$ we obviously have a functor $F_a : B\to C$.
  This gives a function $A_0 \to (C^B)_0$.
  Next, for any $f:\hom_A(a,a')$, we have for any $b:B$ the morphism $F_{(a,b),(a',b)}(f,1_b):F_a(b) \to F_{a'}(b)$.
  These are the components of a natural transformation $F_a \to F_{a'}$.
  Functoriality in $a$ is easy to check, so we have a functor $\hat{F}:A\to C^B$.

  Conversely, suppose given $G:A\to C^B$.
  Then for any $a:A$ and $b:B$ we have the object $G(a)(b):C$, giving a function $A_0 \times B_0 \to C_0$.
  And for $f:\hom_A(a,a')$ and $g:\hom_B(b,b')$, we have the morphism
  \begin{equation*}
     G(a')_{b,b'}(g)\circ G_{a,a'}(f)_b = G_{a,a'}(f)_{b'} \circ  G(a)_{b,b'}(g)
  \end{equation*}
  in $\hom_C(G(a)(b), G(a')(b'))$.
  Functoriality is again easy to check, so we have a functor $\check{G}:A\times B \to C$.

  Finally, it is also clear that these operations are inverses.
\end{proof}

Now for any precategory $A$, we have a hom-functor
\indexdef{hom-functor}%
\[\hom_A : A\op \times A \to \uset.\]
It takes a pair $(a,b): (A\op)_0 \times A_0 \jdeq A_0 \times A_0$ to the set $\hom_A(a,b)$.
For a morphism $(f,f') : \hom_{A\op\times A}((a,b),(a',b'))$, by definition we have $f:\hom_A(a',a)$ and $f':\hom_A(b,b')$, so we can define 
\begin{align*}
  (\hom_A)_{(a,b),(a',b')}(f,f')
  &\defeq (g \mapsto (f'gf))\\
  &: \hom_A(a,b) \to \hom_A(a',b').
\end{align*}
Functoriality is easy to check.

By \autoref{ct:functorexpadj}, therefore, we have an induced functor $\y:A\to \uset^{A\op}$, which we call the \define{Yoneda embedding}.
\indexdef{Yoneda!embedding}%
\indexdef{embedding!Yoneda}%

\begin{thm}[The Yoneda lemma]\label{ct:yoneda}
  \indexdef{Yoneda!lemma}
  For any precategory $A$, any $a:A$, and any functor $F:\uset^{A\op}$, we have an isomorphism
  \begin{equation}\label{eq:yoneda}
    \hom_{\uset^{A\op}}(\y a, F) \cong Fa.
  \end{equation}
  Moreover, this is natural in both $a$ and $F$.
\end{thm}
\begin{proof}
  Given a natural transformation $\alpha:\y a \to F$, we can consider the component $\alpha_a : \y a(a) \to F a$.
  Since $\y a(a)\jdeq \hom_A(a,a)$, we have $1_a : \y a(a)$, so that $\alpha_a(1_a) : F a$.
  This gives a function $(\alpha \mapsto \alpha_a(1_a))$ from left to right in~\eqref{eq:yoneda}.

  In the other direction, given $x:F a$, we define $\alpha:\y a \to F$ by
  \[\alpha_{a'}(f) \defeq F_{a',a}(f)(x). \]
  Naturality is easy to check, so this gives a function from right to left in~\eqref{eq:yoneda}.

  To show that these are inverses, first suppose given $x:F a$.
  Then with $\alpha$ defined as above, we have $\alpha_a(1_a) = F_{a,a}(1_a)(x) = 1_{F a}(x) = x$.
  On the other hand, if we suppose given $\alpha:\y a \to F$ and define $x$ as above, then for any $f:\hom_A(a',a)$ we have
  \begin{align*}
    \alpha_{a'}(f)
    &= \alpha_{a'} (\y a_{a',a}(f))\\
    &= (\alpha_{a'}\circ \y a_{a',a}(f))(1_a)\\
    &= (F_{a',a}(f)\circ \alpha_a)(1_a)\\
    &= F_{a',a}(f)(\alpha_a(1_a))\\
    &= F_{a',a}(f)(x).
  \end{align*}
  Thus, both composites are equal to identities.
  We leave the proof of naturality to the reader.
\end{proof}

\begin{cor}\label{ct:yoneda-embedding}
  The Yoneda embedding $\y :A\to \uset^{A\op}$ is fully faithful.
\end{cor}
\begin{proof}
  By \autoref{ct:yoneda}, we have
  \[ \hom_{\uset^{A\op}}(\y a, \y b) \cong \y b(a) \jdeq \hom_A(a,b). \]
  It is easy to check that this isomorphism is in fact the action of \y on hom-sets.
\end{proof}

\begin{cor}\label{ct:yoneda-mono}
  If $A$ is a category, then $\y_0 : A_0 \to (\uset^{A\op})_0$ is an embedding.
  In particular, if $\y a = \y b$, then $a=b$.
\end{cor}
\begin{proof}
  By \autoref{ct:yoneda-embedding}, \y induces an isomorphism on sets of isomorphisms.
  But as $A$ and $\uset^{A\op}$ are categories and \y is a functor, this is equivalently an isomorphism on identity types, which is the definition of being an embedding.
\end{proof}

\begin{defn}\label{ct:representable}
  A functor $F:\uset^{A\op}$ is said to be \define{representable}
  \indexdef{functor!representable}%
  \indexdef{representable functor}%
  if there exists $a:A$ and an isomorphism $\y a \cong F$.
\end{defn}

\begin{thm}\label{ct:representable-prop}
  If $A$ is a category, then the type ``$F$ is representable'' is a mere proposition.
\end{thm}
\begin{proof}
  By definition ``$F$ is representable'' is just the fiber of $\y_0$ over $F$.
  Since $\y_0$ is an embedding by \autoref{ct:yoneda-mono}, this fiber is a mere proposition.
\end{proof}

In particular, in a category, any two representations of the same functor are equal.
We can use this to give a different proof of \autoref{ct:adjprop}.
First we give a characterization of adjunctions in terms of representability.

\begin{lem}\label{ct:adj-repr}
  For any precategories $A$ and $B$ and a functor $F:A\to B$, the following types are equivalent.
  \begin{enumerate}
  \item $F$ is a left adjoint\index{adjoint!functor}.\label{item:ct:ar1}
  \item For each $b:B$, the functor $(a \mapsto \hom_B(Fa,b))$ from $A\op$ to \uset is representable\index{representable functor}.\label{item:ct:ar2}
  \end{enumerate}
\end{lem}
\begin{proof}
  An element of the type~\ref{item:ct:ar2} consists of a function $G_0:B_0 \to A_0$ together with, for every $a:A$ and $b:B$ an isomorphism
  \[ \gamma_{a,b}:\hom_B(Fa,b) \cong \hom_A(a,G_0 b) \]
  such that $\gamma_{a,b}(g \circ Ff) = \gamma_{a',b}(g)\circ f$ for $f:\hom_{A}(a,a')$.
  
  Given this, for $a:A$ we define $\eta_a \defeq \gamma_{a,Fa}(1_{Fa})$, and for $b:B$ we define $\epsilon_b \defeq \inv{(\gamma_{Gb,b})}(1_{Gb})$.
  Now for $g:\hom_B(b,b')$ we define
  \[ G_{b,b'}(g) \defeq \gamma_{G b, b'}(g \circ \epsilon_b) \]
  The verifications that $G$ is a functor and $\eta$ and $\epsilon$ are natural transformations satisfying the triangle identities are exactly as in the classical case, and as they are all mere propositions we will not care about their values.
  Thus, we have a function~\ref{item:ct:ar2}$\to$\ref{item:ct:ar1}.

  In the other direction, if $F$ is a left adjoint, we of course have $G_0$ specified, and we can take $\gamma_{a,b}$ to be the composite
  \[ \hom_B(Fa,b)
  \xrightarrow{G_{Fa,b}} \hom_A(GFa,Gb)
  \xrightarrow{(\blank\circ \eta_a)} \hom_A(a,Gb).
  \]
  This is clearly natural since $\eta$ is, and it has an inverse given by
  \[ \hom_A(a,Gb)
  \xrightarrow{F_{a,Gb}} \hom_B(Fa,FGb)
  \xrightarrow{(\epsilon_b \circ \blank )} \hom_A(Fa,b)
  \]
  (by the triangle identities).
  Thus we also have~\ref{item:ct:ar1}$\to$~\ref{item:ct:ar2}.

  For the composite~\ref{item:ct:ar2}$\to$\ref{item:ct:ar1}$\to$~\ref{item:ct:ar2}, clearly the function $G_0$ is preserved, so it suffices to check that we get back $\gamma$.
  But the new $\gamma$ is defined to take $f:\hom_B(Fa,b)$ to
  \begin{align*}
    G(f) \circ \eta_a
    &\jdeq \gamma_{G Fa, b}(f \circ \epsilon_{Fa}) \circ \eta_a\\
    &= \gamma_{G Fa, b}(f \circ \epsilon_{Fa} \circ F\eta_a)\\
    &= \gamma_{G Fa, b}(f)
  \end{align*}
  so it agrees with the old one.

  Finally, for~\ref{item:ct:ar1}$\to$\ref{item:ct:ar2}$\to$~\ref{item:ct:ar1}, we certainly get back the functor $G$ on objects.
  The new $G_{b,b'}:\hom_B(b,b') \to \hom_A(Gb,Gb')$ is defined to take $g$ to
  \begin{align*}
    \gamma_{G b, b'}(g \circ \epsilon_b)
    &\jdeq G(g \circ \epsilon_b) \circ \eta_{Gb}\\
    &= G(g) \circ G\epsilon_b \circ \eta_{Gb}\\
    &= G(g)
  \end{align*}
  so it agrees with the old one.
  The new $\eta_a$ is defined to be $\gamma_{a,Fa}(1_{Fa}) \jdeq G(1_{Fa}) \circ \eta_a$, so it equals the old $\eta_a$.
  And finally, the new $\epsilon_b$ is defined to be $\inv{(\gamma_{Gb,b})}(1_{Gb}) \jdeq \epsilon_b \circ F(1_{Gb})$, which also equals the old $\epsilon_b$.
\end{proof}

\begin{cor}\label{ct:adjprop2}[\autoref{ct:adjprop}]
  If $A$ is a category and $F:A\to B$, then the type ``$F$ is a left adjoint'' is a mere proposition.
\end{cor}
\begin{proof}
  By \autoref{ct:representable-prop}, if $A$ is a category then the type in \autoref{ct:adj-repr}\ref{item:ct:ar2} is a mere proposition.
\end{proof}
\index{Yoneda!lemma|)}

\section{Strict categories}
\label{sec:strict-categories}

\index{bargaining|(}%

\begin{defn}
  A \define{strict category}
  \indexdef{category!strict}%
  \indexdef{strict!category}%
  is a precategory whose type of objects is a set.
\end{defn}

In accordance with the mathematical red herring principle,\index{red herring principle} a strict category is not necessarily a category.
In fact, a category is a strict category precisely when it is gaunt (\autoref{ct:gaunt}).
\index{gaunt category}%
\index{category!gaunt}%
Most of the time, category theory is about categories, not strict ones, but sometimes one wants to consider strict categories.
The main advantage of this is that strict categories have a stricter notion of ``sameness'' than equivalence, namely isomorphism (or equivalently, by \autoref{ct:cat-eq-iso}, equality).

Here is one origin of strict categories.

\begin{eg}
  Let $A$ be a precategory and $x:A$ an object.
  Then there is a precategory $\mathsf{mono}(A,x)$ as follows:
  \index{monomorphism}
  \indexsee{mono}{monomorphism}
  \indexsee{monic}{monomorphism}
  \begin{itemize}
  \item Its objects consist of an object $y:A$ and a monomorphism $m:\hom_A(y,x)$.
    (As usual, $m:\hom_A(y,x)$ is a \define{monomorphism} (or is \define{monic}) if $(m\circ f = m\circ g) \Rightarrow (f=g)$.)
  \item Its morphisms from $(y,m)$ to $(z,n)$ are arbitrary morphisms from $y$ to $z$ in $A$ (not necessarily respecting $m$ and $n$).
  \end{itemize}
  An equality $(y,m)=(z,n)$ of objects in $\mathsf{mono}(A,x)$ consists of an equality $p:y=z$ and an equality $\trans{p}{m}=n$, which by \autoref{ct:idtoiso-trans} is equivalently an equality $m=n\circ \idtoiso(p)$.
  Since hom-sets are sets, the type of such equalities is a mere proposition.
  But since $m$ and $n$ are monomorphisms, the type of morphisms $f$ such that $m = n\circ f$ is also a mere proposition.
  Thus, if $A$ is a category, then $(y,m)=(z,n)$ is a mere proposition, and hence $\mathsf{mono}(A,x)$ is a strict category.
\end{eg}

This example can be dualized, and generalized in various ways.
Here is an interesting application of strict categories.

\begin{eg}\label{ct:galois}
  Let $E/F$ be a finite Galois extension
  \index{Galois!extension}%
  of fields, and $G$ its Galois group.
  \index{Galois!group}%
  Then there is a strict category whose objects are intermediate fields $F\subseteq K\subseteq E$, and whose morphisms are field homomorphisms\index{homomorphism!field} which fix $F$ pointwise (but need not commute with the inclusions into $E$).
  There is another strict category whose objects are subgroups $H\subseteq G$, and whose morphisms are morphisms of $G$-sets $G/H \to G/K$.
  The fundamental theorem of Galois theory
  \index{fundamental!theorem of Galois theory}%
  says that these two precategories are isomorphic (not merely equivalent).
\end{eg}

\index{bargaining|)}%

\section{\texorpdfstring{$\dagger$}{†}-categories}
\label{sec:dagger-categories}

It is also worth mentioning a useful kind of precategory whose type of objects is not a set, but which is not a category either.

\begin{defn}
  A \define{$\dagger$-precategory}
  \indexdef{.dagger-precategory@$\dagger$-precategory}%
  \indexdef{precategory!.dagger-@$\dagger$-}%
  is a precategory $A$ together with the following.
  \begin{enumerate}
  \item For each $x,y:A$, a function $\dgr{(-)}:\hom_A(x,y) \to \hom_A(y,x)$.
  \item For all $x:A$, we have $\dgr{(1_x)} = 1_x$.
  \item For all $f,g$ we have $\dgr{(g\circ f)} = \dgr f \circ \dgr g$.
  \item For all $f$ we have $\dgr{(\dgr f)} = f$.
  \end{enumerate}
\end{defn}

\begin{defn}\label{ct:unitary}
  A morphism $f:\hom_A(x,y)$ in a $\dagger$-precategory is \define{unitary}
  \indexdef{.dagger-precategory@$\dagger$-precategory!unitary morphism in}%
  \indexdef{unitary morphism}%
  \indexdef{morphism!unitary}%
  \indexdef{isomorphism!unitary}%
  if $\dgr f \circ f = 1_x$ and $f\circ \dgr f = 1_y$.
\end{defn}

Of course, every unitary morphism is an isomorphism, and being unitary is a mere proposition.
Thus for each $x,y:A$ we have a set of unitary isomorphisms from $x$ to $y$, which we denote $(x\unitaryiso y)$.

\begin{lem}\label{ct:idtounitary}
  If $p:(x=y)$, then $\idtoiso(p)$ is unitary.
\end{lem}
\begin{proof}
  By induction, we may assume $p$ is $\refl x$.
  But then $\dgr{(1_x)} \circ 1_x = 1_x\circ 1_x = 1_x$ and similarly.
\end{proof}

\begin{defn}
  A \define{$\dagger$-category}
  \indexdef{.dagger-category@$\dagger$-category}%
  is a $\dagger$-precategory such that for all $x,y:A$, the function
  \[ (x=y) \to (x \unitaryiso y) \]
  from \autoref{ct:idtounitary} is an equivalence.
\end{defn}

\begin{eg}
  The category \urel from \autoref{ct:rel} becomes a $\dagger$-pre\-cat\-e\-go\-ry if we define $(\dgr R)(y,x) \defeq R(x,y)$.
  The proof that \urel is a category actually shows that every isomorphism is unitary; hence \urel is also a $\dagger$-category.
\end{eg}

\begin{eg}
  Any groupoid becomes a $\dagger$-category if we define $\dgr f \defeq \inv{f}$.
\end{eg}

\begin{eg}\label{ct:hilb}
  Let \uhilb be the following precategory.
  \begin{itemize}
  \item Its objects are finite-dimensional \index{finite!-dimensional vector space} vector spaces\index{vector!space} equipped with an inner product $\langle \blank,\blank\rangle$.
  \item Its morphisms are arbitrary linear maps.
    \index{function!linear}%
    \indexsee{linear map}{function, linear}%
  \end{itemize}
  By standard linear algebra, any linear map $f:V\to W$ between finite
  dimensional inner product spaces has a uniquely defined adjoint\index{adjoint!linear map} $\dgr f:W\to V$, characterized by $\langle f v,w\rangle = \langle v,\dgr f w\rangle$.
  In this way, \uhilb becomes a $\dagger$-precategory.
  Moreover, a linear isomorphism is unitary precisely when it is an \define{isometry},
  \indexdef{isometry}%
  i.e.\ $\langle fv,fw\rangle = \langle v,w\rangle$.
  It follows from this that \uhilb is a $\dagger$-category, though it is not a category (not every linear isomorphism is unitary).
\end{eg}

There has been a good deal of general theory developed for $\dagger$-cat\-e\-gor\-ies under classical\index{mathematics!classical}\index{classical!category theory} foundations.
It was observed early on that the unitary isomorphisms, not arbitrary isomorphisms, are the correct notion of ``sameness'' for objects of a $\dagger$-category, which has caused some consternation among category theorists.
Homotopy type theory resolves this issue by identifying $\dagger$-categories, like strict categories, as simply a different kind of precategory.

\section{The structure identity principle}
\label{sec:sip}
 \index{structure!identity principle|(}

The \emph{structure identity principle} is an informal principle
that expresses that isomorphic structures are identical.  We aim to
prove a general abstract result which can be applied to a wide family
of notions of structure, where structures may be many-sorted or even
dependently-sorted, infinitary, or even higher order.

The simplest kind of single-sorted structure consists of a type with
no additional structure.  The univalence axiom expresses the structure identity principle for that
notion of structure in a strong form: for types $A,B$, the
canonical function $(A=B)\to (\eqv A B)$ is an equivalence.

We start with a precategory $X$.  In our application to
single-sorted first order structures, $X$ will be the category 
of $\bbU$-small sets, where $\bbU$ is a univalent type universe.

\begin{defn}\label{ct:sig}
  A \define{notion of structure} 
  \indexdef{structure!notion of}%
  $(P,H)$ over $X$ consists of the following.
  \begin{enumerate}
  \item A type family $P:X_0 \to \type$.
    For each $x:X_0$ the elements of $Px$ are called \define{$(P,H)$-structures}
    \indexsee{PH-structure@$(P,H)$-structure}{structure}%
    \indexdef{structure!PH@$(P,H)$-}%
    on $x$.
  \item For $x,y:X_0$ and $\alpha:Px$, $\;\beta:Py$, to each $f:\hom_X(x,y)$ a mere proposition 
  \[ H_{\alpha\beta}(f).\]
    If $H_{\alpha\beta}(f)$ is true, we say that $f$ is a \define{$(P,H)$-homomorphism}
    \indexdef{homomorphism!of structures}%
    \indexdef{structure!homomorphism of}%
    from $\alpha$ to $\beta$.
  \item For $x:X_0$ and $\alpha:Px$, we have $H_{\alpha\alpha}(1_x)$.\label{item:sigid}
  \item For $x,y,z:X_0$ and $\alpha:Px$, $\;\beta:Py$, $\;\gamma:Pz$, 
if $f:\hom_X(x,y)$, we have\label{item:sigcmp}
  \[ H_{\alpha\beta}(f)\to H_{\beta\gamma}(g)\to H_{\alpha\gamma}(g\circ   f).\]
   \end{enumerate}
  When $(P,H)$ is a notion of structure, for $\alpha,\beta:Px$ we define
  \[ (\alpha\leq_x\beta) \defeq H_{\alpha\beta}(1_x).\]
  By~\ref{item:sigid} and~\ref{item:sigcmp}, this is a preorder (\autoref{ct:orders}) with $Px$ its type of objects.
  We say that $(P,H)$ is a \define{standard notion of structure}
  \indexdef{structure!standard notion of}%
  if this preorder is in fact a partial order, for all $x:X$.
\end{defn}

Note that for a standard notion of structure, each type $Px$ must actually be a set.
We now define, for any notion of structure $(P,H)$, a \define{precategory of $(P,H)$-structures},
\indexdef{precategory!of PH-structures@of $(P,H)$-structures}%
\indexdef{structure!precategory of PH@precategory of $(P,H)$-}%
$A = \mathsf{Str}_{(P,H)}(X)$.
\begin{itemize}
\item The type of objects of $A$ is the type $A_0 \defeq \sm{x:X} Px$.
  If $a\jdeq (x,\alpha):A_0$, we may write $|a| \defeq x$.
\item For $(x,\alpha):A_0$ and $(y,\beta):A_0$, we define
  \[\hom_A((x,\alpha),(y,\beta)) \defeq \setof{ f:x \to y | H_{\alpha\beta}(f)}.\]
\end{itemize}
The composition and identities are inherited from $X$; conditions~\ref{item:sigid} and \ref{item:sigcmp} ensure that these lift to $A$.

\begin{thm}[Structure identity principle]\label{thm:sip}
  \indexdef{structure!identity principle}%
  If $X$ is a category and $(P,H)$ is a standard notion of structure over $X$, then the precategory $\mathsf{Str}_{(P,H)}(X)$ is a category.
\end{thm}
\begin{proof}
  By the definition of equality in dependent pair types, to give an equality $(x,\alpha)=(y,\beta)$ consists of
  \begin{itemize}
  \item An equality $p:x=y$, and
  \item An equality $\trans{p}{\alpha}=\beta$.
  \end{itemize}
  Since $P$ is set-valued, the latter is a mere proposition.
  On the other hand, it is easy to see that an isomorphism $(x,\alpha)\cong (y,\beta)$ in $\mathsf{Str}_{(P,H)}(X)$ consists of
  \begin{itemize}
  \item An isomorphism $f:x\cong y$ in $X$, such that
  \item $H_{\alpha\beta}(f)$ and $H_{\beta\alpha}(\inv f)$.
  \end{itemize}
  Of course, the second of these is also a mere proposition.
  And since $X$ is a category, the function $(x=y) \to (x\cong y)$ is an equivalence.
  Thus, it will suffice to show that for any $p:x=y$ and for any $(\alpha:Px)$, $(\beta:Py)$, we have $\trans{p}{\alpha}=\beta$ if and only if both  $H_{\alpha\beta}(\idtoiso (p))$ and $H_{\beta\alpha}(\inv{\idtoiso(p)})$.

  The ``only if'' direction is just the existence of the function $\idtoiso$ for the category $\mathsf{Str}_{(P,H)}(X)$.
  For the ``if'' direction, by induction on $p$ we may assume that $y\jdeq x$ and $p\jdeq\refl x$.
  However, in this case $\idtoiso (p)\jdeq 1_x$ and therefore $\inv{\idtoiso(p)}=1_x$.
  Thus, $\alpha\leq_x \beta$ and $\beta\leq_x \alpha$, which implies $\alpha=\beta$ since $(P,H)$ is a standard notion of structure.
\end{proof}

As an example, this methodology gives an alternative way to express the proof of \autoref{ct:functor-cat}.

\begin{eg}\label{ct:sip-functor-cat}
  Let $A$ be a precategory and $B$ a category.
  There is a precategory $B^{A_0}$ whose objects are functions $A_0 \to B_0$, and whose set of morphisms from $F_0:A_0 \to B_0$ to $G_0:A_0 \to B_0$ is $\prd{a:A_0} \hom_B(F_0 a, G_0 a)$.
  Composition and identities are inherited directly from those in $B$.
  It is easy to show that $\gamma:\hom_{B^{A_0}}(F_0, G_0)$ is an isomorphism exactly when each component $\gamma_a$ is an isomorphism, so that we have $\eqv{(F_0 \cong G_0)}{\prd{a:A_0} (F_0 a \cong G_0 a)}$.
  Moreover, the map $\idtoiso : (F_0 = G_0) \to (F_0 \cong G_0)$ of $B^{A_0}$ is equal to the composite
  \[ (F_0 = G_0) \longrightarrow \prd{a:A_0} (F_0 a  = G_0 a) \longrightarrow \prd{a:A_0} (F_0 a \cong G_0 a) \longrightarrow (F_0 \cong G_0) \]
  in which the first map is an equivalence by function extensionality, the second because it is a dependent product of equivalences (since $B$ is a category), and the third as remarked above.
  Thus, $B^{A_0}$ is a category.

  Now we define a notion of structure on $B^{A_0}$ for which $P(F_0)$ is the type of operations $F:\prd{a,a':A_0} \hom_A(a,a') \to \hom_B(F_0 a,F_0 a')$ which extend $F_0$ to a functor (i.e.\ preserve composition and identities).
  This is a set since each $\hom_B(\blank,\blank)$ is so.
  Given such $F$ and $G$, we define $\gamma:\hom_{B^{A_0}}(F_0, G_0)$ to be a homomorphism if it forms a natural transformation.\index{natural!transformation}
  In \autoref{ct:functor-precat} we essentially verified that this is a notion of structure.
  Moreover, if $F$ and $F'$ are both structures on $F_0$ and the identity is a natural transformation from $F$ to $F'$, then for any $f:\hom_A(a,a')$ we have $F'f = F'f \circ 1_{F_0 a} = 1_{F_0 a}\circ F f = F f$.
  Applying function extensionality, we conclude $F = F'$.
  Thus, we have a \emph{standard} notion of structure, and so by \autoref{thm:sip}, the precategory $B^A$ is a category.
\end{eg}

As another example, we consider categories of structures for a first-order signature.
We define a \define{first-order signature},
\indexdef{first-order!signature}%
\indexdef{signature!first-order}%
$\Omega$, to consist of sets $\Omega_0$ and $\Omega_1$ of function symbols, $\omega:\Omega_0$, and relation symbols, $\omega:\Omega_1$, each having an arity\index{arity} $|\omega|$ that is a set.
An \define{$\Omega$-structure}
\indexdef{structure!Omega@$\Omega$-}%
\indexsee{omega-structure@$\Omega$-structure}{structure}%
$a$ consists of a set $|a|$ together with an assignment of an $|\omega|$-ary function $\omega^a:|a|^{|\omega|}\to |a|$ on $|a|$ to each function symbol, $\omega$, and an assignment of an $|\omega|$-ary relation $\omega^a$ on $|a|$, assigning a mere proposition $\omega^ax$ to each $x:|a|^{|\omega|}$, to each relation symbol.
And given $\Omega$-structures $a,b$, a function $f:|a|\to |b|$ is a \define{homomorphism $a\to b$}
\indexdef{homomorphism!of Omega-structures@of $\Omega$-structures}%
\indexdef{structure!homomorphism of Omega@homomorphism of $\Omega$-}%
if it preserves the structure; i.e.\ if for each symbol $\omega$ of the signature and each $x:|a|^{|\omega|}$,
\begin{enumerate}
\item $f(\omega^ax) = \omega^b(f\circ x)$ if $\omega:\Omega_0$, and
\item $\omega^ax\to\omega^b(f\circ x)$ if $\omega:\Omega_1$.
\end{enumerate}
Note that each $x:|a|^{|\omega|}$ is a function $x:|\omega|\to |a|$ so that $f\circ x : b^\omega$.

Now we assume given a (univalent) universe $\bbU$ and a $\bbU$-small signature $\Omega$; i.e. $|\Omega|$ is a $\bbU$-small set and, for each $\omega:|\Omega|$, the set $|\omega|$ is $\bbU$-small.
Then we have the category $\uset_\bbU$ of $\bbU$-small sets.  We want to define the precategory of $\bbU$-small $\Omega$-structures over $\uset_\bbU$ and use \autoref{thm:sip} to show that it is a category.

We use the first order signature $\Omega$ to give us a standard notion of structure $(P,H)$ over $\uset_\bbU$.  

\begin{defn}\label{defn:fo-notion-of-structure}
\mbox{}
\begin{enumerate}
\item For each $\bbU$-small set $x$ define 
  \[ Px \defeq P_0x\times P_1x.\]
  Here
  \begin{align*}
    P_0x &\defeq \prd{\omega:\Omega_0} x^{|\omega|}\to x, \mbox{ and } \\
    P_1x &\defeq \prd{\omega:\Omega_1} x^{|\omega|}\to \propU,
  \end{align*}
\item For $\bbU$-small sets $x,y$ and 
  $\alpha:P^\omega x,\;\beta:P^\omega y,\; f:x\to y$, define
  \[ H_{\alpha\beta}(f) \defeq H_{0,\alpha\beta}(f)\wedge H_{1,\alpha\beta}(f).\]
  Here
  \begin{align*}
    H_{0,\alpha\beta}(f) &\defeq
    \fall{\omega:\Omega_0}{u:x^{|\omega|}} f(\alpha u)=\;\beta(f\circ u),
    \mbox{ and }\\
    H_{1,\alpha\beta}(f) &\defeq
    \fall{\omega:\Omega_1}{u:x^{|\omega|}} \alpha u\to\beta(f\circ u).
  \end{align*}
\end{enumerate}
\end{defn}

It is now routine to check that $(P,H)$ is a standard notion of structure over $\uset_\bbU$ and hence we may use \autoref{thm:sip} to get that the precategory $Str_{(P,H)}(\uset_\bbU)$ is a category.  It only remains to observe that this is essentially the same as the precategory of $\bbU$-small $\Omega$-structures over $\uset_\bbU$.
 \index{structure!identity principle|)}

\section{The Rezk completion}
\label{sec:rezk}
 
In this section we will give a universal way to replace a precategory by a category.
In fact, we will give two.
Both rely on the fact that ``categories see weak equivalences as equivalences''.

To prove this, we begin with a couple of lemmas which are completely standard category theory, phrased carefully so as to make sure we are using the eliminator for $\truncf{-1}$ correctly.
One would have to be similarly careful in classical\index{mathematics!classical}\index{classical!category theory} category theory if one wanted to avoid the axiom of choice: any time we want to define a function, we need to characterize its values uniquely somehow.

\begin{lem}
  If $A,B,C$ are precategories and $H:A\to B$ is an essentially surjective functor, then $(\blank\circ H):C^B \to C^A$ is faithful.
\end{lem}
\begin{proof}
  Let $F,G:B\to C$, and $\gamma,\delta:F\to G$ be such that $\gamma H = \delta H$; we must show $\gamma=\delta$.
  Thus let $b:B$; we want to show $\gamma_b=\delta_b$.
  This is a mere proposition, so since $H$ is essentially surjective, we may assume given an $a:A$ and an isomorphism $f:Ha\cong b$.
  But now we have
  \[ \gamma_b = G(f) \circ \gamma_{Ha} \circ F(\inv{f}) 
  = G(f) \circ \delta_{Ha} \circ F(\inv{f})
  = \delta_b.\qedhere
  \]
\end{proof}

\begin{lem}\label{ct:esofull-precomp-ff}
  If $A,B,C$ are precategories and $H:A\to B$ is essentially surjective and full, then $(\blank\circ H):C^B \to C^A$ is fully faithful.
\end{lem}
\begin{proof}
  It remains to show fullness.
  Thus, let $F,G:B\to C$ and $\gamma:FH \to GH$.
  We claim that for any $b:B$, the type
  \begin{equation}\label{eq:fullprop}
    \sm{g:\hom_C(Fb,Gb)} \prd{a:A}{f:Ha\cong b} (\gamma_a =  \inv{Gf}\circ g\circ Ff)
  \end{equation}
  is contractible.
  Since contractibility is a mere property, and $H$ is essentially surjective, we may assume given $a_0:A$ and $h:Ha_0\cong b$.

  Now take $g\defeq Gh \circ \gamma_{a_0} \circ \inv{Fh}$.
  Then given any other $a:A$ and $f:Ha\cong b$, we must show $\gamma_a =  \inv{Gf}\circ g\circ Ff$.
  Since $H$ is full, there merely exists a morphism $k:\hom_A(a,a_0)$ such that $Hk = \inv{h}\circ f$.
  And since our goal is a mere proposition, we may assume given some such $k$.
  Then we have
  \begin{align*}
    \gamma_a &= \inv{GHk}\circ \gamma_{a_0} \circ FHk\\
    &= \inv{Gf} \circ Gh \circ \gamma_{a_0} \circ \inv{Fh} \circ Ff\\
    &= \inv{Gf}\circ g\circ Ff.
  \end{align*}
  Thus,~\eqref{eq:fullprop} is inhabited.
  It remains to show it is a mere proposition.
  Let $g,g':\hom_C(Fb, Gb)$ be such that for all $a:A$ and $f:Ha\cong b$, we have both $(\gamma_a =  \inv{Gf}\circ g\circ Ff)$ and $(\gamma_a =  \inv{Gf}\circ g'\circ Ff)$.
  The dependent product types are mere propositions, so all we have to prove is $g=g'$.
  But this is a mere proposition, so we may assume $a_0:A$ and $h:Ha_0\cong b$, in which case we have
  \[ g = Gh \circ \gamma_{a_0} \circ \inv{Fh} = g'.\]
  This proves that~\eqref{eq:fullprop} is contractible for all $b:B$.
  Now we define $\delta:F\to G$ by taking $\delta_b$ to be the unique $g$ in~\eqref{eq:fullprop} for that $b$.
  To see that this is natural, suppose given $f:\hom_B(b,b')$; we must show $Gf \circ \delta_b = \delta_{b'}\circ Ff$.
  As before, we may assume $a:A$ and $h:Ha\cong b$, and likewise $a':A$ and $h':Ha'\cong b'$.
  Since $H$ is full as well as essentially surjective, we may also assume $k:\hom_A(a,a')$ with $Hk = \inv{h'}\circ f\circ h$.

  Since $\gamma$ is natural, $GHk\circ \gamma_a = \gamma_{a'} \circ FHk$.
  Using the definition of $\delta$, we have
  \begin{align*}
    Gf \circ \delta_b
    &= Gf \circ Gh \circ \gamma_a \circ \inv{Fh}\\
    &= Gh' \circ GHk\circ \gamma_a \circ \inv{Fh}\\
    &= Gh' \circ \gamma_{a'} \circ FHk \circ \inv{Fh}\\
    &= Gh' \circ \gamma_{a'} \circ \inv{Fh'} \circ Ff\\
    &= \delta_{b'} \circ Ff.
  \end{align*}
  Thus, $\delta$ is natural.
  Finally, for any $a:A$, applying the definition of $\delta_{Ha}$ to $a$ and $1_a$, we obtain $\gamma_a = \delta_{Ha}$.
  Hence, $\delta \circ H = \gamma$.
\end{proof}

The rest of the theorem follows almost exactly the same lines, with the category-ness of $C$ inserted in one crucial step, which we have italicized below for emphasis.
This is the point at which we are trying to define a function into \emph{objects} without using choice, and so we must be careful about what it means for an object to be ``uniquely specified''.
In classical\index{mathematics!classical}\index{classical!category theory} category theory, all one can say is that this object is specified up to unique isomorphism, but in set-theoretic foundations this is not a sufficient amount of uniqueness to give us a function without invoking \choice{}.
In univalent foundations, however, if $C$ is a category, then isomorphism is equality, and we have the appropriate sort of uniqueness (namely, living in a contractible space).

\index{weak equivalence!of precategories|(}%

\begin{thm}\label{ct:cat-weq-eq}
  If $A,B$ are precategories, $C$ is a category, and $H:A\to B$ is a weak equivalence, then $(\blank\circ H):C^B \to C^A$ is an isomorphism.
\end{thm}
\begin{proof}
  By \autoref{ct:functor-cat}, $C^B$ and $C^A$ are categories.
  Thus, by \autoref{ct:eqv-levelwise} it will suffice to show that $(\blank\circ H)$ is an equivalence.
  But since we know from the preceding two lemmas that it is fully faithful, by \autoref{ct:catweq} it will suffice to show that it is essentially surjective.
  Thus, suppose $F:A\to C$; we want there to merely exist a $G:B\to C$ such that $GH\cong F$.

  For each $b:B$, let $X_b$ be the type whose elements consist of:
  \begin{enumerate}
  \item An element $c:C$; and
  \item For each $a:A$ and $h:Ha\cong b$, an isomorphism $k_{a,h}:Fa\cong c$; such that\label{item:eqvprop2}
  \item For each $(a,h)$ and $(a',h')$ as in~\ref{item:eqvprop2} and each $f:\hom_A(a,a')$ such that $h'\circ Hf = h$, we have $k_{a',h'}\circ Ff = k_{a,h}$.\label{item:eqvprop3}
  \end{enumerate}
  We claim that for any $b:B$, the type $X_b$ is contractible.
  As this is a mere proposition, we may assume given $a_0:A$ and $h_0:Ha_0 \cong b$.
  Let $c^0\defeq Fa_0$.
  Next, given $a:A$ and $h:Ha\cong b$, since $H$ is fully faithful there is a unique isomorphism $g_{a,h}:a\to a_0$ with $Hg_{a,h} = \inv{h_0}\circ h$; define $k^0_{a,h} \defeq Fg_{a,h}$.
  Finally, if $h'\circ Hf = h$, then $\inv{h_0}\circ h'\circ Hf = \inv{h_0}\circ h$, hence $g_{a',h'} \circ f = g_{a,h}$ and thus $k^0_{a',h'}\circ Ff = k^0_{a,h}$.
  Therefore, $X_b$ is inhabited.

  Now suppose given another $(c^1,k^1): X_b$.
  Then $k^1_{a_0,h_0}:c^0 \jdeq Fa_0 \cong c^1$.
  \emph{Since $C$ is a category, we have $p:c^0=c^1$ with $\idtoiso(p) = k^1_{a_0,h_0}$.}
  And for any $a:A$ and $h:Ha\cong b$, by~\ref{item:eqvprop3} for $(c^1,k^1)$ with $f\defeq g_{a,h}$, we have
  \[k^1_{a,h} = k^1_{a_0,h_0} \circ k^0_{a,h} = \trans{p}{k^0_{a,h}}\]
  This gives the requisite data for an equality $(c^0,k^0)=(c^1,k^1)$, completing the proof that $X_b$ is contractible.

  Now since $X_b$ is contractible for each $b$, the type $\prd{b:B} X_b$ is also contractible.
  In particular, it is inhabited, so we have a function assigning to each $b:B$ a $c$ and a $k$.
  Define $G_0(b)$ to be this $c$; this gives a function $G_0 :B_0 \to C_0$.

  Next we need to define the action of $G$ on morphisms.
  For each $b,b':B$ and $f:\hom_B(b,b')$, let $Y_f$ be the type whose elements consist of:
  \begin{enumerate}[resume]
  \item A morphism $g:\hom_C(Gb,Gb')$, such that
  \item For each $a:A$ and $h:Ha\cong b$, and each $a':A$ and $h':Ha'\cong b'$, and any $\ell:\hom_A(a,a')$, we have\label{item:eqvprop5}
    \[ (h' \circ H\ell = f \circ h)
    \to
    (k_{a',h'} \circ F\ell = g\circ k_{a,h}). \]
  \end{enumerate}
  We claim that for any $b,b'$ and $f$, the type $Y_f$ is contractible.
  As this is a mere proposition, we may assume given $a_0:A$ and $h_0:Ha_0\cong b$, and each $a'_0:A$ and $h'_0:Ha'_0\cong b'$.
  Then since $H$ is fully faithful, there is a unique $\ell_0:\hom_A(a_0,a_0')$ such that $h'_0 \circ H\ell_0 = f \circ h_0$.
  Define $g_0 \defeq k_{a_0',h_0'} \circ F \ell_0 \circ \inv{(k_{a_0,h_0})}$.

  Now for any $a,h,a',h'$, and $\ell$ such that $(h' \circ H\ell = f \circ h)$, we have $\inv{h}\circ h_0:Ha_0\cong Ha$, hence there is a unique $m:a_0\cong a$ with $Hm = \inv{h}\circ h_0$ and hence $h\circ Hm = h_0$.
  Similarly, we have a unique $m':a_0'\cong a'$ with $h'\circ Hm' = h_0'$.
  Now by~\ref{item:eqvprop3}, we have $k_{a,h}\circ Fm = k_{a_0,h_0}$ and $k_{a',h'}\circ Fm' = k_{a_0',h_0'}$.
  We also have
  \begin{align*}
    Hm' \circ H\ell_0 
    &= \inv{(h')} \circ h_0' \circ H\ell_0\\
    &= \inv{(h')} \circ f \circ h_0\\
    &= \inv{(h')} \circ f \circ h \circ \inv{h} \circ h_0\\
    &= H\ell \circ Hm
  \end{align*}
  and hence $m'\circ \ell_0 = \ell\circ m$ since $H$ is fully faithful.
  Finally, we can compute
  \begin{align*}
    g_0 \circ k_{a,h}
    &= k_{a_0',h_0'} \circ F \ell_0 \circ \inv{(k_{a_0,h_0})} \circ k_{a,h}\\
    &= k_{a_0',h_0'} \circ F \ell_0 \circ \inv{Fm}\\
    &= k_{a_0',h_0'} \circ \inv{(Fm')} \circ F\ell\\
    &= k_{a',h'}\circ F\ell.
  \end{align*}
  This completes the proof that $Y_f$ is inhabited.
  To show it is contractible, since hom-sets are sets, it suffices to take another $g_1:\hom_C(Gb,Gb')$ satisfying~\ref{item:eqvprop5} and show $g_0=g_1$.
  However, we still have our specified $a_0,h_0,a_0',h_0',\ell_0$ around, and~\ref{item:eqvprop5} implies both $g_0$ and $g_1$ must be equal to $k_{a_0',h_0'} \circ F \ell_0 \circ \inv{(k_{a_0,h_0})}$.

  This completes the proof that $Y_f$ is contractible for each $b,b':B$ and $f:\hom_B(b,b')$.
  Therefore, there is a function assigning to each such $f$ its unique inhabitant; denote this function $G_{b,b'}:\hom_B(b,b') \to \hom_C(Gb,Gb')$.
  The proof that $G$ is a functor is straightforward; in each case we can choose $a,h$ and apply~\ref{item:eqvprop5}.

  Finally, for any $a_0:A$, defining $c\defeq Fa_0$ and $k_{a,h}\defeq F g$, where $g:\hom_A(a,a_0)$ is the unique isomorphism with $Hg = h$, gives an element of $X_{Ha_0}$.
  Thus, it is equal to the specified one; hence $GHa=Fa$.
  Similarly, for $f:\hom_A(a_0,a_0')$ we can define an element of $Y_{Hf}$ by transporting along these equalities, which must therefore be equal to the specified one.
  Hence, we have $GH=F$, and thus $GH\cong F$ as desired.
\end{proof}

\index{universal!property!of Rezk completion}%
Therefore, if a precategory $A$ admits a weak equivalence functor $A\to \widehat{A}$, then that is its ``reflection'' into categories: any functor from $A$ into a category will factor essentially uniquely through $\widehat{A}$.
We now give two constructions of such a weak equivalence.

\indexsee{Rezk completion}{completion, Rezk}%
\index{completion!Rezk|(defstyle}%

\begin{thm}
  For any precategory $A$, there is a category $\widehat A$ and a weak equivalence $A\to\widehat{A}$.
\end{thm}

\begin{proof}[First proof]
  Let $\widehat{A}_0 \defeq \setof{ F:\uset^{A\op} | \exis{a:A} (\y a \cong F)}$, with hom-sets inherited from $\uset^{A\op}$.
  Then the inclusion $\widehat{A} \to \uset^{A\op}$ is fully faithful and an embedding on objects.
  Since $\uset^{A\op}$ is a category (by \autoref{ct:functor-cat}, since \uset is so by univalence), $\widehat A$ is also a category.

  Let $A\to\widehat A$ be the Yoneda embedding.
  This is fully faithful by \autoref{ct:yoneda-embedding}, and essentially surjective by definition of $\widehat{A}_0$.
  Thus it is a weak equivalence.
\end{proof}

This proof is very slick, but it has the drawback that it increases universe level.
If $A$ is a category in a universe \bbU, then in this proof \uset must be at least as large as $\uset_\bbU$.
Then $\uset_\bbU$ and $(\uset_\bbU)^{A\op}$ are not themselves categories in \bbU, but only in a higher universe, and \emph{a priori} the same is true of $\widehat A$.
One could imagine a resizing axiom that could deal with this, but it is also possible to give a direct construction using higher inductive types.

\begin{proof}[Second proof]
  We define a higher inductive 1-type $\widehat A_0$ with the following constructors:
  \begin{itemize}
  \item A function $i:A_0 \to \widehat A_0$.
  \item For each $a,b:A$ and $e:a\cong b$, an equality $je:\id{ia}{ib}$.
  \item For each $a,b:A$ and $p:\id a b$, an equality $\id{j(\idtoiso(p))}{\map i p}$.
  \item For each $a:A$, an equality $\id{j(1_a)}{\refl{ia}}$.
  \item For each $(a,b,c:A)$, $(f:a\cong b)$, and $(g:b\cong c)$, an equality $\id{j(g \circ f)}{j(g)\ct j(f)}$.
  \end{itemize}
  This will be the type of objects of $\widehat A$; we now build all the rest of the structure.
  (The following proof is of the sort that can benefit a lot from the help of a computer proof assistant:\index{proof!assistant} it is wide and shallow with many short cases to consider, and a large part of the work consists of writing down what needs to be checked.)

  \mentalpause

  \emph{Step 1:} We define a family $\hom_{\widehat A}:\widehat A_0\to \widehat A_0 \to \set$ by double induction on $\widehat A_0$, which is possible since \set is a 1-type.
  When $x$ and $y$ are of the form $ia$ and $ib$, we take $\hom_{\widehat A}(ia,ib) \defeq \hom_A(a,b)$.
  It remains to consider all the other possible pairs of constructors.

  Let us keep $x=ia$ fixed at first.
  If $y$ varies along the identity $je:\id{ib}{ib'}$, for some $e:b\cong b'$, we require an identity $\id{\hom_A(a,b)}{\hom_A(a,b')}$.
  By univalence, it suffices to give an equivalence $\eqv{\hom_A(a,b)}{\hom_A(a,b')}$.
  We take this to be the function $(e\circ \blank ):\hom_A(a,b)\to \hom_A(a,b')$.
  To see that this is an equivalence, we give its inverse as $(\inv e\circ \blank )$, with witnesses to inversion coming from the fact that $\inv e$ is the inverse of $e$ in $A$.
  
  Next, as $y$ varies along the identity $\id{j(\idtoiso(p))}{\map ip}$, for $p:\id{b}{b'}$, we require an identity $\id{(\idtoiso(p)\circ \blank )}{\map{\hom_A(a,\blank)}{p}}$.
  This is immediate by induction on $p$.

  As $y$ varies along the identity $\id{j(1_b)}{\refl{ib}}$, we require an identity $\id{(1_b\circ \blank )}{\refl{\hom_A(a,b)}}$; this follows from the identity axiom $\id{1_b\circ g}{g}$ of a precategory.
  Similarly, as $y$ varies along the identity $\id{j(g\circ f)}{j(g)\ct j(f)}$, we require an identity $\id{((g\circ f)\circ \blank )}{(g\circ (f\circ \blank ))}$, which follows from associativity.

  Now we consider the other constructors for $x$.
  Say that $x$ varies along the identity $j(e):\id{ia}{ia'}$, for some $e:a \cong a'$; we again must deal with all the constructors for $y$.
  If $y$ is $ib$, then we require an identity $\id{\hom_A(a,b)}{\hom_A(a',b)}$.
  By univalence, this may come from an equivalence, and for this we can use $(\blank\circ \inv e)$, with inverse $(\blank\circ e)$.

  Still with $x$ varying along $j(e)$, suppose now that $y$ also varies along $j(f)$ for some $f:b\cong b'$.
  Then we need to know that the two concatenated identities
  \begin{gather*}
    \hom_A(a,b) = \hom_A(a',b) = \hom_A(a',b') \mathrlap{\qquad\text{and}}\\
    \hom_A(a,b) = \hom_A(a,b') = \hom_A(a',b')
  \end{gather*}
  are identical.
  This follows from associativity: $(f\circ \blank)\circ \inv e = f\circ (\blank\circ \inv e)$.
  The rest of the constructors for $y$ are trivial, since they are 2-fold equalities in sets.

  For the last three constructors of $x$, all but the first constructor for $y$ is likewise trivial.
  When $x$ varies along the equality $\id{j(\idtoiso(p))}{\map i p}$ for $p:a=a'$ and $y$ is $ib$, we require $(\blank\circ \idtoiso(p)) = \map{\hom_A(\blank,b)}{p}$, which follows by induction on $p$.
  Finally, when $x$ varies along $j(1_a)=\refl{ia}$, we use the identity axiom again, and when $x$ varies along $\id{j(g\circ f)}{j(g)\ct j(f)}$, we use associativity again.
  This completes the construction of $\hom_{\widehat A}:\widehat A_0 \to \widehat A_0 \to \set$.

  \mentalpause

  \emph{Step 2:} We give the precategory structure on $\widehat A$, always by induction on $\widehat A_0$.
  We are now eliminating into sets (the hom-sets of $\widehat A$), so all but the first two constructors are trivial to deal with.

  For identities, if $x$ is $ia$ then we have $\hom_{\widehat A}(x,x) \jdeq \hom_A(a,a)$ and we define $1_x \defeq 1_{ia}$.
  If $x$ varies along $je$ for $e:a\cong a'$, we must show that $\trans{je}{1_{ia}} = 1_{ia'}$.
  Here the transport is with respect to the type family $x\mapsto \hom_{\widehat A}(x,x)$.
  But by definition of $\hom_{\widehat A}$, transporting along $je$ is given by composing with $e$ and $\inv e$, and we have $e\circ 1_{ia} \circ \inv{e} = 1_{ia'}$.

  For composition, if $x,y,z$ are $ia,ib,ic$ respectively, then $\hom_{\widehat A}$ reduces to $\hom_A$ and we can define composition in $\widehat A$ to be composition in $A$.
  And when $x$, $y$, or $z$ varies along $je$, then we verify the following equalities:
  \begin{align*}
    e \circ (g\circ f) &= (e\circ g) \circ f,\\
    g\circ f &= (g\circ \inv e) \circ (e\circ f),\\
    (g\circ f) \circ \inv e &= g \circ (f\circ \inv e).
  \end{align*}
  Finally, the associativity and unitality axioms are mere propositions, so all constructors except the first are trivial.
  But in that case, we have the corresponding axioms in $A$.

  \mentalpause

  \emph{Step 3}: We show that $\widehat A$ is a category.
  That is, we must show that for all $x,y:\widehat A$, the function $\idtoiso:(x=y) \to (x\cong y)$ is an equivalence.
  First we define, for all $x,y:\widehat A$, a function $k_{x,y}:(x\cong y) \to (x=y)$ by induction.
  As before, since our goal is a set, it suffices to deal with the first two constructors.

  When $x$ and $y$ are $ia$ and $ib$ respectively, we have $\hom_{\widehat A}(ia,ib)\jdeq \hom_A(a,b)$, with composition and identities inherited as well, so that $(ia\cong ib)$ is equivalent to $(a\cong b)$.
  But now we have the constructor $j:(a\cong b) \to (ia=ib)$.

  Next, if $y$ varies along $j(e)$ for some $e:b\cong b'$, we must show that for $f:a\cong b$ we have $j(\trans{j(e)}{f}) = j(e) \ct j(f)$.
  But by definition of $\hom_{\widehat A}$ on equalities, transporting along $j(e)$ is equivalent to post-composing with $e$, so this equality follows from the last constructor of $\widehat A_0$.
  The remaining case when $x$ varies along $j(e)$ for $e:a\cong a'$ is similar.
  This completes the definition of $k:\prd{x,y:\widehat A_0} (x\cong y) \to (x=y)$.

  Now one thing we must show is that if $p:x=y$, then $k(\idtoiso(p))=p$.
  By induction on $p$, we may assume it is $\refl x$, and hence $\idtoiso(p)\jdeq 1_x$.
  Now we argue by induction on $x:\widehat A_0$, and since our goal is a mere proposition (since $\widehat A_0$ is a 1-type), all constructors except the first are trivial.
  But if $x$ is $ia$, then $k(1_{ia}) \jdeq j(1_a)$, which is equal to $\refl{ia}$ by the penultimate constructor of $\widehat A_0$.

  To complete the proof that $\widehat A$ is a category, we must show that if $f:x\cong y$, then $\idtoiso(k(f))=f$.
  By induction we may assume that $x$ and $y$ are $ia$ and $ib$ respectively, in which case $f$ must arise from an isomorphism $g:a\cong b$ and we have $k(f)\jdeq j(g)$.
  However, for any $p$ we have $\idtoiso(p) = \trans{p}{1}$, so in particular $\idtoiso (j(g)) = \trans{j(g)}{1_{ia}}$.
  And by definition of $\hom_{\widehat A}$ on equalities, this is given by composing $1_{ia}$ with the equivalence $g$, hence is equal to $g$.

  \index{encode-decode method}%
  Note the similarity of this step to the encode-decode method\index{encode-decode method} used in \autoref{sec:compute-coprod,sec:compute-nat,cha:homotopy}.
  Once again we are characterizing the identity types of a higher inductive type (here, $\widehat A_0$) by defining recursively a family of codes (here, $(x,y)\mapsto (x\cong y)$) and encoding and decoding functions by induction on $\widehat A_0$ and on paths.

  \mentalpause

  \emph{Step 4}: We define a weak equivalence $I:A \to \widehat A$.
  We take $I_0 \defeq i : A_0 \to \widehat A_0$, and by construction of $\hom_{\widehat A}$ we have functions $I_{a,b}:\hom_A(a,b) \to \hom_{\widehat A}(Ia,Ib)$ forming a functor $I:A \to \widehat A$.
  This functor is fully faithful by construction, so it remains to show it is essentially surjective.
  That is, for all $x:\widehat A$ we want there to merely exist an $a:A$ such that $Ia\cong x$.
  As always, we argue by induction on $x$, and since the goal is a mere proposition, all but the first constructor are trivial.
  But if $x$ is $ia$, then of course we have $a:A$ and $Ia\jdeq ia$, hence $Ia \cong ia$.
  (Note that if we were trying to prove $I$ to be \emph{split} essentially surjective, we would be stuck, because we know nothing about equalities in $A_0$ and thus have no way to deal with any further constructors.)
\end{proof}

We call the construction $A\mapsto \widehat A$ the \define{Rezk completion},
although there is also an argument (coming from higher topos semantics)
\index{.infinity1-topos@$(\infty,1)$-topos}%
for calling it the \define{stack completion}.
\index{stack}%
\index{completion!Rezk|)}%

We have seen that most precategories arising in practice are categories, since they are constructed from \uset, which is a category by the univalence axiom.
However, there are a few cases in which the Rezk completion is necessary to obtain a category.

\begin{eg}\label{ct:rezk-fundgpd-trunc1}
  Recall from \autoref{ct:fundgpd} that for any type $X$ there is a pregroupoid with $X$ as its type of objects and $\hom(x,y) \defeq \pizero{x=y}$.
  \indexdef{fundamental!groupoid}%
  \index{fundamental!pregroupoid}%
  \indexsee{groupoid!fundamental}{fundamental group\-oid}%
  Its Rezk completion is the \emph{fundamental groupoid} of $X$.
  Recalling that groupoids are equivalent to 1-types, it is not hard to identify this groupoid with $\trunc1X$.
\end{eg}

\begin{eg}\label{ct:hocat}
  Recall from \autoref{ct:hoprecat} that there is a precategory whose type of objects is \type and with $\hom(X,Y) \defeq \pizero{X\to Y}$.
  Its Rezk completion may be called the \define{homotopy category of types}.
  \index{category!of types}%
  \index{homotopy!category of types@(pre)category of types}%
  Its type of objects can be identified with $\trunc1\type$ (see \autoref{ct:ex:hocat}).
\end{eg}

The Rezk completion also allows us to show that the notion of ``category'' is determined by the notion of ``weak equivalence of precategories''.
Thus, insofar as the latter is inevitable, so is the former.

\begin{thm}
  A precategory $C$ is a category if and only if for every weak equivalence of precategories $H:A\to B$, the induced functor $(\blank\circ H):C^B \to C^A$ is an isomorphism of precategories.
\end{thm}
\begin{proof}
  ``Only if'' is \autoref{ct:cat-weq-eq}.
  In the other direction, let $H$ be $I:A\to\widehat A$.
  Then since $(\blank\circ I)_0$ is an equivalence, there exists $R:\widehat A\to A$ such that $RI=1_A$.
  Hence $IRI=I$, but again since $(\blank\circ I)_0$ is an equivalence, this implies $IR =1_{\widehat A}$.
  By \autoref{ct:isoprecat}\ref{item:ct:ipc3}, $I$ is an isomorphism of precategories.
  But then since $\widehat A$ is a category, so is $A$.
\end{proof}

\index{weak equivalence!of precategories|)}%

\newpage

\sectionNotes

The original definition of categories, of course, was in set-theoretic foundations, so that the collection of objects of a category formed a set (or, for large categories, a class).
Over time, it became clear that all ``category-theoretic'' properties of objects were invariant under isomorphism, and that equality of objects in a category was not usually a very useful notion.
Numerous authors~\cite{blanc:eqv-log,freyd:invar-eqv,makkai:folds,makkai:comparing} discovered that a dependently typed logic enabled formulating the definition of category without invoking any notion of equality for objects, and that the statements provable in this logic are precisely the ``category-theoretic'' ones that are invariant under isomorphism.
\index{evil}%

Although most of category theory appears to be invariant under isomorphism of objects and under equivalence of categories, there are some interesting exceptions, which have led to philosophical discussions about what it means to be ``category-theoretic''.
For instance, \autoref{ct:galois} was brought up by Peter May on the categories mailing list in May 2010, as a case where it matters that two categories (defined as usual in set theory) are isomorphic rather than only equivalent.
The case of $\dagger$-categories was also somewhat confounding to those advocating an isomorphism-invariant version of category theory, since the ``correct'' notion of sameness between objects of a $\dagger$-category is not ordinary isomorphism but \emph{unitary} isomorphism.
\index{isomorphism!invariance under}%

The fact that categories satisfying the ``saturation'' or ``univalence'' principle as in \autoref{ct:category} are a good notion of category in univalent foundations occurred independently to Voevodsky, Shulman, and perhaps others around the same time, and was formalized by Ahrens and Kapulkin~\cite{aks:rezk}.
This framework puts all the above examples in a unified context: some precategories are categories, others are strict categories, and so on.
A general theorem that ``isomorphism implies equality'' for a large class of algebraic structures (assuming the univalence axiom) was proven by Coquand and Danielsson; the formulation of the structure identity principle in \autoref{sec:sip} is due to Aczel.

Independently of philosophical considerations about category theory, Rezk~\cite{rezk01css} discovered that when defining a notion of $(\infty,1)$-cat\-e\-go\-ry,
\index{.infinity1-category@$(\infty,1)$-category}%
it was very convenient to use not merely a \emph{set} of objects with spaces of morphisms between them, but a \emph{space} of objects incorporating all the equivalences and homotopies between them.
This yields a very well-behaved sort of model for $(\infty,1)$-categories as particular simplicial spaces, which Rezk called \emph{complete Segal spaces}.
\index{complete!Segal space}%
\index{Segal!space}%
One especially good aspect of this model is the analogue of \autoref{ct:eqv-levelwise}: a map of complete Segal spaces is an equivalence just when it is a levelwise equivalence of simplicial spaces.

When interpreted in Voevodsky's simplicial\index{simplicial!sets} set model of univalent foundations, our precategories are similar to a truncated analogue of Rezk's ``Segal spaces'', while our categories correspond to his ``complete Segal spaces''.
\index{Segal!category}%
Strict categories correspond instead to (a weakened and truncated version of) what are called ``Segal categories''.
It is known that Segal categories and complete Segal spaces are equivalent models for $(\infty,1)$-categories (see e.g.~\cite{bergner:infty-one}), so that in the simplicial set model, categories and strict categories yield ``equivalent'' category theories---although as we have seen, the former still have many advantages.
However, in the more general categorical semantics of a higher topos,
\index{.infinity1-topos@$(\infty,1)$-topos}%
a strict category corresponds to an internal category (in the traditional sense) in the corresponding 1-topos\index{topos} of sheaves, while a category corresponds to a \emph{stack}.
\index{stack}%
The latter are generally a more appropriate sort of ``category'' relative to a topos.

In Rezk's context, what we have called the ``Rezk completion'' corresponds to fibrant replacement
\index{fibrant replacement}
in the model category for complete Segal spaces.
Since this is built using a transfinite induction argument, it most closely matches our second construction as a higher inductive type.
However, in higher topos models of homotopy type theory, the Rezk completion corresponds to \emph{stack completion},\index{completion!stack}\index{stack!completion} which can be constructed either with a transfinite induction~\cite{jt:strong-stacks} or using a Yoneda embedding \cite{bunge:stacks-morita-internal}.

\sectionExercises

\begin{ex}
  For a precategory $A$ and $a:A$, define the \define{slice precategory} $A/a$.
  \indexsee{precategory!slice}{category, slice}%
  \indexsee{slice (pre)category}{category, slice}%
  Show that if $A$ is a category, so is $A/a$.
  \indexdef{category!slice}%
\end{ex}

\begin{ex}
  For any set $X$, prove that the slice category $\uset/X$ is equivalent to the functor category $\uset^X$, where in the latter case we regard $X$ as a discrete category.
\end{ex}

\begin{ex}
  \index{adjoint!functor}%
  \index{adjoint!equivalence}%
  Prove that a functor is an equivalence of categories if and only if it is a \emph{right} adjoint whose unit and counit are isomorphisms.
\end{ex}

\begin{ex}\label{ct:pre2cat}
  Define a \define{pre-2-category}
  \indexdef{pre-2-category}%
  to consist of the structure formed by precategories, functors, and natural transformations in \autoref{sec:transfors}.
  Similarly, define a \define{pre-bicategory}
  \indexdef{pre-bicategory}%
  by replacing the equalities in \autoref{ct:functor-assoc,ct:units} with natural isomorphisms satisfying analogous coherence conditions.
  Define a function from pre-2-categories to pre-bicategories, and show that it becomes an equivalence when restricted and corestricted to those whose hom-pre\-cat\-egories are categories.
\end{ex}

\begin{ex}\label{ct:2cat}
  Define a \define{2-category}
  \indexdef{2-category}%
  to be a pre-2-category satisfying a condition analogous to that of \autoref{ct:category}.
  Verify that the pre-2-category of categories \ucat is a 2-category.
  How much of this chapter can be done internally to an arbitrary 2-category?
\end{ex}

\begin{ex}\label{ct:groupoids}
  Define a 2-category whose objects are 1-types, whose morphisms are functions, and whose 2-morphisms are homotopies.
  Prove that it is equivalent, in an appropriate sense, to the full sub-2-category of \ucat spanned by the \emph{groupoids} (categories in which every arrow is an isomorphism).
\end{ex}

\begin{ex}
  \index{strict!category}%
  Recall that a \emph{strict category} is a precategory whose type of objects is a set.
  Prove that the pre-2-category of strict categories is equivalent to the following pre-2-category.
  \begin{itemize}
  \item Its objects are categories $A$ equipped with a surjection
    $p_A:A_0'\to A_0$, where $A_0'$ is a set.
  \item Its morphisms are functors $F:A\to B$ equipped with a function $F_0':A_0' \to B_0'$ such that $p_B \circ F_0' = F_0 \circ p_A$.
  \item Its 2-morphisms are simply natural transformations.
  \end{itemize}
\end{ex}

\begin{ex}
  Define the pre-2-category of $\dagger$-categories, which has $\dagger$-struc\-tures on its hom-pre\-cat\-egories.
  Show that two $\dagger$-categories are equal precisely when they are ``unitarily equivalent'' in a suitable sense.
\end{ex}

\begin{ex}\label{ct:ex:hocat}
  Prove that a function $X\to Y$ is an equivalence if and only if its image in the homotopy category of \autoref{ct:hocat} is an isomorphism.
  Show that the type of objects of this category is $\trunc1\type$.
\end{ex}

\begin{ex}
  Construct the $\dagger$-Rezk completion of a $\dagger$-precategory into a $\dagger$-category, and give it an appropriate universal property.
\end{ex}

\begin{ex}\label{ex:rezk-vankampen}
  \index{van Kampen theorem}%
  \index{theorem!van Kampen}%
  \index{fundamental!groupoid}%
  \index{fundamental!pregroupoid}%
  Using fundamental (pre)groupoids from \autoref{ct:fundgpd,ct:rezk-fundgpd-trunc1} and the Rezk completion from \autoref{sec:rezk}, give a different proof of van Kampen's theorem (\autoref{sec:van-kampen}).
\end{ex}

\begin{ex}\label{ex:stack}
  Let $X$ and $Y$ be sets and $p:Y\to X$ a surjection.
  \begin{enumerate}
  \item Define, for any precategory $A$, the category $\mathrm{Desc}(A,p)$ of \define{descent data}
    \indexdef{descent data}%
    in $A$ relative to $p$.
  \item Show that any precategory $A$ is a \define{prestack}
    \indexdef{prestack}%
    for $p$, i.e.\ the canonical functor $A^X \to \mathrm{Desc}(A,p)$ is fully faithful.
  \item Show that if $A$ is a category, then it is a \define{stack}
    \indexdef{stack}%
    for $p$, i.e.\ $A^X \to \mathrm{Desc}(A,p)$ is an equivalence.
  \item Show that the statement ``every strict category is a stack for every surjection of sets'' is equivalent to the axiom of choice.
    \index{axiom!of choice}%
    \index{strict!category}%
  \end{enumerate}
\end{ex}


\chapter{Set theory}
\label{cha:set-math}

\index{set|(}%

Our conception of sets as types with particularly simple homotopical character, cf.\
\autoref{sec:basics-sets}, is quite different from the sets of Zermelo--Fraenkel\index{set theory!Zermelo--Fraenkel} set theory, which form a
cumulative hierarchy with an intricate nested membership structure.
For many mathematical purposes, the homotopy-the\-o\-ret\-ic sets are just as good as
the Zermelo--Fraenkel ones, but there are important differences.

We begin this chapter in \autoref{sec:piw-pretopos} by showing that the category $\uset$ has (most of) the usual properties of the category of sets.
\index{mathematics!constructive}%
\index{mathematics!predicative}%
In constructive, predicative, univalent foundations, it is a ``$\Pi\mathsf{W}$-pretopos''; whereas if we assume propositional resizing
\index{propositional!resizing}%
(\autoref{subsec:prop-subsets}) it is an elementary topos,\index{topos} and if we assume \LEM{} and \choice{} then it is a model of Lawvere's \emph{Elementary Theory of the Category of Sets}\index{Lawvere}.
\index{Elementary Theory of the Category of Sets}%
This is sufficient to ensure that the sets in homotopy type theory behave like sets as used by most mathematicians outside of set theory.

In the rest of the chapter, we investigate some subjects that traditionally belong to ``set theory''.
In \autoref{sec:cardinals,sec:ordinals,sec:wellorderings} we study cardinal and ordinal numbers.
These are traditionally defined in set theory using the global membership relation, but we will see that the univalence axiom enables an equally convenient, more ``structural'' approach.

Finally, in \autoref{sec:cumulative-hierarchy} we consider the possibility of constructing \emph{inside} of homotopy type theory a cumulative hierarchy of sets, equipped with a binary membership relation akin to that of Zermelo--Fraenkel set theory.
This combines higher inductive types with ideas from the field of algebraic set theory.
\index{algebraic set theory}%
\index{set theory!algebraic}%

In this chapter we will often use the traditional logical notation described in \autoref{subsec:prop-trunc}.
In addition to the basic theory of \autoref{cha:basics,cha:logic}, we use higher inductive types for colimits and quotients as in \autoref{sec:colimits,sec:set-quotients}, as well as some of the theory of truncation from \autoref{cha:hlevels}, particularly the factorization system of \autoref{sec:image-factorization} in the case $n=-1$.
In \autoref{sec:ordinals} we use an inductive family (\autoref{sec:generalizations}) to describe well-foundedness, and in \autoref{sec:cumulative-hierarchy} we use a more complicated higher inductive type to present the cumulative hierarchy.

\section{The category of sets}
\label{sec:piw-pretopos}

Recall that in \autoref{cha:category-theory} we defined the category \uset to consist of all $0$-types (in some universe \UU) and maps between them, and observed that it is a category (not just a precategory).
We consider successively the levels of structure which \uset possesses.

\subsection{Limits and colimits}
\label{subsec:limits-sets}

\index{limit!of sets}%
\index{colimit!of sets}%

Since sets are closed under products, the universal property of products in \autoref{thm:prod-ump} shows immediately that \uset has finite products.
In fact, infinite products follow just as easily from the equivalence
\[ \Parens{X\to \prd{a:A} B(a)} \eqvsym \Parens{\prd{a:A} (X\to B(a))}.\]
And we saw in \autoref{ex:pullback}\index{pullback} that the pullback of $f:A\to C$ and $g:B\to C$ can be defined as $\sm{a:A}{b:B} f(a)=g(b)$; this is a set if $A,B,C$ are and inherits the correct universal property.
Thus, \uset is a \emph{complete} category in the obvious sense.
\index{category!complete}%
\index{complete!category}%

Since sets are closed under $+$ and contain \emptyt, \uset has finite coproducts.
Similarly, since $\sm{a:A}B(a)$ is a set whenever $A$ and each $B(a)$ are, it yields a coproduct of the family $B$ in \uset.
Finally, we showed in \autoref{sec:pushouts} that pushouts exist in $n$-types, which includes \uset in particular.
Thus, \uset is also \emph{cocomplete}.
\index{category!cocomplete}%
\index{cocomplete category}%

\subsection{Images}
\label{sec:image}

Next, we show that $\uset$ is a \define{regular category}, i.e.:
\indexdef{category!regular}%
\indexdef{regular!category}%
\begin{enumerate}
\item $\uset$ is finitely complete.\label{item:reg1}
\item The kernel pair $\proj1,\proj2: (\sm{x,y:A} f(x)= f(y)) \to A$ of any
  function $f : A \to B$ has a coequalizer.\label{item:reg2}
  \indexdef{kernel!pair}
\item Pullbacks of regular epimorphisms are again regular epimorphisms.\label{item:reg3}
\end{enumerate}
Recall that a \define{regular epimorphism}
\indexdef{epimorphism!regular}%
\indexdef{regular!epimorphism}%
is a morphism that is the coequalizer of \emph{some} pair of maps.
Thus in~\ref{item:reg3} the pullback of a coequalizer is required to again be a coequalizer, but not necessarily of the pulled-back pair.

\index{set-coequalizer}%
\index{image}
The obvious candidate for the coequalizer of the kernel pair of $f:A\to B$ is the \emph{image} of $f$, as defined in \autoref{sec:image-factorization}.
Recall that we defined $\im(f)\defeq \sm{b:B} \brck{\hfib f b}$, with functions 
$\tilde{f}:A\to\im(f)$ and $i_f:\im(f)\to B$ defined by
\begin{align*}
  \tilde{f} & \defeq \lam{a} \Pairr{f(a),\,\bproj{\pairr{a,\refl{f(a)}}}}\\
i_f & \defeq \proj1
\end{align*}
fitting into a diagram:
\begin{equation*}
  \xymatrix{
    **[l]{\sm{x,y:A} f(x)= f(y)}
    \ar@<0.25em>[r]^{\proj1}
    \ar@<-0.25em>[r]_{\proj2}
    &
    {A}
    \ar[r]^(0.4){\tilde{f}}
    \ar[rd]_{f}
    &
    {\im(f)}
    \ar@{..>}[d]^{i_f}
    \\ & &
    B
  }
\end{equation*}

Recall that a function $f:A\to B$ is called \emph{surjective} if
\index{function!surjective}%
\narrowequation{\fall{b:B}\brck{\hfib f b},}
or equivalently $\fall{b:B} \exis{a:A} f(a)=b$.
We have also said that a function $f:A\to B$ between sets is called \emph{injective} if
\index{function!injective}%
$\fall{a,a':A} (f(a) = f(a')) \Rightarrow (a=a')$, or equivalently if each of its fibers is a mere proposition.
Since these are the $(-1)$-connected and $(-1)$-truncated maps in the sense of \autoref{cha:hlevels}, the general theory there implies that $\tilde f$ above is surjective and $i_f$ is injective, and that this factorization is stable under pullback.

We now identify surjectivity and injectivity with the appropriate cat\-e\-go\-ry-theoretic notions.
First we observe that categorical monomorphisms and epimorphisms have a slightly stronger equivalent formulation.

\begin{lem}\label{thm:mono}
  For a morphism $f:\hom_A(a,b)$ in a category $A$, the following are equivalent.
  \begin{enumerate}
  \item $f$ is a \define{monomorphism}:
    \indexdef{monomorphism}%
    for all $x:A$ and ${g,h:\hom_A(x,a)}$, if $f\circ g = f\circ h$ then $g=h$.\label{item:mono1}
  \item (If $A$ has pullbacks) the diagonal map $a\to a\times_b a$ is an isomorphism.\label{item:mono4}
  \item For all $x:A$ and $k:\hom_A(x,b)$, the type $\sm{h:\hom_A(x,a)} (k = f\circ h)$ is a mere proposition.\label{item:mono2}
  \item For all $x:A$ and ${g:\hom_A(x,a)}$, the type $\sm{h:\hom_A(x,a)} (f\circ g = f\circ h)$ is contractible.\label{item:mono3}
  \end{enumerate}
\end{lem}
\begin{proof}
  The equivalence of conditions~\ref{item:mono1} and~\ref{item:mono4} is standard category theory.
  Now consider the function $(f\circ \blank ):\hom_A(x,a) \to \hom_A(x,b)$ between sets.
  Condition~\ref{item:mono1} says that it is injective, while~\ref{item:mono2} says that its fibers are mere propositions; hence they are equivalent.
  And~\ref{item:mono2} implies~\ref{item:mono3} by taking $k\defeq f\circ g$ and recalling that an inhabited mere proposition is contractible.
  Finally,~\ref{item:mono3} implies~\ref{item:mono1} since if $p:f\circ g= f\circ h$, then $(g,\refl{})$ and $(h,p)$ both inhabit the type in~\ref{item:mono3}, hence are equal and so $g=h$.
\end{proof}

\begin{lem}
  A function $f:A\to B$ between sets is injective if and only if it is a monomorphism\index{monomorphism} in \uset.
\end{lem}
\begin{proof}
  Left to the reader.
\end{proof}

Of course, an \define{epimorphism}
\indexdef{epimorphism}
\indexsee{epi}{epimorphism}
is a monomorphism in the opposite category.
We now show that in \uset, the epimorphisms are precisely the surjections, and also precisely the coequalizers (regular epimorphisms).

The coequalizer of a pair of maps $f,g:A\to B$ in $\uset$ is defined as the 0-truncation of a general (homotopy) coequalizer.
For clarity, we may call this the \define{set-coequalizer}.
\indexdef{set-coequalizer}%
\indexsee{coequalizer!of sets}{set-coequalizer}%
It is convenient to express its universal property as follows.

\begin{lem}
\index{universal!property!of set-coequalizer}%
Let $f,g:A\to B$ be functions between sets $A$ and $B$. The 
{set-co}equalizer $c_{f,g}:B\to Q$ has the property that, for any set $C$ and any $h:B\to C$ with $h\circ f = h\circ g$, the type
\begin{equation*}
\sm{k:Q\to C} (k\circ c_{f,g} = h)
\end{equation*}
is contractible.
\end{lem}

\begin{lem}\label{epis-surj}
For any function $f:A\to B$ between sets, the following are equivalent:
\begin{enumerate}
\item $f$ is an epimorphism.
\item Consider the pushout diagram
\begin{equation*}
  \xymatrix{
    {A}
    \ar[r]^{f}
    \ar[d]
    &
    {B}
    \ar[d]^{\iota}
    \\
    {\unit}
    \ar[r]_{t}
    &
    {C_f}
  }
\end{equation*}
in $\uset$ defining the mapping cone\index{cone!of a function}. Then the type $C_f$ is contractible.
\item $f$ is surjective.
\end{enumerate}
\end{lem}

\begin{proof}
Let $f:A\to B$ be a function between sets, and suppose it to be an epimorphism; we show $C_f$ is contractible.
The constructor $\unit\to C_f$ of $C_f$ gives us an element $t:C_f$.
We have to show that
\begin{equation*}
\prd{x:C_f} x= t.
\end{equation*}
Note that $x= t$ is a mere proposition, hence we can use induction on $C_f$.
Of course when $x$ is $t$ we have $\refl{t}:t=t$, so it suffices to find
\begin{align*}
I_0 & : \prd{b:B} \iota(b)= t\\
I_1 & : \prd{a:A} \opp{\alpha_1(a)} \ct I_0(f(a))=\refl{t}.
\end{align*}
where $\iota:B\to C_f$ and $\alpha_1:\prd{a:A} \iota(f(a))= t$ are the other constructors
of $C_f$. Note that $\alpha_1$ is a homotopy from $\iota\circ f$ to
$\mathsf{const}_t\circ f$, so we find the elements
\begin{equation*}
\pairr{\iota,\refl{\iota\circ f}},\pairr{\mathsf{const}_t,\alpha_1}:
\sm{h:B\to C_f} \iota\circ f \htpy h\circ f.
\end{equation*}
By the dual of \autoref{thm:mono}\ref{item:mono3} (and function extensionality), there is a path
\begin{equation*}
\gamma:\pairr{\iota,\refl{\iota\circ f}}=\pairr{\mathsf{const}_t,\alpha_1}.
\end{equation*}
Hence, we may define $I_0(b)\defeq \happly(\projpath1(\gamma),b):\iota(b)=t$.
We also have
\[\projpath2(\gamma) : \trans{\projpath1(\gamma)}{\refl{\iota\circ f}} = \alpha_1. \]
This transport involves precomposition with $f$, which commutes with $\happly$.
Thus, from transport in path types we obtain $I_0(f(a)) = \alpha_1(a)$ for any $a:A$, which gives us $I_1$.

Now suppose $C_f$ is contractible; we show $f$ is surjective.
We first construct a type family $P:C_f\to\prop$ by recursion on $C_f$, which is valid since \prop is a set.
On the point constructors, we define
\begin{align*}
P(t) & \defeq \unit\\
P(\iota(b)) & \defeq \brck{\hfiber{f}b}.
\end{align*}
To complete the construction of $P$, it remains to give a path
\narrowequation{\brck{\hfiber{f}{f(a)}} =_\prop \unit}
for all $a:A$.
However, $\brck{\hfiber{f}{f(a)}}$ is inhabited by $(f(a),\refl{f(a)})$.
Since it is a mere proposition, this means it is contractible --- and thus equivalent, hence equal, to \unit.
This completes the definition of $P$.
Now, since $C_f$ is assumed to be contractible, it follows that $P(x)$ is equivalent to $P(t)$ for any $x:C_f$.
In particular, $P(\iota(b))\jdeq \brck{\hfiber{f}b}$ is equivalent to $P(t)\jdeq \unit$ for each $b:B$, and hence contractible.
Thus, $f$ is surjective.

Finally, suppose $f:A\to B$ to be surjective, and consider a set $C$ and two functions
$g,h:B\to C$ with the property that $g\circ f = h\circ f$. Since $f$ 
is assumed to be surjective, for all $b:B$ the type $\brck{\hfib f b}$ is contractible.
Thus we have the following equivalences:
\begin{align*}
\prd{b:B} (g(b)= h(b))
& \eqvsym \prd{b:B} \Parens{\brck{\hfib f b} \to (g(b)= h(b))}\\
& \eqvsym \prd{b:B} \Parens{\hfib f b \to (g(b)= h(b))}\\
& \eqvsym \prd{b:B}{a:A}{p:f(a)= b} g(b)= h(b)\\
& \eqvsym \prd{a:A} g(f(a))= h(f(a))
\end{align*}
using on the second line the fact that $g(b)=h(b)$ is a mere proposition, since $C$ is a set.
But by assumption, there is an element of the latter type.
\end{proof}


\begin{thm}\label{thm:set_regular}\label{lem:images_are_coequalizers}
The category $\uset$ is regular. Moreover, surjective functions between sets are regular epimorphisms.
\end{thm}

\begin{proof}
It is a standard lemma in category theory that a category is regular as soon as it admits finite limits and a pullback-stable orthogonal 
factorization system\index{orthogonal factorization system} $(\mathcal{E},\mathcal{M})$ with $\mathcal{M}$ the monomorphisms, in which case $\mathcal{E}$ consists automatically of 
the regular epimorphisms.
(See e.g.~\cite[A1.3.4]{elephant}.)
The existence of the factorization system was proved in \autoref{thm:orth-fact}.
\end{proof}

\begin{lem}\label{lem:pb_of_coeq_is_coeq}
Pullbacks of regular epis in \uset are regular epis.
\end{lem}
\begin{proof}
  We showed in \autoref{thm:stable-images} that pullbacks of $n$-connected functions are $n$-connected.
  By \autoref{lem:images_are_coequalizers}, it suffices to apply this when $n=-1$.
\end{proof}

\indexdef{image!of a subset}
One of the consequences of \uset being a regular category is that we have an ``image'' operation on subsets.
That is, given $f:A\to B$, any subset $P:\power A$ (i.e.\ a predicate $P:A\to \prop$) has an \define{image} which is a subset of $B$.
This can be defined directly as $\setof{ y:B | \exis{x:A} f(x)=y \land P(x)}$, or indirectly as the image (in the previous sense) of the composite function
\[ \setof{ x:A | P(x) } \to A \xrightarrow{f} B.\]
\symlabel{subset-image}
We will also sometimes use the common notation $\setof{f(x) | P(x)}$ for the image of $P$.

\subsection{Quotients}\label{subsec:quotients}

\index{set-quotient|(}%
Now that we know that $\uset$ is regular, to show that $\uset$ is exact, we need to show that every
equivalence relation is effective.
\index{effective!equivalence relation|(}%
\index{relation!effective equivalence|(}%
In other words, given an equivalence
relation $R:A\to A\to\prop$, there is a coequalizer $c_R$ of the pair
$\proj1,\proj2:\sm{x,y:A} R(x,y)\to A$ and, moreover, the $\proj1$ and $\proj2$
form the kernel\index{kernel!pair} pair of $c_R$.

We have already seen, in \autoref{sec:set-quotients}, two general ways to construct the quotient of a set by an equivalence relation $R:A\to A\to\prop$.
The first can be described as the set-coequalizer of the two projections
\[\proj1,\proj2:\Parens{\sm{x,y:A} R(x,y)} \to A.\]
The important property of such a quotient is the following.

\begin{defn}
  A relation $R:A\to A\to\prop$ is said to be \define{effective}
  \indexdef{effective!relation}
  \indexdef{effective!equivalence relation}%
  \indexdef{relation!effective equivalence}%
  if the square
\begin{equation*}
  \xymatrix{
    {\sm{x,y:A} R (x,y)}
    \ar[r]^(0.7){\proj1}
    \ar[d]_{\proj2}
    &
    {A}
    \ar[d]^{c_R}
    \\
    {A}
    \ar[r]_{c_R}
    &
    {A/R}
    }
\end{equation*}
is a pullback. 
\end{defn}

Since the standard pullback of $c_R$ and itself is $\sm{x,y:A} (c_R(x)=c_R(y))$, by \autoref{thm:total-fiber-equiv} this is equivalent to asking that the canonical transformation $\prd{x,y:A} R(x,y) \to (c_R(x)=c_R(y))$ be a fiberwise equivalence.

\begin{lem}\label{lem:sets_exact}
Suppose $\pairr{A,R}$ is an equivalence relation. Then there is an equivalence
\begin{equation*}
(c_R(x)= c_R(y))\eqvsym R(x,y)
\end{equation*}
for any $x,y:A$. In other words, equivalence relations are effective.
\end{lem}

\begin{proof}
We begin by extending $R$ to a relation $\widetilde{R}:A/R\to A/R\to\prop$, which we will then show is equivalent
to the identity type on $A/R$. We define $\widetilde{R}$ by double induction on
$A/R$ (note that $\prop$ is a set by univalence for mere propositions). We
define $\widetilde{R}(c_R(x),c_R(y)) \defeq R(x,y)$. For $r:R(x,x')$ and $s:R(y,y')$,
the transitivity and symmetry 
of $R$ gives an equivalence from $R(x,y)$ to $R(x',y')$. This completes the
definition of $\widetilde{R}$.

It remains to show that $\widetilde{R}(w,w')\eqvsym (w= w')$ for every $w,w':A/R$.
The direction $(w=w')\to \widetilde{R}(w,w')$ follows by transport once we show that $\widetilde{R}$ is reflexive, which is an easy induction.
The other direction $\widetilde{R}(w,w')\to (w= w')$ is a mere proposition, so since $c_R:A\to A/R$ is surjective, it suffices to assume that $w$ and $w'$ are of the form $c_R(x)$ and $c_R(y)$.
But in this case, we have the canonical map $\widetilde{R}(c_R(x),c_R(y)) \defeq R(x,y) \to (c_R(x)=c_R(y))$.
(Note again the appearance of the encode-decode method.\index{encode-decode method})
\end{proof}

The second construction of quotients is as the set of equivalence classes of $R$ (a subset
of its power set\index{power set}):
\[ A\sslash R \defeq \setof{ P:A\to\prop | P \text{ is an equivalence class of } R} \]
This requires propositional resizing\index{propositional resizing}\index{impredicative!quotient}\index{resizing} in order to remain in the same universe as $A$ and $R$.

Note that if we regard $R$ as a function from $A$ to $A\to \prop$, then $A\sslash R$ is equivalent to $\im(R)$, as constructed in \autoref{sec:image}.
Now in \autoref{lem:images_are_coequalizers} we have shown that images are
coequalizers. In particular, we immediately get the coequalizer diagram
\begin{equation*}
  \xymatrix{
    **[l]{\sm{x,y:A} R (x)= R (y)}
    \ar@<0.25em>[r]^{\proj1}
    \ar@<-0.25em>[r]_{\proj2}
    &
    {A}
    \ar[r]
    &
    {A \sslash R.}
  }
\end{equation*}
We can use this to give an alternative proof that any equivalence relation is effective and that the two definitions of quotients agree.

\begin{thm}\label{prop:kernels_are_effective}
For any function $f:A\to B$ between any two sets, 
the relation $\ker(f):A\to A\to\prop$ given by 
$\ker(f,x,y)\defeq (f(x)= f(y))$ is effective.
\end{thm}

\begin{proof}
We will use that $\im(f)$ is the coequalizer of $\proj1,\proj2:
(\sm{x,y:A} f(x)= f(y))\to A$.
Note that the kernel pair of the function 
\[c_f\defeq\lam{a} \Parens{f(a),\brck{\pairr{a,\refl{f(a)}}}}
: A \to \im(f)
\]
consists of the two projections
\begin{equation*}
\proj1,\proj2:\Parens{\sm{x,y:A} c_f(x)= c_f(y)}\to A.
\end{equation*}
For any $x,y:A$, we have equivalences
\begin{align*}
  (c_f(x)= c_f(y))
  & \eqvsym \Parens{\sm{p:f(x)= f(y)} \trans{p}{\brck{\pairr{x,\refl{f(x)}}}} =\brck{\pairr{y,\refl{f(x)}}}}\\
  & \eqvsym (f(x)= f(y)),
\end{align*}
where the last equivalence holds because 
$\brck{\hfiber{f}b}$ is a mere proposition for any $b:B$. 
Therefore, we get that
\begin{equation*}
\Parens{\sm{x,y:A} c_f(x)= c_f(y)}\eqvsym \Parens{\sm{x,y:A} f(x)= f(y)}
\end{equation*}
and hence we may conclude that $\ker f$ is an effective relation 
for any function $f$.
\end{proof}

\begin{thm}
Equivalence relations are effective and there is an equivalence $A/R \eqvsym A\sslash  R $. 
\end{thm}

\begin{proof}
We need to analyze the coequalizer diagram
\begin{equation*}
  \xymatrix{
    **[l]{\sm{x,y:A} R (x)= R (y)}
    \ar@<0.25em>[r]^{\proj1}
    \ar@<-0.25em>[r]_{\proj2}
    &
    {A}
    \ar[r]
    &
    {A \sslash R}
  }
\end{equation*}
By the univalence axiom, the type $R(x) = R(y)$ is equivalent to the type of homotopies from $R(x)$ to $R(y)$, which is
equivalent to
\narrowequation{\prd{z:A} R (x,z)\eqvsym R (y,z).}
Since $R$ is an equivalence relation, the latter space is equivalent to $R(x,y)$. To
summarize, we get that $(R(x) = R(y)) \eqvsym R(x,y)$, so $R $ is effective since it is equivalent to an effective relation. Also,
the diagram
\begin{equation*}
  \xymatrix{
    **[l]{\sm{x,y:A} R(x, y)}
    \ar@<0.25em>[r]^{\proj1}
    \ar@<-0.25em>[r]_{\proj2}
    &
    {A}
    \ar[r]
    &
    {A \sslash R.}
  }
\end{equation*}
is a coequalizer diagram. Since coequalizers are unique up to equivalence, it follows that $A/R \eqvsym A\sslash  R $.
\end{proof}

We finish this section by mentioning a possible third construction of the quotient of a set $A$ by an equivalence relation $R$.
Consider the precategory with objects $A$ and hom-sets $R$; the type of objects of the Rezk completion
\index{completion!Rezk}%
(see \autoref{sec:rezk}) of this precategory will then be the
quotient. The reader is invited to check the details.

\index{effective!equivalence relation|)}%
\index{relation!effective equivalence|)}%
\index{set-quotient|)}%

\subsection{\texorpdfstring{$\uset$}{Set} is a \texorpdfstring{$\Pi\mathsf{W}$}{ΠW}-pretopos}
\label{subsec:piw}

\index{structural!set theory|(}%

The notion of a \emph{$\Pi\mathsf{W}$-pretopos}
\index{PiW-pretopos@$\Pi\mathsf{W}$-pretopos}%
\indexsee{pretopos}{$\Pi\mathsf{W}$-pretopos}
--- that is, a locally cartesian closed category
\index{locally cartesian closed category}%
\index{category!locally cartesian closed}%
with disjoint finite coproducts, effective equivalence relations, and initial algebras for polynomial endofunctors --- is intended as a ``predicative''
\index{mathematics!predicative}%
notion of topos, i.e.\ a category of ``predicative sets'', which can serve the purpose for constructive mathematics
\index{mathematics!constructive}%
that the usual category of sets does for classical
\index{mathematics!classical}%
mathematics.

Typically, in constructive type theory, one resorts to an external construction of ``setoids'' --- an exact completion --- to obtain a category with such closure properties.
\index{setoid}\index{completion!exact}%
  In particular, the well-behaved quotients are required for many constructions in mathematics that usually involve (non-constructive) power sets.  It is noteworthy that univalent foundations provides these constructions \emph{internally} (via higher inductive types), without requiring such external constructions.  This represents a powerful advantage of our approach, as we shall see in subsequent examples.

\begin{thm}
  \index{PiW-pretopos@$\Pi\mathsf{W}$-pretopos}
  The category $\uset$ is a $\Pi\mathsf{W}$-pretopos.
\end{thm}
\begin{proof}
  We have an initial object
  \index{initial!set}%
  $\emptyt$ and finite, disjoint sums $A+B$.  These are stable under pullback, simply because pullback has a right adjoint\index{adjoint!functor}.  Indeed, $\uset$ is locally cartesian closed, since for any map $f:A\to B$ between sets, the ``fibrant replacement'' \index{fibrant replacement} $\sm{a:A}f(a)=b$ is equivalent to $A$ (over $B$), and we have dependent function types for the replacement.
We've just shown that $\uset$ is regular (\autoref{thm:set_regular}) and that quotients are effective (\autoref{lem:sets_exact}). We thus have a locally cartesian closed pretopos. Finally, since the $n$-types are closed under the formation of $W$-types by \autoref{ex:ntypes-closed-under-wtypes}, and by \autoref{thm:w-hinit} $W$-types are initial algebras for polynomial endofunctors, we see that $\uset$ is a $\Pi\mathsf{W}$-pretopos.
\end{proof}

\index{topos|(}
One naturally wonders what, if anything, prevents $\uset$ from being an (elementary) topos?
In addition to the structure already mentioned, a topos has a
\emph{subobject classifier}:
\indexdef{subobject classifier}%
\index{classifier!subobject}%
\index{power set}%
a pointed object classifying (equivalence classes of) monomorphisms\index{monomorphism}.  (In fact, in the presence of a subobject
classifier, things become somewhat simpler: one merely needs cartesian closure in order to get the colimits.)
In homotopy type theory, univalence implies that the type $\prop \defeq \sm{X:\UU}\isprop(X)$ does classify monomorphisms (by an argument similar to \autoref{sec:object-classification}), but in general it is as large as the ambient universe $\UU$.
Thus, it is a ``set'' in the sense of being a $0$-type, but it is not ``small'' in the sense of being an object of $\UU$, hence not an object of the category \uset.
However, if we assume an appropriate form of propositional resizing (see \autoref{subsec:prop-subsets}), then we can find a small version of $\prop$, so that \uset becomes an elementary topos.

\begin{thm}\label{thm:settopos}
  \index{propositional!resizing}%
  If there is a type $\Omega:\UU$ of all mere propositions, then the category $\uset_\UU$ is an elementary topos.
\end{thm}
\index{topos|)}

A sufficient condition for this is the law of excluded middle, in the ``mere-propositional'' form that we have called \LEM{}; for then we have $\prop = \bool$, which \emph{is} small, and which then also classifies all mere propositions.
Moreover, in topos theory a well-known sufficient condition for \LEM{} is the axiom of choice, which is of course often assumed as an axiom in classical\index{mathematics!classical} set theory.
In the next section, we briefly investigate the relation between these conditions in our setting.

\index{structural!set theory|)}%

\subsection{The axiom of choice implies excluded middle}
\label{subsec:emacinsets}


We begin with the following lemma.

\begin{lem}\label{prop:trunc_of_prop_is_set}
If $A$ is a mere proposition then its suspension $\susp(A)$ is a set,
and $A$ is equivalent to $\id[\susp(A)]{\north}{\south}$.
\end{lem}

\begin{proof}
To show that $\susp(A)$ is a set, we define a
family $P:\susp(A)\to\susp(A)\to\type$ with the 
property that $P(x,y)$ is a mere proposition for each $x,y:\susp(A)$,
and which is equivalent to its identity type $\idtypevar{\susp(A)}$.
We make the following definitions:
\begin{align*}
P(\north,\north) & \defeq \unit &
P(\south,\north) & \defeq A\\
P(\north,\south) & \defeq A &
P(\south,\south) & \defeq \unit.
\end{align*}
We have to check that the definition preserves paths.
Given any $a : A$, there is a meridian $\merid(a) : \north = \south$,
so we should also have
\begin{equation*}
  P(\north, \north) = P(\north, \south) = P(\south, \north) = P(\south, \south).
\end{equation*}
But since $A$ is inhabited by $a$, it is equivalent to $\unit$, so we have
\begin{equation*}
  P(\north, \north) \eqvsym P(\north, \south) \eqvsym P(\south, \north) \eqvsym P(\south, \south).
\end{equation*}
The univalence axiom turns these into the desired equalities. Also, $P(x,y)$ is a mere
proposition for all $x, y : \susp(A)$, which is proved by induction on $x$ and $y$, and
using the fact that being a mere proposition is a mere proposition.

Note that $P$ is a reflexive relation.
Therefore we may apply \autoref{thm:h-set-refrel-in-paths-sets}, so it suffices to
construct $\tau : \prd{x,y:\susp(A)}P(x,y)\to(x=y)$. We do this by a double induction.
When $x$ is $\north$, we define $\tau(\north)$ by
\begin{equation*}
  \tau(\north,\north,u) \defeq \refl{\north}
  \qquad\text{and}\qquad
  \tau(\north,\south,a) \defeq \merid(a).
\end{equation*}
If $A$ is inhabited by $a$ then $\merid(a) : \north = \south$ so we also need
\narrowequation{
  \trans{\merid(a)}{\tau(\north, \north)} = \tau(\north, \south).
}
This we get by function extensionality using the fact that, for all $x : A$,
\begin{multline*}
  \trans{\merid(a)}{\tau(\north,\north,x)} =
  \tau(\north,\north,x) \ct \opp{\merid(a)} \jdeq \\
  \refl{\north} \ct \merid(a) =
  \merid(a) =
  \merid(x) \jdeq
  \tau(\north, \south, x).
\end{multline*}
In a symmetric fashion we may define $\tau(\south)$ by
\begin{equation*}
  \tau(\south,\north, a) \defeq \opp{\merid(a)}
  \qquad\text{and}\qquad
  \tau(\south,\south, u) \defeq \refl{\south}.
\end{equation*}
To complete the construction of $\tau$, we need to check $\trans{\merid(a)}{\tau(\north)} = \tau(\south)$,
given any $a : A$. The verification proceeds much along the same lines by induction on the
second argument of $\tau$.

Thus, by \autoref{thm:h-set-refrel-in-paths-sets} we have that $\susp(A)$ is a set and that $P(x,y) \eqvsym (\id{x}{y})$ for all $x,y:\susp(A)$.
Taking $x\defeq \north$ and $y\defeq \south$ yields $\eqv{A}{(\id[\susp(A)]\north\south)}$ as desired.
\end{proof}

\begin{thm}[Diaconescu]\label{thm:1surj_to_surj_to_pem}
  \index{axiom!of choice}%
  \index{excluded middle}%
  \index{Diaconescu's theorem}\index{theorem!Diaconescu's}%
  The axiom of choice implies the law of excluded middle.
\end{thm}

\begin{proof}
We use the equivalent form of choice given in \autoref{thm:ac-epis-split}.
Consider a mere proposition $A$.
The function $f:\bool\to\susp(A)$ defined by
$f(\bfalse) \defeq \north$ and $f(\btrue) \defeq \south$
is surjective.
Indeed, we have
$\pairr{\bfalse,\refl{\north}} : \hfiber{f}{\north}$
and $\pairr{\btrue,\refl{\south}} :\hfiber{f}{\south}$.
Since $\bbrck{\hfiber{f}{x}}$ is a mere proposition, by induction the claimed surjectivity follows.

By \autoref{prop:trunc_of_prop_is_set} the suspension $\susp(A)$
is a set, so by the axiom of choice there merely exists a
section $g: \susp(A) \to \bool$ of $f$.
As equality on $\bool$ is decidable we get
\begin{equation*}
 (g(f(\bfalse))= g(f(\btrue))) +
 \lnot (g(f(\bfalse))= g(f(\btrue))),
\end{equation*}
and, since $g$ is a section of $f$, hence injective,
\begin{equation*}
(f(\bfalse) = f(\btrue)) +
\lnot (f(\bfalse) = f(\btrue)).
\end{equation*}
Finally, since $(f(\bfalse)=f(\btrue)) = (\north=\south) = A$ by \autoref{prop:trunc_of_prop_is_set}, we have $A+\neg A$.
\end{proof}




\index{denial}
\begin{thm}\label{thm:ETCS}
  \index{Elementary Theory of the Category of Sets}%
  \index{category!well-pointed}%
  If the axiom of choice holds then the category $\uset$ is a well-pointed boolean\index{topos!boolean}\index{boolean!topos} elementary topos\index{topos} with choice.  
\end{thm}

\begin{proof}
  Since \choice{} implies \LEM{}, we have a boolean elementary topos with choice by \autoref{thm:settopos} and the remark following it.  We leave the proof of well-pointedness as
an exercise for the reader (\autoref{ex:well-pointed}).
\end{proof}

\begin{rmk}
  The conditions on a category mentioned in the theorem are known as Lawvere's\index{Lawvere}
  axioms for the \emph{Elementary Theory of the Category of Sets}~\cite{lawvere:etcs-long}.
\end{rmk}

\section{Cardinal numbers}
\label{sec:cardinals}

\begin{defn}\label{defn:card}
  The \define{type of cardinal numbers}
  \indexdef{type!of cardinal numbers}%
  \indexdef{cardinal number}%
  \indexsee{number!cardinal}{cardinal number}%
  is the 0-truncation of the type \set of sets:
  \[ \card \defeq \pizero{\set} \]
  Thus, a \define{cardinal number}, or \define{cardinal}, is an inhabitant of $\card\jdeq \pizero\set$.
  Technically, of course, there is a separate type $\card_\UU$ associated to each universe \type.
\end{defn}



As usual for truncations, if $A$ is a set, then $\cd{A}$ denotes its image under the canonical projection $\set \to \trunc0\set \jdeq \card$; we call $\cd{A}$ the \define{cardinality}\indexdef{cardinality} of $A$.
By definition, \card is a set.
It also inherits the structure of a semiring from \set.

\begin{defn}
  The operation of \define{cardinal addition}
  \indexdef{addition!of cardinal numbers}%
  \index{cardinal number!addition of}%
  \[ (\blank+\blank) : \card \to \card \to \card \]
  is defined by induction on truncation:
  \[ \cd{A} + \cd{B} \defeq \cd{A+B}.\]
\end{defn}
\begin{proof}
  Since $\card\to\card$ is a set, to define $(\alpha+\blank):\card\to\card$ for all $\alpha:\card$, by induction it suffices to assume that $\alpha$ is $\cd{A}$ for some $A:\set$.
  Now we want to define $(\cd{A}+\blank) :\card\to\card$, i.e.\ we want to define $\cd{A}+\beta :\card$ for all $\beta:\card$.
  However, since $\card$ is a set, by induction it suffices to assume that $\beta$ is $\cd{B}$ for some $B:\set$.
  But now we can define $\cd{A}+\cd{B}$ to be $\cd{A+B}$.
\end{proof}

\begin{defn}
  Similarly, the operation of \define{cardinal multiplication}
  \indexdef{multiplication!of cardinal numbers}%
  \index{cardinal number!multiplication of}%
  \[ (\blank\cdot\blank) : \card \to \card \to \card \]
  is defined by induction on truncation:
  \[ \cd{A} \cdot \cd{B} \defeq \cd{A\times B} \]
\end{defn}

\begin{lem}\label{card:semiring}
  \card is a commutative semiring\index{semiring}, i.e.\ for $\alpha,\beta,\gamma:\card$ we have the following.
  \begin{align*}
    (\alpha+\beta)+\gamma &= \alpha+(\beta+\gamma)\\
    \alpha+0 &= \alpha\\
    \alpha + \beta &= \beta + \alpha\\
    (\alpha \cdot \beta) \cdot \gamma &= \alpha \cdot (\beta\cdot\gamma)\\
    \alpha \cdot 1 &= \alpha\\
    \alpha\cdot\beta &= \beta\cdot\alpha\\
    \alpha\cdot(\beta+\gamma) &= \alpha\cdot\beta + \alpha\cdot\gamma
  \end{align*}
  where $0 \defeq \cd{\emptyt}$ and $1\defeq\cd{\unit}$.
\end{lem}
\begin{proof}
  We prove the commutativity of multiplication, $\alpha\cdot\beta = \beta\cdot\alpha$; the others are exactly analogous.
  Since \card is a set, the type $\alpha\cdot\beta = \beta\cdot\alpha$ is a mere proposition, and in particular a set.
  Thus, by induction it suffices to assume $\alpha$ and $\beta$ are of the form $\cd{A}$ and $\cd{B}$ respectively, for some $A,B:\set$.
  Now $\cd{A}\cdot \cd{B} \jdeq \cd{A\times B}$ and $\cd{B}\times\cd{A} \jdeq \cd{B\times A}$, so it suffices to show $A\times B = B\times A$.
  Finally, by univalence, it suffices to give an equivalence $A\times B \eqvsym B\times A$.
  But this is easy: take $(a,b) \mapsto (b,a)$ and its obvious inverse.
\end{proof}

\begin{defn}
  The operation of \define{cardinal exponentiation} is also defined by induction on truncation:
  \indexdef{exponentiation, of cardinal numbers}%
  \index{cardinal number!exponentiation of}%
  \[ \cd{A}^{\cd{B}} \defeq \cd{B\to A}. \]
\end{defn}

\begin{lem}\label{card:exp}
  For $\alpha,\beta,\gamma:\card$ we have
  \begin{align*}
    \alpha^0 &= 1\\
    1^\alpha &= 1\\
    \alpha^1 &= \alpha\\
    \alpha^{\beta+\gamma} &= \alpha^\beta \cdot \alpha^\gamma\\
    \alpha^{\beta\cdot \gamma} &= (\alpha^{\beta})^\gamma\\
    (\alpha\cdot\beta)^\gamma &= \alpha^\gamma \cdot \beta^\gamma
  \end{align*}
\end{lem}
\begin{proof}
  Exactly like \autoref{card:semiring}.
\end{proof}

\begin{defn}
  The relation of \define{cardinal inequality}
  \index{order!non-strict}%
  \index{cardinal number!inequality of}%
  \[ (\blank\le\blank) : \card\to\card\to\prop \]
  is defined by induction on truncation:
  \symlabel{inj}
  \[ \cd{A} \le \cd{B} \defeq \brck{\inj(A,B)} \]
  where $\inj(A,B)$ is the type of injections from $A$ to $B$.
  \index{function!injective}%
  In other words, $\cd{A} \le \cd{B}$ means that there merely exists an injection from $A$ to $B$.
\end{defn}

\begin{lem}
  Cardinal inequality is a preorder, i.e.\ for $\alpha,\beta:\card$ we have
  \index{preorder!of cardinal numbers}%
  \begin{gather*}
    \alpha \le \alpha\\
    (\alpha \le \beta) \to (\beta\le\gamma) \to (\alpha\le\gamma)
  \end{gather*}
\end{lem}
\begin{proof}
  As before, by induction on truncation.
  For instance, since $(\alpha \le \beta) \to (\beta\le\gamma) \to (\alpha\le\gamma)$ is a mere proposition, by induction on 0-truncation we may assume $\alpha$, $\beta$, and $\gamma$ are $\cd{A}$, $\cd{B}$, and $\cd{C}$ respectively.
  Now since $\cd{A} \le \cd{C}$ is a mere proposition, by induction on $(-1)$-truncation we may assume given injections $f:A\to B$ and $g:B\to C$.
  But then $g\circ f$ is an injection from $A$ to $C$, so $\cd{A} \le \cd{C}$ holds.
  Reflexivity is even easier.
\end{proof}

We may likewise show that cardinal inequality is compatible with the semiring operations.

\begin{lem}\label{thm:injsurj}
  \index{function!injective}%
  \index{function!surjective}%
  Consider the following statements:
  \begin{enumerate}
  \item There is an injection $A\to B$.\label{item:cle-inj}
  \item There is a surjection $B\to A$.\label{item:cle-surj}
  \end{enumerate}
  Then, assuming excluded middle:
  \index{excluded middle}%
  \index{axiom!of choice}%
  \begin{itemize}
  \item Given $a_0:A$, we have~\ref{item:cle-inj}$\to$\ref{item:cle-surj}.
  \item Therefore, if $A$ is merely inhabited, we have~\ref{item:cle-inj} $\to$ merely \ref{item:cle-surj}.
  \item Assuming the axiom of choice, we have~\ref{item:cle-surj} $\to$ merely \ref{item:cle-inj}.
  \end{itemize}
\end{lem}
\begin{proof}
  If $f:A\to B$ is an injection, define $g:B\to A$ at $b:B$ as follows.
  Since $f$ is injective, the fiber of $f$ at $b$ is a mere proposition.
  Therefore, by excluded middle, either there is an $a:A$ with $f(a)=b$, or not.
  In the first case, define $g(b)\defeq a$; otherwise set $g(b)\defeq a_0$.
  Then for any $a:A$, we have $a = g(f(a))$, so $g$ is surjective.

  The second statement follows from this by induction on truncation.
  For the third, if $g:B\to A$ is surjective, then by the axiom of choice, there merely exists a function $f:A\to B$ with $g(f(a)) = a$ for all $a$.
  But then $f$ must be injective.
\end{proof}

\begin{thm}[Schroeder--Bernstein]
  \index{theorem!Schroeder--Bernstein}%
  \index{Schroeder--Bernstein theorem}%
  Assuming excluded middle, for sets $A$ and $B$ we have
  \[ \inj(A,B) \to \inj(B,A) \to (A\cong B) \]
\end{thm}
\begin{proof}
  The usual ``back-and-forth'' argument applies without significant changes.
  Note that it actually constructs an isomorphism $A\cong B$ (assuming excluded middle so that we can decide whether a given element belongs to a cycle, an infinite chain, a chain beginning in $A$, or a chain beginning in $B$).
\end{proof}

\begin{cor}
  Assuming excluded middle, cardinal inequality is a partial order, i.e.\ for $\alpha,\beta:\card$ we have
  \[ (\alpha\le\beta) \to (\beta\le\alpha) \to (\alpha=\beta). \]
\end{cor}
\begin{proof}
  Since $\alpha=\beta$ is a mere proposition, by induction on truncation we may assume $\alpha$ and $\beta$ are $\cd{A}$ and $\cd{B}$, respectively, and that we have injections $f:A\to B$ and $g:B\to A$.
  But then the Schroeder--Bernstein theorem gives an isomorphism $A\eqvsym B$, hence an equality $\cd{A}=\cd{B}$.
\end{proof}

Finally, we can reproduce Cantor's theorem, showing that for every cardinal there is a greater one.

\begin{thm}[Cantor]
  \index{Cantor's theorem}%
  \index{theorem!Cantor's}%
  For $A:\set$, there is no surjection $A \to (A\to \bool)$.
\end{thm}
\begin{proof}
  Suppose $f:A \to (A\to \bool)$ is any function, and define $g:A\to \bool$ by $g(a) \defeq \neg f(a)(a)$.
  If $g = f(a_0)$, then $g(a_0) = f(a_0)(a_0)$ but $g(a_0) = \neg f(a_0)(a_0)$, a contradiction.
  Thus, $f$ is not surjective.
\end{proof}

\begin{cor}
  Assuming excluded middle, for any $\alpha:\card$, there is a cardinal $\beta$ such that $\alpha\le\beta$ and $\alpha\neq\beta$.
\end{cor}
\begin{proof}
  Let $\beta = 2^\alpha$.
  Now we want to show a mere proposition, so by induction we may assume $\alpha$ is $\cd{A}$, so that $\beta\jdeq \cd{A\to \bool}$.
  Using excluded middle, we have a function $f:A\to (A\to \bool)$ defined by
  \[f(a)(a') \defeq
  \begin{cases}
    \btrue &\quad a=a'\\
    \bfalse &\quad a\neq a'.
  \end{cases}
  \]
  And if $f(a)=f(a')$, then $f(a')(a) = f(a)(a) = \btrue$, so $a=a'$; hence $f$ is injective.
  Thus, $\alpha \jdeq \cd{A} \le \cd{A\to \bool} \jdeq 2^\alpha$.

  On the other hand, if $2^\alpha \le \alpha$, then we would have an injection $(A\to\bool)\to A$.
  By \autoref{thm:injsurj}, since we have $(\lam{x} \bfalse):A\to \bool$ and excluded middle, there would then be a surjection $A \to (A\to \bool)$, contradicting Cantor's theorem.
\end{proof}

\section{Ordinal numbers}
\label{sec:ordinals}

\index{ordinal|(}%

\begin{defn}\label{defn:accessibility}
  Let $A$ be a set and
  \[(\blank<\blank):A\to A\to \prop\]
  a binary relation on $A$.
  We define by induction what it means for an element $a:A$ to be \define{accessible}
  \indexdef{accessibility}%
  \indexsee{accessible}{accessibility}%
  by $<$:
  \begin{itemize}
  \item If $b$ is accessible for every $b<a$, then $a$ is accessible.
  \end{itemize}
  We write $\acc(a)$ to mean that $a$ is accessible.
\end{defn}

It may seem that such an inductive definition can never get off the ground, but of course if $a$ has the property that there are \emph{no} $b$ such that $b<a$, then $a$ is vacuously accessible.

Note that this is an inductive definition of a family of types, like the type of vectors considered in \autoref{sec:generalizations}.
More precisely, it has one constructor, say $\acc_<$, with type
\[ \acc_< : \prd{a:A} \Parens{\prd{b:A} (b<a) \to \acc(b)} \to \acc(a). \]
\index{induction principle!for accessibility}%
The induction principle for $\acc$ says that for any $P:\prd{a:A} \acc(a) \to \type$, if we have
\[f:\prd{a:A}{h:\prd{b:A} (b<a) \to \acc(b)}
\Parens{\prd{b:A}{l:b<a} P(b,h(b,l))} \to
P(a,\acc_<(a,h)),
\]
then we have $g:\prd{a:A}{c:\acc(a)} P(a,c)$ defined by induction, with
\[g(a,\acc_<(a,h)) \jdeq f(a,\,h,\,\lam{b}{l} g(b,h(b,l))).\]
This is a mouthful, but generally we apply it only in the simpler case where $P:A\to\type$ depends only on $A$.
In this case the second and third arguments of $f$ may be combined, so that what we have to prove is
\[f:\prd{a:A} \Parens{\prd{b:A} (b<a) \to \acc(b) \times P(b)}
\to P(a).
\]
That is, we assume every $b<a$ is accessible and $g(b):P(b)$ is defined, and from these define $g(a):P(a)$.

The omission of the second argument of $P$ is justified by the following lemma, whose proof is the only place where we use the more general form of the induction principle.

\begin{lem}
  Accessibility\index{accessibility} is a mere property.
\end{lem}
\begin{proof}
  We must show that for any $a:A$ and $s_1,s_2:\acc(a)$ we have $s_1=s_2$.
  We prove this by induction on $s_1$, with
  \[P_1(a,s_1) \defeq \prd{s_2:\acc(a)} (s_1=s_2). \]
  Thus, we must show that for any $a:A$ and ${h_1:\prd{b:A} (b<a) \to \acc(b)}$ and
  \[ k_1:{\prd{b:A}{l:b<a}{t:\acc(b)} h_1(b,l) = t},\]
  we have $\acc_<(a,h) = s_2$ for any $s_2:\acc(a)$.
  We regard this statement as $\prd{a:A}{s_2:\acc(a)} P_2(a,s_2)$, where
  \[P_2(a,s_2) \defeq
  \prd{h_1 : \cdots } 
  {k_1 : \cdots} 
  (\acc_<(a,h_1) = s_2);
  \]
  thus we may prove it by induction on $s_2$.
  Therefore, we assume $h_2 : \prd{b:A} (b<a) \to \acc(b)$, and $k_2$ with a monstrous but irrelevant type,
  and must show that for any $h_1$ and $k_1$ with types as above,
  we have $\acc_<(a,h_1) = \acc_<(a,h_2)$.
  By function extensionality, it suffices to show $h_1(b,l) = h_2(b,l)$ for all $b:A$ and $l:b<a$.
  This follows from $k_1$.
\end{proof}

\begin{defn}
  A binary relation $<$ on a set $A$ is \define{well-founded}
  \indexdef{relation!well-founded}%
  \indexdef{well-founded!relation}%
  if every element of $A$ is accessible.
\end{defn}

The point of well-foundedness is that for $P:A\to \type$, we can use the induction principle of $\acc$ to conclude $\prd{a:A} \acc(a) \to P(a)$, and then apply well-foundedness to conclude $\prd{a:A} P(a)$.
In other words, if from $\fall{b:A} (b<a) \to P(b)$ we can prove $P(a)$, then $\fall{a:A} P(a)$.
This is called \define{well-founded induction}\indexdef{well-founded!induction}.

\begin{lem}
  Well-foundedness is a mere property.
\end{lem}
\begin{proof}
  Well-foundedness of $<$ is the type $\prd{a:A} \acc(a)$, which is a mere proposition since each $\acc(a)$ is.
\end{proof}

\begin{eg}\label{thm:nat-wf}
  Perhaps the most familiar well-founded relation is the usual strict ordering on \nat.
  To show that this is well-founded, we must show that $n$ is accessible for each $n:\nat$.
  \index{strong!induction}%
  This is just the usual proof of ``strong induction'' from ordinary induction on \nat.

  Specifically, we prove by induction on $n:\nat$ that $k$ is accessible for all $k\le n$.
  The base case is just that $0$ is accessible, which is vacuously true since nothing is strictly less than $0$.
  For the inductive step, we assume that $k$ is accessible for all $k\le n$, which is to say for all $k<n+1$; hence by definition $n+1$ is also accessible.

  A different relation on \nat which is also well-founded is obtained by setting only $n < \suc(n)$ for all $n:\nat$.
  Well-foundedness of this relation is almost exactly the ordinary induction principle of \nat.
\end{eg}

\begin{eg}\label{thm:wtype-wf}
  Let $A:\set$ and $B : A \to \set$ be a family of sets.
  Recall from \autoref{sec:w-types} that the $W$-type $\wtype{a:A} B(a)$ is inductively generated by the single constructor
  \begin{itemize}
  \item $\supp : \prd{a:A} (B(a) \to \wtype{x:A} B(x)) \to \wtype{x:A} B(x)$
  \end{itemize}
  We define the relation $<$ on $\wtype{x:A} B(x)$ by recursion on its second argument:
  \begin{itemize}
  \item For any $a:A$ and $f:B(a) \to \wtype{x:A} B(x)$, we define $w<\supp(a,f)$ to mean that there merely exists a $b:B(a)$ such that $w = f(b)$.
  \end{itemize}
  Now we prove that every $w:\wtype{x:A} B(x)$ is accessible for this relation, using the usual induction principle for $\wtype{x:A}B(x)$.
  This means we assume given $a:A$ and $f:B(a) \to \wtype{x:A} B(x)$, and also a lifting $f' : \prd{b:B(a)} \acc(f(b))$.
  But then by definition of $<$, we have $\acc(w)$ for all $w<\supp(a,f)$; hence $\supp(a,f)$ is accessible.
\end{eg}

Well-foundedness allows us to define functions by recursion and prove statements by induction, such as for instance the following.
Recall from \autoref{subsec:prop-subsets} that $\power B$ denotes the \emph{power set}\index{power set} $\power B \defeq (B\to\prop)$.

\begin{lem}\label{thm:wfrec}
  Suppose $B$ is a set and we have a function
  \[ g : \power B \to B \]
  Then if $<$ is a well-founded relation on $A$, there is a function $f:A\to B$ such that for all $a:A$ we have
  \begin{equation*}
    f(a) = g\Big(\setof{ f(a') | a'<a }\Big).
  \end{equation*}
\end{lem}
\noindent
(We are using the notation for images of subsets from \autoref{sec:image}.)
\begin{proof}
  We first define, for every $a:A$ and $s:\acc(a)$, an element $\bar f(a,s):B$.
  By induction, it suffices to assume that $s$ is a function assigning to each $a'<a$ a witness $s(a'):\acc(a')$, and that moreover for each such $a'$ we have an element $\bar f(a',s(a')):B$.
  In this case, we define
  \begin{equation*}
    \bar f(a,s) \defeq g\Big(\setof{ \bar f(a',s(a')) | a'<a }\Big).
  \end{equation*}

  Now since $<$ is well-founded, we have a function $w:\prd{a:A} \acc(a)$.
  Thus, we can define $f(a)\defeq \bar f (a,w(a))$.
\end{proof}

In classical\index{mathematics!classical} logic, well-foundedness has a more well-known reformulation.
In the following, we say that a subset $B: \power A$ is \define{nonempty}
\indexdef{nonempty subset}
if it is unequal to the empty subset $(\lam{x}\bot) : \power X$.
We leave it to the reader to verify that assuming excluded middle, this is equivalent to mere inhabitation, i.e.\ to the condition $\exis{x:A} x\in B$.

\begin{lem}\label{thm:wfmin}
  \index{excluded middle}%
  Assuming excluded middle, $<$ is well-founded if and only if every nonempty subset $B: \power A$ merely has a minimal element.
\end{lem}
\begin{proof}
  Suppose first $<$ is well-founded, and suppose $B\subseteq A$ is a subset with no minimal element.
  That is, for any $a:A$ with $a\in B$, there merely exists a $b:A$ with $b<a$ and $b\in B$.

  We claim that for any $a:A$ and $s:\acc(a)$, we have $a\notin B$.
  By induction, we may assume $s$ is a function assigning to each $a'<a$ a proof $s(a'):\acc(a)$, and that moreover for each such $a'$ we have $a'\notin B$.
  If $a\in B$, then by assumption, there would merely exist a $b<a$ with $b\in B$, which contradicts this assumption.
  Thus, $a\notin B$; this completes the induction.
  Since $<$ is well-founded, we have $a\notin B$ for all $a:A$, i.e. $B$ is empty.

  Now suppose each nonempty subset merely has a minimal element.
  Let $B = \setof{ a:A | \neg \acc(a) }$.
  Then if $B$ is nonempty, it merely has a minimal element.
  Thus there merely exists an $a:A$ with $a\in B$ such that for all $b<a$, we have $\acc(b)$.
  But then by definition (and induction on truncation), $a$ is merely accessible, and hence accessible, contradicting $a\in B$.
  Thus, $B$ is empty, so $<$ is well-founded.
\end{proof}

\begin{defn}
  A well-founded relation $<$ on a set $A$ is \define{extensional}
  \indexdef{relation!extensional}%
  \indexdef{extensional!relation}%
  if for any $a,b:A$, we have
  \[ \Parens{\fall{c:A} (c<a) \Leftrightarrow (c<b)} \to (a=b). \]
\end{defn}

Note that since $A$ is a set, extensionality is a mere proposition.
This notion of ``extensionality'' is unrelated to function extensionality, and also unrelated to the extensionality of identity types.
\index{axiom!of extensionality}%
Rather, it is a ``local'' counterpart of the axiom of extensionality in classical set theory.

\begin{thm}
  The type of extensional well-founded relations is a set.
\end{thm}
\begin{proof}
  By the univalence axiom, it suffices to show that if $(A,<)$ is extensional and well-founded and $f:(A,<) \cong (A,<)$, then $f=\idfunc[A]$.
  \index{automorphism!of extensional well-founded relations}%
  We prove by induction on $<$ that $f(a)=a$ for all $a:A$.
  The inductive hypothesis is that for all $a'<a$, we have $f(a')=a'$.

  Now since $A$ is extensional, to conclude $f(a)=a$ it is sufficient to show
  \[\fall{c:A}(c<f(a)) \Leftrightarrow (c<a).\]
  However, since $f$ is an automorphism, we have $(c<a) \Leftrightarrow (f(c)<f(a))$.
  But $c<a$ implies $f(c)=c$ by the inductive hypothesis, so $(c<a) \to (c<f(a))$.
  On the other hand, if $c<f(a)$, then $f^{-1}(c)<a$, and so $c = f(f^{-1}(c)) = f^{-1}(c)$ by the inductive hypothesis again; thus $c<a$.
  Therefore, we have $(c<a) \Leftrightarrow (c<f(a))$ for any $c:A$, so $f(a)=a$.
\end{proof}

\begin{defn}\label{def:simulation}
  If $(A,<)$ and $(B,<)$ are extensional and well-founded, a \define{simulation}
  \indexdef{simulation}%
  \indexsee{function!simulation}{simulation}%
  is a function $f:A\to B$ such that
  \begin{enumerate}
  \item if $a<a'$, then $f(a)<f(a')$, and\label{item:sim1}
  \item for all $a:A$ and $b:B$, if $b<f(a)$, then there merely exists an $a'<a$ with $f(a')=b$.\label{item:sim2}
  \end{enumerate}
\end{defn}

\begin{lem}
  Any simulation is injective.
\end{lem}
\begin{proof}
  We prove by double well-founded induction that for any $a,b:A$, if $f(a)=f(b)$ then $a=b$.
  The inductive hypothesis for $a:A$ says that for any $a'<a$, and any $b:B$, if $f(a')=f(b)$ then $a=b$.
  The inner inductive hypothesis for $b:A$ says that for any $b'<b$, if $f(a)=f(b')$ then $a=b'$.

  Suppose $f(a)=f(b)$; we must show $a=b$.
  By extensionality, it suffices to show that for any $c:A$ we have $(c<a)\Leftrightarrow (c<b)$.
  If $c<a$, then $f(c)<f(a)$ by \autoref{def:simulation}\ref{item:sim1}.
  Hence $f(c)<f(b)$, so by \autoref{def:simulation}\ref{item:sim2} there merely exists $c':A$ with $c'<b$ and $f(c)=f(c')$.
  By the inductive hypothesis for $a$, we have $c=c'$, hence $c<b$.
  The dual argument is symmetrical.
\end{proof}

In particular, this implies that in \autoref{def:simulation}\ref{item:sim2} the word ``merely'' could be omitted without change of sense.

\begin{cor}
  If $f:A\to B$ is a simulation, then for all $a:A$ and $b:B$, if $b<f(a)$, there \emph{purely} exists an $a'<a$ with $f(a')=b$.
\end{cor}
\begin{proof}
  Since $f$ is injective, $\sm{a:A} (f(a)=b)$ is a mere proposition.
\end{proof}

We say that a subset $C :\power B$ is an \define{initial segment}
\indexdef{initial!segment}%
\indexsee{segment, initial}{initial segment}%
if $c\in C$ and $b<c$ imply $b\in C$.
The image of a simulation must be an initial segment, while the inclusion of any initial segment is a simulation.
Thus, by univalence, every simulation $A\to B$ is \emph{equal} to the inclusion of some initial segment of $B$.

\begin{thm}
  For a set $A$, let $P(A)$ be the type of extensional well-founded relations on $A$.
  If $\mathord{<_A} : P(A)$ and $\mathord{<_B} : P(B)$ and $f:A\to B$, let $H_{\mathord{<_A}\mathord{<_B}}(f)$ be the mere proposition that $f$ is a simulation.
  Then $(P,H)$ is a standard notion of structure over \uset in the sense of \autoref{sec:sip}.
\end{thm}
\begin{proof}
  We leave it to the reader to verify that identities are simulations, and that composites of simulations are simulations.
  Thus, we have a notion of structure.
  For standardness, we must show that if $<$ and $\prec$ are two extensional well-founded relations on $A$, and $\idfunc[A]$ is a simulation in both directions, then $<$ and $\prec$ are equal.
  Since extensionality and well-foundedness are mere propositions, for this it suffices to have $\fall{a,b:A} (a<b) \Leftrightarrow (a\prec b)$.
  But this follows from \autoref{def:simulation}\ref{item:sim1} for $\idfunc[A]$.
\end{proof}

\begin{cor}\label{thm:wfcat}
  There is a category whose objects are sets equipped with extensional well-founded relations, and whose morphisms are simulations.
\end{cor}

In fact, this category is a poset.

\begin{lem}
  For extensional and well-founded $(A,<)$ and $(B,<)$, there is at most one simulation $f:A\to B$.
\end{lem}
\begin{proof}
  Suppose $f,g:A\to B$ are simulations.
  Since being a simulation is a mere property, it suffices to show $\fall{a:A}(f(a)=g(a))$.
  By induction on $<$, we may suppose $f(a')=g(a')$ for all $a'<a$.
  And by extensionality of $B$, to have $f(a)=g(a)$ it suffices to have $\fall{b:B}(b<f(a)) \Leftrightarrow (b<g(a))$.

  But since $f$ is a simulation, if $b<f(a)$, then we have $a'<a$ with $f(a')=b$.
  By the inductive hypothesis, we have also $g(a')=b$, hence $b<g(a)$.
  The dual argument is symmetrical.
\end{proof}

Thus, if $A$ and $B$ are equipped with extensional and well-founded relations, we may write $A\le B$ to mean there exists a simulation $f:A\to B$.
\autoref{thm:wfcat} implies that if $A\le B$ and $B\le A$, then $A=B$.

\begin{defn}
  An \define{ordinal}
  \indexdef{ordinal}%
  \indexsee{number!ordinal}{ordinal}%
  is a set $A$ with an extensional well-founded relation which is \emph{transitive}, i.e.\ satisfies $\fall{a,b,c:A}(a<b)\to (b<c) \to (a<c)$.
\end{defn}

\begin{eg}
  Of course, the usual strict order on \nat is transitive.
  It is easily seen to be extensional as well; thus it is an ordinal.
  As usual, we denote this ordinal by $\omega$.
\end{eg}

\symlabel{ord}
Let \ord denote the type of ordinals.
By the previous results, \ord is a set and has a natural partial order.
We now show that \ord also admits a well-founded relation.

\symlabel{initial-segment}
If $A$ is an ordinal and $a:A$, let $\ordsl A a \defeq \setof{ b:A | b<a}$ denote the initial segment.
\index{initial!segment}%
Note that if $\ordsl A a = \ordsl A b$ as ordinals, then that isomorphism must respect their inclusions into $A$ (since simulations form a poset), and hence they are equal as subsets of $A$.
Therefore, since $A$ is extensional, $a=b$.
Thus the function $a\mapsto \ordsl A a$ is an injection $A\to \ord$.

\begin{defn}
  For ordinals $A$ and $B$, a simulation $f:A\to B$ is said to be \define{bounded}
  \indexdef{simulation!bounded}%
  \indexdef{bounded!simulation}%
  if there exists $b:B$ such that $A = \ordsl B b$.
\end{defn}

The remarks above imply that such a $b$ is unique when it exists, so that boundedness is a mere property.

We write $A<B$ if there exists a bounded simulation from $A$ to $B$.
Since simulations are unique, $A<B$ is also a mere proposition.

\begin{thm}\label{thm:ordord}
  $(\ord,<)$ is an ordinal.
\end{thm}

\noindent
More precisely, this theorem says that the type $\ord_{\UU_i}$ of ordinals in one universe\index{universe level} is itself an ordinal in the next higher universe, i.e.\ $(\ord_{\UU_i},<):\ord_{\UU_{i+1}}$.

\begin{proof}
  Let $A$ be an ordinal; we first show that $\ordsl A a$ is accessible (in \ord) for all $a:A$.
  By well-founded induction on $A$, suppose $\ordsl A b$ is accessible for all $b<a$.
  By definition of accessibility, we must show that $B$ is accessible in \ord for all $B<\ordsl A a$.
  However, if $B<\ordsl A a$ then there is some $b<a$ such that $B = \ordsl{(\ordsl A a)}{b} = \ordsl A b$, which is accessible by the inductive hypothesis.
  Thus, $\ordsl A a$ is accessible for all $a:A$.

  Now to show that $A$ is accessible in \ord, by definition we must show $B$ is accessible for all $B<A$.
  But as before, $B<A$ means $B=\ordsl A a$ for some $a:A$, which is accessible as we just proved.
  Thus, \ord is well-founded.

  For extensionality, suppose $A$ and $B$ are ordinals such that
  \narrowequation{\prd{C:\ord} (C<A) \Leftrightarrow (C<B).}
  Then for every $a:A$, since $\ordsl A a<A$, we have $\ordsl A a<B$, hence there is $b:B$ with $\ordsl A a = \ordsl B b$.
  Define $f:A\to B$ to take each $a$ to the corresponding $b$; it is straightforward to verify that $f$ is an isomorphism.
  Thus $A\cong B$, hence $A=B$ by univalence.

  Finally, it is easy to see that $<$ is transitive.
\end{proof}

Treating \ord as an ordinal is often very convenient, but it has its pitfalls as well.
For instance, consider the following lemma, where we pay attention to how universes are used.

\begin{lem}\label{thm:ordsucc}
  Let \bbU be a universe.
  For any $A:\ord_\bbU$, there is a $B:\ord_\bbU$ such that $A<B$.
\end{lem}
\begin{proof}
  Let $B=A+\unit$, with the element $\ttt:\unit$ being greater than all elements of $A$.
  Then $B$ is an ordinal and it is easy to see that $A\cong \ordsl B \ttt$.
\end{proof}

This lemma illustrates a potential pitfall of the ``typically ambiguous''\index{typical ambiguity} style of using \UU to denote an arbitrary, unspecified universe.
Consider the following alternative proof of it.

\begin{proof}[Another putative proof of \autoref{thm:ordsucc}]
  Note that $C<A$ if and only if $C=\ordsl A a$ for some $a:A$.
  This gives an isomorphism $A \cong \ordsl \ord A$, so that $A<\ord$.
  Thus we may take $B\defeq\ord$.
\end{proof}

The second proof would be valid if we had stated \autoref{thm:ordsucc} in a typically ambiguous style.
But the resulting lemma would be less useful, because the second proof would constrain the second ``\ord'' in the lemma statement to refer to a higher universe level than the first one.
The first proof allows both universes to be the same.

Similar remarks apply to the next lemma, which could be proved in a less useful way by observing that $A\le \ord$ for any $A:\ord$.

\begin{lem}\label{thm:ordunion}
  Let \bbU be a universe.
  For any $X:\type$ and $F:X\to \ord_\bbU$, there exists $B:\ord_\bbU$ such that $Fx\le B$ for all $x:X$.
\end{lem}
\begin{proof}
  Let $B$ be the quotient of the equivalence relation $\eqr$ on $\sm{x:X} Fx$ defined as follows:
  \[ (x,y) \eqr (x',y')
  \;\defeq\;
  \Big(\ordsl{(Fx)}{y} \cong \ordsl{(Fx')}{y'}\Big).
  \]
  Define $(x,y)<(x',y')$ if $\ordsl{(Fx)}{y} < \ordsl{(Fx')}{y'}$.
  This clearly descends to the quotient, and can be seen to make $B$ into an ordinal.
  Moreover, for each $x:X$ the induced map $Fx\to B$ is a simulation.
\end{proof}

\section{Classical well-orderings}
\label{sec:wellorderings}

\index{denial|(}%
We now show the equivalence of our ordinals with the more familiar classical\index{mathematics!classical} well-orderings.

\begin{lem}
  \index{excluded middle}%
  Assuming excluded middle, every ordinal is trichotomous:
  \index{trichotomy of ordinals}%
  \index{ordinal!trichotomy of}%
  \[ \fall{a,b:A} (a<b) \vee (a=b) \vee (b<a). \]
\end{lem}
\begin{proof}
  By induction on $a$, we may assume that for every $a'<a$ and every $b':A$, we have $(a'<b') \vee (a'=b') \vee (b'<a')$.
  Now by induction on $b$, we may assume that for every $b'<b$, we have $(a<b') \vee (a=b') \vee (b'<a)$.

  By excluded middle, either there merely exists a $b'<b$ such that $a<b'$, or there merely exists a $b'<b$ such that $a=b'$, or for every $b'<b$ we have $b'<a$.
  In the first case, merely $a<b$ by transitivity, hence $a<b$ as it is a mere proposition.
  Similarly, in the second case, $a<b$ by transport.
  Thus, suppose $\fall{b':A}(b'<b)\to (b'<a)$.

  Now analogously, either there merely exists $a'<a$ such that $b<a'$, or there merely exists $a'<a$ such that $a'=b$, or for every $a'<a$ we have $a'<b$.
  In the first and second cases, $b<a$, so we may suppose $\fall{a':A}(a'<a)\to (a'<b)$.
  However, by extensionality, our two suppositions now imply $a=b$.
\end{proof}

\begin{lem}
  A well-founded relation contains no cycles, i.e.\
  \[ \fall{n:\mathbb{N}}{a:\mathbb{N}_n\to A} \neg\Big((a_0<a_1) \wedge \dots \wedge (a_{n-1}<a_n)\wedge (a_n<a_0)\Big). \]
\end{lem}
\begin{proof}
  We prove by induction on $a:A$ that there is no cycle containing $a$.
  Thus, suppose by induction that for all $a'<a$, there is no cycle containing $a'$.
  But in any cycle containing $a$, there is some element less than $a$ and contained in the same cycle.
\end{proof}

\indexdef{relation!irreflexive}%
\index{irreflexivity!of well-founded relation}%
In particular, a well-founded relation must be \define{irreflexive}, i.e.\ $\neg(a<a)$ for all $a$.

\begin{thm}\label{thm:wellorder}
  Assuming excluded middle, $(A,<)$ is an ordinal if and only if every nonempty subset $B\subseteq A$ has a least element.
\end{thm}
\begin{proof}
  If $A$ is an ordinal, then by \autoref{thm:wfmin} every nonempty subset merely has a minimal element.
  But trichotomy implies that any minimal element is a least element.
  Moreover, least elements are unique when they exist, so merely having one is as good as having one.

  Conversely, if every nonempty subset has a least element, then by \autoref{thm:wfmin}, $A$ is well-founded.
  We also have trichotomy, since for any $a,b$ the subset
  $ \setof{a,b} \defeq \setof{x:A | x=a \lor x=b} $
  merely has a least element, which must be either $a$ or $b$.
  This implies transitivity, since if $a<b$ and $b<c$, then either $a=c$ or $c<a$ would produce a cycle.
  Similarly, it implies extensionality, for if $\fall{c:A}(c<a)\Leftrightarrow (c<b)$, then $a<b$ implies (letting $c$ be $a$) that $a<a$, which is a cycle, and similarly if $b<a$; hence $a=b$.
\end{proof}

In classical\index{mathematics!classical} mathematics, the characterization of \autoref{thm:wellorder} is taken as the definition of a \emph{well-ordering}, with the \emph{ordinals} being a canonical set of representatives of isomorphism classes for well-orderings.
In our context, the structure identity principle means that there is no need to look for such representatives: any well-ordering is as good as any other.

We now move on to consider consequences of the axiom of choice.
For any set $X$, let $\powerp X$ denote the type of merely inhabited subsets of $X$:
\symlabel{inhabited-powerset}
\[ \powerp X \defeq \setof{ Y : \power X | \exis{x:X} x\in Y}. \]
Assuming excluded middle, this is equivalently the type of \emph{nonempty}\index{nonempty subset} subsets of $X$, and we have $\power X \eqvsym (\powerp X) + \unit$.

\begin{thm}\label{thm:wop}
  \index{axiom!of choice}%
  \index{excluded middle}%
  Assuming excluded middle, the following are equivalent.
  \begin{enumerate}
  \item For every set $X$, there merely exists a function
    $ f: \powerp X \to X $
    such that $f(Y)\in Y$ for all $Y:\power X$.\label{item:wop1}
  \item Every set merely admits the structure of an ordinal.\label{item:wop2}
  \end{enumerate}
\end{thm}

\noindent
Of course,~\ref{item:wop1} is a standard classical\index{mathematics!classical} version of the axiom of choice; see \autoref{ex:choice-function}.

\begin{proof}
  One direction is easy: suppose~\ref{item:wop2}.
  Since we aim to prove the mere proposition~\ref{item:wop1}, we may assume $A$ is an ordinal.
  But then we can define $f(B)$ to be the least element of $B$.

  Now suppose~\ref{item:wop1}.
  As before, since~\ref{item:wop2} is a mere proposition, we may assume given such an $f$.
  We extend $f$ to a function
  \[ \bar f:\power X \eqvsym (\powerp X) + \unit \longrightarrow X+\unit
  \]
  in the obvious way.
  Now for any ordinal $A$, we can define $g_A:A\to X+\unit$ by well-founded recursion:
  \[ g_A(a) \defeq 
    \bar f\Big(X \setminus \setof{ g_A(b) | \strut (b<a) \wedge (g_A(b) \in X) }\Big)
  \]
  (regarding $X$ as a subset of $X+\unit$ in the obvious way).

  Let $A'\defeq \setof{a:A | g_A(a) \in X}$ be the preimage of $X\subseteq X+\unit$; then we claim the restriction $g_A':A' \to X$ is injective.
  For if $a,a':A$ with $a\neq a'$, then by trichotomy and without loss of generality, we may assume $a'<a$.
  Thus $g_A(a') \in \setof{ g_A(b) | b<a }$, so since $f(Y)\in Y$ for all $Y$ we have $g_A(a) \neq g_A(a')$.

  Moreover, $A'$ is an initial segment of $A$.
  For $g_A(a)$ lies in \unit if and only if $\setof{g_A(b)|b<a} = X$, and if this holds then it also holds for any $a'>a$.
  Thus, $A'$ is itself an ordinal.

  Finally, since \ord is an ordinal, we can take $A\defeq\ord$.
  Let $X'$ be the image of $g_\ord':\ord' \to X$; then the inverse of $g_\ord'$ yields an injection $H:X'\to \ord$.
  By \autoref{thm:ordunion}, there is an ordinal $C$ such that $Hx\le C$ for all $x:X'$.
  Then by \autoref{thm:ordsucc}, there is a further ordinal $D$ such that $C<D$, hence $Hx<D$ for all $x:X'$.
  Now we have
  \begin{align*}
    g_{\ord}(D) &= \bar f\Big( X \setminus \setof{ g_\ord(B) | \rule{0pt}{1em} B<D \wedge (g_\ord(B) \in X)} \Big)\\
    &=\bar f\Big( X \setminus \setof{ g_\ord(B) | \rule{0pt}{1em} g_\ord(B) \in X} \Big)
  \end{align*}
  since if $B:\ord$ and $(g_\ord(B) \in X)$, then $B = Hx$ for some $x:X'$, hence $B<D$.
  Now if
  \[\setof{ g_\ord(B) | \rule{0pt}{1em} g_\ord(B) \in X}\]
  is not all of $X$, then $g_\ord(D)$ would lie in $X$ but not in this subset, which would be a contradiction since $D$ is itself a potential value for $B$.
  So this set must be all of $X$, and hence $g_\ord'$ is surjective as well as injective.
  Thus, we can transport the ordinal structure on $\ord'$ to $X$.
\end{proof}

\begin{rmk}
  If we had given the wrong proof of \autoref{thm:ordsucc} or \autoref{thm:ordunion}, then the resulting proof of \autoref{thm:wop} would be invalid: there would be no way to consistently assign universe levels\index{universe level}.
  As it is, we require propositional resizing (which follows from \LEM{}) to ensure that $X'$ lives in the same universe as $X$ (up to equivalence).
\end{rmk}

\begin{cor}
  Assuming the axiom of choice, the function $\ord\to\set$ (which forgets the order structure) is a surjection.
\end{cor}

Note that \ord is a set, while \set is a 1-type.
In general, there is no reason for a 1-type to admit any surjective function from a set.
Even the axiom of choice does not appear to imply that \emph{every} 1-type does so (although see \autoref{ex:acnm-surjset}), but it readily implies that this is so for 1-types constructed out of \set, such as the types of objects of categories of structures as in \autoref{sec:sip}.
The following corollary also applies to such categories.

\begin{cor}
  \index{weak equivalence!of precategories}%
  Assuming \choice{}, \uset admits a weak equivalence functor from a strict category.
\end{cor}
\begin{proof}
  Let $X_0\defeq \ord$, and for $A,B:X_0$ let $\hom_X(A,B) \defeq (A\to B)$.
  Then $X$ is a strict category, since \ord is a set, and the above surjection $X_0 \to \set$ extends to a weak equivalence functor $X\to \uset$.
\end{proof}

Now recall from \autoref{sec:cardinals} that we have a further surjection $\cd{\nameless}:\set\to\card$, and hence a composite surjection $\ord\to\card$ which sends each ordinal to its cardinality.

\begin{thm}
  Assuming \choice{}, the surjection $\ord\to\card$ has a section.
\end{thm}
\begin{proof}
  There is an easy and wrong proof of this: since \ord and \card are both sets, \choice{} implies that any surjection between them \emph{merely} has a section.
  However, we actually have a canonical \emph{specified} section: because \ord is an ordinal, every nonempty subset of it has a uniquely specified least element.
  Thus, we can map each cardinal to the least element in the corresponding fiber.
\end{proof}

It is traditional in set theory to identify cardinals with their image in \ord: the least ordinal having that cardinality.

It follows that \card also canonically admits the structure of an ordinal: in fact, one isomorphic to \ord.
Specifically, we define by well-founded recursion a function $\aleph:\ord\to\ord$, such that $\aleph(A)$ is the least ordinal having cardinality greater than $\aleph({\ordsl A a})$ for all $a:A$.
Then (assuming \choice{}) the image of $\aleph$ is exactly the image of \card.

\index{denial|)}%

\index{ordinal|)}%

\section{The cumulative hierarchy}
\label{sec:cumulative-hierarchy}

\index{bargaining|(}%
We can define a cumulative hierarchy $V$ of all sets in a given universe $\UU$ as a higher inductive type, in such a way that $V$ is again a set (in a larger universe $\UU'$), equipped with a binary ``membership'' relation $x\in y$ which satisfies the usual laws of set theory.

\begin{defn}\label{defn:V}
  The \define{cumulative hierarchy}
  \indexdef{cumulative!hierarchy, set-theoretic}%
  \indexdef{hierarchy!cumulative, set-theoretic}%
  $V$ relative to a type universe $\UU$ is the
  higher inductive type generated by the following constructors.
  \begin{enumerate}
  \item For every $A : \UU$ and $f : A \to V$, there is an element $\vset(A, f)$ : V.
  \item For all $A, B : \UU$, $f : A \to V$ and $g : B \to V$ such that
    \begin{narrowmultline} \label{eq:V-path}
      \big(\fall{a:A} \exis{b:B} \id[V]{f(a)}{g(b)}\big) \land \narrowbreak
      \big(\fall{b:B} \exis{a:A} \id[V]{f(a)}{g(b)}\big)
    \end{narrowmultline}
    there is a path $\id[V]{\vset(A,f)}{\vset(B,g)}$.
  \item The 0-truncation constructor: for all $x,y:V$ and $p,q:x=y$, we have $p=q$.
  \end{enumerate}
\end{defn}

In set-theoretic language, $\vset(A,f)$ can be understood as the set (in the sense of classical set theory) that is the image of $A$ under $f$, i.e.\ $\setof{ f(a) | a \in A }$.
However, we will avoid this notation, since it would clash with our notation for subtypes (but see~\eqref{eq:class-notation} and \autoref{def:TypeOfElements} below).

The hierarchy $V$ is
bootstrapped from the empty map $\rec\emptyt(V) : \emptyt \to V$, which gives the empty set as $\emptyset = \vset(\emptyt,\rec\emptyt(V))$.
Then the singleton $\{\emptyset\}$ enters $V$ through $\unit \to V$, defined as $\ttt \mapsto \emptyset$, and so
on. The type $V$ lives in the same universe as the base universe $\UU$.

The second constructor of $V$ has a form unlike any we have seen before: it involves not only paths in $V$ (which in \autoref{sec:hittruncations} we claimed were slightly fishy) but truncations of sums of them.
It certainly does not fit the general scheme described in \autoref{sec:naturality}, and thus it may not be obvious what its induction principle should be.
Fortunately, like our first definition of the 0-truncation in \autoref{sec:hittruncations}, it can be re-expressed using auxiliary higher inductive types.
We leave it to the reader to work out the details (see \autoref{ex:cumhierhit}).

\index{induction principle!for cumulative hierarchy}%
At the end of the day, the induction principle for $V$ (written in pattern matching language) says that given $P:V\to \set$, in order to construct $h:\prd{x:V} P(x)$, it suffices to give the following.
\begin{enumerate}
\item For any $f:A\to V$, construct $h(\vset(A,f))$, assuming as given $h(f(a))$ for all $a:A$.
\item Verify that if $f : A \to V$ and $g : B \to V$ satisfy~\eqref{eq:V-path}, then $h(\vset(A,f)) = h(\vset(B,g))$, assuming inductively that $h(f(a)) = h(g(b))$ whenever $f(a)=g(b)$.
\end{enumerate}
The second clause checks that the map being defined must respect the paths introduced in \eqref{eq:V-path}.
As usual when we state higher induction principles using pattern matching, it may seem tautologous, but is not.
The point is that ``$h(f(a))$'' is essentially a formal symbol which we cannot peek inside of, which $h(\vset(A,f))$ must be defined in terms of. Thus, in the second clause, we assume equality of these formal symbols when appropriate, and verify that the elements resulting from the construction of the first clause are also equal.  
Of course, if $P$ is a family of mere propositions, then the second clause is automatic.

Observe that, by induction, for each $v:V$ there merely exist $A:\UU$ and $f:A\to V$ such that $v=\vset(A,f)$.
Thus, it is reasonable to try to define the \define{membership relation}
\indexdef{membership, for cumulative hierarchy}%
$x\in v$ on $V$ by setting:
%
%
\symlabel{V-membership}
\begin{equation*}
  (x \in \vset(A,f)) \defeq (\exis{a : A} x = f(a)).
\end{equation*}
To see that the definition is valid, we must use the recursion principle of $V$.  Thus, suppose we have a path $\vset(A, f) = \vset(B, g)$
constructed through~\eqref{eq:V-path}. If $x \in \vset(A,f)$ then there merely is $a : A$ such
that $x = f(a)$, but by~\eqref{eq:V-path} there merely is $b : B$ such that $f(a) = g(b)$, hence
$x = g(b)$ and $x \in \vset(B,f)$. The converse is symmetric.

The \define{subset relation}
\indexdef{subset!relation on the cumulative hierarchy}%
$x\subseteq y$ is defined on $V$ as usual by
\begin{equation*}
  (x \subseteq y) \defeq \fall{z : V} z \in x \Rightarrow z \in y.
\end{equation*}

A \define{class}
\indexdef{class}%
may be taken to be a mere predicate on~$V$. We can say that a class $C : V \to \prop$ is a
\define{$V$-set}
\indexdef{set!in the cumulative hierarchy}%
if there merely exists $v\in V$ such that
\begin{equation*}
  \fall{x : V} C(x) \Leftrightarrow x \in v.
\end{equation*}
We may also use the conventional notation for classes, which matches our standard notation for subtypes:
\begin{equation}
  \setof{ x | C(x) } \defeq \lam{x}C(x).\label{eq:class-notation}
\end{equation}
A class $C: V\to \prop$ will be called \define{$\UU$-small}
\indexdef{class!small}%
\indexdef{small!class}%
if all of its values $C(x)$ lie in $\UU$, specifically $C: V\to \prop_{\UU}$.
Since $V$ lives in the same universe $\UU'$ as does the base universe $\UU$ from which it is built, the same is true for the identity types $v=_V w$ for any $v,w:V$. To obtain a well-behaved theory in the absence of propositional resizing,
\index{propositional!resizing}%
\index{resizing}%
therefore, it will be convenient to have a $\UU$-small ``resizing'' of the identity relation, which we can define by induction as follows.

\begin{defn}\label{def:bisimulation}
  Define the \define{bisimulation}
  \indexdef{bisimulation}%
  relation
  \begin{equation*}
    \mathord\bisim : V \times V \longrightarrow \prop_{\UU}
  \end{equation*}
  by double induction over $V$, where for $\vset(A,f)$ and $\vset(B,g)$ we let:
  \begin{narrowmultline*}
    \vset(A,f)  \bisim \vset(B,g) \defeq \narrowbreak
    \big(\fall{a:A}\exis{b:B} f(a)  \bisim g(b)\big) \land
    \big(\fall{b:B}\exis{a:A} f(a) \bisim g(b)\big).
  \end{narrowmultline*}
\end{defn}
To verify that the definition is correct, we just need to check that it respects paths $\vset(A, f) = \vset(B, g)$ constructed through~\eqref{eq:V-path}, but this is obvious, and that $\prop_{\UU}$ is a set, which it is.  Note that $u \bisim v$ is in $\propU$ by construction.

\begin{lem}\label{lem:BisimEqualsId}
For any $u,v:V$ we have $(u=_V v) = (u \bisim v)$.
\end{lem}

\begin{proof}
An easy induction shows that $\bisim$ is reflexive, so by transport we have $(u=_V v)\to (u \bisim v)$.
Thus, it remains to show that $(u \bisim v)\to (u=_V v)$.
By induction on $u$ and $v$, we may assume they are $\vset(A,f)$ and $\vset(B,g)$ respectively.
Then by definition, $\vset(A,f)\bisim\vset(B,g)$ implies $(\fall{a:A}\exis{b:B}f(a)  \bisim g(b))$ and conversely.
But the inductive hypothesis then tells us that $(\fall{a:A}\exis{b:B}f(a) = g(b))$ and conversely.
So by the path-con\-struc\-tor for $V$ we have $\vset(A,f) =_V \vset(B,g)$.
\end{proof}

Now we can use the resized identity relation to get the following useful principle.

\begin{lem}\label{lem:MonicSetPresent}
For every $u:V$ there is a given $A_u:\UU$ and monic $m_u: A_u \mono V$ such that $u = \vset(A_u, m_u)$.
\end{lem}

\begin{proof}
  Take any presentation $u = \vset(A,f)$ and factor $f:A\to V$ as a surjection followed by an injection:
  \begin{equation*}
    f = m_u\circ e_u : A \epi A_u \mono V.    
  \end{equation*}
  Clearly $u = \vset(A_u, m_u)$ if only $A_u$ is still in $\UU$, which holds if the kernel of $e_u : A \epi A_u$ is in $\UU$.  But the kernel of $e_u : A \epi A_u$ is the pullback along $f : A\to V$ of the identity on $V$, which we just showed to be $\UU$-small, up to equivalence.  Now, this construction of the pair $(A_u, m_u)$ with $m_u :A_u \mono V$ and $u = \vset(A_u, m_u)$ from $u:V$ is unique up to equivalence over $V$, and hence up to identity by univalence.  Thus by the principle of unique choice \eqref{cor:UC} there is a map $c : V\to\sm{A:\UU}(A\to V)$ such that $c(u) = (A_u, m_u)$, with $m_u :A_u \mono V$ and $u = \vset(c(u))$, as claimed.
\end{proof}

\begin{defn}\label{def:TypeOfElements}
  For $u:V$, the just constructed monic presentation $m_u: A_u \mono V$ such that $u = \vset(A_u, m_u)$ may be called the \define{type of members}
  \indexdef{type!of members}%
  of $u$ and denoted $m_u : [u] \mono V$, or even $[u] \mono V$.  We can think of $[u]$ as the ``subclass of $V$ consisting of members of $u$''.
\end{defn}

\begin{thm}\label{thm:VisCST}
  \index{axiom!of set theory, for the cumulative hierarchy}%
  The following hold for $(V, {\in})$:
  \begin{enumerate}
  \item \emph{extensionality:}
    \begin{equation*}
      \fall{x, y : V} x \subseteq y \land y \subseteq x \Leftrightarrow x = y.
    \end{equation*}
     \item \emph{empty set:} for all $x:V$, we have $\neg (x\in \emptyset)$.
    \item \emph{pairing:} for all $u, v:V$, the class $u\cup v \defeq \setof{ x | x = u \vee x = v}$ is a $V$-set.
    \item \emph{infinity:}\index{axiom!of infinity}  there is a $v:V$ with $\emptyset\in v$ and $x\in v$ implies $x\cup \{x\}\in v$.
  \item \emph{union:} for all $v:V$, the class $\cup v\defeq \setof{ x | \exis{u:V} x \in u \in v}$ is a $V$-set.
    \item \emph{function set:} for all $u, v:V$, the class $v^u \defeq \setof{ x | x : u\to v}$ is a $V$-set.%
      \footnote{Here $x:u\to v$ means that $x$ is an appropriate set of ordered pairs, according to the usual way of encoding functions in set theory.}
   \item \emph{$\in$-induction:} if $C : V \to \prop$ is a class such that $C(a)$ holds whenever $C(x)$ for all $x\in a$, then $C(v)$ for all $v:V$.
     \item \emph{replacement:}\index{axiom!of replacement} given any $r : V \to V$ and $a : V$, the class 
       \begin{equation*}
         \setof{ x | \exis{y : V} y \in a \land x = r(y)}
       \end{equation*}
       is a $V$-set.
   \item \emph{separation:}\index{axiom!of separation}  given any $a : V$ and $\UU$-small $C : V \to \propU$, the class
     \begin{equation*}
       \setof{ x | x \in a \land C(x)}
     \end{equation*}
     is a $V$-set.
  \end{enumerate}
\end{thm}

\begin{proof}[Sketch of proof]
  \mbox{}
  \begin{enumerate}
  \item Extensionality: if $\vset(A,f) \subseteq \vset(B, g)$ then $f(a) \in \vset(B, g)$
    for every $a : A$, therefore for every $a : A$ there merely exists $b : B$ such that
    $f(a) = g(b)$. The assumption $\vset(B, g) \subseteq \vset(A, f)$ gives the other half
    of~\eqref{eq:V-path}, therefore $\vset(A,f) = \vset(B,g)$.
    
  \item Empty set: suppose $x\in \emptyset = \vset(\emptyt,\rec\emptyt(V))$.  Then $\exis{a:\emptyt}x=\, \rec\emptyt(V,a)$, which is absurd.
  
  \item Pairing: given $u$ and $v$, let $w=\vset(\bool,\rec\bool(V,u,v))$.
    \index{pair!unordered}

  \item Infinity: take $w = \vset(\nat,I)$, where $I: \nat \to V$ is given by the recursion $I(0) \defeq \emptyset$ and $I(n+1) \defeq I(n)\cup \{I(n)\}$.
 
  \item Union: Take any $v:V$ and any presentation $f :A\to V$ with $v=\vset(A,f)$.  Then let $\tilde{A} \defeq \sm{a:A}[fa]$, where $m_{fa} : [fa] \mono V$ is the type of members from \autoref{def:TypeOfElements}.  $\tilde{A}$ is plainly $\UU$-small, and we have $\cup v \defeq \vset(\tilde{A}, \lam{x} m_{f(\proj1(x))}(\proj2(x)))$.
  
  \item Function set: given $u, v:V$, take the types of elements $[u] \mono V$ and $[u] \mono V$, and the function type $[u]\to [v]$.  We want to define a map
  \[
 r: ([u]\to [v])\ \longrightarrow\ V
  \] 
   with ``$r(f) = \setof{ \pairr{x, f(x)} | x : [u] }$'', but in order for this to make sense we must first define the ordered pair $\pairr{x, y}$, and then we take the map $r': x \mapsto \pairr{x, f(x)}$, and then we can put $r(f)\defeq \vset([u], r')$.  But the ordered pair can be defined in terms of unordered pairing as usual.
   
  \item $\in$-induction: let $C : V \to \prop$ be a class such that $C(a)$ holds whenever $C(x)$ for all $x\in a$, and take any $v=\vset(B,g)$.  To show that $C(v)$ by induction, assume that $C(g(b))$ for all $b:B$.  For every $x\in v$ there merely exists some $b:B$ with $x = g(b)$, and so $C(x)$.  Thus $C(v)$.

  \item Replacement: the statement ``$C$ is a $V$-set'' is a mere proposition, so we may
    proceed by induction as follows. Supposing $x$ is $\vset(A, f)$, we claim that $w
    \defeq \vset(A, r \circ f)$ is the set we are looking for.  If $C(y)$ then there merely exists
    $z : V$ and $a : A$ such that $z = f(a)$ and $y = r(z)$, therefore $y \in w$.
    Conversely, if $y \in w$ then there merely exists $a : A$ such that $y = r(f(a))$, so
    if we take $z \defeq f(a)$ we see that $C(y)$ holds.

  \item Let us say that a class $C: V\to\prop$ is \define{separable}
    \indexdef{class!separable}%
    \indexdef{separable class}%
    if for any $a:V$ the class
  \symlabel{class-intersection}
  \begin{equation*}
    a \cap C \defeq\setof{x | x\in a \wedge C(x)}
  \end{equation*}
  is a $V$-set.
We need to show that any $\UU$-small  $C: V \to \propU$ is separable. Indeed, given $a=\vset(A,f)$, let $A' = \exis{x:A}C(fx)$, and take $f' = f\circ i$, where $i : A' \to A$ is the obvious inclusion.  Then we can take $a' = \vset(A',f')$ and we have $x\in a\wedge C(x) \Leftrightarrow x\in a'$ as claimed.  We needed the assumption that $C$ lands in $\UU$ in order for $A' = \exis{x:A}C(fx)$ to be in $\UU$.\qedhere
\end{enumerate}
\end{proof}

It is also convenient to have a strictly syntactic criterion of separability, so that one can read off from the expression for a class that it produces a $V$-set.  One such familiar condition is being ``$\Delta_0$'', which means that the expression is built up from equality $x=_V y$ and membership $x\in y$, using only mere-propositional connectives $\neg$, $\land$, $\lor$, $\Rightarrow$ and quantifiers $\forall$, $\exists$ over particular sets, i.e.\ of the form $\exists(x\in a)$ and $\forall(y\in b)$ (these are called \define{bounded} quantifiers\index{bounded!quantifier}\index{quantifier!bounded}).\indexdef{separation!.Delta0@$\Delta_0$}%

\begin{cor}\label{cor:Delta0sep}
If the class $C: V \to \prop$ is $\Delta_0$ in the above sense, then it is separable.
\end{cor}
\index{axiom!of $\Delta_0$-separation}%

\begin{proof}
Recall that we have a $\UU$-small resizing $x \bisim y$ of identity $x = y$. Since $x\in y$ is defined in terms of $x=y$, we also have a $\UU$-small resizing of membership
\symlabel{resized-membership}
\begin{equation*}
  x\bin\vset(A,f) \defeq \exis{a:A} x \bisim f(a).
\end{equation*}
Now, let $\Phi$ be a $\Delta_0$ expression for $C$, so that as classes $\Phi = C$ (strictly speaking, we should distinguish expressions from their meanings, but we will blur the difference). Let $\widetilde{\Phi}$ be the result of replacing all occurrences of $=$ and $\in$ by their resized equivalents $\bisim$ and $\bin$.  Clearly then $\widetilde{\Phi}$ also expresses $C$, in the sense that for all $x:V$, $\widetilde{\Phi}(x) \Leftrightarrow C(x)$, and hence $\widetilde{\Phi}=C$ by univalence.  It now suffices to show that $\widetilde{\Phi}$ is $\UU$-small, for then it will be separable by the theorem.  

We show that  $\widetilde{\Phi}$ is $\UU$-small by induction on the construction of the expression.  The base cases are $x \bisim y$ and $x\bin y$, which have already been resized into $\UU$.  It is also clear that $\UU$ is closed under the mere-propositional operations (and $(-1)$-truncation), so it just remains to check the bounded quantifiers $\exists(x\in a)$ and $\forall(y\in b)$.  By definition,  
\begin{align*}
\exists(x\in a) P(x) &\defeq \Brck {\sm{x:V}(x\bin a \land P(x))},\\
\forall(y\in b) P(x) &\defeq  \prd{x:V}(x\bin a \to P(x)).
\end{align*} 
Let us consider $\brck {\sm{x:V}(x\bin a \land P(x))}$.  Although the body $(x\bin a \land P(x))$ is $\UU$-small since $P(x)$ is so by the inductive hypothesis, the quantification over $V$ need not stay inside $\UU$.  However, in the present case we can replace this with a quantification over the type $[a]\mono V$ of members of $a$, and easily show that
\begin{equation*}
  \sm{x:V}(x\bin a \land P(x)) = \sm{x:[a]} P(x).
\end{equation*}
The right-hand side does remain in $\UU$, since both $[a]$ and $P(x)$ are in $\UU$.  The case of $\prd{x:V}(x\bin a \to P(x))$ is analogous, using $\prd{x:V}(x\bin a \to P(x)) = \prd{x:[a]}P(x)$.
\end{proof}

We have shown that in type theory with a universe $\UU$, the cumulative hierarchy $V$ is a model of a ``constructive set theory''
\index{constructive!set theory}%
with many of the standard axioms.
However, as far as we know, it lacks the \emph{strong collection}
\index{axiom!strong collection}%
\index{collection!strong}%
\index{strong!collection}%
and \emph{subset collection}
\index{axiom!subset collection}%
\index{collection!subset}%
\index{subset!collection}%
axioms which are included in \CZF{}~\cite{AczelCZF}.
In the usual interpretation of this set theory into type theory, these two axioms are consequences of the setoid-like definition of equality; while in other constructed models of set theory, strong collection may hold for other reasons.
We do not know whether either of these axioms holds in our model $(V,\in)$, but it seems unlikely.
Since $V$ is a higher inductive type \emph{inside} the system, rather than being an \emph{external} construction, it is not surprising that it differs in some ways from prior interpretations.

Finally, consider the result of adding the axiom of choice for sets to our type theory, in the form  $\choice{}$ from \autoref{subsec:emacinsets} above.  This has the consequence that $\LEM{}$ then also holds, by \autoref{thm:1surj_to_surj_to_pem}, and so $\set$ is a topos\index{topos} with subobject classifier $\bool$, by \autoref{thm:settopos}.  In this case, we have $\prop = \bool:\UU$, and so \emph{all classes are separable}.
Thus we have shown:

\begin{lem}\label{lem:fullsep}
  In type theory with $\choice{}$, the law of \define{(full) separation}
  \indexdef{separation!full}%
  holds for $V$: given \emph{any} class $C : V \to \prop$ and $a : V$, the class $a \cap C$ is a $V$-set.
\end{lem}

\begin{thm}\label{thm:zfc}
In type theory with $\choice{}$ and a universe $\UU$, the cumulative hierarchy $V$ is a model of Zermelo--Fraenkel\index{set theory!Zermelo--Fraenkel} set theory with choice, ZFC.
\end{thm}

\begin{proof}
We have all the axioms listed in \autoref{thm:VisCST}, plus full separation, so we just need to show that there are power sets\index{power set} $\power a:V$ for all $a:V$.  But since we have $\LEM{}$ these are simply function types $\power a = (a\to\bool)$.  Thus $V$ is a model of Zermelo--Fraenkel set theory ZF. We leave the verification of the set-theoretic axiom of choice from $\choice{}$ as an easy exercise. 
\end{proof}

\index{bargaining|)}%

\sectionNotes

The basic properties one expects of the category of sets date back to the early days of elementary topos theory.
The \emph{Elementary theory of the category of sets} referred to in \autoref{subsec:emacinsets} was introduced by Lawvere\index{Lawvere} in
\cite{lawvere:etcs-long}, as a category-theoretic axiomatization of set theory.
\index{Elementary Theory of the Category of Sets}%
The notion of $\Pi W$-pretopos, regarded as a predicative version of an elementary topos, was introduced in~\cite{MoerdijkPalmgren2002}; see also~\cite{palmgren:cetcs}.

The treatment of the category of sets in \autoref{sec:piw-pretopos} roughly follows that in~\cite{RijkeSpitters}.
The fact that epimorphisms are surjective (\autoref{epis-surj}) is well known in classical mathematics, but is not as trivial as it may seem to prove \emph{predicatively}.
\index{mathematics!predicative}%
The proof in~\cite{Mines/R/R:1988} uses the power set operation (which is impredicative), although it can also be seen as a predicative proof of the weaker statement that a map in a universe $\UU_i$ is surjective if it is an epimorphism in the next universe $\UU_{i+1}$.
A predicative proof for setoids was given by Wilander~\cite{Wilander2010}. 
Our proof is similar to Wilander's, but avoids setoids by using pushouts and univalence.

The implication in \autoref{thm:1surj_to_surj_to_pem} from $\choice{}$ to $\LEM{}$ is an adaptation to homotopy type
theory of a theorem from topos theory due to Diaconescu~\cite{Diaconescu}; it was posed as a problem already by Bishop~\cite[Problem~2]{Bishop1967}.

For the intuitionistic theory of ordinal numbers, see~\cite{taylor:ordinals} and also \cite{JoyalMoerdijk1995}.
Definitions of well-foundedness in type theory by an induction principle, including the inductive predicate of accessibility\index{accessibility}, were studied in~\cite{Huet80,Paulson86,Nordstrom88}, although the idea dates back to Getzen's proof of the consistency\index{consistency!of arithmetic} of arithmetic~\cite{Gentzen36}.

The idea of algebraic set theory, which informs our development in \autoref{sec:cumulative-hierarchy} of the cumulative hierarchy, is due to~\cite{JoyalMoerdijk1995}, but it derives from earlier work by~\cite{AczelCZF}.
\index{algebraic set theory}%
\index{set theory!algebraic}%

\sectionExercises

\begin{ex}
  Following the pattern of $\uset$, we would like to make a category $\utype$ of all types and maps between them (in a given universe $\UU$).  In order for this to be a category in the sense of \autoref{sec:cats}, however, we must first declare $\hom(X,Y) \defeq \pizero{X\to Y}$, with composition defined by induction on truncation from ordinary composition $(Y\to Z) \to (X\to Y) \to (X\to Z)$.  This was defined as the \emph{homotopy precategory of types} in \autoref{ct:hoprecat}.  It is still not a category, however, but only a precategory (its type of objects $\UU$ is not even a $0$-type).  It becomes a category by Rezk completion
  \index{completion!Rezk}%
  (see \autoref{ct:hocat}), and its type of objects can be identified with $\trunc1\type$ by \autoref{ct:ex:hocat}.  Show that the resulting category $\utype$, unlike $\uset$, is not a pretopos.
\end{ex}

\begin{ex}
  Show that if every surjection has a section in the category $\uset$, then the axiom of choice holds.
\end{ex}

\begin{ex}\label{ex:well-pointed}
  Show that with $\LEM{}$, the category $\uset$ is well-pointed,
  \indexdef{category!well-pointed}%
  in the sense that the following statement holds: for any $f, g : A\to B$, if $f \neq g$ then there is a function $a : 1\to A$ such that $f(a) \neq g(a)$.
  Show that the slice category
  \index{category!slice}%
  $\uset/\bool$ consisting of functions $A\to \bool$ and commutative triangles does not have this property.
  (Hint: the terminal object in $\uset/\bool$ is the identity function $\bool \to \bool$, so in this category, there are objects $X$ that have no elements $1\to X$.)
\end{ex}

\begin{ex}
  \index{addition!of ordinal numbers}%
  Prove that if $(A,<_A)$ and $(B,<_B)$ are well-founded, extensional, or ordinals, then so is $A+B$, with $<$ defined by
  \begin{align*}
    (a<a') &\defeq (a<_A a') & \text{for }& a,a':A\\
    (b<b') &\defeq (b<_B b') & \text{for }& b,b':B\\
    (a<b) &\defeq \unit      & \text{for }& (a:A),(b:B)\\
    (b<a) &\defeq \emptyt    & \text{for }& (a:A),(b:B).
  \end{align*}
\end{ex}

\begin{ex}
  \index{multiplication!of ordinal numbers}%
  Prove that if $(A,<_A)$ and $(B,<_B)$ are well-founded, extensional, or ordinals, then so is $A\times B$, with $<$ defined by
  \[ ((a,b) <(a',b')) \defeq (a<_A a') \vee ((a=a') \wedge (b<_B b')). \]
\end{ex}

\begin{ex}
  Define the usual algebraic operations on ordinals, and prove that they satisfy the usual properties.
\end{ex}

\begin{ex}\label{ex:prop-ord}
  Note that $\bool$ is an ordinal, under the obvious relation $<$ such that $\bfalse<\btrue$ only.
  \begin{enumerate}
  \item Define a relation $<$ on $\prop$ which makes it into an ordinal.
  \item Show that $\id[\ord]\bool\prop$ if and only if \LEM{} holds.
  \end{enumerate}
\end{ex}

\begin{ex}\label{ex:ninf-ord}
  Recall that we denote \nat by $\omega$ when regarding it as an ordinal; thus we have also the ordinal $\omega+1$.
  On the other hand, let us define
  \[ \nat_\infty \defeq \setof{a:\nat\to\bool | \fall{n:\nat} (a_n \le a_{\suc(n)}) } \]
  where $\le$ denotes the obvious partial order on $\bool$, with $\bfalse\le\btrue$.
  \begin{enumerate}
  \item Define a relation $<$ on $\nat_\infty$ which makes it into an ordinal.
  \item Show that $\id[\ord]{\omega+1}{\nat_\infty}$ if and only if the limited principle of omniscience~\eqref{eq:lpo} holds.%
    \index{limited principle of omniscience}%
  \end{enumerate}
\end{ex}

\begin{ex}
  Show that if $(A,<)$ is well-founded and extensional and $A:\UU$, then there is a simulation $A\to V$, where $(V,\in)$ is the cumulative hierarchy from \autoref{sec:cumulative-hierarchy} built from the universe~\UU.
\end{ex}

\begin{ex}\label{ex:choice-function}
  Show that \autoref{thm:wop}\ref{item:wop1} is equivalent to the axiom of choice~\eqref{eq:ac}.
\end{ex}

\begin{ex}\label{ex:cumhierhit}
  Given types $A$ and $B$, define a \define{bitotal relation}
  \indexsee{bitotal relation}{relation, bitotal}%
  \indexdef{relation!bitotal}%
  to be $R:A\to B\to \prop$ such that
  \[ \Big(\fall{a:A}\exis{b:B} R(a,b) \Big) \land \Big(\fall{b:B}\exis{a:A} R(a,b) \Big). \]
  For such $A,B,R$, let $A\sqcup^R B$ be the higher inductive type generated by
  \begin{itemize}
  \item $i:A\to A\sqcup^R B$
  \item $j:B\to A\sqcup^R B$
  \item For each $a:A$ and $b:B$ such that $R(a,b)$, a path $i(a)=j(b)$.
  \end{itemize}
  Show that the cumulative hierarchy $V$ can be defined by the following more straightforward list of constructors, and that the resulting induction principle is the one given in \autoref{sec:cumulative-hierarchy}.
  \begin{itemize}
  \item For every $A : \UU$ and $f : A \to V$, there is an element $\vset(A, f) : V$.
  \item For any $A,B:\UU$ and bitotal relation
    \index{relation!bitotal}%
    $R:A\to B\to \prop$, and any map $h:A\sqcup^R B \to V$, there is a path $\id{\vset(A,h\circ i)}{\vset(B,h\circ j)}$.
  \item The 0-truncation constructor.
  \end{itemize}
\end{ex}

\begin{ex}
  In \CZF, the \define{axiom of strong collection}
  \indexdef{axiom!strong collection}%
  \indexdef{collection!strong}%
  \indexdef{strong!collection}%
  has the form:
   \begin{multline*}
   \fall{x\in v}\exis{y} R(x,y) \Rightarrow \\
   \exis{w}\big[\big(\fall{x\in v}\exis{y\in w}R(x,y)\big)\land \big(\fall{y\in w}\exis{x\in v}R(x,y) \big)\big]
   \end{multline*}
   Does it hold in the cumulative hierarchy $V$?  (We do not know the answer to this.)
\end{ex}

\begin{ex}
Verify that, if we assume $\choice{}$, then the cumulative hierarchy $V$ satisfies the usual set-theoretic axiom of choice, which may be stated in the form:
  \[
   \fall{x\in V} \fall{y\in x}\exis{z\in V} z\in y \Rightarrow  \exis{c\in(\cup x)^x}\fall{y\in x} c(y)\in y
   \]
\end{ex}

\index{set|)}%


\chapter{Real numbers}
\label{cha:real-numbers}

\index{real numbers|(}%
Any foundation of mathematics worthy of its name must eventually address the construction of real numbers as understood by mathematical analysis, namely as a complete archimedean ordered field.
\index{ordered field}%
There are two notions of completeness. The one by Cauchy requires that the reals be closed under limits of Cauchy sequences\index{Cauchy!sequence}, while the stronger one by Dedekind requires closure under Dedekind cuts.\index{cut!Dedekind}
These lead to two ways of constructing reals, which we study in \autoref{sec:dedekind-reals} and \autoref{sec:cauchy-reals}, respectively. In \autoref{RD-final-field,RC-initial-Cauchy-complete} we characterize the two constructions in terms of universal properties: the Dedekind reals are the final archimedean ordered field, and the Cauchy reals the initial Cauchy complete archimedean ordered field.

In traditional constructive mathematics,
\index{mathematics!constructive}%
real numbers always seem to require certain compromises. For example, the Dedekind reals work better with power sets or some other form of impredicativity, while Cauchy reals work well in the presence of countable choice.
\index{axiom!of choice!countable}%
However, we give a new construction of the Cauchy reals as a higher inductive-inductive type that seems to be a third possibility, which requires neither power sets nor countable choice.

In~\autoref{sec:comp-cauchy-dedek} we compare the two constructions of reals. The Cauchy reals are included in the Dedekind reals. They coincide if excluded middle or countable choice holds, but in general the inclusion might be proper.

In~\autoref{sec:compactness-interval} we consider three notions of compactness of the closed interval~$[0,1]$. We first show that $[0,1]$ is metrically compact\indexdef{metrically compact}\indexdef{compactness!metric} in the sense that it is complete and totally bounded, and that uniformly continuous maps on metrically compact spaces behave as expected. In contrast, the Bolzano--Weierstra\ss{} property that every sequence has a convergent subsequence implies the limited principle of omniscience, which is an instance of excluded middle. Finally, we discuss Heine-Borel compactness. A naive formulation of the finite subcover property does not work, but a proof relevant notion of inductive covers does.
This section is basically standard constructive analysis.

The development of real numbers and analysis in homotopy type theory can be easily made compatible with classical mathematics. By assuming excluded middle and the axiom of choice we get standard classical analysis:\index{mathematics!classical}\index{classical!analysis} the Dedekind and Cauchy reals coincide, foundational questions about the impredicative nature of the Dedekind reals disappear, and the interval is as compact as it could be.

We close the chapter by constructing Conway's surreals as a higher inductive-inductive type in \autoref{sec:surreals};
the construction is more natural in univalent type theory than in  classical set theory.

In addition to the basic theory of \autoref{cha:basics,cha:logic}, as noted above we use ``higher inductive-inductive types'' for the Cauchy reals and the surreals: these combine the ideas of \autoref{cha:hits} with the notion of inductive-inductive type mentioned in \autoref{sec:generalizations}.
We will also frequently use the traditional logical notation described in \autoref{subsec:prop-trunc}, and the fact (proven in \autoref{sec:piw-pretopos}) that our ``sets'' behave the way we would expect.

Note that the total space of the universal cover of the circle, which
in \autoref{subsec:pi1s1-homotopy-theory} played a role similar to ``the real numbers'' in
classical algebraic topology, is \emph{not} the type of reals we are looking for. That
type is contractible, and thus equivalent to the singleton type, so it cannot be equipped
with a non-trivial algebraic structure.

\section{The field of rational numbers}
\label{sec:field-rati-numb}

\indexdef{rational numbers}%
\indexsee{number!rational}{rational numbers}%
We first construct the rational numbers \Q, as the reals can then be seen as a completion
of~\Q. An expert will point out that \Q could be replaced by any approximate field,
\indexdef{field!approximate}%
i.e., a subring of \Q in which arbitrarily precise approximate inverses
\index{inverse!approximate}%
exist. An example is the
ring of dyadic rationals,
\index{rational numbers!dyadic}%
which are those of the form $n/2^k$. 
If we were implementing constructive mathematics on a computer,
an approximate field would be more suitable, but we leave such finesse for those
who care about the digits of~$\pi$.

We constructed the integers \Z in \autoref{sec:set-quotients} as a quotient of $\N\times
\N$, and observed that this quotient is generated by an idempotent. In
\autoref{sec:free-algebras} we saw that \Z is the free group on \unit; we could similarly
show that it is the free commutative ring\index{ring} on \emptyt. The field of rationals \Q is
constructed along the same lines as well, namely as the quotient
\[ \Q \defeq (\Z \times \N)/{\approx} \]
where
\[ (u,a) \approx (v,b) \defeq (u (b + 1) = v (a + 1)). \]
In other words, a pair $(u, a)$ represents the rational number $u / (1 + a)$. There can be
no division by zero because we cunningly added one to the denominator~$a$. Here too we
have a canonical choice of representatives, namely fractions in lowest terms. Thus we may
apply \autoref{lem:quotient-when-canonical-representatives} to obtain a set \Q, which
again has a decidable equality.
\index{decidable!equality}%

We do not bother to write down the arithmetical operations on \Q as we trust our readers
know how to compute with fractions even in the case when one is added to the denominator.
Let us just record the conclusion that there is an entirely unproblematic construction of
the ordered field of rational numbers \Q, with a decidable equality and decidable order.
It can also be characterized as the initial ordered field.
\index{initial!ordered field}%

\symlabel{positive-rationals}
\indexdef{rational numbers!positive}%
\indexdef{positive!rational numbers}%
Let $\Qp = \setof{ q : \Q | q > 0 }$ be the type of positive rational numbers.

\section{Dedekind reals}
\label{sec:dedekind-reals}

\index{real numbers!Dedekind|(}%
Let us first recall the basic idea of Dedekind's construction. We use two-sided Dedekind
cuts, as opposed to an often used one-sided version, because the symmetry makes
constructions more elegant, and it works constructively as well as classically.
\index{mathematics!constructive}%
A \emph{Dedekind cut}\index{cut!Dedekind} consists of a pair $(L, U)$ of subsets $L, U \subseteq \Q$, called the
\emph{lower} and \emph{upper cut} respectively, which are:
\begin{enumerate}
\item \emph{inhabited:} there are $q \in L$ and $r \in U$,
\item \emph{rounded:} $q \in L \Leftrightarrow \exis {r \in \Q} q < r \land r \in L$
  and $r \in U \Leftrightarrow \exis {q \in \Q} q \in U \land q < r$,
  \index{rounded!Dedekind cut}
\item \emph{disjoint:} $\lnot (q \in L \land q \in U)$, and
\item \emph{located:} $q < r \Rightarrow q \in L \lor r \in U$.
  \index{locatedness}%
\end{enumerate}
Reading the roundedness condition from left to right tells us that cuts are \emph{open},
\index{open!cut}%
and from right to left that they are \emph{lower}, respectively \emph{upper}, sets. The
locatedness condition states that there is no large gap between $L$ and $U$. Because cuts
are always open, they never include the ``point in between'', even when it is rational. A
typical Dedekind cut looks like this:
\begin{center}
  \begin{tikzpicture}[x=\textwidth]
    \draw[<-),line width=0.75pt] (0,0) -- (0.297,0) node[anchor=south east]{$L\ $};
    \draw[(->,line width=0.75pt] (0.300, 0) node[anchor=south west]{$\ U$} -- (0.9, 0) ;
  \end{tikzpicture}
\end{center}
We might naively translate the informal definition into type theory by saying that a cut
is a pair of maps $L, U : \Q \to \prop$. But we saw in \autoref{subsec:prop-subsets} that
$\prop$ is an ambiguous\index{typical ambiguity} notation for $\prop_{\UU_i}$ where~$\UU_i$ is a universe. Once we
use a particular $\UU_i$ to define cuts, the type of reals will reside in the next
universe $\UU_{i+1}$, a property of reals two levels higher in $\UU_{i+2}$, a property of
subsets of reals in $\UU_{i+3}$, etc. In principle we should be able to keep track of the
universe levels\index{universe level}, especially with the help of a proof assistant, but doing so here would
just burden us with bureaucracy that we prefer to avoid. We shall therefore make a
simplifying assumption that a single type of propositions $\Omega$ is sufficient for all
our purposes.

In fact, the construction of the Dedekind reals is quite resilient to logical
manipulations. There are several ways in which we can make sense of using a single type
$\Omega$:
\begin{enumerate}

\item We could identify $\Omega$ with the ambiguous $\prop$ and track all the universes
  that appear in definitions and constructions.

\item We could assume the propositional resizing axiom,
  \index{propositional!resizing}%
  as in \autoref{subsec:prop-subsets}, which essentially collapses the $\prop_{\UU_i}$'s to the
  lowest level\index{universe level}, which we call $\Omega$.

\item A classical mathematician who is not interested in the intricacies of type-theoretic
  universes or computation may simply assume the law of excluded middle~\eqref{eq:lem} for
  mere propositions so that $\Omega \jdeq \bool$.
  \index{excluded middle}
  This not only eradicates questions about
  levels\index{universe level} of $\prop$, but also turns everything we do into the standard classical\index{mathematics!classical}
  construction of real numbers.

\item On the other end of the spectrum one might ask for a minimal requirement that makes
  the constructions work. The condition that a mere predicate be a Dedekind cut is
  expressible using only conjunctions, disjunctions, and existential quantifiers\index{quantifier!existential} over~\Q, which
  is a countable set. Thus we could take $\Omega$ to be the initial \emph{$\sigma$-frame},
  \index{initial!sigma-frame@$\sigma$-frame}%
  \index{sigma-frame@$\sigma$-frame!initial|defstyle}%
  i.e., a lattice\index{lattice} with countable joins\index{join!in a lattice} in which binary meets distribute over countable
  joins. (The initial $\sigma$-frame cannot be the two-point lattice $\bool$ because
  $\bool$ is not closed under countable joins, unless we assume excluded middle.) This
  would lead to a construction of~$\Omega$ as a higher inductive-inductive type, but one
  experiment of this kind in \autoref{sec:cauchy-reals} is enough.
\end{enumerate}

In all of the above cases $\Omega$ is a set.
Without further ado, we translate the informal definition into type theory.
Throughout this chapter, we use the
logical notation from \autoref{defn:logical-notation}.

\begin{defn} \label{defn:dedekind-reals}
  A \define{Dedekind cut}
  \indexsee{Dedekind!cut}{cut, Dedekind}%
  \indexdef{cut!Dedekind}%
  is a pair $(L, U)$ of mere predicates $L : \Q \to \Omega$ and $U
  : \Q \to \Omega$ which is:
  \begin{enumerate}
  \item \emph{inhabited:} $\exis{q : \Q} L(q)$ and $\exis{r : Q} U(r)$,
  \item \emph{rounded:} for all $q, r : \Q$,
    \index{rounded!Dedekind cut}
    \begin{align*}
      L(q) &\Leftrightarrow \exis{r : \Q} (q < r) \land L(r)
      \qquad\text{and}\\
      U(r) &\Leftrightarrow \exis{q : \Q} (q < r) \land U(q),
    \end{align*}
  \item \emph{disjoint:} $\lnot (L(q) \land U(q))$ for all $q : \Q$,
  \item \emph{located:} $(q < r) \Rightarrow L(q) \lor U(r)$ for all $q, r : \Q$.
  \index{locatedness}%
  \end{enumerate}
  We let $\dcut(L, U)$ denote the conjunction of these conditions. The type of
  \define{Dedekind reals} is
  \indexsee{Dedekind!real numbers}{real numbers, De\-de\-kind}%
  \indexdef{real numbers!Dedekind}%
  \begin{equation*}
    \RD \defeq \setof{ (L, U) : (\Q \to \Omega) \times (\Q \to \Omega) | \dcut(L,U)}.
  \end{equation*}
\end{defn}

It is apparent that $\dcut(L, U)$ is a mere proposition, and since $\Q \to \Omega$ is a
set the Dedekind reals form a set too. See
\autoref{ex:RD-extended-reals,ex:RD-lower-cuts,ex:RD-interval-arithmetic} for variants of
Dedekind cuts which lead to extended reals, lower and upper reals, and the interval
domain.

There is an embedding $\Q \to \RD$ which associates with each rational $q : \Q$ the cut
$(L_q, U_q)$ where
\begin{equation*}
  L_q(r) \defeq (r < q)
  \qquad\text{and}\qquad
  U_q(r) \defeq (q < r).
\end{equation*}
We shall simply write $q$ for the cut $(L_q, U_q)$ associated with a rational number.

\subsection{The algebraic structure of Dedekind reals}
\label{sec:algebr-struct-dedek}

The construction of the algebraic and order-theoretic structure of Dedekind reals proceeds
as usual in intuitionistic logic. Rather than dwelling on details we point out the
differences between the classical\index{mathematics!classical} and intuitionistic setup. Writing $L_x$ and $U_x$ for
the lower and upper cut of a real number $x : \RD$, we define addition as%
\indexdef{addition!of Dedekind reals}%
\begin{align*}
  L_{x + y}(q) &\defeq \exis{r, s : \Q} L_x(r) \land L_y(s) \land q = r + s, \\
  U_{x + y}(q) &\defeq \exis{r, s : \Q} U_x(r) \land U_y(s) \land q = r + s,
\end{align*}
and the additive inverse by
\begin{align*}
  L_{-x}(q) &\defeq \exis{r : \Q} U_x(r) \land q = - r, \\
  U_{-x}(q) &\defeq \exis{r : \Q} L_x(r) \land q = - r.
\end{align*}
With these operations $(\RD, 0, {+}, {-})$ is an abelian\index{group!abelian} group. Multiplication is a bit
more cumbersome:
\indexdef{multiplication!of Dedekind reals}%
\begin{align*}
  L_{x \cdot y}(q) &\defeq
  \begin{aligned}[t]
    \exis{a, b, c, d : \Q} & L_x(a) \land U_x(b) \land L_y(c) \land U_y(d) \land {}\\
                           & \qquad q < \min (a \cdot c, a \cdot d, b \cdot c, b \cdot d),
  \end{aligned} \\
  U_{x \cdot y}(q) &\defeq
  \begin{aligned}[t]
    \exis{a, b, c, d : \Q} & L_x(a) \land U_x(b) \land L_y(c) \land U_y(d) \land {}\\
                           & \qquad \max (a \cdot c, a \cdot d, b \cdot c, b \cdot d) < q.
  \end{aligned}
\end{align*}
\index{interval!arithmetic}%
These formulas are related to multiplication of intervals in interval arithmetic, where
intervals $[a,b]$ and $[c,d]$ with rational endpoints multiply to the interval
\begin{equation*}
  [a,b] \cdot [c,d] =
  [\min(a c, a d, b c, b d), \max(a c, a d, b c, b d)].
\end{equation*}
For instance, the formula for the lower cut can be read as saying that $q < x \cdot y$
when there are intervals $[a,b]$ and $[c,d]$ containing $x$ and $y$, respectively, such
that $q$ is to the left of $[a,b] \cdot [c,d]$. It is generally useful to think of an
interval $[a,b]$ such that $L_x(a)$ and $U_x(b)$ as an approximation of~$x$, see
\autoref{ex:RD-interval-arithmetic}.

We now have a commutative ring\index{ring} with unit
\index{unit!of a ring}%
$(\RD, 0, 1, {+}, {-}, {\cdot})$. To treat
multiplicative inverses, we must first introduce order. Define $\leq$ and $<$ as
\begin{align*}
  (x \leq y) &\ \defeq \ \fall{q : \Q} L_x(q) \Rightarrow L_y(q), \\
  (x < y)    &\ \defeq \ \exis{q : \Q} U_x(q) \land L_y(q).
\end{align*}

\begin{lem} \label{dedekind-in-cut-as-le}
  For all $x : \RD$ and $q : \Q$, $L_x(q) \Leftrightarrow (q < x)$ and $U_x(q)
  \Leftrightarrow (x < q)$.
\end{lem}

\begin{proof}
  If $L_x(q)$ then by roundedness there merely is $r > q$ such that $L_x(r)$, and since
  $U_q(r)$ it follows that $q < x$. Conversely, if $q < x$ then there is $r : \Q$ such
  that $U_q(r)$ and $L_x(r)$, hence $L_x(q)$ because $L_x$ is a lower set. The other half
  of the proof is symmetric.
\end{proof}

\index{partial order}%
\index{transitivity!of . for reals@of $<$ for reals}
\index{transitivity!of . for reals@of $\leq$ for reals}
\index{relation!irreflexive}
\index{irreflexivity!of . for reals@of $<$ for reals}
The relation $\leq$ is a partial order, and $<$ is transitive and irreflexive. Linearity
\index{order!linear}%
\index{linear order}%
\begin{equation*}
  (x < y) \lor (y \leq x)
\end{equation*}
is valid if we assume excluded middle, but without it we get weak linearity
\index{order!weakly linear}
\index{weakly linear order}
\begin{equation} \label{eq:RD-linear-order}
  (x < y) \Rightarrow (x < z) \lor (z < y).
\end{equation}
At first sight it might not be clear what~\eqref{eq:RD-linear-order} has to do with
linear order. But if we take $x \jdeq u - \epsilon$ and $y \jdeq u + \epsilon$ for
$\epsilon > 0$, then we get
\begin{equation*}
  (u - \epsilon < z) \lor (z < u + \epsilon).
\end{equation*}
This is linearity ``up to a small numerical error'', i.e., since it is unreasonable to
expect that we can actually compute with infinite precision, we should not be surprised
that we can decide~$<$ only up to whatever finite precision we have computed.

To see that~\eqref{eq:RD-linear-order} holds, suppose $x < y$. Then there merely exists $q : \Q$ such that $U_x(q)$ and
$L_y(q)$. By roundedness there merely exist $r, s : \Q$ such that $r < q < s$, $U_x(r)$
and $L_y(s)$. Then, by locatedness $L_z(r)$ or $U_z(s)$. In the first case we get $x < z$
and in the second $z < y$. 

Classically, multiplicative inverses exist for all numbers which are different from zero.
However, without excluded middle, a stronger condition is required. Say that $x, y : \RD$
are \define{apart}
\indexdef{apartness}%
from each other, written $x \apart y$, when $(x < y) \lor (y < x)$:
\symlabel{apart}
\begin{equation*}
  (x \apart y) \defeq (x < y) \lor (y < x).
\end{equation*}
If $x \apart y$, then $\lnot (x = y)$.
The converse is true if we assume excluded middle, but is not provable constructively.
\index{mathematics!constructive}%
Indeed, if $\lnot (x = y)$ implies $x\apart y$, then a little bit of excluded middle follows; see \autoref{ex:reals-apart-neq-MP}.

\begin{thm} \label{RD-inverse-apart-0}
  A real is invertible if, and only if, it is apart from $0$.
\end{thm}

\begin{rmk}
  We observe that a real is invertible if, and only if, it is merely
  invertible.  Indeed, the same is true in any ring,\index{ring} since a ring is a set, and
  multiplicative inverses are unique if they exist.  See the discussion
  following \autoref{cor:UC}.
\end{rmk}

\begin{proof}
  Suppose $x \cdot y = 1$. Then there merely exist $a, b, c, d : \Q$ such that
  $a < x < b$, $c < y < d$ and $0 < \min (a c, a d, b c, b d)$. From $0 < a c$ and $0 < b c$ it follows
  that $a$, $b$, and $c$ are either all positive or all negative.
  Hence either $0 < a < x$ or $x < b < 0$, so that $x \apart 0$.

  Conversely, if $x \apart 0$ then
  \begin{align*}
    L_{x^{-1}}(q) &\defeq
    \exis{r : \Q} U_x(r) \land ((0 < r \land q r < 1) \lor (r < 0 \land 1 < q r))
    \\
    U_{x^{-1}}(q) &\defeq
    \exis{r : \Q} L_x(r) \land ((0 < r \land q r > 1) \lor (r < 0 \land 1 > q r))
  \end{align*}
  defines the desired inverse. Indeed, $L_{x^{-1}}$ and $U_{x^{-1}}$ are inhabited because
  $x \apart 0$.
\end{proof}

\index{ordered field!archimedean}%
\index{dense}%
\indexsee{order-dense}{dense}%
The archimedean principle can be stated in several ways. We find it most illuminating in the
form which says that $\Q$ is dense in $\RD$.

\begin{thm}[Archimedean principle for $\RD$] \label{RD-archimedean}
  For all $x, y : \RD$ if $x < y$ then there merely exists $q : \Q$ such that
  $x < q < y$.
\end{thm}

\begin{proof}
  By definition of $<$.
\end{proof}

Before tackling completeness of Dedekind reals, let us state precisely what algebraic
structure they possess. In the following definition we are not aiming at a minimal
axiomatization, but rather at a useful amount of structure and properties.

\begin{defn} \label{ordered-field} An \define{ordered field}
  \indexdef{ordered field}%
  \indexsee{field!ordered}{ordered field}%
  is a set $F$ together with
  constants $0$, $1$, operations $+$, $-$, $\cdot$, $\min$, $\max$, and mere relations
  $\leq$, $<$, $\apart$ such that:
  \begin{enumerate}
  \item $(F, 0, 1, {+}, {-}, {\cdot})$ is a commutative ring with unit;
    \index{unit!of a ring}%
    \index{ring}%
  \item $x : F$ is invertible if, and only if, $x \apart 0$;
  \item $(F, {\leq}, {\min}, {\max})$ is a lattice;
  \item the strict order $<$ is transitive, irreflexive,
    \index{relation!irreflexive}
    \index{irreflexivity!of . in a field@of $<$ in a field}%
    and weakly linear ($x < y \Rightarrow x < z \lor z < y$);\index{transitivity!of . in a field@of $<$ in a field}
    \index{order!weakly linear}
    \index{weakly linear order}
    \index{strict!order}%
    \index{order!strict}%
  \item apartness $\apart$ is irreflexive, symmetric and cotransitive ($x \apart y \Rightarrow x \apart z \lor y \apart z$);
    \index{relation!irreflexive}
    \index{irreflexivity!of apartness}%
    \indexdef{relation!cotransitive}%
    \index{cotransitivity of apartness}%
  \item for all $x, y, z : F$:
    \begin{align*}
      x \leq y &\Leftrightarrow \lnot (y < x), &
      x < y \leq z &\Rightarrow x < z, \\
      x \apart y &\Leftrightarrow (x < y) \lor (y < x), &
      x \leq y < z &\Rightarrow x < z, \\
      x \leq y &\Leftrightarrow x + z \leq y + z, &
      x \leq y \land 0 \leq z &\Rightarrow x z \leq y z, \\
      x < y &\Leftrightarrow x + z < y + z, &
      0 < z \Rightarrow (x < y &\Leftrightarrow x z < y z), \\
      0 < x + y &\Rightarrow 0 < x \lor 0 < y, &
      0 &< 1.
    \end{align*}
  \end{enumerate}
  Every such field has a canonical embedding $\Q \to F$. An ordered field is
  \define{archimedean}
  \indexdef{ordered field!archimedean}%
  \indexsee{archimedean property}{ordered field, archi\-mede\-an}%
  when for all $x, y : F$, if $x < y$ then there merely exists $q :
  \Q$ such that $x < q < y$.
\end{defn}

\begin{thm} \label{RD-archimedean-ordered-field}
  The Dedekind reals form an ordered archimedean field.
\end{thm}

\begin{proof}
  We omit the proof in the hope that what we have demonstrated so far makes the theorem
  plausible.
\end{proof}

\subsection{Dedekind reals are Cauchy complete}
\label{sec:RD-cauchy-complete}

Recall that $x : \N \to \Q$ is a \emph{Cauchy sequence}\indexdef{Cauchy!sequence} when it satisfies
\begin{equation} \label{eq:cauchy-sequence}
  \prd{\epsilon : \Qp} \sm{n : \N} \prd{m, k \geq n} |x_m - x_k| < \epsilon.
\end{equation}
Note that we did \emph{not} truncate the inner existential because we actually want to
compute rates of convergence---an approximation without an error estimate carries little
useful information. By \autoref{thm:ttac}, \eqref{eq:cauchy-sequence} yields a function $M
: \Qp \to \N$, called the \emph{modulus of convergence}\indexdef{modulus!of convergence}, such that $m, k \geq M(\epsilon)$
implies $|x_m - x_k| < \epsilon$. From this we get $|x_{M(\delta/2)} - x_{M(\epsilon/2)}|<
\delta + \epsilon$ for all $\epsilon : \Qp$. In fact, the map $(\epsilon \mapsto
x_{M(\epsilon/2)}) : \Qp \to \Q$ carries the same information about the limit as the
original Cauchy condition~\eqref{eq:cauchy-sequence}. We shall work with these
approximation functions rather than with Cauchy sequences.

\begin{defn} \label{defn:cauchy-approximation}
  A \define{Cauchy approximation}
  \indexdef{Cauchy!approximation}%
  is a map $x : \Qp \to \RD$ which satisfies
  \begin{equation}
    \label{eq:cauchy-approx}
    \fall{\delta, \epsilon :\Qp} |x_\delta - x_\epsilon| < \delta + \epsilon.
  \end{equation}
  The \define{limit}
  \index{limit!of a Cauchy approximation}%
  of a Cauchy approximation $x : \Qp \to \RD$ is a number $\ell : \RD$ such
  that
  \begin{equation*}
    \fall{\epsilon, \theta : \Qp} |x_\epsilon - \ell| < \epsilon + \theta.
  \end{equation*}
\end{defn}

\begin{thm} \label{RD-cauchy-complete}
  Every Cauchy approximation in $\RD$ has a limit.
\end{thm}

\begin{proof}
  Note that we are showing existence, not mere existence, of the limit.
  Given a Cauchy approximation $x : \Qp \to \RD$, define
  \begin{align*}
    L_y(q) &\defeq \exis{\epsilon, \theta : \Qp} L_{x_\epsilon}(q + \epsilon + \theta),\\
    U_y(q) &\defeq \exis{\epsilon, \theta : \Qp} U_{x_\epsilon}(q - \epsilon - \theta).
  \end{align*}
  It is clear that $L_y$ and $U_y$ are inhabited, rounded, and disjoint. To establish
  locatedness, consider any $q, r : \Q$ such that $q < r$. There is $\epsilon : \Qp$ such
  that $5 \epsilon < r - q$. Since $q + 2 \epsilon < r - 2 \epsilon$ merely
  $L_{x_\epsilon}(q + 2 \epsilon)$ or $U_{x_\epsilon}(r - 2 \epsilon)$. In the first case
  we have $L_y(q)$ and in the second $U_y(r)$.

  To show that $y$ is the limit of $x$, consider any $\epsilon, \theta : \Qp$. Because
  $\Q$ is dense in $\RD$ there merely exist $q, r : \Q$ such that
  \begin{narrowmultline*}
    x_\epsilon - \epsilon - \theta/2 < q < x_\epsilon - \epsilon - \theta/4
    < x_\epsilon < \\
    x_\epsilon + \epsilon + \theta/4 < r < x_\epsilon + \epsilon + \theta/2,
  \end{narrowmultline*}
  and thus $q < y < r$. Now either $y < x_\epsilon + \theta/2$ or $x_\epsilon - \theta/2 < y$.
  In the first case we have
  \begin{equation*}
    x_\epsilon - \epsilon - \theta/2 < q < y < x_\epsilon + \theta/2,
  \end{equation*}
  and in the second
  \begin{equation*}
    x_\epsilon - \theta/2 < y < r < x_\epsilon + \epsilon + \theta/2.
  \end{equation*}
  In either case it follows that $|y - x_\epsilon| < \epsilon + \theta$.
\end{proof}

For sake of completeness we record the classic formulation as well.

\begin{cor}
  Suppose $x : \N \to \RD$ satisfies the Cauchy condition~\eqref{eq:cauchy-sequence}. Then
  there exists $y : \RD$ such that
  \begin{equation*}
    \prd{\epsilon : \Qp} \sm{n : \N} \prd{m \geq n} |x_m - y| < \epsilon.
  \end{equation*}
\end{cor}

\begin{proof}
  By \autoref{thm:ttac} there is $M : \Qp \to \N$ such that $\bar{x}(\epsilon) \defeq
  x_{M(\epsilon/2)}$ is a Cauchy approximation. Let $y$ be its limit, which exists by
  \autoref{RD-cauchy-complete}. Given any $\epsilon : \Qp$, let $n \defeq M(\epsilon/4)$
  and observe that, for any $m \geq n$,
  \begin{narrowmultline*}
    |x_m - y| \leq |x_m - x_n| + |x_n - y| =
    |x_m - x_n| + |\bar{x}(\epsilon/2) - y| < \narrowbreak
    \epsilon/4 + \epsilon/2 + \epsilon/4 = \epsilon.\qedhere
  \end{narrowmultline*}
\end{proof}

\subsection{Dedekind reals are Dedekind complete}
\label{sec:RD-dedekind-complete}

We obtained $\RD$ as the type of Dedekind cuts on $\Q$. But we could have instead started
with any archimedean ordered field $F$ and constructed Dedekind cuts\index{cut!Dedekind} on $F$. These would
again form an archimedean ordered field $\bar{F}$, the \define{Dedekind completion of $F$},%
\index{completion!Dedekind}%
\indexsee{Dedekind!completion}{completion, Dedekind}
with $F$ contained as a subfield. What happens if we apply this construction to
$\RD$, do we get even more real numbers? The answer is negative. In fact, we shall prove a
stronger result: $\RD$ is final.

Say that an ordered field~$F$ is \define{admissible for $\Omega$}
\indexsee{admissible!ordered field}{ordered field, admissible}%
\indexdef{ordered field!admissible}%
when the strict order
$<$ on~$F$ is a map ${<} : F \to F \to \Omega$.

\begin{thm} \label{RD-final-field}
  Every archimedean ordered field which is admissible for $\Omega$ is a subfield of~$\RD$.
\end{thm}

\begin{proof}
  Let $F$ be an archimedean ordered field. For every $x : F$ define $L, U : \Q \to
  \Omega$ by
  \begin{equation*}
    L_x(q) \defeq (q < x)
    \qquad\text{and}\qquad
    U_x(q) \defeq (x < q).
  \end{equation*}
  (We have just used the assumption that $F$ is admissible for $\Omega$.)
  Then $(L_x, U_x)$ is a Dedekind cut.\index{cut!Dedekind} Indeed, the cuts are inhabited and rounded because
  $F$ is archimedean and $<$ is transitive, disjoint because $<$ is irreflexive, and
  located because $<$ is a weak linear order. Let $e : F \to \RD$ be the map $e(x) \defeq (L_x,
  U_x)$.

  We claim that $e$ is a field embedding which preserves and reflects the order. First of
  all, notice that $e(q) = q$ for a rational number $q$. Next we have the equivalences,
  for all $x, y : F$,
  \begin{narrowmultline*}
    x < y \Leftrightarrow
    (\exis{q : \Q} x < q < y) \Leftrightarrow \narrowbreak
    (\exis{q : \Q} U_x(q) \land L_y(q)) \Leftrightarrow
    e(x) < e(y),
  \end{narrowmultline*}
  so $e$ indeed preserves and reflects the order. That $e(x + y) = e(x) + e(y)$ holds
  because, for all $q : \Q$,
  \begin{equation*}
    q < x + y \Leftrightarrow
    \exis{r, s : \Q} r < x \land s < y \land q = r + s.
  \end{equation*}
  The implication from right to left is obvious. For the other direction, if $q < x +
  y$ then there merely exists $r : \Q$ such that $q - y < r < x$, and by taking $s \defeq
  q - r$ we get the desired $r$ and $s$. We leave preservation of multiplication by $e$ as
  an exercise.
\end{proof}

To establish that the Dedekind cuts on $\RD$ do not give us anything new, we need just one
more lemma.

\begin{lem} \label{lem:cuts-preserve-admissibility}
  If $F$ is admissible for $\Omega$ then so is its Dedekind completion.
  \index{completion!Dedekind}%
\end{lem}

\begin{proof}
  Let $\bar{F}$ be the Dedekind completion of $F$. The strict order on $\bar{F}$ is
  defined by
  \begin{equation*}
    ((L,U) < (L',U')) \defeq \exis{q : \Q} U(q) \land L'(q).
  \end{equation*}
  Since $U(q)$ and $L'(q)$ are elements of $\Omega$, the lemma holds as long as $\Omega$
  is closed under conjunctions and countable existentials, which we assumed from the outset.
\end{proof}

\begin{cor} \label{RD-dedekind-complete}
  \indexdef{complete!ordered field, Dedekind}%
  \indexdef{Dedekind!completeness}%
  The Dedekind reals are Dedekind complete: for every real-valued Dedekind cut $(L, U)$
  there is a unique $x : \RD$ such that $L(y) = (y < x)$ and $U(y) = (x < y)$ for all $y :
  \RD$.
\end{cor}

\begin{proof}
  By \autoref{lem:cuts-preserve-admissibility} the Dedekind completion $\barRD$ of $\RD$
  is admissible for $\Omega$, so by \autoref{RD-final-field} we have an embedding $\barRD
  \to \RD$, as well as an embedding $\RD \to \barRD$. But these embeddings must be
  isomorphisms, because their compositions are order-preserving field homomorphisms\index{homomorphism!field} which
  fix the dense subfield~$\Q$, which means that they are the identity. The corollary now
  follows immediately from the fact that $\barRD \to \RD$ is an isomorphism.
\end{proof}

\index{real numbers!Dedekind|)}%

\section{Cauchy reals}
\label{sec:cauchy-reals}

\index{real numbers!Cauchy|(}%
\index{completion!Cauchy|(}%
\indexsee{Cauchy!completion}{completion, Cauchy}%
The Cauchy reals are, by intent, the completion of \Q under limits of Cauchy sequences.\index{Cauchy!sequence}
In the classical construction of the Cauchy reals, we consider the set $\mathcal{C}$ of all Cauchy sequences in \Q and then form a suitable quotient $\mathcal{C}/{\approx}$.
Then, to show that $\mathcal{C}/{\approx}$ is Cauchy complete, we consider a Cauchy sequence $x : \N \to \mathcal{C}/{\approx}$, lift it to a sequence of sequences $\bar{x} : \N \to \mathcal{C}$, and construct the limit of $x$ using $\bar{x}$. However, the lifting of~$x$ to $\bar{x}$ uses
the axiom of countable choice (the instance of~\eqref{eq:ac} where $X=\N$) or the law of excluded middle, which we may wish to avoid.
\indexdef{axiom!of choice!countable}%
Every construction of reals whose last step is a quotient suffers from this deficiency.
There are three common ways out of the conundrum in constructive mathematics:
\index{mathematics!constructive}%
\index{bargaining}%
\begin{enumerate}
\item Pretend that the reals are a setoid $(\mathcal{C}, {\approx})$, i.e., the type of
  Cauchy sequences $\mathcal{C}$ with a coincidence\index{coincidence, of Cauchy approximations} relation attached to it by
  administrative decree. A sequence of reals then simply \emph{is} a sequence of Cauchy
  sequences representing them.
\item Give in to temptation and accept the axiom of countable choice. After all, the axiom
  is valid in most models of constructive mathematics based on a computational viewpoint,
  such as realizability models.
\item Declare the Cauchy reals unworthy and construct the Dedekind reals instead.
  Such a verdict is perfectly valid in certain contexts, such as in sheaf-theoretic models of constructive mathematics.
  However, as we saw in \autoref{sec:dedekind-reals}, the constructive Dedekind reals have their own problems.
\end{enumerate}

Using higher inductive types, however, there is a fourth solution, which we believe to be preferable to any of the above, and interesting even to a classical mathematician.
The idea is that the Cauchy real numbers should be the \emph{free complete metric space}\index{free!complete metric space} generated by~\Q.
In general, the construction of a free gadget of any sort requires applying the gadget operations repeatedly many times to the generators.
For instance, the elements of the free group on a set $X$ are not just binary products and inverses of elements of $X$, but words obtained by iterating the product and inverse constructions.
Thus, we might naturally expect the same to be true for Cauchy completion, with the relevant ``operation'' being ``take the limit of a Cauchy sequence''.
(In this case, the iteration would have to take place transfinitely, since even after infinitely many steps there will be new Cauchy sequences to take the limit of.)

The argument referred to above shows that if excluded middle or countable choice hold, then Cauchy completion is very special: when building the completion of a space, it suffices to stop applying the operation after \emph{one step}.
This may be regarded as analogous to the fact that free monoids and free groups can be given explicit descriptions in terms of (reduced) words.
However, we saw in \autoref{sec:free-algebras} that higher inductive types allow us to construct free gadgets \emph{directly}, whether or not there is also an explicit description available.
In this section we show that the same is true for the Cauchy reals (a similar technique would construct the Cauchy completion of any metric space; see \autoref{ex:metric-completion}).
Specifically, higher inductive types allow us to \emph{simultaneously} add limits of Cauchy sequences and quotient by the coincidence relation, so that we can avoid the problem of lifting a sequence of reals to a sequence of representatives.
\index{completion!Cauchy|)}%

\subsection{Construction of Cauchy reals}
\label{sec:constr-cauchy-reals}

The construction of the Cauchy reals $\RC$ as a higher inductive type is a bit more subtle than that of the free algebraic structures considered in \autoref{sec:free-algebras}.
We intend to include a ``take the limit'' constructor whose input is a Cauchy sequence of reals, but the notion of ``Cauchy sequence of reals'' depends on having some way to measure the ``distance'' between real numbers.
In general, of course, the distance between two real numbers will be another real number, leading to a potentially problematic circularity.

However, what we actually need for the notion of Cauchy sequence of reals is not the general notion of ``distance'', but a way to say that ``the distance\index{distance} between two real numbers is less than $\epsilon$'' for any $\epsilon:\Qp$.
This can be represented by a family of binary relations, which we will denote $\mathord{\close\epsilon} : \RC\to\RC\to \prop$.
The intended meaning of $x \close\epsilon y$ is $|x - y| < \epsilon$, but since we do not have notions of subtraction, absolute value, or inequality available yet (we are only just defining $\RC$, after all), we will have to define these relations $\close\epsilon$ at the same time as we define $\RC$ itself.
And since $\close\epsilon$ is a type family indexed by two copies of $\RC$, we cannot do this with an ordinary mutual (higher) inductive definition; instead we have to use a \emph{higher inductive-inductive definition}.
\index{inductive-inductive type!higher}

Recall from \autoref{sec:generalizations} that the ordinary notion of inductive-inductive definition allows us to define a type and a type family indexed by it by simultaneous induction.
Of course, the ``higher'' version of this allows both the type and the family to have path constructors as well as point constructors.
We will not attempt to formulate any general theory of higher inductive-inductive definitions, but hopefully the description we will give of $\RC$ and $\close\epsilon$ will make the idea transparent.

\begin{rmk}
  We might also consider a \emph{higher inductive-recursive definition}, in which $\close\epsilon$ is defined using the \emph{recursion} principle of $\RC$, simultaneously with the \emph{inductive} definition of $\RC$.
  We choose the inductive-inductive route instead for two reasons.
  Firstly, higher inductive-re\-cur\-sive definitions seem to be more difficult to justify in homotopical semantics.
  Secondly, and more importantly, the inductive-inductive definition yields a more powerful induction principle, which we will need in order to develop even the basic theory of Cauchy reals.
\end{rmk}

Finally, as we did for the discussion of Cauchy completeness of the Dedekind reals in \autoref{sec:RD-cauchy-complete}, we will work with \emph{Cauchy approximations} (\autoref{defn:cauchy-approximation}) instead of Cauchy sequences.
Of course, our Cauchy approximations will now consist of Cauchy reals, rather than Dedekind reals or rational numbers.

\begin{defn}\label{defn:cauchy-reals}
  Let $\RC$ and the relation $\closesym:\Qp \times \RC \times \RC \to \type$ be the following higher inductive-inductive type family.
  The type $\RC$ of \define{Cauchy reals}
  \indexdef{real numbers!Cauchy}%
  \indexsee{Cauchy!real numbers}{real numbers, Cau\-chy}%
  is generated by the following constructors:
  \begin{itemize}
  \item \emph{rational points:} 
    for any $q : \Q$ there is a real $\rcrat(q)$.
    \index{rational numbers!as Cauchy real numbers}%
  \item \emph{limit points}:
    for any $x : \Qp \to \RC$ such that
    \begin{equation}
      \label{eq:RC-cauchy}
      \fall{\delta, \epsilon : \Qp} x_\delta \close{\delta + \epsilon} x_\epsilon
    \end{equation}
    there is a point $\rclim(x) : \RC$. We call $x$ a \define{Cauchy approximation}.
    \indexdef{Cauchy!approximation}%
    \index{limit!of a Cauchy approximation}%
  \item \emph{paths:}
    for $u, v : \RC$ such that
    \begin{equation}
      \label{eq:RC-path}
      \fall{\epsilon : \Qp} u \close\epsilon v
    \end{equation}
    then there is a path $\rceq(u, v) : \id[\RC]{u}{v}$.
  \end{itemize}
  Simultaneously, the type family $\closesym:\RC\to\RC\to\Qp \to\type$ is generated by the following constructors.
  Here $q$ and $r$ denote rational numbers; $\delta$, $\epsilon$, and $\eta$ denote positive rationals; $u$ and $v$ denote Cauchy reals; and $x$ and $y$ denote Cauchy approximations:
  \begin{itemize}
  \item for any $q,r,\epsilon$, if $-\epsilon < q - r < \epsilon$, then $\rcrat(q) \close\epsilon \rcrat(r)$,
  \item for any $q,y,\epsilon,\delta$, if $\rcrat(q) \close{\epsilon - \delta} y_\delta$, then $\rcrat(q) \close{\epsilon} \rclim(y)$,
  \item for any $x,r,\epsilon,\delta$, if $x_\delta \close{\epsilon - \delta} \rcrat(r)$, then $\rclim(x) \close\epsilon \rcrat(r)$,
  \item for any $x,y,\epsilon,\delta,\eta$, if $x_\delta \close{\epsilon - \delta - \eta} y_\eta$, then $\rclim(x) \close\epsilon \rclim(y)$,
  \item for any $u,v,\epsilon$, if $\xi,\zeta : u \close{\epsilon} v$, then $\xi=\zeta$ (propositional truncation).
  \end{itemize}
\end{defn}

\mentalpause

The first constructor of $\RC$ says that any rational number can be regarded as a real number.
The second says that from any Cauchy approximation to a real number, we can obtain a new real number called its ``limit''.
And the third expresses the idea that if two Cauchy approximations coincide, then their limits are equal.

The first four constructors of $\closesym$ specify when two rational numbers are close, when a rational is close to a limit, and when two limits are close.
In the case of two rational numbers, this is just the usual notion of $\epsilon$-closeness for rational numbers, whereas the other cases can be derived by noting that each approximant $x_\delta$ is supposed to be within $\delta$ of the limit $\rclim(x)$.

We remind ourselves of proof-relevance: a real number obtained from $\rclim$ is represented not
just by a Cauchy approximation $x$, but also a proof $p$ of~\eqref{eq:RC-cauchy}, so we
should technically have written $\rclim(x,p)$ instead of just $\rclim(x)$.
A similar observation also applies to $\rceq$ and~\eqref{eq:RC-path}, but we shall write just
$\rceq : u = v$ instead of $\rceq(u, v, p) : u = v$. These abuses of notation are
mitigated by the fact that we are omitting mere propositions and information that is
readily guessed.
Likewise, the last constructor of $\mathord{\close\epsilon}$ justifies our leaving the other four nameless.

We are immediately able to populate $\RC$ with many real numbers. For suppose $x : \N \to
\Q$ is a traditional Cauchy sequence\index{Cauchy!sequence} of rational numbers, and let $M : \Qp \to \N$ be its
modulus of convergence. Then $\rcrat \circ x \circ M : \Qp \to \RC$ is a Cauchy
approximation, using the first constructor of $\closesym$ to produce the necessary witness.
Thus, $\rclim(\rcrat \circ x \circ m)$ is a real number. Various famous
real numbers $\sqrt{2}$, $\pi$, $e$, \dots{} are all limits of such Cauchy sequences of
rationals.

\subsection{Induction and recursion on Cauchy reals}
\label{sec:induct-recurs-cauchy}

In order to do anything useful with $\RC$, of course, we need to give its induction principle.
As is the case whenever we inductively define two or more objects at once, the basic induction principle for $\RC$ and $\closesym$ requires a simultaneous induction over both at once.
Thus, we should expect it to say that assuming two type families over $\RC$ and $\closesym$, respectively, together with data corresponding to each constructor, there exist sections of both of these families.
However, since $\closesym$ is indexed on two copies of $\RC$, the precise dependencies of these families is a bit subtle.
The induction principle will apply to any pair of type families:
\begin{align*}
A&:\RC\to\type\\
B&:\prd{x,y:\RC} A(x) \to A(y) \to \prd{\epsilon:\Qp} (x\close\epsilon y) \to \type.
\end{align*}
The type of $A$ is obvious, but the type of $B$ requires a little thought.
Since $B$ must depend on $\closesym$, but $\closesym$ in turn depends on two copies of $\RC$ and one copy of $\Qp$, it is fairly obvious that $B$ must also depend on the variables $x,y:\RC$ and $\epsilon:\Qp$ as well as an element of $(x\close\epsilon y)$.
What is slightly less obvious is that $B$ must also depend on $A(x)$ and $A(y)$.

This may be more evident if we consider the non-dependent case (the recursion principle), where $A$ is a simple type (rather than a type family).
In this case we would expect $B$ not to depend on $x,y:\RC$ or $x\close\epsilon y$.
But the recursion principle (along with its associated uniqueness principle) is supposed to say that $\RC$ with $\close\epsilon$ is an ``initial object'' in some category, so in this case the dependency structure of $A$ and $B$ should mirror that of $\RC$ and $\close\epsilon$: that is, we should have $B:A\to A\to \Qp \to \type$.
Combining this observation with the fact that, in the dependent case, $B$ must also depend on $x,y:\RC$ and $x\close\epsilon y$, leads inevitably to the type given above for $B$.

\symlabel{RC-recursion}
It is helpful to think of $B$ as an $\epsilon$-indexed family of relations between the types $A(x)$ and $A(y)$.
With this in mind, we may write $B(x,y,a,b,\epsilon,\xi)$ as $(x,a) \bsim_\epsilon^\xi (y,b)$.
Since $\xi:x\close\epsilon y$ is unique when it exists, we generally omit it from the notation and write $(x,a) \bsim_\epsilon (y,b)$; this is harmless as long as we keep in mind that this relation is only defined when $x\close\epsilon y$.
We may also sometimes simplify further and write $a\bsim_\epsilon b$, with $x$ and $y$ inferred from the types of $a$ and $b$, but sometimes it will be necessary to include them for clarity.

\index{induction principle!for Cauchy reals}%
Now, given a type family $A:\RC\to\type$ and a family of relations $\bsim$ as above, the hypotheses of the induction principle consist of the following data, one for each constructor of $\RC$ or $\closesym$:
\begin{itemize}
\item For any $q : \Q$, an element $f_q:A(\rcrat(q))$.
\item For any Cauchy approximation $x$, and any $a:\prd{\epsilon:\Qp} A(x_\epsilon)$ such that
  \begin{equation}
    \fall{\delta, \epsilon : \Qp}
    (x_\delta,a_\delta) \bsim_{\delta+\epsilon} (x_\epsilon,a_\epsilon),
    \label{eq:depCauchyappx}
  \end{equation}
  an element $f_{x,a}:A(\rclim(x))$.
  We call such $a$ a \define{dependent Cauchy approximation}
  \indexdef{Cauchy!approximation!dependent}%
  \indexsee{approximation, Cauchy}{Cauchy approximation}%
  \indexdef{dependent!Cauchy approximation}%
  over $x$.
\item For $u, v : \RC$ such that $h:\fall{\epsilon : \Qp} u \close\epsilon v$, and all $a:A(u)$ and $b:A(v)$ such that
  $\fall{\epsilon:\Qp} (u,a) \bsim_\epsilon (v,b)$,
  a dependent path $\dpath{A}{\rceq(u,v)}{a}{b}$.
\item For $q,r:\Q$ and $\epsilon:\Qp$, if $-\epsilon < q - r < \epsilon$, we have
  \narrowequation{(\rcrat(q),f_q) \bsim_\epsilon (\rcrat(r),f_r).}
\item For $q:\Q$ and $\delta,\epsilon:\Qp$ and $y$ a Cauchy approximation, and $b$ a dependent Cauchy approximation over $y$, if $\rcrat(q) \close{\epsilon - \delta} y_\delta$, then
  \[(\rcrat(q),f_q) \bsim_{\epsilon-\delta} (y_\delta,b_\delta)
  \;\Rightarrow\;
  (\rcrat(q),f_q) \bsim_\epsilon (\rclim(y),f_{y,b}).\]
\item Similarly, for $r:\Q$ and $\delta,\epsilon:\Qp$ and $x$ a Cauchy approximation, and $a$ a dependent Cauchy approximation over $x$, if $x_\delta \close{\epsilon - \delta} \rcrat(r)$, then
  \[(x_\delta,a_\delta) \bsim_{\epsilon-\delta} (\rcrat(r),f_r)
  \;\Rightarrow\;
  (\rclim(x),f_{x,a}) \bsim_\epsilon (\rcrat(q),f_r).
  \]
\item For $\epsilon,\delta,\eta:\Qp$ and $x,y$ Cauchy approximations, and $a$ and $b$ dependent Cauchy approximations over $x$ and $y$ respectively, if we have $x_\delta \close{\epsilon - \delta - \eta} y_\eta$, then
  \[ (x_\delta,a_\delta) \bsim_{\epsilon - \delta - \eta} (y_\eta,b_\eta)
  \;\Rightarrow\;
  (\rclim(x),f_{x,a}) \bsim_\epsilon (\rclim(y),f_{y,b}).\]
\item For $\epsilon:\Qp$ and $x,y:\RC$ and $\xi,\zeta:x\close{\epsilon} y$, and $a:A(x)$ and $b:A(y)$, any two elements of $(x,a) \bsim_\epsilon^\xi (y,b)$ and $(x,a) \bsim_\epsilon^\zeta (y,b)$ are dependently equal over $\xi=\zeta$.
  Note that as usual, this is equivalent to asking that $\bsim$ takes values in mere propositions.
\end{itemize}
Under these hypotheses, we deduce functions
\begin{align*}
  f&:\prd{x:\RC} A(x)\\
  g&:\prd{x,y:\RC}{\epsilon:\Qp}{\xi:x\close{\epsilon} y}
  (x,f(x)) \bsim_\epsilon^\xi (y,f(y))
\end{align*}
which compute as expected:
\begin{align}
  f(\rcrat(q)) &\defeq f_q, \label{eq:rcsimind1}\\
  f(\rclim(x)) &\defeq f_{x,(f,g)[x]}. \label{eq:rcsimind2}
\end{align}
Here $(f,g)[x]$ denotes the result of applying $f$ and $g$ to a Cauchy approximation $x$ to obtain a dependent Cauchy approximation over $x$.
That is, we define $(f,g)[x]_\epsilon \defeq f(x_\epsilon) : A(x_\epsilon)$, and then for any $\epsilon,\delta:\Qp$ we have $g(x_\epsilon,x_\delta,\epsilon+\delta,\xi)$ to witness the fact that $(f,g)[x]$ is a dependent Cauchy approximation, where $\xi: x_\epsilon \close{\epsilon+\delta} x_\delta$ arises from the assumption that $x$ is a Cauchy approximation.

We will never use this notation again, so don't worry about remembering it.
Generally we use the pattern-matching convention, where $f$ is defined by equations such as~\eqref{eq:rcsimind1} and~\eqref{eq:rcsimind2} in which the right-hand side of~\eqref{eq:rcsimind2} may involve the symbols $f(x_\epsilon)$ and an assumption that they form a dependent Cauchy approximation.

However, this induction principle is admittedly still quite a mouthful.
To help make sense of it, we observe that it contains as special cases two separate induction principles for~$\RC$ and for~$\closesym$.
Firstly, suppose given only a type family $A:\RC\to\type$, and define $\bsim$ to be constant at \unit.
Then much of the required data becomes trivial, and we are left with:
\begin{itemize}
\item for any $q : \Q$, an element $f_q:A(\rcrat(q))$,
\item for any Cauchy approximation $x$, and any $a:\prd{\epsilon:\Qp} A(x_\epsilon)$, an element $f_{x,a}:A(\rclim(x))$,
\item for $u, v : \RC$ and $h:\fall{\epsilon : \Qp} u \close\epsilon v$, and $a:A(u)$ and $b:A(v)$, we have $\dpath{A}{\rceq(u,v)}{a}{b}$.
\end{itemize}
Given these data, the induction principle yields a function $f:\prd{x:\RC} A(x)$ such that
\begin{align*}
  f(\rcrat(q)) &\defeq f_q,\\
  f(\rclim(x)) &\defeq f_{x,f(x)}.
\end{align*}
We call this principle \define{$\RC$-induction}; it says essentially that if we take $\close\epsilon$ as given, then $\RC$ is inductively generated by its constructors.

In particular, if $A$ is a mere property, the third hypothesis in $\RC$-induction is trivial.
Thus, we may prove mere properties of real numbers by simply proving them for rationals and for limits of Cauchy approximations.
Here is an example.

\begin{lem}
  For any $u:\RC$ and $\epsilon:\Qp$, we have $u\close\epsilon u$.
\end{lem}
\begin{proof}
  Define $A(u) \defeq \fall{\epsilon:\Qp} (u\close\epsilon u)$.
  Since this is a mere proposition (by the last constructor of $\closesym$), by $\RC$-induction, it suffices to prove it when $u$ is $\rcrat(q)$ and when $u$ is $\rclim(x)$.
  In the first case, we obviously have $|q-q|<\epsilon$ for any $\epsilon$, hence $\rcrat(q) \close\epsilon \rcrat(q)$ by the first constructor of $\closesym$.
  And in the second case, we may assume inductively that $x_\delta \close\epsilon x_\delta$ for all $\delta,\epsilon:\Qp$.
  Then in particular, we have $x_{\epsilon/3} \close{\epsilon/3} x_{\epsilon/3}$, whence $\rclim(x) \close{\epsilon} \rclim(x)$ by the fourth constructor of $\closesym$.
\end{proof}

\begin{thm}\label{thm:Cauchy-reals-are-a-set}
  $\RC$ is a set.
\end{thm}
\begin{proof}
  We have just shown that the mere relation
  \narrowequation{P(u,v) \defeq \fall{\epsilon:\Qp} (u\close\epsilon v)}
  is reflexive.
  Since it implies identity, by the path constructor of $\RC$, the result follows from \autoref{thm:h-set-refrel-in-paths-sets}.
\end{proof}

We can also show that although $\RC$ may not be a quotient of the set of Cauchy sequences of \emph{rationals}, it is nevertheless a quotient of the set of Cauchy sequences of \emph{reals}.
(Of course, this is not a valid \emph{definition} of $\RC$, but it is a useful property.)
We define the type of Cauchy approximations to be
\symlabel{cauchy-approximations}%
\index{Cauchy!approximation!type of}%
\begin{equation*}
  \CAP \defeq
  \setof{ x : \Qp \to \RC |
    \fall{\epsilon, \delta : \Qp} x_\delta \close{\delta + \epsilon} x_\epsilon
  }.
\end{equation*}
The second constructor of $\RC$ gives a function $\rclim:\CAP\to\RC$.

\begin{lem} \label{RC-lim-onto}
  Every real merely is a limit point: $\fall{u : \RC} \exis{x : \CAP} u = \rclim(x)$.
  In other words, $\rclim:\CAP\to\RC$ is surjective.
\end{lem}
\begin{proof}
  By $\RC$-induction, we may divide into cases on $u$.
  Of course, if $u$ is a limit $\rclim(x)$, the statement is trivial.
  So suppose $u$ is a rational point $\rcrat(q)$; we claim $u$ is equal to $\rclim(\lam{\epsilon} \rcrat(q))$.
  By the path constructor of $\RC$, it suffices to show $\rcrat(q) \close\epsilon \rclim(\lam{\epsilon} \rcrat(q))$ for all $\epsilon:\Qp$.
  And by the second constructor of $\closesym$, for this it suffices to find $\delta:\Qp$ such that $\rcrat(q)\close{\epsilon-\delta} \rcrat(q)$.
  But by the first constructor of $\closesym$, we may take any $\delta:\Qp$ with $\delta<\epsilon$.
\end{proof}

\begin{lem} \label{RC-lim-factor}
  If $A$ is a set and $f : \CAP \to A$ respects coincidence\index{coincidence!of Cauchy approximations} of Cauchy approximations, in the sense that
  \begin{equation*}
    \fall{x, y : \CAP} \rclim(x) = \rclim(y) \Rightarrow f(x) = f(y),
  \end{equation*}
  then $f$ factors uniquely through $\rclim : \CAP \to \RC$.
\end{lem}
\begin{proof}
  Since $\rclim$ is surjective, by \autoref{lem:images_are_coequalizers}, $\RC$ is the quotient of $\CAP$ by the kernel pair\index{kernel!pair} of $\rclim$.
  But this is exactly the statement of the lemma.
\end{proof}

For the second special case of the induction principle, suppose instead that we take $A$ to be constant at $\unit$.
In this case, $\bsim$ is simply an $\epsilon$-indexed family of relations on $\epsilon$-close pairs of real numbers, so we may write $u\bsim_\epsilon v$ instead of $(u,\ttt)\bsim_\epsilon (v,\ttt)$.
Then the required data reduces to the following, where $q, r$ denote rational numbers, $\epsilon, \delta, \eta$ positive rational numbers, and $x, y$ Cauchy approximations:
\begin{itemize}
\item if $-\epsilon < q - r < \epsilon$, then
  $\rcrat(q) \bsim_\epsilon \rcrat(r)$,
\item if $\rcrat(q) \close{\epsilon - \delta} y_\delta$ and
  $\rcrat(q)\bsim_{\epsilon-\delta} y_\delta$,
  then $\rcrat(q) \bsim_\epsilon \rclim(y)$,
\item if $x_\delta \close{\epsilon - \delta} \rcrat(r)$ and
  $x_\delta \bsim_{\epsilon-\delta} \rcrat(r)$,
  then $\rclim(y) \bsim_\epsilon \rcrat(q)$,
\item if $x_\delta \close{\epsilon - \delta - \eta} y_\eta$ and
  $x_\delta\bsim_{\epsilon - \delta - \eta} y_\eta$,
  then $\rclim(x) \bsim_\epsilon \rclim(y)$.
\end{itemize}
The resulting conclusion is $\fall{u,v:\RC}{\epsilon:\Qp} (u\close\epsilon v) \to (u \bsim_\epsilon v)$.
We call this principle \define{$\closesym$-induction}; it says essentially that if we take $\RC$ as given, then $\close\epsilon$ is inductively generated (as a family of types) by \emph{its} constructors.
For example, we can use this to show that $\closesym$ is symmetric.

\begin{lem}\label{thm:RCsim-symmetric}
  For any $u,v:\RC$ and $\epsilon:\Qp$, we have $(u\close\epsilon v) = (v\close\epsilon u)$.
\end{lem}
\begin{proof}
  Since both are mere propositions, by symmetry it suffices to show one implication.
  Thus, let $(u\bsim_\epsilon v) \defeq (v\close\epsilon u)$.
  By $\closesym$-induction, we may reduce to the case that $u\close\epsilon v$ is derived from one of the four interesting constructors of $\closesym$.
  In the first case when $u$ and $v$ are both rational, the result is trivial (we can apply the first constructor again).
  In the other three cases, the inductive hypothesis (together with commutativity of addition in $\Q$) yields exactly the input to another of the constructors of $\closesym$ (the second and third constructors switch, while the fourth stays put).
\end{proof}

The general induction principle, which we may call \define{$(\RC,\closesym)$-induction}, is therefore a sort of joint $\RC$-induction and $\closesym$-induction.
Consider, for instance, its non-dependent version, which we call \define{$(\RC,\closesym)$-recursion}, which is the one that we will have the most use for.
\index{recursion principle!for Cauchy reals}%
Ordinary $\RC$-recursion tells us that to define a function $f : \RC \to A$ it suffices to:
\begin{enumerate}
\item for every $q : \Q$ construct $f(\rcrat(q)) : A$,
\item for every Cauchy approximation $x : \Qp \to \RC$, construct $f(x) : A$,
  assuming that $f(x_\epsilon)$ has already been defined for all $\epsilon : \Qp$,
\item prove $f(u) = f(v)$ for all $u, v : \RC$ satisfying $\fall{\epsilon:\Qp} u\close\epsilon v$.\label{item:rcrec3}
\end{enumerate}
However, it is generally quite difficult to show~\ref{item:rcrec3} without knowing something about how $f$ acts on $\epsilon$-close Cauchy reals.
The enhanced principle of $(\RC,\closesym)$-recursion remedies this deficiency, allowing us to specify an \emph{arbitrary} ``way in which $f$ acts on $\epsilon$-close Cauchy reals'', which we can then prove to be the case by a simultaneous induction with the definition of $f$.
This is the family of relations $\bsim$.
Since $A$ is independent of $\RC$, we may assume for simplicity that $\bsim$ depends only on $A$ and $\Qp$, and thus there is no ambiguity in writing $a\bsim_\epsilon b$ instead of $(u,a) \bsim_\epsilon (v,b)$.
In this case, defining a function $f:\RC\to A$ by $(\RC,\closesym)$-recursion requires the following cases (which we now write using the pattern-matching convention).
\begin{itemize}
\item For every $q : \Q$, construct $f(\rcrat(q)) : A$.
\item For every Cauchy approximation $x : \Qp \to \RC$, construct $f(x) : A$, assuming inductively that $f(x_\epsilon)$ has already been defined for all $\epsilon : \Qp$ and form a ``Cauchy approximation with respect to $\bsim$'', i.e.\ that $\fall{\epsilon,\delta:\Qp} (f(x_\epsilon) \bsim_{\epsilon+\delta} f(x_\delta))$.
\item Prove that the relations $\bsim$ are \emph{separated}, i.e.\ that, for any $a,b:A$,
  \indexdef{relation!separated family of}%
  \indexdef{separated family of relations}%
\narrowequation{(\fall{\epsilon:\Qp} a\bsim_\epsilon b) \Rightarrow (a=b).}
\item Prove that if $-\epsilon< q-r <\epsilon$ for $q,r:\Q$, then $f(\rcrat(q))\bsim_\epsilon f(\rcrat(r))$.
\item For any $q:\Q$ and any Cauchy approximation $y$, prove that
\narrowequation{f(\rcrat(q)) \bsim_\epsilon f(\rclim(y)),} assuming inductively that $\rcrat(q)\close{\epsilon-\delta} y_\delta$ and $f(\rcrat(q)) \bsim_{\epsilon-\delta} f(y_\delta)$ for some $\delta:\Qp$, and that $\eta \mapsto f(x_\eta)$ is a Cauchy approximation with respect to $\bsim$.
\item For any Cauchy approximation $x$ and any $r:\Q$, prove that
\narrowequation{f(\rclim(x)) \bsim_\epsilon f(\rcrat(r)),}
assuming inductively that $x_\delta \close{\epsilon-\delta} \rcrat(r)$ and $f(x_\delta) \bsim_{\epsilon-\delta} f(\rcrat(r))$ for some $\delta:\Qp$, and that $\eta\mapsto f(x_\eta)$ is a Cauchy approximation with respect to $\bsim$.
\item For any Cauchy approximations $x,y$, prove that
\narrowequation{f(\rclim(x)) \bsim_\epsilon f(\rclim(y)),}
assuming inductively that $x_\delta \close{\epsilon-\delta-\eta} y_\eta$ and $f(x_\delta) \bsim_{\epsilon-\delta-\eta} f(y_\eta)$ for some $\delta,\eta:\Qp$, and that $\theta\mapsto f(x_\theta)$ and $\theta\mapsto f(y_\theta)$ are Cauchy approximations with respect to $\bsim$.
\end{itemize}
Note that in the last four proofs, we are free to use the specific definitions of $f(\rcrat(q))$ and $f(\rclim(x))$ given in the first two data.
However, the proof of separatedness must apply to \emph{any} two elements of $A$, without any relation to $f$: it is a sort of ``admissibility'' condition on the family of relations $\bsim$.
Thus, we often verify it first, immediately after defining $\bsim$, before going on to define $f(\rcrat(q))$ and $f(\rclim(x))$.

Under the above hypotheses, $(\RC,\closesym)$-recursion yields a function $f:\RC\to A$ such that $f(\rcrat(q))$ and $f(\rclim(x))$ are judgmentally equal to the definitions given for them in the first two clauses.
Moreover, we may also conclude
\begin{equation}
  \fall{u,v:\RC}{\epsilon:\Qp} (u\close\epsilon v) \to (f(u) \bsim_\epsilon f(v)).\label{eq:RC-sim-recursion-extra}
\end{equation}

As a paradigmatic example, $(\RC,\closesym)$-recursion allows us to extend functions defined on $\Q$ to all of $\RC$, as long as they are sufficiently continuous.
\index{function!continuous}%

\begin{defn}\label{defn:lipschitz}
  A function $f:\Q\to\RC$ is \define{Lipschitz}
  \indexdef{function!Lipschitz}%
  \indexdef{Lipschitz!function}%
  \indexdef{Lipschitz!constant}%
  \indexdef{constant!Lipschitz}%
  if there exists $L:\Qp$ (the \define{Lipschitz constant}) such that
  \[ |q - r|<\epsilon \Rightarrow (f(q) \close{L\epsilon} f(r)) \]
  for all $\epsilon:\Qp$ and $q,r:\Q$.
  Similarly, $g:\RC\to\RC$ is \define{Lipschitz} if there exists $L:\Qp$ such that
  \[ (u\close\epsilon v) \Rightarrow (g(u) \close{L\epsilon} g(v)) \]
  for all $\epsilon:\Qp$ and $u,v:\RC$..
\end{defn}

In particular, note that by the first constructor of $\closesym$, if $f:\Q\to\Q$ is Lipschitz in the obvious sense, then so is the composite $\Q\xrightarrow{f} \Q \to \RC$.

\begin{lem}\label{RC-extend-Q-Lipschitz}
  Suppose $f : \Q \to \RC$ is Lipschitz with constant $L : \Qp$.
  Then there exists a Lipschitz map $\bar{f} : \RC \to \RC$, also with constant $L$, such that $\bar{f}(\rcrat(q)) \jdeq f(q)$ for all $q:\Q$.
\end{lem}

\begin{proof}
  We define $\bar{f}$ by $(\RC,\closesym)$-recursion, with codomain $A\defeq \RC$.
  We define the relation $\mathord{\bsim}: \RC \to \RC \to \Qp \to \prop$ to be
  \begin{align*}
    (u \bsim_\epsilon v) &\defeq (u \close{L\epsilon} v).
  \end{align*}
  For $q : \Q$, we define
  \begin{equation*}
    \bar{f}(\rcrat(q)) \defeq \rcrat(f(q)).
  \end{equation*}
  For a Cauchy approximation $x : \Qp \to \RC$, we define
  \begin{equation*}
    \bar{f}(\rclim(x)) \defeq \rclim(\lamu{\epsilon : \Qp} \bar{f}(x_{\epsilon/L})).
  \end{equation*}
  For this to make sense, we must verify that $y \defeq \lamu{\epsilon : \Qp} \bar{f}(x_{\epsilon/L})$ is a Cauchy approximation.
  However, the inductive hypothesis for this step is that for any $\delta,\epsilon:\Qp$ we have $\bar{f}(x_\delta) \bsim_{\delta+\epsilon} \bar{f}(x_\epsilon)$, i.e.\ $\bar{f}(x_\delta) \close{L\delta+L\epsilon} \bar{f}(x_\epsilon)$.
  Thus we have
  \[y_\delta \jdeq f(x_{\delta/L}) \close{\delta + \epsilon} f(x_{\epsilon/L})   \jdeq y_\epsilon. \]
  
  For proving separatedness, we simply observe that $\fall{\epsilon:\Qp} a\bsim_\epsilon b$ means $\fall{\epsilon:\Qp} a\close{L\epsilon} b$, which implies $\fall{\epsilon:\Qp}a\close\epsilon b$ and thus $a=b$.

  To complete the $(\RC,\closesym)$-recursion, it remains to verify the four conditions on $\bsim$.
  This basically amounts to proving that $\bar f$ is Lipschitz for all the four constructors of $\closesym$.
  \begin{enumerate}
  \item When $u$ is $\rcrat(q)$ and $v$ is $\rcrat(r)$ with $-\epsilon < |q-r| <\epsilon$, the assumption that $f$ is Lipschitz yields $f(q) \close{L\epsilon} f(r)$, hence $\bar{f}(\rcrat(q)) \bsim_\epsilon \bar{f}(\rcrat(r))$ by definition.
  \item When $u$ is $\rclim(x)$ and $v$ is $\rcrat(q)$ with $x_\eta \close{\epsilon - \eta} \rcrat(q)$, then the
      inductive hypothesis is $\bar{f}(x_\eta) \close{L \epsilon - L \eta} \rcrat(f(q))$, which proves
      \narrowequation{\bar{f}(\rclim(x)) \close{L \epsilon} \bar{f}(\rcrat(q))}
      by the third constructor of $\closesym$.
  \item The symmetric case when $u$ is rational and $v$ is a limit is essentially identical.
  \item When $u$ is $\rclim(x)$ and $v$ is $\rclim(y)$, with $\delta, \eta : \Qp$ such that $x_\delta \close{\epsilon - \delta - \eta} y_\eta$,
      the inductive hypothesis is $\bar{f}(x_\delta) \close{L \epsilon - L \delta - L \eta} \bar{f}(y_\eta)$, which proves $\bar{f}(\rclim(x)) \close{L
        \epsilon} \bar{f}(\rclim(y))$ by the fourth constructor of $\closesym$.
  \end{enumerate}
  This completes the $(\RC,\closesym)$-recursion, and hence the construction of $\bar f$.
  The desired equality $\bar f(\rcrat(q))\jdeq f(q)$ is exactly the first computation rule for $(\RC,\closesym)$-recursion, and the additional condition~\eqref{eq:RC-sim-recursion-extra} says exactly that $\bar f$ is Lipschitz with constant $L$.
\end{proof}

At this point we have gone about as far as we can without a better characterization of $\closesym$.
We have specified, in the constructors of $\closesym$, the conditions under which we want Cauchy reals of the two different forms to be $\epsilon$-close.
However, how do we know that in the resulting inductive-inductive type family, these are the \emph{only} witnesses to this fact?
We have seen that inductive type families (such as identity types, see \autoref{sec:identity-systems}) and higher inductive types have a tendency to contain ``more than was put into them'', so this is not an idle question.

In order to characterize $\closesym$ more precisely, we will define a family of relations $\approx_\epsilon$ on $\RC$ \emph{recursively}, so that they will compute on constructors, and prove that this family is equivalent to $\close\epsilon$.

\begin{thm}\label{defn:RC-approx}
  There is a family of mere relations $\mathord\approx:\RC\to\RC\to\Qp\to\prop$ such that
  \begin{align}
    (\rcrat(q) \approx_\epsilon \rcrat(r))  &\defeq
    (-\epsilon < q - r < \epsilon)\label{eq:RCappx1}\\
    (\rcrat(q) \approx_\epsilon \rclim(y)) &\defeq
    \exis{\delta : \Qp} \rcrat(q) \approx_{\epsilon - \delta} y_\delta\label{eq:RCappx2}\\
    (\rclim(x) \approx_\epsilon \rcrat(r)) &\defeq
    \exis{\delta : \Qp} x_\delta \approx_{\epsilon - \delta} \rcrat(r)\label{eq:RCappx3}\\
    (\rclim(x) \approx_\epsilon \rclim(y)) &\defeq
    \exis{\delta, \eta : \Qp} x_\delta \approx_{\epsilon - \delta - \eta} y_\eta.\label{eq:RCappx4}
  \end{align}
  Moreover, we have
  \begin{gather}
    (u \approx_\epsilon v) \Leftrightarrow \exis{\theta:\Qp} (u \approx_{\epsilon-\theta} v) \label{RC-sim-rounded}\\
    (u \approx_\epsilon v) \to (v\close\delta w) \to (u\approx_{\epsilon+\delta} w)\label{eq:RC-sim-rtri}\\ 
    (u \close\epsilon v) \to (v\approx_\delta w) \to (u\approx_{\epsilon+\delta} w)\label{eq:RC-sim-ltri}.
  \end{gather}
\end{thm}

The additional conditions~\eqref{RC-sim-rounded}--\eqref{eq:RC-sim-ltri} turn out to be required in order to make the inductive definition go through.
Condition~\eqref{RC-sim-rounded} is called being \define{rounded}.
\indexsee{relation!rounded}{rounded relation}%
\indexdef{rounded!relation}%
Reading it from right to left gives \define{monotonicity} of $\approx$,
\index{monotonicity}%
\index{relation!monotonic}%
\begin{equation*}
  (\delta < \epsilon) \land (u \approx_\delta v) \Rightarrow (u \approx_\epsilon v)
\end{equation*}
while reading it left to right to \define{openness} of $\approx$,
\index{open!relation}%
\index{relation!open}%
\begin{equation*}
  (u \approx_\epsilon v) \Rightarrow \exis{\epsilon : \Qp} (\delta < \epsilon) \land (u \approx_\delta v).
\end{equation*}
Conditions~\eqref{eq:RC-sim-rtri} and~\eqref{eq:RC-sim-ltri} are forms of the triangle inequality, which say that $\approx$ is a ``module'' over $\closesym$ on both sides.

\begin{proof}
  We will define $\mathord\approx:\RC\to\RC\to\Qp\to\prop$ by double $(\RC,\closesym)$-recursion.
  First we will apply $(\RC,\closesym)$-recursion with codomain the subset of $\RC\to\Qp\to\prop$ consisting of those families of predicates which are rounded and satisfy the one appropriate form of the triangle inequality.
  Thinking of these predicates as half of a binary relation, we will write them as $(u,\epsilon) \mapsto (\hapx_\epsilon u)$, with the symbol $\hapname$ referring to the whole relation.
  Now we can write $A$ precisely as
  \begin{multline*}
    A \defeq\; \Bigg\{ \hapname :\RC\to\Qp\to\prop \;\bigg|\; \\
    \Big(\fall{u:\RC}{\epsilon:\Qp}
    \big((\hapx_\epsilon u) \Leftrightarrow \exis{\theta:\Qp} (\hapx_{\epsilon-\theta} u)\big)\Big)  \\
    \land \Big(\fall{u,v:\RC}{\eta,\epsilon:\Qp} (u\close\epsilon v) \to\\
    \big((\hapx_\eta u) \to (\hapx_{\eta+\epsilon} v) \big) \land \big((\hapx_\eta v) \to (\hapx_{\eta+\epsilon} u) \big)\Big)\Bigg\}
  \end{multline*}
  As usual with subsets, we will use the same notation for an inhabitant of $A$ and its first component $\hapname$.
  As the family of relations required for $(\RC,\closesym)$-recursion, we consider the following, which will ensure the other form of the triangle inequality:
  \begin{narrowmultline*}
    (\hapname \bsim_\epsilon \hapbname ) \defeq \narrowbreak
    \fall{u:\RC}{\eta:\Qp} ((\hapx_\eta u) \to (\hapxb_{\epsilon+\eta} u))
    \land \narrowbreak
    ((\hapxb_\eta u) \to (\hapx_{\epsilon+\eta} u)).
  \end{narrowmultline*}
  We observe that these relations are separated.
  For assuming
  \narrowequation{\fall{\epsilon:\Qp} (\hapname \bsim_\epsilon \hapbname),}
  to show $\hapname = \hapbname$ it suffices to show $(\hapx_\epsilon u) \Leftrightarrow (\hapxb_\epsilon u)$ for all $u:\RC$.
  But $\hapx_\epsilon u$ implies $\hapx_{\epsilon-\theta} u$ for some $\theta$, by roundedness, which together with $\hapname \bsim_\epsilon \hapbname$ implies $\hapxb_\epsilon u$; and the converse is identical.

  Now the first two data the recursion principle requires are the following.
  \begin{itemize}
  \item For any $q:\Q$, we must give an element of $A$, which we denote $(\rcrat(q)\approx_{(\blank)} \blank)$.
  \item For any Cauchy approximation $x$, if we assume defined a function $\Qp \to A$, which we will denote by $\epsilon \mapsto (x_\epsilon \approx_{(\blank)} \blank)$, with the property that 
    \begin{equation}
      \fall{u:\RC}{\delta,\epsilon,\eta:\Qp} (x_\delta \approx_\eta u) \to (x_\epsilon \approx_{\eta+\delta+\epsilon} u),\label{eq:appxrec2}
    \end{equation}
    we must give an element of $A$, which we write as $(\rclim(x)\approx_{(\blank)} \blank)$.
  \end{itemize}
  In both cases, we give the required definition by using a nested $(\RC,\closesym)$-recursion, with codomain the subset of $\Qp\to\prop$ consisting of rounded families of mere propositions.
  Thinking of these propositions as zero halves of a binary relation, we will write them as $\epsilon \mapsto (\tap{\epsilon})$, with the symbol $\tapname$ referring to the whole family.
  Now we can write the codomain of these inner recursions precisely as
  \begin{narrowmultline*}
    C \defeq
    \bigg\{ \tapname :\Qp\to\prop \;\;\Big|\;\; \narrowbreak
    \fall{\epsilon:\Qp} \Big((\tap\epsilon) \Leftrightarrow \exis{\theta:\Qp} (\tap{\epsilon-\theta})\Big)\bigg\}
  \end{narrowmultline*}
  We take the required family of relations to be the remnant of the triangle inequality:
  \begin{narrowmultline*}
    (\tapname \bbsim_\epsilon \tapbname) \defeq
    \fall{\eta:\Qp} ((\tap\eta) \to (\tapb{\epsilon+\eta})) \land
    \narrowbreak
    ((\tapb\eta) \to (\tap{\epsilon+\eta})).
  \end{narrowmultline*}
  These relations are separated by the same argument as for $\bsim$, using roundedness of all elements of $C$.

  Note that if such an inner recursion succeeds, it will yield a family of predicates $\hapname : \RC\to\Qp\to \prop$ which are rounded
(since their image in $\Qp\to\prop$ lies in $C$) and satisfy
  \[ \fall{u,v:\RC}{\epsilon:\Qp} (u\close\epsilon v) \to \big((\hapx_{(\blank)} u) \bbsim_\epsilon (\hapx_{(\blank)} u)\big). \]
  Expanding out the definition of $\bbsim$, this yields precisely the third condition for $\hapname$ to belong to $A$; thus it is exactly what we need.

  It is at this point that we can give the definitions~\eqref{eq:RCappx1}--\eqref{eq:RCappx4}, as the first two clauses of each of the two inner recursions, corresponding to rational points and limits.
  In each case, we must verify that the relation is rounded and hence lies in $C$.
  In the rational-rational case~\eqref{eq:RCappx1} this is clear, while in the other cases it follows from an inductive hypothesis.
  (In~\eqref{eq:RCappx2} the relevant inductive hypothesis is that $(\rcrat(q) \approx_{(\blank)} y_\delta) : C$, while in~\eqref{eq:RCappx3} and~\eqref{eq:RCappx4} it is that $(x_\delta \approx_{(\blank)} \blank) : A$.)

  The remaining data of the sub-recursions consist of showing that \eqref{eq:RCappx1}--\eqref{eq:RCappx4} satisfy the triangle inequality on the right with respect to the constructors of $\closesym$.
  There are eight cases --- four in each sub-recursion --- corresponding to the eight possible ways that $u$, $v$, and $w$ in~\eqref{eq:RC-sim-rtri} can be chosen to be rational points or limits.
  First we consider the cases when $u$ is $\rcrat(q)$.
  \begin{enumerate}
  \item Assuming $\rcrat(q)\approx_\phi \rcrat(r)$ and $-\epsilon<|r-s|<\epsilon$, we must show $\rcrat(q)\approx_{\phi+\epsilon} \rcrat(s)$.
    But by definition of $\approx$, this reduces to the triangle inequality for rational numbers.
  \item We assume $\phi,\epsilon,\delta:\Qp$ such that $\rcrat(q)\approx_\phi \rcrat(r)$ and $\rcrat(r) \close{\epsilon-\delta} y_\delta$, and inductively that
    \begin{equation}
      \fall{\psi:\Qp}(\rcrat(q) \approx_{\psi} \rcrat(r)) \to (\rcrat(q) \approx_{\psi+\epsilon-\delta} y_\delta).\label{eq:RCappx-rtri-rrl1}
    \end{equation}
    We assume also that $\psi,\delta\mapsto (\rcrat(q) \approx_{\psi} y_\delta)$ is a Cauchy approximation with respect to $\bbsim$, i.e.\
    \begin{equation}
      \fall{\psi,\xi,\zeta:\Qp} (\rcrat(q) \approx_{\psi} y_\xi) \to (\rcrat(q) \approx_{\psi+\xi+\zeta} y_\zeta),\label{eq:RCappx-rtri-rrl2}
    \end{equation}
    although we do not need this assumption in this case.
    Indeed, \eqref{eq:RCappx-rtri-rrl1} with $\psi\defeq \phi$ yields immediately $\rcrat(q) \approx_{\phi+\epsilon-\delta} y_\delta$, and hence $\rcrat(q) \approx_{\phi+\epsilon} \rclim(y)$ by definition of $\approx$.
  \item We assume $\phi,\epsilon,\delta:\Qp$ such that $\rcrat(q)\approx_\phi \rclim(y)$ and $y_\delta \close{\epsilon-\delta} \rcrat(r)$, and inductively that
    \begin{gather}
      \fall{\psi:\Qp}(\rcrat(q) \approx_{\psi} y_\delta) \to (\rcrat(q) \approx_{\psi+\epsilon-\delta} \rcrat(r)).\label{eq:RCappx-rtri-rlr1}\\
      \fall{\psi,\xi,\zeta:\Qp} (\rcrat(q) \approx_{\psi} y_\xi) \to (\rcrat(q) \approx_{\psi+\xi+\zeta} y_\zeta).\label{eq:RCappx-rtri-rlr2}
    \end{gather}
    By definition, $\rcrat(q)\approx_\phi \rclim(y)$ means that we have $\theta:\Qp$ with $\rcrat(q) \approx_{\phi-\theta} y_\theta$.
    By assumption~\eqref{eq:RCappx-rtri-rlr2}, therefore, we have also $\rcrat(q) \approx_{\phi+\delta} y_\delta$, and then by~\eqref{eq:RCappx-rtri-rlr1} it follows that $\rcrat(q) \approx_{\phi+\epsilon} \rcrat(r)$, as desired.
  \item We assume $\phi,\epsilon,\delta,\eta:\Qp$ such that $\rcrat(q)\approx_\phi \rclim(y)$ and $y_\delta \close{\epsilon-\delta-\eta} z_\eta$, and inductively that 
    \begin{gather}
      \fall{\psi:\Qp}(\rcrat(q) \approx_{\psi} y_\delta) \to (\rcrat(q) \approx_{\psi+\epsilon-\delta-\eta} z_\eta), \label{eq:RCappx-rtri-rll1}\\
      \fall{\psi,\xi,\zeta:\Qp} (\rcrat(q) \approx_{\psi} y_\xi) \to (\rcrat(q) \approx_{\psi+\xi+\zeta} y_\zeta), \label{eq:RCappx-rtri-rll2}\\
      \fall{\psi,\xi,\zeta:\Qp} (\rcrat(q) \approx_{\psi} z_\xi) \to (\rcrat(q) \approx_{\psi+\xi+\zeta} z_\zeta). \label{eq:RCappx-rtri-rll3}
    \end{gather}
    Again, $\rcrat(q)\approx_\phi \rclim(y)$ means we have $\xi:\Qp$ with $\rcrat(q) \approx_{\phi-\xi} y_\xi$, while~\eqref{eq:RCappx-rtri-rll2} then implies $\rcrat(q) \approx_{\phi+\delta} y_\delta$ and~\eqref{eq:RCappx-rtri-rll1} implies $\rcrat(q) \approx_{\phi+\epsilon-\eta} z_\eta$.
    But by definition of $\approx$, this implies $\rcrat(q) \approx_{\phi+\epsilon} \rclim(z)$ as desired.
  \end{enumerate}
  Now we move on to the cases when $u$ is $\rclim(x)$, with $x$ a Cauchy approximation.
  In this case, the ambient inductive hypothesis of the definition of $(\rclim(x) \approx_{(\blank)} {\blank}) : A$ is that we have ${(x_\delta \approx_{(\blank)} {\blank})}: A$, so that in addition to being rounded they satisfy the triangle inequality on the right.
  \begin{enumerate}\setcounter{enumi}{4}
  \item Assuming $\rclim(x)\approx_\phi \rcrat(r)$ and $-\epsilon<|r-s|<\epsilon$, we must show $\rclim(x)\approx_{\phi+\epsilon} \rcrat(s)$.
    By definition of $\approx$, the former means $x_\delta \approx_{\phi-\delta} \rcrat(r)$, so that above triangle inequality implies $x_\delta \approx_{\epsilon+\phi-\delta} \rcrat(s)$, hence $\rclim(x)\approx_{\phi+\epsilon} \rcrat(s)$ as desired.
  \item We assume $\phi,\epsilon,\delta:\Qp$ such that $\rclim(x)\approx_\phi \rcrat(r)$ and $\rcrat(r) \close{\epsilon-\delta} y_\delta$, and two unneeded inductive hypotheses.
    By definition, we have $\eta:\Qp$ such that $x_\eta \approx_{\phi-\eta} \rcrat(r)$, so the inductive triangle inequality gives $x_\eta \approx_{\phi+\epsilon-\eta-\delta} y_\delta$.
    The definition of $\approx$ then immediately yields $\rclim(x) \approx_{\phi+\epsilon} \rclim(y)$.
  \item We assume $\phi,\epsilon,\delta:\Qp$ such that $\rclim(x)\approx_\phi \rclim(y)$ and $y_\delta \close{\epsilon-\delta} \rcrat(r)$, and two unneeded inductive hypotheses.
    By definition we have $\xi,\theta:\Qp$ such that $x_\xi \approx_{\phi-\xi-\theta} y_\theta$.
    Since $y$ is a Cauchy approximation, we have $y_\theta \close{\theta+\delta} y_\delta$, so the inductive triangle inequality gives $x_\xi \approx_{\phi+\delta-\xi} y_\delta$ and then $x_\xi \close{\phi+\epsilon-\xi} \rcrat(r)$.
    The definition of $\approx$ then gives $\rclim(x) \approx_{\phi+\epsilon}\rcrat(r)$, as desired.
  \item Finally, we assume $\phi,\epsilon,\delta,\eta:\Qp$ such that $\rclim(x)\approx_\phi \rclim(y)$ and $y_\delta \close{\epsilon-\delta-\eta} z_\eta$.
    Then as before we have $\xi,\theta:\Qp$ with $x_\xi \approx_{\phi-\xi-\theta} y_\theta$, and two applications of the triangle inequality suffices as before.
  \end{enumerate}

  This completes the two inner recursions, and thus the definitions of the families of relations $(\rcrat(q)\approx_{(\blank)}\blank)$ and $(\rclim(x)\approx_{(\blank)}\blank)$.
  Since all are elements of $A$, they are rounded and satisfy the triangle inequality on the right with respect to $\closesym$.
  What remains is to verify the conditions relating to $\bsim$, which is to say that these relations satisfy the triangle inequality on the \emph{left} with respect to the constructors of $\closesym$.
  The four cases correspond to the four choices of rational or limit points for $u$ and $v$ in~\eqref{eq:RC-sim-ltri}, and since they are all mere propositions, we may apply $\RC$-induction and assume that $w$ is also either rational or a limit.
  This yields another eight cases, whose proofs are essentially identical to those just given; so we will not subject the reader to them.
\end{proof}

We can now prove:

\begin{thm}\label{thm:RC-sim-characterization}
  For any $u,v:\RC$ and $\epsilon:\Qp$ we have $(u\close\epsilon v) = (u\approx_\epsilon v)$.
\end{thm}
\begin{proof}
  Since both are mere propositions, it suffices to prove bidirectional implication.
  For the left-to-right direction, we use $\closesym$-induction applied to $C(u,v,\epsilon)\defeq (u\approx_\epsilon v)$.
  Thus, it suffices to consider the four constructors of $\closesym$.
  In each case, $u$ and $v$ are specialized to either rational points or limits, so that the definition of $\approx$ evaluates, and the inductive hypothesis always applies.

  For the right-to-left direction, we use $\RC$-induction to assume that $u$ and $v$ are rational points or limits, allowing $\approx$ to evaluate.
  But now the definitions of $\approx$, and the inductive hypotheses, supply exactly the data required for the relevant constructors of $\closesym$.
\end{proof}

\index{encode-decode method}%
Stretching a point, one might call $\approx$ a fibration of ``codes'' for $\closesym$, with the two directions of the above proof being \encode and \decode respectively.
By the definition of $\approx$, from \autoref{thm:RC-sim-characterization} we get equivalences
\begin{align*}
  (\rcrat(q) \close\epsilon \rcrat(r))  &=
  (-\epsilon < q - r < \epsilon)\\
  (\rcrat(q) \close\epsilon \rclim(y)) &=
  \exis{\delta : \Qp} \rcrat(q) \close{\epsilon - \delta} y_\delta\\
  (\rclim(x) \close\epsilon \rcrat(r)) &=
  \exis{\delta : \Qp} x_\delta \close{\epsilon - \delta} \rcrat(r)\\
  (\rclim(x) \close\epsilon \rclim(y)) &=
  \exis{\delta, \eta : \Qp} x_\delta \close{\epsilon - \delta - \eta} y_\eta.
\end{align*}
Our proof also provides the following additional information.

\begin{cor}
  \index{triangle!inequality for R@inequality for $\RC$}%
  \indexsee{inequality!triangle}{triangle inequality}%
  $\closesym$ is rounded\index{rounded!relation} and satisfies the triangle inequality:
    \begin{gather}
      \eqvspaced{
        (u \close\epsilon v)
      }{
        \exis{\theta : \Qp} u \close{\epsilon - \theta} v
      }\\
      (u\close\epsilon v) \to (v\close\delta w) \to (u\close{\epsilon+\delta} w). \label{item:RC-sim-triangle}
    \end{gather}
\end{cor}

With the triangle inequality in hand, we can show that ``limits'' of Cauchy approximations actually behave like limits.

\begin{lem}\label{thm:RC-sim-lim}
  For any $u:\RC$, Cauchy approximation $y$, and $\epsilon,\delta:\Qp$, if $u\close\epsilon y_\delta$ then $u\close{\epsilon+\delta} \rclim(y)$.
\end{lem}
\begin{proof}
  We use $\RC$-induction on $u$.
  If $u$ is $\rcrat(q)$, then this is exactly the second constructor of $\closesym$.
  Now suppose $u$ is $\rclim(x)$, and that each $x_\eta$ has the property that for any $y,\epsilon,\delta$, if $x_\eta\close\epsilon y_\delta$ then $x_\eta \close{\epsilon+\delta} \rclim(y)$.
  In particular, taking $y\defeq x$ and $\delta\defeq\eta$ in this assumption, we conclude that $x_\eta \close{\eta+\theta} \rclim(x)$ for any $\eta,\theta:\Qp$.

  Now let $y,\epsilon,\delta$ be arbitrary and assume $\rclim(x) \close\epsilon y_\delta$.
  By roundedness, there is a $\theta$ such that $\rclim(x) \close{\epsilon-\theta} y_\delta$.
  Then by the above observation, for any $\eta$ we have $x_\eta \close{\eta+\theta/2} \rclim(x)$, and hence $x_\eta \close{\epsilon+\eta-\theta/2} y_\delta$ by the triangle inequality.
  Hence, the fourth constructor of $\closesym$ yields $\rclim(x) \close{\epsilon+2\eta+\delta-\theta/2} \rclim(y)$.
  Thus, if we choose $\eta \defeq \theta/4$, the result follows.
\end{proof}

\begin{lem}\label{thm:RC-sim-lim-term}
  For any Cauchy approximation $y$ and any $\delta,\eta:\Qp$ we have $y_\delta \close{\delta+\eta} \rclim(y)$.
\end{lem}
\begin{proof}
  Take $u\defeq y_\delta$ and $\epsilon\defeq \eta$ in the previous lemma.
\end{proof}

\begin{rmk}
  We might have expected to have $y_\delta \close{\delta} \rclim(y)$, but this fails in examples.
  For instance, consider $x$ defined by $x_\epsilon \defeq \epsilon$.
  Its limit is clearly $0$, but we do not have $|\epsilon - 0 |<\epsilon$, only $\le$.
\end{rmk}

As an application, \autoref{thm:RC-sim-lim-term} enables us to show that the extensions of Lipschitz functions from \autoref{RC-extend-Q-Lipschitz} are unique.

\begin{lem}\label{RC-continuous-eq}
  \index{function!continuous}%
  Let $f,g:\RC\to\RC$ be continuous, in the sense that
  \[ \fall{u:\RC}{\epsilon:\Qp}\exis{\delta:\Qp}\fall{v:\RC} (u\close\delta v) \to (f(u) \close\epsilon f(v)) \]
  and analogously for $g$.
  If $f(\rcrat(q))=g(\rcrat(q))$ for all $q:\Q$, then $f=g$.
\end{lem}
\begin{proof}
  We prove $f(u)=g(u)$ for all $u$ by $\RC$-induction.
  The rational case is just the hypothesis.
  Thus, suppose $f(x_\delta)=g(x_\delta)$ for all $\delta$.
  We will show that $f(\rclim(x))\close\epsilon g(\rclim(x))$ for all $\epsilon$, so that the path constructor of $\RC$ applies.

  Since $f$ and $g$ are continuous, there exist $\theta,\eta$ such that for all $v$, we have
  \begin{align*}
    (\rclim(x)\close\theta v) &\to (f(\rclim(x)) \close{\epsilon/2} f(v))\\
    (\rclim(x)\close\eta v) &\to (g(\rclim(x)) \close{\epsilon/2} g(v)).
  \end{align*}
  Choosing $\delta < \min(\theta,\eta)$, by \autoref{thm:RC-sim-lim-term} we have both $\rclim(x)\close\theta y_\delta$ and $\rclim(x)\close\eta y_\delta$.
  Hence
  \[ f(\rclim(x)) \close{\epsilon/2} f(y_\delta) = g(y_\delta) \close{\epsilon/2} g(\rclim(x))\]
  and thus $f(\rclim(x))\close\epsilon g(\rclim(x))$ by the triangle inequality.
\end{proof}

\subsection{The algebraic structure of Cauchy reals}
\label{sec:algebr-struct-cauchy}

We first define the additive structure $(\RC, 0, {+}, {-})$. Clearly, the additive unit element
$0$ is just $\rcrat(0)$, while the additive inverse ${-} : \RC \to \RC$ is obtained as the
extension of the additive inverse ${-} : \Q \to \Q$, using \autoref{RC-extend-Q-Lipschitz}
with Lipschitz constant~$1$. We have to work a bit harder for addition.

\begin{lem} \label{RC-binary-nonexpanding-extension}
  Suppose $f : \Q \times \Q \to \Q$ satisfies, for all $q, r, s : \Q$,
  \begin{equation*}
    |f(q, s) - f(r, s)| \leq |q - r|
    \qquad\text{and}\qquad
    |f(q, r) - f(q, s)| \leq |r - s|.
  \end{equation*}
  Then there is a function $\bar{f} : \RC \times \RC \to \RC$ such that
  $\bar{f}(\rcrat(q), \rcrat(r)) = f(q,r)$ for all $q, r : \Q$. Furthermore,
  for all $u, v, w : \RC$ and $q : \Qp$,
  \begin{equation*}
    u \close\epsilon v \Rightarrow \bar{f}(u,w) \close\epsilon \bar{f}(v,w)
    \quad\text{and}\quad
    v \close\epsilon w \Rightarrow \bar{f}(u,v) \close\epsilon \bar{f}(u,w).
  \end{equation*}
\end{lem}

\begin{proof}
  We use $(\RC, {\closesym})$-recursion to construct the curried form of $\bar{f}$ as a map
  $\RC \to A$ where $A$ is the space of non-expanding\index{function!non-expanding}\index{non-expanding function} real-valued
  functions:
  \begin{equation*}
    A \defeq
    \setof{ h : \RC \to \RC |
      \fall{\epsilon : \Qp} \fall{u, v : \RC}
      u \close\epsilon v \Rightarrow h(u) \close\epsilon h(v)
    }.
  \end{equation*}
  We shall also need a suitable $\bsim_\epsilon$ on $A$, which we define as
  \begin{equation*}
    (h \bsim_\epsilon k) \defeq \fall{u : \RC} h(u) \close\epsilon k(u).
  \end{equation*}
  Clearly, if $\fall{\epsilon : \Qp} h \bsim_\epsilon k$ then $h(u) = k(u)$ for all $u :
  \RC$, so $\bsim$ is separated.

  For the base case we define $\bar{f}(\rcrat(q)) : A$, where $q : \Q$, as the
  extension of the Lipschitz map $\lam{r} f(q,r)$ from $\Q \to \Q$ to $\RC \to \RC$, as
  constructed in \autoref{RC-extend-Q-Lipschitz} with Lipschitz constant~$1$. Next, for a
  Cauchy approximation $x$, we define $\bar{f}(\rclim(x)) : \RC \to \RC$ as
  \begin{equation*}
    \bar{f}(\rclim(x))(v) \defeq \rclim (\lam{\epsilon} \bar{f}(x_\epsilon)(v)).
  \end{equation*}
  For this to be a valid definition, $\lam{\epsilon} \bar{f}(x_\epsilon)(v)$ should be a
  Cauchy approximation, so consider any $\delta, \epsilon : \Q$. Then by assumption
  $\bar{f}(x_\delta) \bsim_{\delta + \epsilon} \bar{f}(x_\epsilon)$, hence
  $\bar{f}(x_\delta)(v) \close{\delta + \epsilon} \bar{f}(x_\epsilon)(v)$. Furthermore,
  $\bar{f}(\rclim(x))$ is non-expanding because $\bar{f}(x_\epsilon)$ is such by induction
  hypothesis. Indeed, if $u \close\epsilon v$ then, for all $\epsilon : \Q$,
  \begin{equation*}
    \bar{f}(x_{\epsilon/3})(u) \close{\epsilon/3} \bar{f}(x_{\epsilon/3})(v),
  \end{equation*}
  therefore $\bar{f}(\rclim(x))(u) \close\epsilon \bar{f}(\rclim(x))(v)$ by the fourth constructor of $\closesym$.

  We still have to check four more conditions, let us illustrate just one. Suppose
  $\epsilon : \Qp$ and for some $\delta : \Qp$ we have $\rcrat(q) \close{\epsilon - \delta}
  y_\delta$ and $\bar{f}(\rcrat(q)) \bsim_{\epsilon - \delta} \bar{f}(y_\delta)$. To show
  $\bar{f}(\rcrat(q)) \bsim_\epsilon \bar{f}(\rclim(y))$, consider any $v : \RC$ and observe that
  \begin{equation*}
    \bar{f}(\rcrat(q))(v) \close{\epsilon - \delta} \bar{f}(y_\delta)(v).
  \end{equation*}
  Therefore, by the second constructor of $\closesym$, we have
  \narrowequation{\bar{f}(\rcrat(q))(v) \close\epsilon \bar{f}(\rclim(y))(v)}
  as required.
\end{proof}

We may apply \autoref{RC-binary-nonexpanding-extension} to any bivariate rational function
which is non-expanding separately in each variable. Addition is such a function, therefore
we get ${+} : \RC \times \RC \to \RC$.
\indexdef{addition!of Cauchy reals}%
Furthermore, the extension is unique as long as we
require it to be non-expanding in each variable, and just as in the univariate case,
identities on rationals extend to identities on reals. Since composition of non-expanding
maps is again non-expanding, we may conclude that addition satisfies the usual properties,
such as commutativity and associativity.
\index{associativity!of addition!of Cauchy reals}%
Therefore, $(\RC, 0, {+}, {-})$ is a commutative
group.

We may also apply \autoref{RC-binary-nonexpanding-extension} to the functions $\min : \Q \times
\Q \to \Q$ and $\max : \Q \times \Q \to \Q$, which turns $\RC$ into a lattice. The partial
order $\leq$ on $\RC$ is defined in terms of $\max$ as
\symlabel{leq-RC}
\index{order!non-strict}%
\index{non-strict order}%
\begin{equation*}
  (u \leq v) \defeq (\max(u, v) = v).
\end{equation*}
The relation $\leq$ is a partial order because it is such on $\Q$, and the axioms of a
partial order are expressible as equations in terms of $\min$ and $\max$, so they transfer
to $\RC$.

\index{absolute value}%
Another function which extends to $\RC$ by the same method is the absolute value $|{\blank}|$.
Again, it has the expected properties because they transfer from $\Q$ to $\RC$.

\symlabel{lt-RC}
From $\leq$ we get the strict order $<$ by
\index{strict!order}%
\index{order!strict}%
\begin{equation*}
  (u < v) \defeq \exis{q, r : \Q} (u \leq \rcrat(q)) \land (q < r) \land (\rcrat(r) \leq v).
\end{equation*}
That is, $u < v$ holds when there merely exists a pair of rational numbers $q < r$ such that $x \leq
\rcrat(q)$ and $\rcrat(r) \leq v$. It is not hard to check that $<$ is irreflexive and
transitive, and has other properties that are expected for an ordered field.
The archimedean principle follows directly from the definition of~$<$.

\index{ordered field!archimedean}%
\begin{thm}[Archimedean principle for $\RC$] \label{RC-archimedean}
  For every $u, v : \RC$ such that $u < v$ there merely exists $q : \Q$ such that $u < q < v$.
\end{thm}

\begin{proof}
  From $u < v$ we merely get $r, s : \Q$ such that $u \leq r < s \leq v$, and we may take $q
  \defeq (r + s) / 2$.
\end{proof}

We now have enough structure on $\RC$ to express $u \close\epsilon v$ with standard concepts.

\begin{lem}\label{thm:RC-le-grow}
  If $q:\Q$ and $u:\RC$ satisfy $u\le \rcrat(q)$, then for any $v:\RC$ and $\epsilon:\Qp$, if $u\close\epsilon v$ then $v\le \rcrat(q+\epsilon)$.
\end{lem}
\begin{proof}
  Note that the function $\max(\rcrat(q),\blank):\RC\to\RC$ is Lipschitz with constant $1$.
  First consider the case when $u=\rcrat(r)$ is rational.
  For this we use induction on $v$.
  If $v$ is rational, then the statement is obvious.
  If $v$ is $\rclim(y)$, we assume inductively that for any $\epsilon,\delta$, if $\rcrat(r)\close\epsilon y_\delta$ then $y_\delta \le \rcrat(q+\epsilon)$, i.e.\ $\max(\rcrat(q+\epsilon),y_\delta)=\rcrat(q+\epsilon)$.

  Now assuming $\epsilon$ and $\rcrat(r)\close\epsilon \rclim(y)$, we have $\theta$ such that $\rcrat(r)\close{\epsilon-\theta} \rclim(y)$, hence $\rcrat(r)\close\epsilon y_\delta$ whenever $\delta<\theta$.
  Thus, the inductive hypothesis gives $\max(\rcrat(q+\epsilon),y_\delta)=\rcrat(q+\epsilon)$ for such $\delta$.
  But by definition,
  \[\max(\rcrat(q+\epsilon),\rclim(y)) \jdeq \rclim(\lam{\delta} \max(\rcrat(q+\epsilon),y_\delta)).\]
  Since the limit of an eventually constant Cauchy approximation is that constant, we have 
  \[\max(\rcrat(q+\epsilon),\rclim(y)) = \rcrat(q+\epsilon),\] hence $\rclim(y)\le \rcrat(q+\epsilon)$.
  
  Now consider a general $u:\RC$.
  Since $u\le \rcrat(q)$ means $\max(\rcrat(q),u)=\rcrat(q)$, the assumption $u\close\epsilon v$ and the Lipschitz property of $\max(\rcrat(q),-)$ imply $\max(\rcrat(q),v) \close\epsilon \rcrat(q)$.
  Thus, since $\rcrat(q)\le \rcrat(q)$, the first case implies $\max(\rcrat(q),v) \le \rcrat(q+\epsilon)$, and hence $v\le \rcrat(q+\epsilon)$ by transitivity of $\le$.
\end{proof}

\begin{lem}\label{thm:RC-lt-open}
  Suppose $q:\Q$ and $u:\RC$ satisfy $u<\rcrat(q)$.  Then:
  \begin{enumerate}
  \item For any $v:\RC$ and $\epsilon:\Qp$, if $u\close\epsilon v$ then $v< \rcrat(q+\epsilon)$.\label{item:RCltopen1}
  \item There exists $\epsilon:\Qp$ such that for any $v:\RC$, if $u\close\epsilon v$ we have $v<\rcrat(q)$.\label{item:RCltopen2}
  \end{enumerate}
\end{lem}
\begin{proof}
  By definition, $u<\rcrat(q)$ means there is $r:\Q$ with $r<q$ and $u\le \rcrat(r)$.
  Then by \autoref{thm:RC-le-grow}, for any $\epsilon$, if $u\close\epsilon v$ then $v\le \rcrat(r+\epsilon)$.
  Conclusion~\ref{item:RCltopen1} follows immediately since $r+\epsilon<q+\epsilon$, while for~\ref{item:RCltopen2} we can take any $\epsilon <q-r$.
\end{proof}

We are now able to show that the auxiliary relation $\closesym$ is what we think it is.

\begin{thm} \label{RC-sim-eqv-le}
  \index{distance}%
  $\eqv{(u \close\epsilon v)}{(|u - v| < \rcrat(\epsilon))}$
  for all $u, v : \RC$ and $\epsilon : \Qp$.
\end{thm}
\begin{proof}
  The Lipschitz properties of subtraction and absolute value imply that if $u\close\epsilon v$, then $|u-v| \close\epsilon |u-u| = 0$.
  Thus, for the left-to-right direction, it will suffice to show that if $u\close\epsilon 0$, then $|u|<\rcrat(\epsilon)$.
  We proceed by $\RC$-induction on $u$.

  If $u$ is rational, the statement follows immediately since absolute value and order extend the standard ones on $\Qp$.
  If $u$ is $\rclim(x)$, then by roundedness we have $\theta:\Qp$ with $\rclim(x)\close{\epsilon-\theta} 0$.
  By the triangle inequality, therefore, we have $x_{\theta/3} \close{\epsilon-2\theta/3} 0$, so the inductive hypothesis yields $|x_{\theta/3}|<\rcrat(\epsilon-2\theta/3)$.
  But $x_{\theta/3} \close{2\theta/3} \rclim(x)$, hence $|x_{\theta/3}| \close{2\theta/3} |\rclim(x)|$ by the Lipschitz property, so \autoref{thm:RC-lt-open}\ref{item:RCltopen1} implies $|\rclim(x)|<\rcrat(\epsilon)$.

  In the other direction, we use $\RC$-induction on $u$ and $v$.
  If both are rational, this is the first constructor of $\closesym$.

  If $u$ is $\rcrat(q)$ and $v$ is $\rclim(y)$, we assume inductively that for any $\epsilon,\delta$, if $|\rcrat(q)-y_\delta|<\rcrat(\epsilon)$ then $\rcrat(q) \close{\epsilon} y_\delta$.
  Fix an $\epsilon$ such that $|\rcrat(q) - \rclim(y)|<\rcrat(\epsilon)$.
  Since $\Q$ is order-dense in $\RC$, there exists $\theta<\epsilon$ with $|\rcrat(q) - \rclim(y)|<\rcrat(\theta)$.
  Now for any $\delta,\eta$ we have $\rclim(y)\close{2\delta} y_\delta$, hence by the Lipschitz property
  \[ |\rcrat(q) - \rclim(y)| \close{\delta+\eta} |\rcrat(q) - y_\delta|. \]
  Thus, by \autoref{thm:RC-lt-open}\ref{item:RCltopen1}, we have $|\rcrat(q) - y_\delta| < \rcrat(\theta+2\delta)$.
  So by the inductive hypothesis, $\rcrat(q) \close{\theta+2\delta} y_\delta$, and thus $\rcrat(q)\close{\theta+4\delta} \rclim(y)$ by the triangle inequality.
  Thus, it suffices to choose $\delta \defeq (\epsilon-\theta)/4$.

  The remaining two cases are entirely analogous.
\end{proof}

\indexdef{multiplication!of Cauchy reals}%
Next, we would like to equip $\RC$ with multiplicative structure. For each $q : \Q$ the
map $r \mapsto q \cdot r$ is Lipschitz with constant\footnote{We defined Lipschitz
  constants as \emph{positive} rational numbers.} $|q| + 1$, and so we can extend it to
multiplication by $q$ on the real numbers. Therefore $\RC$ is a vector space\index{vector!space} over $\Q$.
In general, we can define multiplication of real numbers as
\begin{equation}
  u \cdot v \defeq
  {\textstyle \frac{1}{2}} \cdot ((u + v)^2 - u^2 - v^2),\label{mult-from-square}
\end{equation}
so we just need squaring\index{squaring function} $u \mapsto u^2$ as a map $\RC \to \RC$. Squaring is not a
Lipschitz map, but it is Lipschitz on every bounded domain, which allows us to patch it
together. Define the open and closed intervals
\indexdef{interval!open and closed}%
\indexdef{open!interval}%
\indexdef{closed!interval}%
\begin{equation*}
  [u,v] \defeq \setof{ x : \RC | u \leq x \leq v }
  \qquad\text{and}\qquad
  (u,v) \defeq \setof{ x : \RC | u < x < v }.
\end{equation*}
Although technically an element of $[u,v]$ or $(u,v)$ is a Cauchy real number together with a proof, since the latter inhabits a mere proposition it is uninteresting.
Thus, as is common with subset types, we generally write simply $x:[u,v]$ whenever $x:\RC$ is such that $u\leq x \leq v$, and similarly.

\begin{thm} \label{RC-squaring}
  There exists a unique function ${(\blank)}^2 : \RC \to \RC$ which extends squaring $q \mapsto
  q^2$ of rational numbers and satisfies
  \begin{equation*}
    \fall{n : \N}
    \fall{u, v : [-n, n]}
    |u^2 - v^2| \leq 2 \cdot n \cdot |u - v|.
  \end{equation*}
\end{thm}

\begin{proof}
  We first observe that for every $u : \RC$ there merely exists $n : \N$ such that $-n
  \leq u \leq n$, see \autoref{ex:traditional-archimedean}, so the map
  \begin{equation*}
    e : \Parens{\sm{n : \N} [-n, n]} \to \RC
    \qquad\text{defined by}\qquad
    e(n, x) \defeq x
  \end{equation*}
  is surjective. Next, for each $n : \N$, the squaring map
  \begin{equation*}
    s_n : \setof{ q : \Q | -n \leq q \leq n } \to \Q
    \qquad\text{defined by}\qquad
    s_n(q) \defeq q^2
  \end{equation*}
  is Lipschitz with constant $2 n$, so we can use \autoref{RC-extend-Q-Lipschitz} to
  extend it to a map $\bar{s}_n : [-n, n] \to \RC$ with Lipschitz constant $2 n$, see
  \autoref{RC-Lipschitz-on-interval} for details. The maps $\bar{s}_n$ are compatible: if
  $m < n$ for some $m, n : \N$ then $s_n$ restricted to $[-m, m]$ must agree with $s_m$
  because both are Lipschitz, and therefore continuous in the sense
  of~\autoref{RC-continuous-eq}. Therefore, by \autoref{lem:images_are_coequalizers} the map
  \begin{equation*}
    \Parens{\sm{n : \N} [-n, n]} \to \RC,
    \qquad\text{given by}\qquad
    (n, x) \mapsto s_n(x)
  \end{equation*}
  factors uniquely through $\RC$ to give us the desired function.
\end{proof}

At this point we have the ring structure of the reals and the archimedean order. To
establish $\RC$ as an archimedean ordered field, we still need inverses.

\begin{thm}
  \index{apartness}%
  A Cauchy real is invertible if, and only if, it is apart from zero.
\end{thm}

\begin{proof}
  First, suppose $u : \RC$ has an inverse $v : \RC$ By the archimedean principle there is $q :
  \Q$ such that $|v| < q$. Then $1 = |u v| < |u| \cdot v < |u| \cdot q$ and hence $|u| >
  1/q$, which is to say that $u \apart 0$.

  For the converse we construct the inverse map
  \begin{equation*}
    ({\blank})^{-1} : \setof{ u : \RC | u \apart 0 } \to \RC
  \end{equation*}
  by patching together functions, similarly to the construction of squaring in
  \autoref{RC-squaring}. We only outline the main steps. For every $q : \Q$ let
  \begin{equation*}
    [q, \infty) \defeq \setof{u : \RC | q \leq u}
    \qquad\text{and}\qquad
    (-\infty, q] \defeq \setof{u : \RC | u \leq -q}.
  \end{equation*}
  Then, as $q$ ranges over $\Qp$, the types $(-\infty, q]$ and $[q, \infty)$ jointly cover
  $\setof{u : \RC | u \apart 0}$. On each such $[q, \infty)$ and $(-\infty, q]$ the
  inverse function is obtained by an application of \autoref{RC-extend-Q-Lipschitz}
  with Lipschitz constant $1/q^2$. Finally, \autoref{lem:images_are_coequalizers}
  guarantees that the inverse function factors uniquely through $\setof{ u : \RC | u
    \apart 0 }$.
\end{proof}

We summarize the algebraic structure of $\RC$ with a theorem.

\begin{thm} \label{RC-archimedean-ordered-field}
  The Cauchy reals form an archimedean ordered field.
\end{thm}

\subsection{Cauchy reals are Cauchy complete}
\label{sec:cauchy-reals-cauchy-complete}

We constructed $\RC$ by closing $\Q$ under limits of Cauchy approximations, so it better
be the case that $\RC$ is Cauchy complete. Thanks to \autoref{RC-sim-eqv-le} there is no
difference between a Cauchy approximation $x : \Qp \to \RC$ as defined in the construction
of $\RC$, and a Cauchy approximation in the sense of \autoref{defn:cauchy-approximation}
(adapted to $\RC$).

Thus, given a Cauchy approximation $x : \Qp \to \RC$ it is quite natural to expect that
$\rclim(x)$ is its limit, where the notion of limit is defined as in
\autoref{defn:cauchy-approximation}. But this is so by \autoref{RC-sim-eqv-le} and
\autoref{thm:RC-sim-lim-term}. We have proved:

\begin{thm}
  Every Cauchy approximation in $\RC$ has a limit.
\end{thm}

An archimedean ordered field in which every Cauchy approximation has a limit is called
\define{Cauchy complete}.
\indexdef{Cauchy!completeness}%
\indexdef{complete!ordered field, Cauchy}%
\index{ordered field}%
The Cauchy reals are the least such field.

\begin{thm} \label{RC-initial-Cauchy-complete}
  The Cauchy reals embed into every Cauchy complete archimedean ordered field.
\end{thm}

\begin{proof}
  \index{limit!of a Cauchy approximation}%
  Suppose $F$ is a Cauchy complete archimedean ordered field. Because limits are unique,
  there is an operator $\lim$ which takes Cauchy approximations in $F$ to their limits. We
  define the embedding $e : \RC \to F$ by $(\RC, {\closesym})$-recursion as
  \begin{equation*}
    e(\rcrat(q)) \defeq q
    \qquad\text{and}\qquad
    e(\rclim(x)) \defeq \lim (e \circ x).
  \end{equation*}
  A suitable $\bsim$ on $F$ is
  \begin{equation*}
    (a \bsim_\epsilon b) \defeq |a - b| < \epsilon.
  \end{equation*}
  This is a separated relation because $F$ is archimedean. The rest of the clauses for
  $(\RC, {\closesym})$-recursion are easily checked. One would also have to check that $e$ is
  an embedding of ordered fields which fixes the rationals.
\end{proof}

\index{real numbers!Cauchy|)}%

\section{Comparison of Cauchy and Dedekind reals}
\label{sec:comp-cauchy-dedek}

\index{real numbers!Dedekind|(}%
\index{real numbers!Cauchy|(}%
\index{depression|(}

Let us also say something about the relationship between the Cauchy and Dedekind reals. By
\autoref{RC-archimedean-ordered-field}, $\RC$ is an archimedean ordered field. It is also
admissible\index{ordered field!admissible} for $\Omega$, as can be easily checked. (In case $\Omega$ is the initial
$\sigma$-frame
\index{initial!sigma-frame@$\sigma$-frame}%
\index{sigma-frame@$\sigma$-frame!initial}%
it takes a simple induction, while in other cases it is immediate.)
Therefore, by \autoref{RD-final-field} there is an embedding of ordered fields
\begin{equation*}
  \RC \to \RD
\end{equation*}
which fixes the rational numbers.
(We could also obtain this from \autoref{RC-initial-Cauchy-complete,RD-cauchy-complete}.)
In general we do not expect $\RC$ and $\RD$ to coincide
without further assumptions. 

\begin{lem} \label{lem:untruncated-linearity-reals-coincide}
  If for every $x : \RD$ there merely exists
  \begin{equation}
    \label{eq:untruncated-linearity}
    c : \prd{q, r : \Q} (q < r) \to (q < x) + (x < r)
  \end{equation}
  then the Cauchy and Dedekind reals coincide.
\end{lem}

\begin{proof}
  Note that the type in~\eqref{eq:untruncated-linearity} is an untruncated variant
  of~\eqref{eq:RD-linear-order}, which states that~$<$ is a weak linear order.
  We already know that $\RC$ embeds into $\RD$, so it suffices to show that every Dedekind
  real merely is the limit of a Cauchy sequence\index{Cauchy!sequence} of rational numbers.

  Consider any $x : \RD$. By assumption there merely exists $c$ as in the statement of the
  lemma, and by inhabitation of cuts\index{cut!Dedekind} there merely exist $a, b : \Q$ such that $a < x < b$.
  We construct a sequence\index{sequence} $f : \N \to \setof{ \pairr{q, r} \in \Q \times \Q | q < r }$ by
  recursion:
  \begin{enumerate}
  \item Set $f(0) \defeq \pairr{a, b}$.
  \item Suppose $f(n)$ is already defined as $\pairr{q_n, r_n}$ such that $q_n < r_n$.
    Define $s \defeq (2 q_n + r_n)/3$ and $t \defeq (q_n + 2 r_n)/3$. Then $c(s,t)$
    decides between $s < x$ and $x < t$. If it decides $s < x$ then we set $f(n+1) \defeq
    \pairr{s, r_n}$, otherwise $f(n+1) \defeq \pairr{q_n, t}$.
  \end{enumerate}
  Let us write $\pairr{q_n, r_n}$ for the $n$-th term of the sequence~$f$. Then it is easy
  to see that $q_n < x < r_n$ and $|q_n - r_n| \leq (2/3)^n \cdot |q_0 - r_0|$ for all $n
  : \N$. Therefore $q_0, q_1, \ldots$ and $r_0, r_1, \ldots$ are both Cauchy sequences
  converging to the Dedekind cut~$x$. We have shown that for every $x : \RD$ there merely
  exists a Cauchy sequence converging to $x$.
\end{proof}

The lemma implies that either countable choice or excluded middle suffice for coincidence
of $\RC$ and $\RD$.

\begin{cor} \label{when-reals-coincide}
  \index{axiom!of choice!countable}%
  \index{excluded middle}%
  If excluded middle or countable choice holds then $\RC$ and $\RD$ are equivalent.
\end{cor}

\begin{proof}
  If excluded middle holds then $(x < y) \to (x < z) + (z < y)$ can be proved: either $x <
  z$ or $\lnot (x < z)$. In the former case we are done, while in the latter we get $z <
  y$ because $z \leq x < y$. Therefore, we get~\eqref{eq:untruncated-linearity} so that we
  can apply \autoref{lem:untruncated-linearity-reals-coincide}.

  Suppose countable choice holds. The set $S = \setof{ \pairr{q, r} \in \Q \times \Q | q <
    r }$ is equivalent to $\N$, so we may apply countable choice to the statement that $x$
  is located,
  \begin{equation*}
    \fall{\pairr{q, r} : S} (q < x) \lor (x < r).
  \end{equation*}
  Note that $(q < x) \lor (x < r)$ is expressible as an existential statement $\exis{b :
    \bool} (b = \bfalse \to q < x) \land (b = \btrue \to x < r)$. The (curried form) of
  the choice function is then precisely~\eqref{eq:untruncated-linearity} so that
  \autoref{lem:untruncated-linearity-reals-coincide} is applicable again.
\end{proof}

\index{real numbers!Dedekind|)}%
\index{real numbers!Cauchy|)}%
\index{real numbers!agree}%

\index{depression|)}

\section{Compactness of the interval}
\label{sec:compactness-interval}

\index{mathematics!classical|(}%
\index{mathematics!constructive|(}%

We already pointed out that our constructions of reals are entirely compatible with
classical logic. Thus, by assuming the law of excluded middle~\eqref{eq:lem} and the axiom
of choice~\eqref{eq:ac} we could develop classical analysis,\index{classical!analysis}\index{analysis!classical} which would essentially
amount to copying any standard book on analysis.

\index{analysis!constructive}%
\index{constructive!analysis}%
Nevertheless, anyone interested in computation, for example a numerical analyst, ought to
be curious about developing analysis in a computationally meaningful setting. That
analysis in a constructive setting is even possible was demonstrated by~\cite{Bishop1967}.
As a sample of the differences and similarities between classical and constructive
analysis we shall briefly discuss just one topic---compactness of the closed interval
$[0,1]$ and a couple of theorems surrounding the concept.

Compactness is no exception to the common phenomenon in constructive mathematics that
classically equivalent notions bifurcate. The three most frequently used notions of
compactness are:
\indexdef{compactness}%
\begin{enumerate}
\item \define{metrically compact:} ``Cauchy complete and totally bounded'',
  \indexdef{metrically compact}%
  \indexdef{compactness!metric}%
\item \define{Bolzano--Weierstra\ss{} compact:} ``every sequence has a convergent subsequence'',
  \index{compactness!Bolzano--Weierstrass@Bolzano--Weierstra\ss{}}%
  \indexsee{Bolzano--Weierstrass@Bolzano--Weierstra\ss{}}{compactness}%
  \index{sequence}%
\item \define{Heine-Borel compact:} ``every open cover has a finite subcover''.
  \index{compactness!Heine-Borel}%
  \indexsee{Heine-Borel}{compactness}%
\end{enumerate}
These are all equivalent in classical mathematics.
Let us see how they fare in homotopy type theory. We can use either the Dedekind or the
Cauchy reals, so we shall denote the reals just as~$\R$. We first recall several basic
definitions.

\indexsee{space!metric}{metric space}
\index{metric space|(}%

\begin{defn} \label{defn:metric-space}
  A \define{metric space}
  \indexdef{metric space}%
  $(M, d)$ is a set $M$ with a map $d : M \times M \to \R$
  satisfying, for all $x, y, z : M$,
  \begin{align*}
    d(x,y) &\geq 0, &
    d(x,y) &= d(y,x), \\
    d(x,y) &= 0 \Leftrightarrow x = y, &
    d(x,z) &\leq d(x,y) + d(y,z).
  \end{align*}
\end{defn}

\begin{defn} \label{defn:complete-metric-space}
  A \define{Cauchy approximation}
  \index{Cauchy!approximation}%
  in $M$ is a sequence $x : \Qp \to M$ satisfying
  \begin{equation*}
    \fall{\delta, \epsilon} d(x_\delta, x_\epsilon) < \delta + \epsilon.
  \end{equation*}
  \index{limit!of a Cauchy approximation}%
  The \define{limit} of a Cauchy approximation $x : \Qp \to M$ is a point $\ell : M$
  satisfying
  \begin{equation*}
    \fall{\epsilon, \theta : \Qp} d(x_\epsilon, \ell) < \epsilon + \theta.
  \end{equation*}
  \indexdef{metric space!complete}%
  \indexdef{complete!metric space}%
  A \define{complete metric space} is one in which every Cauchy approximation has a limit.
\end{defn}

\begin{defn} \label{defn:total-bounded-metric-space}
  For a positive rational $\epsilon$, an \define{$\epsilon$-net}
  \indexdef{epsilon-net@$\epsilon$-net}%
  in a metric space $(M,
  d)$ is an element of
  \begin{equation*}
    \sm{n : \N}{x_1, \ldots, x_n : M}
    \fall{y : M} \exis{k \leq n} d(x_k, y) < \epsilon.
  \end{equation*}
  In words, this is a finite sequence of points $x_1, \ldots, x_n$ such that every point
  in $M$ merely is within $\epsilon$ of some~$x_k$.

  A metric space $(M, d)$ is \define{totally bounded}
  \indexdef{totally bounded metric space}%
  \indexdef{metric space!totally bounded}%
  when it has $\epsilon$-nets of all
  sizes:
  \begin{equation*}
    \prd{\epsilon : \Qp} 
    \sm{n : \N}{x_1, \ldots, x_n : M}
    \fall{y : M} \exis{k \leq n} d(x_k, y) < \epsilon.
  \end{equation*}
\end{defn}

\begin{rmk}
  In the definition of total boundedness we used sloppy notation $\sm{n : \N}{x_1, \ldots, x_n : M}$. Formally, we should have written $\sm{x : \lst{M}}$ instead,
  where $\lst{M}$ is the inductive type of finite lists\index{type!of lists} from \autoref{sec:bool-nat}.
  However, that would make the rest of the statement a bit more cumbersome to express.
\end{rmk}

Note that in the definition of total boundedness we require pure existence of an
$\epsilon$-net, not mere existence. This way we obtain a function which assigns to each
$\epsilon : \Qp$ a specific $\epsilon$-net. Such a function might be called a ``modulus of
total boundedness''. In general, when porting classical metric notions to homotopy type
theory, we should use propositional truncation sparingly, typically so that we avoid
asking for a non-constant map from $\R$ to $\Q$ or $\N$. For instance, here is the
``correct'' definition of uniform continuity.

\begin{defn} \label{defn:uniformly-continuous}
  A map $f : M \to \R$ on a metric space is \define{uniformly continuous}
  \indexdef{function!uniformly continuous}%
  \indexdef{uniformly continuous function}%
  when
  \begin{equation*}
    \prd{\epsilon : \Qp}
    \sm{\delta : \Qp}
    \fall{x, y : M}
    d(x,y) < \delta \Rightarrow |f(x) - f(y)| < \epsilon.
  \end{equation*}
  In particular, a uniformly continuous map has a modulus of uniform continuity\indexdef{modulus!of uniform continuity},
  which is a function that assigns to each $\epsilon$ a corresponding $\delta$.
\end{defn}

Let us show that $[0,1]$ is compact in the first sense.

\begin{thm} \label{analysis-interval-ctb}
  \index{compactness!metric}%
  \index{interval!open and closed}%
  The closed interval $[0,1]$ is complete and totally bounded.
\end{thm}

\begin{proof}
  Given $\epsilon : \Qp$, there is $n : \N$ such that $2/k < \epsilon$, so we may take the
  $\epsilon$-net $x_i = i/k$ for $i = 0, \ldots, k-1$. This is an $\epsilon$-net because,
  for every $y : [0,1]$ there merely exists $i$ such that $0 \leq i < k$ and $(i -
  1)/k < y < (i+1)/k$, and so $|y - x_i| < 2/k < \epsilon$.

  For completeness of $[0,1]$, consider a Cauchy approximation $x : \Qp \to
  [0,1]$ and let $\ell$ be its limit in $\R$. Since $\max$ and $\min$ are Lipschitz maps,
  the retraction $r : \R \to [0,1]$ defined by $r(x) \defeq \max(0, \min(1, x))$ commutes
  with limits of Cauchy approximations, therefore
  \begin{equation*}
    r(\ell) =
    r (\lim x) =
    \lim (r \circ x) =
    r (\lim x) =
    \ell,
  \end{equation*}
  which means that $0 \leq \ell \leq 1$, as required.
\end{proof}

We thus have at least one good notion of compactness in homotopy type theory.
Unfortunately, it is limited to metric spaces because total boundedness is a metric
notion. We shall consider the other two notions shortly, but first we prove that a
uniformly continuous map on a totally bounded space has a \define{supremum},
\indexsee{least upper bound}{supremum}%
i.e.\ an upper bound which is less than or equal to all other upper bounds.

\begin{thm} \label{ctb-uniformly-continuous-sup}
  \indexdef{supremum!of uniformly continuous function}%
  A uniformly continuous map $f : M \to \R$ on a totally bounded metric space
  $(M, d)$ has a supremum $m : \R$. For every $\epsilon : \Qp$ there exists $u : M$ such
  that $|m - f(u)| < \epsilon$.
\end{thm}

\begin{proof}
  Let $h : \Qp \to \Qp$ be the modulus of uniform continuity of~$f$.
  We define an approximation $x : \Qp \to \R$ as follows: for any $\epsilon : \Q$ total
  boundedness of $M$ gives a $h(\epsilon)$-net $y_0, \ldots, y_n$. Define
  \begin{equation*}
    x_\epsilon \defeq \max (f(y_0), \ldots, f(y_n)).
  \end{equation*}
  We claim that $x$ is a Cauchy approximation. Consider any $\epsilon, \eta : \Q$, so that
  \begin{equation*}
    x_\epsilon \jdeq \max (f(y_0), \ldots, f(y_n))
    \quad\text{and}\quad
    x_\eta \jdeq \max (f(z_0), \ldots, f(z_m))
  \end{equation*}
  for some $h(\epsilon)$-net $y_0, \ldots, y_n$ and $h(\eta)$-net $z_0, \ldots, z_m$.
  Every $z_i$ is merely $h(\epsilon)$-close to some $y_j$, therefore $|f(z_i) - f(y_j)| <
  \epsilon$, from which we may conclude that
  \begin{equation*}
    f(z_i) < \epsilon + f(y_j) \leq \epsilon + x_\epsilon,
  \end{equation*}
  therefore $x_\eta < \epsilon + x_\epsilon$. Symmetrically we obtain $x_\eta < \eta +
  x_\eta$, therefore $|x_\eta - x_\epsilon| < \eta + \epsilon$.

  We claim that $m \defeq \lim x$ is the supremum of~$f$. To prove that $f(x) \leq m$ for
  all $x : M$ it suffices to show $\lnot (m < f(x))$. So suppose to the contrary that $m <
  f(x)$. There is $\epsilon : \Qp$ such that $m + \epsilon < f(x)$. But now merely for
  some $y_i$ participating in the definition of $x_\epsilon$ we get $|f(x) - f(y_i) <
  \epsilon$, therefore $m < f(x) - \epsilon < f(y_i) \leq m$, a contradiction.

  We finish the proof by showing that $m$ satisfies the second part of the theorem, because
  it is then automatically a least upper bound. Given any $\epsilon : \Qp$, on one hand
  $|m - f(x_{\epsilon/2})| < 3 \epsilon/4$, and on the other $|f(x_{\epsilon/2}) - f(y_i)| <
  \epsilon/4$ merely for some $y_i$ participating in the definition of $x_{\epsilon/2}$,
  therefore by taking $u \defeq y_i$ we obtain $|m - f(u)| < \epsilon$ by triangle
  inequality.
\end{proof}

Now, if in \autoref{ctb-uniformly-continuous-sup} we also knew that $M$ were complete, we
could hope to weaken the assumption of uniform continuity to continuity, and strengthen
the conclusion to existence of a point at which the supremum is attained. The usual proofs
of these improvements rely on the the facts that in a complete totally bounded space
\begin{enumerate}
\item continuity implies uniform continuity, and
\item every sequence has a convergent subsequence.
\end{enumerate}
The first statement follows easily from Heine-Borel compactness, and the second is just
Bolzano--Weierstra\ss{} compactness.
\index{compactness!Bolzano--Weierstrass@Bolzano--Weierstra\ss{}}%
Unfortunately, these are both somewhat problematic. Let
us first show that Bolzano--Weierstra\ss{} compactness implies an instance of excluded middle
known as the \define{limited principle of omniscience}:
\indexsee{axiom!limited principle of omniscience}{limited principle of omniscience}%
\indexdef{limited principle of omniscience}%
for every $\alpha : \N \to \bool$,
\begin{equation} \label{eq:lpo}
  \Parens{\sm{n : \N} \alpha(n) = \btrue} +
  \Parens{\prd{n : \N} \alpha(n) = \bfalse}.
\end{equation}
Computationally speaking, we would not expect this principle to hold, because it asks us to decide
whether infinitely many values of a function are~$\bfalse$.
  
\begin{thm} \label{analysis-bw-lpo}
  Bolzano--Weierstra\ss{} compactness of $[0,1]$ implies the limited principle of omniscience.
  \index{compactness!Bolzano--Weierstrass@Bolzano--Weierstra\ss{}}%
\end{thm}

\begin{proof}
  Given any $\alpha : \N \to \bool$, define the sequence\index{sequence} $x : \N \to [0,1]$ by
  \begin{equation*}
    x_n \defeq
    \begin{cases}
      0 & \text{if $\alpha(k) = \bfalse$ for all $k < n$,}\\
      1 & \text{if $\alpha(k) = \btrue$ for some $k < n$}.
    \end{cases}
  \end{equation*}
  If the Bolzano--Weierstra\ss{} property holds, there exists a strictly increasing $f : \N \to
  \N$ such that $x \circ f$ is a Cauchy sequence\index{Cauchy!sequence}. For a sufficiently large $n :
  \N$ the $n$-th term $x_{f(n)}$ is within $1/6$ of its limit. Either $x_{f(n)} < 2/3$ or
  $x_{f(n)} > 1/3$. If $x_{f(n)} < 2/3$ then~$x_n$ converges to $0$ and so $\prd{n : \N}
  \alpha(n) = \bfalse$. If $x_{f(n)} > 1/3$ then $x_{f(n)} = 1$, therefore $\sm{n : \N}
  \alpha(n) = \btrue$.
\end{proof}

While we might not mourn Bolzano--Weierstra\ss{} compactness too much, it seems harder to live
without Heine--Borel compactness, as attested by the fact that both classical mathematics
and Brouwer's Intuitionism accepted it. As we do not want to wade too deeply into general
topology, we shall work with basic open sets. In the case of $\R$ these are the open
intervals with rational endpoints. A family of such intervals, indexed by a type~$I$,
would be a map
\begin{equation*}
  \mathcal{F} : I \to \setof{(q, r) : \Q \times \Q | q < r},
\end{equation*}
with the idea that a pair of rationals $(q, r)$ with $q < r$ determines the type $\setof{ x : \R | q < x < r}$. It is slightly more convenient to allow degenerate intervals as well, so we take a
\define{family of basic intervals}
\indexdef{family!of basic intervals}%
\indexdef{interval!family of basic}%
to be a map
\begin{equation*}
  \mathcal{F} : I \to \Q \times \Q.
\end{equation*}
To be quite precise, a family is a dependent pair $(I, \mathcal{F})$, not just
$\mathcal{F}$. A \define{finite family of basic intervals} is one indexed by $\setof{ m :
  \N | m < n}$ for some $n : \N$. We usually present it by a finite list $[(q_0, r_0), \ldots,
(q_{n-1}, r_{n-1})]$. Finally, a \define{finite subfamily}\indexdef{subfamily, finite, of intervals} of $(I, \mathcal{F})$ is given
by a list of indices $[i_1, \ldots, i_n]$ which then determine the finite family
$[\mathcal{F}(i_1), \ldots, \mathcal{F}(i_n)]$.

As long as we are aware of the distinction between a pair $(q, r)$ and the corresponding
interval $\setof{ x : \R | q < x < r}$, we may safely use the same notation $(q, r)$ for
both. Intersections\indexdef{intersection!of intervals} and inclusions\indexdef{inclusion!of intervals}\indexdef{containment!of intervals} of intervals are expressible in terms of their
endpoints:
\symlabel{interval-intersection}
\symlabel{interval-subset}
\begin{align*}
  (q, r) \cap (s, t) &\ \defeq\  (\max(q, s), \min(r, t)),\\
  (q, r) \subseteq (s, t) &\ \defeq\ (q < r \Rightarrow s \leq q < r \leq t).
\end{align*}
We say that $\intfam{i}{I}{(q_i, r_i)}$ \define{(pointwise) covers $[a,b]$}
\indexdef{interval!pointwise cover}%
\indexdef{cover!pointwise}%
\indexdef{pointwise!cover}%
when
\begin{equation} \label{eq:cover-pointwise-truncated}
  \fall{x : [a,b]} \exis{i : I} q_i < x < r_i.
\end{equation}
The \define{Heine-Borel compactness for $[0,1]$}
\indexdef{compactness!Heine-Borel}%
states that every covering family of $[0,1]$
merely has a finite subfamily which still covers $[0,1]$.

\index{depression}
\begin{thm} \label{classical-Heine-Borel}
  \index{excluded middle}%
  If excluded middle holds then $[0,1]$ is Heine-Borel compact.
\end{thm}

\begin{proof}
  Assume for the purpose of reaching a contradiction that a family $\intfam{i}{I}{(a_i,
    b_i)}$ covers $[0,1]$ but no finite subfamily does. We construct a sequence of closed
  intervals $[q_n, r_n]$ which are nested, their sizes shrink to~$0$, and none of them is covered
  by a finite subfamily of $\intfam{i}{I}{(a_i, b_i)}$.

  We set $[q_0, r_0] \defeq [0,1]$. Assuming $[q_n, r_n]$ has been constructed, let $s
  \defeq (2 q_n + r_n)/3$ and $t \defeq (q_n + 2 r_n)/3$. Both $[q_n, t]$ and $[s, r_n]$
  are covered by $\intfam{i}{I}{(a_i, b_i)}$, but they cannot both have a finite subcover,
  or else so would $[q_n, r_n]$. Either $[q_n, t]$ has a finite subcover or it does not.
  If it does we set $[q_{n+1}, r_{n+1}] \defeq [s, r_n]$, otherwise we set $[q_{n+1},
  r_{n+1}] \defeq [q_n, t]$.

  The sequences $q_0, q_1, \ldots$ and $r_0, r_1, \ldots$ are both Cauchy and they
  converge to a point $x : [0,1]$ which is contained in every $[q_n, r_n]$.
  There merely exists $i : I$ such that $a_i < x < b_i$. Because the sizes of the
  intervals $[q_n, r_n]$ shrink to zero, there is $n : \N$ such that $a_i < q_n \leq x
  \leq r_n < b_i$, but this means that $[q_n, r_n]$ is covered by a single interval $(a_i,
  b_i)$, while at the same time it has no finite subcover. A contradiction.
\end{proof}

Without excluded middle, or a pinch of Brouwerian Intuitionism, we seem to be stuck.
Nevertheless, Heine-Borel compactness of $[0,1]$ \emph{can} be recovered in a constructive
setting, in a fashion that is still compatible with classical mathematics! For this to be
done, we need to revisit the notion of cover. The trouble with
\eqref{eq:cover-pointwise-truncated} is that the truncated existential allows a space to
be covered in any haphazard way, and so computationally speaking, we stand no chance of
merely extracting a finite subcover. By removing the truncation we get
\begin{equation} \label{eq:cover-pointwise}
  \prd{x : [0,1]} \sm{i : I} q_i < x < r_i,
\end{equation}
which might help, were it not too demanding of covers. With this definition we
could not even show that $(0,3)$ and $(2,5)$ cover $[1,4]$ because that would amount
to exhibiting a non-constant map $[1,4] \to \bool$, see
\autoref{ex:reals-non-constant-into-Z}.  Here we can take a lesson from ``pointfree topology''
\index{pointfree topology}%
\index{topology!pointfree}%
(i.e.\ locale theory):
\index{locale}%
the notion of cover ought to be expressed in terms of open sets, without
reference to points. Such a ``holistic'' view of space will then allow us to analyze the
notion of cover, and we shall be able to recover Heine-Borel compactness.  Locale
theory uses power sets,
\index{power set}%
which we could obtain by assuming propositional resizing;
\index{propositional!resizing}%
but instead we can steal ideas from the predicative cousin of locale theory,
\index{mathematics!predicative}%
which is called ``formal topology''.
\index{formal!topology}%

\index{acceptance|(}

Suppose that we have a family $\pairr{I, \mathcal{F}}$ and an interval $(a, b)$. How might
we express the fact that $(a,b)$ is covered by the family, without referring to points?
Here is one: if $(a, b)$ equals some $\mathcal{F}(i)$ then it is covered by the family.
And another one: if $(a,b)$ is covered by some other family $(J, \mathcal{G})$, and in
turn each $\mathcal{G}(j)$ is covered by $\pairr{I, \mathcal{F}}$, then $(a,b)$ is covered
$\pairr{I, \mathcal{F}}$. Notice that we are listing \emph{rules} which can be used to
\emph{deduce} that $\pairr{I, \mathcal{F}}$ covers $(a,b)$. We should find sufficiently
good rules and turn them into an inductive definition.

\begin{defn} \label{defn:inductive-cover}
  The \define{inductive cover $\cover$}
  \indexdef{inductive!cover}%
  \indexdef{cover!inductive}%
  is a mere relation
  \begin{equation*}
    {\cover} : (\Q \times \Q) \to \Parens{\sm{I : \type} (I \to \Q \times \Q)} \to \prop
  \end{equation*}
  defined inductively by the following rules, where $q, r, s, t$ are rational numbers and
  $\pairr{I, \mathcal{F}}$, $\pairr{J, \mathcal{G}}$ are families of basic intervals:
  \begin{enumerate}

  \item \emph{reflexivity:}
    \index{reflexivity!of inductive cover}%
    $\mathcal{F}(i) \cover \pairr{I, \mathcal{F}}$ for all $i : I$,
      
  \item \emph{transitivity:}
    \index{transitivity!of inductive cover}%
    if $(q, r) \cover \pairr{J, \mathcal{G}}$ and $\fall{j : J} \mathcal{G}(j) \cover \pairr{I,\mathcal{F}}$
    then $(q, r) \cover \pairr{I, \mathcal{F}}$,

  \item \emph{monotonicity:}
    \index{monotonicity!of inductive cover}%
    if $(q, r) \subseteq (s, t)$ and $(s,t) \cover \pairr{I, \mathcal{F}}$ then $(q, r) \cover
    \pairr{I, \mathcal{F}}$,

  \item \emph{localization:}
    \index{localization of inductive cover}%
    if $(q, r) \cover (I, \mathcal{F})$ then $(q, r) \cap (s, t) \cover
    \intfam{i}{I}{(\mathcal{F}(i) \cap (s, t))}$.

  \item \label{defn:inductive-cover-interval-1}
    if $q < s < t < r$ then $(q, r) \cover [(q, t), (r, s)]$,

  \item \label{defn:inductive-cover-interval-2}
    $(q, r) \cover \intfam{u}{\setof{ (s,t) : \Q \times \Q | q < s < t < r}}{u}$.
  \end{enumerate}
\end{defn}

The definition should be read as a higher-inductive type in which the listed rules are
point constructors, and the type is $(-1)$-truncated. The first four clauses are of a
general nature and should be intuitively clear. The last two clauses are specific to the
real line: one says that an interval may be covered by two intervals if they overlap,
while the other one says that an interval may be covered from within. Incidentally, if $r
\leq q$ then $(q, r)$ is covered by the empty family by the last clause.

Inductive covers enjoy the Heine-Borel property, the proof of which requires a lemma.

\begin{lem} \label{reals-formal-topology-locally-compact}
  Suppose $q < s < t < r$ and $(q, r) \cover \pairr{I, \mathcal{F}}$. Then there merely
  exists a finite subfamily of $\pairr{I, \mathcal{F}}$ which inductively covers $(s, t)$.
\end{lem}

\begin{proof}
  We prove the statement by induction on $(q, r) \cover \pairr{I, \mathcal{F}}$. There are
  six cases:
  \begin{enumerate}

  \item Reflexivity: if $(q, r) = \mathcal{F}(i)$ then by monotonicity $(s, t)$ is covered
    by the finite subfamily $[\mathcal{F}(i)]$.

  \item Transitivity:
    suppose $(q, r) \cover \pairr{J, \mathcal{G}}$ and $\fall{j : J} \mathcal{G}(j) \cover
    \pairr{I, \mathcal{F}}$. By the inductive hypothesis there merely exists
    $[\mathcal{G}(j_1), \ldots, \mathcal{G}(j_n)]$ which covers $(s, t)$.
    Again by the inductive hypothesis, each of $\mathcal{G}(j_k)$ is covered by a finite
    subfamily of $\pairr{I, \mathcal{F}}$, and we can collect these into a finite
    subfamily which covers $(s, t)$.

  \item Monotonicity:
    if $(q, r) \subseteq (u, v)$ and $(u, v) \cover \pairr{I, \mathcal{F}}$ then we may
    apply the inductive hypothesis to $(u, v) \cover \pairr{I, \mathcal{F}}$ because $u <
    s < t < v$.

  \item Localization:
    suppose $(q', r') \cover \pairr{I, \mathcal{F}}$ and $(q, r) = (q', r') \cap (a, b)$.
    Because $q' < s < t < r'$, by the inductive hypothesis there is a finite subcover
    $[\mathcal{F}(i_1), \ldots, \mathcal{F}(i_n)]$ of $(s, t)$. We also know that $a < s <
    t < b$, therefore $(s, t) = (s, t) \cap (a, b)$ is covered by
    $[\mathcal{F}(i_1) \cap (a,b), \ldots, \mathcal{F}(i_n) \cap (a,b)]$, which is a
    finite subfamily of $\intfam{i}{I}{(\mathcal{F}(i) \cap (a, b))}$.

  \item If $(q, r) \cover [(q, v), (u, r)]$ for some $q < u < v < r$ then by monotonicity
    $(s, t) \cover [(q, v), (u, r)]$.

  \item Finally, $(s, t) \cover \intfam{z}{\setof{ (u,v):\Q \times \Q | q < u < v < r}}{z}$ by
    reflexivity. \qedhere
  \end{enumerate}
\end{proof}

Say that \define{$\pairr{I, \mathcal{F}}$ inductively covers
  $[a, b]$} when there merely exists $\epsilon : \Qp$ such that $(a - \epsilon, b +
\epsilon) \cover \pairr{I, \mathcal{F}}$.

\begin{cor} \label{interval-Heine-Borel}
  \index{compactness!Heine-Borel}%
  \index{interval!open and closed}%
  A closed interval is Heine-Borel compact for inductive covers.
\end{cor}

\begin{proof}
  Suppose $[a, b]$ is inductively covered by $\pairr{I, \mathcal{F}}$, so there merely is
  $\epsilon : \Qp$ such that $(a - \epsilon, b + \epsilon) \cover \pairr{I, \mathcal{F}}$.
  By \autoref{reals-formal-topology-locally-compact} there is a finite subcover of
  $(a - \epsilon/2, b + \epsilon/2)$, which is therefore a finite subcover of $[a, b]$.
\end{proof}

Experience from formal topology\index{topology!formal} shows that the rules for inductive covers are sufficient
for a constructive development of pointfree topology. But we can also provide our own
evidence that they are a reasonable notion.

\begin{thm} \label{inductive-cover-classical}
  \mbox{}
  \begin{enumerate}
  \item An inductive cover is also a pointwise cover.
  \item Assuming excluded middle, a pointwise cover is also an inductive cover.
  \end{enumerate}
\end{thm}

\begin{proof}
  \mbox{}
  \begin{enumerate}

  \item 
    Consider a family of basic intervals $\pairr{I, \mathcal{F}}$, where we write $(q_i,
    r_i) \defeq \mathcal{F}(i)$, an interval $(a,b)$ inductively covered by $\pairr{I,
      \mathcal{F}}$, and $x$ such that $a < x < b$.
    We prove by induction on $(a,b) \cover \pairr{I, \mathcal{F}}$ that there merely
    exists $i : I$ such that $q_i < x < r_i$. Most cases are pretty obvious, so we show
    just two. If $(a,b) \cover \pairr{I, \mathcal{F}}$ by reflexivity, then there merely
    is some $i : I$ such that $(a,b) = (q_i, r_i)$ and so $q_i < x < r_i$. If $(a,b)
    \cover \pairr{I, \mathcal{F}}$ by transitivity via $\intfam{j}{J}{(s_j, t_j)}$ then by
    the inductive hypothesis there merely is $j : J$ such that $s_j < x < t_j$, and then since
    $(s_j, t_j) \cover \pairr{I, \mathcal{F}}$ again by the inductive hypothesis there merely
    exists $i : I$ such that $q_i < x < r_i$. Other cases are just as exciting.

  \item Suppose $\intfam{i}{I}{(q_i, r_i)}$ pointwise covers $(a, b)$. By
    \autoref{defn:inductive-cover-interval-2} of \autoref{defn:inductive-cover} it
    suffices to show that $\intfam{i}{I}{(q_i, r_i)}$ inductively covers $(c, d)$ whenever
    $a < c < d < b$, so consider such $c$ and $d$. By \autoref{classical-Heine-Borel}
    there is a finite subfamily $[i_1, \ldots, i_n]$ which already pointwise covers $[c,
    d]$, and hence $(c,d)$. Let $\epsilon : \Qp$ be a Lebesgue number
    \index{Lebesgue number}
    for $(q_{i_1}, r_{i_1}), \ldots, (q_{i_n}, r_{i_n})$ as in
    \autoref{ex:finite-cover-lebesgue-number}. There is a positive $k : \N$ such that $2 (d - c)/k
    < \min(1, \epsilon)$. For $0 \leq i \leq k$ let
    \begin{equation*}
      c_k \defeq ((k - i) c + i d) / k.
    \end{equation*}
    The intervals $(c_0, c_2)$, $(c_1, c_3)$, \dots, $(c_{k-2}, c_k)$ inductively cover
    $(c,d)$ by repeated use of transitivity and~\autoref{defn:inductive-cover-interval-1}
    in \autoref{defn:inductive-cover}. Because their widths are below $\epsilon$ each of
    them is contained in some $(q_i, r_i)$, and we may use transitivity and monotonicity to
    conclude that $\intfam{i}{I}{(q_i, r_i)}$ inductively cover $(c, d)$. \qedhere
  \end{enumerate}
\end{proof}

The upshot of the previous theorem is that, as far as classical mathematics is concerned,
there is no difference between a pointwise and an inductive cover. In particular, since it
is consistent to assume excluded middle in homotopy type theory, we cannot exhibit an
inductive cover which fails to be a pointwise cover. Or to put it in a different way, the
difference between pointwise and inductive covers is not what they cover but in the
\emph{proofs} that they cover. 

We could write another book by going on like this, but let us stop here and hope that we
have provided ample justification for the claim that analysis can be developed in homotopy
type theory. The curious reader should consult \autoref{ex:mean-value-theorem} for
constructive versions of the mean value theorem.

\index{acceptance|)}

\index{mathematics!classical|)}%
\index{mathematics!constructive|)}%

\section{The surreal numbers}
\label{sec:surreals}

\index{surreal numbers|(}%

In this section we consider another example of a higher inductive-in\-duc\-tive type, which draws together many of our threads: Conway's field \NO of \emph{surreal numbers}~\cite{conway:onag}.
The surreal numbers are the natural common generalization of the (Dedekind) real numbers (\autoref{sec:dedekind-reals}) and the ordinal numbers (\autoref{sec:ordinals}).
Conway, working in classical\index{mathematics!classical} mathematics with excluded middle and Choice, defines a surreal number to be a pair of \emph{sets} of surreal numbers, written $\surr L R$, such that every element of $L$ is strictly less than every element of $R$.
This obviously looks like an inductive definition, but there are three issues with regarding it as such.

Firstly, the definition requires the relation of (strict) inequality between surreals, so that relation must be defined simultaneously with the type \NO of surreals.
(Conway avoids this issue by first defining \emph{games}\index{game!Conway}, which are like surreals but omit the compatibility condition on $L$ and $R$.)
As with the relation $\closesym$ for the Cauchy reals, this simultaneous definition could \emph{a priori} be either inductive-inductive or inductive-recursive.
We will choose to make it inductive-inductive, for the same reasons we made that choice for $\closesym$.

Moreover, we will define strict inequality $<$ and non-strict inequality $\le$ for surreals separately (and mutually inductively).
Conway defines $<$ in terms of $\le$, in a way which is sensible classically but not constructively.
\index{mathematics!constructive}%
Furthermore, a negative definition of $<$ would make it unacceptable as a hypothesis of the constructor of a higher inductive type (see \autoref{sec:strictly-positive}).

Secondly, Conway says that $L$ and $R$ in $\surr L R$ should be ``sets of surreal numbers'', but the naive meaning of this as a predicate $\NO\to\prop$ is not positive, hence cannot be used as input to an inductive constructor.
However, this would not be a good type-theoretic translation of what Conway means anyway, because in set theory the surreal numbers form a proper class, whereas the sets $L$ and $R$ are true (small) sets, not arbitrary subclasses of \NO.
In type theory, this means that \NO will be defined relative to a universe \UU, but will itself belong to the next higher universe $\UU'$, like the sets \ord and \card of ordinals and cardinals, the cumulative hierarchy $V$, or even the Dedekind reals in the absence of propositional resizing.
\index{propositional!resizing}%
We will then require the ``sets'' $L$ and $R$ of surreals to be \UU-small, and so it is natural to represent them by \emph{families} of surreals indexed by some \UU-small type.
(This is all exactly the same as what we did with the cumulative hierarchy in \autoref{sec:cumulative-hierarchy}.)
That is, the constructor of surreals will have type
\[ \prd{\LL,\RR:\UU} (\LL\to\NO) \to (\RR\to \NO) \to (\text{some condition}) \to \NO \]
which is indeed strictly positive.\index{strict!positivity}

Finally, after giving the mutual definitions of \NO and its ordering, Conway declares two surreal numbers $x$ and $y$ to be \emph{equal} if $x\le y$ and $y\le x$.
This is naturally read as passing to a quotient of the set of ``pre-surreals'' by an equivalence relation.
However, in the absence of the axiom of choice, such a quotient presents the same problem as the quotient in the usual construction of Cauchy reals: it will no longer be the case that a pair of families \emph{of surreals} yield a new surreal $\surr L R$, since we cannot necessarily ``lift'' $L$ and $R$ to families of pre-surreals.
Of course, we can solve this problem in the same way we did for Cauchy reals, by using a \emph{higher} inductive-inductive definition.

\begin{defn}\label{defn:surreals}
  The type \NO of \define{surreal numbers},
  \indexdef{surreal numbers}%
  \indexsee{number!surreal}{surreal numbers}%
  along with the relations $\mathord<:\NO\to\NO\to\type$ and $\mathord\le:\NO\to\NO\to\type$, are defined higher inductive-inductively as follows.
  The type \NO has the following constructors.
  \begin{itemize}
  \item For any $\LL,\RR:\UU$ and functions $\LL\to \NO$ and $\RR\to \NO$, whose values we write as $x^L$ and $x^R$ for $L:\LL$ and $R:\RR$ respectively, if $\fall{L:\LL}{R:\RR} x^L<x^R$, then there is a surreal number $x$.
  \item For any $x,y:\NO$ such that $x\le y$ and $y\le x$, we have $\noeq(x,y):x=y$.
  \end{itemize}
  We will refer to the inputs of the first constructor as a \define{cut}.
  \indexdef{cut!of surreal numbers}%
  If $x$ is the surreal number constructed from a cut, then the notation $x^L$ will implicitly assume $L:\LL$, and similarly $x^R$ will assume $R:\RR$.
  In this way we can usually avoid naming the indexing types $\LL$ and $\RR$, which is convenient when there are many different cuts under discussion.
  Following Conway, we call $x^L$ a \emph{left option}\indexdef{option of a surreal number} of $x$ and $x^R$ a \emph{right option}.

  The path constructor implies that different cuts can define the same surreal number.
  Thus, it does not make sense to speak of the left or right options of an arbitrary surreal number $x$, unless we also know that $x$ is defined by a particular cut.
  Thus in what follows we will say, for instance, ``given a cut defining a surreal number $x$'' in contrast to ``given a surreal number $x$''.

  The relation $\le$ has the following constructors.
  \index{non-strict order}%
  \index{order!non-strict}%
  \begin{itemize}
  \item Given cuts defining two surreal numbers $x$ and $y$, if $x^L<y$ for all $L$, and $x<y^R$ for all $R$, then $x\le y$.
  \item Propositional truncation:
    for any $x,y:\NO$, if $p,q:x\le y$, then $p=q$.
  \end{itemize}
  And the relation $<$ has the following constructors.
  \index{strict!order}%
  \index{order!strict}%
  \begin{itemize}
  \item Given cuts defining two surreal numbers $x$ and $y$, if there is an $L$ such that $x\le y^L$, then $x<y$.
  \item Given cuts defining two surreal numbers $x$ and $y$, if there is an $R$ such that $x^R\le y$, then $x<y$.
  \item Propositional truncation: for any $x,y:\NO$, if $p,q:x<y$, then $p=q$.
  \end{itemize}
\end{defn}

\noindent
We compare this with Conway's definitions:
\begin{itemize}\footnotesize
\item[-] If $L,R$ are any two sets of numbers, and no member of $L$ is $\ge$ any member of $R$, then there is a number $\surr L R$.
  All numbers are constructed in this way.
\item[-] $x\ge y$ iff (no $x^R\le y$ and $x\le$ no $y^L$).
\item[-] $x=y$ iff ($x \ge y$ and $y\ge x$).
\item[-] $x>y$ iff ($x\ge y$ and $y\not\ge x$).
\end{itemize}
The inclusion of $x\ge y$ in the definition of $x>y$ is unnecessary if all objects are [surreal] numbers rather than ``games''\index{game!Conway}.
Thus, Conway's $<$ is just the negation of his $\ge$, so that his condition for $\surr L R$ to be a surreal is the same as ours.
Negating Conway's $\le$ and canceling double negations, we arrive at our definition of $<$, and we can then reformulate his $\le$ in terms of $<$ without negations.

We can immediately populate $\NO$ with many surreal numbers.
Like Conway, we write
\symlabel{surreal-cut}
\[\surr{x,y,z,\dots}{u,v,w,\dots}\]
for the surreal number defined by a cut where $\LL\to\NO$ and $\RR\to\NO$ are families described by $x,y,z,\dots$ and $u,v,w,\dots$.
Of course, if $\LL$ or $\RR$ are $\emptyt$, we leave the corresponding part of the notation empty.
There is an unfortunate clash with the standard notation $\setof{x:A | P(x)}$ for subsets, but we will not use the latter in this section.
\begin{itemize}
\item We define $\iota_{\nat}:\nat\to\NO$ recursively by
  \begin{align*}
    \iota_{\nat}(0) &\defeq \surr{}{},\\
    \iota_\nat(\suc(n)) &\defeq \surr{\iota_\nat(n)}{}.
  \end{align*}
  That is, $\iota_\nat(0)$ is defined by the cut consisting of $\emptyt\to\NO$ and $\emptyt\to\NO$.
  Similarly, $\iota_\nat(\suc(n))$ is defined by $\unit\to\NO$ (picking out $\iota_\nat(n)$) and $\emptyt\to\NO$.
\item Similarly, we define $\iota_{\Z}:\Z\to\NO$ using the sign-case recursion principle (\autoref{thm:sign-induction}):
  \begin{align*}
    \iota_{\Z}(0) &\defeq \surr{}{},\\
    \iota_\Z(n+1) &\defeq \surr{\iota_\Z(n)}{} & &\text{$n\ge 0$,}\\
    \iota_\Z(n-1) &\defeq \surr{}{\iota_\Z(n)} & &\text{$n\le 0$.}
  \end{align*}
\item By a \define{dyadic rational}
  \indexdef{rational numbers!dyadic}%
  \indexsee{dyadic rational}{rational numbers, dyadic}%
  we mean a pair $(a,n)$ where $a:\Z$ and $n:\nat$, and such that if $n>0$ then $a$ is odd.
  We will write it as $a/2^n$, and identify it with the corresponding rational number.
  If $\Q_D$ denotes the set of dyadic rationals, we define $\iota_{\Q_D}:\Q_D\to\NO$ by induction on $n$:
  \begin{align*}
    \iota_{\Q_D}(a/2^0) &\defeq \iota_\Z(a),\\
    \iota_{\Q_D}(a/2^n) &\defeq \surr{a/2^n - 1/2^n}{a/2^n + 1/2^n},
    \quad \text{for $n>0$.}
  \end{align*}
  Here we use the fact that if $n>0$ and $a$ is odd, then $a/2^n \pm 1/2^n$ is a dyadic rational with a smaller denominator than $a/2^n$.
\item We define $\iota_{\RD}:\RD\to\NO$, where $\RD$ is (any version of) the Dedekind reals from \autoref{sec:dedekind-reals}, by
  \begin{align*}
    \iota_{\RD}(x) &\defeq
    \surr{q\in\Q_D \text{ such that } q<x}{q\in\Q_D \text{ such that } x<q}.
  \end{align*}
  Unlike in the previous cases, it is not obvious that this extends $\iota_{\Q_D}$ when we regard dyadic rationals as Dedekind reals.
  This follows from the simplicity theorem (\autoref{thm:NO-simplicity}).
\item Recall the type \ord of \emph{ordinals}\index{ordinal} from \autoref{sec:ordinals}, which is well-ordered by the relation $<$, where $A<B$ means that $A = \ordsl B b$ for some $b:B$.
  We define $\iota_{\ord}:\ord\to\NO$ by well-founded recursion (\autoref{thm:wfrec}) on $\ord$:
  \begin{equation*}
    \iota_{\ord}(A) \defeq
    \surr{\iota_\ord(\ordsl A a) \text{ for all } a:A}{}.
  \end{equation*}
  It will also follow from the simplicity theorem that $\iota_\ord$ restricted to finite ordinals agrees with $\iota_\nat$.
\item A few more interesting examples taken from Conway:
  \begin{align*}
    \omega &\defeq \surr{0,1,2,3,\dots}{} \qquad\text{(also an ordinal)}\\
    -\omega &\defeq \surr{}{\dots,-3,-2,-1,0}\\
    1/\omega &\defeq \textstyle\surr{0}{1,\frac12,\frac14,\frac18,\dots}\\
    \omega-1 &\defeq \surr{0,1,2,3,\dots}{\omega}\\
    \omega/2 &\defeq \surr{0,1,2,3,\dots}{\dots,\omega-2,\omega-1,\omega}.
  \end{align*}
\end{itemize}

In identifying surreal numbers presented by different cuts, the following simple observation is useful.

\begin{thm}[Conway's simplicity theorem]\label{thm:NO-simplicity}
  \index{simplicity theorem}%
  \index{theorem!Conway's simplicity}%
  Suppose $x$ and $z$ are surreal numbers defined by cuts, and that the following hold.
  \begin{itemize}
  \item $x^L < z < x^R$ for all $L$ and $R$.
  \item For every left option $z^L$ of $z$, there exists a left option $x^{L'}$ with $z^L\le x^{L'}$.
  \item For every right option $z^R$ of $z$, there exists a right option $x^{R'}$ with $x^{R'}\le z^R$.
  \end{itemize}
  Then $x=z$.
\end{thm}
\begin{proof}
  Applying the path constructor of $\NO$, we must show $x\le z$ and $z\le x$.
  The first entails showing $x^L<z$ for all $L$, which we assumed, and $x<z^R$ for all $R$.
  But by assumption, for any $z^R$ there is an $x^{R'}$ with $x^{R'}\le z^R$ hence $x<z^R$ as desired.
  Thus $x\le z$; the proof of $z\le x$ is symmetric.
\end{proof}

\index{induction principle!for surreal numbers}
In order to say much more about surreal numbers, however, we need their induction principle.
The mutual induction principle for $(\NO,\le,<)$ applies to three families of types:
\begin{align*}
  A &: \NO\to\type\\
  B &: \prd{x,y:\NO}{a:A(x)}{b:A(y)} (x\le y) \to \type\\
  C &: \prd{x,y:\NO}{a:A(x)}{b:A(y)} (x<y) \to \type.
\end{align*}
As with the induction principle for Cauchy reals, it is helpful to think of $B$ and $C$ as families of relations between the types $A(x)$ and $A(y)$.
\symlabel{NO-recursion}
Thus we write $B(x,y,a,b,\xi)$ as $(x,a) \ble^\xi (y,b)$ and $C(x,y,a,b,\xi)$ as $(x,a) \blt^\xi (y,b)$.
Similarly, we usually omit the $\xi$ since it inhabits a mere proposition and so is uninteresting, and we may often omit $x$ and $y$ as well, writing simply $a\ble b$ or $a\blt b$.
With these notations, the hypotheses of the induction principle are the following.
\begin{itemize}
\item For any cut defining a surreal number $x$, together with
  \begin{enumerate}
  \item for each $L$, an element $a^L:A(x^L)$, and
  \item for each $R$, an element $a^R:A(x^R)$, such that
  \item for all $L$ and $R$ we have $(x^L,a^L) \blt (x^R,a^R)$
  \end{enumerate}
  there is a specified element $f_a:A(x)$.
  We call such data a \define{dependent cut}
  \indexdef{cut!of surreal numbers!dependent}%
  \indexdef{dependent!cut}%
  over the cut defining~$x$.
\item For any $x,y:\NO$ with $a:A(x)$ and $b:A(y)$, if $x\le y$ and $y\le x$ and also $(x,a) \ble (y,b)$
  and $(y,b) \ble (x,a)$,
  then $\dpath{A}{\noeq}{a}{b}$.
\item Given cuts defining two surreal numbers $x$ and $y$, and dependent cuts $a$ over $x$ and $b$ over $y$, such that for all $L$ we have $x^L<y$ and $(x^L,a^L)\blt (y,f_b)$,
  and for all $R$ we have $x<y^R$ and $(x,f_a) \blt (y^R,b^R)$,
  then $(x,f_a) \ble (y,f_b)$.
\item $\ble$ takes values in mere propositions.
\item Given cuts defining two surreal numbers $x$ and $y$, dependent cuts $a$ over $x$ and $b$ over $y$, and an $L_0$ such that $x\le y^{L_0}$ and $(x,f_a) \ble (y^{L_0},b^{L_0})$,
  we have $(x,f_a) \blt (y,f_b)$.
\item Given cuts defining two surreal numbers $x$ and $y$, dependent cuts $a$ over $x$ and $b$ over $y$, and an ${R_0}$ such that $x^{R_0}\le y$ together with $(x^{R_0},a^{R_0}),\ble (y,f_b)$,
  we have $(x,f_a) \blt (y,f_b)$.
\item $\blt$ takes values in mere propositions.
\end{itemize}
Under these hypotheses we deduce a function $f:\prd{x:\NO} A(x)$ such that
\begin{align}
  f(x) &\;\jdeq\; f_{f[x]} \label{eq:noind1}\\
  (x\le y) &\;\Rightarrow\; (x,f(x)) \ble (y,f(y)) \notag\\
  (x< y) &\;\Rightarrow\; (x,f(x)) \blt (y,f(y)). \notag
\end{align}
In the computation rule~\eqref{eq:noind1} for the point constructor, $x$ is a surreal number defined by a cut, and $f[x]$ denotes the dependent cut over $x$ defined by applying $f$ (and using the fact that $f$ takes $<$ to $\blt$).
As usual, we will generally use pattern-matching notation, where the definition of $f$ on a cut $\surr{x^L}{x^R}$ may use the symbols $f(x^L)$ and $f(x^R)$ and the assumption that they form a dependent cut.

As with the Cauchy reals, we have special cases resulting from trivializing some of $A$, $\ble$, and~$\blt$.
Taking $\ble$ and $\blt$ to be constant at \unit, we have \define{\NO-induction}, which for simplicity we state only for mere properties:
\begin{itemize}
\item Given $P:\NO\to\prop$, if $P(x)$ holds whenever $x$ is a surreal number defined by a cut such that $P(x^L)$ and $P(x^R)$ hold for all
$L$ and $R$, then $P(x)$ holds for all $x:\NO$.
\end{itemize}
This should be compared with Conway's remark:
\begin{quote}\footnotesize
  In general when we wish to establish a proposition $P(x)$ for all numbers $x$, we will prove it inductively by deducing $P(x)$ from the truth of all the propositions $P(x^L)$ and $P(x^R)$.
  We regard the phrase ``all numbers are constructed in this way'' as justifying the legitimacy of this procedure.
\end{quote}
With $\NO$-induction, we can prove

\begin{thm}[Conway's Theorem 0]\label{thm:NO-refl-opt}\ 
  \index{theorem!Conway's 0}%
  \begin{enumerate}
  \item For any $x:\NO$, we have $x\le x$.\label{item:NO-le-refl}
  \item For any $x:\NO$ defined by a cut, we have $x^L <x$ and $x<x^R$ for all $L$ and $R$.\label{item:NO-lt-opt}
  \end{enumerate}
\end{thm}
\begin{proof}
  Note first that if $x\le x$, then whenever $x$ occurs as a left option of some cut $y$, we have $x<y$ by the first constructor of $<$, and similarly whenever $x$ occurs as a right option of a cut $y$, we have $y<x$ by the second constructor of $<$.
  In particular,~\ref{item:NO-le-refl}$\Rightarrow$\ref{item:NO-lt-opt}.

  We prove~\ref{item:NO-le-refl} by $\NO$-induction on $x$.
  Thus, assume $x$ is defined by a cut such that $x^L\le x^L$ and $x^R \le x^R$ for all $L$ and $R$.
  But by our observation above, these assumptions imply $x^L<x$ and $x<x^R$ for all $L$ and $R$, yielding $x\le x$ by the constructor of $\le$.
\end{proof}

\begin{cor}\label{thm:NO-set}
  \NO is a 0-type.
\end{cor}
\begin{proof}
  The mere relation $R(x,y)\defeq (x\le y) \land (y\le x)$ implies identity by the path constructor of $\NO$, and contains the diagonal by \autoref{thm:NO-refl-opt}\ref{item:NO-le-refl}.
  Thus, \autoref{thm:h-set-refrel-in-paths-sets} applies.
\end{proof}

By contrast, Conway's Theorem 1 (transitivity of $\le$) is somewhat harder to establish with our definition; see \autoref{thm:NO-unstrict-transitive}.



We will also need the joint recursion principle, \define{$(\NO,\le,<)$-recursion}, which it is convenient to state as follows.
Suppose $A$ is a type equipped with relations $\mathord\ble:A\to A\to\prop$ and $\mathord\blt:A\to A\to\prop$.
Then we can define $f:\NO\to A$ by doing the following.
\begin{enumerate}
\item For any $x$ defined by a cut, assuming $f(x^L)$ and $f(x^R)$ to be defined such that $f(x^L)\blt f(x^R)$ for all $L$ and $R$, we must define $f(x)$.  (We call this the \emph{primary clause} of the recursion.)\label{item:NO-rec-primary}
\item Prove that $\ble$ is \emph{antisymmetric}\index{relation!antisymmetric}: if $a\ble b$ and $b\ble a$, then $a=b$.
\item For $x,y$ defined by cuts such that $x^L<y$ for all $L$ and $x<y^R$ for all $R$, and assuming inductively that $f(x^L)\blt f(y)$ for all $L$, $f(x)\blt f(y^R)$ for all $R$, and also that $f(x^L)\blt f(x^R)$ and $f(y^L)\blt f(y^R)$ for all $L$ and $R$, we must prove $f(x)\ble f(y)$.
\item For $x,y$ defined by cuts and an $L_0$ such that $x\le y^{L_0}$, and assuming inductively that $f(x)\ble f(y^{L_0})$, and also that $f(x^L)\blt f(x^R)$ and $f(y^L)\blt f(y^R)$ for all $L$ and $R$, we must prove $f(x)\blt f(y)$.
\item For $x,y$ defined by cuts and an $R_0$ such that $x^{R_0}\le y$, and assuming inductively that $f(x^{R_0})\ble f(y)$, and also that $f(x^L)\blt f(x^R)$ and $f(y^L)\blt f(y^R)$ for all $L$ and $R$, we must prove $f(x)\blt f(y)$.\label{item:NO-rec-last}
\end{enumerate}
The last three clauses can be more concisely described by saying we must prove that $f$ (as defined in the first clause) takes $\le$ to $\ble$ and $<$ to $\blt$.
We will refer to these properties by saying that \emph{$f$ preserves inequalities}.
Moreover, in proving that $f$ preserves inequalities, we may assume the particular instance of $\le$ or $<$ to be obtained from one of its constructors, and we may also use inductive hypotheses that $f$ preserves all inequalities appearing in the input to that constructor.

If we succeed at~\ref{item:NO-rec-primary}--\ref{item:NO-rec-last} above, then we obtain $f:\NO\to A$, which computes on cuts as specified by~\ref{item:NO-rec-primary}, and which preserves all inequalities:
\begin{narrowmultline*}
  \fall{x,y:\NO}\Big((x\le y) \to (f(x)\ble f(y))\Big) \land
  \narrowbreak
  \Big((x< y) \to (f(x)\blt f(y))\Big).  
\end{narrowmultline*}
Like $(\RC,\closesym)$-recursion for the Cauchy reals, this recursion principle is essential for defining functions on $\NO$, since we cannot first define a function on ``pre-surreals'' and only later prove that it respects the notion of equality.

\begin{eg}
  Let us define the \emph{negation} function $\NO\to\NO$.
  We apply the joint recursion principle with $A\defeq\NO$, with $(x\ble y)\defeq (y\le x)$, and $(x\blt y)\defeq (y< x)$.
  Clearly this $\ble$ is antisymmetric.

  For the main clause in the definition, we assume $x$ defined by a cut, with $-x^L$ and $-x^R$ defined such that $-x^L \blt -x^R$ for all $L$ and $R$.
  By definition, this means $-x^R< -x^L$ for all $L$ and $R$, so we can define $-x$ by the cut $\surr{-x^R}{-x^L}$.
  This notation, which follows Conway, refers to the cut whose left options are indexed by the type $\RR$ indexing the right options of $x$, and whose right options are indexed by the type $\LL$ indexing the left options of $x$, with the corresponding families $\RR\to\NO$ and $\LL\to\NO$ defined by composing those for $x$ with negation.

  We now have to verify that $f$ preserves inequalities.
  \begin{itemize}
  \item For $x\le y$, we may assume $x^L<y$ for all $L$ and $x < y^R$ for all $R$, and show $-y\le -x$.
    But inductively, we may assume $-y <-x^L$ and $-y^R<-x$, which gives the desired result, by definition of $-y$, $-x$, and the constructor of $\le$.
  \item For $x<y$, in the first case when it arises from some $x\le y^{L_0}$, we may inductively assume $-y^{L_0} \le -x$, in which case $-y<-x$ follows by the constructor of $<$.
  \item Similarly, if $x<y$ arises from $x^{R_0}\le y$, the inductive hypothesis is $-y \le -x^R$, yielding $-y<-x$ again.
  \end{itemize}
\end{eg}

To do much more than this, however, we will need to characterize the relations $\le$ and $<$ more explicitly, as we did for the Cauchy reals in \autoref{thm:RC-sim-characterization}.
Also as there, we will have to simultaneously prove a couple of essential properties of these relations, in order for the induction to go through.

\begin{thm}\label{defn:No-codes}
  There are relations $\mathord\preceq:\NO\to\NO\to\prop$ and $\mathord\prec:\NO\to\NO\to\prop$ such that if $x$ and $y$ are surreals defined by cuts, then
  \begin{align*}
    (x\preceq y) &\defeq
    \big(\fall{L} x^L\prec y\big) \land \big(\fall{R} x\prec y^R\big)\\
    (x\prec y) &\defeq
    \big(\exis{L} x\preceq y^L\big) \lor \big(\exis{R} x^R \preceq y\big).
  \end{align*}
  Moreover, we have
  \begin{equation}\label{eq:NO-codes-unstrict}
    (x\prec y) \to (x\preceq y)
  \end{equation}
  and all the reasonable transitivity properties making $\prec$ and $\preceq$ into a ``bimodule''\index{bimodule} over $\le$ and $<$:
  \begin{equation}\label{eq:NO-codes-transitivity}
    \begin{array}{c@{\hspace{1cm}}c}
      (x \le y) \to (y\preceq z) \to (x\preceq z) &
      (x \preceq y) \to (y\le z) \to (x\preceq z) \\
      (x \le y) \to (y\prec z) \to (x\prec z) &
      (x \preceq y) \to (y< z) \to (x\prec z) \\
      (x < y) \to (y\preceq z) \to (x\prec z) &
      (x \prec y) \to (y\le z) \to (x\prec z).
  \end{array}
  \end{equation}
\end{thm}

\begin{proof}
  We define $\preceq$ and $\prec$ by double $(\NO,\le,<)$-induction on $x,y$.
  The first induction is a simple recursion, whose codomain is the subset $A$ of $(\NO\to\prop)\times (\NO\to\prop)$ consisting of pairs of predicates of which one implies the other and which satisfy ``transitivity on the right'', i.e.~\eqref{eq:NO-codes-unstrict} and the right column of~\eqref{eq:NO-codes-transitivity} with $(x\preceq \blank)$ and $(x\prec \blank)$ replaced by the two given predicates.
  As in the proof of \autoref{defn:RC-approx}, we regard these predicates as half of binary relations, writing them as $y\mapsto (\hle y)$ and $y\mapsto (\hlt y)$, with $\hlname$ denoting the pair of relations.
  We equip $A$ with the following two relations:
  \begin{align*}
    (\hlname \ble \hlbname) &\defeq
    \fall{y:\NO} \Big( (\hleb y) \to (\hle y) \Big) \land
    \Big( (\hltb y) \to (\hlt y) \Big),\\
    (\hlname \blt \hlbname) &\defeq
    \fall{y:\NO} \Big( (\hleb y) \to (\hlt y) \Big).
  \end{align*}
  Note that $\ble$ is antisymmetric, since if $\hlname \ble \hlbname$ and $\hlbname \ble \hlname$, then $(\hleb y) \Leftrightarrow (\hle y)$ and $(\hltb y) \Leftrightarrow (\hlt y)$ for all $y$, hence $\hlname=\hlbname$ by univalence for mere propositions and function extensionality.
  Moreover, to say that a function $\NO\to A$ preserves inequalities is exactly to say that, when regarded as a pair of binary relations on $\NO$, it satisfies ``transitivity on the left'' (the left column of~\eqref{eq:NO-codes-transitivity}).

  Now for the primary clause of the recursion, we assume given $x$ defined by a cut, and relations $(x^L \prec \blank)$, $(x^R \prec \blank)$, $(x^L \preceq \blank)$, and $(x^R \preceq \blank)$ for all $L$ and $R$, of which the strict ones imply the non-strict ones, which satisfy transitivity on the right, and such that
  \begin{equation}\label{eq:NO-prec-outer-IH}
    \fall{L,R}{y:\NO}\Big( (x^R\preceq y) \to (x^L \prec y) \Big).
  \end{equation}
  We now have to define $(x\prec y)$ and $(x\preceq y)$ for all $y$.
  Here in contrast to \autoref{defn:RC-approx}, rather than a nested recursion, we use a nested induction, in order to be able to inductively use transitivity on the left with respect to the inequalities $x^L<x$ and $x<x^R$.
  Define $A':\NO\to\type$ by taking $A'(y)$ to be the subset $A'$ of $\prop\times\prop$ consisting of two mere propositions, denoted $\tle y$ and $\tlt y$ (with $\tlname:A'(y)$), such that
  \begin{gather}
    (\tlt y) \to (\tle y)\\
    \fall{L} (\tle y)\to (x^L\prec y) \label{eq:NO-prec-IHL}\\
    \fall{R} (x^R \preceq y) \to (\tlt y) \label{eq:NO-prec-IHR}.
  \end{gather}
  Using notation analogous to $\ble$ and $\blt$, we equip $A'$ with the two relations defined for $\tlname:A'(y)$ and $\tlbname:A'(z)$ by
  \begin{align*}
    (\tlname \bble \tlbname) &\defeq
    \Big((\tle y) \to (\tleb z)\Big) \land \Big((\tlt y) \to (\tltb z)\Big)\\
    (\tlname \bblt \tlbname) &\defeq
    \Big((\tle y) \to (\tltb z)\Big). 
  \end{align*}
  Again, $\bble$ is evidently antisymmetric in the appropriate sense.
  Moreover, a function $\prd{y:\NO} A'(y)$ which preserves inequalities is precisely a pair of predicates of which one implies the other, which satisfy transitivity on the right, and transitivity on the left with respect to the inequalities $x^L<x$ and $x<x^R$.
  Thus, this inner induction will provide what we need to complete the primary clause of the outer recursion.

  For the primary clause of the inner induction, we assume also given $y$ defined by a cut, and properties $(x\prec y^L)$, $(x\prec y^R)$, $(x\preceq y^L)$, and $(x\preceq y^R)$ for all $L$ and $R$, with the strict ones implying the non-strict ones, transitivity on the left with respect to $x^L<x$ and $x<x^R$, and on the right with respect to $y^L<y^R$.
  We can now give the definitions specified in the theorem statement:
  \begin{align}
    (x\preceq y) &\defeq
    (\fall{L} x^L\prec y) \land (\fall{R} x\prec y^R), \label{eq:NO-preceq-def}\\
    (x\prec y) &\defeq
    (\exis{L} x\preceq y^L) \lor (\exis{R} x^R \preceq y).\label{eq:NO-prec-def}
  \end{align}
  For this to define an element of $A'(y)$, we must show first that $(x\prec y) \to (x\preceq y)$.
  The assumption $x\prec y$ has two cases.
  On one hand, if there is $L_0$ with $x\preceq y^{L_0}$, then by transitivity on the right with respect to $y^{L_0}<y^R$, we have $x\prec y^R$ for all $R$.
  Moreover, by transitivity on the left with respect to $x^L<x$, we have $x^L \prec y^{L_0}$ for any $L$, hence $x^L\prec y$ by transitivity on the right.
  Thus, $x\preceq y$.

  On the other hand, if there is $R_0$ with $x^{R_0}\preceq y$, then by transitivity on the left with respect to $x^L<x^{R_0}$ we have $x^L \prec y$ for all $L$.
  And by transitivity on the left and right with respect to $x<x^{R_0}$ and $y<y^R$, we have $x\prec y^R$ for any $R$.
  Thus, $x\preceq y$.

  We also need to show that these definitions are transitive on the left with respect to $x^L<x$ and $x<x^R$.
  But if $x\preceq y$, then $x^L\prec y$ for all $L$ by definition; while if $x^R\preceq y$, then $x\prec y$ also by definition.

  Thus,~\eqref{eq:NO-preceq-def} and~\eqref{eq:NO-prec-def} do define an element of $A'(y)$.
  We now have to verify that this definition preserves inequalities, as a dependent function into $A'$, i.e.\ that these relations are transitive on the right.
  Remember that in each case, we may assume inductively that they are transitive on the right with respect to all inequalities arising in the inequality constructor.
  \begin{itemize}
  \item Suppose $x\preceq y$ and $y\le z$, the latter arising from $y^L<z$ and $y<z^R$ for all $L$ and $R$.
    Then the inductive hypothesis (of the inner recursion) applied to $y<z^R$ yields $x\prec z^R$ for any $R$.
    Moreover, by definition $x\preceq y$ implies that $x^L \prec y$ for any $L$, so by the inductive hypothesis of the outer recursion we have $x^L \prec z$.
    Thus, $x\preceq z$.
  \item Suppose $x\preceq y$ and $y<z$.
    First, suppose $y<z$ arises from $y\le z^{L_0}$.
    Then the inner inductive hypothesis applied to $y\le z^{L_0}$ yields $x \preceq z^{L_0}$, hence $x\prec z$.

    Second, suppose $y<z$ arises from $y^{R_0}\le z$.
    Then by definition, $x\preceq y$ implies $x\prec y^{R_0}$, and then the inner inductive hypothesis for $y^{R_0}\le z$ yields $x\prec z$.
  \item Suppose $x\prec y$ and $y\le z$, the latter arising from $y^L<z$ and $y<z^R$ for all $L$ and $R$.
    By definition, $x\prec y$ implies there merely exists $R_0$ with $x^{R_0}\preceq y$ or $L_0$ with $x\preceq y^{L_0}$.
    If $x^{R_0}\preceq y$, then the outer inductive hypothesis yields $x^{R_0}\preceq z$, hence $x\prec z$.
    If $x\preceq y^{L_0}$, then the inner inductive hypothesis for $y^{L_0}<z$ (which holds by the constructor of $y\le z$) yields $x\prec z$.

  \end{itemize}
  This completes the inner induction.
  Thus, for any $x$ defined by a cut, we have $(x\prec \blank)$ and $(x\preceq \blank)$ defined by~\eqref{eq:NO-preceq-def} and~\eqref{eq:NO-prec-def}, and transitive on the right.

  To complete the outer recursion, we need to verify these definitions are transitive on the left.
  After a $\NO$-induction on $z$, we end up with three cases that are essentially identical to those just described above for transitivity on the right.
  Hence, we omit them.
\end{proof}

\begin{thm}\label{thm:NO-encode-decode}
  For any $x,y:\NO$ we have $(x<y)=(x\prec y)$ and $(x\le y)=(x\preceq y)$.
\end{thm}
\begin{proof}
  From left to right, we use $(\NO,\le,<)$-induction where $A(x)\defeq\unit$, with $\preceq$ and $\prec$ supplying the relations $\ble$ and $\blt$.
  In all the constructor cases, $x$ and $y$ are defined by cuts, so the definitions of $\preceq$ and $\prec$ evaluate, and the inductive hypotheses apply.

  From right to left, we use $\NO$-induction to assume that $x$ and $y$ are defined by cuts.
  But now the definitions of $\preceq$ and $\prec$, and the inductive hypotheses, supply exactly the data required for the relevant constructors of $\le$ and $<$.
\end{proof}

\begin{cor}\label{thm:NO-unstrict-transitive}
  The relations $\le$ and $<$ on $\NO$ satisfy
  \[ \fall{x,y:\NO} (x<y) \to (x\le y) \]
  and are transitive:
  \index{transitivity!of . for surreals@of $<$ for surreals}
  \index{transitivity!of . for surreals@of $\leq$ for surreals}
  \begin{gather*}
    (x\le y) \to (y\le z) \to (x\le z)\\
    (x\le y) \to (y< z) \to (x< z)\\
    (x< y) \to (y\le z) \to (x< z).
  \end{gather*}
\end{cor}

As with the Cauchy reals, the joint $(\NO,\le,<)$-recursion principle remains essential when defining all operations on $\NO$.

\begin{eg}
\index{addition!of surreal numbers}%
We define $\mathord+:\NO\to\NO\to\NO$ by a double recursion.
For the outer recursion, we take the codomain to be the subset of $\NO\to\NO$ consisting of functions $g$ such that $(x<y) \to (g(x)<g(x))$ and $(x\le y) \to (g(x)\le g(y))$ for all $x,y$.
For such $g,h$ we define $(g\ble h)\defeq \fall{x:\NO} g(x)\le h(x)$ and $(g\blt h)\defeq \fall{x:\NO} g(x)< h(x)$.
Clearly $\ble$ is antisymmetric.

For the primary clause of the recursion, we suppose $x$ defined by a cut, and we define $(x+\blank)$ by an inner recursion on $\NO$ with codomain $\NO$, with relations $\bble$ and $\bblt$ coinciding with $\le$ and $<$.
For the primary clause of the inner recursion, we suppose also $y$ defined by a cut, and give Conway's definition:
\[ x+y \defeq \surr{x^L+y, x+y^L}{x^R+y,x+y^R}. \]
In other words, the left options of $x+y$ are all numbers of the form $x^L+y$ for some left option $x^L$, or $x+y^L$ for some left option $y^L$.
Now we verify that this definition preserves inequality:
\begin{itemize}
\item If $y\le z$ arises from knowing that $y^L<z$ and $y<z^R$ for all $L$ and $R$, then the inner inductive hypothesis gives $x+y^L<x+z$ and $x+y < x+z^R$, while the outer inductive hypotheses give $x^L+y < x^L+z$ and $x^R+ y < x^R+z$.
  And since each $x^L+z$ is by definition a left option of $x+z$, we have $x^L+z < x+z$, and similarly $x+y < x^R+y$.
  Thus, using transitivity, $x^L+y < x+z$ and $x+y < x^R+z$, and so we may conclude $x+y \le x+z$ by the constructor of $\le$.
\item If $y<z$ arises from an $L_0$ with $y\le z^{L_0}$, then inductively $x+y \le x+z^{L_0}$, hence $x+y<x+z$ since $x+z^{L_0}$ is a right option of $x+z$.
\item Similarly, if $y<z$ arises from $y^{R_0}\le z$, then $x+y<x+z$ since $x+y^{R_0}\le x+z$.
\end{itemize}
This completes the inner recursion.
For the outer recursion, we have to verify that $+$ preserves inequality on the left as well.
After an $\NO$-induction, this proceeds in exactly the same way.
\end{eg}

\index{acceptance|(}%
\index{mathematics!formalized}%
In the Appendix to Part Zero of~\cite{conway:onag}, Conway discusses how the surreal numbers may be formalized in ZFC set theory: by iterating along the ordinals and passing to sets of representatives of lowest rank for each equivalence class, or by representing numbers with ``sign-expansions''.
He then remarks that
\begin{quote}\footnotesize
  The curiously complicated nature of these constructions tells us more about the nature of formalizations within ZF than about our system of numbers\dots
\end{quote}
and goes on to advocate for a general theory of ``permissible kinds of construction'' which should include
\begin{enumerate}\footnotesize
\item Objects may be created from earlier objects in any reasonably constructive fashion.\label{item:conway1}
\item Equality among the created objects can be any desired equivalence relation.\label{item:conway2}
\end{enumerate}
\noindent
Condition~\ref{item:conway1} can be naturally read as justifying general principles of \emph{inductive definition}, such as those presented in \autoref{sec:strictly-positive,sec:generalizations}.
In particular, the condition of strict positivity for constructors can be regarded as a formalization of what it means to be ``reasonably constructive''.
Condition~\ref{item:conway2} then suggests we should extend this to \emph{higher} inductive definitions of all sorts, in which we can impose path constructors making objects equal in any reasonable way.
For instance, in the next paragraph Conway says:
\begin{quote}\footnotesize
  \dots we could also, for instance, freely create a new object $(x,y)$ and call it the ordered pair of $x$ and $y$.
  We could also create an ordered pair $[x,y]$ different from $(x,y)$ but co-existing with it\dots
  If instead we wanted to make $(x,y)$ into an unordered pair, we could define equality by means of the equivalence relation $(x,y)=(z,t)$ if and only if $x=z,y=t$ \emph{or} $x=t,y=z$.
\end{quote}
The freedom to introduce new objects with new names, generated by certain forms of constructors, is precisely what we have in the theory of inductive definitions.
Just as with our two copies of the natural numbers $\nat$ and $\nat'$ in \autoref{sec:appetizer-univalence}, if we wrote down an identical definition to the cartesian product type $A\times B$, we would obtain a distinct product type $A\times' B$ whose canonical elements we could freely write as $[x,y]$.
And we could make one of these a type of unordered pairs by adding a suitable path constructor. 

To be sure, Conway's point was not to complain about ZF in particular, but to argue against all foundational theories at once:
\begin{quote}\footnotesize
  \dots this proposal is not of any particular theory as an alternative to ZF\dots{}
  What is proposed is instead that we give ourselves the freedom to create arbitrary mathematical theories of these kinds, but prove a metatheorem which ensures once and for all that any such theory could be formalized in terms of any of the standard foundational theories.
\end{quote}
One might respond that, in fact, univalent foundations is not one of the ``standard foundational theories'' which Conway had in mind, but rather the \emph{metatheory} in which we may express our ability to create new theories, and about which we may prove Conway's metatheorem.
For instance, the surreal numbers are one of the ``mathematical theories'' Conway has in mind, and we have seen that they can be constructed and justified inside univalent foundations.
Similarly, Conway remarked earlier that
\begin{quote}\footnotesize
  \dots set theory would be such a theory, sets being constructed from earlier ones by processes corresponding to the usual axioms, and the equality relation being that of having the same members.
\end{quote}
This description closely matches the higher-inductive construction of the cumulative hierarchy of set theory in \autoref{sec:cumulative-hierarchy}.
Conway's metatheorem would then correspond to the fact we have referred to several times that we can construct a model of univalent foundations inside ZFC (which is outside the scope of this book).

However, univalent foundations is so rich and powerful in its own right that it would be foolish to relegate it to only a metatheory in which to construct set-like theories.
We have seen that even at the level of sets (0-types), the higher inductive types in univalent foundations yield direct constructions of objects by their universal properties (\autoref{sec:free-algebras}), such as a constructive theory of Cauchy completion (\autoref{sec:cauchy-reals}).
But most importantly, the potential to model homotopy theory and category theory directly in the foundational system (\autoref{cha:homotopy,cha:category-theory}) gives univalent foundations an advantage which no set-theoretic foundation can match.
\index{acceptance|)}%

\index{surreal numbers|)}%

\sectionNotes

Defining algebraic operations on Dedekind reals, especially multiplication, is both somewhat tricky and tedious.
There are several ways to get arithmetic going: each has its own advantages, but they all seem to require some technical work.
For instance, Richman~\cite{Richman:reals} defines multiplication on the Dedekind reals first on the positive cuts and then extends it algebraically to all Dedekind cuts, while Conway~\cite{conway:onag} has observed that the definition of multiplication for surreal numbers works well for Dedekind reals.

Our treatment of the Dedekind reals borrows many ideas from~\cite{BauerTaylor09} where the Dedekind reals are constructed in the context of Abstract Stone Duality.
\index{Abstract Stone Duality}%
This is a (restricted) form of simply typed $\lambda$-calculus with a distinguished object $\Sigma$ which classifies open sets, and by duality also the closed ones. In~\cite{BauerTaylor09} you can also find detailed proofs of the basic properties of arithmetical operations.

The fact that $\RC$ is the least Cauchy complete archimedean ordered field, as was proved in \autoref{RC-initial-Cauchy-complete}, indicates that our Cauchy reals probably coincide with the Escard{\'o}-Simpson reals~\cite{EscardoSimpson:01}.
\index{real numbers!Escardo-Simpson@Escard\'o-Simpson}%
It would be interesting to check\index{open!problem} whether this is really the case. The notion of Escard{\'o}-Simpson reals, or more precisely the corresponding closed interval, is interesting because it can be stated in any category with finite products.

In constructive set theory augmented by the ``regular extension axiom'', one may also try to define Cauchy completion by closing under limits of Cauchy sequences with a transfinite iteration.
It would also be interesting to check whether this construction agrees with ours.

It is constructive folklore that coincidence of Cauchy and Dedekind reals requires dependent choice but it is less well known that countable choice suffices. Recall that \define{dependent choice}
\indexdef{axiom!of choice!dependent}%
\index{axiom!of choice!countable}%
\index{total!relation}%
states that for a total relation $R$ on $A$, by which we mean $\fall{x : A} \exis{y : A} R(x,y)$, and for any $a : A$ there merely exists $f : \N \to A$ such that $f(0) = a$ and $R(f(n), f(n+1))$ for all $n : \N$. Our \autoref{when-reals-coincide} uses the typical trick for converting an application of dependent choice to one using countable choice. Namely, we use countable choice once to make in advance all the choices that could come up, and then use the choice function to avoid the dependent choices.

The intricate relationship between various notions of compactness in a constructive
setting is discussed in \cite{bridges2002compactness}. Palmgren~\cite{Palmgren:FT} has a 
good comparison between pointwise analysis and 
pointfree topology.

The surreal numbers were defined by~\cite{conway:onag}, using a sort of inductive definition but without justifying it explicitly in terms of any foundational system.
For this reason, some later authors have tended to use sign-expansions or other more explicit presentations which can be coded more obviously into set theory.
The idea of representing them in type theory was first considered by Hancock, while
Setzer and Forsberg~\cite{forsbergfinite} noted that the surreals and their inequality relations $<$ and $\le$ naturally form an inductive-inductive definition.
The \emph{higher} inductive-inductive version presented here, which builds in the correct notion of equality for surreals, is new.

\sectionExercises

\begin{ex}
 Give an alternative definition of the Dedekind reals by first defining the square and then use \autoref{mult-from-square}.
 Check that one obtains a commutative ring.
\end{ex}

\begin{ex} \label{ex:RD-extended-reals}
  Suppose we remove the boundedness condition in \autoref{defn:dedekind-reals}.
  Then we obtain the \define{extended reals}
  \indexdef{real numbers!extended}%
  \indexdef{extended real numbers}%
  which contain $-\infty \defeq
  (\emptyt, \Q)$ and $\infty \defeq (\Q, \emptyt)$. Which definitions of arithmetical
  operations on cuts still make sense for extended reals? What algebraic structure do we
  get?
\end{ex}

\begin{ex} \label{ex:RD-lower-cuts}
  By considering one-sided cuts we obtain \define{lower} and \define{upper} Dedekind reals,
  \indexdef{real numbers!Dedekind!upper}%
  \indexdef{real numbers!Dedekind!lower}%
  \indexdef{lower Dedekind reals}%
  \indexdef{upper Dedekind reals}%
  \index{cut!Dedekind}%
  respectively. For example, a lower real is given by a predicate $L : \Q \to \Omega$
  which is
  \begin{enumerate}
  \item \emph{inhabited:} $\exis{q : \Q} L(q)$ and
  \item \emph{rounded:} $L(q) = \exis{r : \Q} q < r \land L(r)$.
    \index{rounded!Dedekind cut}
  \end{enumerate}
  (We could also require $\exis{r : \Q} \lnot L(r)$ to exclude the cut $\infty \defeq
  \Q$.) Which arithmetical operations can you define on the lower reals? In particular,
  what happens with the additive inverse?
\end{ex}

\begin{ex} \label{ex:RD-interval-arithmetic}
  \index{interval!arithmetic}%
  Suppose we remove the locatedness condition in \autoref{defn:dedekind-reals}.
  Then we obtain the \define{interval domain}
  \indexdef{interval!domain}%
  $\mathbb{I}$ because cuts are allowed
  to have ``gaps'', which are just intervals. Define the partial order $\sqsubseteq$ on
  $\mathbb{I}$ by
  \begin{narrowmultline*}
    ((L, U) \sqsubseteq (L', U'))
    \defeq \narrowbreak
    (\fall{q : \Q} L(q) \Rightarrow L'(q)) \land
    (\fall{q : \Q} U(q) \Rightarrow U'(q)).
  \end{narrowmultline*}
  What are the maximal elements of $\mathbb{I}$ with respect to $\mathbb{I}$? Define the
  ``endpoint'' operations which assign to an element of the interval domain its lower and
  upper endpoints. Are the endpoints reals, lower reals, or upper reals (see
  \autoref{ex:RD-lower-cuts})? Which definitions of arithmetical operations on cuts still
  make sense for the interval domain?
\end{ex}

\begin{ex} \label{ex:RD-lt-vs-le}
  Show that, for all $x, y : \RD$,
  \begin{equation*}
    \lnot (x < y) \Rightarrow y \leq x
  \end{equation*}
  and
  \begin{equation*}
    \eqv{(x \leq y)}{\Parens{\prd{\epsilon : \Qp} x < y + \epsilon}}.
  \end{equation*}
  Does $\lnot (x \leq y)$ imply $y < x$?
\end{ex}

\begin{ex} \label{ex:reals-non-constant-into-Z}
  \mbox{}
  \begin{enumerate}
  \item 
    Assuming excluded middle, construct a non-constant map $\RD \to \Z$.
  \item 
    Suppose $f : \RD \to \Z$ is a map such that $f(0) = 0$ and $f(x) \neq 0$ for all $x >
    0$. Derive from this the limited principle of omniscience~\eqref{eq:lpo}.
\index{limited principle of omniscience}%
  \end{enumerate}
\end{ex}

\begin{ex} \label{ex:traditional-archimedean}
  \index{ordered field!archimedean}%
  Show that in an ordered field $F$, density of $\Q$ and the traditional archimedean axiom
  are equivalent:
  \begin{equation*}
    (\fall{x, y : F} x < y \Rightarrow \exis{q : \Q} x < q < y)
    \Leftrightarrow
    (\fall{x : F} \exis{k : \Z} x < k).
  \end{equation*}  
\end{ex}

\begin{ex} \label{RC-Lipschitz-on-interval} Suppose $a, b : \Q$ and $f : \setof{ q : \Q |
    a \leq q \leq b } \to \RC$ is Lipschitz with constant~$L$. Show that there exists a unique
  extension $\bar{f} : [a,b] \to \RC$ of $f$ which is Lipschitz with
  constant~$L$. Hint: rather than redoing \autoref{RC-extend-Q-Lipschitz} for closed
  intervals, observe that there is a retraction $r : \RC \to [-n,n]$ and apply
  \autoref{RC-extend-Q-Lipschitz} to $f \circ r$.
\end{ex}

\begin{ex} \label{ex:metric-completion}
  \index{completion!of a metric space}%
  Generalize the construction of $\RC$ to construct the Cauchy completion of any metric space. First, think about which notion of real numbers is most natural as the codomain for the distance\index{distance} function of a metric space. Does it matter? Next, work out the details of two constructions:
  \begin{enumerate}
  \item Follow the construction of Cauchy reals to define the completion of a metric space as an inductive-inductive type closed under limits of Cauchy sequences.\index{Cauchy!sequence}
  \item Use the following construction due to Lawvere~\cite{lawvere:metric-spaces}\index{Lawvere} and Richman~\cite{Richman00thefundamental}, where the completion of a metric space $(M, d)$ is given as the type of \define{locations}.
    \indexdef{location}%
    A location is a function $f : M \to \R$ such that
    \begin{enumerate}
    \item $f(x) \geq |f(y) - d(x,y)|$ for all $x, y : M$, and
    \item $\inf_{x \in M} f(x) = 0$, by which we mean $\fall{\epsilon : \Qp} \exis{x : M} |f(x)| < \epsilon$ and $\fall{x : M} f(x) \geq 0$.
    \end{enumerate}
    The idea is that $f$ looks like it is measuring the distance from a point.
  \end{enumerate}
  \index{universal!property!of metric completion}%
  Finally, prove the following universal property of metric completions: a locally uniformly continuous map from a metric space to a Cauchy complete metric space extends uniquely to a locally uniformly continuous map on the completion. (We say that a map is \define{locally uniformly continuous}
  \indexdef{function!locally uniformly continuous}%
  \indexdef{locally uniformly continuous map}%
  if it is uniformly continuous on open balls.)
\end{ex}

\index{metric space|)}%

\begin{ex} \label{ex:reals-apart-neq-MP}
  \define{Markov's principle}
  \indexdef{axiom!Markov's principle}%
  \indexdef{Markov's principle}%
  says that for all $f : \nat \to \bool$,
  \begin{equation*}
    (\lnot \lnot \exis{n : \nat} f(n) = \btrue)
    \Rightarrow
    \exis{n : \nat} f(n) = \btrue.
  \end{equation*}
  This is a particular instance of the law of double negation~\eqref{eq:ldn}. Show that
  $\fall{x, y: \RD} x \neq y \Rightarrow x \apart y$ implies Markov's principle. Does the
  converse hold as well?
\end{ex}

\begin{ex} \label{ex:reals-apart-zero-divisors}
  \index{apartness}%
  Verify that the following ``no zero divisors'' property holds for the real numbers:
  $x y \apart 0 \Leftrightarrow x \apart 0 \land y \apart 0$.
\end{ex}

\begin{ex} \label{ex:finite-cover-lebesgue-number}
  Suppose $(q_1, r_1), \ldots, (q_n, r_n)$ pointwise cover $(a, b)$. Then there is
  $\epsilon : \Qp$ such that whenever $a < x < y < b$ and $|x - y| < \epsilon$
  then there merely exists $i$ such that $q_i < x < r_i$ and $q_i < y < r_i$. Such an
  $\epsilon$ is called a \define{Lebesgue number}
  \indexdef{Lebesgue number}%
  for the given cover.
\end{ex}

\begin{ex} \label{ex:mean-value-theorem}
  Prove the following approximate version of the intermediate value theorem:
  \begin{quote}
    \emph{
      If $f : [0,1] \to \R$ is uniformly continuous and $f(0) < 0 < f(1)$ then
      for every $\epsilon : \Qp$ there merely exists $x : [0,1]$ such that $|f(x)| <
      \epsilon$.
    }
  \end{quote}
  Hint: do not try to use the bisection method because it leads to the axiom of choice.
  Instead, approximate $f$ with a piecewise linear map. How do you construct a piecewise
  linear map?
\end{ex}

\begin{ex}
  Check whether everything in~\cite{knuth74:_surreal_number} can be done using the higher
  inductive-inductive surreals of \autoref{sec:surreals}.
\end{ex}

\index{real numbers|)}%



\cleartooddpage[\thispagestyle{empty}] 
\phantomsection 
\part*{Appendix}


\appendix

\renewcommand{\chaptermark}[1]{\markboth{\textsc{Appendix. \thechapter. #1}}{}}
\renewcommand{\sectionmark}[1]{\markright{\textsc{\thesection\ #1}}}


\titleformat{\chapter}[display]{\fontsize{23}{25}\fontseries{m}\fontshape{it}\selectfont}{\chaptertitlename}{20pt}{\fontsize{35}{35}\fontseries{b}\fontshape{n}\selectfont}
\chapter{Formal type theory}
\label{cha:rules}

\index{formal!type theory|(}%
\index{type theory!formal|(}%
\index{rules of type theory|(}%

Just as one can develop mathematics in set theory without explicitly using the axioms of Zermelo--Fraenkel set theory, 
in this book we have developed mathematics in univalent foundations without explicitly referring to a formal
system of homotopy type theory. Nevertheless, it is important to \emph{have} a
precise description of homotopy type theory as a formal system in order to, for example,
\begin{itemize}
\item state and prove its metatheoretic properties, including logical
consistency,
\item construct models, e.g.\  in simplicial sets, model categories, higher toposes,
etc., and
\item implement it in proof assistants like \Coq or \Agda.
  \index{proof!assistant}
\end{itemize}
Even the logical consistency\index{consistency} of homotopy type theory, namely that in the empty context there is no term $a:\emptyt$, is not obvious: if we had erroneously
chosen a definition of equivalence for which $\eqv{\emptyt}{\unit}$, then
univalence would imply that $\emptyt$ has an element, since $\unit$ does.
Nor is it obvious that, for example, our definition of $\Sn^1$ as a higher
inductive type yields a type which behaves like the ordinary circle.

There are two aspects of type theory which we must pin down before addressing
such questions. Recall from the Introduction that type theory
comprises a set of rules specifying when the judgments $a:A$ and $a\jdeq a':A$
hold---for example, products are characterized by the rule that whenever $a:A$
and $b:B$, $(a,b):A\times B$. To make this precise, we must first define
precisely the syntax of terms---the objects $a,a',A,\dots$ which these judgments
relate; then, we must define precisely the judgments and their rules of
inference---the manner in which judgments can be derived from other judgments.

In this appendix, we present two formulations of Martin-L\"{o}f type
theory, and of the extensions that constitute homotopy type theory.
The first presentation (\autoref{sec:syntax-informally}) describes the syntax of
terms and the forms of judgments as an extension of the untyped
$\lambda$-calculus, while leaving the rules of inference informal.
The second (\autoref{sec:syntax-more-formally}) defines the terms, judgments,
and rules of inference inductively in the style of natural deduction, as
is customary in much type-theoretic literature.

\section*{Preliminaries}
\label{sec:formal-prelim}

In \autoref{cha:typetheory}, we presented the two basic \define{judgments}
\index{judgment}
of type theory. The first, $a:A$, asserts that a term $a$ has type $A$. The second,
$a\jdeq b:A$, states that the two terms $a$ and $b$ are \define{judgmentally
equal}
\index{equality!judgmental}
\index{judgmental equality}
at type $A$. These judgments are inductively defined by a set of
inference rules described in \autoref{sec:syntax-more-formally}.

To construct an element $a$ of a type $A$ is to derive $a:A$; in the book, we
give informal arguments which describe the construction of $a$, but formally,
one must specify a precise term $a$ and a full derivation that $a:A$.

However, the main difference between the presentation of type theory in the book
and in this appendix is that here judgments are explicitly
formulated in an ambient \define{context},
\index{context}
or list of assumptions, of the form
\[
  x_1:A_1, x_2:A_2,\dots,x_n:A_n.
\]
An element $x_i : A_i$ of the context expresses the assumption that the
variable
\index{variable}%
$x_i$ has type $A_i$. The variables $x_1, \ldots, x_n$ appearing in
the context must be distinct. We abbreviate contexts with the letters $\Gamma$
and $\Delta$.

The judgment $a:A$ in context $\Gamma$ is written 
\[ \oftp\Gamma aA \]
and means that $a:A$ under the assumptions listed in $\Gamma$. When the list of
assumptions is empty, we write simply
\[ \oftp{}aA \]
or
\[ \oftp\emptyctx aA \]
where $\emptyctx$ denotes the empty context. The same applies to the equality
judgment
\[
  \jdeqtp\Gamma{a}{b}{A}
\]

However, such judgments are sensible only for \define{well-formed} contexts,
\index{context!well-formed}%
a notion captured by our third and final judgment
\[
  \wfctx{(x_1:A_1, x_2:A_2,\dots,x_n:A_n)}
\]
expressing that each $A_i$ is a type in the context $x_1:A_1,
x_2:A_2,\dots,x_{i-1}:A_{i-1}$.  In particular, therefore, if $\oftp\Gamma aA$ and
$\wfctx\Gamma$, then we know that each $A_i$ contains only the variables
$x_1,\dots,x_{i-1}$, and that $a$ and $A$ contain only the variables
$x_1,\dots,x_n$.
\index{variable!in context}

In informal mathematical presentations, the context is
implicit. At each point in a proof, the mathematician knows which
variables are available and what types they have, either by historical
convention ($n$ is usually a number, $f$ is a function, etc.) or
because variables are explicitly introduced with sentences such as
``let $x$ be a real number''. We discuss some benefits of using explicit
contexts in \autoref{sec:more-formal-pi,sec:more-formal-sigma}.

We write $B[a/x]$ for the \define{substitution}
\index{substitution}%
of a term $a$ for free occurrences of
the variable~$x$ in the term $B$, with possible capture-avoiding
renaming of bound variables,
\index{variable!and substitution}%
as discussed in
\autoref{sec:function-types}. The general form of substitution
\[
   B[a_1,\dots,a_n/x_1,\dots,x_n]
\]
substitutes expressions $a_1,\dots,a_n$ for the variables
$x_1,\dots,x_n$ simultaneously.

To \define{bind a variable $x$ in an expression $B$}
\indexdef{variable!bound}%
means to incorporate both of them into a larger expression, called an \define{abstraction},
\indexdef{abstraction}%
whose purpose is to express the fact that $x$ is ``local'' to $B$, i.e., it
is not to be confused with other occurrences of $x$ appearing
elsewhere. Bound variables are familiar to programmers, but less so to mathematicians.
Various notations are used for binding, such as $x \mapsto B$,
$\lam x B$, and $x \,.\, B$, depending on the situation. We may write $C[a]$ for the
substitution of a term $a$ for the variable in the abstracted expression, i.e.,
we may define $(x.B)[a]$ to be $B[a/x]$. As discussed in
\autoref{sec:function-types}, changing the name of a bound variable everywhere within an expression (``$\alpha$-conversion'')
\index{alpha-conversion@$\alpha $-conversion}%
does not change the expression. Thus, to be very
precise, an expression is an equivalence class of syntactic forms
which differ in names of bound variables.

One may also regard each variable $x_i$ of a judgment
\[
  x_1:A_1, x_2:A_2,\dots,x_n:A_n \vdash a : A
\]
to be bound in its \define{scope},
\indexdef{variable!scope of}%
\index{scope}%
consisting of the expressions $A_{i+1},
\ldots, A_n$, $a$, and $A$.

\section{The first presentation}
\label{sec:syntax-informally}

The objects and types of our type theory may be written as terms using
the following syntax, which is an extension of $\lambda$-calculus with
\emph{variables} $x, x',\dots$,
\index{variable}%
\emph{primitive constants}
\index{primitive!constant}%
\index{constant!primitive}%
$c,c',\dots$, \emph{defined constants}\index{constant!defined} $f,f',\dots$, and term forming
operations
\[
  t \production x \mid \lam{x} t \mid t(t') \mid c \mid f
\]
The notation used here means that a term $t$ is either a variable $x$, or it
has the form $\lam{x} t$ where $x$ is a variable and $t$ is a term, or it has
the form $t(t')$ where $t$ and $t'$ are terms, or it is a primitive constant
$c$, or it is a defined constant $f$. The syntactic markers '$\lambda$', '(',
')', and '.' are punctuation for guiding the human eye.

We use $t(t_1,\dots,t_n)$ as an abbreviation for the repeated application
$t(t_1)(t_2)\dots (t_n)$. We may also use \emph{infix}\index{infix notation} notation, writing $t_1\;
\star\; t_2$ for $\star(t_1,t_2)$ when $\star$ is a primitive or defined
constant.

Each defined constant has zero, one or more \define{defining equations}.
\index{equation, defining}%
\index{defining equation}%
There are two kinds of defined constant. An \emph{explicit}
\index{constant!explicit}
defined constant $f$ has a single defining equation
  \[ f(x_1,\dots,x_n)\defeq t,\]
where $t$ does not involve $f$. 
For example, we might introduce the explicit defined constant $\circ$ with defining equation
  \[ \circ (x,y)(z) \defeq x(y(z)),\]
and use infix notation $x\circ y$ for $\circ(x,y)$. This of course is just composition of functions.

The second kind of defined constant is used to specify a (parameterized) mapping
$f(x_1,\dots,x_n,x)$, where $x$ ranges over a type whose elements are generated
by zero or more primitive constants.  For each such primitive constant $c$ there
is a defining equation of the form
\[
  f(x_1,\dots,x_n,c(y_1,\dots,y_m)) \defeq t,
\]
where $f$ may occur in $t$, but only in such a way that it is clear that the
equations determine a totally defined function. The paradigm examples of such
defined functions are the functions defined by primitive recursion on the
natural numbers. We may call this kind of definition of a function a \emph{total
  recursive definition}.
\index{total!recursive definition}%
In computer science and logic this kind of definition
of a function on a recursive data type has been called a \define{definition by
  structural recursion}.
\index{definition!by structural recursion}%
\index{structural!recursion}%
\index{recursion!structural}%

\define{Convertibility}
\index{convertibility of terms}%
\index{term!convertibility of}%
$t \conv t'$ between terms $t$
and $t'$ is the equivalence relation generated by the defining equations for constants,
the computation rule\index{computation rule!for function types}
\[
  (\lam{x} t)(u) \defeq t[u/x],
\]
and the rules which make it a \emph{congruence} with respect to application and $\lambda$-abstraction\index{lambda abstraction@$\lambda$-abstraction}:
\begin{itemize}
\item if $t \conv t'$ and $s \conv s'$ then $t(s) \conv t'(s')$, and
\item if $t \conv t'$ then $(\lam{x} t) \conv (\lam{x} t')$.
\end{itemize}
\noindent
The equality judgment $t \jdeq u : A$ is then derived by the following single rule:
\begin{itemize}
\item if $t:A$, $u:A$, and $t \conv u$, then $t \jdeq u : A$.
\end{itemize}
Judgmental equality is an equivalence relation.

\subsection{Type universes}

We postulate a hierarchy of \define{universes} denoted by primitive constants
\index{type!universe}
\begin{equation*}
  \UU_0, \quad \UU_1, \quad  \UU_2, \quad \ldots
\end{equation*}
The first two rules for universes say that they form a cumulative hierarchy of types:
\begin{itemize}
\item $\UU_m : \UU_n$ for $m < n$,
\item if $A:\UU_m$ and $m \le n$, then $A:\UU_n$,
\end{itemize}
and the third expresses the idea that an object of a universe can serve as a type and stand to the
right of a colon in judgments:
\begin{itemize}
\item if $\Gamma \vdash A : \UU_n$, and $x$ is a new variable,%
\footnote{By ``new'' we mean that it does not appear in $\Gamma$ or $A$.}
then $\vdash (\Gamma, x:A)\; \ctx$.
\end{itemize}
In the body of the book, an equality judgment $A \jdeq B : \UU_n$ between types
$A$ and $B$ is usually abbreviated to $A \jdeq B$. This is an instance of
typical ambiguity\index{typical ambiguity}, as we can always switch to a larger universe, which however does not affect the validity of the judgment.

The following conversion rule allows us to replace a type by one equal to it in a typing judgment:
\begin{itemize}
\item if $a:A$ and $A \jdeq B$ then $a:B$.
\end{itemize}

\subsection{Dependent function types (\texorpdfstring{$\Pi$}{Π}-types)}

We introduce a primitive constant $c_\Pi$, but write
$c_\Pi(A,\lam{x} B)$ as $\tprd{x:A}B$. Judgments concerning
such expressions and expressions of the form $\lam{x} b$ are introduced by the following rules:
\begin{itemize}
\item if $\Gamma \vdash A:\UU_n$ and $\Gamma,x:A \vdash B:\UU_n$, then $\Gamma \vdash \tprd{x:A}B : \UU_n$
\item if $\Gamma, x:A \vdash b:B$ then $\Gamma \vdash (\lam{x} b) : (\tprd{x:A} B)$
\item if $\Gamma\vdash g:\tprd{x:A} B$ and $\Gamma\vdash t:A$ then $\Gamma\vdash g(t):B[t/x]$
\end{itemize}
If $x$ does not occur freely in $B$, we abbreviate $\tprd{x:A} B$ as the non-dependent function type 
$A\rightarrow B$ and derive the following rule:
\begin{itemize}
\item if $\Gamma\vdash g:A \rightarrow B$ and $\Gamma\vdash t:A$ then $\Gamma\vdash g(t):B$
\end{itemize}
Using non-dependent function types and leaving implicit the context $\Gamma$, the rules above can be written in the following alternative style that we use in the rest of this section of the appendix.
\begin{itemize}
\item if $A:\UU_n$ and $B:A\to\UU_n$, then $\tprd{x:A}B(x) : \UU_n$
\item if $x:A \vdash b:B$ then $ \lam{x} b : \tprd{x:A} B(x)$
\item if $g:\tprd{x:A} B(x)$ and $t:A$ then $g(t):B(t)$
\end{itemize}

\subsection{Dependent pair types (\texorpdfstring{$\Sigma$}{Σ}-types)}

We introduce primitive constants $c_\Sigma$ and $c_{\mathsf{pair}}$. An
expression of the form $c_\Sigma(A,\lam{a} B)$ is written as $\sm{a:A}B$,
and an expression of the form $c_{\mathsf{pair}}(a,b)$ is written as $\tup
a b$. We write $A\times B$ instead of $\sm{x:A} B$ if $x$ is not free in $B$.

Judgments concerning such expressions are introduced by the following
rules:
\begin{itemize}
\item if $A:\UU_n$ and $B: A \rightarrow \UU_n$, then $\sm{x:A}B(x) : \UU_n$
\item if, in addition, $a:A$ and $b:B(a)$, then $\tup a b:\sm{x:A}B(x)$
\end{itemize}
If we have $A$ and $B$ as above, $C : \sm{x:A}B(x) \rightarrow \UU_m$, and
\[
  d:\tprd{x:A}{y:B(x)} C(\tup x y)
\]
we can introduce a defined constant 
\[
  f:\tprd{p:\sm{x:A}B(x)} C(p)
\]
with the defining equation
\[
  f(\tup x y)\defeq d(x,y).
\]
Note that $C$, $d$, $x$, and $y$ may contain extra implicit parameters $x_1,\ldots,x_n$ if they were obtained in some non-empty context; therefore, the fully explicit recursion schema is
\begin{narrowmultline*}
 f(x_1,\dots,x_n,\tup{x(x_1,\dots,x_n)}{y(x_1,\dots,x_n)}) \defeq
 \narrowbreak
 d(x_1,\dots,x_n,\tup{x(x_1,\dots,x_n)}{y(x_1,\dots,x_n)}).
\end{narrowmultline*}

\subsection{Coproduct types}

We introduce primitive constants $c_+$, $c_\inlsym$, and $c_\inrsym$.
We write $A+B$ instead of $c_+(A,B)$, $\inl(a)$ instead of
$c_\inlsym(a)$, and $\inr(a)$ instead of $c_\inrsym(a)$:
\begin{itemize}
\item if $A,B : \UU_n$ then $A + B : \UU_n$
\item moreover, $\inl: A \rightarrow A+B$ and $\inr: B \rightarrow A+B$
\end{itemize}
If we have $A$ and $B$ as above, $C : A+B \rightarrow \UU_m$, 
$d:\tprd{x:A} C(\inl(x))$, and $e:\tprd{y:B} C(\inr(y))$,
then we can introduce a defined constant $f:\tprd{z:A+B}C(z)$ with the defining equations
\begin{equation*}
  f(\inl(x)) \defeq d(x)
  \qquad\text{and}\qquad
  f(\inr(y)) \defeq e(y).
\end{equation*}

\subsection{The finite types}

We introduce primitive constants $\ttt$, $\emptyt$, $\unit$, satisfying the following rules:
\begin{itemize}
\item $\emptyt : \UU_0$, $\unit : \UU_0$
\item $\ttt:\unit$
\end{itemize}

Given $C : \emptyt \rightarrow \UU_n$ we can introduce a defined constant $f:\tprd{x:\emptyt} C(x)$, with no defining equations.

Given $C : \unit \rightarrow \UU_n$ and $d : C(\ttt)$ we can introduce a defined constant $f:\tprd{x:\unit} C(x)$, with defining equation $f(\ttt) \defeq d$.

\subsection{Natural numbers}

The type of natural numbers is obtained by introducing primitive constants
$\N$, $0$, and $\suc$ with the following rules:
\begin{itemize}
  \item $\N : \UU_0$,
  \item $0:\N$,
  \item $\suc:\N\rightarrow \N$.
\end{itemize}
Furthermore, we can define functions by primitive recursion. If we have
$C : \N \rightarrow \UU_k $ we can introduce a defined constant $f:\tprd{x:\N}C(x)$ whenever we have
\begin{align*}
  d & : C(0) \\
  e & : \tprd{x:\N}(C(x)\rightarrow C(\suc (x)))
\end{align*}
with the defining equations
\begin{equation*}
  f(0) \defeq d
  \qquad\text{and}\qquad
  f(\suc (x)) \defeq e(x,f(x)).
\end{equation*}

\subsection{\texorpdfstring{$W$}{W}-types}

For $W$-types we introduce primitive constants $c_\wtypesym$ and $c_\suppsym$.
An expression of the form $c_\wtypesym(A,\lam{x} B)$ is written as
$\wtype{x:A}B$, and an expression of the form $c_\suppsym(x,u)$ is written
as $\supp(x,u)$:
\begin{itemize}
\item if $A:\UU_n$ and $B: A \rightarrow \UU_n$, then $\wtype{x:A}B(x) : \UU_n$
\item if moreover, $a:A$ and $g:B(a)\rightarrow \wtype{x:A}B(x)$ then $\supp(a,g):\wtype{x:A}B(x)$.
\end{itemize}
Here also we can define functions by total recursion. If we have $A$ and $B$
as above and $C : \wtype{x:A}B(x) \rightarrow \UU_m$, then we can introduce a defined constant
$f:\tprd{z:\wtype{x:A}B(x)} C(z)$ whenever we have
\[
  d:\tprd{x:A}{u:B(x) \rightarrow \wtype{x:A}B(x)}((\tprd{y:B(x)}C(u(y))) \rightarrow C(\supp(x,u)))
\]
with the defining equation
\[
  f(\supp(x,u)) \defeq d(x,u,f\circ u).
\]

\subsection{Identity types}

We introduce primitive constants $c_\idsym$ and $c_{\refl{}}$. We write
$\id[A] a b$ for $c_\idsym(A,a,b)$ and $\refl a$ for $c_{\refl{}}(A,a)$, when
$a:A$ is understood:
\begin{itemize}
\item If $A : \UU_n$, $a:A$, and $b:A$ then $\id[A] a b : \UU_n$.
\item If $a:A$ then $\refl a :\id[A] a a $.
\end{itemize}
Given $a:A$, if $y:A, z:\id[A] a y \vdash C : \UU_m$ and 
$\vdash d:C[a,\refl{a}/y,z]$ then we can introduce a defined constant 
\[
  f:\tprd{y:A}{z:\id[A] a y} C
\]
with defining equation
\[
  f(a,\refl{a})\defeq d.
\]

\section{The second presentation}
\label{sec:syntax-more-formally}

In this section, there are three kinds of judgments 
\begin{mathpar}
\wfctx\Gamma
\and
\oftp\Gamma{a}{A}
\and
\jdeqtp\Gamma{a}{a'}{A}
\end{mathpar}
which we specify by providing inference rules for deriving them. A typical \define{inference rule}
\indexsee{inference rule}{rule}%
\indexdef{rule}%
has the form
\begin{equation*}
  \inferrule*[right=\textsc{Name}]
  {\mathcal{J}_1 \\ \cdots \\ \mathcal{J}_k}
  {\mathcal{J}}
\end{equation*}
It says that we may derive the \define{conclusion} $\mathcal{J}$, provided that we have
already derived the \define{hypotheses} $\mathcal{J}_1, \ldots, \mathcal{J}_k$.
(Note that, being judgments rather than types, these are not hypotheses \emph{internal} to the type theory in the sense of \autoref{sec:types-vs-sets}; they are instead hypotheses in the deductive system, i.e.\ the metatheory.)
On the
right we write the \textsc{Name} of the rule, and there may be extra side conditions that
need to be checked before the rule is applicable.

A \define{derivation}
\index{derivation}%
of a judgment is a tree constructed from such inference
rules, with the judgment at the root of the tree. For example, with the rules given below, the following is a derivation of
$\oftp{\emptyctx}{\lamu{x:\unit} x}{\unit\to\unit}$.
\begin{mathpar}
\inferrule*[right=$\Pi$-\intro]
  {\inferrule*[right=$\Vble$]
    {\inferrule*[right=\ctx-\textsc{ext}]
      {\inferrule*[right=$\unit$-\form]
        {\inferrule*[right=\ctx-\textsc{emp}]
          {\ }
          {\wfctx {\emptyctx}}}
        {\oftp{}{\unit}{\UU_0}}}
      {\wfctx {\tmtp x\unit}}}
   {\oftp{\tmtp x\unit}{x}{\unit}}}
 {\oftp{\emptyctx}{\lamu{x:\unit} x}{\unit\to\unit}}
\end{mathpar}

\subsection{Contexts}
\label{subsec:contexts}

\index{context}%
A context is a list
\begin{equation*}
  \tmtp{x_1}{A_1}, \tmtp{x_2}{A_2}, \ldots, \tmtp{x_n}{A_n}
\end{equation*}
which indicates that the distinct variables
\index{variable}%
$x_1, \ldots, x_n$ are assumed to have types $A_1, \ldots, A_n$, respectively. The list may be empty. We abbreviate contexts with the letters $\Gamma$ and $\Delta$, and we may juxtapose them to form larger contexts.

The judgment $\wfctx{\Gamma}$ formally expresses the fact that $\Gamma$ is a well-formed context, and is governed by the rules of inference
\begin{mathpar}
  \inferrule*[right=\ctx-\textsc{emp}]
  {\ }
  {\wfctx\emptyctx}
\and
  \inferrule*[right=\ctx-\textsc{ext}]
  {\oftp{\tmtp{x_1}{A_1}, \ldots, \tmtp{x_{n-1}}{A_{n-1}}}{A_n}{\UU_i}}
  {\wfctx{(\tmtp{x_1}{A_1}, \ldots, \tmtp{x_n}{A_n})}}
\end{mathpar}
with a side condition for the second rule: the variable $x_n$ must be distinct from the variables $x_1, \ldots, x_{n-1}$.
Note that the hypothesis and conclusion of $\ctx$-\textsc{ext} are judgments of different forms: the hypothesis says that in the context of variables $x_1, \ldots, x_{n-1}$, the expression $A_n$ has type $\UU_i$; while the conclusion says that the extended context $(\tmtp{x_1}{A_1}, \ldots, \tmtp{x_n}{A_n})$ is well-formed.
(It is a meta-theoretic property of the system that if $\oftp{\tmtp{x_1}{A_1}, \ldots, \tmtp{x_{n}}{A_{n}}}{b}{B}$ is derivable, then the context $(\tmtp{x_1}{A_1}, \ldots, \tmtp{x_{n}}{A_{n}})$ must be well-formed; thus $\ctx$-\textsc{ext} does not need to hypothesize well-formedness of the context to the left of $x_n$.)

\subsection{Structural rules}

\index{structural!rules|(}%
\index{rule!structural|(}%

The fact that the context holds assumptions is expressed by the rule which says that we may derive those typing judgments which are listed in the context:
\begin{mathpar}
  \inferrule*[right=$\Vble$]
  {\wfctx {(\tmtp{x_1}{A_1}, \ldots, \tmtp{x_n}{A_n})} }
  {\oftp{\tmtp{x_1}{A_1}, \ldots, \tmtp{x_n}{A_n}}{x_i}{A_i}}
\end{mathpar}
As with $\ctx$-\textsc{ext}, the hypothesis and conclusion of the rule $\Vble$ are judgments of different forms, only now they are reversed: we start with a well-formed context and derive a typing judgment.

The following important principles, called \define{substitution}
\indexdef{rule!of substitution}%
and
\define{weakening},
\indexdef{rule!of weakening}%
need not be explicitly assumed. Rather, it is possible to
show, by induction on the structure of all possible derivations, that whenever
the hypotheses of these rules are derivable, their conclusion is also
derivable.\footnote{Such rules are called \define{admissible}\indexdef{rule!admissible}\indexsee{admissible!rule}{rule, admissible}.}
For the typing judgments these principles are manifested as
\begin{mathpar}
  \inferrule*[right=$\Subst_1$]
  {\oftp\Gamma{a}{A} \\ \oftp{\Gamma,\tmtp xA,\Delta}{b}{B}}
  {\oftp{\Gamma,\Delta[a/x]}{b[a/x]}{B[a/x]}}
\and
  \inferrule*[right=$\Weak_1$]
  {\oftp\Gamma{A}{\UU_i} \\ \oftp{\Gamma,\Delta}{b}{B}}
  {\oftp{\Gamma,\tmtp xA,\Delta}{b}{B}}
\end{mathpar}
and for judgmental equalities they become
\begin{mathpar}
  \inferrule*[right=$\Subst_2$]
  {\oftp\Gamma{a}{A} \\ \jdeqtp{\Gamma,\tmtp xA,\Delta}{b}{c}{B}}
  {\jdeqtp{\Gamma,\Delta[a/x]}{b[a/x]}{c[a/x]}{B[a/x]}}
\and
  \inferrule*[right=$\Weak_2$]
  {\oftp\Gamma{A}{\UU_i} \\ \jdeqtp{\Gamma,\Delta}{b}{c}{B}}
  {\jdeqtp{\Gamma,\tmtp xA,\Delta}{b}{c}{B}}
\end{mathpar}
In addition to the judgmental equality rules given for each type former, we also
assume that judgmental equality is an equivalence relation respected by typing.
\begin{mathparpagebreakable}
  \inferrule*{\oftp\Gamma{a}{A}}{\jdeqtp\Gamma{a}{a}{A}}
\and
  \inferrule*{\jdeqtp\Gamma{a}{b}{A}}{\jdeqtp\Gamma{b}{a}{A}}
\and
  \inferrule*{\jdeqtp\Gamma{a}{b}{A} \\ \jdeqtp\Gamma{b}{c}{A}}{\jdeqtp\Gamma{a}{c}{A}}
\and
  \inferrule*{\oftp\Gamma{a}{A} \\ \jdeqtp\Gamma{A}{B}{\UU_i}}{\oftp\Gamma{a}{B}}
\and
  \inferrule*{\jdeqtp\Gamma{a}{b}{A} \\ \jdeqtp\Gamma{A}{B}{\UU_i}}{\jdeqtp\Gamma{a}{b}{B}}
\end{mathparpagebreakable}
Additionally, for all the type formers below, we assume rules stating that each constructor preserves definitional equality in each of its arguments; for instance, along with the $\Pi$-\intro\ rule, we assume the rule
\[
  \inferrule*[right=$\Pi$-\intro-eq]
  {\jdeqtp\Gamma{A}{A'}{\UU_i} \\
   \jdeqtp{\Gamma,\tmtp xA}{B}{B'}{\UU_i} \\
   \jdeqtp{\Gamma,\tmtp xA}{b}{b'}{B}}
  {\jdeqtp\Gamma{\lamu{x:A} b}{\lamu{x:A'} b'}{\tprd{x:A} B}}
\]
However, we omit these rules for brevity.

\index{rule!structural|)}%
\index{structural!rules|)}%

\subsection{Type universes}

\index{type!universe}%

We postulate an infinite hierarchy of type universes
\begin{equation*}
  \UU_0, \quad \UU_1, \quad  \UU_2, \quad \ldots
\end{equation*}
Each universe is contained in the next, and any type in $\UU_i$ is also in $\UU_{i+1}$:
\begin{mathpar}
\inferrule*[right=\UU-\textsc{intro}]
  {\wfctx \Gamma }
  {\oftp\Gamma{\UU_i}{\UU_{i+1}}}
\and
\inferrule*[right=\UU-\textsc{cumul}]
  {\oftp\Gamma{A}{\UU_i}}
  {\oftp\Gamma{A}{\UU_{i+1}}}
\end{mathpar}
We shall set up the rules of type theory in such a way that $\oftp\Gamma{a}{A}$
implies $\oftp\Gamma{A}{\UU_i}$ for some $i$. In other words, if $A$ plays the role of a type then it is in some universe. Another property of our type system is that $\jdeqtp\Gamma{a}{b}{A}$
implies $\oftp\Gamma{a}{A}$ and $\oftp\Gamma{b}{A}$.

\subsection{Dependent function types (\texorpdfstring{$\Pi$}{Π}-types)}
\label{sec:more-formal-pi}

\index{type!dependent function}%
\index{type!function}%

In \autoref{sec:function-types}, we introduced non-dependent functions $A\to B$ in
order to define a family of types as a function $\lam{x:A} B:A\to\UU_i$, which
then gives rise to a type of dependent functions $\tprd{x:A} B$. But with explicit contexts
we may replace $\lam{x:A} B:A\to\UU_i$ with the judgment
\begin{equation*}
  \oftp{\tmtp xA}{B}{\UU_i}.
\end{equation*}
Consequently, we may define dependent functions directly, without reference to non-dependent ones. This way we follow the general principle that each type former, with its constants and rules, should be introduced independently of all other type formers.
In fact, henceforth each type former is introduced systematically by:
\begin{itemize}
\item a \define{formation rule}, stating when the type former can be applied;\index{formation rule}\index{rule!formation}
\item some \define{introduction rules}, stating how to inhabit the type;\index{introduction rule}\index{rule!introduction}
\item \define{elimination rules}, or an induction principle, stating how to use an
  element of the type;
  \index{induction principle}\index{eliminator}
\item \define{computation rules}, which are judgmental equalities explaining what happens when elimination rules are applied to results of introduction rules;
  \index{computation rule}
  \indexsee{rule!computation}{computation rule}
\item optional \define{uniqueness principles}, which are judgmental equalities explaining how every element of the type is uniquely determined by the results of elimination rules applied to it.
  \index{uniqueness!principle}
  \indexsee{principle!uniqueness}{uniqueness principle}
\end{itemize}
(See also \autoref{rmk:introducing-new-concepts}.)

For the dependent function type these rules are:
\begin{mathparpagebreakable}
  \def\premise{\oftp{\Gamma}{A}{\UU_i} \and \oftp{\Gamma,\tmtp xA}{B}{\UU_i}}
  \inferrule*[right=$\Pi$-\form]
    \premise
    {\oftp\Gamma{\tprd{x:A}B}{\UU_i}}
\and
  \inferrule*[right=$\Pi$-\intro]
  {\oftp{\Gamma,\tmtp xA}{b}{B}}
  {\oftp\Gamma{\lam{x:A} b}{\tprd{x:A} B}}
\and
  \inferrule*[right=$\Pi$-\elim]
  {\oftp\Gamma{f}{\tprd{x:A} B} \\ \oftp\Gamma{a}{A}}
  {\oftp\Gamma{f(a)}{B[a/x]}}
\and
  \inferrule*[right=$\Pi$-\comp]
  {\oftp{\Gamma,\tmtp xA}{b}{B} \\ \oftp\Gamma{a}{A}}
  {\jdeqtp\Gamma{(\lam{x:A} b)(a)}{b[a/x]}{B[a/x]}}
\and
  \inferrule*[right=$\Pi$-\uniq]
  {\oftp\Gamma{f}{\tprd{x:A} B}}
  {\jdeqtp\Gamma{f}{(\lamu{x:A}f(x))}{\tprd{x:A} B}}
\end{mathparpagebreakable}

The expression $\lam{x:A} b$ binds free occurrences of $x$ in $b$, as does $\tprd{x:A} B$ for
$B$.

When $x$ does not occur freely in $B$ so that $B$ does not depend on $A$, we obtain as a
special case the ordinary function type $A\to B \defeq \tprd{x:A} B$. We take this as the \emph{definition} of $\to$.

We may abbreviate an expression $\lam{x:A} b$ as $\lamu{x:A} b$, with the understanding
that the omitted type $A$ should be filled in appropriately before type-checking.

\subsection{Dependent pair types (\texorpdfstring{$\Sigma$}{Σ}-types)}
\label{sec:more-formal-sigma}

\index{type!dependent pair}%
\index{type!product}%

In \autoref{sec:sigma-types}, we needed $\to$ and $\prdsym$ types in order to
define the introduction and elimination rules for $\smsym$; as with $\prdsym$, contexts allow us to state the rules for $\smsym$ independently:
\begin{mathparpagebreakable}
  \def\premise{\oftp{\Gamma}{A}{\UU_i} \and \oftp{\Gamma,\tmtp xA}{B}{\UU_i}}
  \inferrule*[right=$\Sigma$-\form]
    \premise
    {\oftp\Gamma{\tsm{x:A} B}{\UU_i}}
  \and
  \inferrule*[right=$\Sigma$-\intro]
    {\oftp{\Gamma, \tmtp x A}{B}{\UU_i} \\
     \oftp\Gamma{a}{A} \\ \oftp\Gamma{b}{B[a/x]}}
    {\oftp\Gamma{\tup ab}{\tsm{x:A} B}}
  \and
  \inferrule*[right=$\Sigma$-\elim]
    {\oftp{\Gamma, \tmtp z {\tsm{x:A} B}}{C}{\UU_i} \\
     \oftp{\Gamma,\tmtp x A,\tmtp y B}{g}{C[\tup x y/z]} \\
     \oftp\Gamma{p}{\tsm{x:A} B}}
    {\oftp\Gamma{\ind{\tsm{x:A} B}(z.C,x.y.g,p)}{C[p/z]}}
  \and
  \inferrule*[right=$\Sigma$-\comp]
    {\oftp{\Gamma, \tmtp z {\tsm{x:A} B}}{C}{\UU_i} \\
     \oftp{\Gamma, \tmtp x A, \tmtp y B}{g}{C[\tup x y/z]} \\\\
     \oftp\Gamma{a'}{A} \\ \oftp\Gamma{b'}{B[a'/x]}}
    {\jdeqtp\Gamma{\ind{\tsm{x:A} B}(z.C,x.y.g,\tup{a'}{b'})}{g[a',b'/x,y]}{C[\tup {a'} {b'}/z]}}
\end{mathparpagebreakable}
The expression $\tsm{x:A} B$ binds free occurrences of $x$ in $B$. Furthermore, because
$\ind{\tsm{x:A} B}$ has some arguments with free variables beyond those in $\Gamma$,
we bind (following the variable names above) $z$ in $C$, and $x$ and $y$ in $g$.
These bindings are written as $z.C$ and $x.y.g$, to indicate the names of the bound
variables.
\index{variable!bound}%
In particular, we treat $\ind{\tsm{x:A} B}$ as a primitive,
two of whose arguments contain binders; this is superficially similar to, but
different from, $\ind{\tsm{x:A} B}$ being a function that takes functions as
arguments.

When $B$ does not contain free occurrences of $x$, we obtain as a special case
the cartesian product $A \times B \defeq \tsm{x:A} B$. We take this
as the \emph{definition} of the cartesian product.

Notice that we don't postulate a judgmental uniqueness principle for $\Sigma$-types, even
though we could have; see \autoref{thm:eta-sigma} for a proof of the corresponding
propositional uniqueness principle.

\subsection{Coproduct types}

\index{type!coproduct}%

\begin{mathparpagebreakable}
  \inferrule*[right=$+$-\form]
  {\oftp\Gamma{A}{\UU_i} \\ \oftp\Gamma{B}{\UU_i}}
  {\oftp\Gamma{A+B}{\UU_i}}
\\
  \inferrule*[right=$+$-\intro${}_1$]
  {\oftp\Gamma{A}{\UU_i} \\ \oftp\Gamma{B}{\UU_i} \\\\ \oftp\Gamma{a}{A}}
  {\oftp\Gamma{\inl(a)}{A+B}}
\and
  \inferrule*[right=$+$-\intro${}_2$]
  {\oftp\Gamma{A}{\UU_i} \\ \oftp\Gamma{B}{\UU_i} \\\\ \oftp\Gamma{b}{B}}
  {\oftp\Gamma{\inr(b)}{A+B}}
\\
  \inferrule*[right=$+$-\elim]
  {\oftp{\Gamma,\tmtp z{(A+B)}}{C}{\UU_i} \\\\
   \oftp{\Gamma,\tmtp xA}{c}{C[\inl(x)/z]} \\
   \oftp{\Gamma,\tmtp yB}{d}{C[\inr(y)/z]} \\\\
   \oftp\Gamma{e}{A+B}}
  {\oftp\Gamma{\ind{A+B}(z.C,x.c,y.d,e)}{C[e/z]}}
\and
  \inferrule*[right=$+$-\comp${}_1$]
  {\oftp{\Gamma,\tmtp z{(A+B)}}{C}{\UU_i} \\
   \oftp{\Gamma,\tmtp xA}{c}{C[\inl(x)/z]} \\
   \oftp{\Gamma,\tmtp yB}{d}{C[\inr(y)/z]} \\\\
   \oftp\Gamma{a}{A}}
  {\jdeqtp\Gamma{\ind{A+B}(z.C,x.c,y.d,\inl(a))}{c[a/x]}{C[\inl(a)/z]}}
\and
  \inferrule*[right=$+$-\comp${}_2$]
  {\oftp{\Gamma,\tmtp z{(A+B)}}{C}{\UU_i} \\
   \oftp{\Gamma,\tmtp xA}{c}{C[\inl(x)/z]} \\
   \oftp{\Gamma,\tmtp yB}{d}{C[\inr(y)/z]} \\\\
   \oftp\Gamma{b}{B}}
  {\jdeqtp\Gamma{\ind{A+B}(z.C,x.c,y.d,\inr(b))}{d[b/y]}{C[\inr(b)/z]}}
\end{mathparpagebreakable}
In $\ind{A+B}$, $z$ is bound in $C$, $x$ is bound in $c$, and $y$ is bound in
$d$.

\subsection{The empty type \texorpdfstring{$\emptyt$}{0}}

\index{type!empty|(}%

\begin{mathparpagebreakable}
  \inferrule*[right=$\emptyt$-\form]
  {\wfctx\Gamma}
  {\oftp\Gamma\emptyt{\UU_i}}
\and
  \inferrule*[right=$\emptyt$-\elim]
  {\oftp{\Gamma,\tmtp x\emptyt}{C}{\UU_i} \\ \oftp\Gamma{a}{\emptyt}}
  {\oftp\Gamma{\ind{\emptyt}(x.C,a)}{C[a/x]}}
\end{mathparpagebreakable}
In $\ind{\emptyt}$, $x$ is bound in $C$. The empty type has no introduction rule and no computation rule.

\index{type!empty|)}%

\subsection{The unit type \texorpdfstring{$\unit$}{1}}

\index{type!unit|(}%

\begin{mathparpagebreakable}
  \inferrule*[right=$\unit$-\form]
  {\wfctx\Gamma}
  {\oftp\Gamma\unit{\UU_i}}
\and
  \inferrule*[right=$\unit$-\intro]
  {\wfctx\Gamma}
  {\oftp\Gamma{\ttt}{\unit}}
\and
  \inferrule*[right=$\unit$-\elim]
  {\oftp{\Gamma,\tmtp x\unit}{C}{\UU_i} \\
   \oftp{\Gamma,\tmtp y\unit}{c}{C[y/x]} \\
   \oftp\Gamma{a}{\unit}}
  {\oftp\Gamma{\ind{\unit}(x.C,y.c,a)}{C[a/x]}}
\and
  \inferrule*[right=$\unit$-\comp]
  {\oftp{\Gamma,\tmtp x\unit}{C}{\UU_i} \\
   \oftp{\Gamma,\tmtp y\unit}{c}{C[y/x]}}
  {\jdeqtp\Gamma{\ind{\unit}(x.C,y.c,\ttt)}{c[\ttt/y]}{C[\ttt/x]}}
\end{mathparpagebreakable}
In $\ind{\unit}$, $x$ is bound in $C$, and $y$ is bound in $c$.

Notice that we don't postulate a judgmental uniqueness principle for the unit
type; see \autoref{sec:finite-product-types} for a proof of the corresponding
propositional uniqueness statement.

\index{type!unit|)}%

\subsection{The natural number type}

\index{natural numbers|(}%

We give the rules for natural numbers, following \autoref{sec:inductive-types}.

\begin{mathparpagebreakable}
  \def\premise{
     \oftp{\Gamma,\tmtp x{\N}}{C}{\UU_i} \\
     \oftp\Gamma{c_0}{C[0/x]} \\
     \oftp{\Gamma,\tmtp{x}\N,\tmtp y C}{c_s}{C[\suc(x)/x]}}
  \inferrule*[right=$\N$-\form]
  {\wfctx\Gamma}
  {\oftp\Gamma{\N}{\UU_i}}
\and
  \inferrule*[right=$\N$-\intro${}_1$]
  {\wfctx\Gamma}
  {\oftp\Gamma{0}{\N}}
\and
  \inferrule*[right=$\N$-\intro${}_2$]
  {\oftp\Gamma{n}{\N}}
  {\oftp\Gamma{\suc(n)}{\N}}
\and
  \inferrule*[right=$\N$-\elim]
  {\premise \\ \oftp\Gamma{n}{\N}}
  {\oftp\Gamma{\ind{\N}(x.C,c_0,x.y.c_s,n)}{C[n/x]}}
\and
  \inferrule*[right=$\N$-\comp${}_1$]
  {\premise}
  {\jdeqtp\Gamma{\ind{\N}(x.C,c_0,x.y.c_s,0)}{c_0}{C[0/x]}}
\and
  \inferrule*[right=$\N$-\comp${}_2$]
  {\premise \\ \oftp\Gamma{n}{\N}}
  {\Gamma\vdash
    {\begin{aligned}[t]
      &\ind{\N}(x.C,c_0,x.y.c_s,\suc(n)) \\
      &\quad \jdeq c_s[n,\ind{\N}(x.C,c_0,x.y.c_s,n)/x,y] : C[\suc(n)/x]
    \end{aligned}}}
\end{mathparpagebreakable}
In $\ind{\N}$, $x$ is bound in $C$, and $x$ and $y$ are bound in $c_s$.

Other inductively defined types follow the same general scheme.

\index{natural numbers|)}%

\subsection{Identity types}

\index{type!identity|(}%

The presentation here corresponds to the (unbased) path induction principle for identity types in
\autoref{sec:identity-types}.

\begin{mathparpagebreakable}
  \inferrule*[right=$\idsym$-\form]
  {\oftp\Gamma{A}{\UU_i} \\ \oftp\Gamma{a}{A} \\ \oftp\Gamma{b}{A}}
  {\oftp\Gamma{\id[A]{a}{b}}{\UU_i}}
\and
  \inferrule*[right=$\idsym$-\intro]
  {\oftp\Gamma{A}{\UU_i} \\ \oftp\Gamma{a}{A}}
  {\oftp\Gamma{\refl a}{\id[A]aa}}
\and
  \inferrule*[right=$\idsym$-\elim]
  {\oftp{\Gamma,\tmtp xA,\tmtp yA,\tmtp p{\id[A]xy}}{C}{\UU_i} \\
   \oftp{\Gamma,\tmtp zA}{c}{C[z,z,\refl z/x,y,p]} \\
   \oftp\Gamma{a}{A} \\ \oftp\Gamma{b}{A} \\ \oftp\Gamma{p'}{\id[A]ab}}
  {\oftp\Gamma{\indidb{A}(x.y.p.C,z.c,a,b,p')}{C[a,b,p'/x,y,p]}}
\and
  \inferrule*[right=$\idsym$-\comp]
  {\oftp{\Gamma,\tmtp xA,\tmtp yA,\tmtp p{\id[A]xy}}{C}{\UU_i} \\
   \oftp{\Gamma,\tmtp zA}{c}{C[z,z,\refl z/x,y,p]} \\
   \oftp\Gamma{a}{A}}
  {\jdeqtp\Gamma{\indidb{A}(x.y.p.C,z.c,a,a,\refl a)}{c[a/z]}{C[a,a,\refl a/x,y,p]}}
\end{mathparpagebreakable}
In $\indidb{A}$, $x$, $y$, and $p$ are bound in $C$, and $z$ is bound in
$c$.

\index{type!identity|)}%

\subsection{Definitions}

\index{definition}%

Although the rules we listed so far allows us to construct everything we need directly, we
would still like to be able to use named constants, such as $\isequiv$, as a matter of
convenience. Informally, we can think of these constants simply as
abbreviations, but the situation is a bit subtler in the formalization.

For example, consider function composition, which we takes $f:A\to B$ and
$g:B\to C$ to $g\circ f:A\to C$. Somewhat unexpectedly, to make this work formally, $\circ$ must take as arguments not only $f$ and $g$, but also their types $A$, $B$, $C$:
\begin{narrowmultline*}
  {\circ} \defeq \lam{A:\UU_i}{B:\UU_i}{C:\UU_i}
  \narrowbreak
  \lam{g:B\to C}{f:A\to B}{x:A} g(f(x)).
\end{narrowmultline*}
From a practical perspective, we do not want to annotate each application of
$\circ$ with $A$, $B$ and $C$, as the are usually quite easily guessed from surrounding information. We would like to simply write $g\circ f$.
Then, strictly speaking, $g \circ f$ is not an abbreviation for $\lam{x : A} g(f(x))$,
because it involves additional \define{implicit arguments} which we want to suppress.
\index{implicit argument}

Inference of implicit arguments, typical ambiguity\index{typical ambiguity} (\autoref{sec:universes}),
ensuring that symbols are only defined once, etc., are collectively called
\define{elaboration}. \index{elaboration, in type theory}
Elaboration must take place prior to checking a derivation, and is
thus not usually presented as part of the core type theory. However, it is
essentially impossible to use any implementation of type theory which does not
perform elaboration; see \cite{Coq,norell2007towards} for further discussion.

\section{Homotopy type theory}
\label{sec:hott-features}

In this section we state the additional axioms of homotopy type theory which distinguish it from standard Martin-L\"{o}f type theory: function extensionality, the
univalence axiom, and higher inductive types. We state them in the style
of the second presentation \autoref{sec:syntax-more-formally}, although the first presentation \autoref{sec:syntax-informally} could be used just as well.

\subsection{Function extensionality and univalence}

There are two basic ways of introducing axioms which do not introduce new syntax or judgmental equalities (function extensionality and univalence are of this form):
either add a primitive constant to inhabit the axiom, or prove all theorems which depend on the axiom by hypothesizing a variable that inhabits the axiom, cf.\ \autoref{sec:axioms}.
While these are essentially equivalent, we opt for the former approach because we feel that the axioms of homotopy type theory are an essential part of the core theory.

\index{function extensionality}%
\autoref{axiom:funext} is formalized by introduction of a constant $\funext$ which
asserts that $\happly$ is an equivalence:
\begin{mathparpagebreakable}
  \inferrule*[right=$\Pi$-\textsc{ext}]
  {\oftp\Gamma{f}{\tprd{x:A} B} \\
   \oftp\Gamma{g}{\tprd{x:A} B}}
  {\oftp\Gamma{\funext(f,g)}{\isequiv(\happly_{f,g})}}
\end{mathparpagebreakable}
The definitions of $\happly$ and $\isequiv$ can be found in~\eqref{eq:happly} and
\autoref{sec:concluding-remarks}, respectively.

\index{univalence axiom}%
\autoref{axiom:univalence} is formalized in a similar fashion, too:
\begin{mathparpagebreakable}
  \inferrule*[right=$\UU_i$-\textsc{univ}]
  {\oftp\Gamma{A}{\UU_i} \\
   \oftp\Gamma{B}{\UU_i}}
  {\oftp\Gamma{\univalence(A,B)}{\isequiv(\idtoeqv_{A,B})}}
\end{mathparpagebreakable}
The definition of $\idtoeqv$ can be found in~\eqref{eq:uidtoeqv}.

\subsection{The circle}

\index{type!circle}%

Here we give an example of a basic higher inductive type; others follow the same
general scheme, albeit with elaborations.

Note that the rules below do not precisely follow the pattern of the ordinary
inductive types in \autoref{sec:syntax-more-formally}: the rules refer to the
notions of transport and functoriality of maps (\autoref{sec:functors}), and the
second computation rule is a propositional, not judgmental, equality. These
differences are discussed in \autoref{sec:dependent-paths}.

\begin{mathparpagebreakable}
  \inferrule*[right=$\Sn^1$-\form]
  {\wfctx\Gamma}
  {\oftp\Gamma{\Sn^1}{\UU_i}}
\and
  \inferrule*[right=$\Sn^1$-\intro${}_1$]
  {\wfctx\Gamma}
  {\oftp\Gamma{\base}{\Sn^1}}
\and
  \inferrule*[right=$\Sn^1$-\intro${}_2$]
  {\wfctx\Gamma}
  {\oftp\Gamma{\lloop}{\id[\Sn^1]{\base}{\base}}}
\and
  \inferrule*[right=$\Sn^1$-\elim]
  {\oftp{\Gamma,\tmtp x{\Sn^1}}{C}{\UU_i} \\
   \oftp{\Gamma}{b}{C[\base/x]} \\
   \oftp{\Gamma}{\ell}{\dpath C \lloop b b} \\
   \oftp\Gamma{p}{\Sn^1}}
  {\oftp\Gamma{\ind{\Sn^1}(x.C,b,\ell,p)}{C[p/x]}}
\and
  \inferrule*[right=$\Sn^1$-\comp${}_1$]
  {\oftp{\Gamma,\tmtp x{\Sn^1}}{C}{\UU_i} \\
   \oftp{\Gamma}{b}{C[\base/x]} \\
   \oftp{\Gamma}{\ell}{\dpath C \lloop b b}}
  {\jdeqtp\Gamma{\ind{\Sn^1}(x.C,b,\ell,\base)}{b}{C[\base/x]}}
\and
  \inferrule*[right=$\Sn^1$-\comp${}_2$]
  {\oftp{\Gamma,\tmtp x{\Sn^1}}{C}{\UU_i} \\
   \oftp{\Gamma}{b}{C[\base/x]} \\
   \oftp{\Gamma}{\ell}{\dpath C \lloop b b}}
  {\oftp\Gamma{\Sn^1\text{-}\mathsf{loopcomp}}
    {\id {\apd{(\lamu{y:\Sn^1} \ind{\Sn^1}(x.C,b,\ell,y))}{\lloop}} {\ell}}}
\end{mathparpagebreakable}
In $\ind{\Sn^1}$, $x$ is bound in $C$. The notation ${\dpath C \lloop b b}$ for dependent paths was introduced in \autoref{sec:dependent-paths}.
\index{rules of type theory|)}%

\section{Basic metatheory}
\index{metatheory|(}%

This section discusses the meta-theoretic properties of the type theory presented in 
\autoref{sec:syntax-informally}, and similar results hold for \autoref{sec:syntax-more-formally}. Figuring out which of these still hold when we add the features from \autoref{sec:hott-features} quickly leads to open questions,\index{open!problem} as discussed at the end of this section.

Recall that \autoref{sec:syntax-informally} defines the terms of type theory as
an extension of the untyped $\lambda$-calculus. The $\lambda$-calculus 
has its own notion of computation, namely the computation rule\index{computation rule!for function types}: 
\[
  (\lam{x} t)(u) \defeq t[u/x]
\]
This rule, together with the defining equations for the defined constants form
\emph{rewriting rules}\index{rewriting rule}\index{rule!rewriting} that determine reduction steps for a rewriting 
system. These steps yield a notion of computation in the sense that each rule
has a natural direction: one simplifies $(\lam{x} t)(u)$ by evaluating the
function at its argument.

Moreover, this system is \emph{confluent}\index{confluence}, that is, if $a$ simplifies in some
number of steps to both $a'$ and $a''$, there is some $b$ to which both $a'$ and
$a''$ eventually simplify. Thus we can define $t\conv u$ to mean that $t$ and
$u$ simplify to the same term.

(The situation is similar in \autoref{sec:syntax-more-formally}: Although there
we presented the computation rules as undirected equalities $\jdeq$, we can give
an operational semantics by saying that the application of an eliminator to an
introductory form simplifies to its equal, not the other way around.)

It is straightforward to show that the system in \autoref{sec:syntax-informally}
has the following properties:

\begin{thm}\label{thm:conversion-preserves-typing}
If $A : \UU$ and $A \conv A'$ then $A' : \UU$.
If $t:A$ and $t \conv t'$ then $t':A$.
\end{thm}

We say that a term is \define{normalizable}
\indexdef{term!normalizable}%
\index{normalization}%
\indexdef{normalizable term}%
(respectively, \define{strongly
normalizable})
\indexdef{term!strongly normalizable}%
\index{normalization!strong}%
\index{strong!normalization}%
if some (respectively, every), sequence of rewriting steps from the term
terminates.

\begin{thm}\label{thm:strong-normalization}
If $A : \UU$ then $A$ is strongly normalizable.
If $t:A$ then $A$ and $t$ are strongly normalizable.
\end{thm}

We say that a term is in \define{normal form}
\index{normal form}%
\index{term!normal form of}%
if it cannot be further
simplified, and that a term is \define{closed}
\index{closed!term}%
\index{term!closed}%
if no variable occurs freely in
it. A closed normal type has to be a primitive type, i.e., of the form
$c(\vec{v})$ for some primitive constant $c$ (where the list $\vec{v}$ of closed
normal terms may be omitted if empty, for instance, as with $\N$). In fact, we
can explicitly describe all normal forms:

\begin{lem}\label{lem:normal-forms}
  The terms in normal form can be described by the following syntax:
  \begin{align*}
    v & \production  k \mid \lam{x} v \mid c(\vec{v}) \mid f(\vec{v}), \\
    k &\production x \mid k(v) \mid f(\vec{v})(k),
  \end{align*}
  where $f(\vec{v})$ represents a partial application of the defined function $f$.
  In particular, a type in normal form is of the form $k$ or $c(\vec{v})$.
\end{lem}

\begin{thm}
  If $A$ is in normal form then the 
  judgment $A : \UU$ is decidable. If $A : \UU$ and $t$ is in normal form then the judgment
  $t:A$ is decidable.
\end{thm}

Logical consistency\index{consistency} (of the system in \autoref{sec:syntax-informally}) follows
immediately: if we had $a:\emptyt$ in the empty context, then by
\autoref{thm:conversion-preserves-typing,thm:strong-normalization}, $a$
simplifies to a normal term $a':\emptyt$. But by
\autoref{lem:normal-forms} no such term exists.

\begin{cor}
 The system in \autoref{sec:syntax-informally} is logically consistent.
\end{cor}

Similarly, we have the \emph{canonicity}\indexdef{canonicity} property that if $a:\N$ in the empty
context, then $a$ simplifies to a normal term $\suc^k(0)$ for some numeral $k$.

\begin{cor}
 The system in \autoref{sec:syntax-informally} has the canonicity property.
\end{cor}

Finally, if $a,A$ are in normal form, it is \emph{decidable} whether $a:A$; in
other words, because type-checking amounts to verifying the correctness of a
proof, this means we can always ``recognize a correct proof when we see one.''

\begin{cor}
The property of being a proof in the system in \autoref{sec:syntax-informally} is decidable.
\end{cor}

\mentalpause

The above results do not apply to the extended system of homotopy type
theory (i.e., the above system extended by \autoref{sec:hott-features}), since
occurrences of the univalence axiom and constructors of higher inductive types
never simplify, breaking \autoref{lem:normal-forms}. It is an open question\index{open!problem}
whether one can simplify applications of these constants in order to restore
canonicity. We also do not have a schema describing all permissible higher
inductive types, nor are we certain how to correctly formulate their rules
(e.g., whether the computation rules on higher constructors should be judgmental
equalities).

The consistency\index{consistency} of Martin-L\"{o}f type theory extended with univalence and higher
inductive types could be shown by inventing an appropriate normalization procedure, but currently
the only proofs that these systems are consistent are via semantic models---for
univalence, a model in Kan\index{Kan complex} complexes due to Voevodsky \cite{klv:ssetmodel}, and
for higher inductive types, a model due to Lumsdaine and Shulman \cite{ls:hits}.

Other metatheoretic issues, and a summary of our current results, are discussed
in greater length in the ``Constructivity'' and ``Open problems'' sections of
the introduction to this book.

\index{metatheory|)}%

\sectionNotes\label{subsec:general-remarks}


The system of rules with introduction (primitive constants) and elimination
and computation rules (defined constant) is inspired by Gentzen natural
deduction. The possibility of strengthening the elimination rule for
existential quantification was indicated in \cite{howard:pat}. The
strengthening of the axioms for disjunction appears in \cite{Martin-Lof-1972},
and for absurdity elimination and identity type in \cite{Martin-Lof-1973}. The
$W$-types were introduced in \cite{Martin-Lof-1979}. They generalize a notion
of trees introduced by \cite{Tait-1968}.
\index{Martin-L\"of}%


The generalized form of primitive recursion for natural numbers and ordinals
appear in \cite{Hilbert-1925}. This motivated G\"odel's system $T$,
\cite{Goedel-T-1958}, which was analyzed by \cite{Tait-1966}, who used,
following \cite{Goedel-T-1958}, the terminology ``definitional equality'' for
conversion: two terms are \emph{judgmentally equal} if they reduce to a
common term by means of a sequence of applications of the reduction
rules. This terminology was also used by de Bruijn \cite{deBruijn-1973} in his
presentation of \emph{AUTOMATH}.\index{AUTOMATH}

Streicher \cite[Theorem 4.13]{Streicher-1991}, explains how to give the
semantics in a contextual category of terms in normal form using a simple syntax
similar to the one we have presented.

Our second presentation comprises fairly standard presentation of
intensional Martin-L\"{o}f type theory, with some additional features needed in
homotopy type theory. Compared to a reference presentation of
\cite{hofmann:syntax-and-semantics}, the type theory of this book has a few
non-critical differences:
\begin{itemize}
\item universes \`{a} la Russell, in the sense of
\cite{martin-lof:bibliopolis}; and
\item judgmental $\eta$ and function extensionality for $\Pi$ types;
\end{itemize}
and a few features essential for homotopy type theory:
\begin{itemize}
\item the univalence axiom; and
\item higher inductive types.
\end{itemize}
As a matter of convenience, the book primarily defines functions by induction
using definition by \emph{pattern matching}.
\index{pattern matching}%
\index{definition!by pattern matching}%
It is possible to formalize the
notion of pattern matching, as done in \autoref{sec:syntax-informally}. However, the
standard type-theoretic presentation, adopted in \autoref{sec:syntax-more-formally}, is to introduce a single \emph{dependent
eliminator} for each type former, from which functions out of that type must be
defined. This approach is easier to formalize both syntactically and
semantically, as it amounts to the universal property of the type former.
The two approaches are equivalent; see \autoref{sec:pattern-matching} for a
longer discussion.

\index{type theory!formal|)}%
\index{formal!type theory|)}%


\nocite{Angiuli13}
\nocite{BauerAcceptanceVideo}

\bibliographystyle{halpha}
\phantomsection 
\addcontentsline{toc}{part}{\bibname}
\markboth{}{\textsc{Bibliography}}
{\renewcommand{\markboth}[2]{} 
\OPTbibliographyfont
\bibliography{references}}

\cleartooddpage[\thispagestyle{empty}]


\phantomsection 
\markboth{}{\textsc{Index of symbols}}
\addcontentsline{toc}{part}{Index of symbols}
\chapter*{Index of symbols}

\newcommand{\pg}[1]{p.~\pageref{#1}}
\newcommand{\symbolindex}[2]{\hbox{\makebox[0.2\textwidth][s]{#1\hfill}\hspace*{2.5em}\parbox[t]{0.65\textwidth}{#2\hfill}}\\[1.5pt]}


{\OPTindexfont 

\noindent
\symbolindex{$x \defeq a$	}{definition, \pg{defn:defeq}}
\symbolindex{$a \jdeq b$	}{judgmental equality, \pg{defn:judgmental-equality}}
\symbolindex{$a =_A b$	}{identity type, \pg{sec:identity-types}}
\symbolindex{$a = b$	}{identity type, \pg{sec:identity-types}}
\symbolindex{$x \defid b$	}{propositional equality by definition, \pg{rmk:defid}}
\symbolindex{$\idtypevar{A}(a,b)$	}{identity type, \pg{sec:identity-types}}
\symbolindex{$\dpath{P}{p}{a}{b}$	}{dependent path type, \pg{eq:dpath}}
\symbolindex{$a \neq b$	}{disequality, \pg{sec:disequality}}
\symbolindex{$\refl{x}$	}{reflexivity path at $x$, \pg{sec:identity-types}}
\symbolindex{$\opp{p}$	}{path reversal, \pg{lem:opp}}
\symbolindex{$p\ct q$	}{path concatenation, \pg{lem:concat}}
\symbolindex{$p\leftwhisker r$	}{left whiskering, \pg{thm:EckmannHilton}}
\symbolindex{$r\rightwhisker q$	}{right whiskering, \pg{thm:EckmannHilton}}
\symbolindex{$r\hct s$	}{horizontal concatenation of 2-paths, \pg{thm:EckmannHilton}}
\symbolindex{$g\circ f$	}{composite of functions, \pg{ex:composition}}
\symbolindex{$g\circ f$	}{composite of morphisms in a precategory, \pg{ct:precategory}}
\symbolindex{$\inv{f}$	}{quasi-inverse of an equivalence, \pg{thm:equiv-eqrel}}
\symbolindex{$\inv{f}$	}{inverse of an isomorphism in a precategory, \pg{ct:inv}}
\symbolindex{$\emptyt$	}{empty type, \pg{sec:coproduct-types}}
\symbolindex{$\unit$	}{unit type, \pg{sec:finite-product-types}}
\symbolindex{$\ttt$	}{canonical inhabitant of $\unit$, \pg{sec:finite-product-types}}
\symbolindex{$\bool$	}{type of booleans, \pg{sec:type-booleans}}
\symbolindex{$\btrue$, $\bfalse$	}{constructors of $\bool$, \pg{sec:type-booleans}}
\symbolindex{$\izero$, $\ione$	}{point-constructors of the interval $\interval$, \pg{sec:interval}}
\symbolindex{$\choice{}$	}{axiom of choice, \pg{eq:ac}}
\symbolindex{$\choice{\infty}$	}{``type-theoretic axiom of choice'', \pg{thm:ttac}}
\symbolindex{$\acc(a)$	}{accessibility predicate, \pg{defn:accessibility}}
\symbolindex{$P \land Q$	}{logical conjunction (``and''), \pg{defn:logical-notation}}
\symbolindex{$\apfunc{f}(p)$ or $\ap{f}{p}$	}{application of $f:A\to B$ to $p:\id[A]xy$, \pg{lem:map}}
\symbolindex{$\apd{f}{p}$	}{application of $f:\prd{a:A} B(a)$ to $p:\id[A]xy$, \pg{lem:mapdep}}
\symbolindex{$\apdtwo{f}{p}$	}{two-dimensional dependent $\apfunc{}$, \pg{thm:apd2}}
\symbolindex{$x\apart y$	}{apartness of real numbers, \pg{apart}}
\symbolindex{$\base$	}{basepoint of $\Sn^1$, \pg{sec:intro-hits}}
\symbolindex{$\base$            }{basepoint of $\Sn^2$, \pg{s2a} and \pg{s2b}}
\symbolindex{$\biinv(f)$        }{proposition that $f$ is bi-invertible, \pg{defn:biinv}}
\symbolindex{$x \bisim y$       }{bisimulation, \pg{def:bisimulation}}
\symbolindex{$\blank$           }{blank used for implicit $\lambda$-abstractions, \pg{blank}}
\symbolindex{$\CAP$             }{type of Cauchy approximations, \pg{cauchy-approximations}}
\symbolindex{$\card$            }{type of cardinal numbers, \pg{defn:card}}
\symbolindex{$\modal A$         }{reflector or modality applied to $A$, \pg{defn:reflective-subuniverse} and \pg{defn:modality}}
\symbolindex{$\cocone{X}{Y}$    }{type of cocones, \pg{defn:cocone}}
\symbolindex{$\code$            }{family of codes for paths, \pg{sec:compute-coprod}, \pg{S1-universal-cover}, \pg{sec:general-encode-decode}}
\symbolindex{$A \setminus B$    }{subset complement, \pg{complement}}
\symbolindex{$\cons(x,\ell)$    }{concatenation constructor for lists, \pg{lst} and \pg{lst-freemonoid}}
\symbolindex{$\contr_x$         }{path to the center of contraction, \pg{defn:contractible}}
\symbolindex{$\mathcal{F}\cover\pairr{J, \mathcal{G}}$          }{inductive cover, \pg{defn:inductive-cover}}
\symbolindex{$\dcut(L,U)$       }{the property of being a Dedekind cut, \pg{defn:dedekind-reals}}
\symbolindex{$\surr{L}{R}$      }{cut defining a surreal number, \pg{surreal-cut}}
\symbolindex{$\dgr{X}$          }{morphism reversal in a $\dagger$-category, \pg{sec:dagger-categories}}
\symbolindex{$\decode$          }{decoding function for paths, \pg{sec:compute-coprod}, \pg{S1-universal-cover}, \pg{sec:general-encode-decode}}
\symbolindex{$\encode$	}{encoding function for paths, \pg{sec:compute-coprod}, \pg{S1-universal-cover}, \pg{sec:general-encode-decode}}
\symbolindex{$\mreturn^\modal_A$ or $\mreturn_A$	}{the function $A\to\modal A$, \pg{defn:reflective-subuniverse} and \pg{defn:modality}}
\symbolindex{$A \epi B$	}{epimorphism or surjection}
\symbolindex{$\noeq(x,y)$	}{path-constructor of the surreals, \pg{defn:surreals}}
\symbolindex{$\rceq(u,v)$	}{path-constructor of the Cauchy reals, \pg{defn:cauchy-reals}}
\symbolindex{$a \eqr b$	}{an equivalence relation, \pg{equivalencerelation}}
\symbolindex{$\eqv{X}{Y}$	}{type of equivalences, \pg{eq:eqv}}
\symbolindex{$\texteqv{X}{Y}$	}{type of equivalences (same as $\eqv{X}{Y}$)}
\symbolindex{$\cteqv{A}{B}$	}{type of equivalences of categories, \pg{ct:equiv}}
\symbolindex{$P \Leftrightarrow Q$	}{logical equivalence, \pg{defn:logical-notation}}
\symbolindex{$\exis{x:A} B(x)$	}{logical notation for mere existential, \pg{defn:logical-notation}}
\symbolindex{$\extend f$	}{extension of $f:A\to B$ along $\eta_A$, \pg{extend}}
\symbolindex{$\bot$	}{logical falsity, \pg{defn:logical-notation}}
\symbolindex{$\hfib{f}{b}$	}{fiber of $f:A\to B$ at $b:B$, \pg{defn:homotopy-fiber}}
\symbolindex{$\Fin(n)$  }{standard finite type, \pg{fin}}
\symbolindex{$\fall{x:A} B(x)$	}{logical notation for dependent function type, \pg{defn:logical-notation}}
\symbolindex{$\funext$	}{function extensionality, \pg{axiom:funext}}
\symbolindex{$A\to B$	}{function type, \pg{sec:function-types}}
\symbolindex{$\glue$	}{path constructor of $A \sqcup^C B$, \pg{sec:colimits}}
\symbolindex{$\happly$	}{function making a path of functions into a homotopy, \pg{eq:happly}}
\symbolindex{$\hom_A(a,b)$	}{hom-set in a precategory, \pg{ct:precategory}}
\symbolindex{$f \htpy g$	}{homotopy between functions, \pg{defn:homotopy}}
\symbolindex{$\interval$	}{the interval type, \pg{sec:interval}}
\symbolindex{$\idfunc[A]$	}{the identity function of $A$, \pg{idfunc}}
\symbolindex{$1_a$	}{identity morphism in a precategory, \pg{ct:precategory}}
\symbolindex{$\idtoeqv$	}{function $(A=B)\to(\eqv A B)$ which univalence inverts, \pg{eq:uidtoeqv}}
\symbolindex{$\idtoiso$	}{function $(a=b) \to (a\cong b)$ in a precategory, \pg{ct:idtoiso}}
\symbolindex{$\im(f)$	}{image of map $f$, \pg{defn:modal-image}}
\symbolindex{$\im_n(f)$	}{$n$-image of map $f$, \pg{defn:modal-image}}
\symbolindex{$P \Rightarrow Q$	}{logical implication (``implies''), \pg{defn:logical-notation}}
\symbolindex{$a \in P$	}{membership in a subset or subtype, \pg{membership}}
\symbolindex{$x\in v$	}{membership in the cumulative hierarchy, \pg{V-membership}}
\symbolindex{$x\bin v$	}{resized membership, \pg{resized-membership}}
\symbolindex{$\ind{\emptyt}$	}{induction for ${\emptyt}$, \pg{defn:induction-emptyt},}
\symbolindex{$\ind{\unit}$	}{induction for ${\unit}$, \pg{defn:induction-unit},}
\symbolindex{$\ind{\bool}$	}{induction for ${\bool}$, \pg{defn:induction-bool},}
\symbolindex{$\ind{\nat}$	}{induction for ${\nat}$, \pg{defn:induction-nat}, and}
\symbolindex{$\indid{A}$	}{path induction for $=_A$, \pg{defn:induction-ML-id},}
\symbolindex{$\indidb{A}$	}{based path induction for $=_A$, \pg{defn:induction-PM-id},}
\symbolindex{$\ind{A \times B}$	}{induction for ${A \times B}$, \pg{defn:induction-times},}
\symbolindex{$\ind{\sm{x:A} B(x)}$	}{induction for ${\sm{x:A} B}$, \pg{defn:induction-sm},}
\symbolindex{$\ind{A + B}$	}{induction for ${A + B}$, \pg{defn:induction-plus},}
\symbolindex{$\ind{\wtype{x:A} B(x)}$	}{induction for ${\wtype{x:A} B}$, \pg{defn:induction-wtype}}
\symbolindex{$\ordsl{A}{a}$	}{initial segment of an ordinal, \pg{initial-segment}}
\symbolindex{$\inj(A,B)$	}{type of injections, \pg{inj}}
\symbolindex{$\inl$	}{first injection into a coproduct, \pg{sec:coproduct-types}}
\symbolindex{$\inr$	}{second injection into a coproduct, \pg{sec:coproduct-types}}
\symbolindex{$A \cap B$	}{intersection of subsets, \pg{intersection}, classes, \pg{class-intersection}, or intervals, \pg{interval-intersection}}
\symbolindex{$\iscontr(A)$	}{proposition that $A$ is contractible, \pg{defn:contractible}}
\symbolindex{$\isequiv(f)$	}{proposition that $f$ is an equivalence, \pg{basics-isequiv}, \pg{cha:equivalences}, and \pg{sec:concluding-remarks}}
\symbolindex{$\ishae(f)$	}{proposition that $f$ is a half-adjoint equivalence, \pg{defn:ishae}}
\symbolindex{$a\cong b$	}{type of isomorphisms in a (pre)category, \pg{ct:isomorphism}}
\symbolindex{$A\cong B$	}{type of isomorphisms between precategories, \pg{ct:isocat}}
\symbolindex{$a\unitaryiso b$	}{type of unitary isomorphisms, \pg{ct:unitary}}
\symbolindex{$\isotoid$	}{inverse of $\idtoiso$ in a category, \pg{isotoid}}
\symbolindex{$\istype{n}(X)$	}{proposition that $X$ is an $n$-type, \pg{def:hlevel}}
\symbolindex{$\isprop(A)$	}{proposition that $A$ is a mere proposition, \pg{defn:isprop}}
\symbolindex{$\isset(A)$	}{proposition that $A$ is a set, \pg{defn:set}}
\symbolindex{$A*B$	}{join of $A$ and $B$, \pg{join}}
\symbolindex{$\ker(f)$	}{kernel of a map of pointed sets, \pg{kernel}}
\symbolindex{$\lam{x} b(x)$	}{$\lambda$-abstraction, \pg{eq:lambda-abstraction}}
\symbolindex{$\lcoh{f}{g}{\eta}$	}{type of left adjoint coherence data, \pg{defn:lcoh-rcoh}}
\symbolindex{$\LEM{}$	}{law of excluded middle, \pg{eq:lem}}
\symbolindex{$\LEM{\infty}$	}{inconsistent propositions-as-types \LEM{}, \pg{thm:not-lem} and \pg{lem-infty}}
\symbolindex{$x < y$	}{strict inequality on natural numbers, \pg{leq-nat}, ordinals, \pg{sec:ordinals}, Cauchy reals, \pg{lt-RC}, surreals, \pg{defn:surreals}, etc.}
\symbolindex{$x \le y$	}{non-strict inequality on natural numbers, \pg{leq-nat}, Cauchy reals, \pg{leq-RC}, surreals, \pg{defn:surreals}, etc.}
\symbolindex{$\preceq$, $\prec$	}{recursive versions of $\le$ and $<$ for surreals, \pg{defn:No-codes}}
\symbolindex{$\ble$, $\blt$, $\bble$, $\bblt$	}{orderings on codomain of $\NO$-recursion, \pg{NO-recursion}}
\symbolindex{$\rclim(x)$	}{limit of a Cauchy approximation, \pg{defn:cauchy-reals}}
\symbolindex{$\linv(f)$	}{type of left inverses to $f$, \pg{defn:linv-rinv}}
\symbolindex{$\lst{X}$	}{type of lists of elements of $X$, \pg{lst} and \pg{lst-freemonoid}}
\symbolindex{$\lloop$	}{path-constructor of $\Sn^1$, \pg{sec:intro-hits}}
\symbolindex{$\Map_*(A,B)$	}{type of based maps, \pg{based-maps}}
\symbolindex{$x\mapsto b$	}{alternative notation for $\lambda$-abstraction, \pg{mapsto}}
\symbolindex{$\max(x,y)$	}{maximum in some ordering, e.g.\ \pg{ordered-field} and \pg{leq-RC}}
\symbolindex{$\merid(a)$	}{meridian of $\susp A$ at $a:A$, \pg{sec:suspension}}
\symbolindex{$\min(x,y)$	}{minimum in some ordering, e.g.\ \pg{ordered-field} and \pg{leq-RC}}
\symbolindex{$A \mono B$	}{monomorphism or embedding}
\symbolindex{$\N$	}{type of natural numbers, \pg{sec:inductive-types}}
\symbolindex{$\north$	}{north pole of $\susp A$, \pg{sec:suspension}}
\symbolindex{$\natw$, $\zerow$, $\sucw$	}{natural numbers encoded as a $W$-type, \pg{natw}}
\symbolindex{$\nalg$	}{type of $\nat$-algebras, \pg{defn:nalg}}
\symbolindex{$\nhom(C,D)$	}{type of $\nat$-homomorphisms, \pg{defn:nhom}}
\symbolindex{$\nil$	}{empty list, \pg{lst} and \pg{lst-freemonoid}}
\symbolindex{$\NO$	}{type of surreal numbers, \pg{defn:surreals}}
\symbolindex{$\neg P$	}{logical negation (``not''), \pg{defn:logical-notation}}
\symbolindex{$\typele{n}$, $\typeleU{n}$	}{universe of $n$-types, \pg{universe-of-ntypes}}
\symbolindex{$\Omega(A,a)$, $\Omega A$		}{loop space of a pointed type, \pg{def:loopspace}}
\symbolindex{$\Omega^k(A,a)$, $\Omega^k A$		}{iterated loop space, \pg{def:loopspace}}
\symbolindex{$A\op$	}{opposite precategory, \pg{ct:opposite-category}}
\symbolindex{$P \lor Q$	}{logical disjunction (``or''), \pg{defn:logical-notation}}
\symbolindex{$\ord$	}{type of ordinal numbers, \pg{ord}}
\symbolindex{$\pairr{a,b}$	}{(dependent) pair, \pg{sec:finite-product-types} and \pg{defn:dependent-pair}}
\symbolindex{$\pairpath$	}{constructor for $=_{A \times B}$, \pg{defn:pairpath}}
\symbolindex{$\pi_n(A)$	}{$n^{\mathrm{th}}$ homotopy group of $A$, \pg{thm:homotopy-groups} and \pg{def-of-homotopy-groups}}
\symbolindex{$\power A$	}{power set, \pg{powerset}}
\symbolindex{$\powerp A$	}{merely-inhabited power set, \pg{inhabited-powerset}}
\symbolindex{$\Zpred$	}{predecessor function $\Z\to\Z$, \pg{subsec:pi1s1-encode-decode}}
\symbolindex{$A\times B$	}{cartesian product type, \pg{sec:finite-product-types}}
\symbolindex{$\prd{x:A} B(x)$	}{dependent function type, \pg{sec:pi-types}}
\symbolindex{$\proj1(t)$	}{the first projection from a pair, \pg{defn:proj} and \pg{defn:dependent-proj1}}
\symbolindex{$\proj2(t)$	}{the second projection from a pair, \pg{defn:proj} and \pg{defn:dependent-proj1}}
\symbolindex{$\prop$, $\propU$	}{universe of mere propositions, \pg{propU}}
\symbolindex{$A \times_C B$	}{pullback of $A$ and $B$ over $C$, \pg{eq:defn-pullback}}
\symbolindex{$A \sqcup^C B$	}{pushout of $A$ and $B$ under $C$, \pg{sec:colimits}}
\symbolindex{$\Q$	}{type of rational numbers, \pg{sec:field-rati-numb}}
\symbolindex{$\Qp$	}{type of positive rational numbers, \pg{positive-rationals}}
\symbolindex{$\qinv(f)$	}{type of quasi-inverses to $f$, \pg{qinv}}
\symbolindex{$A/R$	}{quotient of a set by an equivalence relation, \pg{sec:set-quotients}}
\symbolindex{$A\sslash R$	}{alternative definition of quotient, \pg{def:VVquotient}}
\symbolindex{$\R$	}{type of real numbers (either), \pg{sec:compactness-interval}}
\symbolindex{$\RC$	}{type of Cauchy real numbers, \pg{defn:cauchy-reals}}
\symbolindex{$\RD$	}{type of Dedekind real numbers, \pg{defn:dedekind-reals}}
\symbolindex{$\rcrat(q)$	}{rational number regarded as a Cauchy real, \pg{defn:cauchy-reals}}
\symbolindex{$\rcoh{f}{g}{\epsilon}$	}{type of right adjoint coherence data, \pg{defn:lcoh-rcoh}}
\symbolindex{$\rec{\emptyt}$	}{recursor for ${\emptyt}$, \pg{defn:recursor-emptyt}}
\symbolindex{$\rec{\unit}$	}{recursor for ${\unit}$, \pg{defn:recursor-unit}}
\symbolindex{$\rec{\bool}$	}{recursor for ${\bool}$, \pg{defn:recursor-bool}}
\symbolindex{$\rec{\nat}$	}{recursor for ${\nat}$, \pg{defn:recursor-nat}}
\symbolindex{$\rec{A \times B}$	}{recursor for ${A \times B}$, \pg{defn:recursor-times}}
\symbolindex{$\rec{\sm{x:A} B(x)}$	}{recursor for ${\sm{x:A} B}$, \pg{defn:recursor-sm}}
\symbolindex{$\rec{A + B}$	}{recursor for ${A + B}$, \pg{defn:recursor-plus}}
\symbolindex{$\rec{\wtype{x:A} B(x)}$	}{recursor for ${\wtype{x:A} B}$, \pg{defn:recursor-wtype}}
\symbolindex{$\rinv$	}{type of right inverses to $f$, \pg{defn:linv-rinv}}
\symbolindex{$\south$	}{south pole of $\susp A$, \pg{sec:suspension}}
\symbolindex{$\Sn^n$	}{$n$-dimensional sphere, \pg{sec:circle}}
\symbolindex{$\seg$	}{path-constructor of the interval $\interval$, \pg{sec:interval}}
\symbolindex{$\set$, $\setU$	}{universe of sets, \pg{setU}}
\symbolindex{$\uset$	}{category of sets, \pg{ct:precatset}}
\symbolindex{$\vset(A,f)$	}{constructor of the cumulative hierarchy, \pg{defn:V}}
\symbolindex{$x\sim_\epsilon y$	}{relation of $\epsilon$-closeness for $\RC$, \pg{defn:cauchy-reals}}
\symbolindex{$x\approx_\epsilon y$	}{recursive version of $\sim_\epsilon$, \pg{defn:RC-approx}}
\symbolindex{$\bsim_\epsilon$ or $\bbsim_\epsilon$	}{closeness relations on codomain of $\RC$-recursion, \pg{RC-recursion}}
\symbolindex{$A\wedge B$	}{smash product of $A$ and $B$, \pg{smash}}
\symbolindex{$\setof{x : A | P(x)}$	}{subset type, \pg{defn:setof}}
\symbolindex{$\setof{ f(x) | P(x)}$	}{image of a subset, \pg{subset-image}}
\symbolindex{$B \subseteq C$	}{containment of subset types, \pg{subset}}
\symbolindex{$(q,r)\subseteq (s,t)$	}{inclusion of intervals, \pg{interval-subset}}
\symbolindex{$\suc$	}{successor function $\N\to\N$, \pg{sec:inductive-types}}
\symbolindex{$\Zsuc$	}{successor function $\Z\to\Z$, \pg{sec:pi1s1-initial-thoughts}}
\symbolindex{$A+B$	}{coproduct type, \pg{sec:coproduct-types}}
\symbolindex{$\sm{x:A} B(x)$	}{dependent pair type, \pg{sec:sigma-types}}
\symbolindex{$\supp(a, f)$	}{constructor for $W$-type, \pg{defn:supp}}
\symbolindex{$\surf$	}{2-path constructor of $\Sn^2$, \pg{s2a} and \pg{s2b}}
\symbolindex{$\susp A$	}{suspension of $A$, \pg{sec:suspension}}
\symbolindex{$\total{f}$	}{induced map on total spaces, \pg{defn:total-map}}
\symbolindex{$\trans{p}{u}$	}{transport of $u:P(x)$ along $p:x=y$, \pg{lem:transport}}
\symbolindex{$\transfib{P}{p}{u}$	}{transport of $u:P(x)$ along $p:x=y$, \pg{lem:transport}}
\symbolindex{$\transtwo{X}{Y}$	}{two-dimensional transport, \pg{thm:transport2}}
\symbolindex{$\transconst{X}{Y}{Z}$	}{transporting in a constant family, \pg{thm:trans-trivial}}
\symbolindex{$\trunc{n}{A}$	}{$n$-truncation of $A$, \pg{sec:truncations}}
\symbolindex{$\tproj[A]{n}{a}$, $\tproj{n}{a}$	}{image of $a:A$ in $\trunc n A$, \pg{sec:truncations}}
\symbolindex{$\brck{A}$	}{propositional truncation of $A$, \pg{subsec:prop-trunc} and \pg{sec:hittruncations}}
\symbolindex{$\bproj{a}$	}{image of $a:A$ in $\brck A$, \pg{subsec:prop-trunc} and \pg{sec:hittruncations}}
\symbolindex{$\top$	}{logical truth, \pg{defn:logical-notation}}
\symbolindex{$\nameless$	}{an unnamed object or variable}
\symbolindex{$A \cup B$	}{union of subsets, \pg{union}}
\symbolindex{$\UU$	}{universe type, \pg{sec:universes}}
\symbolindex{$\modaltype$	}{universe of modal types, \pg{eq:modaltype}}
\symbolindex{$\pointed\type$	}{universe of pointed types, \pg{def:pointedtype}}
\symbolindex{$\ua$	}{inverse to $\idtoeqv$ from univalence, \pg{ua}}
\symbolindex{$V$	}{cumulative hierarchy, \pg{defn:V}}
\symbolindex{$\walg(A,B)$	}{type of $w$-algebras, \pg{walg}}
\symbolindex{$\whom_{A,B}(C,D)$	}{type of $\w$-homomorphisms, \pg{whom}}
\symbolindex{$\wtype{x:A} B(x)$	}{$W$-type (inductive type), \pg{sec:w-types}}
\symbolindex{$A\vee B$	}{wedge of $A$ and $B$, \pg{wedge}}
\symbolindex{$\y$	}{Yoneda embedding, \pg{ct:yoneda}}
\symbolindex{$\Z$	}{type of integers, \pg{defn-Z}}}


\cleartooddpage[\thispagestyle{empty}]


\indexsee{principle}{axiom}
\indexsee{number!real}{real numbers}
\indexsee{abelian group}{group, abelian}
\indexsee{sequence!Cauchy}{Cauchy sequence}
\indexsee{adjunction}{adjoint functor}
\indexsee{higher topos}{$(\infty,1)$-topos}
\indexsee{topos!higher}{$(\infty,1)$-topos}
\indexsee{source!of a function}{domain}
\indexsee{target!of a function}{codomain}
\indexsee{type!truncation of}{truncation}
\indexsee{propositional!truncation}{truncation}
\indexsee{proof-relevant mathematics}{mathematics, proof-relevant}
\indexsee{classical!mathematics}{mathematics, classical}
\indexsee{classical!logic}{logic}
\indexsee{constructive!mathematics}{mathematics, constructive}
\indexsee{constructive!logic}{logic}
\indexsee{intuitionistic logic}{logic}
\indexsee{definition!inductive}{type, inductive}
\indexsee{inductive!definition}{type, inductive}
\indexsee{bounded!totally}{totally bounded}
\indexsee{sum!of numbers}{addition}
\indexsee{continuous map}{function, continuous}
\indexsee{function!continuity of@``continuity'' of}{``continuity''}
\indexsee{function!functoriality of@``functoriality'' of}{``functoriality''}
\indexsee{codes}{encode-decode method}
\indexsee{inequality}{order}
\indexsee{Coq@\Coq}{proof assistant}
\indexsee{Agda@\Agda}{proof assistant}
\indexsee{NuPRL@\NuPRL}{proof assistant}
\indexsee{generator!of an inductive type}{constructor}
\indexsee{groupoid!.infinity-@$\infty$-}{$\infty$-groupoid}
\indexsee{higher groupoid}{$\infty$-groupoid}
\indexsee{hierarchy!of n-types@of $n$-types}{$n$-type}
\indexsee{homotopy!n-type@$n$-type}{$n$-type}
\indexsee{homotopy!theory, classical}{classical homotopy theory}
\indexsee{homotopy!fiber}{fiber}
\indexsee{homotopy!limit}{limit of types}
\indexsee{homotopy!colimit}{colimit of types}
\indexsee{implementation}{proof assistant}
\indexsee{notation, abuse of}{abuse of notation}
\indexsee{language, abuse of}{abuse of language}
\indexsee{operator!induction}{induction principle}
\indexsee{operator!modal}{modality}
\indexsee{commutative!group}{group, abelian}
\indexsee{countable axiom of choice}{axiom of choice, countable}

\phantomsection 
\addcontentsline{toc}{part}{Index}
\markboth{}{\textsc{Index}}
\renewcommand{\markboth}[2]{}
{\OPTindexfont
\setlength{\columnsep}{\OPTindexcolumnsep}
\printindex}

\ifOPTcover
\cleartooddpage[\thispagestyle{empty}]
\pagestyle{empty}
\cleartoevenpage

\ThisLRCornerWallPaper{0.7}{\OPTbackimage}
\pagecolor{covercolor}
\color{covertext}

{
\parindent=0pt
\parskip=\baselineskip
{\OPTbacktitlefont
\textit{From the Introduction:}}
\OPTbackfont

\emph{Homotopy type theory} is a new branch of mathematics that combines aspects of several different fields in a surprising way. It is based on a recently discovered connection between \emph{homotopy theory} and \emph{type theory}.
It touches on topics as seemingly distant as the homotopy groups of spheres, the algorithms for type checking, and the definition of weak $\infty$-groupoids.

Homotopy type theory brings new ideas into the very foundation of mathematics.
On the one hand, there is Voevodsky's subtle and beautiful \emph{univalence axiom}.
The univalence axiom implies, in particular, that isomorphic structures can be identified, a principle that mathematicians have been happily using on workdays, despite its incompatibility with the ``official'' doctrines of conventional foundations.
On the other hand, we have \emph{higher inductive types}, which provide direct, logical descriptions of some of the basic spaces and constructions of homotopy theory: spheres, cylinders, truncations, localizations, etc.
Both ideas are impossible to capture directly in classical set-theoretic foundations, but when combined in homotopy type theory, they permit an entirely new kind of ``logic of homotopy types''.

This suggests a new conception of foundations of mathematics, with intrinsic homotopical content, an ``invariant'' conception of the objects of mathematics --- and convenient machine implementations, which can serve as a practical aid to the working mathematician.
This is the \emph{Univalent Foundations} program.

The present book is intended as a first systematic exposition of the basics of univalent foundations, and a collection of examples of this new style of reasoning --- but without requiring the reader to know or learn any formal logic, or to use any computer proof assistant.
We believe that univalent foundations will eventually become a viable alternative to set theory as the ``implicit foundation'' for the unformalized mathematics done by most mathematicians.

\bigskip

\begin{center}
  {\Large
  \textit{Get a free copy of the book at HomotopyTypeTheory.org.}}
\end{center}
}

\else
\fi


\end{document}